\documentclass[12pt]{book}
\title{Lectures on Lagrangian torus fibrations}
\author{Jonny Evans}
\date{}
\usepackage{amsthm,amsmath,amsfonts,amscd,amssymb}
\usepackage[T1]{fontenc}
\usepackage{makeidx}
\usepackage{tikz}
\usepackage{caption}
\usepackage{pgfplots}
\usepackage{bm}
\usepackage{scalerel}
\usepackage{fancyhdr}
\fancypagestyle{cup}{\fancyhf{}\fancyfoot[L]{\footnotesize This material has been/will be published by Cambridge University
Press as {\em Lectures on Lagrangian Torus Fibrations} by Jonny
Evans. This pre-publication version is free to view and download
for personal use only. Not for re-distribution, re-sale, or use
in derivative works.\\ {\copyright} Jonny Evans 2018--2023.}}

\usetikzlibrary{patterns}
\usetikzlibrary{decorations.pathreplacing}
\usetikzlibrary{decorations.markings}
\usetikzlibrary{decorations.pathmorphing}
\tikzset{snake it/.style={decorate, decoration=snake}}
\tikzset{->-/.style={decoration={
              markings,
              mark=at position .5 with {\arrow{>}}},postaction={decorate}}}
\tikzset{->>-/.style={decoration={
              markings,
              mark=at position .5 with {\arrow{>>}}},postaction={decorate}}}
\tikzset{->>>-/.style={decoration={
              markings,
              mark=at position .5 with {\arrow{>>>}}},postaction={decorate}}}
\tikzset{->>>>-/.style={decoration={
              markings,
              mark=at position .5 with {\arrow{>>>>}}},postaction={decorate}}}
\newcommand{\SmoothConic}[2]{
\begin{scope}[shift={(#1-4,#2-1.5)}]
\filldraw[fill=white,draw=none] (3.5,1) to[out=45,in=-45] (3.5,2) -- (4.5,2) to[out=-135,in=135] (4.5,1) -- (3.5,1);
\filldraw[fill=white] (4,2) circle [x radius=0.5,y radius=0.25];
\filldraw[fill=white,dashed] (3.5,1) arc [x radius=0.5,y radius=0.25,start angle=180,end angle=0];
\filldraw[fill=white] (3.5,1) arc [x radius=0.5,y radius=0.25,start angle=180,end angle=360];
\draw (3.5,1) to[out=45,in=-45] (3.5,2);
\draw (4.5,1) to[out=135,in=-135] (4.5,2);
\draw[dashed] (3.66,1.25) arc [x radius=0.34,y radius=0.17,start angle=180,end angle=0];
\draw (3.66,1.25) arc [x radius=0.34,y radius=0.17,start angle=180,end angle=360];
\draw[dashed] (3.66,1.75) arc [x radius=0.34,y radius=0.17,start angle=180,end angle=0];
\draw (3.66,1.75) arc [x radius=0.34,y radius=0.17,start angle=180,end angle=360];
\draw[dashed] (3.71,1.5) arc [x radius=0.29,y radius=0.145,start angle=180,end angle=0];
\draw (3.71,1.5) arc [x radius=0.29,y radius=0.145,start angle=180,end angle=360];
\end{scope}
}

\newcommand{\SingularConic}[2]{
\begin{scope}[shift={(#1,#2)}]
\filldraw[fill=white,draw=none] (-0.5,-0.5) -- (0,0.5) -- (0.5,-0.5) -- cycle;
\filldraw[fill=white,draw=none] (-0.5,0.5) -- (0,0.5) -- (0.5,0.5) -- cycle;
\filldraw[fill=white,draw=black,thick] (0,0.5) circle [x radius=0.5,y radius=0.25];
\filldraw[fill=white,dashed] (-0.5,-0.5) arc [x radius=0.5,y radius=0.25,start angle=180,end angle=0];
\filldraw[fill=white] (-0.5,-0.5) arc [x radius=0.5,y radius=0.25,start angle=180,end angle=360];
\draw (-0.5,-0.5) -- (0.5,0.5);
\draw (-0.5,0.5) -- (0.5,-0.5);
\draw[dashed] (-0.25,-0.25) arc [x radius=0.25,y radius=0.125,start angle=180,end angle=0];
\draw (-0.25,-0.25) arc [x radius=0.25,y radius=0.125,start angle=180,end angle=360];
\draw[dashed] (-0.25,0.25) arc [x radius=0.25,y radius=0.125,start angle=180,end angle=0];
\draw (-0.25,0.25) arc [x radius=0.25,y radius=0.125,start angle=180,end angle=360];
\node at (0,0) {\(\bullet\)};
\end{scope}
}

\usepackage{pinlabel}
\usepackage{parskip}
\usepackage{ragged2e}
\usepackage{float}
\usepackage{subfig}
\usepackage[utf8]{inputenc}
\newcommand{\CC}{\mathbb{C}}
\newcommand{\QQ}{\mathbb{Q}}
\newcommand{\RR}{\mathbb{R}}

\newcommand{\VV}{\mathbb{V}}

\newcommand{\ZZ}{\mathbb{Z}}
\newcommand{\II}{I}
\newcommand{\JJ}{J}
\newcommand{\cp}[1]{\mathbb{CP}^{#1}}
\newcommand{\rp}[1]{\mathbb{RP}^{#1}}
\newcommand{\OP}[1]{\mathrm{#1}}
\newcommand{\Lie}{\mathcal{L}}
\newcommand{\Orb}{\Omega}
\newcommand{\vtx}{\Gamma}
\newcommand{\Vect}{\OP{vect}}
\newcommand{\oo}{\mathfrak{o}}
\newcommand{\matr}[4]{\left(\begin{array}{cc}#1 & #2\\ #3 & #4\end{array}\right)}

\newcommand{\vect}[2]{\left(\begin{array}{c}#1\\#2\end{array}\right)}

\newenvironment{Proof}{\begin{proof}}{\end{proof}\ignorespacesafterend}
\newenvironment{Solution}{\begin{proof}[Solution]}{\end{proof}\ignorespacesafterend}
\begingroup
\makeatletter
\@for\theoremstyle:=definition,remark,plain\do{%
  \expandafter\g@addto@macro\csname th@\theoremstyle\endcsname{%
    \addtolength\thm@preskip\parskip
  }%
}
\endgroup

\renewcommand{\subset}{\subseteq}
\usepackage{graphicx}
\newtheorem{Theorem}{Theorem}[chapter]
\newtheorem{Lemma}[Theorem]{Lemma}
\newtheorem{Exercise}[Theorem]{Exercise}
\newtheorem{Corollary}[Theorem]{Corollary}

\newtheorem{Proposition}[Theorem]{Proposition}
\theoremstyle{remark}
\newtheorem{Remark}[Theorem]{Remark}
\theoremstyle{definition}
\newtheorem{Definition}[Theorem]{Definition}
\newtheorem{Example}[Theorem]{Example}
\newtheorem{Problem}[Theorem]{Problem}
\newtheorem{Recipe}[Theorem]{Recipe}
\newcommand{\lmatrix}{\left\langle\begin{matrix}}
\newcommand{\rmatrix}{\end{matrix}\right)}
\makeindex
\begin{document}
\frontmatter
\maketitle
\tableofcontents\thispagestyle{empty}\addtocontents{toc}{\protect\thispagestyle{empty}}

\chapter*{Preface}
\addcontentsline{toc}{chapter}{Preface}
\thispagestyle{cup}

This book is aimed at graduate students and researchers in
symplectic geometry. The primary message of the book is that
when a symplectic manifold \(X\) admits a Lagrangian torus
fibration \(f\colon X\to B\), the base \(B\) inherits an
integral affine structure from which we can ``read off'' a lot
of information about \(X\).

The book is based on a 10-hour lecture series I gave in 2019 for
graduate students at the London Taught Course Centre. It also
draws on sessions on toric geometry and symplectic reduction
which I taught between 2014--2017 for the Geometry Topics Course
at the London School of Geometry and Number Theory. It is
heavily expanded from both of these. It could be used as the
basis for a one-semester graduate-level course: the core content
is the foundational material in Chapters
\ref{ch:arnold_liouville}--\ref{ch:lag_fib}, the examples and
constructions in Chapters \ref{ch:glob}--\ref{ch:symp_cut}, and
the material on almost toric geometry in Chapters
\ref{ch:focusfocus}--\ref{ch:almost_toric_manifolds}. The
lecturer could then choose whether to include more about
Lagrangian submanifolds (Chapter \ref{ch:visible_Lagrangian} and
Appendix \ref{ch:tropical_lag}), or about connections to
low-dimensional topology or algebraic geometry (Chapters
\ref{ch:surgery}--\ref{ch:cusps} and Appendix
\ref{ch:markov_triples}).

There are many good books and papers which cover similar ground
to this book, including: Arnold's book \cite{Arnold} on
classical mechanics; Audin's book \cite{Audin} on torus actions;
Duistermaat's paper \cite{Duistermaat} on action-angle
coordinates; Symington's groundbreaking paper \cite{Symington}
on almost toric geometry and her follow-up paper with Leung
\cite{SymingtonLeung}; Auroux's survey \cite{Auroux} on mirror
symmetry, in which almost toric fibrations play a crucial role;
Zung's papers \cite{Zung1,Zung2} on the geometry and topology of
Lagrangian fibrations; Vianna's papers
\cite{Vianna1,Vianna2,Vianna3} on exotic tori and almost toric
geometry, and his paper with Cheung \cite{CheungVianna} on the
appearance of mutations in a variety of contexts; Mikhalkin
\cite{Mikhalkin} and Matessi's papers \cite{Matessi1,Matessi2}
on tropical Lagrangian submanifolds. Where there is common
ground, I have tried to give a different perspective.

We will not discuss {\em special} Lagrangian torus fibrations,
or much about the connection to {\em mirror symmetry}. For the
reader who is interested in this, there are many good places to
start, including Kontsevich and Soibelman's influential paper on
homological mirror symmetry and torus fibrations \cite{KS},
Gross's series of papers \cite{Gross1,Gross2,Gross3}, and much
of the early work of Joyce (see for example
\cite{JoyceSpLagSYZ}). We will also not get as far as discussing
the piecewise-smooth torus fibrations of Casta\~{n}o-Bernard and
Matessi \cite{CastanoBernardMatessi1,CastanoBernardMatessi2}, or
the far-reaching and highly technical constructions of
W.-D. Ruan \cite{Ruan1,Ruan2,Ruan3}.

Whilst reading, you will see that some lemmas are left as
exercises. This is because the proof is either (a) easy, (b)
fun, or (c) too much of a distraction from the main
narrative\footnote{In case (c), you shouldn't feel too bad if
you can't figure out the proof for yourself!}. You will find the
proofs of these in the sections called ``solutions to inline
exercises'' at the end of each chapter. There are also extensive
appendices: some to provide background and make the book more
self-contained, some to discuss in more detail matters which are
mentioned in the main text at a point where a full discussion
would distract.

Starting in Chapter 1, I will not assume you already know about
symplectic geometry and Lagrangian submanifolds (though it
wouldn't hurt). I will assume that you know:
\begin{itemize}
\item Differential forms and De Rham cohomology (and occasionally
singular homology, though only in passing).
\item Lie derivatives, though I have included an appendix (Appendix
\ref{ch:lie_ds}) which gives a high-level overview of this,
including a proof of Cartan's ``magic formulas'' for taking
Lie derivatives of differential forms.
\item Some basic notions from differential topology like
submersions, and critical or regular values.
\item The fundamental group and the theory of covering spaces.
\end{itemize}
There will probably be other things that I assume in passing,
but these are the most important ingredients. In the remainder
of the preface, I will assume familiarity with much more, so
that I can put this book in context.

Let \(X\) be a symplectic manifold. Roughly speaking, a
Lagrangian torus fibration on \(X\) is a map \(f\colon X\to B\)
with Lagrangian fibres. We usually call the target space \(B\)
the {\em base} of the fibration. We will see very early on
(Theorem \ref{thm:littlearnoldliouville} and Corollary
\ref{cor:bigarnoldliouville}) that the regular fibres must be
tori, and that we can use the symplectic structure to get a
natural local coordinate system on \(B\) whose transition maps
are integral affine transformations. Moreover, under nice
conditions, one can reconstruct \(X\) starting from this
integral affine manifold \(B\) (Theorem
\ref{thm:uniqueness}). Since the base has only half as many
dimensions as the total space, Lagrangian torus fibrations give
us a way of compressing information in a way that helps us to
visualise and understand 4- or 6-dimensional spaces using 2- or
3-dimensional integral affine geometry.

If we restrict to {\em regular} Lagrangian fibrations (with only
regular fibres) then we can only study a very restricted class
of symplectic manifolds (total spaces of torus bundles over a
flat base). For this reason, over the course of the book, we
gradually expand the class of critical points that \(f\) is
allowed to have. In Chapter \ref{ch:glob}, we introduce {\em
toric critical points}, which naturally appear in the theory of
toric varieties. This gives us a wealth of interesting examples
like \(X=\cp{n}\) where the integral affine base is simply a
polytope in \(\RR^n\), and we start to use the integral affine
geometry of this polytope to understand the symplectic geometry
of \(X\) (for example using {\em visible Lagrangian
submanifolds} in Chapter \ref{ch:visible_Lagrangian}). In Chapter
\ref{ch:symp_cut}, we introduce the {\em symplectic cut}
operation: this widens our class of examples to include
things like resolutions of singularities.

In Chapters \ref{ch:focusfocus}-\ref{ch:almost_toric_manifolds},
we allow ourselves another type of critical point: the {\em
focus-focus critical point}. This was intensively studied by San
V\~{u} Ng\d{o}c \cite{Ngoc}, who understood the asymptotic
behaviour of action coordinates as you approach a focus-focus
point; understanding V\~{u} Ng\d{o}c's calculation is the aim of
Chapter \ref{ch:focusfocus}. Margaret Symington \cite{Symington}
developed a general theory of Lagrangian torus fibrations with
at worst toric and focus-focus critical points, which she called
{\em almost toric fibrations}. In Chapter
\ref{ch:focusfocus_examples}, we find many examples, including
Milnor fibres of cyclic quotient singularities (Chapter
\ref{ch:focusfocus_examples}). In Chapter
\ref{ch:almost_toric_manifolds}, we explain Symington's
operations for modifying almost toric fibrations (nodal trades,
nodal slides, mutations).

Symington's ideas will allow us to get to our first real
highlight: the almost toric fibrations on \(\cp{2}\) discovered
by Vianna in 2013 \cite{Vianna1,Vianna2}. In these papers,
Vianna discovered infinitely many non-Hamiltonian-isotopic
Lagrangian tori in \(\cp{2}\). These tori are very hard to see
in our ``usual'' pictures of \(\cp{2}\), but become very easy to
construct and study using almost toric fibrations. We will not
develop any of the Floer theory required to distinguish these
tori, and refer the interested reader to Auroux's paper
\cite{Auroux} for an introduction, to Vianna's papers
\cite{Vianna1,Vianna2,Vianna3} for details, and
Pascaleff-Tonkonog \cite{TonkonogPascaleff} for later
developments. Instead, we content ourselves with the
construction of the tori; in general, the methods developed in
this book are useful for constructing and visualising, but not
so useful for proving constraints.

In Chapter \ref{ch:surgery}, we explain some of the most useful
surgery constructions that behave well with respect to almost
toric fibrations: non-toric blow-up, and rational
blow-up/blow-down:
\begin{itemize}
\item If you blow-up a toric variety at a toric fixed point then the
result is again toric, and the moment polytope is obtained
from the original moment polytope by truncating at the vertex
corresponding to the fixed point (see Example
\ref{exm:symp_cut_blow_up}). Non-toric blow-up allows us to
blow-up a point in the toric boundary which is not a toric
fixed-point and obtain an almost toric fibration on the
result. This operation was discovered by Zung \cite{Zung2},
and further elaborated by Symington \cite{Symington}.
\item Rational blow-up/blow-down is a family of operations which
allow us to replace a {\em chain} of symplectically embedded
spheres with a symplectically embedded rational homology
ball. The simplest example replaces a single sphere of
self-intersection \(-4\) with an open neighbourhood of the
zero-section in \(T^*\rp{2}\). This has proved useful in
low-dimensional topology for constructing {\em small exotic
4-manifolds}.
\end{itemize}
We will use both non-toric blow-up and rational blow-down to
understand Lisca's classification of symplectic fillings of lens
spaces. Again, we will give an almost toric construction of all
of Lisca's fillings, but shy away from proving the
classification, as this would require nontrivial input from
pseudoholomorphic curve theory.

Finally, in Chapter \ref{ch:cusps}, we will study integral
affine cones and see that these correspond to symplectic
manifolds with singularities modelled on elliptic and cusp
singularities. This will allow us to understand the minimal
resolutions of cusp singularities and provide us with an almost
toric fibration on a K3 surface. The pictures from this chapter
will aid the reader who is interested in reading Engel's
beautiful paper \cite{Engel} on the Looijenga cusp conjecture.

Appendices \ref{ch:sla}-\ref{ch:moser} provide some background
material on symplectic linear algebra, complex projective
geometry, cotangent bundles, and Moser isotopy, in an effort to
make the book more self-contained. Appendix
\ref{ch:toric_varieties} gives a construction of a toric variety
associated to a convex polytope with vertices at integer lattice
points, as a more algebro-geometric alternative to the
construction using symplectic cuts from Chapter
\ref{ch:symp_cut}. Appendix \ref{ch:visible_contact} discusses
the contact geometry and Reeb dynamics of hypersurfaces which
are fibred with respect to a Lagrangian torus
fibration. Appendix \ref{ch:tropical_lag} gives a brief
exposition of Mikhalkin's theory of tropical Lagrangian
submanifolds. Appendix \ref{ch:markov_triples} explains some of
the integral affine geometry behind the Diophantine Markov
equation, which underlies Vianna's constructions of almost toric
fibrations on \(\cp{2}\).

My goal in writing this book is to provide you with the tools
necessary for you to make your own investigations, to probe
hitherto unexplored regions of our most cherished and familiar
symplectic manifolds, and to bring back and show me the new
things that you find. Appendix \ref{ch:open_problems}, the final
chapter of the book, gives a few open problems as inspiration.

\section*{Acknowledgements}

The aforementioned papers by Denis Auroux, Margaret Symington,
and Renato Vianna have been enormously influential on my
thinking and geometric intuition, and this book has grown out of
my attempts to spread the appreciation of these papers in the
wider geometry community. I have also shared many formative
conversations and correspondence on this topic with people
including: Denis Auroux, Daniel Cavey, Georgios Dimitroglou
Rizell, Paul Hacking, Ailsa Keating, Jarek K\k{e}dra, Momchil
Konstantinov, Yank{\i} Lekili, Diego Matessi, Mirko Mauri, Emily
Maw, Mark McLean, Jie Min, Martin Schwingenheuer, Daniele Sepe,
Ivan Smith, Jack Smith, Tobias Sodoge, Dmitry Tonkonog,
Giancarlo Urz\'{u}a, Renato Vianna, and Chris Wendl. Thanks also
to Matt Buck, Yank{\i} Lekili, Patrick Ramsey, and the anonymous
referees for careful reading and corrections; to Leo Digiosia
for spotting a gap in an earlier attempted proof of Theorem
\ref{thm:ngoc}, and several more typos.

I would like to thank the 2014, 2015, 2016 and 2017 cohorts of
graduate students at the London School of Geometry and Number Theory,
for their insightful comments and questions during my topics sessions
on symplectic reduction and/or toric varieties from which the early
parts of these notes developed. Thanks also to the audience for my
2019 lectures at the London Taught Course Centre, whose patience and
endurance was tested by listening to this material in five 2-hour
blocks, and whose nonwithstanding cheerful engagement and repartee
helped me to improve these notes immeasurably.

\section*{Notation}

One point of confusion will be the fact that I often take
vectors to be row vectors and matrices to ``act'' from the
right. Apart from the typographical convenience of writing row
vectors versus column vectors, this is because my vectors are
usually momenta and hence naturally transform as covectors. To
remind the reader when I am doing this, I use the convention
\[\lmatrix a & b \\ c & d\rmatrix\] to emphasise that a matrix
will be acting from the right.

\vspace{1cm}

\begin{flushright}
{\em Jonny Evans

Lancaster, 2021}
\end{flushright}

\mainmatter
\part{Lagrangian torus fibrations}
\chapter{The Arnold-Liouville theorem}
\label{ch:arnold_liouville}
\thispagestyle{cup}

\section{Hamilton's equations in 2D}

Let \((p,q)\) be coordinates on \(\RR^2\) and \(H(p,q)\) be a smooth
function. A smooth path \((p(t),q(t))\) is said to satisfy Hamilton's
equations \index{Hamilton's equations}
for the Hamiltonian\index{Hamiltonian|(} \(H\) if\footnote{A dot over a variable
stands for differentiation with respect to time,
e.g.\ \(\dot{p}=\frac{dp}{dt}\).}
\begin{equation}
\label{eq:hamilton}\dot{p}=-\frac{\partial H}{\partial q}, \qquad\dot{q}=\frac{\partial H}{\partial p}.
\end{equation}
This can be used to describe the classical motion of a particle
moving on a one-dimensional line. We think of \(q(t)\) as the
position of the particle on the line at time \(t\), \(p(t)\) as
its momentum, and \(H\) as its energy. For example, if \(H(p,q)
= \frac{p^2}{2m}\) (the usual expression for kinetic energy of a
particle with mass \(m\)) then Hamilton's equations become
\[\dot{p}=0,\qquad\dot{q}=p/m,\] which are the statements that
(a) there is no force acting and (b) momentum is mass times
velocity. You can add in external (conservative) forces by
adding potential energy terms to \(H\). Observe that
\[\dot{H}=\frac{\partial H}{\partial p}\dot{p}+\frac{\partial
H}{\partial q}\dot{q}=\dot{q}\dot{p}-\dot{p}\dot{q}=0,\] so
energy is conserved.\index{conservation of energy}

From a purely mathematical point of view, Equation
\eqref{eq:hamilton} is a machine for turning the
Hamiltonian\footnote{Any function can be used as a Hamiltonian,
not only ones with physical relevance. The adjective {\em
Hamiltonian} \index{Hamiltonian|)} is just here to indicate the
way we're using the function \(H\), not that there is anything
special about \(H\).} function \(H(p,q)\) into a one-parameter
family of maps \(\phi^H_t\colon\RR^2\to\RR^2\) called the
associated {\em Hamiltonian flow}).\index{flow!Hamiltonian} The
flow is defined as follows: \[\phi^H_T(p_0,q_0)=(p(T),q(T)),\]
where \((p(t),q(t))\) is the solution to the differential
equation \eqref{eq:hamilton} with \(p(0)=p_0\) and
\(q(0)=q_0\). Conservation of \(H\) means that the flow
satisfies \(H(\phi^H_t(p,q))=H(p,q)\).

\begin{Remark}
In conclusion, given a function \(H\) we get a flow
\(\phi^H_t\) conserving \(H\). This is a simple instance of
{\em Noether's theorem}.\index{Noether's theorem} See Section
\ref{sct:noether} for a full discussion.

\end{Remark}
\begin{Example}\label{exm:rotation1}
If \(H_1=\frac{1}{2}(p^2+q^2)\) then \(\dot{p}=-q\), \(\dot{q}=p\),
so \[\vect{p(t)}{q(t)}=\matr{\cos t}{-\sin t}{\sin t}{\cos
t}\vect{p(0)}{q(0)}.\] This corresponds to a rotation of the plane
with constant angular speed. Conservation of \(H_1\) means that
points stay a fixed distance from the origin.

\end{Example}
\begin{Example}\label{exm:rotation2}
If \(H_2=\sqrt{p^2+q^2}\) then \(\dot{p}=-q/H_2\),
\(\dot{q}=p/H_2\). Since \(\dot{H}_2=0\), we can treat \(H_2\) as a
constant, so the solution is
\[\vect{p(t)}{q(t)}=\matr{\cos(t/H_2)}{-\sin(t/H_2)}{\sin(t/H_2)}{\cos(t/H_2)}\vect{p(0)}{q(0)}.\]
This flow has the same orbits (circles of radius \(H_2\)), but now
the orbit at radius \(H_2\) has period \(2\pi H_2\).

\end{Example}
\begin{figure}[htb]
\begin{center}
\begin{tikzpicture}
\draw (0,0) circle [radius = 1];
\draw (0,0) circle [radius = 2];
\node at (1,0) {\(\circ\)};
\node at (2,0) {\(\circ\)};
\node (b1) at (-1,0) {\(\bullet\)};
\node (b2) at (-2,0) {\(\bullet\)};
\node (l2) at (-2,1.25) [left] {\(\phi^{H_1}_{\pi}(2,0)\)};
\node (l1) at (-2,-1.25) [left] {\(\phi^{H_1}_{\pi}(1,0)\)};
\draw[->] (l1.east) -- (b1);
\draw[->] (l2.south) -- (b2);
\begin{scope}[shift={(7,0)}]
\draw (0,0) circle [radius = 1];
\draw (0,0) circle [radius = 2];
\node at (1,0) {\(\circ\)};
\node at (2,0) {\(\circ\)};
\node (c1) at (-1,0) {\(\bullet\)};
\node (c2) at (0,2) {\(\bullet\)};
\node (m2) at (-2,1.25) [left] {\(\phi^{H_2}_{\pi}(2,0)\)};
\node (m1) at (-2,-1.25) [left] {\(\phi^{H_2}_{\pi}(1,0)\)};
\draw[->] (m1.east) -- (c1);
\draw[->] (m2.east) -- (c2);
\end{scope}
\end{tikzpicture}
\end{center}
\caption{Above: Snapshots at \(t=\pi\) of the Hamiltonian systems in Examples \ref{exm:rotation1} (left) and \ref{exm:rotation2} (right) showing the orbits and positions of \(\phi^H_{\pi}(1,0)\) and \(\phi^H_{\pi}(2,0)\).}
\label{fig:simplest_ham_systems}
\end{figure}
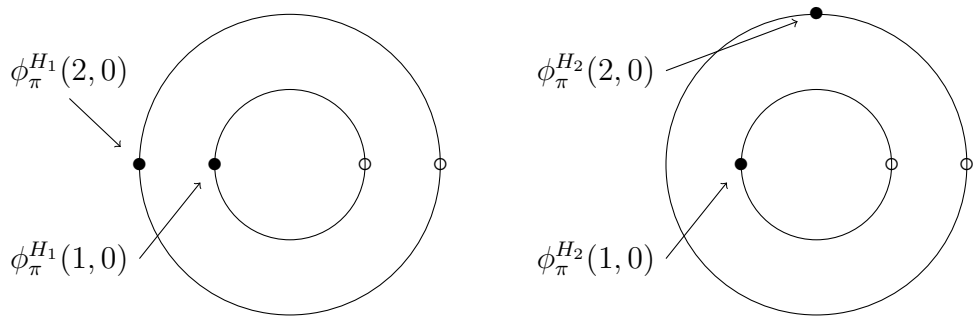

\begin{Theorem}\label{thm:2dactionangle}
If all level sets of \(H\) are closed (circular) orbits then
there exists a diffeomorphism \(\alpha\colon\RR\to\RR\) such
that, for the Hamiltonian \(\alpha\circ H\), all orbits have
period \(2\pi\).
\end{Theorem}
\begin{Proof}
By the chain rule, Hamilton's equations for \(\alpha\circ H\) are
\[\dot{p}=-\frac{\partial(\alpha\circ H)}{\partial
q}=-\alpha'(H)\frac{\partial H}{\partial
q}\quad\mbox{and}\quad \dot{q}=\frac{\partial(\alpha\circ
H)}{\partial p}=\alpha'(H)\frac{\partial H}{\partial p},\] so
the effect of postcomposing \(H\) with \(\alpha\) is to
rescale \((\dot{p},\dot{q})\) by \(\alpha'(H)\). Since \(H\)
is conserved along orbits, \(\alpha'(H)\) is constant along
orbits. This means that the orbits of the Hamiltonian flow for
\(\alpha\circ H\) are just the orbits of the flow for \(H\),
traversed at \(\alpha'(H)\) times the speed. Let
\(\Orb:=H^{-1}(b)\) be one of the orbits. If the period of the
orbit \(\Orb\) of the flow \(\phi^H_t\) is \(T(b)\) then its
period under the flow \(\phi^{\alpha\circ H}_t\) is
\(T(b)/\alpha'(b)\). To ensure that all periods are \(2\pi\),
we should therefore use
\(\alpha(b)=\frac{1}{2\pi}\int_0^bT(c)\,dc\). \qedhere

\end{Proof}
\begin{Example}
Let us revisit Example \ref{exm:rotation2}. The period of the
orbit \(H_2^{-1}(b)\) is \(2\pi b\), so the proof of Theorem
\ref{thm:2dactionangle} gives us
\(\alpha(b)=\frac{1}{2\pi}\int_0^b2\pi c\, dc=b^2/2\). This
tells us that to give all orbits the same period, we should
use the Hamiltonian \(\frac{1}{2}H_2^2\), which is precisely
the Hamiltonian \(H_1\) from Example \ref{exm:rotation1}.

\end{Example}
Periods are usually hard to find explicitly; for example, to
calculate the period of a simple pendulum in terms of its
length, initial displacement and the gravitational constant, you
need to use elliptic functions (see, for example, {\cite[Chapter
1]{Greenhill}} or {\cite[\S 44]{Whittaker}}). Similarly, the map
\(\alpha\) is difficult to write down explicitly in
examples. The following theorem gives a useful formula

\begin{Theorem}\label{thm:pdq}
In a 1-parameter family of closed orbits \(\Orb_b\),
\(b\in\RR\), of a Hamiltonian system, the period of \(\Orb_b\)
is \(\frac{d}{db}\int_{\Orb_b}p\,dq\).

\end{Theorem}
\begin{Proof}\belowdisplayskip=-12pt
Assume for simplicity\footnote{One can always find coordinates
\((p,q)\) in which the orbits have this form.} that we have
coordinates \((p,q)\), with \(q\in\RR/2\pi\ZZ\), such that the
orbits have the form \(\Orb_b:=\{(p_b(q),q)\ :\
q\in\RR/2\pi\ZZ\}\) for some functions \(p_b\).

\begin{center}
\begin{tikzpicture}
\draw[->-] (-1,-1) -- (-1,2);
\draw (1,-1) -- (1,2);
\draw[->-] (-1,0) to[out=-45,in=-135] (1,0);
\draw[dotted] (-1,0) to[out=45,in=135] (1,0);
\node at (0,-0.5) [below] {\(q\)};
\node at (-1,0.5) [left] {\(p\)};
\draw (-1,0.75) to[out=0,in=180] (-0.5,1.25) to[out=0,in=180] (0.5,0.45) to[out=0,in=180] (1,1.25);
\draw[dotted] (1,1.25) to[out=180,in=0] (-1,0.75);
\node at (1,1.25) [right] {\(\Orb_b\)};

\end{tikzpicture}
\end{center}
Then
\begin{align*}T(b)&=\int_0^{2\pi}\frac{dt}{dq}\,dq=\int_0^{2\pi}\frac{dq}{\dot{q}}\\
&=\int_0^{2\pi}\frac{dq}{\partial H/\partial
p_b}=\int_0^{2\pi}\frac{\partial p_b}{\partial
H}\,dq=\frac{d}{db}\int_0^{2\pi}p\,dq.\end{align*}\qedhere
\end{Proof}

\begin{Remark}
This means that \(\alpha(b)=\frac{1}{2\pi}\int_{\Orb_b}p\,dq\)
is another way of writing the function we found in Theorem
\ref{thm:2dactionangle}. Note that \[\alpha(b_1) - \alpha(b_0)
= \frac{1}{2\pi}\int_{\Orb_{b_1}-\Orb_{b_0}}p\,dq =
\frac{1}{2\pi}\int_C dp\wedge dq,\] by Stokes's theorem, where
\(C\) is the cylinder of orbits
\(\bigcup_{b\in[b_0,b_1]}\Orb_b\). Therefore if we choose
\(\alpha(b_0)=0\), the function \(\alpha(b)\) is just the {\em
\(dp\wedge dq\)-area} of the cylinder connecting \(\Orb_b\) to
\(\Orb_{b_0}\).

\end{Remark}
\begin{center}
\begin{tikzpicture}
\draw[->-] (-1,-1) -- (-1,2);
\draw (1,-1) -- (1,2);
\filldraw[draw=none,fill=lightgray,opacity=0.5] (-1,-0.25) to[out=0,in=180] (-0.5,-0.5) to[out=0,in=180] (0.5,-0.1) to[out=0,in=180] (1,-0.5) -- (1,1.25) to[out=180,in=0] (0.5,0.45) to[out=180,in=0] (-0.5,1.25) to[out=180,in=0] (-1,0.75) -- cycle;
\draw (-1,0.75) to[out=0,in=180] (-0.5,1.25) to[out=0,in=180] (0.5,0.45) to[out=0,in=180] (1,1.25);
\draw[dotted] (1,1.25) to[out=180,in=0] (-1,0.75);
\node at (1,1.25) [right] {\(\Orb_{b_1}\)};
\draw (-1,-0.25) to[out=0,in=180] (-0.5,-0.5) to[out=0,in=180] (0.5,-0.1) to[out=0,in=180] (1,-0.5);
\draw[dotted] (1,-0.5) to[out=180,in=0] (-1,-0.25);
\node at (1,-0.5) [right] {\(\Orb_{b_0}\)};
\node at (0,0.25) {\(C\)};

\end{tikzpicture}
\end{center}
Our goal in this first lecture is to generalise these
observations to Hamiltonian systems in higher dimensions. It
will be convenient to introduce the language of symplectic
geometry.

\section{Symplectic geometry}

This section uses Lie derivatives, Lie brackets, and the magic
formulas that relate these to exterior derivative and interior
product; we refer to Appendix \ref{ch:lie_ds} for a quick
overview of these concepts and a proof of the magic formulas.

\begin{Definition}
Let \(X\) be a manifold and \(\omega\) a 2-form. Let
\(\Vect(X)\) denote the space of vector fields on \(X\) and
\(\Omega^1(X)\) the space of 1-forms. Define a map
\(\Vect(X)\to\Omega^1(X)\) by \(V\mapsto \iota_V\omega\). We
say that \(\omega\) is {\em nondegenerate}
\index{nondegenerate 2-form} if this map is an isomorphism. A
{\em symplectic form}\index{symplectic form} is a closed,
nondegenerate 2-form.

\end{Definition}
\begin{Definition}
Let \(\omega\) be a symplectic form on a manifold
\(X\). Suppose we are given a smooth function \(H\colon
X\to\RR\). By nondegeneracy of \(\omega\), there is a unique
vector field \(V_H\) such that \(\iota_{V_H}\omega=-dH\). We
call vector fields arising in this way {\em Hamiltonian vector
fields}.\index{vector field!Hamiltonian} The flow \(\phi^H_t\)
along \(V_H\) is called a {\em Hamiltonian
flow}.\index{flow!Hamiltonian}\index{Hamiltonian flow|see
{flow, Hamiltonian}}

\end{Definition}
\begin{Example}
Let \(\omega=dp\wedge dq\) on \(X=\RR^2\) and pick a
Hamiltonian function \(H(p,q)\). Recall that the Hamiltonian
flow is defined by \((p(t),q(t))=\phi^H_t(p(0),q(0))\) and the
Hamiltonian vector field is \(V_H=(\dot{p},\dot{q})\). Using
the explicit formula for \(\omega\), we have
\(\iota_{V_H}\omega=\dot{p}\,dq-\dot{q}\,dp\). By definition,
\(\iota_{V_H}\omega=-dH=-\frac{\partial H}{\partial
p}\,dp-\frac{\partial H}{\partial q}\,dq\). Comparing
components, we recover Hamilton's equations\index{Hamilton's
equations}:
\[\dot{p}=-\frac{\partial H}{\partial q},\qquad
\dot{q}=\frac{\partial H}{\partial p}.\]

\end{Example}
\begin{Lemma}\label{lma:conserved}
A Hamiltonian flow \(\phi^H_t\) satisfies\index{conservation of energy}
\[(\phi^H_t)^*\omega=\omega\quad\mbox{and}\quad(\phi^H_t)^*H=H.\]
\end{Lemma}
\begin{Proof}
We have
\[\frac{d}{dt}((\phi^H_t)^*\omega)=(\phi^H_t)^*\Lie_{V_H}\omega
\quad\mbox{and}\quad
\frac{d}{dt}((\phi^H_t)^*H)=(\phi^H_t)^*\Lie_{V_H}H,\] so it
suffices to show that that the Lie derivatives
\(\Lie_{V_H}\omega\) and \(\Lie_{V_H}H\) vanish. For this, we
use Cartan's formula (Equation \eqref{eq:cartan_mf})
\(\Lie_V\eta=\iota_Vd\eta+d\iota_V\eta\) for the Lie
derivative of a differential form \(\eta\) along a vector
field \(V\).

We have \[\Lie_{V_H}\omega=d\iota_{V_H}\omega +
\iota_{V_H}d\omega\] Since \(d\omega=0\) the second term
vanishes. Since \(\iota_{V_H}\omega=-dH\), we get
\[\Lie_{V_H}\omega=-ddH=0.\] Finally, we have
\(\Lie_{V_H}H=\iota_{V_H}dH=-\omega(V_H,V_H)=0\), as
\(\omega\) is antisymmetric.\qedhere

\end{Proof}
\begin{Remark}
Note that if \(H\) is also allowed to depend\footnote{We call
a Hamiltonian {\em autonomous} it does not depend on \(t\) and
{\em non-autonomous} otherwise. You should imagine that if
\(H\) is autonomous then the system is just getting on by
itself, whereas if \(H\) depends on \(t\) then there is some
external input changing the system.} explicitly on \(t\) then
the previous argument for conservation of energy
(\((\phi^H_t)^*H=H\)) breaks down; an extra \(dH_t/dt\) term
appears in \(d((\phi^{H_t}_t)^*H_t)/dt\). Nonetheless, the
flow preserves the symplectic form. For example, consider the
Hamiltonian \(H_t = t\). We have \(\phi^{H_t}_t(x)=x\) for all
\(t\), which certainly preserves the symplectic form, but
energy changes over time.

\end{Remark}
\begin{Lemma}\label{lma:liebrackethamiltonian}
The Lie bracket of two Hamiltonian vector fields\index{vector
field!Hamiltonian!Lie bracket}\index{Lie bracket!of
Hamiltonian vector fields} \(V_F\) and \(V_G\) is the
Hamiltonian vector field \(V_{\{F,G\}}\), where
\(\{F,G\}=\omega(V_F,V_G)\).
\end{Lemma}
\begin{Proof}
By Equation \eqref{eq:other_mf} in Appendix \ref{ch:lie_ds},
we have
\(\iota_{[V_F,V_G]}\omega=[\Lie_{V_F},\iota_{V_G}]\omega\). Since
\(V_F\) is Hamiltonian, \(\Lie_{V_F}\omega=0\). Therefore
\[\iota_{[V_F,V_G]}\omega = \Lie_{V_F}\iota_{V_G}\omega =
d\iota_{V_F}\iota_{V_G}\omega +
\iota_{V_F}d\iota_{V_G}\omega.\] Since
\(d\iota_{V_G}\omega=-ddG=0\), we get
\(\iota_{[V_F,V_G]}\omega=d\iota_{V_F}\iota_{V_G}\omega\). Since
\(\iota_{V_F}\iota_{V_G}\omega=-\omega(V_F,V_G)\) this tells
us that \([V_F,V_G]=V_{\omega(V_F,V_G)}\) as
required. \qedhere
\end{Proof}
\begin{Definition}
The quantity \(\{F,G\}=\omega(V_F,V_G)\) is called the {\em
Poisson bracket}\index{Poisson!bracket} of \(F\) and \(G\). We
say that \(F\) and \(G\) {\em Poisson
commute}\index{Poisson!commute} if \(\{F,G\}=0\).

\end{Definition}
\begin{Remark}[Exercise \ref{exr:comm_exr}]\label{rmk:comm_exr}
Recall that the flows along two vector fields commute if and
only if the Lie bracket of the vector fields vanishes. Lemma
\ref{lma:liebrackethamiltonian} shows that two Hamiltonian
flows \(\phi^F_t\) and \(\phi^G_t\) commute if and only if the
Poisson bracket \(\{F,G\}\) is locally constant.

\end{Remark}
\begin{Lemma}[Exercise \ref{exr:heisenberg}]\label{lma:heisenberg}
Let \(F\) and \(G\) be smooth functions. Define
\(F_t(x):=F(\phi^G_t(x))\). Then \(\frac{dF_t}{dt}=\{G,F_t\}\).

\end{Lemma}
\begin{Remark}
Lemma \ref{lma:heisenberg} should look familiar to readers who know
some quantum mechanics; it is the classical counterpart of
Heisenberg's equation of motion for a quantum observable \(\hat{F}\)
evolving under the quantum Hamiltonian \(\hat{G}\).

\end{Remark}
\section{Integrable Hamiltonian systems}

\begin{Definition}[Hamiltonian \(\RR^n\)-actions]\label{dfn:hamrn}
Suppose we have a symplectic manifold \((X,\omega)\) and a map
\[\bm{H}=(H_1,\ldots,H_n)\colon X\to\RR^n\] for which the
components \(H_1,\ldots,H_n\) satisfy \(\{H_i,H_j\}=0\) for
all pairs \(i,j\). In what follows, we will assume that the
vector fields \(V_{H_i}\) can
be integrated for all time, so that the flows \(\phi^{H_i}_t\)
are defined for all \(t\in\RR\). By Remark \ref{rmk:comm_exr},
the flows \(\phi^{H_1}_{t_1},\ldots,\phi^{H_n}_{t_n}\) commute
with one another and hence define an action of the group
\(\RR^n\) on \(X\). We call this a {\em Hamiltonian
\(\RR^n\)-action}.\index{Hamiltonian!Rn
action@$\mathbb{R}^n$-action} We write
\(\phi^{\bm{H}}_{\bm{t}} :=
\phi^{H_1}_{t_1}\cdots\phi^{H_n}_{t_n}\) for this
\(\RR^n\)-action and \(\Orb(x)\) for its orbit through \(x\in
X\).

\end{Definition}
\begin{Example}[Not a Hamiltonian \(\RR^n\)-action]\label{exm:countereg}
Consider the Hamiltonians \(x\) and \(y\) on \(\RR^2\). These
generate an \(\RR^2\)-action on \(\RR^2\) where \((s,t)\) acts
by \(\phi^x_t\phi^y_s(x_0,y_0) = (x_0+s,y_0+t)\). This example
is {\em not} a Hamiltonian \(\RR^2\)-action because the
Poisson bracket \(\{x,y\}=1\) is not zero (i.e.\ the
Hamiltonians do not Poisson-commute even though the flows
commute).

\end{Example}
\begin{Remark}
More generally, for a Lie group \(G\) with Lie algebra
\(\mathfrak{g}\), a Hamiltonian
\(G\)-action\index{Hamiltonian!G-action@$G$-action} is a
\(G\)-action in which every one-parameter subgroup
\(\exp(t\xi)\), \(\xi\in\mathfrak{g}\), acts as a Hamiltonian
flow \(\phi^{H_\xi}_t\), and the assignment \(\xi\mapsto
H_\xi\) is a Lie algebra map
(i.e.\ \(H_{[\xi_1,\xi_2]}=\{H_{\xi_1},H_{\xi_2}\}\) for all
\(\xi_1,\xi_2\in\mathfrak{g}\)).

\end{Remark}
\begin{Definition}
A submanifold \(L\) of a symplectic manifold \((X,\omega)\) is
called {\em isotropic}\index{isotropic!submanifold} if
\(\omega\) vanishes on vectors tangent to \(L\) and {\em
Lagrangian}\index{Lagrangian!submanifold} if it is isotropic
and \(2\dim(L)=\dim(X)\).

\end{Definition}
\begin{Lemma}[Exercise \ref{exr:isotropic_orbits}]\label{lma:isotropic_orbits}
If \(L\) is an isotropic submanifold of the symplectic manifold
\((X,\omega)\) then \(2\dim(L)\leq \dim(X)\).

\end{Lemma}
\begin{Lemma}
Suppose that \(\bm{H}\colon X\to\RR^n\) generates a
Hamiltonian \(\RR^n\)-action. The orbits of this action are
isotropic. As a consequence, if \(X\) contains a
regular\index{regular!point} \index{regular!value}
\index{regular!fibre} point\footnote{Recall if \(\bm{H}\colon
X\to \RR^n\) is a smooth map then a point \(x\in X\) is called
{\em regular} if \(d\bm{H}\) is surjective at \(x\) and a
point \(\bm{b}\in\RR^n\) is called a {\em regular value} if
the fibre \(\bm{H}^{-1}(\bm{b})\) consists entirely of regular
points; in this case we call \(\bm{H}^{-1}(\bm{b})\) a {\em
regular fibre}.} of \(\bm{H}\) then \(n\leq \frac{1}{2}\dim
X\).
\end{Lemma}
\begin{Proof}
The tangent space to an orbit is spanned by the vectors
\(V_{H_1},\ldots,V_{H_n}\), which satisfy
\(\omega(V_{H_i},V_{H_j})=\{H_i,H_j\}=0\), so the orbits are
isotropic. If \(x\in X\) is a regular point then the
differentials \(dH_1,\ldots,dH_n\) are linearly independent at
\(x\), so the vectors \(V_{H_1}(x),\ldots,V_{H_n}(x)\) span an
\(n\)-dimensional isotropic space, which can have dimension at
most \(\frac{1}{2}\dim X\). \qedhere

\end{Proof}
\begin{Corollary}\label{cor:reg_fibres}
If \(\dim X=2n\) and \(\bm{H}\colon X\to\RR^n\) is a smooth
map with connected fibres whose components satisfy
\(\{H_i,H_j\}=0\), then the regular fibres are Lagrangian
orbits of the \(\RR^n\)-action.
\end{Corollary}
\begin{Proof}
Since \(\{H_i,H_j\}=0\), Lemma \ref{lma:heisenberg} implies
that \(H_j\) is constant along the flow of \(V_{H_i}\). In
particular, this means that if \(x\in \bm{H}^{-1}(\bm{b})\)
then its orbit \(\Orb(x)\) is contained in the fibre
\(\bm{H}^{-1}(\bm{b})\). If \(\bm{b}\) is a regular value then
the fibre \(\bm{H}^{-1}(\bm{b})\) is \(n\)-dimensional, and
the orbit of each point in the fibre is a \(n\)-dimensional
isotropic (i.e.\ Lagrangian) submanifold, so the fibre is a
union of Lagrangian submanifolds. These orbits are open
submanifolds of the fibre: if \(\Orb(x)\subset
\bm{H}^{-1}(\bm{b})\) then for any open neighbourhood
\(T\subset\RR^n\) of \(0\), the subset
\(\{\phi^{\bm{H}}_{\bm{t}}(x)\ :\ \bm{t}\in T\}\) is an open
neighbourhood of \(x\in \bm{H}^{-1}(\bm{b})\) contained in
\(\Orb(x)\). If the fibre is connected then it cannot be a
union of more than one open submanifold, so the
\(\RR^n\)-action is transitive on connected regular fibres, as
required. \qedhere

\end{Proof}
\begin{Definition}
Let \((X,\omega)\) be a \(2n\)-dimensional symplectic
manifold. We say that a smooth map \(\bm{H}\colon X\to\RR^n\)
is a {\em complete commuting Hamiltonian
system}\index{Hamiltonian!system!complete commuting} if the
components \(H_1,\ldots,H_n\) satisfy \(\{H_i,H_j\}=0\) for
all \(i,j\). We say that a complete commuting Hamiltonian
system \(\bm{H}\) is an {\em integrable Hamiltonian
system}\index{Hamiltonian!system!integrable} if
\begin{itemize}
\item \(\bm{H}(X)\) contains a dense open set of regular values,
\item \(\bm{H}\) is proper (preimages of compact sets are compact) and has
connected fibres.
\end{itemize}
The first assumption rules out trivial examples; the
properness condition ensures that the flows of the vector
fields \(V_{H_1},\ldots,V_{H_n}\) exist for all time.

\end{Definition}
\section{Period lattices}

We want to generalise the idea that all orbits have the same period,
but now we have \(n\) Hamiltonians.

\begin{figure}[htb]
\begin{center}
\begin{tikzpicture}
\draw (0,0) -- (2.5,0) node [right] {\(\sigma(B)\)};
\draw (0,0) circle [radius = 1];
\draw (0,0) circle [radius = 2];
\node at (1,0) {\(\circ\)};
\node at (2,0) {\(\circ\)};
\node (b1) at (-1,0) {\(\bullet\)};
\node (b2) at (-2,0) {\(\bullet\)};
\node (l2) at (-2,1.25) [left] {\(\phi^{H_1}_{\pi}(2,0)\)};
\node (l1) at (-2,-1.25) [left] {\(\phi^{H_1}_{\pi}(1,0)\)};
\draw[->] (l1.east) -- (b1);
\draw[->] (l2.south) -- (b2);
\begin{scope}[shift={(7,0)}]
\draw (0,0) -- (2.5,0) node [right] {\(\sigma(B)\)};
\draw (0,0) circle [radius = 1];
\draw (0,0) circle [radius = 2];
\node at (1,0) {\(\circ\)};
\node at (2,0) {\(\circ\)};
\node (c1) at (-1,0) {\(\bullet\)};
\node (c2) at (0,2) {\(\bullet\)};
\node (m2) at (-2,1.25) [left] {\(\phi^{H_2}_{\pi}(2,0)\)};
\node (m1) at (-2,-1.25) [left] {\(\phi^{H_2}_{\pi}(1,0)\)};
\draw[->] (m1.east) -- (c1);
\draw[->] (m2.east) -- (c2);
\end{scope}
\begin{scope}[shift={(0,-4.5)}]
\draw[->-] (0,0) -- (3,0) node [midway,above] {\(b\)};
\draw[->] (0,-1.5) -- (0,1.5) node [midway,above=15pt,sloped] {period of \(\Orb(\sigma(b))\)};
\draw (0,1) -- (3,1);
\draw (0,-1) -- (3,-1);
\node at (0,0) [left] {\(0\)};
\node at (0,1) [left] {\(2\pi\)};
\node at (0,-1) [left] {\(-2\pi\)};
\end{scope}
\begin{scope}[shift={(7,-4.5)}]
\draw[->-] (0,0) -- (3,0) node [midway,above] {\(b\)};
\draw[->] (0,-1.5) -- (0,1.5) node [midway,above,sloped] {period of \(\Orb(\sigma(b))\)};
\draw (0,0) -- (3,1.5) node [right] {\(2\pi b\)};
\draw (0,0) -- (1.5,1.5) node [right] {\(4\pi b\)};
\draw (0,0) -- (1,1.5) node [above] {\(6\pi b\)};
\draw (0,0) -- (3,-1.5) node [right] {\(2\pi b\)};
\draw (0,0) -- (1.5,-1.5) node [right] {\(-4\pi b\)};
\draw (0,0) -- (1,-1.5) node [below] {\(-6\pi b\)};
\node at (3,0) [right] {\(0\)};
\end{scope}
\end{tikzpicture}
\end{center}
\caption{Above: Snapshots at \(t=\pi\) of the Hamiltonian systems in Examples \ref{exm:rotation1_rev} (left) and \ref{exm:rotation2_rev} (right) showing the orbits and positions of \(\phi^H_{\pi}(1,0)\) and \(\phi^H_{\pi}(2,0)\). Below: The period lattices from Example \ref{exm:rotation1_rev} (left: standard) and Example \ref{exm:rotation2_rev} (right: non-standard).}
\label{fig:periods_rotations}
\end{figure}
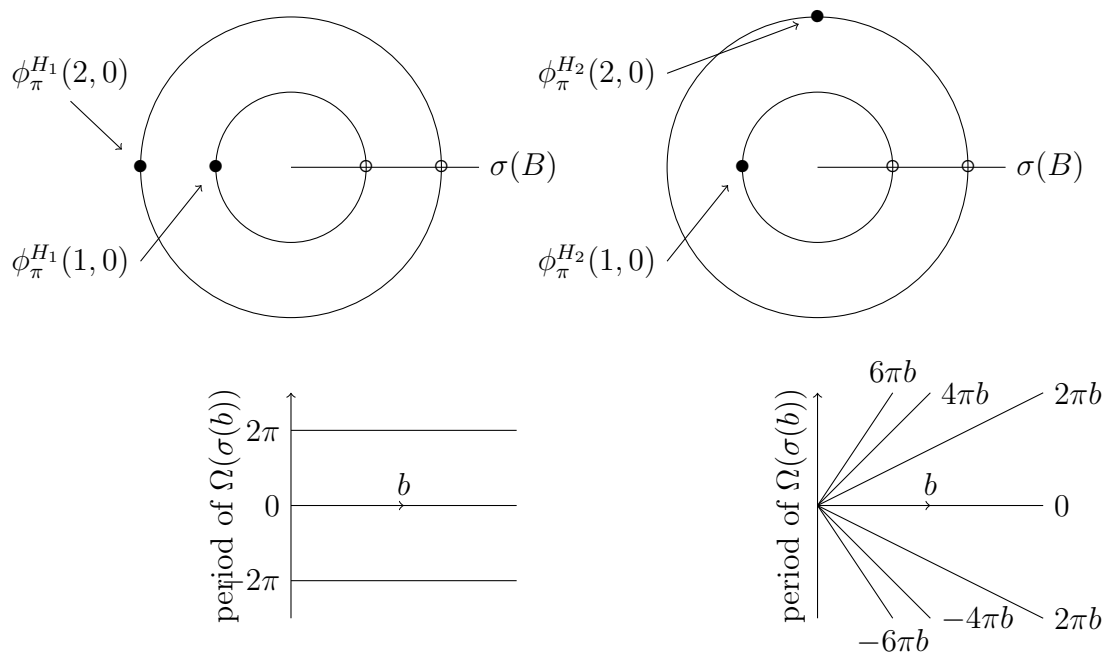

\begin{Definition}
Suppose we have an integrable Hamiltonian system
\(\bm{H}\colon X\to\RR^n\). Let
\(B\subset\bm{H}(X)\subset\RR^n\) be an open subset of the
image of \(\bm{H}\). A {\em local section over \(B\)} is a map
\(\sigma\colon B\to X\) such that
\(\bm{H}\circ\sigma=\OP{id}\).

\end{Definition}
\begin{Remark}
Note that if \(\sigma\) is a local section over \(B\) then
\(\sigma(\bm{b})\) is necessarily a regular point of
\(\bm{H}\) for every \(\bm{b}\in B\) because
\(d\bm{H}(d\sigma(T_{\bm{b}}B)) = \OP{id}(T_{\bm{b}}B) =
T_{\bm{b}}B\).

\end{Remark}
\begin{Definition}
Given an integrable Hamiltonian system \(\bm{H}\colon
X\to\RR^n\) and a local section \(\sigma\colon B\to X\), over
a subset \(B\subset \bm{H}(X)\), the {\em period lattice at
\(\bm{b}\in B\)} is defined to be:
\[\Lambda^{\bm{H}}_{\bm{b}}:=\{\bm{t}\in
\RR^n\,:\,\phi^{\bm{H}}_{\bm{t}}(\sigma(\bm{b}))=\sigma(\bm{b})\},\]
and the {\em period lattice}\index{lattice, period|see {period
lattice}}\index{standard period lattice|see {period lattice,
standard}} \index{period lattice} is \[\Lambda^{\bm{H}} :=
\{(\bm{b}, \bm{t})\in
B\times\RR^n\,:\,\bm{t}\in\Lambda^{\bm{H}}_{\bm{b}}\}.\] We
will often omit the superscript \(\bm{H}\) if \(\bm{H}\) is
clear from the context. We say that the period lattice is
standard\index{period lattice!standard} if \(\Lambda = B\times
(2\pi\ZZ)^n\).

\end{Definition}
\begin{Lemma}
\(\Lambda^{\bm{H}}_{\bm{b}}\) consists of tuples
\(\bm{t}\in\RR^n\) such that \(\phi^{\bm{H}}_{\bm{t}}\) fixes
every point of
the orbit \(\Orb(\sigma(\bm{b}))\).
\end{Lemma}
\begin{Proof}
By definition, \(\bm{t}\in\Lambda^{\bm{H}}_{\bm{b}}\) if and
only if \(\phi^{\bm{H}}_{\bm{t}}\) fixes
\(\sigma(\bm{b})\). Any other point in this orbit can be
written as \(\phi^{\bm{H}}_{\bm{t}'}(\sigma(\bm{b}))\) for
some \(\bm{t}'\). Therefore if
\(\bm{t}\in\Lambda^{\bm{H}}_{\bm{b}}\), we have
\[\phi^{\bm{H}}_{\bm{t}}(\phi^{\bm{H}}_{\bm{t}'}(\sigma(\bm{b})))
=
\phi^{\bm{H}}_{\bm{t}'}(\phi^{\bm{H}}_{\bm{t}}(\sigma(\bm{b})))
= \phi^{\bm{H}}_{\bm{t}'}(\sigma(\bm{b})),\] so
\(\phi^{\bm{H}}_{\bm{t}}\) fixes every point in the
orbit.\qedhere

\end{Proof}
\begin{Remark}
If the orbit \(\Orb(\sigma(\bm{b}))\) is dense in
\(\bm{H}^{-1}(\bm{b})\), this means
\[\Lambda^{\bm{H}}_{\bm{b}} =
\{{\bm{t}}\in\RR^n\,:\,\phi^{\bm{H}}_{\bm{t}}|_{\bm{H}^{-1}(\bm{b})}
= \OP{id}_{{\bm{H}}^{-1}(\bm{b})}\}.\]

\end{Remark}
\begin{Example}\label{exm:rotation1_rev}
In Example \ref{exm:rotation1}, the Hamiltonian is
\(H_1(p,q)=\frac{1}{2}(p^2+q^2)\) on \(\RR^2\). If we take
\(B=\RR_{>0}\) and choose the section
\(\sigma(b)=(\sqrt{2b},0)\) then \(\phi^{H_1}_t(\sigma(b)) =
(\sqrt{2b} \cos t, \sqrt{2b} \sin t)\) and the period lattice
is standard\index{period lattice!standard}: every point
\(\sigma(b)\) returns to itself after time \(2\pi\). See
Figure \ref{fig:periods_rotations} (left).

\end{Example}
\begin{Example}\label{exm:rotation2_rev}
In Example \ref{exm:rotation2}, the Hamiltonian is
\(H_2(p,q)=\sqrt{p^2+q^2}\) on \(\RR^2\). If we take the
section \(\sigma(b)=(b,0)\) then we have
\(\phi^{H_2}_t(\sigma(b))=(b\cos(t/b), b\sin(t/b))\) so the
point \(\sigma(b)\) returns to itself after time \(2\pi
b\). The period lattice is therefore \(\{(b,2\pi bn)\ :\ b>0,\
n\in\ZZ\}\). See Figure \ref{fig:periods_rotations} (right).

\end{Example}

\begin{Example}\label{exm:periodlattice}
Consider a Hamiltonian system on \(\RR^2\) whose level sets
are shown in Figure \ref{fig:periodlattice}. This Hamiltonian
generates an \(\RR\)-action which has three types of orbits:
the fixed points (marked \(\bullet\) in the figure); the two
separatrices (arcs connecting the central fixed point to
itself); the remaining orbits are closed loops either inside
or outside the separatrices. The separatrices have infinite
period (it takes infinitely long to flow around them). If we
take as Lagrangian section the wiggly line segment on the left
then the period lattice looks like the figure on the right.

\end{Example}
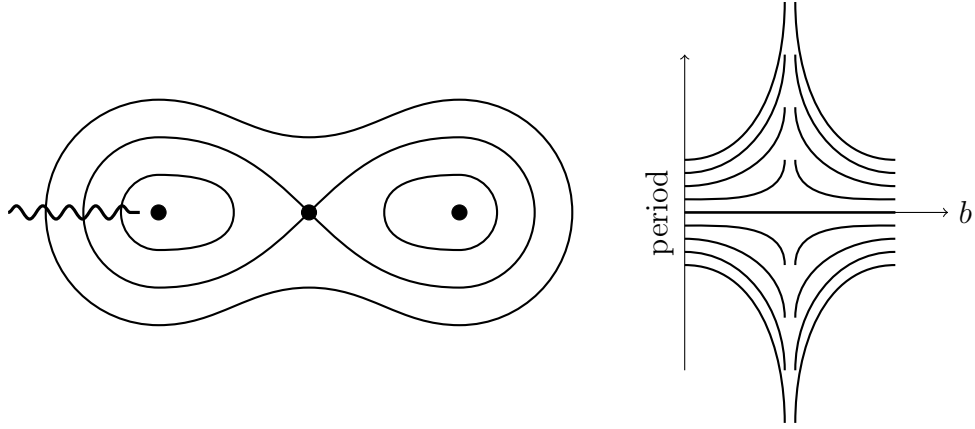
\begin{figure}[htb]
\begin{center}
\begin{tikzpicture}[baseline=0]
\draw[thick] (0,0) to[out=45,in=180] (2,1) to[out=0,in=90] (3,0) to[out=-90,in=0] (2,-1) to[out=180,in=-45] (0,0);
\draw[thick] (0,0) to[out=135,in=0] (-2,1) to[out=180,in=90] (-3,0) to[out=-90,in=180] (-2,-1) to[out=0,in=225] (0,0);
\filldraw (-2,0) circle [radius=1mm];
\filldraw (2,0) circle [radius=1mm];
\filldraw (0,0) circle [radius=1mm];
\draw[thick] (1,0) to[out=90,in=180] (2,0.5) to[out=0,in=90] (2.5,0) to[out=-90,in=0] (2,-0.5) to[out=180,in=-90] (1,0);
\draw[thick] (-1,0) to[out=90,in=0] (-2,0.5) to[out=180,in=90] (-2.5,0) to[out=-90,in=180] (-2,-0.5) to[out=0,in=-90] (-1,0);
\draw[thick] (3.5,0) to[out=90,in=0] (2,1.5) to[out=180,in=0] (0,1) to[out=180,in=0] (-2,1.5) to[out=180,in=90] (-3.5,0)
to[out=-90,in=180] (-2,-1.5) to[out=0,in=180] (0,-1) to[out=0,in=180] (2,-1.5) to[out=0,in=-90] (3.5,0);
\draw[snake it,very thick] (-4,0) -- (-2.25,0);
\end{tikzpicture}\qquad
\begin{tikzpicture}[scale=0.7,baseline=0]
\draw[thick] (-2,0) -- (2,0);
\draw[thick] (-2,0.5) to[out=0,in=-90] (-0.1,2);
\draw[thick] (2,0.5) to[out=180,in=-90] (0.1,2);
\draw[thick] (-2,0.75) to[out=0,in=-90] (-0.1,3);
\draw[thick] (2,0.75) to[out=180,in=-90] (0.1,3);
\draw[thick] (-2,0.25) to[out=0,in=-90] (-0.1,1);
\draw[thick] (2,0.25) to[out=180,in=-90] (0.1,1);
\draw[thick] (-2,1) to[out=0,in=-90] (-0.1,4);
\draw[thick] (2,1) to[out=180,in=-90] (0.1,4);
\draw[->] (-2,-3) -- (-2,3) node [midway,above,sloped] {period};
\draw[->] (-2,0) -- (3,0) node [right] {\(b\)};
\begin{scope}[yscale=-1,xscale=1]
\draw[thick] (-2,0) -- (2,0);
\draw[thick] (-2,0.5) to[out=0,in=-90] (-0.1,2);
\draw[thick] (2,0.5) to[out=180,in=-90] (0.1,2);
\draw[thick] (-2,0.75) to[out=0,in=-90] (-0.1,3);
\draw[thick] (2,0.75) to[out=180,in=-90] (0.1,3);
\draw[thick] (-2,0.25) to[out=0,in=-90] (-0.1,1);
\draw[thick] (2,0.25) to[out=180,in=-90] (0.1,1);
\draw[thick] (-2,1) to[out=0,in=-90] (-0.1,4);
\draw[thick] (2,1) to[out=180,in=-90] (0.1,4);
\end{scope}
\end{tikzpicture}
\end{center}
\caption{The Hamiltonian system (left) and period lattice (right) for Example \ref{exm:periodlattice}. The wiggly line is a Lagrangian section. The infinite period of the separatrix is what gives rise to the vertical asymptotes of the period lattice.}
\label{fig:periodlattice}
\end{figure}

The justification for ``lattice'' in the name period lattice
comes from the following result:

\begin{Lemma}[Exercise \ref{exr:discrete}]\label{lma:discrete}
For each \(\bm{b}\in B\), the intersection
\(\Lambda_{\bm{b}}=\Lambda\cap(\{\bm{b}\}\times\RR^n)\) is a
lattice\index{lattice|(} in \(\RR^n\), that is a discrete
subgroup of \(\RR^n\). The rank of the lattice is lower
semicontinuous as a function of \(\bm{b}\), that is,
\(\bm{b}\) has a neighbourhood \(V\) such that
\(rank(\Lambda_{\bm{b}'})\geq rank(\Lambda_{\bm{b}})\) for all
\(\bm{b}'\in V\).

\end{Lemma}
\begin{Example}
In Example \ref{exm:periodlattice}, the period lattice for
most orbits is isomorphic to \(\ZZ\), but where \(\sigma(B)\)
intersects the separatrix orbit the period lattice is the zero
lattice; this corresponds to the vertical asymptote in Figure
\ref{fig:periodlattice}.

\end{Example}
The following result can be found in Arnold's book {\cite[Lemma
3, p.276]{Arnold}}, and tells us that lattices are what we think
they are. We will use it below to explain why compact orbits are
diffeomorphic to tori.

\begin{Lemma}\label{lma:latticeclassification}
If \(\Lambda\subset\RR^n\) is a lattice\index{lattice|)} then
there is a basis \(e_1,\ldots,e_n\) of \(\RR^n\) such that
\(\Lambda\) is the \(\ZZ\)-linear span of the vectors
\(e_1,\ldots,e_k\) for some \(k\leq n\).

\end{Lemma}
\section{Liouville coordinates}

In what follows, we will usually use {\em Lagrangian sections}
to define the period lattice, i.e.\, sections whose image is a
Lagrangian submanifold.\index{Lagrangian!section} These always
exist locally:

\begin{Lemma}[Exercise \ref{exr:locallagrangian}]\label{lma:locallagrangian}
Let \(\bm{H}\colon X\to\RR^n\) be an integrable Hamiltonian
system. There exists a local {\em Lagrangian} section through
any regular point \(x\).

\end{Lemma}
\begin{Theorem}[Liouville coordinates]\label{thm:liouvillecoordinates}
Let \(\bm{H}\colon X\to\RR^n\) be an integrable Hamiltonian
system, let \(B\subset\RR^n\) be an open set, and let
\(\sigma\colon B\to X\) be a local Lagrangian section. Define
\[\Psi\colon B\times\RR^n\to X,\qquad\Psi(\bm{b},\bm{t}) =
\phi^{\bm{H}}_{\bm{t}}(\sigma(\bm{b})).\] Then \(\Psi\) is
both an immersion and a submersion and \(\Psi^*\omega=\sum
db_i\wedge dt_i\), where \((b_1,\ldots,b_n)\) are the standard
coordinates on \(B\subset\RR^n\). This means that
\((b_1,\ldots,b_n,t_1,\ldots,t_n)\) provide local symplectic
coordinates on a neighbourhood of \(\sigma(B)\); we call these
{\em Liouville coordinates}.\index{Liouville coordinates}
\end{Theorem}
\begin{Proof}
We first verify that \(\Psi^*\omega=\sum_{i=1}^n db_i\wedge
dt_i\) on pairs of basis vectors \(\partial_{b_i}\) and
\(\partial_{t_i}\). First, observe that, by definition of
\(\Psi\), we have \[\Psi_*\partial_{b_i} =
(\phi^{\bm{H}}_{\bm{t}})_*\sigma_*(\partial_{b_i}),\qquad
\Psi_*\partial_{t_i}V_{H_i}.\]

The vectors \(\Psi_*\partial_{b_i}\) and
\(\Psi_*\partial_{b_j}\) are tangent to
\(\phi^{\bm{H}}_{\bm{t}}(\sigma(B))\), which is the image of a
Lagrangian submanifold under a series of Hamiltonian flows,
hence Lagrangian. Therefore
\(\omega(\Psi_*\partial_{b_i},\Psi_*\partial_{b_j})=0\).

Since \(\Psi_*\partial_{t_i}=V_{H_i}\), we have
\(\omega(\Psi_*\partial_{t_i},\Psi_*\partial_{t_j}) =
\omega(V_{H_i},V_{H_j}) = \{H_i,H_j\} = 0\).

Finally, we have
\(\omega(\Psi_*\partial_{b_i},\Psi_*\partial_{t_j}) = -
(\iota_{V_{H_j}}\omega)(\Psi_*\partial_{b_i})
dH_j(\Psi_*\partial_{b_i})\). Since the flow along
\(\phi^{\bm{H}}_{\bm{t}}\) preserves the level sets of
\(H_j\), we have
\((H_j\circ\Psi)(\bm{b},\bm{t})=b_j\). Therefore
\(dH_j(\Psi_*(\partial_{b_i})) = db_j(\partial_{b_i}) =
\delta_{ij}\). This completes the verification that
\(\Psi^*\omega = \sum db_i\wedge dt_i\).

This implies that \(\Psi\) is both an immersion and a
submersion: if this failed at some point then \(\Psi^*\omega\)
would be degenerate there.\qedhere

\end{Proof}
\begin{Remark}
Note that the period lattice is given by
\(\Lambda^{\bm{H}}=\Psi^{-1}(\sigma(B))\). Since \(\Psi\) is a
symplectic map and \(\sigma\) is a Lagrangian section, the
period lattice is a Lagrangian submanifold of \(B\times
\RR^n\) with respect to \(\sum db_i\wedge dt_i\).

\end{Remark}
\section{The Arnold-Liouville theorem}

\begin{Theorem}[Little Arnold-Liouville theorem]\label{thm:littlearnoldliouville}
Let\index{Arnold-Liouville theorem!little} \(\bm{H}\colon
X\to\RR^n\) be an integrable Hamiltonian system and
\(\sigma\colon B\to X\) be a local section. Each orbit
\(\Orb(\sigma(\bm{b}))\) is diffeomorphic to
\(\left(\RR^k/\ZZ^k\right)\times\RR^{n-k}\) for some \(k\). In
particular, if \(\Orb(\sigma(\bm{b}))\) is compact then it is
a torus\index{Lagrangian!torus}.
\end{Theorem}
\begin{Proof}
The action of \(\RR^n\) defines a diffeomorphism
\(\RR^n/\Lambda_{\bm{b}}\to \Orb(\sigma(\bm{b}))\). Since
\(\Lambda_{\bm{b}}\) is a lattice, the result follows from the
classification of lattices in Lemma
\ref{lma:latticeclassification}. \qedhere

\end{Proof}
We now focus attention on a neighbourhood of a regular fibre
(i.e.\ one containing no critical points). By Corollary
\ref{cor:reg_fibres}, a regular fibre is an orbit of the
\(\RR^n\)-action. Since \(\bm{H}\) is proper, its fibres are
compact, so by Theorem \ref{thm:littlearnoldliouville}, a
regular fibre is a torus; this is the analogue of assuming that
our fibres are circles in Theorem \ref{thm:2dactionangle}. Since
the set of regular values is open, we can shrink the domain
\(B\) of our local Lagrangian section so that it is a disc
consisting entirely of regular values. Our goal is to find a map
\(\alpha\colon B\to \RR^n\) such that \(\alpha\circ \bm{H}\) has
standard period lattice.

\begin{Lemma}[Exercise \ref{exr:ham_postcompose}]\label{lma:ham_postcompose}
Let \(\bm{H}\colon X\to B\subset \RR^n\) be an integrable
Hamiltonian system over a disc with only regular fibres, let
\(\alpha\colon B\to C\subset\RR^n\) be a diffeomorphism, and
let \(\bm{G} := \alpha\circ \bm{H}\). Let \(A(b)\) be the
matrix with \(ij\)th entry\footnote{i.e.\ \(i\)th row, \(j\)th
column.} \(A_{ij}(b)=\frac{\partial \alpha_i}{\partial
b_j}(b)\) (the Jacobian of \(\alpha\)). Then:
\begin{itemize}
\item [(i)] the Hamiltonian vector fields of \(\bm{G}\) and
\(\bm{H}\) are related by \(V_{G_i} = \sum_j A_{ij}
V_{H_j}\),
\item [(ii)] the Hamiltonian flows of \(\bm{G}\) and \(\bm{H}\)
are related by \(\phi^{\bm{G}}_{\bm{t}} =
\phi^{\bm{H}}_{A^T\bm{t}}\), and
\item [(iii)] the period lattices \(\Lambda^{\bm{G}}\) and
\(\Lambda^{\bm{H}}\) are related by
\(A^T\Lambda_{\alpha(\bm{b})}^{\bm{G}} =
\Lambda^{\bm{H}}_{\bm{b}}\).

\end{itemize}
\end{Lemma}
\begin{Theorem}[Action-angle coordinates]\label{thm:actionangle}
Let\index{action-angle coordinates} \(\bm{H}\colon X\to
B\subset\RR^n\) be an integrable Hamiltonian system over the
disc with only regular fibres and pick a local Lagrangian
section \(\sigma\). There is a local change of coordinates
\(\alpha\colon B\to C\subset \RR^n\) such that
\(\bm{G}:=\alpha\circ \bm{H}\) generates a Hamiltonian torus
action\index{Hamiltonian!torus action} on \(X\). In other
words, the period lattice \(\Lambda^{\bm{G}}\) is standard and
the map \((\bm{c}, \bm{t}) \mapsto
\phi^{\bm{G}}_{\bm{t}}(\sigma(\alpha^{-1}(\bm{c})))\) defined
in Theorem \ref{thm:liouvillecoordinates} descends to give a
symplectomorphism \(C\times(\RR/2\pi\ZZ)^n\to
\bm{G}^{-1}(C)=\bm{H}^{-1}(B)\).
\end{Theorem}
\begin{Proof}
The following proof is due to Duistermaat \cite{Duistermaat}.

For each \(\bm{b}\in B\), let \(2\pi
\bm{\tau}_1(\bm{b}),\ldots,2\pi \bm{\tau}_n(\bm{b})\in\RR^n\)
be a collection of vectors (smoothly varying in \(\bm{b}\))
which span the lattice of periods
\(\Lambda^{\bm{H}}_{\bm{b}}\). This is possible because \(B\)
is contractible so there is no obstruction to picking sections
of the projection \(\Lambda\to B\). This means that
\(\phi^{\bm{H}}_{\bm{\tau}_i(\bm{b})} = \OP{id}\) for
\(i=1,\ldots,n\). Let us write \(\bm{\tau}_i(\bm{b}) =
(A_{i1}(\bm{b}),\ldots,A_{in}(\bm{b}))\). Let \(A\) be the
matrix with \(ij\)th entry \(A_{ij}(\bm{b})\). Then
\(\Lambda^{\bm{H}}_{\bm{b}}=2\pi A^T\ZZ^n\). By Lemma
\ref{lma:ham_postcompose}(iii), it is sufficient to find a map
\(\alpha=(\alpha_1,\ldots,\alpha_n)\colon B\to\RR^n\) whose
Jacobian \(\partial \alpha_i/\partial b_j\) is \(A_{ij}\).

By the Poincar\'{e} lemma, we can find such functions
\(\alpha_i\) provided
\begin{equation}
\label{eq:poincare}\frac{\partial A_{ij}}{\partial
b_k}=\frac{\partial A_{ik}}{\partial b_j},
\end{equation}
so it remains to check this identity.

Let \(\Psi\colon B\times\RR^n\to X\) be the Liouville
coordinates associated to our choice of Lagrangian section and
\(\Lambda=\Psi^{-1}(\sigma(B))\) be the period lattice. Since
\(\Psi\) is symplectic and \(\sigma(B)\) is Lagrangian,
\(\Lambda\) is Lagrangian. Moreover, \(\Lambda\) is a union of
sheets, each traced out by a single lattice point. For
example, \(\{(\bm{b},\bm{\tau}_i(\bm{b}))\ :\ \bm{b}\in B\}\)
traces out a Lagrangian sheet for each \(i\). In coordinates,
this is \(\{(b_1, \ldots, b_n,A_{i1}(\bm{b}), \ldots,
A_{in}(\bm{b}))\ :\ \bm{b}\in B\}\), which is Lagrangian if
and only if Equation \eqref{eq:poincare} holds (Exercise
\ref{exr:lag_section}). \qedhere

\end{Proof}
\begin{Definition}
The Liouville coordinates associated to the new, periodic
Hamiltonian system are called {\em action-angle
coordinates}\index{action-angle coordinates}. More precisely,
the new Hamiltonians \(\alpha_1\circ \bm{H}, \ldots,
\alpha_n\circ \bm{H}\) are called {\em action coordinates} and
the new \(2\pi\)-periodic conjugate coordinates
\(t_1,\ldots,t_n\) are called {\em angle coordinates}.

\end{Definition}
\begin{Corollary}[Big Arnold-Liouville theorem]\label{cor:bigarnoldliouville}
If\index{Arnold-Liouville theorem!big} \(\bm{H}\colon
M\to\RR^n\) is an integrable Hamiltonian system then any
regular fibre is a torus and admits a neighbourhood
symplectomorphic to \(B\times T^n\), where \(B\subset\RR^n\)
is an open ball and the symplectic form is given by
\(\sum_{i=1}^ndb_i\wedge dt_i\). Under this symplectomorphism,
the orbits of the original system are sent to the tori
\(\{\bm{b}\}\times T^n\).

\end{Corollary}
\section{Solutions to inline exercises}

\begin{Exercise}[Remark \ref{rmk:comm_exr}]\label{exr:comm_exr}
Recall that the flows along two vector fields commute if and
only if the Lie bracket of the vector fields vanishes. Show
that two Hamiltonian flows \(\phi^F_t\) and \(\phi^G_t\)
commute if and only if the Poisson
bracket\index{Poisson!bracket} \(\{F,G\}\) is locally
constant.
\end{Exercise}
\begin{Solution}
We have \([V_F,V_G]=V_{\{F,G\}}\) by Lemma
\ref{lma:liebrackethamiltonian}. Since
\(\iota_{V_{\{F,G\}}}\omega=-d\{F,G\}\), we see that the Lie bracket
vanishes if and only if \(d\{F,G\}=0\) so that all partial
derivatives of \(\{F,G\}\) vanish. This happens if and only if
\(\{F,G\}\) is locally constant. \qedhere

\end{Solution}
\begin{Exercise}[Lemma \ref{lma:heisenberg}]\label{exr:heisenberg}
Let \(F\) and \(G\) be smooth functions. Define
\(F_t(x):=F(\phi^G_t(x))\). Then \(\frac{dF_t}{dt}=\{G,F_t\}\).
\end{Exercise}
\begin{Proof}
We have \(\frac{dF_t}{dt} = dF(V_G) =
-\omega(V_F,V_G)=\{G,F\}\). \qedhere

\end{Proof}
\begin{Exercise}[Lemma \ref{lma:isotropic_orbits}]\label{exr:isotropic_orbits}
If \(L\) is an isotropic\index{isotropic!submanifold}
submanifold of the symplectic manifold \((X,\omega)\) then
\(2\dim(L)\leq \dim(X)\).
\end{Exercise}
\begin{Proof}
For any point \(x\in L\), the tangent space \(T_xL\) is an
isotropic subspace of the symplectic vector space
\(T_xX\). The claim now follows from Lemma
\ref{lma:dimension_bound} in the appendix on symplectic linear
algebra.\qedhere

\end{Proof}
\begin{Exercise}[Lemma \ref{lma:locallagrangian}]\label{exr:locallagrangian}
Let \(\bm{H}\colon X\to \RR^n\) be an integrable Hamiltonian
system. There exists a local Lagrangian section through any
regular point \(x\).
\end{Exercise}
\begin{Proof}
It is a theorem of Darboux (see {\cite[Section 43.B]{Arnold}},
{\cite[Corollary I.1.11]{Audin}}, {\cite[Theorem
3.15]{McDuffSalamon}}) that any point \(x\) in a symplectic
manifold is the centre of a coordinate chart
\((p_1,\ldots,p_n,q_1,\ldots,q_n)\) where the symplectic form
is \(\sum_i dp_i\wedge dq_i\). Let us work locally in these
coordinates. We treat this local chart as a symplectic vector
space and use some notions from the appendix on symplectic
linear algebra. If we define \(J\colon\RR^{2n}\to\RR^{2n}\) to
be the linear map \(J(\bm{p},\bm{q})= (-\bm{q},\bm{p})\) then
\(J\) is an \(\omega\)-compatible complex structure (see
Definition \ref{dfn:cpt_j}). Thus if \(L\) is a Lagrangian
subspace in \(\RR^{2n}\), the subspace \(JL\) is a
complementary Lagrangian subspace (Lemma
\ref{lma:lag_orth}). Since \(T_x\Orb(x)\) is the tangent space
to the orbit \(\Orb(x)\), its image \(JT_x\Orb(x)\) is a
Lagrangian complement. The subspace \(JT_x\Orb(x)\) is the
tangent space of a linear Lagrangian submanifold \(L\) of the
Darboux ball, which is transverse to \(\Orb(x)\) at \(x\). The
differentials \(dH_1,\ldots,dH_n\) are linearly independent at
\(x\) but vanish on \(T_x\Orb(x)\) because \(\bm{H}\) is
constant on \(\Orb(x)\). Therefore these differentials
restrict to linearly independent forms on \(L\) near
\(x\). This implies that the map \(\bm{H}|_L\colon L\to\RR^n\)
is a local diffeomorphism in a neighbourhood of \(x\), so that
its local inverse is a local section of \(\bm{H}\) near \(x\)
whose image is contained in \(L\) and hence
Lagrangian.\qedhere

\end{Proof}
\begin{Exercise}[Lemma \ref{lma:discrete}]\label{exr:discrete}
Let \(\bm{H}\colon X\to\RR^n\) be an integrable Hamiltonian
system, let \(B\subset\RR^n\) be an open set of regular values
and let \(\sigma\colon B\to X\) be a local Lagrangian section;
write \(\Lambda\) for the period lattice. For each \(\bm{b}\in
B\), the intersection \(\Lambda_{\bm{b}} =
\Lambda\cap(\{\bm{b}\}\times\RR^n)\) is a
lattice in \(\RR^n\), that is a discrete
subgroup of \(\RR^n\). The rank of the lattice is lower
semicontinuous as a function of \(\bm{b}\), that is,
\(\bm{b}\) has a neighbourhood \(V\) such that
\(rank(\Lambda_{\bm{b}'})\geq rank(\Lambda_{\bm{b}})\) for all
\(\bm{b}'\in V\).
\end{Exercise}
\begin{Proof}
Let \(\sigma\colon B\to X\) be a local Lagrangian section of
\(\bm{H}\) such that \(\sigma(\bm{b})\) is a regular point of
\(\bm{H}\) for all \(\bm{b}\in B\). We will first show that,
for all \(\bm{b}\in B\), the period lattice
\(\Lambda_{\bm{b}}\) is a discrete subgroup of \(\RR^n\).

The subset \(\Lambda_{\bm{b}}\) is the stabiliser of
\(\bm{b}\) under the action of \(\RR^n\), so it is a subgroup
of \(\RR^n\). To prove discreteness, we need to show that
there is an open set \(W\subset\RR^n\) such that
\(W\cap\Lambda_{\bm{b}} = \{0\}\). Since \(\Psi\colon
B\times\RR^n\to X\) is a local diffeomorphism, there is an
open set \(W'\subset B\times\RR^n\) (containing \((\bm{b},
0)\)) such that \(\Psi\colon W'\to \Psi(W')\) is a
diffeomorphism. There exist open sets \(\bm{b}\in W_1\subset
B\) and \(0\in W_2\subset\RR^n\) such that \(W_1\times
W_2\subset W'\) as these product sets form a basis for the
product topology. In particular, \(0\) is the only point
\(\bm{t}\) in \(W_2\) such that \(\Psi(\bm{b},\bm{t}) =
\bm{b}\). We may therefore take \(W=W_2\) to see that
\(\Lambda_{\bm{b}}\) is discrete.

To see that the rank of the lattice is lower semicontinuous,
we need to show, for each \(\bm{b}\in B\), there is a
neighbourhood \(V\) of \(\bm{b}\) such that
\(rank(\Lambda_{\bm{b}'})\geq rank(\Lambda_{\bm{b}})\) for
\(\bm{b}'\in V\).

Let \(\bm{\lambda}_1(\bm{b}),\ldots,\bm{\lambda}_k(\bm{b})\)
be a \(\ZZ\)-basis for \(\Lambda_{\bm{b}}=\{\bm{t}\in\RR^n\ :
\phi^{\bm{H}}_{\bm{t}}(\bm{b}) = \bm{b}\}\). Then, since
\(\Psi\) is an immersion (Theorem
\ref{thm:liouvillecoordinates}), there is an open
neighbourhood of \(\bm{b}\in B\) such that, for \(\bm{b}'\) in
this open neighbourhood, there are solutions
\(\bm{t}=\bm{\lambda}_1(\bm{b}'), \ldots, \bm{t} =
\bm{\lambda}_k(\bm{b}')\) to the equation
\(\phi^{\bm{H}}_{\bm{t}}(\bm{b}')=\bm{b}'\) which vary
continuously in \(\bm{b}'\). Since the condition of being
linearly independent is an open condition, the points
\(\bm{\lambda}_1(\bm{b}'),\ldots,\bm{\lambda}_k(\bm{b}')\) are
linearly independent for \(\bm{b}'\) in a, possibly smaller,
neighbourhood of \(\bm{b}\), so the rank of the lattice
\(\Lambda_{\bm{b}'}\) is at least \(k\) for \(\bm{b}'\) in a
neighbourhood of \(\bm{b}\). \qedhere

\end{Proof}
\begin{Exercise}[Lemma \ref{lma:ham_postcompose}]\label{exr:ham_postcompose}
Let \(\bm{H}\colon X\to B\subset \RR^n\) be an integrable
Hamiltonian system over a disc with only regular fibres, let
\(\alpha\colon B\to C\subset\RR^n\) be a diffeomorphism, and
let \(\bm{G} := \alpha\circ \bm{H}\). Let \(A(b)\) be the
matrix with \(ij\)th entry\footnote{i.e.\ \(i\)th row, \(j\)th
column.} \(A_{ij}(b)=\frac{\partial \alpha_i}{\partial
b_j}(b)\) (the Jacobian of \(\alpha\)). Then:
\begin{itemize}
\item [(i)] the Hamiltonian vector fields of \(\bm{G}\) and
\(\bm{H}\) are related by \(V_{G_i} = \sum_j A_{ij}
V_{H_j}\),
\item [(ii)] the Hamiltonian flows of \(\bm{G}\) and \(\bm{H}\)
are related by \(\phi^{\bm{G}}_{\bm{t}} =
\phi^{\bm{H}}_{A^T\bm{t}}\), and
\item [(iii)] the period lattices \(\Lambda^{\bm{G}}\) and
\(\Lambda^{\bm{H}}\) are related by
\(A^T\Lambda_{\alpha(\bm{b})}^{\bm{G}} =
\Lambda^{\bm{H}}_{\bm{b}}\).
\end{itemize}
\end{Exercise}
\begin{Solution}
Let us write \(A_{ij} = \frac{\partial \alpha_i}{\partial
b_j}\). We have
\begin{align*}
\iota_{\sum_jA_{ij}V_{H_j}}\omega &= \sum_j\frac{\partial \alpha_i}{\partial b_j}\iota_{V_{H_j}}\omega\\
&= -\sum_j\frac{\partial \alpha_i}{\partial b_j}dH_j\\
&= -d(\alpha_i\circ \bm{H}) = -dG_i.
\end{align*}
This proves (i): \(V_{G_i} = \sum_j A_{ij} V_{H_j}\). Thus, if
\(\bm{t}\) is the row vector \((t_1,\ldots,t_n)\), then
\[\sum_i t_i V_{G_i} = \sum_{i,j}t_iA_{ij}V_{H_j},\] where the
matrix \(A_{ij}\) is constant on each orbit. Therefore we
obtain (ii):
\(\phi^{\bm{G}}_{\bm{t}}=\phi^{\bm{H}}_{A^T\bm{t}}\).

The lattice \(\Lambda^{\bm{G}}_{\alpha(\bm{b})}\) of
\(\bm{G}\) on the consists of tuples
\(\bm{t}=(t_1,\ldots,t_n)\) such that
\(\phi^{\bm{G}}_{\bm{t}}=\OP{id}\) on
\(\bm{G}^{-1}(\alpha(\bm{b}))\). By (ii), this is equivalent
to \(\phi^{\bm{H}}_{A^T\bm{t}}=\OP{id}\) on the orbit
\(\bm{H}^{-1}(\bm{b})\), so
\(\bm{t}\in\Lambda^{\bm{G}}_{\alpha(\bm{b})}\) if and only if
\(A^T\bm{t}\in\Lambda^{\bm{H}}_{\bm{b}}\), which gives (iii):
\[A^T\Lambda^{\bm{G}}_{\alpha(\bm{b})} =
\Lambda^{\bm{H}}_{\bm{b}}.\qedhere\]

\end{Solution}
\begin{Exercise}[From proof of Theorem \ref{thm:actionangle}]\label{exr:lag_section}
Show that a section \(\sigma(\bm{b})=(\bm{b},\bm{t}(\bm{b}))\)
is Lagrangian with respect to the symplectic form \(\omega =
\sum db_i\wedge dt_i\) if and only if \(\partial t_i/\partial
b_j = \partial t_j/\partial b_i\) for all \(i,j\).
\end{Exercise}
\begin{Solution}
The tangent space to the section \(\sigma\) is spanned by the
vectors \(\sigma_*(\partial_{b_i})\) so it suffices to check
that \(\omega(\sigma_*(\partial_{b_i}),
\sigma_*(\partial_{b_j})) = 0\) for all \(i,j\). We have
\(\sigma_*(\partial_{b_i}) = \partial_{b_i} + \sum_k(\partial
t_i/\partial b_k)\partial_{b_k}\), which gives
\[\omega(\sigma_*(\partial_{b_i}), \sigma_*(\partial_{b_j})) =
\partial t_i/\partial b_j - \partial t_j/\partial
b_i.\qedhere\]

\end{Solution}
\chapter{Lagrangian fibrations}
\label{ch:lag_fib}
\thispagestyle{cup}

We have seen that an integrable Hamiltonian system is a map
\(X\to\RR^n\) whose regular fibres are Lagrangian
submanifolds. This structure, called a {\em Lagrangian
fibration}\footnote{The word {\em fibration} also appears in
algebraic topology (e.g.\ {\em Serre fibrations}) where it
describes maps with a homotopy lifting property. Lagrangian
torus fibrations are not fibrations in that sense: though they
are fibre-bundles over the regular locus, homotopy lifting fails
near the critical points. This is an unfortunate accident of
history.} turns out to be very useful for studying the geometry
and topology of symplectic manifolds.

In this chapter, we introduce a general definition of Lagrangian
fibration. We then discuss the {\em regular Lagrangian
fibrations}: those with no critical points, i.e.\ proper
submersions \(X\to B\) with connected Lagrangian fibres. We will
see that these are locally the same as integrable Hamiltonian
systems (Remark
\ref{rmk:locally_modelled_on_integrable_system}). In particular,
the fibres are tori (Corollary \ref{cor:fibre_tori}). For this
reason, we often use the name {\em Lagrangian torus fibration}
instead of Lagrangian fibration. Next, we will see that local
action coordinates equip the image \(B\) with a geometric
structure called an {\em integral affine structure}, which can
also be understood in terms of the symplectic areas of cylinders
connecting fibres. Finally, we will show that under certain
assumptions (existence of a global Lagrangian section), the
integral affine manifold \(B\) is enough information to
reconstruct the Lagrangian fibration \(X\to B\) completely.

As the book progresses, we will allow our fibrations to have
progressively worse critical points.

\section{Lagrangian fibrations}

\begin{Definition}
Recall that a stratification\index{stratification} of a
topological space \(B\) is a filtration
\[\emptyset=:B_{-1}\subset B_0\subset\cdots\subset B_d\subset
B_{d+1}\subset\cdots\subset B,\] where each \(B_d\) is a
closed subset such that, for each \(d\), the \(d\)-stratum
\(S_d(B):=B_d\setminus B_{d-1}\) is a smooth \(d\)-dimensional
manifold (possibly empty) and \(B=\bigcup_{d\geq 0}B_d\). We
say that \(B\) is finite-dimensional if the \(d\)-stratum is
empty for sufficiently large \(d\), and we say that \(B\) is
\(n\)-dimensional if \(B\) is finite-dimensional and \(n\) is
maximal such that \(S_n(B)\) is nonempty (in this case we call
\(S_n(B)\) the {\em top stratum}).

\end{Definition}
We adopt the following working definition of a Lagrangian torus
fibration, given in {\cite[Definition 2.5]{EvansMauri}}. It is
extremely weak because it places no restrictions on the critical
points of the fibration.

\begin{Definition}
Let\index{fibration!Lagrangian torus, general definition}
\((X,\omega)\) be a \(2n\)-dimensional symplectic manifold and
\(B\) be an \(n\)-dimensional stratified space. A Lagrangian
torus fibration \(f\colon X\to B\) is a proper continuous map
such that \(f\) is a smooth submersion over the top stratum
with connected Lagrangian fibres and the other fibres are
themselves connected stratified spaces with isotropic
strata. We call \(B^{reg}:= S_n(B)\) the {\em regular locus}
of \(H\) and \(B^{sing} := B\setminus S_n(B)\) the {\em
discriminant locus}.\index{discriminant locus}

\end{Definition}
\begin{Remark}
Throughout Chapter 1, \(B\) denoted an open subset of
\(\RR^n\). This is no longer the case. However, it is still
the target (``base'') of the fibration, hence the choice of
letter.

\end{Remark}
\section{Regular Lagrangian fibrations}

We first study Lagrangian fibrations with no critical points. It
turns out (Lemma \ref{lma:lagfibreg}) that these are locally
equivalent to integrable Hamiltonian systems.

\begin{Definition}
We say that a Lagrangian fibration \(f\colon X\to B\) is {\em
regular}\index{fibration!regular
Lagrangian|(}\index{Lagrangian!fibration|see {fibration,
Lagrangian}} if \(B = B^{reg}\), that is if \(f\) is a smooth
proper submersion with connected Lagrangian fibres.

\end{Definition}
\begin{Lemma}\label{lma:invariantlagrangians}
Let \((X,\omega)\) be a symplectic manifold. Suppose that
\(H\colon X\to\RR\) is a Hamiltonian function and \(L\subset
X\) is a Lagrangian submanifold such that \(L\subset
H^{-1}(c)\) for some \(c\in\RR\). Then \(\phi^H_t(x)\in L\)
for all \(x\in L\), \(t\in\RR\), i.e.\ \(L\) is invariant
under the Hamiltonian flow of \(H\).
\end{Lemma}
\begin{Proof}
Since \(L\subset H^{-1}(c)\), the function \(H\) is constant
on \(L\), so the directional derivative \(v(H)=dH(v)\)
vanishes whenever \(v\in TL\). We have
\(\iota_{V_H}\omega=-dH\). If \(v\in TL\) then
\[\omega(V_H,v)=-dH(v)=0.\] This means that \(V_H\) is in the
symplectic orthogonal complement\footnote{See Definition
\ref{dfn:symplectic_orthogonal_complement} for the definition
of the symplectic orthogonal complement.}
\((TL)^{\omega}\). Since \(L\) is Lagrangian,
\(TL=(TL)^{\omega}\), so this shows that \(V_H\in TL\). Since
\(V_H\) is tangent to \(L\), the flow of \(V_H\) preserves
\(L\). \qedhere

\end{Proof}
\begin{Lemma}\label{lma:lagfibreg}
Let \((X,\omega)\) be a symplectic \(2n\)-manifold, \(B\) be
an \(n\)-manifold and let \(f\colon X\to B\) be a regular
Lagrangian fibration. Let \((b_1,\ldots,b_n)\) be local
coordinates on \(B\). The functions \(b_1\circ
f,\ldots,b_n\circ f\) Poisson commute\index{Poisson!commute}.
\end{Lemma}
\begin{Proof}
Fix a point \(c\in B\) with \(b_i(c)=c_i\). The Lagrangian
fibre \(f^{-1}(c)\) is contained in all the level sets
\(\{b_i\circ f=c_i\}\), \(i=1,\ldots,n\). By Lemma
\ref{lma:invariantlagrangians}, the Hamiltonian vector field
\(V_{b_i\circ f}\) is tangent to \(L\) (for all
\(i\)). Therefore \[\{b_i\circ f,b_j\circ
f\}=\omega(V_{b_i\circ f},V_{b_j\circ f})=0\] because
\(V_{b_i\circ f},V_{b_j\circ f}\in TL\) and \(L\) is
Lagrangian.\qedhere

\end{Proof}
\begin{Remark}\label{rmk:locally_modelled_on_integrable_system}
In particular, \(f\) is locally modelled on an integrable
Hamiltonian system.

\end{Remark}
\begin{Corollary}[Exercise \ref{exr:fibre_tori}]\label{cor:fibre_tori}
If \(f\colon X\to B\) is a proper submersion with connected
Lagrangian fibres then the fibres are Lagrangian
tori.\index{fibration!regular Lagrangian|)}
\index{Lagrangian!torus}

\end{Corollary}
\section{Integral affine structures}

The big Arnold-Liouville theorem (Corollary
\ref{cor:bigarnoldliouville}) gives us more information than
Corollary \ref{cor:fibre_tori}: we will be able to show that the
base of the Lagrangian fibration has an {\em integral affine
structure}\index{integral affine!structure|(}.

\begin{Definition}
An integral affine transformation\index{integral
affine!transformation} is a map \(T\colon\RR^n\to\RR^n\) of
the form\footnote{We think of \(\RR^n\) as consisting of row
vectors and matrices acting on the right.}
\(T(\bm{b})=\bm{b}A+\bm{C}\) where \(A\in GL(n,\ZZ)\) and
\(\bm{C}\in\RR^n\). An integral affine structure on a manifold
\(B\) is an atlas for \(B\) whose transition functions are
integral affine transformations.

\end{Definition}
\begin{Lemma}\label{lma:zaffine}
Suppose \(\bm{G}\colon X\to\RR^n\) and \(\bm{H}\colon
X\to\RR^n\) are submersions defining integrable Hamiltonian
systems such that the period lattices are both
standard\index{period lattice!standard}. Suppose that
\(\psi\colon \bm{H}(X)\to \bm{G}(X)\) is a diffeomorphism such
that \(\bm{G}=\psi\circ \bm{H}\). Then \(\psi\) is (the
restriction to \(\bm{H}(X)\) of) an integral affine
transformation.
\end{Lemma}
\begin{Proof}
Let \(\phi^{\bm{G}}_{\bm{t}} =
\phi^{G_1}_{t_1}\cdots\phi^{G_n}_{t_n}\) and
\(\phi^{\bm{H}}_{\bm{t}} =
\phi^{H_1}_{t_1}\cdots\phi^{H_n}_{t_n}\) be the Hamiltonian
\(\RR^n\)-actions. Since \(\bm{G}=\psi\circ \bm{H}\), Lemma
\ref{lma:ham_postcompose}(iii) implies
\(A(\bm{b})\Lambda_{\psi(\bm{b})}^{\bm{G}} =
\Lambda_{\bm{b}}^{\bm{H}}\) where
\(A(\bm{b})=d_{\bm{b}}\psi\). Since both period lattices are
assumed to be standard, this means \(A(\bm{b})\in GL(n,\ZZ)\)
for all \(\bm{b}\in \bm{G}(X)\). Since \(GL(n,\ZZ)\) is
discrete, this is only possible if \(d\psi\) is constant. Thus
\(\psi(\bm{b})=\bm{b}A+\bm{C}\) for some \(A\in GL(n,\ZZ)\)
and \(\bm{C}\in\RR^n\). \qedhere

\end{Proof}
\begin{Remark}
This proof contains the first instance of a useful trick we
will use repeatedly in what follows. Namely, by showing that
the derivative of \(\psi\) belongs to some discrete set, we
were able to severely constrain \(\psi\). For further examples
of this trick in action, see Proposition
\ref{prp:straight_lines} (the boundary of the moment polytope
is piecewise linear) and Theorem \ref{thm:visibility}
(``visible Lagrangians'' live over straight lines).

\end{Remark}
\begin{Theorem}\label{thm:zaffine_base}
If \(f\colon X\to B\) is a regular Lagrangian fibration then
\(B\) inherits an integral affine structure.
\end{Theorem}
\begin{Proof}
Suppose we are given a coordinate chart\footnote{We write
partially-defined maps with \(\dashrightarrow\) to save
overburdening the notation with domains and targets.}
\(\varphi\colon B\dashrightarrow \RR^n\). By Lemma
\ref{lma:lagfibreg}, \(\varphi\circ f\) is an integrable
Hamiltonian system. Let
\(\alpha\colon\RR^n\dashrightarrow\RR^n\) be the map
constructed in the proof of Theorem \ref{thm:actionangle} so
that \(\alpha\circ\varphi\circ f\) are action
coordinates. This gives us a modified chart
\(\alpha\circ\varphi\colon B\dashrightarrow \RR^n\). If we
modify a whole atlas in this way, we obtain a new atlas; we
will check that the resulting transition functions are
integral affine transformations. Suppose we have charts
\(\varphi_1\colon B\dashrightarrow \RR^n\) and
\(\varphi_2\colon B\dashrightarrow\RR^n\) which we modify
using \(\alpha_1\colon \RR^n\dashrightarrow \RR^n\),
\(\alpha_2\colon \RR^n\dashrightarrow \RR^n\). The transition
map for the modified atlas is \(\psi_{12} :=
\alpha_2\circ\varphi_2\circ\varphi_1^{-1}\circ
\alpha_1^{-1}\). We know that \(\bm{H} :=
\alpha_1\circ\varphi_1\circ f\) and \(\bm{G} :=
\alpha_2\circ\varphi_2\circ f\) are integrable systems with
standard period lattice, and \(\psi_{12}\circ \bm{H} =
\bm{G}\), so by Lemma \ref{lma:zaffine}, \(\psi_{12}\) is an
integral affine transformation. \qedhere

\end{Proof}
\begin{Remark}\label{rmk:smooth_atlas}
In the construction of this integral affine
structure\index{integral affine!structure|)}, we modified the
atlas and, hence, the smooth structure of \(B\). In other
words, we don't get to pick the smooth structure on \(B\): it
is dictated to us by the geometry of the fibration.

\end{Remark}
\section{Flux map}
\label{sct:flux}

There is a more geometric way to characterise the action
coordinates. Let \(f\colon X\to B\) be a regular Lagrangian
fibration. We assume for simplicity\footnote{Exercise
\ref{exr:non_exact}: Explain how to modify the construction to
get an integral affine structure on \(B\) even if \(\omega\) is
not exact. Disclaimer: This is one of the exercises that
requires a lot of work.} that \(\omega=d\lambda\) for some
1-form \(\lambda\).

Consider the local system \(\xi\to B\) whose fibre over \(b\) is
the abelian group \(H_1(f^{-1}(b);\ZZ)\cong\ZZ^n\). Let
\(p\colon\tilde{B}\to B\) be the universal cover and let
\(\tilde{\xi}=p^*\xi\). Since \(\tilde{B}\) is simply-connected,
\(\tilde{\xi}\) is trivial. Let \(c_1,\ldots,c_n\) be a
\(\ZZ\)-basis of continuous sections of
\(\tilde{\xi}\to\tilde{B}\).

\begin{Definition}[Flux map]\label{dfn:fluxmap}
The {\em flux map}\index{flux|(} is defined to be the map
\(\II\colon\tilde{B}\to\RR^n\) given by
\[\II(\tilde{b})=(I_1(\tilde{b}),\ldots,I_n(\tilde{b})):=
\left(\frac{1}{2\pi}\int_{c_1(\tilde{b})}\lambda,\ldots,\frac{1}{2\pi}\int_{c_n(\tilde{b})}\lambda\right).\]

\end{Definition}
\begin{Lemma}[Flux map = action coordinates]\label{lma:fluxaction}
Suppose that \(\tilde{U}\subset\tilde{B}\) and \(U\subset B\)
are open subsets such that \(p|_{\tilde{U}}\colon\tilde{U}\to
U\) is a diffeomorphism. Then \(\II\circ
(p|_{\tilde{U}})^{-1}\colon U\to\RR^n\) gives action
coordinates on \(U\).\index{action and flux}
\end{Lemma}
\begin{Proof}
By Corollary \ref{cor:bigarnoldliouville}, it is sufficient to
prove this for the local model \((U\times T^n,\omega_0)=\sum
db_i\wedge dt_i)\). In that case, we can pick \(\lambda=\sum
b_idt_i\) and take \(c_1,\ldots,c_n\) to be the standard basis
of \(H_1(T^n;\ZZ)\). Then we get \(I_i(b)=b_i\), which
recovers the action coordinates. \qedhere

\end{Proof}
\begin{Definition}[Fundamental action domain]\label{dfn:fad}
We call \(\II(\tilde{U})\) a {\em fundamental action
domain}\index{action domain, fundamental} for the Lagrangian
fibration.

\end{Definition}
\begin{Remark}
If we pick a different \(\lambda'\) such that
\(d\lambda'=d\lambda\) then \(\lambda-\lambda'\) is closed, so
\(\int_{c_i(b)}(\lambda-\lambda')\) is constant (by Stokes's
theorem) and the flux map changes by an additive constant. If
we pick a different \(\ZZ\)-basis \((c'_1,\ldots,c'_n)\) then
we can express the new integrals as a \(\ZZ\)-linear
combination of \(I_1,\ldots,I_n\). This means that the flux
map is determined up to an integral affine transformation.

\end{Remark}
The integral affine structure\index{integral affine!structure|(}
from Theorem \ref{thm:zaffine_base} can now be understood in the
following way. We pull back the integral affine structure from
\(\RR^n\) along \(\II\) to get an integral affine structure on
\(\tilde{B}\); this integral affine structure on \(\tilde{B}\)
descends to one on \(B\) (it is invariant under the
action\footnote{Conventions: We think of \(\II(b)\) as a row
vector, write concatenation of loops as \(\alpha\cdot\beta\)
meaning ``follow \(\alpha\) then \(\beta\)'', and write the deck
group acting on the right.} of deck transformations). We will
prove this because it introduces an important new idea: the {\em
affine monodromy}.

\begin{Corollary}\label{cor:zaffine}
If we equip \(\tilde{B}\) with the integral affine structure pulled
back from \(\RR^n\) along \(\II\) then it is invariant under the
action of the deck group of the cover \(p\colon \tilde{B}\to B\).
\end{Corollary}
\begin{Proof}
If \(g\colon\tilde{B}\to\tilde{B}\) is a deck transformation of the
cover \(p\) then \(c_1(\tilde{b}),\ldots,c_n(\tilde{b})\) and
\(c_1(\tilde{b}g),\ldots,c_n(\tilde{b}g)\) are both
\(\ZZ\)-bases for the \(\ZZ\)-module
\(H_1(f^{-1}(p(\tilde{b}));\ZZ)\) and therefore they are
related by some change-of-basis matrix \(M(g)\in
GL(n,\ZZ)\). This implies that
\(\II(\tilde{b}g)=\II(\tilde{b})M(g)\). Since \(M(g)\) is an integral
affine transformation, this shows that the integral affine structure
descends to the quotient \(B\). \qedhere

Note that, with our conventions, \(M(g_1g_2)=M(g_1)M(g_2)\). Indeed,
\(M\colon\pi_1(B)\to GL(n,\ZZ)\) is the monodromy of the
local system \(\xi\to A\).

\end{Proof}
\begin{Definition}\label{dfn:affine_monodromy}
We call \(M\colon\pi_1(B)\to GL(n,\ZZ)\) the {\em affine
monodromy}\index{affine monodromy} in what follows. The first
example we will encounter where the affine monodromy is
nontrivial will be the fibrations with focus-focus critical
points in Chapter \ref{ch:focusfocus}.

\end{Definition}
\begin{Remark}
The manifold \(B\) can be reconstructed in the usual way as a
quotient of a closed fundamental domain for the universal
cover \(\tilde{B}\to B\) where the identifications are made
using deck transformations. If we wish to reconstruct the
integral affine structure on \(B\) then we use a fundamental
action domain and the identifications are made using the
integral affine transformations \(M(g)\) corresponding to deck
transformations \(g\).

\end{Remark}
\begin{Remark}
Given any integral affine manifold \(B\), there is a {\em
developing map}\index{developing map}, that is a
(globally-defined) local diffeomorphism
\(\II\colon\tilde{B}\to\RR^n\) from the universal cover into
Euclidean space such that the integral affine structure
inherited by \(\tilde{B}\) from the covering map agrees with
the pullback of the integral affine structure along the
developing map. In our context, the flux map is the developing
map.\index{flux|)}

\end{Remark}
\begin{Remark}
Suppose that \(f\colon X\to B\) is an integrable system with
\(B\subset\RR^n\), so that \(B\) already has an integral
affine structure as open subset of \(\RR^n\). This does not
agree with the integral affine structure\index{integral
affine!structure|)} constructed in Corollary \ref{cor:zaffine}
unless the period lattice is standard.

\end{Remark}
\section{Uniqueness}

\begin{Definition}
Let \(f\colon X\to B\) and \(g\colon Y\to C\) be regular
Lagrangian fibrations. If \(\phi\colon B\to C\) is a
diffeomorphism then a {\em symplectomorphism fibred over
\(\phi\)}\index{fibred symplectomorphism|(} is a
symplectomorphism \(\Phi\colon X\to Y\) such that \(g\circ\Phi
= \phi\circ f\).

\begin{center}
\begin{tikzpicture}
\node (X) at (0,0) {\(X\)};
\node (Y) at (2,0) {\(Y\)};
\node (B) at (0,-2) {\(B\)};
\node (C) at (2,-2) {\(C\)};
\draw[->] (X) -- (B) node [midway,left] {\(f\)};
\draw[->] (Y) -- (C) node [midway,right] {\(g\)};
\draw[->] (X) -- (Y) node [midway,above] {\(\Phi\)};
\draw[->] (B) -- (C) node [midway,below] {\(\phi\)};

\end{tikzpicture}
\end{center}
If \(\phi=\OP{id}\), we will simply call \(\Phi\) a {\em
fibred symplectomorphism} and if moreover \(f=g\) then we call
\(\Phi\) a {\em fibred automorphism of \(f\)}.\index{fibred
automorphism|(}

\end{Definition}
An argument similar to the one which proved Lemma
\ref{lma:zaffine} shows that the map \(\phi\) is an isomorphism
of integral affine manifolds \(B\to C\). We now tackle the
converse question: if there is an integral affine isomorphism
\(\phi\colon B\to C\), is there a symplectomorphism \(X\to Y\)
fibred over \(\phi\)? We first prove some preliminary lemmas.

\begin{Lemma}\label{lma:unique_fibred_symp}
Let \(\Phi\colon X\to X\) be a fibred automorphism of
\(f\colon X\to B\) and suppose there is a Lagrangian section
\(\sigma\colon B\to X\) such that
\(\Phi\circ\sigma=\sigma\). Then \(\Phi=\OP{id}\).
\end{Lemma}
\begin{Proof}
The property that \(\Phi=\OP{id}\) can be checked locally, so
we lose nothing by passing to a small affine coordinate chart
in \(B\). Without loss of generality, therefore, we will
assume that \(f=\bm{H}\colon X\to B\subset\RR^n\) is an
integrable Hamiltonian system with Lagrangian section
\(\sigma\). By Corollary \ref{cor:bigarnoldliouville},
\(X\cong B\times T^n\) with symplectic form \(\omega = \sum
db_i\wedge dt_i\). Since we have used the section \(\sigma\)
to define the Liouville coordinates, the section is given in
these coordinates by \(\sigma(\bm{b})=(\bm{b},0)\). The fact
that \(\Phi\) is fibred means that
\(\Phi(\bm{b},\bm{t})=(\bm{b},\bm{q}(\bm{b},\bm{t}))\) for
some function \(\bm{q}(\bm{b},\bm{t})\). The condition that
\(\Phi\) is symplectic means in particular that
\(\omega(\Phi_*\partial_{b_i},\Phi_*\partial_{t_j}) =
\delta_{ij}\), which becomes \(\partial q_i/\partial
t_j=\delta_{ij}\). Upon integrating, this means
\(\bm{q}(\bm{b},\bm{t}) = \bm{q}(\bm{b},0)+\bm{t}\), so the
condition \(\bm{q}(\bm{b},0)=0\) tells us that
\(\bm{q}(\bm{b},\bm{t})=\bm{t}\), and hence \(\Phi\) is the
identity. \qedhere

\end{Proof}
\begin{Lemma}\label{lma:globcoords}
Assume that \(\bm{F}\colon X\to\RR^n\) and \(\bm{G}\colon
Y\to\RR^n\) are integrable Hamiltonian systems with no
critical points. Assume that the period lattices
\(\Lambda^{\bm{F}}\) and \(\Lambda^{\bm{G}}\) are both
standard, and that we are given global Lagrangian sections
\(\sigma\) of \(\bm{F}\) and \(\tau\) of \(\bm{G}\). Suppose
there is an integral affine transformation
\(\phi\colon\RR^n\to\RR^n\) such that
\(\phi(\bm{F}(X))=\bm{G}(X)\). Then there is a unique
symplectomorphism \(\Phi\colon X\to Y\) fibred over \(\phi\)
satisfying \(\Phi\circ\sigma=\tau\circ\phi\).
\end{Lemma}
\begin{Proof}
Write \(\bm{F}=(F_1, \ldots, F_n)\) and \(\bm{G}=(G_1, \ldots,
G_n)\). Let \((s_1,\ldots,s_n)\) and \((t_1,\ldots,t_n)\) be
the \(2\pi\)-periodic Liouville (angle) coordinates associated
to the Lagrangian sections. Write
\(\phi(\bm{b})=\bm{b}A+\bm{C}\) for some \(A\in GL(n,\ZZ)\)
and \(\bm{C}\in\RR^n\). As usual, we think of \(\bm{F}\) and
\(\bm{G}\) as row vectors and write \(A\) acting on the right.

By Corollary \ref{cor:bigarnoldliouville}, \(X\) is
symplectomorphic to \(\bm{F}(X)\times T^n\) with symplectic
form \(\sum_i dF_i\wedge ds_i\) and \(Y\) is symplectomorphic
to \(\bm{G}(X)\times T^n\) with symplectic form \(\sum_i
dG_i\wedge dt_i\). Under these identifications, we have
\(\sigma(\bm{b})=(\bm{b},0)\) and \(\tau(\bm{c})=(\bm{c},0)\).

Define a map \(\bm{F}(X)\times\RR^n\to\bm{G}(X)\times\RR^n\)
by \[(\bm{c},\bm{t})= \left(\bm{b}A + \bm{C},
A^{-1}\bm{s}\right).\] Because \(A\in GL(n,\ZZ)\), and because
both period lattices \(\Lambda^{\bm{F}}\) and
\(\Lambda^{\bm{G}}\) are standard, the matrix \(A^{-1}\) sends
\(\Lambda^{\bm{F}}\) isomorphically to \(\Lambda^{\bm{G}}\),
and descends to a well-defined diffeomorphism
\(\Phi\colon\bm{F}(X)\times T^n\to \bm{G}(X)\times T^n\). We
need to show \(\Phi\) is symplectic. We have \(dG_j=\sum_i
dF_iA_{ij}\) and \(dt_j=\sum_k A^{-1}_{jk}ds_k\), so \[\sum_j
dG_j\wedge dt_j=\sum_{i,j,k} A_{ij} A^{-1}_{jk}\,dF_i\wedge
ds_k=\sum_{i,k}\delta_{ik}dF_i\wedge ds_k=\sum_i dF_i\wedge
ds_i,\] which shows that \(\Phi\) is a symplectic map.

Note that, by construction, \[\Phi(\sigma(\bm{b})) =
\Phi(\bm{b},0) = (\bm{b}A+\bm{C},0)=\tau(\phi(\bm{b})).\] If
\(\Phi'\) were another symplectomorphism fibred over \(\phi\)
with this property then \(\Phi^{-1}\circ\Phi'\) would be a
fibred automorphism\index{fibred automorphism|)} of \(\bm{F}\)
fixing \(\sigma\), and hence equal to the identity by Lemma
\ref{lma:unique_fibred_symp}. \qedhere

\end{Proof}
From now on, we will suppose for convenience that
\(\phi=\OP{id}\), so that we have two regular Lagrangian
fibrations \(f\colon X\to B\) and \(g\colon Y\to B\) which equip
\(B\) with the same integral affine structure and we ask if
there is a fibred symplectomorphism \(\Phi\colon X\to Y\).

\begin{Theorem}\label{thm:uniqueness}
Suppose that we have regular Lagrangian fibrations \(f\colon
X\to B\) and \(g\colon Y\to B\) over the same integral affine
base. Suppose moreover that both fibrations admit global
Lagrangian sections \(\sigma\) and \(\tau\). Then there is a
unique fibred symplectomorphism \(\Phi\colon X\to Y\) such
that \(\Phi\circ\sigma=\tau\).
\end{Theorem}
\begin{Proof}
Given a sufficiently small \(U\subset B\), Lemma
\ref{lma:globcoords} produces a unique fibred
symplectomorphism \(\Phi_U\colon f^{-1}(U)\to g^{-1}(U)\)
satisfying \(\Phi_U\circ\sigma = \tau\). We would like to
define \(\Phi\) by \(\Phi(x)=\Phi_U(x)\) if \(f(x)\in U\). The
only thing to check is that this prescription is well-defined
independently of the choice of \(U\). In other words, given
subsets \(U,V\subset B\) and \(x\in X\) such that \(f(x)\in
U\cap V\), we want to show that \(\Phi_U(x)=\Phi_V(x)\). Since
\(\Phi_U\circ\sigma=\tau\) and \(\Phi_V\circ\sigma=\tau\), we
see that the restrictions of these fibred symplectomorphisms
to \(f^{-1}(U\cap V)\) must agree by the uniqueness part of
Lemma \ref{lma:globcoords}, so \(\Phi_U(x)=\Phi_V(x)\), as
required. \qedhere

\end{Proof}
The assumption that there is a global Lagrangian section is
necessary, as the following example illustrates.

\begin{Example}
Consider the quotient \(K\) of the product \(\RR\times T^3\)
by the equivalence relation \((t,x,y,z)\sim
(t+1,x,y,y+z)\). The symplectic form \(\omega = dt\wedge
dx+dy\wedge dz\) descends to \(K\) because \(d(t+1)\wedge
dx+dy\wedge d(y+z)=dt\wedge dx+dy\wedge dz\). The symplectic
manifold \((K,\omega)\) is called the Kodaira-Thurston
manifold\index{Kodaira-Thurston manifold} and was the first
known example of a symplectic manifold which does not admit a
compatible K\"{a}hler structure\footnote{You can see this
because, for example, the first Betti number of a K\"{a}hler
manifold must be even, but \(b_1(K)=3\).}; see
\cite{Thurston}.

The projection \((t,x,y,z)\mapsto (t,y)\) is a well-defined
regular Lagrangian fibration \(K\to T^2\). The action of
\((\theta_1,\theta_2)\in T^2\) by \((t,x,y,z)\mapsto
(t,x+\theta_1,y,z+\theta_2)\) has the fibres of \(f\) as its
orbits. If there were a section\footnote{Lagrangian or not.}
\(T^2\to K\), say \((t,y)\mapsto (t,x(t,y),y,z(t,y))\), then
one would get a diffeomorphism \(T^4\to K,\quad
(t,y,\theta_1,\theta_2)\mapsto
(t,x(t,y)+\theta_1,y,z(t,y)+\theta_2)\). There is no such
diffeomorphism because \(K\not\cong T^4\) (for example,
\(b_1(K)=3\neq 4 = b_1(T^4)\)). Therefore there is no section.

The base of this fibration is the torus \(T^2\) with its
product integral affine structure. This same integral affine
manifold arises as the base of a different Lagrangian
fibration: the standard torus fibration \(T^4\to T^2\) where
we equip \(T^4\) with the symplectic form \(d\theta_1\wedge
d\theta_2 + d\theta_3 \wedge d\theta_4\) and the torus
fibration is \(\bm{\theta}\mapsto(\theta_1,\theta_3)\). This
shows that it is possible to have two inequivalent Lagrangian
fibrations over the same integral affine base provided one of
them does not admit a global Lagrangian section.

\end{Example}
\begin{Remark}
In fact, one can also compare two Lagrangian fibrations
\(f\colon X\to B\) and \(g\colon Y\to B\) without assuming the
existence of a global Lagrangian section. Given a subset
\(U\subset B\), consider the set \(\mathcal{S}(U)\) of fibred
symplectomorphisms \(\Phi\colon f^{-1}(U)\to g^{-1}(U)\). This
assignment \(U\mapsto \mathcal{S}(U)\) is a {\em sheaf} over
\(B\). Using the language of sheaf theory, one can formulate
an analogue of Theorem \ref{thm:uniqueness} without mentioning
Lagrangian sections. There is an element
\(\Phi\in\mathcal{S}(B)\) (i.e.\ a fibred
symplectomorphism\index{fibred symplectomorphism|)}) if and
only if a certain characteristic class vanishes. See
{\cite[Section 2]{Duistermaat}} for a full discussion.

\end{Remark}
When we do have global Lagrangian sections, Theorem
\ref{thm:uniqueness} is a wonderful compression of information:
to reconstruct our \(2n\)-dimensional space \(X\), all we need
is an \(n\)-dimensional integral affine manifold. For example,
if \(n=2,3\), this brings 4- and 6-dimensional spaces into the
range of visualisation.

\section{Lagrangian and non-Lagrangian sections}

We now turn to the question of when a Lagrangian fibration
admits a Lagrangian section\index{Lagrangian!section}. First we
see what happens to the symplectic form in Liouville coordinates
when we pick a non-Lagrangian section.

\begin{Lemma}\label{lma:liouville_nonlag}
Let \(\bm{H}\colon X\to\RR^n\) be an integrable Hamiltonian
system, let \(B\subset\bm{H}(X)\subset\RR^n\) be an open set,
and let \(\sigma\colon B\to X\) be a (not necessarily
Lagrangian) section. Define \[\Psi\colon B\times\RR^n\to
X,\qquad\Psi(\bm{b},\bm{t}) =
\phi^{\bm{H}}_{\bm{t}}(\sigma(\bm{b})).\] Let \(\beta\) denote
the pullback of the 2-form \(\sigma^*\omega\) on \(B\) to
\(B\times\RR^n\). Then \(\Psi\) is both an immersion and a
submersion and \(\Psi^*\omega=\sum db_i\wedge dt_i + \beta\),
where \((b_1,\ldots,b_n)\) are the standard coordinates on
\(B\subset\RR^n\).
\end{Lemma}
\begin{Proof}
The only difference with the proof of Theorem
\ref{thm:liouvillecoordinates} is that
\(\omega(\Psi_*\partial_{b_i},\Psi_*\partial_{b_j})\) does not
need to vanish. Instead,
\begin{align*}\omega(\Psi_*\partial_{b_i},
\Psi_*\partial_{b_i}) &=
\omega((\phi^{\bm{H}}_{\bm{t}})_*\sigma_*\partial_{b_i},
(\phi^{\bm{H}}_{\bm{t}})_*\sigma_*\partial_{b_j}) \\
&= \omega(\sigma_*\partial_{b_i},\sigma_*\partial_{b_j})\\
&=\sigma^*\omega(\partial_{b_i},\partial_{b_j}), \end{align*}
which gives the term \(\beta\) in \(\Psi^*\omega\) as
claimed. This 2-form is still nondegenerate (each
\(\partial_{b_i}\) pairs nontrivially with the corresponding
\(\partial_{t_i}\)) so \(\Psi\) is still a submersion and an
immersion.\qedhere

\end{Proof}
\begin{Lemma}\label{lma:lag_section_modify}
In the situation of the previous lemma, if there is a 1-form
\(\eta\) on \(B\) with \(\sigma^*\omega = d\eta\) then there
is a Lagrangian section over \(B\).
\end{Lemma}
\begin{Proof}
If \(\tau(\bm{b}) = (\bm{b}, \bm{t}(\bm{b}))\) is another
section (written with respect to the coordinate system
\(\Psi\)) then we can compute \(\tau^*\omega\) by following
the calculation in Exercise \ref{exr:lag_section}. We get
\[\omega(\tau_*\partial_{b_i}, \tau_*\partial_{b_j}) =
\frac{\partial t_i}{\partial b_j} - \frac{\partial
t_j}{\partial b_i} + \beta(\partial_{b_i}, \partial_{b_j}).\]
By comparing with the formula for the exterior derivative of
the 1-form \(\sum t_i(\bm{b})db_i\), we see that
\(\tau^*\omega = d(\sum t_i(\bm{b})db_i) + \beta\). Now
suppose that \(\beta = d\eta\) for some 1-form \(\eta=\sum
\eta_i(\bm{b})db_i\). Taking \(t_i(\bm{b}) = -\eta_i(\bm{b})\)
we get a section for which \(\tau^*\omega = -\beta + \beta =
0\), i.e.\ a Lagrangian section.\qedhere

\end{Proof}
\begin{Corollary}\label{cor:extend_section_1}
If \(\bm{H}\colon X\to\RR^n\) is an integrable Hamiltonian
system with \(\bm{H}(X)=B\) and \(\sigma\) is a section over
\(B\) with\footnote{\(H^2_{dR}(B)\) denotes the De Rham
cohomology group of closed 2-forms modulo exact 2-forms;
\(H^2_{dR}(B)=0\) is a fancy way of saying ``if \(d\beta=0\)
then \(\beta=d\eta\)''.} \([\sigma^*\omega] = 0\in
H^2_{dR}(B)\) then \(\bm{H}\) admits a Lagrangian section over
\(B\). In fact, if \(\sigma\) is Lagrangian over a subset
\(B'\subset B\) and\footnote{\(H^2_{dR}(B,B')\) denotes the
relative De Rham cohomology. This is again closed forms modulo
exact forms, but where the forms \(\beta\) and \(\eta\) are
required to vanish on \(B'\). This is a slightly different
formulation to the standard setup in, say, the book by Bott
and Tu {\cite[p.78--79]{BottTu}} but equivalent to it (as
explained in the MathOverflow answer \cite{EbertMO} by
Ebert).} \([\sigma^*\omega] = 0\in H^2_{dR}(B,B')\) then
\(\bm{H}\) admits a Lagrangian section which agrees with
\(\sigma\) over \(B'\).
\end{Corollary}
\begin{Proof}
Note that \(\beta:=\sigma^*\omega\) is closed, so defines a de
Rham cohomology class. If \([\beta]=0\) in de Rham cohomology
then there exists a 1-form such that \(\beta = d\eta\). If
\(\beta = 0\) on \(B'\) then it defines a class in relative de
Rham cohomology \(H^2_{dR}(B,B')\), which vanishes if and only
if \(\beta = d\eta\) for a 1-form \(\eta\) which itself
vanishes on \(B'\). Inspecting the proof of Lemma
\ref{lma:liouville_nonlag}, this means that the Lagrangian
section built using \(\eta\) coincides with \(\sigma\) on
\(B'\).\qedhere

\end{Proof}
\begin{Remark}\label{rmk:rel_coh}
We will use the condition on relative
cohomology\index{relative cohomology|(} to find Lagrangian
sections for non-regular Lagrangian fibrations: we will first
construct Lagrangian sections near the critical fibres, then
extend them over the regular locus using this result,
providing the relevant relative cohomology group vanishes.

\end{Remark}
\begin{Corollary}\label{cor:extend_section_2}
Let \(f\colon X\to B\) be a regular Lagrangian fibration and
suppose \(\sigma\) is a section which is Lagrangian over a
(possibly empty) subset \(B'\subset B\). If \([\sigma^*\omega]
= 0 \in H^2_{dR}(B,B')\) then \(f\) admits a Lagrangian
section.
\end{Corollary}
\begin{Proof}
By the cohomological assumption, there exists a 1-form
\(\eta\) on \(B\) such that \(\eta=0\) on \(B'\) and \(d\eta =
\sigma^*\omega\). Cover \(B\) by integral affine coordinate
charts; the Lagrangian fibration is equivalent to an
integrable Hamiltonian system over each of these charts, and
we can apply Lemma \ref{lma:lag_section_modify} (using
\(\eta\)) to modify \(\sigma\) and obtain a Lagrangian
section. Since we are using the same 1-form on different
charts, we modify \(\sigma\) in the same way on overlaps
between charts, so we find a Lagrangian section over the whole
of \(B\).\qedhere

\end{Proof}
\begin{Remark}[Exercise \ref{exr:coh_vanishing}]\label{rmk:coh_vanishing}
We will later apply this when \(B\) is a punctured surface and
\(B'\) is a neighbourhood of a strict subset of the
punctures. This satisfies \(H^2(B,B')=0\).\index{relative
cohomology|)}

\end{Remark}
\section{Solutions to inline exercises}

\begin{Exercise}[Corollary \ref{cor:fibre_tori}]\label{exr:fibre_tori}
If \(f\colon X\to B\) is a proper submersion with connected
Lagrangian fibres then the fibres are Lagrangian
tori.
\end{Exercise}
\begin{Solution}
By Lemma \ref{lma:lagfibreg}, if we pick local coordinates
\((b_1,\ldots,b_n)\) on \(B\) then the functions \(b_1\circ
f,\ldots,b_n\circ f\) form an integrable Hamiltonian system,
so this follows from the little Arnold-Liouville theorem
(Theorem \ref{thm:littlearnoldliouville}). \qedhere

\end{Solution}
\begin{Exercise}\label{exr:non_exact}
If \(\omega\) is not an exact 2-form, how can we construct the
integral affine structure\index{integral affine!structure} on
\(B\)?
\end{Exercise}
\begin{Solution}
We need to define the \index{flux|(}flux map
\(\tilde{B}\to\RR^n\). As before, we fix the universal cover
\(p\colon\tilde{B}\to B\) and write \(\xi\to B\) for the local
system with fibre \(H_1(f^{-1}(b);\ZZ)\) over \(b\in B\). We
pick a \(\ZZ\)-basis of global sections \(c_1,\ldots,c_n\) of
\(\tilde{\xi}=p^*\xi\). Write \(\tilde{f}\colon p^*X\to
\tilde{B}\) for the pullback of \(f\) to the universal cover
(i.e.\ the Lagrangian fibration whose fibre over \(\tilde{b}\)
is \(f^{-1}(p(\tilde{b}))\)). We continue to write \(\omega\)
for the pullback of \(\omega\) to \(p^*X\).

Fix a basepoint \(\tilde{b}_0\in\tilde{B}\). Given a point
\(\tilde{b}\in\tilde{B}\), pick a path \(\gamma\colon[0,1]\to
\tilde{B}\) from \(\tilde{b}_0\) to \(\tilde{b}\). A {\em
family of loops over \(\gamma\)} (see Figure
\ref{fig:family_of_loops}) is a homotopy \(C\colon
S^1\times[0,1]\to p^*X\) satisfying
\(\tilde{f}(C(s,t))=\gamma(t)\), i.e.\ if \(t\) is fixed,
\(C(s,t)\) is a loop in \(\tilde{f}^{-1}(\gamma(t))\).

\begin{figure}[htb]
\begin{center}
\begin{tikzpicture}[scale=0.7]
\draw (0,0) circle [x radius = 2, y radius = 3];
\draw (80:1) arc [x radius = 0.5, y radius = 1, start angle = 80, end angle = 280];
\draw (100:0.9) arc [x radius = 0.5, y radius = 0.9, start angle = 100, end angle = -100];
\draw (0,-3) -- (8,-3);
\draw (0,-1) -- (8,-1);
\draw[dashed] (0,-1.95) circle [x radius = 0.3,y radius = 1.05];
\draw (0,-0.9) arc [x radius = 0.3,y radius = 1.05, start angle = 90, end angle = -90];
\begin{scope}[shift={(8,0)}]
\draw (0,0) circle [x radius = 2, y radius = 3];
\draw (80:1) arc [x radius = 0.5, y radius = 1, start angle = 80, end angle = 280];
\draw (100:0.9) arc [x radius = 0.5, y radius = 0.9, start angle = 100, end angle = -100];
\draw[dashed] (0,-1.95) circle [x radius = 0.3,y radius = 1.05];
\draw (0,-0.9) arc [x radius = 0.3,y radius = 1.05, start angle = 90, end angle = -90];
\end{scope}
\draw[dashed] (3,-2) circle [x radius = 0.3,y radius = 1];
\draw (3,-1) arc [x radius = 0.3,y radius = 1, start angle = 90, end angle = -90];
\draw[dashed] (5,-2) circle [x radius = 0.3,y radius = 1];
\draw (5,-1) arc [x radius = 0.3,y radius = 1, start angle = 90, end angle = -90];
\draw[->] (4,-3.2) -- (4,-6.3);
\node at (4,-4.5) [left] {\(\tilde{f}\)};
\begin{scope}[shift={(0,-6.5)}]
\draw (0,0) -- (8,0);
\node at (4,0) [below] {\(\gamma\)};
\node at (0,0) [left] {\(\tilde{b}_0\)};
\node at (8,0) [right] {\(\tilde{b}\)};
\end{scope}
\node at (4,-2) {\(C\)};
\end{tikzpicture}
\end{center}
\caption{A family of loops over \(\gamma\).}
\label{fig:family_of_loops}
\end{figure}
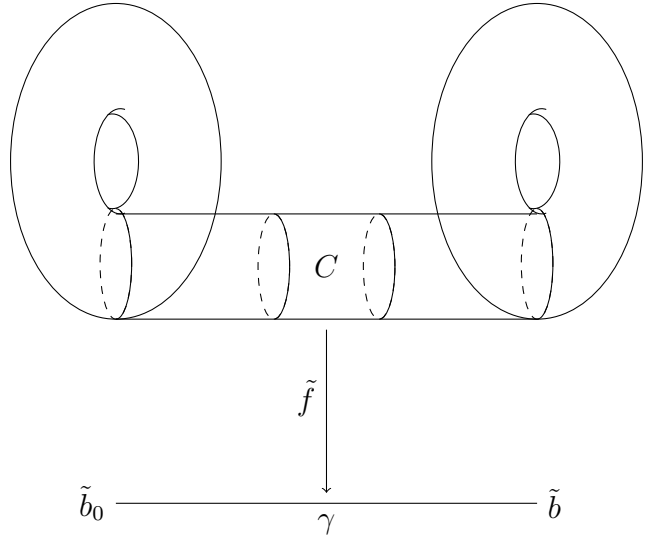

For \(k=1,\ldots,n\), pick a family of loops \(C_k\) over
\(\gamma\) with \(C_k(\cdot,t)\in c_k(\gamma(t))\) for all
\(t\in[0,1]\). Define \[\II(\tilde{b}) =
(I_1(\tilde{b}),\ldots,I_n(\tilde{b})),\qquad I_k(\tilde{b}) =
\int_{C_k}\omega.\]

It remains to understand how this flux map depends on the
choices we made, namely:
\begin{enumerate}
\item a basis \(c_1,\ldots,c_n\) of \(p^*\xi\),
\item a basepoint \(\tilde{b}_0\),
\item a path \(\gamma\) from \(\tilde{b}_0\) to \(\tilde{b}\),
\item a family of loops \(C_k\) over \(\gamma\) for each
\(k\in\{1,\ldots,n\}\).

\end{enumerate}
We deal first with the choice of \(\gamma\) and \(C_k\). Since
\(\tilde{B}\) is simply-connected, a different choice of path
\(\gamma'\) from \(\tilde{b}_0\) to \(\tilde{b}\) will be
homotopic to \(\gamma\) via some homotopy
\(h\colon[0,1]\times[0,1]\to \tilde{B}\). Choose \(C_k\) over
\(\gamma\) and \(C'_k\) over \(\gamma'\). We will show that
\(\int_{C_k}\omega=\int_{C'_k}\omega\).

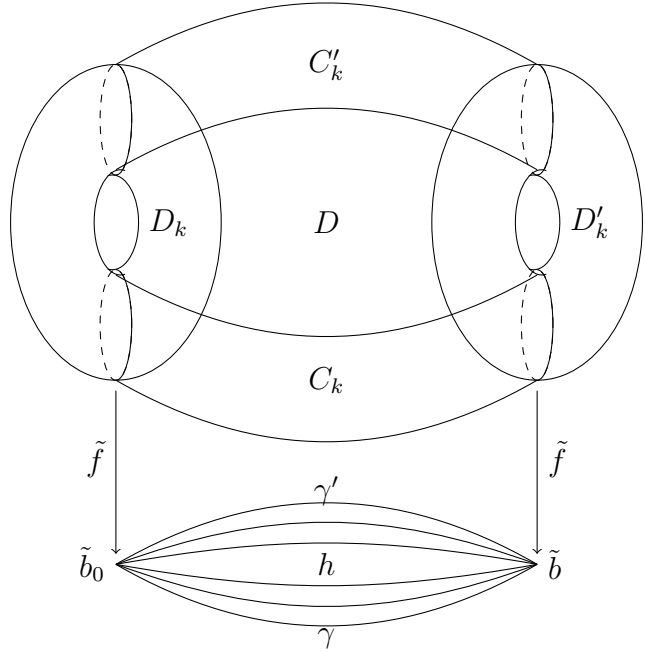
\begin{figure}[htb]
\begin{center}
\begin{tikzpicture}[scale=0.7]
\draw (0,0) circle [x radius = 2, y radius = 3];
\draw (80:1) arc [x radius = 0.5, y radius = 1,start angle = 80, end angle = 280];
\draw (100:0.9) arc [x radius = 0.5, y radius = 0.9,start angle = 100, end angle = -100];
\draw (0,-3) to[out=-30,in=-150] (8,-3);
\draw (0,-1) to[out=-30,in=-150] (8,-1);
\draw (0,3) to[out=30,in=150] (8,3);
\draw (0,1) to[out=30,in=150] (8,1);
\draw[dashed] (0,-1.95) circle [x radius = 0.3,y radius = 1.05];
\draw (0,-0.9) arc [x radius = 0.3,y radius = 1.05,start angle = 90, end angle = -90];
\draw[dashed] (0,1.95) circle [x radius = 0.3,y radius = 1.05];
\draw (0,0.9) arc [x radius = 0.3,y radius = 1.05,start angle = -90, end angle = 90];
\begin{scope}[shift={(8,0)}]
\draw (0,0) circle [x radius = 2, y radius = 3];
\draw (80:1) arc [x radius = 0.5, y radius = 1,start angle = 80, end angle = 280];
\draw (100:0.9) arc [x radius = 0.5, y radius = 0.9,start angle = 100, end angle = -100];
\draw[dashed] (0,-1.95) circle [x radius = 0.3,y radius = 1.05];
\draw (0,-0.9) arc [x radius = 0.3,y radius = 1.05,start angle = 90, end angle = -90];
\draw[dashed] (0,1.95) circle [x radius = 0.3,y radius = 1.05];
\draw (0,0.9) arc [x radius = 0.3,y radius = 1.05,start angle = -90, end angle = 90];
\end{scope}
\draw[->] (0,-3.2) -- (0,-6.3);
\draw[->] (8,-3.2) -- (8,-6.3);
\node at (0,-4.5) [left] {\(\tilde{f}\)};
\node at (8,-4.5) [right] {\(\tilde{f}\)};
\begin{scope}[shift={(0,-6.5)}]
\draw (0,0) to[out=30,in=150] (8,0);
\draw (0,0) to[out=20,in=160] (8,0);
\draw (0,0) to[out=10,in=170] (8,0);
\node at (4,0) {\(h\)};
\draw (0,0) to[out=-30,in=-150] (8,0);
\draw (0,0) to[out=-20,in=-160] (8,0);
\draw (0,0) to[out=-10,in=-170] (8,0);
\node at (4,1.4) {\(\gamma'\)};
\node at (4,-1.4) {\(\gamma\)};
\node at (0,0) [left] {\(\tilde{b}_0\)};
\node at (8,0) [right] {\(\tilde{b}\)};
\end{scope}
\node at (4,3) {\(C'_k\)};
\node at (4,-3) {\(C_k\)};
\node at (4,0) {\(D\)};
\node at (1,0) {\(D_k\)};
\node at (9,0) {\(D'_k\)};
\end{tikzpicture}
\end{center}
\caption{Different choices of paths and homotopies for the solution of Exercise \ref{exr:non_exact}.}
\label{fig:indep_of_C_k}
\end{figure}

The loops \(C_k(\cdot,0)\) and \(C'_k(\cdot,0)\) are
homologous in \(\tilde{f}^{-1}(\gamma(0))\) by assumption, and
therefore freely homotopic because \(\pi_1(T^n)\cong
H_1(T^n;\ZZ)\cong\ZZ^n\). Let \(D_k\colon S^1\times[0,1]\to
\tilde{f}^{-1}(\gamma(0))\) be a free homotopy with
\(D_k(\cdot,0)=C_k(\cdot,0)\) and
\(D_k(\cdot,1)=C'_k(\cdot,0)\). By the homotopy lifting
property of the submersion \(\tilde{f}\), we can find a map
\(D\colon S^1\times[0,1]\times [0,1]\to p^*X\) with
\(\tilde{f}\circ D = h\) and \(D(s,t,0)=D_k(s,t)\). Define
\(D'_k(s,t)=D(s,t,1)\); this defines a cylinder in the
Lagrangian torus \(\tilde{f}^{-1}(\gamma(1))\) (see Figure
\ref{fig:indep_of_C_k}). Consider \(D\) as a 3-chain (in the
sense of singular homology). Because \(d\omega=0\), we have
\[0=\int_D d\omega=\int_{\partial D} \omega\] by Stokes's
theorem. But \(\partial D = C_k + D'_k - C'_k - D_k\), so \[0
= \int_{C_k} \omega + \int_{D'_k} \omega - \int_{C'_k}
\omega - \int_{D_k} \omega.\] Since \(D_k\) and \(D'_k\) are
contained in Lagrangian fibres, the integrals \(\int_{D_k}
\omega\) and \(\int_{D'_k} \omega\) vanish, and we see that
\[\int_{C_k} \omega = \int_{C'_k} \omega\] as required.

If we change basepoint to \(\tilde{b}'_0\), we can choose a
path \(\beta\) from \(\tilde{b}'_0\) to \(\tilde{b}_0\) and
homotopies \(\Gamma_1,\ldots,\Gamma_n\) over \(\beta\). Given
another point \(\tilde{b}\), choose \(\gamma\) from
\(\tilde{b}_0\) to \(\tilde{b}\) and homotopies
\(C_1,\ldots,C_n\) over \(\gamma\) to define the flux map
\(\II(\tilde{b})\). We can then choose the concatenated path
\(\gamma\cdot\beta\) from \(\tilde{b}'_0\) to \(\tilde{b}\)
and the concatenated homotopies \(C_k\cdot\Gamma_k\) to define
the flux map \(\II'(\tilde{b})\). The resulting flux maps
differ by translation: \(\II'(\tilde{b})=\JJ+\II(\tilde{b})\)
with \[\JJ = (J_1,\ldots,J_n),\qquad J_k=\int_{\Gamma_k}
\omega.\]

Finally, if we change the basis of sections \(c_1,\ldots,c_n\)
by an element of \(GL(n,\ZZ)\) then the result is to apply a
\(\ZZ\)-linear transformation to the flux
map\index{flux|)}. The argument that proved Corollary
\ref{cor:zaffine} shows that the integral affine structure on
\(\tilde{B}\) descends to \(B\). \qedhere

\end{Solution}
\begin{Exercise}[Remark \ref{rmk:coh_vanishing}]\label{exr:coh_vanishing}
Suppose that \(B\) is a 2-dimensional surface with a nonempty
set of punctures and that \(B'\subset B\) is a collar
neighbourhood of a strict subset of the punctures. Then
\(H^2(B,B')=0\).
\end{Exercise}
\begin{Proof}
Note first that the second cohomology of a punctured surface
is zero (provided there is at least one puncture). We have
\(H^2(B,B')\cong H^2(B/B')\). The quotient \(B/B'\) is the
result of filling in a strict subset of the punctures, so is
homeomorphic to a surface with fewer (but still some)
punctures. Therefore \(H^2(B,B') = H^2(B/B') = 0\).\qedhere

\end{Proof}
\chapter{Global action-angle coordinates and torus actions}
\label{ch:glob}
\thispagestyle{cup}

\section{Hamiltonian torus actions}

One way of stating the Arnold-Liouville theorem is that, after a
suitable change of coordinates in the target, the
\(\RR^n\)-action generated by the Hamiltonian vector fields
\(V_{H_1},\ldots,V_{H_n}\) actually factors through a
\(T^n\)-action. In this chapter, we work backwards, assuming
that we have a globally-defined torus action, even on the
non-regular fibres, and see what kinds of critical points can
occur.

\begin{Definition}
Let \(\bm{H}\colon X\to\RR^n\) be an integrable Hamiltonian
system such that the Hamiltonian \(\RR^n\)-action
\(\phi^{\bm{H}}_{\bm{t}}\) factors through a Hamiltonian
\(T^n\)-action\index{Hamiltonian!torus action}, that is
\(\phi^{\bm{H}}_{\bm{t}}=\OP{id}\) for any
\(\bm{t}\in(2\pi\ZZ)^n\). Then we call \(\bm{H}\) the {\em
moment map}\index{moment map} for the torus action. It is
conventional to write \(\mu\) rather than \(\bm{H}\) for a
moment map, and we will do this wherever we want to emphasise
the existence of the torus action. We will call a symplectic
\(2n\)-manifold \(X\) a {\em toric
manifold}\index{toric manifold} if it admits a Hamiltonian
\(T^n\)-action.

\end{Definition}
We saw in Lemma \ref{lma:globcoords} that the image of a moment
map determines the Hamiltonian system completely up to fibred
symplectomorphism, at least if there are no critical points and
there is a global Lagrangian section. We therefore concentrate
on the image \(\mu(X)\) of the moment map, which we will call
the {\em moment image}\index{moment image} or {\em moment
polytope}\index{moment polytope}. The Atiyah-Guillemin-Sternberg
convexity theorem\index{convexity theorem}, discussed in Section
\ref{sct:AGSD} below, tells us that \(\mu(X)\) is indeed a
rational convex polytope\index{polytope!convex rational}. We
will not give a full proof of this theorem, as there are many
excellent expositions in the literature (e.g.\ Atiyah
{\cite[Theorem 1]{Atiyah}}, Audin \cite{Audin},
Guillemin-Sternberg {\cite[Theorem 4]{GS}}, McDuff-Salamon
{\cite[Theorem 5.47]{McDuffSalamon}}, amongst others). Instead,
we will prove the much easier Proposition
\ref{prp:straight_lines} below: that under mild conditions, the
boundary of the moment image is piecewise linear. This has the
advantage of being a local result, which will apply in
situations where we only have a torus action on some parts of
the manifold. In particular, it will apply in situations where
there is no sense in which the image of the Hamiltonian system
is convex, like the almost toric setting in Chapter
\ref{ch:almost_toric_manifolds}.

We need the following preliminary lemma:

\begin{Lemma}[Exercise \ref{exr:1_ps_ham}]\label{lma:1_ps_ham}
Let \(\mu\colon X\to\RR^n\) be the moment map of a Hamiltonian
\(T^n\)-action. If \(s\colon\RR^n\to\RR\) is a linear map,
\(s(b_1,\ldots,b_n)=\sum s_ib_i\) then \(s\circ\mu\) generates
the Hamiltonian flow \(\phi^{\mu}_{(s_1t_1,\ldots,s_nt_n)}\).

\end{Lemma}
In this case, the \(\RR\)-action (flow) \(\phi^{s\circ\mu}_{t}\)
can be thought of as a subgroup of the \(T^n\) action, coming
from the homomorphism \[s^T\colon\RR\to\RR^n/(2\pi\ZZ)^n,\qquad
s^T(t)=(s_1t,\ldots,s_nt).\]
Now suppose that \(X\) is a symplectic \(2n\)-manifold, and that
\(\mu\colon X\to\RR^n\) is the moment map for a Hamiltonian
\(T^n\)-action. Write
\(\partial\mu(X)\) for the boundary of the moment image. We will
assume that \(\partial\mu(X)\) is a piecewise smooth
hypersurface; we will show that, under mild assumptions,
\(\partial\mu(X)\) is piecewise {\em linear}. Pick local smooth
embeddings \(\delta_i \colon(0,1)^{n-1} \to \partial\mu(X)
\subset \RR^n\) parametrising the smooth pieces of
\(\partial\mu(X)\) and assume that there are smooth lifts
\(\gamma_i\colon(0,1)^{n-1}\to X\) such that \(\delta_i =
\mu\circ\gamma_i\).

\begin{Proposition}[Piecewise linearity of the toric boundary]\label{prp:straight_lines}
The\index{moment image!has piecewise linear boundary} image of
each \(\delta_i\) is contained in an affine hyperplane
\(\Pi_i\) with rational slopes, that is
\(\Pi_i=\{x\in\RR^n\,:\,\alpha\cdot x=c\}\) for some integer
vector \(\alpha\). If \(z\in\mu^{-1}(\delta_i)\) then the
stabiliser of \(z\) is precisely the 1-dimensional subtorus
\(s_i^T(\RR)\subset T^n\) where \(s_i(x)=\alpha\cdot x\).
\end{Proposition}
\begin{Proof}
Let \(\Pi_i(\bm{t})\) be the tangent hyperplane to
\(\delta_i\) at \(\delta_i(\bm{t})\), with normal vector
\(\alpha=(\alpha_1,\ldots,\alpha_n)\). We say that
\(\delta_i\) has rational slopes at \(\bm{b}\) if \(\alpha\) is
parallel to an integer vector. Otherwise, at least one of the
ratios \(\alpha_k/\alpha_\ell\) is irrational. We will show
that the tangent hyperplane to \(\Pi_i(\bm{t})\) has rational
slopes for all \(\bm{t}\), which is only possible if
\(\Pi_i(\bm{t})\) is independent of \(\bm{t}\) (otherwise the
slopes would need to take irrational values by the
intermediate value theorem). This will imply that \(\delta_i\)
coincides with its tangent hyperplane. The statement about
stabilisers will come up naturally in the proof.

Suppose that \(\delta_i\) has an irrational slope at
\(\bm{b}:=\delta_i(\bm{t})\). Pick a ball \(B\) centred at
\(\bm{b}\) and a smooth function \(S\colon B\to\RR\) such that
\(\mu(X)\cap B=\{\bm{p}\in B\,:\,S(\bm{p})\geq 0\}\) and
\(\partial\mu(X)\cap B=S^{-1}(0)\) is a regular level set. The
function \(S\circ\mu\) has a minimum along \(S^{-1}(0)\), so
if \(z\in\mu^{-1}(\bm{b})\) then \(d_{\bm{b}}S\circ
d_z\mu=0\). Let \(s:=d_{\bm{b}}S\) and consider the
Hamiltonian function \(H:=s\circ \mu\); by Lemma
\ref{lma:1_ps_ham}, this generates the \(\RR\)-action given by
\(s^T(\RR)\subset T^n\). But \(d_zH=s\circ d_z\mu=0\), so this
\(\RR\)-action fixes any point \(z\in \mu^{-1}(\bm{b})\). The
stabiliser of \(z\) is a closed subgroup of \(T^n\) containing
\(s^T(\RR)\subset T^n\); if \(\delta_i\) has an irrational
slope at \(\bm{b}\) then the closure of this subgroup is at
least 2-dimensional, so the stabiliser of \(z\) contains a
2-torus. This means there are two linearly independent
components of \(\mu\) whose Hamiltonian vector fields vanish
at \(z\); in particular, the rank of
\(d_{\gamma_i(\bm{t})}\mu\) is at most \(n-2\). Since
\(\delta_i=\mu\circ\gamma_i\), we have \(d_{\bm{t}}\delta_i =
d_{\gamma_i(\bm{t})}\mu\circ d_{\bm{t}}\gamma_i\), and this
means that the rank of \(d_{\bm{t}}\delta_i\) is at most
\(n-2\). This contradicts the assumption that \(\delta_i\) is
an embedding (\(\delta_i\) fails to be an immersion at
\(\bm{t}\)).

If \(\delta_i\) has rational slopes then we can take
\(S(x)=\alpha\cdot x\) as the function which is constant along
\(\delta_i\) and the same argument gives us the stabiliser as
claimed. \qedhere

\end{Proof}
\begin{Remark}
As this book progresses, we will allow our Lagrangian torus
fibrations \(f\colon X\to B\) to have more and more different
types of critical points. If \(B^{reg}\subset B\) denotes the
set of regular values of \(f\) then we know \(B^{reg}\)
inherits an integral affine structure. We can now allow \(f\)
to have ``toric critical points'', where \(X\) admits a local
Hamiltonian torus action having \(f\) as its moment
map. Proposition \ref{prp:straight_lines} tells us that \(B\)
will have the structure of an integral affine manifold with
piecewise linear boundary and corners, extending the integral
affine structure on \(B^{reg}\).

\end{Remark}
\section{Delzant polytopes and toric manifolds}
\label{sct:AGSD}

\begin{Definition}
A {\em rational convex polytope} \(P\) is a subset of
\(\RR^n\) defined as the intersection of a finite collection
of half-spaces \(S_{\alpha,b}=\{x\in\RR^n\ :\
\alpha_1x_1+\cdots+\alpha_nx_n\leq b\}\) with
\(\alpha_1,\ldots,\alpha_n\in\ZZ\) and \(b\in\RR^n\). We say
that \(P\) is a {\em Delzant\footnote{Audin \cite{Audin} calls
these {\em primitive polytopes}.}
polytope}\index{polytope!Delzant|(}
\index{Delzant!polytope|see {polytope, Delzant}}if it is a
convex rational polytope\index{polytope!convex rational} such
that every point on a \(k\)-dimensional facet has a
neighbourhood isomorphic (via an integral affine
transformation) to a neighbourhood of the origin in the
polytope \([0,\infty)^{n-k}\times\RR^k\). A vertex of a
polytope is called {\em Delzant}\index{Delzant!vertex} if the
germ of the polytope at that vertex is Delzant.

\end{Definition}
\begin{Example}
The polygon in Figure \ref{fig:non_delzant} fails to be Delzant:
there is no integral affine transformation sending the marked vertex
to the origin and sending the two marked edges to the \(x\)- and
\(y\)-axes. Indeed, the primitive integer vectors \((0,1)\) and
\((2,1)\) pointing along these edges span a strict sublattice of
the integer lattice \(\ZZ^2\).

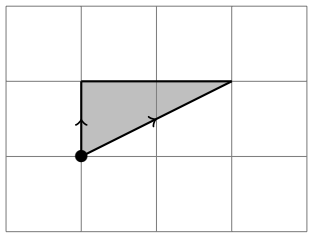
\begin{figure}[htb]
\begin{center}
\begin{tikzpicture}
\draw[step=1,gray,thin] (-1,-1) grid (3,2);
\filldraw[fill=gray,draw=none,opacity=0.5] (0,0) -- (0,1) -- (2,1) -- cycle;
\draw[thick] (0,1) -- (2,1);
\draw[->-,thick] (0,0) -- (0,1);
\draw[->-,thick] (0,0) -- (2,1);
\node at (0,0) {\(\bullet\)};
\end{tikzpicture}
\end{center}
\caption{A non-Delzant polygon.}
\label{fig:non_delzant}
\end{figure}

\end{Example}
\begin{Theorem}\label{thm:convexity}
Let \(X\) be a toric manifold, that is a symplectic
\(2n\)-manifold equipped with a Hamiltonian
\(T^n\)-action with moment map
\(\mu\colon X\to\RR^n\).
\begin{enumerate}
\item (Atiyah-Guillemin-Sternberg convexity theorem {\cite[Theorem
1]{Atiyah}}, {\cite[Theorem
4]{GS}})\index{Atiyah-Guillemin-Sternberg theorem|see {convexity theorem}}\index{convexity theorem}
The moment image \(\Delta:=\mu(X)\) is a Delzant
polytope. If \(X\) is compact, then \(\Delta\) is the convex
hull of \(\{\mu(x)\ :\ x\in Fix(X)\}\), where \(Fix(X)\) is
the set of fixed points of the torus action.
\item (Delzant existence theorem {\cite[Section 3]{Delzant}}) For
\index{Delzant!existence theorem} any compact Delzant
polytope\index{polytope!Delzant|)} \(\Delta\subset\RR^n\)
there exists a symplectic \(2n\)-manifold \(X_\Delta\) and a
map \(\mu\colon X_\Delta\to\RR^n\) with
\(\mu(X_\Delta)=\Delta\) such that \(\mu\) generates a
Hamiltonian \(T^n\)-action. Moreover, \(X_\Delta\) is a
projective variety. Such varieties
are often called {\em projective toric
varieties}\index{toric variety}.
\item (Delzant uniqueness theorem {\cite[Theorem 2.1]{Delzant}})
The\index{Delzant!uniqueness theorem} moment polytope
determines the pair \((X,\mu)\) up to fibred
symplectomorphism.

\end{enumerate}
\end{Theorem}
We will not prove (1) or (3). We will see two constructions of
\(X_\Delta\) later, proving (2). In the remainder of this
chapter, we will focus instead on examples where we can extract
geometric information about \(X\) from the moment polytope.

\section{Examples}

\begin{Example}\label{exm:cn}
Consider\index{moment polytope!standard torus action on
Cn@standard torus action on $\mathbb{C}^n$} the \(n\)-torus
action on \(\CC^n\) given by
\[(z_1,\ldots,z_n)\mapsto(e^{it_1}z_1,\ldots,e^{it_n}z_n).\]
This is Hamiltonian, with moment map \[\mu(z_1,\ldots,z_n) =
\left(\frac{1}{2}|z_1|^2, \ldots, \frac{1}{2}|z_n|^2\right).\]
The image of the moment map is the nonnegative orthant. This
is a manifold with boundary and corners: the \(\mu\)-preimage
of a boundary stratum of codimension \(k\) is an
\((n-k)\)-dimensional torus. For example, the preimage of the
vertex is a single fixed point (the origin), the preimage of a
point on the positive \(b_1\)-axis is a circle with fixed
radius in the \(z_1\)-plane, the preimage of a point on the
interior of the \(b_1b_2\)-plane is a 2-torus, and so forth.

\begin{center}
\begin{tikzpicture}[baseline=0]
\filldraw[fill=lightgray,opacity=0.5,draw=none] (2,0) -- (0,0) -- (0,2) -- (2,2) -- cycle;
\draw[thick,black] (2,0) -- (0,0) -- (0,2);
\node at (-0.5,0.5) {\(\CC^{2}\)};
\end{tikzpicture}
\qquad
\begin{tikzpicture}[baseline=0]
\filldraw[fill=lightgray,draw=gray,opacity=0.5] (1.5,-0.5) -- (3.3,-0.2) -- (1.8,0.3) -- (0,0) -- cycle;
\filldraw[fill=lightgray,draw=gray,opacity=0.5] (1.8,0.3) -- (0,0) -- (0,2) -- (1.8,2.3) -- cycle;
\filldraw[fill=lightgray,draw=gray,opacity=0.5] (1.5,-0.5) -- (0,0) -- (0,2) -- (1.5,1.5)-- cycle;
\filldraw[fill=lightgray,draw=gray,opacity=0.5] (1.5,-0.5) -- (3.3,-0.2) -- (3.3,1.8) -- (1.5,1.5) -- cycle;
\filldraw[fill=lightgray,draw=gray,opacity=0.5] (1.8,2.3) -- (0,2) -- (1.5,1.5) -- (3.3,1.8) -- cycle;
\filldraw[fill=lightgray,draw=gray,opacity=0.5] (1.8,0.3) -- (3.3,-0.2) -- (3.3,1.8) -- (1.8,2.3) -- cycle;
\draw[black,thick] (0,0) -- (0,2);
\draw[black,thick] (0,0) -- (1.5,-0.5);
\draw[black,dotted,thick] (0,0) -- (1.8,0.3);
\node at (3.7,0.5) {\(\CC^{3}\)};
\end{tikzpicture}
\end{center}

\end{Example}
\begin{Remark}
The critical values of \(\mu\) are precisely the boundary
points of the moment polytope. The boundary is stratified into
facets of dimension \(0\) (vertices), \(1\) (edges), \(2\)
(faces), etc, so we can classify the critical values according
to the dimension of the stratum to which they belong. By
definition, any Delzant polytope is locally isomorphic to
\(\RR^k\times[0,\infty)^{n-k}\) in a neighbourhood of a point
in a \(k\)-dimensional facet. In Example \ref{exm:cn}, we have
found a system whose moment image is \([0,\infty)^{n-k}\), so
by Theorem \ref{thm:convexity}(3), this means that the
integrable Hamiltonian system in a neighbourhood of a critical
point living over a \(k\)-dimensional facet is
fibred-symplectomorphic to the
system \begin{align*}\mu\colon\RR^k\times
(S^1)^k\times\CC^{n-k}&\to\RR^n,\\
\mu(\bm{p},\bm{q},z_{k+1},\ldots,z_n)&=
\left(\bm{p},\frac{1}{2}|z_{k+1}|^2,\ldots,\frac{1}{2}|z_n|^2\right).\end{align*}
Such critical points are called {\em
toric}\index{toric critical point|see {critical point,
toric}}\index{critical point!toric}\footnote{In fact, it is a
theorem of Eliasson \cite{Eliasson} and Dufour--Molino
\cite{DufourMolino} that toric critical points can be
characterised purely in terms of the Hessian of the
Hamiltonian system at the critical point. They call such
critical points {\em elliptic}.}\index{elliptic critical
point|see {critical point, toric}}\index{Eliasson normal form
theorem}\index{Hessian} and the set of all toric critical
points is often called the {\em toric
boundary}\index{toric boundary|(} of \(X\). It is not a
boundary in the usual sense: it is a union of submanifolds of
codimension 2. Instead, considering \(X\) as a projective
variety, it is the boundary in the sense of algebraic
geometry: it is a divisor, and is often called the {\em toric
divisor}.\index{divisor!toric|see {toric,
boundary}}\index{toric boundary|)}

\end{Remark}
Here is a nice way to understand the genus 1
Heegaard\index{Heegaard decomposition|(} decomposition of the
3-sphere using the moment map for \(\CC^2\).

\begin{figure}[htb]
\begin{center}
\begin{tikzpicture}[baseline=0]
\filldraw[fill=lightgray,draw=none,opacity=0.5] (2.5,0) -- (0,0) -- (0,2.5) -- (2.5,2.5) -- cycle;
\draw[black,thick] (2.5,0) -- (0,0) -- (0,2.5);
\draw[very thick] (2,0) -- (0,2);
\draw [decorate,decoration={brace,amplitude=7pt}] (0,2) -- (2,0) node [black,midway,xshift=10pt,yshift=10pt] {\(S\)};
\draw [decorate,decoration={brace,amplitude=7pt}] (2,0) -- (1,1) node [black,midway,xshift=-10pt,yshift=-8pt] {\(S_1\)};
\draw [decorate,decoration={brace,amplitude=7pt}] (1,1) -- (0,2) node [black,midway,xshift=-8pt,yshift=-10pt] {\(S_2\)};
\node at (1,1) {\(\bullet\)};
\node at (1,1) [below left] {\(T\)};
\node at (2,0) {\(\bullet\)};
\node at (0,2) {\(\bullet\)};
\node at (2,0) [below] {\(s_1\)};
\node at (0,2) [left] {\(s_2\)};
\end{tikzpicture}
\qquad
\begin{tikzpicture}[baseline=-1.3cm]
\draw (0,0) circle [radius=1];
\draw[shift={(-0.125,0)}] (90:0.25) arc (90:-90:0.25);
\draw[shift={(0.125,0)}] (70:0.35) arc (70:290:0.35);
\draw (0,0) circle [radius=0.7];
\node (si) at (2.5,0) {core circle \(s_i\)};
\draw[->] (si) -- (0:0.7);
\node at (2.15,-1.15) {torus fibre \(T\)};
\draw[->] (1,-1) -- (-45:1);
\node at (2.15,1.15) {solid torus \(S_i\)};
\draw[->] (1,1) -- (45:0.8);
\end{tikzpicture}
\end{center}
\caption{The unit sphere in \(\CC^2\) lives over the slanted line; the fibre \(T\) separates it into two solid tori.}
\label{fig:heegaard}
\end{figure}
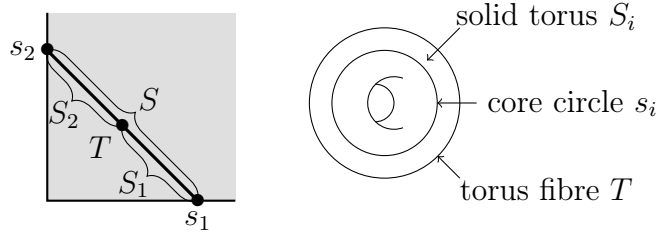

\begin{Example}[Heegaard decomposition of \(S^3\)]\label{exm:heegaard}
Let \(\mu\colon\CC^2\to\RR^2\) be the moment map from Example
\ref{exm:cn}. The preimage of the line segment
\(b_1+b_2=\frac{1}{2}\), \(b_1,b_2\geq 0\), is the subset
\(S:=\{(z_1,z_2)\in\CC^2\ :\ |z_1|^2+|z_2|^2=1\}\), that is
the unit 3-sphere; this is the slanted line segment in Figure
\ref{fig:heegaard}. The fibre
\(T:=\mu^{-1}\left(\tfrac{1}{4},\tfrac{1}{4}\right)\) is a
torus with \(T\subset S\). We can see from Figure
\ref{fig:heegaard} that \(T\) separates \(S\) into two pieces
\(S_1\), \(S_2\), and it is also easy to see that each piece
is homeomorphic to a {\em solid torus} \(S^1\times D^2\): the
``core circles'' of these solid tori are the fibres
\(s_1=\mu^{-1}\left(\tfrac{1}{2},0\right)\),
\(s_2=\mu^{-1}\left(0,\tfrac{1}{2}\right)\) over the points
where the line segment intersects the \(b_1\)- and
\(b_2\)-axes.\index{Heegaard decomposition|)}

\end{Example}
\begin{Example}[Exercise \ref{exr:2_sphere}]\label{exm:2_sphere}
Consider\index{moment polytope!S2@$S^2$} the unit 2-sphere
\((S^2,\omega)\) where \(\omega\) is the area form. By
comparing infinitesimal area elements, one can show that the
projection map from \(S^2\) to a circumscribed cylinder is
area-preserving\index{Archimedes's theorem on
areas}\footnote{If Cicero is to be believed
{\cite[XXIII–64,65]{Cicero}}, a diagram representing this
theorem was engraved on the tomb of Archimedes (who proved
it).}. Let \(\mu\colon S^2\to\RR\) be the height function
\(\mu(x,y,z)=z\) (thinking of \(S^2\) embedded in the standard
way in \(\RR^3\)). Then \(\mu\) is a moment map for the circle
action which rotates around the \(z\)-axis. The moment image
is \([-1,1]\subset\RR\).


\makeatletter
\pgfdeclareradialshading[tikz@ball]{ball}{\pgfqpoint{-20bp}{20bp}}{%
color(0bp)=(tikz@ball!0!white); color(17bp)=(tikz@ball!0!gray);
color(21bp)=(tikz@ball!70!black); color(25bp)=(black!70);
color(30bp)=(black!70)} \makeatother
\begin{figure}[htb]
\begin{center}
\begin{tikzpicture}
\begin{scope}[scale=2]
\draw (0,-1) circle [x radius = 1,y radius = 0.3];
\draw (-1,0) arc [x radius = 1,y radius = 0.3,start angle=180,end angle=360];
\draw[dotted] (-1,0) arc [x radius = 1,y radius = 0.3,start angle=180,end angle=0];
\begin{scope}[opacity=0.8]
\clip (0,0) circle (1);
\shade [ball color=gray] (-0.0,0) ellipse (1.2 and 1);
\end{scope}
\draw (0,1) circle [x radius = 1,y radius = 0.3];
\draw (-1,-1) -- (-1,1);
\draw (1,-1) -- (1,1);
\end{scope}
\draw[thick,->] (2.5,0) -- (4,0) node [midway,above] {\(\mu\)};
\draw (4.5,-2) -- (4.5,2);
\end{tikzpicture}
\end{center}
\end{figure}

\end{Example}
\begin{Example}\label{exm:s2rescaled}
If we take \(S^2\) with the area form \(\lambda\omega\) (where
\(\omega\) is the form giving area \(4\pi\)) then the rescaled
height function \(\lambda z\) is a moment map for the circle
action which rotates around the \(z\)-axis with period
\(2\pi\). The moment image is \([-\lambda,\lambda]\).

\end{Example}
\begin{Example}\label{exm:products_of_spheres}
One\index{moment polytope!product of spheres} can form more
examples by taking products. If we take
\(S^2\times\cdots\times S^2\) with the product symplectic
form\index{symplectic form!on product of spheres} giving the
\(i\)th factor symplectic area \(4\pi\lambda_i\) then we get a
\(T^n\)-action on \((S^2)^n\), whose moment map is
\(\mu((x_1,y_1,z_1), \ldots, (x_n,y_n,z_n)) =
(z_1,\ldots,z_n)\), with image the hypercuboid
\([-\lambda_1,\lambda_1]\times \cdots
\times[-\lambda_n,\lambda_n]\). For example, if we use equal
areas \(\lambda_1=\cdots=\lambda_n=1\) then the moment image
for \(S^2\times S^2\) is a square, whose vertices correspond
to the fixed points \(\{(0,0,\pm 1)\}\times\{(0,0,\pm 1)\}\),
and whose edges correspond to the spheres
\(S^2\times\{(0,0,\pm 1)\}\) and \(\{(0,0,\pm 1)\}\times
S^2\). For \(S^2\times S^2\times S^2\) the moment image is a
cube whose horizontal faces correspond to the submanifolds
\(S^2\times S^2\times\{(0,0,\pm 1)\}\), etc.

\begin{center}
\begin{tikzpicture}[baseline=0]
\filldraw[fill=gray,opacity=0.5,draw=none] (0,0) -- (0,2) -- (2,2) -- (2,0) -- (0,0);
\draw[black,thick] (0,0) -- (0,2) -- (2,2) -- (2,0) -- (0,0);
\node at (1,0) [below] {\(S^2\times S^2\)};
\begin{scope}[shift={(3,0)}]
\filldraw[fill=lightgray,draw=none,opacity=0.5] (1.5,-0.5) -- (3.3,-0.2) -- (1.8,0.3) -- (0,0) -- cycle;
\filldraw[fill=lightgray,draw=none,opacity=0.5] (1.8,0.3) -- (0,0) -- (0,2) -- (1.8,2.3) -- cycle;
\filldraw[fill=lightgray,draw=none,opacity=0.5] (1.5,-0.5) -- (0,0) -- (0,2) -- (1.5,1.5)-- cycle;
\filldraw[fill=lightgray,draw=none,opacity=0.5] (1.5,-0.5) -- (3.3,-0.2) -- (3.3,1.8) -- (1.5,1.5) -- cycle;
\filldraw[fill=lightgray,draw=none,opacity=0.5] (1.8,2.3) -- (0,2) -- (1.5,1.5) -- (3.3,1.8) -- cycle;
\filldraw[fill=lightgray,draw=none,opacity=0.5] (1.8,0.3) -- (3.3,-0.2) -- (3.3,1.8) -- (1.8,2.3) -- cycle;
\draw[black,thick,dotted] (1.5,-0.5) -- (3.3,-0.2) -- (1.8,0.3) -- (0,0) -- cycle;
\draw[black,thick,dotted] (1.8,0.3) -- (0,0) -- (0,2) -- (1.8,2.3) -- cycle;
\draw[black,thick] (1.5,-0.5) -- (0,0) -- (0,2) -- (1.5,1.5)-- cycle;
\draw[black,thick] (1.5,-0.5) -- (3.3,-0.2) -- (3.3,1.8) -- (1.5,1.5) -- cycle;
\draw[black,thick] (1.8,2.3) -- (0,2) -- (1.5,1.5) -- (3.3,1.8) -- cycle;
\draw[black,thick,dotted] (1.8,0.3) -- (3.3,-0.2) -- (3.3,1.8) -- (1.8,2.3) -- cycle;
\node at (1.75,-0.5) [below] {\(S^2\times S^2\times S^2\)};
\end{scope}
\end{tikzpicture}

\end{center}
\end{Example}
\begin{Definition}[Affine length]\label{dfn:affinelength}
If \(\ell\colon[0,L]\to\RR^n\) is a line segment of the form
\(\ell(t)=at+b\) with \(a\in\ZZ^n\) a primitive
vector\footnote{An integer vector \(a\) is called primitive if
it is a shortest integer vector on the line it spans, in other
words if \(\lambda a\in\ZZ^n\) implies \(|\lambda|\geq 1\).}
and \(b\in\RR^n\) then we say \(\ell\) is a {\em rational}
line segment and define the {\em affine length}\index{affine
length|(} of \(\ell\) to be \(L\).

\end{Definition}
\begin{Example}
Consider the triangle in Figure \ref{fig:non_delzant}. The
horizontal edge has affine length \(2\) and the other two edges both
have affine length \(1\).

\end{Example}
\begin{Lemma}\label{lma:sphere_area}
If \(\ell\colon[0,L]\to\RR^n\) is a rational line segment
whose image is an edge of the moment polytope then
\(\mu^{-1}(\ell([0,L]))\) is a symplectic sphere of symplectic
area\index{symplectic area} \(2\pi L\).
\end{Lemma}
\begin{Proof}
By Theorem \ref{thm:convexity}(3), the preimage of an edge is
determined up to fibred symplectomorphism by its moment image
\(\ell([0,L])\). By comparing with Example
\ref{exm:s2rescaled}, we see that the preimage of such an edge
is symplectomorphic to
\(\left(S^2,\frac{L\omega}{2}\right)\).\index{affine length|)}
\qedhere

\end{Proof}
\begin{Example}\label{exm:cpn}
Consider\index{moment polytope!CPn@$\mathbb{CP}^n$} the
complex projective \(n\)-space\index{projective space!complex}
\(\cp{n}\), with homogeneous coordinates
\([z_1:\cdots:z_{n+1}]\) (see Appendix
\ref{ch:complex_projective_spaces}). This has a torus action
\([z_1: \cdots: z_{n+1}] \mapsto [e^{it_1}z_1: \cdots:
e^{it_n}z_n: z_{n+1}]\) which is Hamiltonian, for the
Fubini-Study form\footnote{If you are not familiar with this
symplectic form, we will construct it in Example
\ref{exm:cpn_reduction}.} \(\omega\), with moment map
\[\mu([z_1:\cdots:z_{n+1}] = \left(\frac{|z_1|^2}{|z|^2},
\ldots, \frac{|z_n|^2}{|z|^2}\right),\] where
\(|z|^2=\sum_{i=1}^{n+1}|z_i|^2\). The moment image is the
simplex \[\{(b_1,\ldots,b_n)\in\RR^n\ :\ b_1,\ldots,b_n\geq
0,\ b_1+\cdots+b_n\leq 1\}.\] For example, \(\mu(\cp{2})\) and
\(\mu(\cp{3})\) are drawn in Figure \ref{fig:cp2_cp3}. In each
case, the {\em hyperplane at infinity}
\(\{[z_1:\cdots:z_n:0]\}\) projects via \(\mu\) to the facet
\(b_1+\cdots+b_n=1\) of the simplex.

\begin{figure}
\begin{center}
\begin{tikzpicture}[baseline=0]
\filldraw[fill=lightgray,opacity=0.5,draw=black,thick] (2,0) -- (0,0) -- (0,2) -- cycle;
\draw[thick] (2,0) -- (0,0) -- (0,2) -- cycle;
\node at (0.75,-0.5) {\(\cp{2}\)};
\begin{scope}[shift={(4,0.5)}]
\filldraw[fill=lightgray,draw=black,thick,opacity=0.5] (1.5,-0.5) -- (1.8,0.3) -- (0,0) -- cycle;
\filldraw[fill=lightgray,draw=black,thick,opacity=0.5] (1.8,0.3) -- (0,0) -- (0,2) -- cycle;
\filldraw[fill=lightgray,draw=black,thick,opacity=0.5] (1.5,-0.5) -- (0,0) -- (0,2) -- cycle;
\filldraw[fill=lightgray,draw=black,thick,opacity=0.5] (1.5,-0.5) -- (1.8,0.3) -- (0,2) -- cycle;
\node at (0.75,-1) {\(\cp{3}\)};
\end{scope}
\end{tikzpicture}
\end{center}
\caption{The moment polytopes for \(\cp{2}\) and \(\cp{3}\).}
\label{fig:cp2_cp3}
\end{figure}

\end{Example}
\begin{Example}\label{exm:taut_bun}
The\index{moment polytope!O-1@$\mathcal{O}(-1)$} {\em
tautological bundle}\index{tautological
bundle|(}\index{O-1@$\mathcal{O}(-1)$|(} over \(\cp{1}\) is
the variety
\[\mathcal{O}(-1):=\{(z_1,z_2,[z_3:z_4])\in\CC^2\times\cp{1}\ :\
z_1z_4=z_2z_3\}.\] This has a holomorphic projection
\(\pi\colon\mathcal{O}(-1)\to\cp{1}\),
\(\pi(z_1,z_2,[z_3:z_4])=[z_3:z_4]\), which exhibits it as the
total space of a holomorphic line bundle over \(\cp{1}\). This
is a fancy way of saying that \(\pi^{-1}([z_3:z_4])\) is a
complex line (specifically \(\{(z_1,z_2)\in\CC^2\ :\
z_1z_4=z_2z_3\}\subset\CC^2\)) for all
\([z_3:z_4]\in\cp{1}\). The symplectic form
\(\omega_{\CC^2}\oplus\omega_{\cp{1}}\) on
\(\CC^2\times\cp{1}\) pulls back to a symplectic
form\index{symplectic form!on O-1@on $\mathcal{O}(-1)$} on
\(\mathcal{O}(-1)\), with respect to which the following
\(T^2\)-action is Hamiltonian: \[(z_1,z_2,[z_3:z_4])\mapsto
(e^{it_1}z_1, e^{it_2}z_2, [e^{it_1}z_3: e^{it_2}z_4]).\] The
moment map is the sum of the moment maps for \(\CC^2\) and
\(\cp{1}\): \[\mu(z_1,z_2,[z_3:z_4]) =
\left(\frac{1}{2}|z_1|^2+\frac{|z_3|^2}{|z_3|^2+|z_4|^2},\quad
\frac{1}{2}|z_2|^2+\frac{|z_4|^2}{|z_3|^2+|z_4|^2}\right).\]
The image of the moment map is the subset in Figure
\ref{fig:taut_bun}:
\[\Delta_{\mathcal{O}(-1)}:=\left\{(b_1,b_2)\in\RR^2\ :\
b_1,b_2\geq 0,\ b_1+b_2\geq 1\right\}.\]

\begin{figure}
\begin{center}
\begin{tikzpicture}
\filldraw[fill=lightgray,opacity=0.5,draw=none] (0,3) -- (0,2) -- (2,0) -- (3,0) -- (3,3) -- cycle;
\draw[black,thick] (0,3) -- (0,2) -- (2,0) -- (3,0);
\node at (1,1) [below left] {\(\mu(\cp{1})\)};
\end{tikzpicture}
\end{center}
\caption{The moment polygon \(\Delta_{\mathcal{O}(-1)}\).}
\label{fig:taut_bun}
\end{figure}
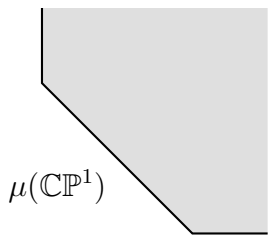

The zero-section
\(\cp{1}=\{z_1=z_2=0\}\subset\mathcal{O}(-1)\) projects down
to the edge \(b_1+b_2=1\). An alternative moment map can be
obtained by postcomposing with the integral affine
transformation\footnote{This is the first instance of the
notation mentioned in the preface: the angle bracket reminds
the reader that our matrix acts from the right. This will be
more important when the matrix appears in isolation.}
\[(b_1,b_2)\mapsto(b_1,b_2)\lmatrix 1 & 1 \\ 0 & 1
\rmatrix+(0,-1),\] which sends the moment polygon to
\[\{(b_1,b_2)\in\RR^2\ :\ b_1,b_2\geq 0,\ b_1-b_2\geq 1\}.\]

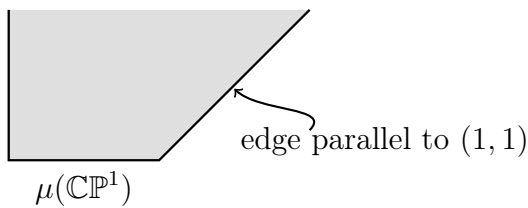
\begin{figure}
\begin{center}
\begin{tikzpicture}
\filldraw[fill=lightgray,opacity=0.5,draw=black,thick] (0,2) -- (0,0) -- (2,0) -- (4,2);
\draw[thick] (0,2) -- (0,0) -- (2,0) -- (4,2);
\node at (1,0) [below] {\(\mu(\cp{1})\)};
\node at (5,0.25) {edge parallel to \((1,1)\)};
\draw[thick,->] (4,0.4) to[out=45,in=-45] (3,0.95);
\end{tikzpicture}
\end{center}
\caption{Alternative moment polygon for \(\mathcal{O}(-1)\).}
\label{fig:taut_bun_2}
\end{figure}

\end{Example}
This is an important example because of the role played by
\(\mathcal{O}(-1)\) in birational geometry. The projection
\(\varpi\colon\mathcal{O}(-1)\to\CC^2\) given by
\(\varpi(z_1,z_2,[z_3:z_4])=(z_1,z_2)\) is a birational map
called the {\em blow-down}\index{blow-up} or {\em contraction}
of a \(-1\)-curve. It is an isomorphism away from
\((0,0)\in\CC^2\), but it contracts the sphere
\(\{(0,0,[z_3:z_4])\ :\ [z_3:z_4]\in\cp{1}\}\) (known as the
{\em exceptional sphere}) to the origin.

When we introduce the symplectic cut operation in Section
\ref{sct:symp_cut}, we will see that if we take a toric variety
\(X_\Delta\) and blow-up a fixed point of the torus action
(living over a vertex \(v\in\Delta\)), we get a new toric
variety \(X_{\Delta'}\) whose moment polytope \(\Delta'\)
differs from the previous one by truncating at the vertex
\(v\). More precisely, we use an integral affine transformation
to put \(\Delta\) in such a position that \(v\) sits at the
origin and \(\Delta\) is locally isomorphic to \([0,\infty)^n\)
near \(v\), then we truncate \(\Delta\) using the hyperplane
\(b_1+\cdots+b_n=c\) for some positive \(c\). Varying the
constant \(c\) will give different symplectic structures (in
particular, for \(n=2\), the symplectic area of the exceptional
sphere will vary).\index{tautological
bundle|)}\index{O-1@$\mathcal{O}(-1)$|)}

\begin{Example}\label{exm:o_n}
The\index{moment polytope!O-n@$\mathcal{O}(n)$} bundle
\(\mathcal{O}(-n)\)\index{O-n@$\mathcal{O}(n)$|(} over
\(\cp{1}\) is the variety\footnote{The discerning reader will
spot that this is the pullback of \(\mathcal{O}(-1)\) along
the degree \(n\) holomorphic map \(\cp{1}\to\cp{1}\),
\([z_3:z_4]\mapsto[z_3^n:z_4^n]\).}
\[\mathcal{O}(-n):=\{(z_1,z_2,[z_3:z_4])\in\CC^2\times\cp{1}\ :\
z_1z_4^n = z_2z_3^n\}\] The Hamiltonians
\[H_1=\frac{1}{2}|z_1|^2+\frac{|z_3|^2}{|z_3|^2+|z_4|^2},\qquad
H_2=\frac{1}{2}|z_2|^2+\frac{|z_4|^2}{|z_3|^2+|z_4|^2}\] still
generate circle actions, but the period lattice for the
\(\RR^2\)-action generated by \((H_1,H_2)\), while constant,
is no longer standard: the element
\(\phi^{H_1}_{2\pi/n}\phi^{H_2}_{2\pi/n}\) now acts as the
identity. This means that the period lattice is spanned by
\(\ZZ\vect{2\pi/n}{2\pi/n}\oplus\ZZ\vect{2\pi}{0}\). If we use
the combination \(\mu=\left(H_1,\frac{H_1+H_2}{n}\right)\)
then we get a standard period lattice, so this is a valid
moment map. This has the effect of applying the affine
transformation \(\lmatrix 1 & 1/n \\ 0 & 1/n\rmatrix\) to the
moment polygon in Figure \ref{fig:taut_bun}; we also translate
by \((0,-1/n)\) so that the horizontal edge \(\mu(\cp{1})\)
sits on the \(b_1\)-axis).

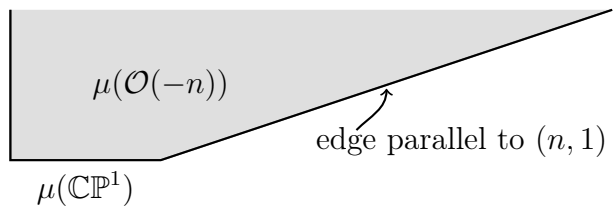
\begin{figure}[htb]
\begin{center}
\begin{tikzpicture}
\filldraw[fill=lightgray,opacity=0.5,draw=black,thick] (0,2) -- (0,0) -- (2,0) -- (8,2);
\draw[thick] (0,2) -- (0,0) -- (2,0) -- (8,2);
\node at (2,1) {\(\mu(\mathcal{O}(-n))\)};
\node at (1,0) [below] {\(\mu(\cp{1})\)};
\node at (6,0.25) {edge parallel to \((n,1)\)};
\draw[thick,->] (4.6,0.4) to[out=45,in=-90] (5,0.95);
\end{tikzpicture}
\caption{The moment polygon for \(\mathcal{O}(-n)\).}
\label{fig:o_n}
\end{center}
\end{figure}

Similarly, one can define the bundles
\(\mathcal{O}(n)\to\cp{1}\), \(n\geq 0\), and these admit
torus actions; the moment map now sends a neighbourhood of the
zero-section in \(\mathcal{O}(n)\) to the region shown in
Figure \ref{fig:o_n_pos}. For example, a complex line in
\(\cp{2}\) has normal bundle \(\mathcal{O}(1)\), and in the
moment image of \(\cp{2}\) we see precisely the \(n=1\)
neighbourhood surrounding the \(b_1\)-axis.

\begin{figure}[htb]
\begin{center}
\begin{tikzpicture}
\filldraw[fill=lightgray,opacity=0.5,draw=black,thick] (0,2) -- (0,0) -- (5,0) -- (1,2);
\draw[thick] (0,2) -- (0,0) -- (5,0) -- (1,2);
\node at (1.5,1) {\(\mu(\mathcal{O}(n))\)};
\node at (2.5,0) [below] {\(\mu(\cp{1})\)};
\node at (5,1.75) {edge parallel to \((-n,1)\)};
\draw[thick,->] (3.6,1.5) to[out=-90,in=45] (2.5,1.3);
\end{tikzpicture}
\caption{\(\mathcal{O}(n)\) for \(n\geq 0\).}
\label{fig:o_n_pos}
\end{center}
\end{figure}
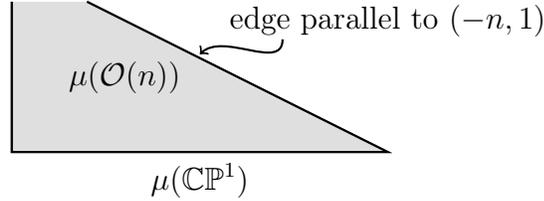

\end{Example}
The following lemma now follows immediately from these examples
and Theorem \ref{thm:convexity}(3).

\begin{Lemma}\label{lma:selfint}
Let\index{self-intersection!of curve!from moment polygon}
\(\Delta\subset\RR^2\) be a moment polygon and
\(e\subset\Delta\) an edge connecting two vertices
\(P,Q\). Assume that this edge is traversed from \(P\) to
\(Q\) as you move anticlockwise around the boundary of
\(\Delta\). Let \(v,w\) be primitive integer vectors pointing
along the other edges emerging from \(P\) and \(Q\)
respectively. Then a neighbourhood of \(\mu^{-1}(e)\) in
\(X_\Delta\) is symplectomorphic to a neighbourhood of the
zero-section in \(\mathcal{O}(n)\) where \(n=\det M\) where
\(M\) is the matrix with rows \(v,w\) (you may also see \(\det
M\) written as \(v\wedge w\)).

\end{Lemma}
\begin{center}
\begin{tikzpicture}
\filldraw[draw=none,fill=lightgray,opacity=0.5] (-1,1) -- (0,0) -- (2,-0.2) -- (3,1.2) -- cycle;
\draw[thick,->] (0,0) -- (-1,1) node [midway,left] {\(v\)};
\draw[thick,->] (2,-0.2) -- (3,1.2) node [midway,right] {\(w\)};
\draw[thick] (0,0) -- (2,-0.2) node [midway,below] {\(e\)};
\node at (1,0.5) {\(\Delta\)};
\node at (0,0) [below left] {\(P\)};
\node at (2,-0.2) [below right] {\(Q\)};

\end{tikzpicture}
\end{center}
\begin{Proof}
This is easily checked for the local models discussed above,
and any edge is integral affine equivalent to one of these
local models. It is therefore enough to check that \(v\wedge
w\) is unchanged by an integral affine transformation. The
determinant is unchanged by orientation-preserving integral
affine transformations. An orientation-reversing
transformation will switch the sign of \(v\wedge w\), but also
switch the order to \(w\wedge v\) because it switches
anticlockwise to clockwise, so these sign effects will
cancel.\index{O-n@$\mathcal{O}(n)$|)} \qedhere

\end{Proof}
\section{Non-Delzant polytopes}

\begin{Example}\label{exm:cycquot}
Consider\index{moment polytope!cyclic quotient singularity}
the group of \(n\)th roots of unity \(\bm{\mu}_n\) acting on
\(\CC^2\) via \((z_1,z_2)\mapsto (\mu z_1,\mu^a z_2)\) where
\(\gcd(a,n)=1\). Let \(X=\CC^2/\bm{\mu}_n\) be the quotient by
this group action. This is a symplectic {\em
orbifold}\index{orbifold}: the origin is a singular point. We
call this kind of singularity a {\em cyclic quotient
singularity of type
\(\frac{1}{n}(1,a)\)}.\index{singularity!cyclic
quotient}\index{singularity!1n1a@$\frac{1}{n}(1,a)$|see
{singularity, cyclic quotient}}

Hamiltonian flows still make perfect sense on \(X\) provided
they fix the origin. Consider the Hamiltonians
\(H_1=\frac{1}{2}|z_1|^2\) and \(H_2=\frac{1}{2}|z_2|^2\);
these are invariant under the action of \(\bm{\mu}_n\) and
hence define functions on \(X\). The flow is simply
\((e^{it_1}z_1,e^{it_2}z_2)\). However, the period lattice is
no longer standard\index{period lattice!standard}; we have
\(\phi_{H_1}^{2\pi/n}\phi_{H_2}^{2\pi a/n}=\OP{id}\). If
instead we use the Hamiltonians
\[\left(H_2,\frac{1}{n}(H_1+aH_2)\right)\] then the lattice of
periods becomes standard. The moment image is a convex wedge
in the plane bounded by the rays emanating from the origin in
directions \((0,1)\) and \((n,a)\); we will denote this
noncompact polygon by \(\pi(n,a)\)\index{pi na@$\pi(n,a)$}:

\begin{center}
\begin{tikzpicture}
\filldraw[fill=lightgray,opacity=0.5] (0,2) -- (0,0) -- (5,2);
\draw[thick,black,->] (0,0) -- (0,2) node [left] {\((0,1)\)};
\draw[thick,black,->] (0,0) -- (5,2) node [right] {\((n,a)\)};
\node at (1,1) {\(\pi(n,a)\)};
\draw[dashed] (0,1.5) -- (5*1.5/2,1.5);
\node at (2.5,1.5) [above] {\(\ell\)};

\end{tikzpicture}
\end{center}
This polygon is not Delzant\index{polytope!non-Delzant} at the
origin, corresponding to the fact that \(X\) is not smooth at
the origin.

\end{Example}
\begin{Remark}
The {\em link}\index{singularity!link} of a singularity is the
boundary of a small Euclidean neighbourhood of the singular
point. In this example, the link of the
\(\frac{1}{n}(1,a)\)-singularity is the preimage of a
horizontal line segment \(\ell\) running across
\(\pi(n,a)\). As in Example \ref{exm:heegaard}, this has a
decomposition as a union of two solid tori; this means it is a
{\em lens space}\index{lens space}.\index{Lna@$L(n,a)$|see
{lens space}} By definition, this is the lens space
\(L(n,a)\).

\end{Remark}
\begin{Lemma}[Exercise \ref{exr:lens_diffeo}]\label{lma:lens_diffeo}
The lens space \(L(n,a+kn)\) is diffeomorphic to \(L(n,a)\)
for all integers \(k\). The lens space\index{lens
space!diffeomorphism between} \(L(n,a)\) is diffeomorphic to
\(L(n,\bar{a})\) where \(a\bar{a} = 1\mod n\).

\end{Lemma}
\section{Solutions to inline exercises}

\begin{Exercise}[Lemma \ref{lma:1_ps_ham}]\label{exr:1_ps_ham}
Let \(\mu\colon X\to\RR^n\) be the moment map of a Hamiltonian
\(T^n\)-action. If \(s\colon\RR^n\to\RR\) is a linear map,
\(s(b_1,\ldots,b_n)=\sum s_ib_i\) then \(s\circ\mu\) generates
the Hamiltonian flow \(\phi^{\mu}_{(s_1t_1,\ldots,s_nt_n)}\).
\end{Exercise}
\begin{Proof}
We have \(\iota_{V_{\sum_i s_i\mu_i}}\omega = -d\left(\sum_i
s_i\mu_i\right) = -\sum_i s_id\mu_i =
\sum_i\iota_{s_iV_{\mu_i}}\omega\), and \(\sum s_iV_{\mu_i}\)
generates the flow \(\Phi^\mu_{(s_1t,\ldots,s_nt)}\). \qedhere

\end{Proof}
\begin{Exercise}[Example \ref{exm:2_sphere}]\label{exr:2_sphere}
Consider the unit 2-sphere \((S^2,\omega)\) where \(\omega\)
is the area form. By comparing infinitesimal area elements,
show that the projection map from \(S^2\) to a circumscribed
cylinder is area-preserving. Let \(H\colon S^2\to\RR\) be the
height function \(H(x,y,z)=z\) (thinking of \(S^2\) embedded
in the standard way in \(\RR^3\)). Show that \(H\) is an
action coordinate.
\end{Exercise}
\begin{Solution}
Let \(\hat{n}=(x,y,z)\) be the unit normal vector field to the unit
sphere. The area element on the unit sphere is given by
\(\sigma=\iota_{\hat{n}}(dx\wedge dy\wedge dz)=x dy\wedge dz + y
dz\wedge dx + z dx\wedge dy\). If we use cylindrical coordinates
\(x=r\cos\theta\), \(y=r\sin\theta\) then \(dx=\cos\theta
dr-r\sin\theta d\theta\) and \(dy=\sin\theta dr+r\cos\theta
d\theta\), so (after some algebra): \[\sigma = r^2d\theta\wedge
dz+rz dr\wedge d\theta.\] The unit sphere is defined by the equation
\(r^2+z^2=1\), which means that \(rdr=-zdz\) on the
sphere. Therefore \[\sigma=(1-z^2)d\theta\wedge dz-z^2dz\wedge
d\theta=d\theta\wedge dz.\] The unit cylinder has area element
\(\tau=d\theta\wedge dz\). The projection map from the sphere to the
cylinder is \(p(r,\theta,z)=(1,\theta,z)\), so
\(p^*\tau=d\theta\wedge dz=\sigma\).

Observe that the Hamiltonian \(H(x,y,z)=z\) gives the
Hamiltonian vector field \(\partial_\theta\), which rotates
the sphere with constant speed so that all orbits have period
\(2\pi\). Therefore \(H\) is an action coordinate (with angle
coordinate \(\theta\)). \qedhere

\end{Solution}
\begin{Exercise}[Lemma \ref{lma:lens_diffeo}]\label{exr:lens_diffeo}
The lens space \(L(n,a+kn)\) is diffeomorphic to \(L(n,a)\)
for all integers \(k\). The lens space \(L(n,a)\) is
diffeomorphic to \(L(n,\bar{a})\) where \(a\bar{a}=-1\mod n\).
\end{Exercise}
\begin{Solution}
Let \(X\) be the \(\frac{1}{n}(1,a)\) singularity and
\(\bm{H}\colon X\to\RR^2\) be the moment map from Example
\ref{exm:cycquot} with image \(\pi(n,a)\). Recall that the
lens space \(L(n,a)\) is the preimage under \(\bm{H}\) of the
horizontal line segment \(\ell\) shown in Figure
\ref{fig:lens_space_exr_1}.

\begin{figure}[htb]
\begin{center}
\begin{tikzpicture}
\filldraw[fill=lightgray,opacity=0.5] (0,2) -- (0,0) -- (5,2);
\draw[thick,black,->] (0,0) -- (0,2) node [left] {\((0,1)\)};
\draw[thick,black,->] (0,0) -- (5,2) node [right] {\((n,a)\)};
\node at (1,1) {\(\pi(n,a)\)};
\draw[dashed] (0,1.5) -- (5*1.5/2,1.5);
\node at (3,1.5) [above] {\(\ell\)};
\end{tikzpicture}
\end{center}
\caption{}
\label{fig:lens_space_exr_1}
\end{figure}

Let \(X'\) be the cyclic quotient singularity
\(\frac{1}{n}(1,a+kn)\), whose moment image is \(\pi(n,a+kn)\)
shown in Figure \ref{fig:lens_space_exr_2}. The integral
affine transformation \(M=\lmatrix 1&k\\0&1 \rmatrix\) relates
these moment polygons: \(\pi(n,a)M=\pi(n,a+kn)\).

\begin{figure}[htb]
\begin{center}
\begin{tikzpicture}
\filldraw[fill=lightgray,opacity=0.5] (0,7) -- (0,0) -- (5,7);
\draw[thick,black,->] (0,0) -- (0,2) node [left] {\((0,1)\)};
\draw[thick,black,->] (0,0) -- (5,7) node [right] {\((n,a+kn)\)};
\draw[dashed] (0,1.5) -- (5*1.5/2,1.5+3.75);
\draw[dashed] (0,1.5) -- (15/14,1.5);
\node at (5*1.5/4,5) {\(\pi(n,a+kn)\)};
\node at (6/14,1.5) [below] {\(\ell'\)};
\node at (5*1.5/4,3.5) [left] {\(\ell M\)};
\end{tikzpicture}
\end{center}
\caption{}
\label{fig:lens_space_exr_2}
\end{figure}

Since the moment polygons are related by \(M\), Lemma
\ref{lma:globcoords} gives us a fibred symplectomorphism
\(X'\to X\). The image of \(L(n,a)\) under this fibred
symplectomorphism lives over the (now slanted) line \(\ell
M\). We can isotope \(\ell M\) until it is a horizontal
segment \(\ell'\). The preimages are isotopic, and hence
diffeomorphic. The preimage of \(\ell'\) is \(L(n,a+kn)\) by
definition. Thus \(L(n,a)\cong L(n,a+kn)\).

If \(a\bar{a}=1\mod n\) then \(a\bar{a}+tn=1\) for some
\(t\). Let \(N=\lmatrix -a & t\\ n & \bar{a}\rmatrix \). We
have \(\pi(n,a)N=\pi(n,\bar{a})\), which in turn shows that
the associated lens spaces \(L(n,a)\) and \(L(n,\bar{a})\) are
diffeomorphic via the fibred symplectomorphism associated to
the integral affine transformation \(N\). If this seems like
magic, the trick to finding \(N\) is first to reflect
\(\pi(n,a)\) in the \(y\)-axis to get the wedge \(\pi(-n,a)\),
and then hunt for a matrix in \(SL(2,\ZZ)\) which sends
\((-n,a)\) to \((0,1)\). The composite is then \(N\in
GL(2,\ZZ)\). \qedhere

\end{Solution}
\chapter{Symplectic reduction}
\label{ch:symp_cut}
\thispagestyle{cup}

We now introduce {\em symplectic reduction}\index{symplectic
reduction|(}\index{symplectic quotient|see {symplectic
reduction}}, an operation which allows us to construct many
interesting symplectic manifolds. A special case of this is {\em
symplectic cut}, which you will use in Exercise
\ref{exr:toric_from_cut} to construct all toric manifolds.

\section{Symplectic reduction}

\begin{Definition}\label{dfn:ham_circ}
Let \((X,\omega)\) be a symplectic manifold and let \(H\colon
X\to\RR\) be a Hamiltonian. Suppose that
\(\phi^H_{2\pi}(x)=x\) for all \(x\in X\). Then the flow
defines an action of the circle \(S^1 = \RR/2\pi\ZZ\) on
\(X\). We this a {\em Hamiltonian circle
action}\index{Hamiltonian!circle action|(}. We will write
\(M_c:=H^{-1}(c)\) for the level sets of \(c\).

\end{Definition}
\begin{Remark}
Recall that a group action is called {\em
effective}\index{effective} if the only group element which
acts as the identity is the identity, and {\em
free}\index{free group action} if every point has trivial stabiliser. The
quotient of a manifold by a free circle action is again a
manifold. If all stabilisers are finite then the quotient is
an orbifold\index{orbifold|(}.

\end{Remark}
Here are some of the key facts about Hamiltonian circle actions.

\begin{Lemma}[Exercise \ref{exr:easy_properties}]\label{lma:easy_properties}
Suppose \(H\colon X\to\RR\) generates a circle action.
\begin{enumerate}
\item [(a)] The critical points of \(H\) are precisely the fixed
points of the circle action.
\item [(b)] The level sets \(M_c\) are preserved by the circle action.
\item [(c)] If \(x\) is a regular point then \(T_xM_c=\ker(dH)\).
\item [(d)] If \(v\in T_xM_c\) satisfies \(\omega(v,w)=0\) for all \(w\in
T_xM_c\) then \(v\in\OP{span}(V_H)\).

\end{enumerate}
\end{Lemma}
By Lemma \ref{lma:easy_properties}(a), the stabiliser of a
Hamiltonian circle action at a critical point of \(H\) is the
whole circle. Since any smooth function on a compact manifold
has critical points, Hamiltonian circle actions on compact
manifolds are never free. For this reason, we restrict attention
to a regular level set.

\begin{Lemma}
If \(c\) is a regular value of \(H\) and \(M_c=H^{-1}(c)\) is
the regular level set over \(c\) then the quotient \(Q_c :=
M_c/S^1\) is an orbifold\index{orbifold|)}.
\end{Lemma}
\begin{Proof}
By Lemma \ref{lma:easy_properties}(b), the level set is
\(M_c\) is preserved by the circle action, so the quotient
makes sense. We need to show that the stabiliser of the circle
action at a point \(x\in M_c\) is finite. Since \(x\) is
regular, it is not a critical point of \(H\), so by Lemma
\ref{lma:easy_properties}(a), the stabiliser at \(x\) is a
proper subgroup of \(S^1\). The circle is compact, and
stabilisers are closed subgroups, hence compact. The only
proper compact subgroups of \(S^1\) are finite. \qedhere

\end{Proof}
\begin{Lemma}[Symplectic reduction]\label{lma:symp_reduction}
Suppose \(H\colon X\to\RR\) generates a circle action and
\(c\) is a regular value of \(H\). Suppose for simplicity that
the action on \(M_c\) is free. Write \(i\colon M_c\to X\) and
\(p\colon M_c\to Q_c:=M_c/S^1\) for the inclusion and quotient
maps respectively. There is a unique symplectic
form\index{symplectic form!on symplectic reduction} \(\sigma\)
on \(Q_c\) such that \(i^*\omega = p^*\sigma\). We call
\((Q_c,\sigma)\) the {\em symplectic quotient} or {\em
symplectic reduction} of \(X\) by the Hamiltonian circle
action {\em at level \(c\)}.

\end{Lemma}
\begin{Remark}
One can drop the assumption that the action is free at the
cost of allowing quotients which are orbifolds.

\end{Remark}
\begin{Proof}
Suppose that \(v,w\in T_xM_c\) are tangent vectors to the
level set. We want to show that \(\omega(v,w)\) depends only
on \(p_*v\) and \(p_*w\in T_{p(x)}Q_c\), so that
\(\omega(v,w)=\sigma(p_*v,p_*w)\) for some 2-form \(\sigma\)
on \(Q_c\). In other words, we want to show that if \(x'\in
M_c\) is another point with \(p(x')=p(x)\) and \(v',w'\) are
vectors in \(T_{x'}M_c\) with \(p_*v' = p_*v\) and \(p_*w' =
p_*w\) then \(\omega(v',w')=\omega(v,w)\).

Since \(p(x')=p(x)\), we have \(\phi^H_t(x)=x'\) for some
\(t\in S^1\). The vectors \(v'\), \((\phi^H_t)_*v\), \(w'\)
and \((\phi^H_t)_*w\) all live in \(T_{x'}M_c\) so we can
add/subtract them. Since \(p\circ\phi^H_t=p\), we have
\(p_*\circ(\phi^H_t)_*=p_*\), so \[p_*(v' - (\phi^H_t)_*v) =
p_*v'-p_*v=0.\] Similarly \(p_*(w'-(\phi^H_t)_*w)=0\). Since
the kernel of \(p_*\) is spanned by \(V_H\), we have \[v' =
(\phi^H_t)_*v + aV_H,\qquad w' = (\phi^H_t)_*w + bV_H\] for
some \(a,b\in\RR\). Thus
\begin{align*}
\omega(v',w')&=\omega\left((\phi^H_t)_*v+aV_H,(\phi^H_t)_*w+bV_H\right)\\
&=\omega(v,w)+a\, dH(v)-b\, dH(w)=\omega(v,w),
\end{align*}
where we have used the fact that \(v\) and \(w\) are tangent
to a level set, so are annihilated by \(dH\). This shows the
existence of \(\sigma\).

We will now show that \(\sigma\) is nondegenerate. Given a
nonzero vector \(u\in T_{p(x)}Q_c\), pick a vector \(v\in
T_xM_c\) with \(p_*v=u\). This is possible since \(p\) is a
submersion. Since the projection \(p_*v=u\) is nonzero, \(v\)
is not a multiple of \(V_H\). By Lemma
\ref{lma:easy_properties}(d), there exists \(w\in T_xM_c\)
such that \(\omega(v,w)\neq 0\). Therefore
\(\sigma(u,p_*w)\neq 0\), showing that \(\sigma\) is
nondegenerate.

The fact that \(\sigma\) is closed follows from Lemma
\ref{lma:closed_submersion} below applied to
\(d\sigma\).\qedhere

\end{Proof}
\begin{Lemma}[Exercise \ref{exr:closed_submersion}]\label{lma:closed_submersion}
If \(p\colon M\to Q\) is a submersion and \(\eta\) is a \(k\)-form
on \(Q\) such that \(p^*\eta=0\) then \(\eta=0\).

\end{Lemma}
We finish this section by proving a lemma that will help us to
construct Hamiltonian circle or torus actions on symplectic
reductions.

\begin{Lemma}\label{lma:symmetries_descend}
Suppose \(G\colon X\to\RR\) Poisson-commutes with \(H\) then:
\begin{enumerate}
\item [(a)] \(G|_{H^{-1}(c)}\) descends to a function
\(\overline{G}\colon H^{-1}(c)/S^1\to\RR\).
\item [(b)] The Hamiltonian vector field \(V_{\overline{G}}\) is equal
to\footnote{Here, we are abusively identifying \(V_G\) with
its restriction to \(H^{-1}(c)\).} \(p_*V_G\).
\item [(c)] If \(G\) generates a Hamiltonian circle action on
\(H^{-1}(c)\) then \(\overline{G}\) generates a circle
action on the symplectic quotient.
\end{enumerate}
\end{Lemma}
\begin{Proof}
In this proof, we will write \(i\colon H^{-1}(c)\to X\) for
the inclusion of the \(c\)-level set.

(a) Since \(G\) Poisson-commutes with \(H\), Lemma
\ref{lma:heisenberg} implies it is constant along
\(H\)-orbits, and hence descends to a function
\(\overline{G}\) on the quotient, i.e.\, \(i^*G =
p^*\overline{G}\).

(b) Since \(G\) Poisson-commutes with \(H\), the vector field
\(V_G\) is tangent to \(H^{-1}(c)\). This means that the
restriction \(v := V_G|_{H^{-1}(c)}\) makes sense as a vector
field on \(H^{-1}(c)\) and \(V_G = i_*v\). We want to show
that \(\iota_{p_*v}\sigma = -d\overline{G}\).

We know that \(\iota_{V_G}\omega = -dG\). Pulling back via
\(i\) gives \(\iota_vi^*\omega = -di^*G\), or
\(\iota_vp^*\sigma = -dp^*\overline{G}\). This implies that
\(p^*\iota_{p_*v}\sigma = p^*(-d\overline{G})\). Since \(p\)
is a submersion, Lemma \ref{lma:closed_submersion} implies
that \(\iota_{p_*v}\sigma = -d\overline{G}\) as required.

(c) Part (b) implies that \(\phi^{\overline{G}}_t(p(x)) =
p(\phi^G_t(x))\), so if \(\phi^G_{2\pi}=\OP{id}\) then
\(\phi^{\overline{G}}_{2\pi}=\OP{id}\).\qedhere

\end{Proof}
\section{Examples}

\begin{Example}[Complex projective spaces]\label{exm:cpn_reduction}
Let\index{symplectic reduction!CPn as a@$\mathbb{CP}^n$ as a}
\(X=\RR^{2n}\) with coordinates \((x_1,y_1,\ldots,x_n,y_n)\)
and symplectic form \(\sum_i dx_i\wedge dy_i\). Consider the
Hamiltonian \(H=\frac{1}{2}\sum_i(x_i^2+y_i^2)\). The
Hamiltonian vector field is \((-y_1,x_1,\ldots,-y_n,x_n)\),
which generates a circle action rotating each \(xy\)-plane at
constant angular speed. The non-empty regular level sets are
the spheres \(M_c\) of radius \(\sqrt{2c}\). If \(c>0\) we get
a nonempty regular level set, whose symplectic quotient is
called {\em complex projective
space}\index{projective space!complex} \(\cp{n-1}\). When
\(c=1\), we call the reduced symplectic form on \(\cp{n-1}\)
the {\em Fubini-Study
form}\index{symplectic form!Fubini-Study}
\(\omega_{FS}\). This is normalised to that the moment image
is a simplex with edges of affine length \(1\), which means
that the complex lines in \(\cp{n-1}\) have symplectic area
\(2\pi\). (Exercise \ref{exr:check_area}: Check directly from
this definition of \(\omega_{FS}\) on \(\cp{1}\) that
\(\int_{\cp{1}}\omega_{FS}=2\pi\).)

\end{Example}
\begin{Remark}
If \(c<0\) then the symplectic quotient is empty. In general,
it is an interesting problem to understand how the topology of
the symplectic quotient varies when the parameter \(c\)
crosses a critical value; this is the subject of {\em
symplectic birational geometry}\index{symplectic birational
geometry} \cite{GuilleminSternbergBirational,TJLiRuan}.

\end{Remark}
\begin{Remark}
We can identify \(\RR^{2n}\) with \(\CC^n\) by introducing
complex coordinates \(z_k=x_k+iy_k\). Then the circle action
generated by \(H\) is precisely the action \(\mathbf{z}\mapsto
e^{it}\mathbf{z}\). The orbits (other than the origin) are
precisely the circles of fixed radius in the complex lines of
\(\CC^n\), so our symplectic quotient coincides with the usual
definition\footnote{If you are unfamiliar with this
description of complex projective space, and with the role it
plays in algebraic geometry, you can read more in Appendix
\ref{ch:complex_projective_spaces}.} of \(\cp{n-1}\) as the
space of complex lines through the origin in \(\CC^n\). The
advantage of defining it as a symplectic quotient is the clean
construction of \(\omega_{FS}\); we did not need to write down
an explicit K\"{a}hler potential.

\end{Remark}
In fact, we also recover the standard torus action on
\(\cp{n-1}\) using Lemma \ref{lma:symmetries_descend}. Consider
the Hamiltonian system \(\bm{G}=(G_1,\ldots,G_n)\colon
\CC^n\to\RR^n\) where \(G_k = \frac{1}{2}|z_k|^2\). These all
commute with the Hamiltonian \(H = \sum_{k=1}^nG_k\) from
Example \ref{exm:cpn_reduction}. Therefore they descend to give
a Hamiltonian system \((\overline{G_1},\ldots,\overline{G_n})
\colon \cp{n-1} \to \RR^n\). The image of this Hamiltonian
system is simply the image of \(\bm{G}|_{H^{-1}(c)}\), which is
the intersection of the hyperplane \(\sum G_k = c\) with the
nonnegative orthant in \(\RR^n\), that is an \((n-1)\)-simplex
(Figure \ref{fig:moment_image_cp2_symp_quotient}).

\begin{figure}[htb]
\begin{center}
\begin{tikzpicture}
\draw[->] (0,0) -- (0,3) node [above] {\(G_3\)};
\draw[->] (0,0) -- (-20:3) node [right] {\(G_1\)};
\draw[dotted] (20:1.2) -- (20:2);
\draw (0,0) -- (20:1.2);
\draw[->] (20:2) -- (20:3) node [right] {\(G_2\)};;
\filldraw[fill=gray,opacity=0.5] (0,2) -- (-20:2) -- (20:2) -- cycle;
\end{tikzpicture}
\end{center}
\caption{The moment image of \(\cp{2}\) in \(\RR^3\) when considered as a symplectic quotient of \(\CC^3\).}
\label{fig:moment_image_cp2_symp_quotient}
\end{figure}

\begin{Example}[Weighted projective spaces]\label{exm:wps}
Let\index{symplectic reduction!weighted projective space as}
\(a_1,\ldots,a_n\) be positive integers and consider the
Hamiltonian \(H\colon\CC^n\to\RR\) given by
\[H(z_1,\ldots,z_n) = \frac{1}{2}\sum_{k=1}^na_k|z_k|^2.\]
This generates the circle action \(z_k\mapsto
e^{ia_kt}z_k\). The symplectic quotient is an
orbifold\index{orbifold} called the weighted projective
space\index{projective space!weighted}
\(\mathbb{P}(a_1,\ldots,a_n)\). If \(a_k\neq 1\) then the
points of the form \((0,\ldots,0,z_k,0,\ldots,0)\) have
nontrivial stabiliser \(\{\mu\in S^1\,:\,\mu^{a_k}=1\}\). If
\(\gcd(a_k,a_{\ell})=1\) then points with \(z_k\) and
\(z_\ell\) both nonzero have trivial stabiliser, so if we
assume that the \(a_k\) are pairwise coprime then the only
singularities are the isolated cyclic quotient singularities
with all but one of the \(z_k\) equal to zero. As in the case
of \(\cp{n-1}\), the functions \(\frac{1}{2}|z_k|^2\) descend
and generate a torus action on \(\mathbb{P}(a_1,\ldots,a_n)\),
whose moment polytope\index{moment polytope!weighted
projective space} is the simplex in \(\RR^n\) with vertices at
\[(c/a_1,0,\ldots,0),\, (0,c/a_2,0,\ldots,0),\, \ldots,\,
(0,\ldots,0,c/a_n).\] In Figure
\ref{fig:moment_image_wps_symp_quotient}(a), you can see the
moment triangle for \(\mathbb{P}(1,2,3)\). If you project it
onto the \(yz\)-plane, you get the triangle shown in Figure
\ref{fig:moment_image_wps_symp_quotient}(b). This has a smooth
point over the Delzant vertex and two singularities modelled
on the \(\frac{1}{2}(1,1)\) and \(\frac{1}{3}(1,2)\)
singularities.

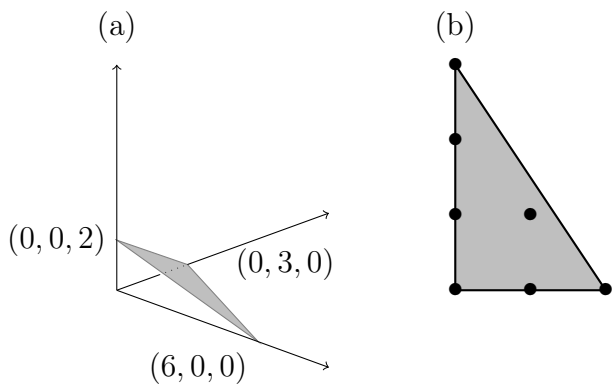
\begin{figure}[htb]
\begin{center}
\begin{tikzpicture}
\node at (0,3.5) {(a)};
\draw[->] (0,0) -- (0,3);
\draw[->] (0,0) -- (-20:3);
\draw[dotted] (20:0.6) -- (20:1);
\draw (0,0) -- (20:0.6);
\draw[->] (20:1) -- (20:3);
\filldraw[fill=gray,opacity=0.5] (0,2/3) -- (-20:2) -- (20:1) -- cycle;
\node at (0,2/3) [left] {\((0,0,2)\)};
\node at (20:1) [right=0.5cm] {\((0,3,0)\)};
\node at (-20:2) [below left] {\((6,0,0)\)};
\begin{scope}[shift={(4.5,0)}]
\node at (0,3.5) {(b)};
\filldraw[fill=lightgray,draw=black,thick] (0,0) -- (2,0) -- (0,3) -- cycle;
\node at (0,0) {\(\bullet\)};
\node at (1,0) {\(\bullet\)};
\node at (2,0) {\(\bullet\)};
\node at (0,1) {\(\bullet\)};
\node at (0,2) {\(\bullet\)};
\node at (0,3) {\(\bullet\)};
\node at (1,1) {\(\bullet\)};
\end{scope}
\end{tikzpicture}
\end{center}
\caption{(a) The moment image of the weighted projective space \(\mathbb{P}(1,2,3)\) in \(\RR^3\) when considered as a symplectic quotient of \(\CC^3\). The coordinates of the vertices of the triangle are given assuming the symplectic reduction is at level \(c=6\). (b) The projection of this triangle to the \(yz\)-plane, with integer lattice points marked.}
\label{fig:moment_image_wps_symp_quotient}
\end{figure}

\end{Example}
\begin{Example}\label{exm:periodic_geodesics}
Consider the round metric on the \(n\)-dimensional unit sphere
\(S^n\); in particular, geodesics are great circles, and if
you traverse a great circle with constant speed 1 then it
returns to its starting point after time \(2\pi\). Let
\(H\colon T^*S^n\to \RR\) be the Hamiltonian
\(\frac{1}{2}|\eta|^2\) generating the cogeodesic
flow\footnote{If you are unfamiliar with cotangent bundles and
cogeodesic flow, you can read more about them in Appendix
\ref{ch:cot_bun}.}. The cogeodesic
flow\index{flow!cogeodesic}\index{cogeodesic flow|see {flow,
cogeodesic}} does not define a circle action on \(T^*S^n\),
but the periodicity of geodesics with speed \(1\) means that
the flow does define a circle action on the level set
\(H^{-1}(1/2)\) (on which geodesics move with speed
\(1\)). This allows us to perform symplectic reduction at that
level. The result is that the space of (oriented) geodesics
parametrised by arc-length on the round \(S^n\) is naturally a
symplectic manifold of dimension \(2n-2\). We can also
identify this manifold of oriented geodesics as a homogeneous
space: the great circles are intersections of
\(S^n\subset\RR^{n+1}\) with oriented 2-planes through the
origin, so the manifold of oriented geodesics coincides with
the Grassmannian\index{Grassmannian} \(\widetilde{Gr}(2,n+1)\)
of oriented 2-planes in \(\RR^{n+1}\).

\end{Example}
\begin{Remark}[Exercise \ref{exr:grassmannian}]\label{rmk:grassmannian}
One can identify the Grassmannian \(\widetilde{Gr}(2,n+1)\)
with the homogeneous space \(O(n+1)/(SO(2)\times O(n-1))\), or
with the quadric\index{quadric hypersurface!and Grassmannian} hypersurface
\(\sum_{i=1}^{n+1}z_i^2=0\) in \(\cp{n}\) with homogeneous
coordinates\footnote{Homogeneous coordinates are introduced in
Appendix \ref{ch:complex_projective_spaces}.}
\([z_1:\cdots:z_{n+1}]\).

\end{Remark}
\begin{Remark}[Exercise \ref{exr:fixing_periods}]\label{rmk:fixing_periods}
The zero-section \(S^n\subset T^*S^n\) consists of fixed
points of the cogeodesic flow, but all the other orbits are
circles. By a similar argument to Theorem
\ref{thm:2dactionangle}, one can modify
\(H=\frac{1}{2}|\eta|^2\) away from the zero-section to get a
Hamiltonian for which all orbits have the same period (i.e.\
giving a circle action away from the zero-section). What
Hamiltonian should you take instead?

\end{Remark}
\section{Symplectic cut}
\label{sct:symp_cut}

The symplectic cut\index{symplectic cut|(} is a particularly
useful case of symplectic reduction, introduced by Lerman
\cite{Lerman}. Here is a simple example to illustrate the
operation.

\begin{Example}
Consider the cylinder \(X = \RR\times S^1\) with symplectic
form \(dp\wedge dq\) (where \(q\) is the angular coordinate on
\(S^1\)). The function \(H = p\) generates the circle action
which rotates the \(S^1\) factor, and the symplectic reduction
at any level \(c\) yields a single point. If we symplectically
reduce at level \(c\) but leave all the other levels alone
then this has the effect of ``pinching'' the cylinder along a
circle. Although the result is singular, it is the union of
two smooth symplectic discs, which we will denote by
\(\overline{X}_{H\leq c}\) and \(\overline{X}_{H\geq
c}\). Each of these contains the point \(H^{-1}(c)/S^1\) at
its centre. We will adopt the convention of calling
\(\overline{X}_{H\geq c}\) {\em the symplectic cut of \(X\) at
level \(c\).}

\end{Example}
\begin{figure}[htb]
\begin{center}
\begin{tikzpicture}
\node (c) at (0,0) {\(\bullet\)};
\node at (c) [below=0.2cm] {\(c\)};
\node at (0.3,1.5) [right] {\(H^{-1}(c)\)};
\draw (-2,0) -- (2,0);
\draw (-2,1) -- (2,1);
\draw (-2,2) -- (2,2);
\draw (-2,1.5) circle [x radius = 0.3, y radius = 0.5];
\draw (2,2) arc [x radius = 0.3, y radius = 0.5, start angle = 90, end angle = -90];
\draw[dotted] (2,1) arc [x radius = 0.3, y radius = 0.5, start angle = 270, end angle = 90];
\draw (0,2) arc [x radius = 0.3, y radius = 0.5, start angle = 90, end angle = -90];
\draw[dotted] (0,1) arc [x radius = 0.3, y radius = 0.5, start angle = 270, end angle = 90];
\draw[->] (0,0.8) -- (0,0.2);
\node at (0,0.5) [right] {\(H\)};
\node at (-1,1.5) {\(X\)};
\begin{scope}[shift={(0,-0.5)}]
\draw (-2,3) to[out=0,in=-90] (0,3.5) to[out=90,in=0] (-2,4);
\draw (2,3) to[out=180,in=-90] (0,3.5) to[out=90,in=180] (2,4);
\draw (-2,3.5) circle [x radius = 0.3, y radius = 0.5];
\draw (2,4) arc [x radius = 0.3, y radius = 0.5, start angle = 90, end angle = -90];
\draw[dotted] (2,3) arc [x radius = 0.3, y radius = 0.5, start angle = 270, end angle = 90];
\node at (-1,3.5) {\(\overline{X}_{H\leq c}\)};
\node at (1,3.5) {\(\overline{X}_{H\geq c}\)};
\node (a) at (0,3.5) {\(\bullet\)};
\node (d) at (0,4.5) {\(H^{-1}(c)/S^1\)};
\draw[->] (0.1,4.3) to[bend left] (0.1,3.6);
\end{scope}
\end{tikzpicture}
\end{center}
\caption{Symplectic cut. It is as if we have cut \(X\) along \(H^{-1}(c)\) and sewn up the wound by collapsing the circle-orbit to a point.}
\label{fig:symp_cut}
\end{figure}
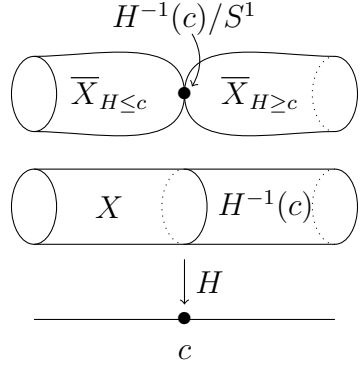

We can perform symplectic cuts in a completely systematic and
general way. Let \((X,\omega)\) be a symplectic manifold and
suppose that \(H\colon X\to\RR\) is a Hamiltonian generating a
circle action \(\phi^H_t\). Pick a regular value \(c\) for which
the circle action on \(H^{-1}(c)\) is free. Consider the product
\(X\times\CC\) with the symplectic form \(\omega + dp\wedge dq\)
(where \(p+iq\in\CC\)) and define the Hamiltonian
\(\tilde{H}_c(x,\xi) = H(x) - c - \frac{1}{2}|\xi|^2\).

\begin{Lemma}[Exercise \ref{exr:circle_action_symp_cut}]\label{lma:circle_action_symp_cut}
The Hamiltonian \(\tilde{H}_c\) generates the circle action
\[\left(x,\xi\right)\mapsto \left(\phi^H_t(x),
e^{-it}\xi\right)\] on \(X\times\CC\).

\end{Lemma}
\begin{Definition}[Symplectic cut at level \(c\)]\label{dfn:symplectic_cut}
The symplectic cut \(\overline{X}_{H\geq c}\) of \(X\) at
level \(c\) is defined to be the symplectic reduction of
\(X\times\RR^2\) with respect to the Hamiltonian
\(\tilde{H}_c\) at level zero. One can also define a
symplectic cut \(\overline{X}_{H\leq c}\) by using the
Hamiltonian \(H-c+\frac{1}{2}|\xi|^2\). We will write
\([x,\xi]\) for the equivalence class of
\((x,\xi)\in\tilde{H}_c^{-1}(0)\) in the symplectic cut.

\end{Definition}
We now try to understand what \(\overline{X}_{H\geq c}\) looks
like.

\begin{Lemma}
Consider the function \(X\times\CC\to\RR\), \((x,\xi)\mapsto
H(x)\). This Poisson-commutes with
\(\tilde{H}_c\) and hence (by Lemma
\ref{lma:symmetries_descend}) descends to a function
\(\overline{H}\colon \overline{X}_{H\geq c}\to \RR\) which
generates a circle action on \(\overline{X}_{H\geq
c}\). Moreover:
\begin{enumerate}
\item [(a)] The image of \(\overline{H}\) is contained in
\([c,\infty)\).
\item [(b)] The preimage \(\overline{H}^{-1}(c)\) is
symplectomorphic to the symplectic reduction
\(H^{-1}(c)/S^1\).
\item [(c)] There is a symplectomorphism \(\Psi\colon
H^{-1}((c,\infty))\to \overline{H}^{-1}((c,\infty))\) which
intertwines the circle actions: \[\Psi(\phi^H_t(x)) =
\phi^{\overline{H}}_t(\Psi(x)).\]

\end{enumerate}
\end{Lemma}
In other words, \(\overline{X}_{H\geq c}\) contains an open set
symplectomorphic to \(H^{-1}((c,\infty))\subset X\); where it
has been cut (along \(H^{-1}(c)\)) the circle-orbits are
collapsed to points, yielding a symplectic submanifold
symplectomorphic to \(H^{-1}(c)/S^2\). Everything below level
\(c\) is thrown away.

\begin{Proof}
The function \((x,\xi)\mapsto H(x)\) generates the Hamiltonian
vector field \((V_H,0)\) on
\(X\times\CC\). Recall that \(\tilde{H}_c = H - c +
\frac{1}{2}|\xi|^2\). Since \(dH(V_H,0)=0\) and
\(d|\xi|^2(V_H,0)=0\) we have \(\{H,\tilde{H}_c\} =
-d\tilde{H}_c(V_H,0)=0\). By Lemma
\ref{lma:symmetries_descend}(a), this function descends to
give a function \(\overline{H}\) on the symplectic
quotient. Since \((V_H,0)\) generates a circle action on
\(X\times\CC\), Lemma \ref{lma:symmetries_descend}(c) implies
that \(\overline{H}\) generates a circle action on
\(\overline{X}_{H\geq c}\).\index{Hamiltonian!circle action|)}

Property (a): If \((x,\xi)\in \tilde{H}_c^{-1}(0)\) then
\(H(x) - c + \frac{1}{2}|\xi|^2=0\), so \[H(x) = c +
\frac{1}{2}|\xi|^2\geq c.\] Therefore
\(\overline{H}([x,\xi])=H(x)\geq c\).

Property (b): The map \(H^{-1}(c)/S^1\to
\overline{H}^{-1}(c)\subset \overline{X}_{H\geq c}\) defined
by \([x]\mapsto [x,0]\) is the required symplectomorphism.

Property (c): The required symplectomorphism is \[\Psi(x) =
\left[x,\sqrt{2(H(x)-c)}\right].\qedhere\]

\end{Proof}
\section{Further examples}

Most of our examples will be based on the following special
case.

\begin{Example}\label{exm:toric_one_cut}
Suppose that \(X\) is a toric manifold\index{toric manifold}
and that \(\mu\colon X\to\RR^n\) is the moment map with moment
image \(\Delta = \mu(X)\). Pick a \(\ZZ\)-linear map
\(s\colon\RR^n\to\RR\) and take \(H=s\circ\mu\). If we
symplectically cut at level \(c\) then the result is still
toric by Lemma \ref{lma:symmetries_descend}, and the moment
image is given by \[\{x\in \Delta\, :\, s(x) \geq c\}.\] In
other words, this is the truncation\index{moment
polytope!truncation}\index{truncation|see {moment polytope,
truncation}} of the moment image by the half-space \(\{s\geq
c\}\).

\end{Example}
\begin{Remark}[Exercise \ref{exr:why_z_linear}]\label{rmk:why_z_linear}
We need \(s\) to be \(\ZZ\)-linear (or at least
\(\QQ\)-linear) for this construction to work. Can you explain
why?

\end{Remark}
\begin{Theorem}[Exercise \ref{exr:toric_from_cut}]\label{thm:toric_from_cut}
Any convex rational polytope\index{polytope!convex rational}
\(\Delta\) occurs as the moment image of a toric Hamiltonian
system (on a possibly singular
space).\index{toric variety!from symplectic cut}

\end{Theorem}
\begin{Example}\label{exm:symp_cut_blow_up}
Take \(X=\CC^n\) and
\(H(\bm{z})=\frac{1}{2}\sum_{k=1}^n|z_k|^2\). We saw in
Example \ref{exm:cpn_reduction} that the symplectic {\em
reduction} at level \(c > 0\) is \(\cp{n-1}\). If instead we
take the symplectic {\em cut} at this level, we get a toric
manifold\index{toric manifold!CPn@$\mathbb{CP}^n$} whose
moment image is \[\{\bm{b}\in\RR_{\geq 0}^n\, :\, \sum b_k
\geq c\}.\] (See Figure \ref{fig:symp_cut_blow_up} for the
\(n=2\) case). This contains the complement of the symplectic
ball \(B_c := \{\bm{z}\in\CC^n\,:\, |z|^2\leq 2c\}\) and it
also contains a copy of \(\cp{n-1}\) living over the newly cut
facet \(\sum b_k = c\). This is known as the {\em symplectic
blow-up}\index{blow-up!symplectic|(} of \(X\) in the ball
\(B_c\). More generally, you can perform this operation on any
Delzant vertex of a moment
polytope\index{moment polytope!truncation} (see Example
\ref{exm:2_ball_3_ball} below). The symplectic area of the
exceptional curve is \(2\pi c\).

\begin{figure}
\begin{center}
\begin{tikzpicture}
\filldraw[fill=lightgray,opacity=0.5,draw=none] (0,3) -- (0,0) -- (3,0) -- (3,3) -- cycle;
\draw[black,thick] (0,3) -- (0,0) -- (3,0);
\draw[dotted,thick] (0,2) -- (2,0);
\draw[snake it,->] (3.5,1.5) -- (4.5,1.5) node [above,midway] {cut};
\begin{scope}[shift={(5,0)}]
\filldraw[fill=lightgray,opacity=0.5,draw=none] (0,3) -- (0,2) -- (2,0) -- (3,0) -- (3,3) -- cycle;
\draw[black,thick] (0,3) -- (0,2) -- (2,0) -- (3,0);
\end{scope}
\end{tikzpicture}
\end{center}
\caption{The moment polygon for symplectic blow-up of \(\CC^2\) at the origin.}
\label{fig:symp_cut_blow_up}
\end{figure}
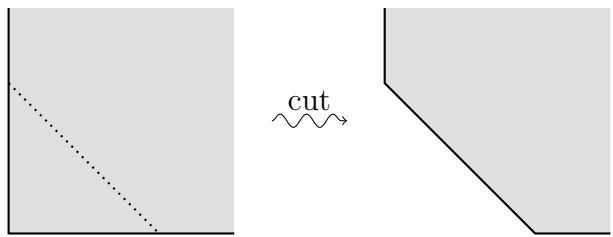

\end{Example}
\begin{Remark}
Compare Figure \ref{fig:symp_cut_blow_up} with Figure
\ref{fig:taut_bun}. This shows that symplectically blowing up
a ball is essentially the same as the complex blow-up of a
smooth point, though of course they happen in different
categories, so one must be careful when relating them in
practice. A paper which very carefully explains (and uses) the
relationship is {\cite[Section 2]{McDuffPolterovich}}.

\end{Remark}
\begin{Example}
More generally, if \(Y\subset X\) is a symplectic submanifold
of real codimension \(2k\) then one can perform this
symplectic blow-up\index{blow-up!along symplectic submanifold}
along \(Y\): each fibre of the normal bundle to \(Y\) is
replaced by its blow-up in a ball centred at the origin. The
result is a symplectic manifold which contains:
\begin{itemize}
\item the complement of a neighbourhood of \(Y\)
\item a copy of the projectivisation of the normal bundle of \(Y\)
in \(X\) (considered as a complex vector bundle).
\end{itemize}
For more details, see Usher's MathOverflow answer
\cite{UsherMO}.

\end{Example}
\begin{Example}[Exercise \ref{exr:2_ball_3_ball}]\label{exm:2_ball_3_ball}
There is a blow-up of \(\cp{1}\times\cp{1}\) in two disjoint
symplectic balls which is symplectomorphic to a blow-up of
\(\cp{2}\) in three disjoint balls.

\end{Example}
\begin{Example}[Exercise \ref{exr:dp_as_slice}]\label{exm:dp_as_slice}
Show that the common blow-up from Example
\ref{exm:2_ball_3_ball} arises as a symplectic reduction of
\(\cp{1}\times\cp{1}\times\cp{1}\).\index{blow-up!symplectic|)}

\end{Example}
\begin{Example}
As in Example \ref{exm:symp_cut_blow_up}, we take \(X=\CC^n\)
and \(H(\bm{z})=\frac{1}{2}\sum_{k=1}^n |z_k|^2\), but this
time we use the symplectic cut \(\overline{X}_{H \leq
c}\). The result is a symplectic toric
manifold\index{toric manifold!CPn@$\mathbb{CP}^n$} whose
moment image\index{moment polytope!CPn@$\mathbb{CP}^n$} is the
\(n\)-simplex \(\{\bm{b}\in\RR^n\,:\, b_k\geq 0,\sum b_k\leq
c\}\). When \(c=1\), we saw this moment polytope arise in
Example \ref{exm:cpn} as the moment image of the standard
torus action on \(\cp{n}\), so \(\overline{X}_{H\leq 1}\) is
symplectomorphic to \(\cp{n}\) with the Fubini-Study
form. More generally, if \(c>0\), the only difference is that
the symplectic form is rescaled by a factor of \(c\).

\end{Example}
The previous example illustrates how the standard operation of
{\em projective
compactification}\index{projective compactification} in
algebraic geometry can be understood using symplectic cut. We
can do something similar starting with Example
\ref{exm:periodic_geodesics}.

\begin{Example}\label{exm:periodic_geodesics_ctd}
Consider the Hamiltonian \(H\colon T^*S^n\to\RR\) from Example
\ref{exm:periodic_geodesics}. The symplectic cut
\(\overline{(T^*S^n)}_{H\leq 1/2}\) is a compact symplectic
manifold containing an open subset symplectomorphic to the
open unit cotangent bundle\index{cotangent bundle} of \(S^n\)
and a ``compactifying divisor''\index{divisor!compactifying}
\(H^{-1}(1/2)/S^1\) symplectomorphic to the
Grassmannian\index{Grassmannian} \(\widetilde{Gr}(2,n+1)\). In
fact, \(\overline{(T^*S^n)}_{H\leq 1/2}\) is symplectomorphic
to a quadric\index{quadric hypersurface!as symplectic cut} hypersurface in
\(\cp{n+1}\) and the compactifying divisor is a hyperplane
section. For an in-depth discussion of this and similar
examples, see \cite{AudinSkeletons}.

\end{Example}
\section{Resolution of singularities}

We have seen that symplectic cut can be used as a symplectic
Ersatz for blowing-up\index{blow-up!resolving singularities|(}
in algebraic geometry. Just as blowing up allows us to resolve
singularities\index{singularity!resolution|(} in algebraic
geometry, symplectic cut can be used to remove singularities in
symplectic geometry.

\begin{Example}\label{exm:minimal_resolution}
Consider\index{singularity!resolution!minimal|(} the cyclic
quotient singularity\index{singularity!cyclic quotient}
\(\frac{1}{a_0}(1,a_1)\) with moment polygon \(\pi(a_0,a_1)\)
from Example \ref{exm:cycquot}. We can symplectically cut
using a horizontal slice at a level above the singularity to
obtain a polygon \(\tilde{\pi}(a_0,a_1)\) which now has two
vertices: a Delzant corner and a
corner\index{polytope!non-Delzant} modelled on
\(\pi(a_1,a_2)\) where \(0 < a_2\leq a_1\) satisfies \(a_2 =
-a_0\mod a_1\), i.e.\ \(a_2 = y_1a_1 - a_0\) for some
\(y_1\). Namely, the matrix \(M_1 = \lmatrix 0 & -1 \\ 1 &
y_1\rmatrix \) satisfies \((-1,0)M_1 = (0,1)\),
\((a_0,a_1)M_1=(a_1,a_2)\). See Figure
\ref{fig:minimal_resolution}.

If \(a_2=a_1\) then this second corner is Delzant, and the
symplectic cut is smooth. Otherwise, we can iterate this
procedure by slicing the polygon \(\pi(a_1,a_2)\) horizontally
(i.e.\ slicing \(\tilde{\pi}(a_0,a_1)\) parallel to
\((1,0)M_1^{-1}\)). We get a decreasing sequence of positive
integers \(a_0,a_1,a_2,\ldots\) and a sequence of positive
integers \(y_1,y_2,\ldots\) with \(a_{k+1} = y_k
a_k-a_{k-1}\). At some point we necessarily find
\(a_{m+1}=a_{m+2}\) because the sequence cannot continue
decreasing forever. The result of all these cuts is a Delzant
polygon. See Figures \ref{fig:minimal_resolution} and
\ref{fig:min_res_3_2}.

Unpacking the recursion formula for \(a_{k+1}\), we see that
\[a_0/a_1 = y_1 - a_2/a_1 = y_1 - \frac{1}{a_1/a_2} = y_1 -
\frac{1}{y_2 - a_3/a_2} = \cdots = y_1 -
\frac{1}{y_2-\frac{1}{\ddots - \frac{1}{y_{m}}}}.\] Thus the
numbers \(y_i\) are the coefficients in the continued
fraction\index{continued
fraction!Hirzebruch-Jung}\index{continued fraction!minimal
resolution|(} expansion\footnote{This is the {\em
Hirzebruch-Jung} convention for taking continued fractions
with minus signs.} of \(a_0/a_1\).

This process has introduced \(m\) compact edges into our
polygon. There are symplectic spheres\index{symplectic chain
of spheres|(} living over these in the iterated symplectic
cut. These spheres have
self-intersection\index{self-intersection!of curve} numbers
\(-y_1,-y_2,\ldots,-y_{m}\). To see this, it suffices to
check it for \(y_1\) because of the iterative nature of the
process. There are two cases to consider:
\begin{itemize}
\item If our process terminated after a single cut then
\(a_0/a_1=y_1\) is an integer, so \(a_1=1\) and
\(a_0=y_1\). In this case, the compact edge of
\(\tilde{\pi}(a_0,a_1)\) has rays emanating in the \((0,1)\)
and \((a_0,a_1)=(y_1,1)\) directions from its endpoints. By
Lemma \ref{lma:selfint}, the corresponding sphere has
self-intersection \(-y_1\).
\item If our process takes two or more steps then the leftmost
compact edge has outgoing rays pointing in the \((0,1)\) and
\((1,0)M_1^{-1} = (y_1,1)\) directions. Again, by Lemma
\ref{lma:selfint}, the corresponding sphere has
self-intersection \(-y_1\).

\end{itemize}
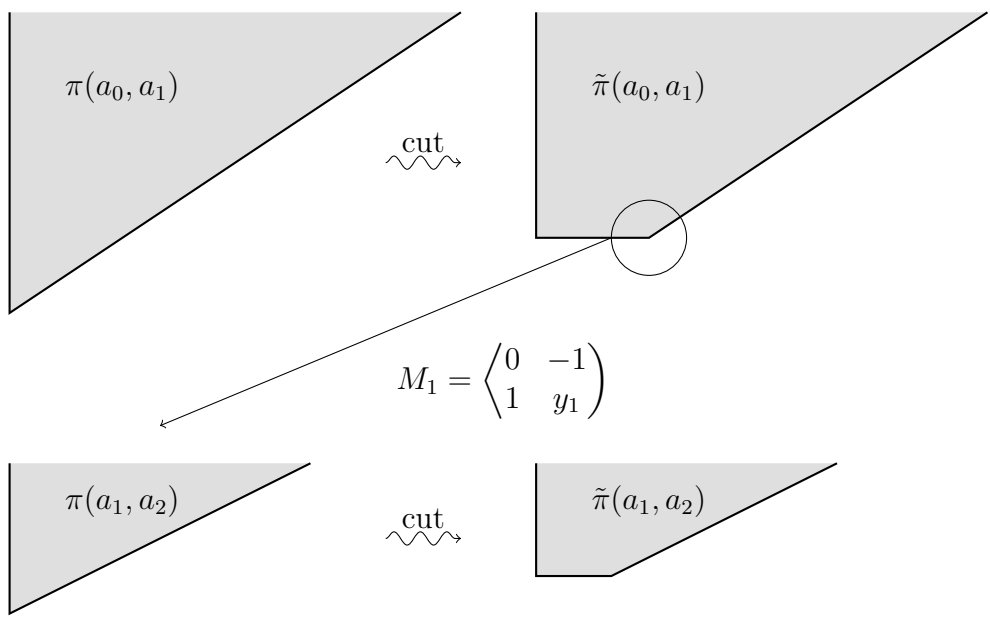
\begin{figure}[htb]
\begin{center}
\begin{tikzpicture}
\filldraw[fill=lightgray,opacity=0.5,draw=none] (0,4) -- (0,0) -- (6,4) -- cycle;
\draw[black,thick] (0,4) -- (0,0) -- (6,4);
\node at (1.5,3) {\(\pi(a_0,a_1)\)};
\draw[snake it,->] (5,2) -- (6,2) node [midway,above] {cut};
\begin{scope}[shift={(7,0)}]
\filldraw[fill=lightgray,opacity=0.5,draw=none] (0,4) -- (0,1) -- (1.5,1) -- (6,4) -- cycle;
\draw[black,thick] (0,4) -- (0,1) -- (1.5,1) -- (6,4);
\node at (1.5,3) {\(\tilde{\pi}(a_0,a_1)\)};
\draw (1.5,1) circle [radius=0.5cm];
\draw[->] (1,1) -- (-5,-1.5) node [midway,below right] {\(M_1=\lmatrix 0 & -1 \\ 1 & y_1\rmatrix \)};
\end{scope}
\begin{scope}[shift={(0,-4)}]
\filldraw[fill=lightgray,opacity=0.5,draw=none] (0,2) -- (0,0) -- (4,2) -- cycle;
\draw[black,thick] (0,2) -- (0,0) -- (4,2);
\node at (1.5,1.5) {\(\pi(a_1,a_2)\)};
\draw[snake it,->] (5,1) -- (6,1) node [midway,above] {cut};
\end{scope}
\begin{scope}[shift={(7,-4)}]
\filldraw[fill=lightgray,opacity=0.5,draw=none] (0,2) -- (0,0.5) -- (1,0.5) -- (4,2) -- cycle;
\draw[black,thick] (0,2) -- (0,0.5) -- (1,0.5) -- (4,2);
\node at (1.5,1.5) {\(\tilde{\pi}(a_1,a_2)\)};
\end{scope}
\end{tikzpicture}
\end{center}
\caption{The minimal resolution of cyclic quotient surface singularities via symplectic cuts. In this specific example we have used \(a_0=3\), \(a_1=2\), \(a_2=1\), \(y_1=y_2=2\) (so the process terminates after two cuts: the result is shown in Figure \ref{fig:min_res_3_2}).}
\label{fig:minimal_resolution}
\end{figure}

\begin{figure}[htb]
\begin{center}
\begin{tikzpicture}
\filldraw[fill=lightgray,opacity=0.5,draw=none] (0,4) -- (0,1) -- (1,1) -- (3,2) -- (6,4) -- cycle;
\draw[black,thick] (0,4) -- (0,1) node {\(\bullet\)} -- (1,1) node {\(\bullet\)} -- (3,2) node {\(\bullet\)} -- (6,4);
\node at (0.5,1) [below] {\((1,0)\)};
\node at (2,1.5) [below right] {\((2,1)\)};
\node at (4.5,3) [below right] {\((3,2)\)};
\end{tikzpicture}
\end{center}
\caption{The combined result of performing both cuts in Figure \ref{fig:minimal_resolution}; vertices are marked with dots for emphasis. In this example (the \(\frac{1}{3}(1,2)\) singularity) the vectors indicate the directions of the edges and using Lemma \ref{lma:selfint} you can check that we have introduced two \(-2\)-spheres in the minimal resolution.}
\label{fig:min_res_3_2}
\end{figure}
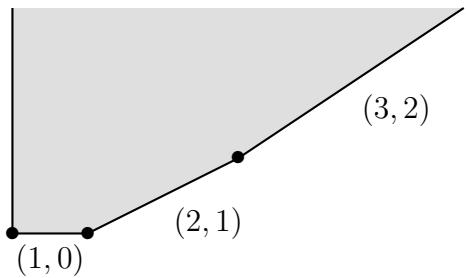

\end{Example}
This process allows us to replace a \(\frac{1}{a_0}(1,a_1)\)
singularity by a chain of symplectic
spheres\index{symplectic chain of spheres|)} with
self-intersections determined by the continued fraction
expansion of \(a_0/a_1\). This is precisely what happens when we
perform the {\em minimal resolution} of this singularity in
complex geometry, so we often refer to this sequence of
symplectic cuts as taking the minimal resolution. Figures
\ref{fig:minimal_resolution} and \ref{fig:min_res_3_2}
illustrate this process in the case \(a_0=3\), \(a_1=2\). Since
\(\frac{3}{2} = 2-\frac{1}{2}\), the minimal resolution replaces
the \(\frac{1}{3}(1,2)\) singularity by two
\(-2\)-spheres\index{continued fraction!minimal
resolution|)}.

\begin{Remark}
Just as symplectic blow-up involves a choice of parameter (the
area of the exceptional sphere), forming the minimal
resolution by sympplectic cuts involves choices of
edge-lengths for the symplectic cuts. This subtlety is absent
from the minimal resolution in complex geometry, but appears
when you try to equip the minimal resolution with an ample
line bundle or K\"{a}hler
form.\index{singularity!resolution!minimal|)}

\end{Remark}
In higher dimensions, resolution of singularities is more
complicated. We consider just one 6-dimensional example to give
a flavour of what can happen.

\begin{Example}[The conifold]\label{exm:odp_polytope}
Consider the affine variety
\(C:=\{z_1z_4=z_2z_3\}\subset\CC^4\). This has a singularity
at the origin which goes by many names ({\em \(A_1\)
singularity}, {\em node}, {\em ordinary double point}, {\em
conifold})\index{singularity!conifold|(}\index{singularity!A1@$A_1$|see
{singularity, conifold}}\index{singularity!ordinary double
point|see {singularity, conifold}}\index{singularity!nodal|see
{singularity, conifold}}. It admits a Hamiltonian
\(T^3\)-action: \[\left(e^{it_1}z_1, e^{it_2}z_2, e^{it_3}z_3,
e^{i(t_2+t_3-t_1)}z_4\right).\] The moment map for this action
is \[\mu(z_1,\ldots,z_4) = \left(\frac{1}{2}(|z_1|^2-|z_4|^2),
\frac{1}{2}(|z_2|^2+|z_4|^2),
\frac{1}{2}(|z_3|^2+|z_4|^2)\right).\] If we write
\(H_1,H_2,H_3\) for these three functions then we have
\[H_2\geq 0,\quad H_3\geq 0,\quad H_1+H_2\geq 0,\quad
H_1+H_3\geq 0.\] These four inequalities cut out a polyhedral
cone \(\Delta\) in \(\RR^3\), spanned by the four rays
\[(r,0,0),\quad (0,r,0),\quad (0,0,r),\quad (-r,r,r),\quad r >
0.\] (See Figure \ref{fig:odp_polytope}). These rays are the
moment images of the curves \[(z,0,0,0),\quad (0,z,0,0),\quad
(0,0,z,0),\quad (0,0,0,z)\] all of which are contained in
\(C\), so \(\mu(C)\) contains all four of these rays and hence
their convex hull, which is the whole of \(\Delta\). Thus
\(\mu(C)=\Delta\)\index{moment polytope!conifold}. Note that
\(\Delta\) is not Delzant\index{polytope!non-Delzant} at the
origin (corresponding to the nodal singularity of \(C\)).

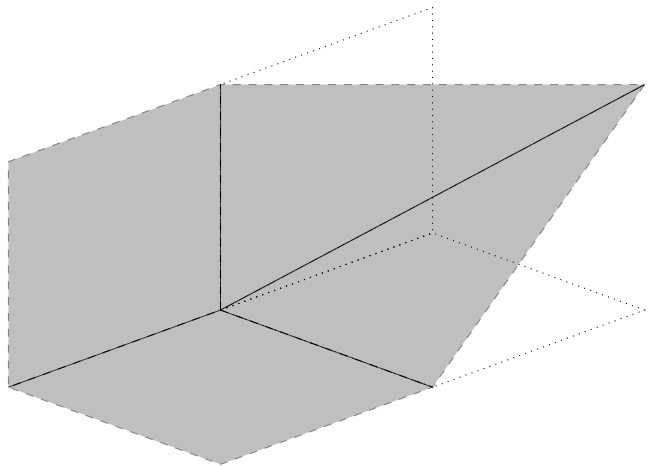
\begin{figure}
\begin{center}
\begin{tikzpicture}
\draw (0,0) -- (0,3);
\draw (0,0) -- (200:3);
\draw (0,0) -- (-20:3);
\draw (0,0) -- (5.63815572471545,3);
\draw[dotted] (0,0) -- (20:3);
\filldraw[dashed,fill=gray,opacity=0.5] (0,0) -- (0,3) -- (5.63815572471545,3) -- cycle;
\filldraw[dashed,fill=gray,opacity=0.5] (0,0) -- (-20:3) -- (5.63815572471545,3) -- cycle;
\filldraw[dashed,fill=gray,opacity=0.5] (0,0) -- (0,3) -- (-2.81907786235773,1.97393957002299) -- (200:3) -- cycle;
\filldraw[dashed,fill=gray,opacity=0.5] (0,0) -- (-20:3) -- (0,-2.05212085995401) -- (200:3) -- cycle;
\draw[dotted] (0,0) -- (20:3) -- (5.63815572471545,0) -- (-20:3) -- cycle;
\draw[dotted] (0,0) -- (20:3) -- (2.81907786235773,4.02606042997701) -- (0,3) -- cycle;
\end{tikzpicture}
\end{center}
\caption{The moment polytope for the nodal 3-fold in Example \ref{exm:odp_polytope} is an infinite cone extending this figure.}
\label{fig:odp_polytope}
\end{figure}

One can find some resolutions of this nodal 3-fold by taking
symplectic cuts.

\end{Example}
\begin{Example}[The fully-resolved conifold]\label{exm:symp_cut_blow_up_node}
Take the symplectic cut of the conifold \(C\) with respect to
the Hamiltonian \(\frac{1}{2}\sum_{k=1}^4|z_k|^2\) at level
\(c=1\). In terms of our functions \(H_1,H_2,H_3\) this is
just \(H_1+H_2+H_3\), so we obtain the polytope\index{moment
polytope!fully-resolved conifold} shown in Figure
\ref{fig:odp_polytope_resolution_1}. We have introduced a new
quadrilateral facet \(H_1+H_2+H_3=1\). This quadrilateral has
vertices at \((1,0,0)\), \((0,1,0)\), \((0,0,1)\) and
\((-1,1,1)\). In fact, this is isomorphic to a square under
the integral affine transformation \((x,y,z)\mapsto (y,z)\),
so the preimage of our new facet is the toric manifold
associated to a square. In Example
\ref{exm:products_of_spheres}, we saw that this is
symplectomorphic to \(\cp{1}\times \cp{1}\). We have replaced
the nodal point of \(C\) with an ``exceptional
divisor''\index{divisor!exceptional}
\(\cp{1}\times\cp{1}\). Note that this is precisely what we
would get if we treated \(z_1,\ldots,z_4\) as homogeneous
coordinates: \(z_1z_4=z_2z_3\) is a smooth projective quadric
surface, and hence biholomorphic to \(\cp{1}\times\cp{1}\)
(see Example \ref{exm:segre}). In the language of algebraic
geometry, we have performed a blow-up of \(\CC^4\) at the
origin and taken the proper transform of \(C\).

\end{Example}
\begin{figure}
\begin{center}
\begin{tikzpicture}
\draw[dotted] (0,0) -- (0,3);
\draw[dotted] (0,0) -- (200:3);
\draw[dotted] (0,0) -- (-20:3);
\draw (5.63815572471545/3,1) -- (5.63815572471545,3);
\draw[dotted] (0,0) -- (20:3);
\filldraw[dashed,fill=gray,opacity=0.5] (0,1) -- (0,3) -- (5.63815572471545,3) -- (5.63815572471545/3,1) -- cycle;
\filldraw[dashed,fill=gray,opacity=0.5] (-20:1) -- (-20:3) -- (5.63815572471545,3) -- (5.63815572471545/3,1) -- cycle;
\filldraw[dashed,fill=gray,opacity=0.5] (0,1) -- (0,3) -- (-2.81907786235773,1.97393957002299) -- (200:3) -- (200:1) -- cycle;
\filldraw[dashed,fill=gray,opacity=0.5] (-20:1) -- (-20:3) -- (0,-2.05212085995401) -- (200:3) -- (200:1) -- cycle;
\filldraw[fill=gray,opacity=0.5] (200:1) -- (-20:1) -- (5.63815572471545/3,1) -- (0,1) -- cycle;
\draw (200:1) -- (-20:1) -- (5.63815572471545/3,1) -- (0,1) -- cycle;
\draw (200:1) -- (200:3);
\draw (-20:1) -- (-20:3);
\draw (0,1) -- (0,3);
\draw[dotted] (0,0) -- (20:3) -- (5.63815572471545,0) -- (-20:3) -- cycle;
\draw[dotted] (0,0) -- (20:3) -- (2.81907786235773,4.02606042997701) -- (0,3) -- cycle;
\end{tikzpicture}
\end{center}
\caption{The moment polytope of the fully-resolved conifold in Example \ref{exm:symp_cut_blow_up_node}. This is obtained by intersecting the polytope from Figure \ref{fig:odp_polytope} with the half-space \(x+y+z\geq 1\).}
\label{fig:odp_polytope_resolution_1}
\end{figure}
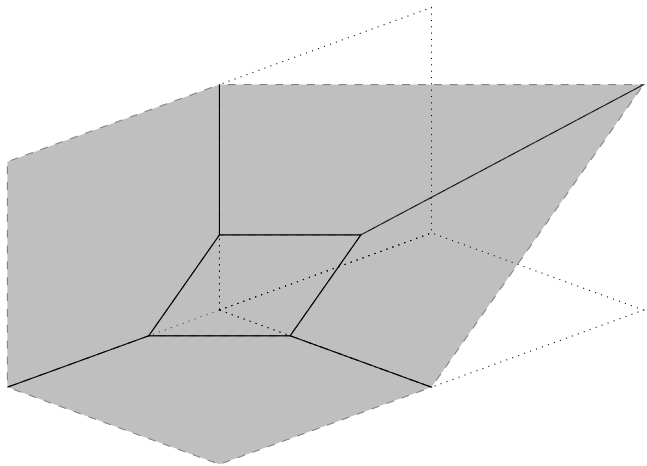

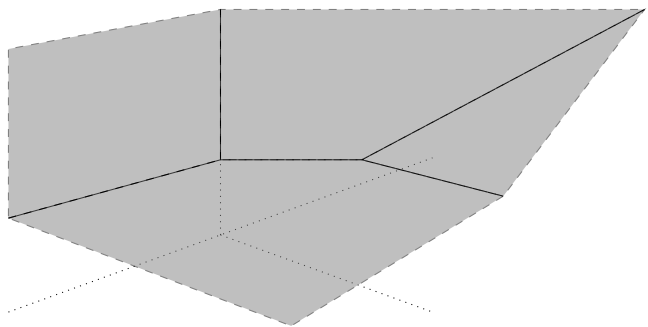
\begin{figure}
\begin{center}
\begin{tikzpicture}
\coordinate (a) at (-2.81907786235773,-0.0260604299770062+0.25);
\coordinate (b) at (0,3);
\coordinate (c) at (0,1);
\coordinate (d) at (5.63815572471545/3,1);
\coordinate (e) at (5.63815572471545,3);
\coordinate (f) at (3.7587704831464,0.315959713348663+0.2);
\coordinate (g) at (-2.81907786235773,1.97393957002299+0.5);
\coordinate (h) at ({-2.81907786235773+3.7587704831464},{-0.0260604299770062+0.315959713348663-1.5});
\draw[dotted] (200:3) -- (20:3);
\draw[dotted] (0,0) -- (c.center);
\draw[dotted] (0,0) -- (-20:3);
\draw (a.center) -- (c.center) -- (b.center);
\draw (f.center) -- (d.center) -- (e.center);
\draw (c.center) -- (d.center);
\filldraw[dashed,fill=gray,opacity=0.5] (a.center) -- (c.center) -- (b.center) -- (g.center) -- cycle;
\filldraw[dashed,fill=gray,opacity=0.5] (d.center) -- (e.center) -- (f.center) -- cycle;
\filldraw[dashed,fill=gray,opacity=0.5] (b.center) -- (c.center) -- (d.center) -- (e.center) -- cycle;
\filldraw[dashed,fill=gray,opacity=0.5] (a.center) -- (c.center) -- (d.center) -- (f.center) -- (h.center) -- cycle;
\end{tikzpicture}
\end{center}
\caption{The moment polytope of the small-resolved conifold in Example \ref{exm:odp_polytope}. This is obtained by intersecting the polytope in Figure \ref{fig:odp_polytope} with the half-space \(z\geq \epsilon\).}
\label{fig:odp_polytope_resolution_2}
\end{figure}

\begin{Example}[Small-resolved conifold]\label{exm:small_resolution_symp_cut}
Take\index{singularity!resolution!small|(} a linear map \(s\)
which vanishes along one of the facets of \(\Delta\), for
example \(s(x,y,z)=z\). Take the symplectic cut
\(\overline{C}_{s\circ \mu\geq \epsilon}\) of the conifold for
some small \(\epsilon>0\). If \(\Delta\) were Delzant, this
truncation would have no effect on the combinatorics of the
moment polytope\index{moment polytope!small-resolution} and
would just change the lengths of some edges. Because
\(\Delta\) is not Delzant, the result is to introduce a new
edge (see Figure \ref{fig:odp_polytope_resolution_2}). This is
an example of a {\em small resolution}: the symplectic cut is
a smooth manifold, but instead of replacing the singularity
with a divisor, we have replaced it with a curve. In
algebro-geometric language, we have blown-up the {\em Weil
divisor}\index{divisor!Weil} \(\{z_3=z_4=0\}\subset
C\). Roughly speaking, a Weil divisor is a complex codimension
1 submanifold; if a Weil divisor can be defined by the
vanishing of a single polynomial then it is called a {\em
Cartier divisor}\index{divisor!Cartier}. Blowing up along a
Cartier divisor has no effect, but we have blown up along a
non-Cartier divisor (requiring both \(z_3\) and \(z_4\) to
vanish).\index{singularity!resolution|)}\index{singularity!resolution!small|)}\index{blow-up!resolving
singularities|)}\index{singularity!conifold|)}

\end{Example}
\section{Solutions to inline exercises}

\begin{Exercise}[Lemma \ref{lma:easy_properties}]\label{exr:easy_properties}
Suppose \(H\colon X\to\RR\) generates a circle
action.
\begin{enumerate}
\item [(a)] The critical points of \(H\) are precisely the fixed
points of the circle action.
\item [(b)] The level sets \(M_c\) are preserved by the circle action.
\item [(c)] If \(x\) is a regular point then \(T_xM_c=\ker(dH)\).
\item [(d)] If \(v\in T_xM_c\) satisfies \(\omega(v,w)=0\) for all
\(w\in T_xM_c\) then \(v\in\OP{span}(V_H)\).
\end{enumerate}
\end{Exercise}
\begin{Solution}
(a) The vector field \(V_H\)
satisfies \(\iota_{V_H}\omega=-dH\). Since \(\omega\) is
nondegenerate, we see that \(V_H=0\) if and only if \(dH=0\),
so the zeros of \(V_H\) coincide with the critical points of
\(H\). A point is fixed under the circle action if and only if
the vector field vanishes there.

(b) This is immediate from Lemma \ref{lma:conserved}.

(c) Suppose that \(\gamma\) is a path in
\(M_c=H^{-1}(c)\). Then the directional derivative of \(H\)
along \(\dot{\gamma}\) vanishes because \(H\circ\gamma=c\), so
\(dH(\dot{\gamma})=0\). Because \(T_xM_c\) consists of tangent
vectors to paths through \(x\) in \(M_c\), we see that
\(T_xM_c\subset\ker(dH)\). Since \(x\) is
regular, both \(T_xM_c\) and \(\ker(dH)\)
have codimension 1, so the containment
\(T_xM_c\subset\ker(dH)\) implies they coincide.

(d) The symplectic orthogonal complement \((T_xM_c)^{\omega}\)
is 1-dimensional and contains the span of \(V_H\), therefore
it equals the span of \(V_H\). This means that if
\(\omega(v,w)=0\) for all \(w\in T_xM_c\) then
\(v\in\OP{span}(V_H)\).\qedhere

\end{Solution}
\begin{Exercise}[Lemma \ref{lma:closed_submersion}]\label{exr:closed_submersion}
If \(p\colon M\to Q\) is a submersion and \(\eta\) is a \(k\)-form
on \(Q\) such that \(p^*\eta=0\) then \(\eta=0\).
\end{Exercise}
\begin{Solution}
If \(\eta\neq 0\) then there exist vectors \(\xi_1,\ldots,\xi_k\in
TQ\) such that \(\eta(\xi_1,\ldots,\xi_k)\neq 0\). Since \(p\) is a
submersion, \(\xi_i=p_*v_i\) for some vectors \(v_1,\ldots,v_k\in
TM\). Therefore \(0\neq
\eta(\xi_1,\ldots,\xi_k)=(p^*\eta)(v_1,\ldots,v_k)=0\), which is a
contradiction.\qedhere

\end{Solution}
\begin{Exercise}[Example \ref{exm:cpn_reduction}]\label{exr:check_area}
Check directly from the definition of \(\omega_{FS}\) on
\(\cp{1}\) that \(\int_{\cp{1}}\omega_{FS}=2\pi\).
\end{Exercise}
\begin{Proof}
Recall that \(\omega_{FS}\) comes from symplectically cutting
\(\CC^2\) at radius \(r = \sqrt{2}\). Consider the affine
coordinate patch \(\{[x+iy:1]\,:\,x+iy\in\CC\}\) on
\(\cp{1}\). This covers all but a point of \(\cp{1}\), and we
can pick a section of the symplectic quotient \(p\colon
S^3(r)\to\cp{1}\) living over this patch, for example
\([x+iy:1]\mapsto (\frac{r(x+iy)}{1+x^2+y^2},
\frac{r}{1+x^2+y^2})\in S^3(r)\subset\CC^2\). The integral of
\(\omega_{FS}\) over this coordinate patch is given by the
integral of \(\omega_{\CC^2}\) over the section (since
\(p^*\omega_{FS} = i^*\omega_{\CC^2}\) where \(i\) is the
inclusion \(S^3(r)\to\CC^2\)). The form \(\omega_{\CC^2}\) is
the sum \(pr_1^*\omega_{\CC}+pr_2^*\omega_{\CC}\) where
\(pr_1\) and \(pr_2\) are the projections to the first and
second factors. The projection of our section to the second
factor is 1-dimensional, \(pr_2^*\omega_{\CC}\) integrates
trivially over the section. The image of the section under
\(pr_1\) is the disc of radius \(r\), so the integral of
\(pr_1^*\omega_{\CC}\) over the section is \(\pi r^2\). Since
\(r=\sqrt{2}\), this gives area \(2\pi\), as required.\qedhere

\end{Proof}
\begin{Exercise}[Remark \ref{rmk:grassmannian}]\label{exr:grassmannian}
One can identify the Grassmannian\index{Grassmannian}
\(\widetilde{Gr}(2,n+1)\) with the homogeneous space
\(O(n+1)/(SO(2)\times O(n-1))\), or with the quadric
hypersurface \(\sum_{i=1}^{n+1}z_i^2=0\) in \(\cp{n}\) with
homogeneous coordinates \([z_1:\cdots:z_{n+1}]\).
\end{Exercise}
\begin{Solution}
There is a transitive action of \(O(n+1)\) on oriented
\(2\)-planes in \(\RR^{n+1}\). By the orbit-stabiliser
theorem\footnote{i.e.\ the theorem which identifies the
\(G\)-orbit of \(x\) with \(G/\OP{Stab}(x)\), in whatever
category you are working, e.g.\ differentiable manifolds and
smooth actions.}, this means that
\(\widetilde{Gr}(2,n+1)=O(n+1)/\OP{Stab}(\RR^2)\), where
\(\OP{Stab}(\RR^2)\) is the subgroup of \(O(n+1)\) stabilising
the standard oriented 2-plane
\(\{(x_1,x_2,0,\ldots,0)\,:\,x_1,x_2\in\RR\}\). This
stabiliser consists of block-matrices \(\begin{pmatrix} A & 0
\\ 0 & B\end{pmatrix}\) with \(A\in SO(2)\) and \(B\in
O(n-1)\).

To identify this homogeneous space with the quadric
hypersurface, observe that the quadric also admits an action
of \(O(n+1)\) (inherited from \(\CC^{n+1}\)) precisely because
the quadratic form \(\sum_{k=1}^{n+1}z_k^2\) is preserved by
orthogonal matrices. This is (a) transitive and (b) the
stabiliser of the point \([1:i:\cdots:0]\) is isomorphic to
\(SO(2)\times O(n-1)\).

(a) To see transitivity, let us prove that the orbit of
\([1:i:0:\cdots:0]\) is the whole quadric. Suppose that
\(\bm{z}=\bm{x}+i\bm{y}\in\CC^{n+1}\) is a nonzero complex
vector with \(\sum_{k=1}^{n+1} z_k^2=0\). The real and
imaginary parts of this condition become
\(|\bm{x}|^2=|\bm{y}|^2\) and \(\bm{x}\cdot\bm{y}=0\). Use the
Gram-Schmidt process to extend
\(\hat{\bm{x}}=\bm{x}/|\bm{x}|\),
\(\hat{\bm{y}}=\bm{y}/|\bm{y}|\) to an orthonormal basis of
\(\RR^{n+1}\) use these basis vectors as the columns of an
orthogonal matrix \(A\). By construction,
\(A(1,i,0,\ldots,0)=\hat{\bm{x}}+i\hat{\bm{y}}\), so
\(A[1:i:0:\cdots:0]=[\hat{\bm{z}}]=[\bm{z}]\). This shows that
\([\bm{z}]\) lies in the \(O(n+1)\)-orbit of
\([1:i:0\cdots:0]\).

(b) To understand the stabiliser, suppose that
\(A[1:i:0:\cdots:0]=[1:i:0:\cdots:0]\). If the columns of
\(A\) are \(\bm{a}_1,\ldots,\bm{a}_{n+1}\) then this condition
becomes \(\bm{a}_1+i\bm{a}_2=re^{i\theta}(1,i,0,\ldots,0)\)
for some \(re^{i\theta}\in\CC\setminus\{0\}\). Since \(A\) is
orthogonal, we get \(r=1\) and \[\bm{a}_1=\begin{pmatrix}
\cos\theta\\ -\sin\theta\\ 0\\ \vdots\\0\end{pmatrix},\qquad
\bm{a}_2=\begin{pmatrix}
\sin\theta\\ \cos\theta\\ 0\\ \vdots\\0\end{pmatrix}\] The
upper-right 2-by-2 block of \(A\) is therefore an element of
\(SO(2)\). Orthogonality of \(A\) now implies that \(A\) is
block-diagonal and the lower-left \((n-1)\)-by-\((n-1)\) block
is orthogonal. Thus \(A\in SO(2)\times O(n-1)\).\qedhere

\end{Solution}
\begin{Exercise}[Remark \ref{rmk:fixing_periods}]\label{exr:fixing_periods}
The zero-section \(S^n\subset T^*S^n\) consists of fixed
points of the cogeodesic flow, but all the other orbits of the
cogeodesic flow on the round sphere are circles. By a similar
argument to Theorem \ref{thm:2dactionangle}, one can modify
\(H=\frac{1}{2}|\eta|^2\) away from the zero-section to get a
Hamiltonian for which all orbits have the same period (i.e.\
giving a circle action away from the zero-section). What
Hamiltonian should you take instead?
\end{Exercise}
\begin{Solution}
The geodesics in the level set \(H^{-1}(c)\) have speed
\(|p|=\sqrt{2c}\), so have period \(T(c)=2\pi/\sqrt{2c}\). The proof
of Theorem \ref{thm:2dactionangle} tells us to use the Hamiltonian
\(\alpha\circ H\) where \(\alpha(b)=\frac{1}{2\pi}\int_0^b
T(c)\,dc=\frac{1}{\sqrt{2}}\int_0^b\frac{dc}{\sqrt{c}}=\sqrt{2b}\). In
other words, use the Hamiltonian \(|p|\). This Hamiltonian is not
smooth at \(p=0\) (i.e.\ along the zero-section) so the new
Hamiltonian flow only makes sense away from the
zero-section. \qedhere

\end{Solution}
\begin{Exercise}[Lemma \ref{lma:circle_action_symp_cut}]\label{exr:circle_action_symp_cut}
The Hamiltonian \(\tilde{H}_c\) generates the circle
action
\[\left(x,\xi\right)\mapsto \left(\phi^H_t(x),
e^{-it}\xi\right)\] on \(X\times\CC\).
\end{Exercise}
\begin{Solution}
The Hamiltonian vector field is \(V_H-\partial_\theta\) where
\(\theta\) is the angular coordinate on \(\CC\). The flowlines
are therefore \((\phi^H_t(x), e^{-it}\xi)\), which all have
period \(2\pi\). \qedhere

\end{Solution}
\begin{Exercise}[Remark \ref{rmk:why_z_linear}]\label{exr:why_z_linear}
Why do we need \(s\) to be \(\ZZ\)-linear (or at least
\(\QQ\)-linear) for the construction in Example
\ref{exm:toric_one_cut} to work?
\end{Exercise}
\begin{Solution}
If \(s\) is \(\ZZ\)-linear then the Hamiltonian
\(H:=s\circ\mu\) generates a circle action: the Hamiltonian
vector field \(V_H\) is a \(\ZZ\)-linear combination of the
periodic vector fields generating the torus action. If \(s\)
is only \(\QQ\)-linear, we can rescale it to clear
denominators and get a \(\ZZ\)-linear map, so the symplectic
cut can still be made to work. If \(s\) is not \(\QQ\)-linear
then the subgroup of \(T^n\) generated by the flow of \(V_H\)
is not closed. If we persist in taking the quotient, the
result will likely fail to be Hausdorff. \qedhere

\end{Solution}
\begin{Exercise}[Theorem \ref{thm:toric_from_cut}]\label{exr:toric_from_cut}
Any convex rational polytope \(\Delta\) occurs as the moment
image of a toric Hamiltonian system (on a possibly singular
space).
\end{Exercise}
\begin{Proof}
Start with the Hamiltonian system \(\mu\colon T^*T^n\to\RR^n\)
given in canonical coordinates by \((\bm{p},\bm{q})\to
\bm{p}\). The moment image is the whole of \(\RR^n\). The
polytope \(\Delta\) is an intersection of a collection of
half-spaces \(s_j(x)\geq c_j\) where \(s_1,\ldots,s_m\) are
\(\QQ\)-linear maps and \(c_1,\ldots,c_m\) are real
numbers. Take the symplectic cut of \(T^*T^n\) by
\(s_1\circ\mu\) at level \(c_1\). Then take the symplectic cut
of the result by \(s_2\circ\mu\) at level \(c_2\), and
continue. Each time you cut, the moment image is intersected
with another half-space. The final result is a Hamiltonian
system whose moment image is \(\Delta\). The total space will
have singularities if \(\Delta\) is not Delzant. This will
happen when we quotient by a non-free circle action. \qedhere

\end{Proof}
\begin{Exercise}[Example \ref{exm:2_ball_3_ball}]\label{exr:2_ball_3_ball}
There is a blow-up of \(\cp{1}\times\cp{1}\) in two disjoint
symplectic balls which is symplectomorphic to a blow-up of
\(\cp{2}\) in three disjoint balls.
\end{Exercise}
\begin{Proof}
Blow up the shaded balls by symplectic
cut\index{symplectic cut|)}. This has the effect of truncating
the moment polygons to obtain the moment
hexagon\index{moment polytope!hexagon} of the common blow-up.

\begin{center}
\begin{tikzpicture}
\begin{scope}[shift={(-3,0)}]
\filldraw[fill=gray,opacity=0.5] (1,-1) -- (1,1) -- (-1,1) -- (-1,-1) -- cycle;
\draw[thick] (1,-1) -- (1,1) -- (-1,1) -- (-1,-1) -- cycle;
\filldraw[dotted,fill=darkgray] (-1,-1) -- (0,-1) -- (-1,0) -- cycle;
\filldraw[dotted,fill=darkgray] (1,1) -- (0,1) -- (1,0) -- cycle;
\end{scope}
\filldraw[fill=lightgray,opacity=0.5] (1,0) -- (0,1) -- (-1,1) -- (-1,0) -- (0,-1) -- (1,-1) -- cycle;
\draw[thick] (1,0) -- (0,1) -- (-1,1) -- (-1,0) -- (0,-1) -- (1,-1) -- cycle;
\begin{scope}[shift={(3,0)}]
\filldraw[fill=lightgray,opacity=0.5] (-1,-1) -- (2,-1) -- (-1,2) -- cycle;
\draw[thick] (-1,-1) -- (2,-1) -- (-1,2) -- cycle;
\filldraw[dotted,fill=darkgray] (-1,-1) -- (0,-1) -- (-1,0) -- cycle;
\filldraw[dotted,fill=darkgray] (2,-1) -- (1,0) -- (1,-1) -- cycle;
\filldraw[dotted,fill=darkgray] (-1,2) -- (-1,1) -- (0,1) -- cycle;
\end{scope}
\end{tikzpicture}
\end{center}

\end{Proof}
\begin{Exercise}[Example \ref{exm:dp_as_slice}]\label{exr:dp_as_slice}
Show that the common blow-up from Example
\ref{exm:2_ball_3_ball} arises as a symplectic
reduction\index{symplectic reduction|)} of
\(\cp{1}\times\cp{1}\times\cp{1}\).
\end{Exercise}
\begin{Proof}
Make the cut as shown.

\begin{center}
\begin{tikzpicture}
\filldraw[fill=lightgray,draw=none,opacity=0.5] (1.5,-0.5) -- (3.3,-0.2) -- (1.8,0.3) -- (0,0) -- cycle;
\filldraw[fill=lightgray,draw=none,opacity=0.5] (1.8,0.3) -- (0,0) -- (0,2) -- (1.8,2.3) -- cycle;
\filldraw[fill=lightgray,draw=none,opacity=0.5] (1.5,-0.5) -- (0,0) -- (0,2) -- (1.5,1.5)-- cycle;
\filldraw[fill=lightgray,draw=none,opacity=0.5] (1.5,-0.5) -- (3.3,-0.2) -- (3.3,1.8) -- (1.5,1.5) -- cycle;
\filldraw[fill=lightgray,draw=none,opacity=0.5] (1.8,2.3) -- (0,2) -- (1.5,1.5) -- (3.3,1.8) -- cycle;
\filldraw[fill=lightgray,draw=none,opacity=0.5] (1.8,0.3) -- (3.3,-0.2) -- (3.3,1.8) -- (1.8,2.3) -- cycle;
\draw[black,thick,dotted] (1.5,-0.5) -- (3.3,-0.2) -- (1.8,0.3) -- (0,0) -- cycle;
\draw[black,thick,dotted] (1.8,0.3) -- (0,0) -- (0,2) -- (1.8,2.3) -- cycle;
\draw[black,thick] (1.5,-0.5) -- (0,0) -- (0,2) -- (1.5,1.5)-- cycle;
\draw[black,thick] (1.5,-0.5) -- (3.3,-0.2) -- (3.3,1.8) -- (1.5,1.5) -- cycle;
\draw[black,thick] (1.8,2.3) -- (0,2) -- (1.5,1.5) -- (3.3,1.8) -- cycle;
\draw[black,thick,dotted] (1.8,0.3) -- (3.3,-0.2) -- (3.3,1.8) -- (1.8,2.3) -- cycle;
\filldraw[dotted,fill=darkgray,opacity=0.5] (0,1) -- (0.75,1.75) --(2.4,1.65) -- (3.3,0.8) -- (2.55,0.05) -- (0.9,0.15) -- cycle;
\end{tikzpicture}
\end{center}

This gives a hexagonal moment polygon which is \(\ZZ\)-affine
equivalent to the previous one (e.g.\ via projection to the
\(xy\)-plane). \qedhere

\end{Proof}
\chapter{Visible Lagrangian submanifolds}
\label{ch:visible_Lagrangian}
\thispagestyle{cup}

We will now study Lagrangian submanifolds of toric manifolds. It
will turn out that if the moment image of Lagrangian submanifold
has codimension \(k\) then it is contained in a affine subspace
of codimension \(k\). This does not quite determine the
Lagrangian completely, but gives severe restrictions. Just as we
have been specifying a symplectic manifold by drawing a
polytope, we will be able to specify a Lagrangian submanifold by
drawing an affine subspace of the moment polytope; these are
called {\em visible
Lagrangians}\index{Lagrangian!submanifold!visible|(}. While most
Lagrangian submanifolds of toric varieties are not visible, the
visible ones are useful to know about, and we discuss the theory
and numerous examples in this chapter. In Appendix
\ref{ch:tropical_lag}, we will see a more versatile construction
due to Mikhalkin and Matessi, which assigns a {\em tropical
Lagrangian} to a {\em tropical curve} in the polytope.

\section{Visible Lagrangian submanifolds}

\begin{Theorem}\label{thm:visibility}
Consider the integrable Hamiltonian system \(\bm{H}\colon
\RR^n\times T^n\to\RR^n\), \(\bm{H}(\bm{p},\bm{q})=\bm{p}\)
where \(q_1,\ldots,q_n\) are taken modulo \(2\pi\) and the
symplectic form is \(\sum dp_i\wedge dq_i\). Let
\(L\subset\RR^n\times T^n\) be a Lagrangian
submanifold. Suppose that \(\bm{H}|_L\colon L\to\RR^n\)
factors as \(\bm{H}|_L=f\circ g\), where \(g\colon L\to K\) is
a bundle over a \(k\)-dimensional manifold \(K\), \(k<n\), and
\(f\colon K\to\RR^n\) is an embedding. Then \(K\) is an affine
linear subspace of \(\RR^n\) which is rational with respect to
the lattice \((2\pi\ZZ)^n\).

\end{Theorem}
\begin{Definition}\label{dfn:visible}
We call Lagrangian submanifolds which project in this way {\em
visible}.

\end{Definition}
\begin{Remark}
Theorem \ref{thm:visibility} was first observed when \(n=2\)
and \(\dim(K)=1\) by Symington {\cite[Corollary
7.9]{Symington}}.

\end{Remark}
\begin{Proof}
Let \(\bm{s}=(s_1,\ldots,s_k)\) be local coordinates on \(K\)
and \(\bm{t}=(t_{k+1},\ldots,t_n)\) be local coordinates on
the fibre of \(g\). By assumption, the inclusion of \(L\)
into \(\RR^n\) has the form \((\bm{s},\bm{t})\mapsto
(\bm{p}(\bm{s}),\bm{q}(\bm{s},\bm{t}))\) for some functions
\(\bm{p},\bm{q}\). The vectors \(\partial_{s_i}\) and
\(\partial_{t_j}\) pushforward to
\((\partial_{s_i}\bm{p},\partial_{s_i}\bm{q})\) and
\((0,\partial_{t_j}\bm{q})\). The Lagrangian condition on
\(L\) is equivalent to
\(\partial_{s_i}\bm{p}\cdot\partial_{t_j}\bm{q}=0\) and
\(\partial_{s_i}\bm{p} \cdot \partial_{s_j}\bm{q} =
\partial_{s_j}\bm{p}\cdot\partial_{s_i}\bm{q}\) for all
\(i,j\). The first of these conditions implies that the
tangent space of the fibre of \(g\) is
orthogonal\footnote{with respect to the Euclidean metric on
\(\RR^n\).} to the \(k\)-dimensional subspace \(f_*(TK)\)
spanned by \(\partial_{s_1}\bm{p}, \ldots,
\partial_{s_k}\bm{p}\). Since the tangent space of the fibre
of \(g\) is \((n-k)\)-dimensional, it must be precisely
\(f_*(TK)^{\perp}\); in other words, for each \(\bm{s}\in K\),
the fibre of \(g\) over \(\bm{s}\) is an integral submanifold
of the distribution on \(T^n\) given by
\(f_*(TK)^{\perp}\). This distribution has an integral
submanifold if and only if \(f_*(TK)\) is a rational subspace
with respect to the lattice \((2\pi\ZZ)^n\). Since \(f_*(TK)\)
varies smoothly in \(\bm{s}\), and must always be rational, it
is necessarily constant. Therefore \(f(K)\) is a rational
affine subspace. \qedhere

\end{Proof}
\begin{Remark}\label{rmk:visible_orthogonal}
As a consequence of the proof, we see that if the visible
Lagrangian projects to an affine subspace \(K\) in the
\(\bm{p}\)-plane then its fibre in the \(\bm{q}\)-torus above
a point in \(K\) is a translate of the subtorus
\(K^\perp/(K^\perp\cap(2\pi\ZZ)^n)\).

\end{Remark}
\begin{Example}\label{exm:visible_orthogonal}
Suppose \(n=2\) and \(K\) is the \(p_1\)-axis. Then
\(L\cap\{(p_1,0)\}\) is a circle
\(\{(q_1,\theta)\,:\,\theta\in[0,2\pi]\}\) for some fixed
\(q_1\). For example, \(L\) could be the cylinder
\(\{p_2=0,\quad q_1=0\}\).

\end{Example}
\begin{Remark}\label{rmk:qsdependence}
Note that the dependence of \(q_i\) on the coordinates \(s_j\) can
be nontrivial.

\end{Remark}
\begin{Example}\label{exm:qsdependence}
Let \((p_1,p_2,q_1,q_2)\) be coordinates on \(X=\RR^2\times
T^2\) with symplectic form \(\sum dp_i\wedge dq_i\). The
Lagrangian embedding \(i\colon\RR\times S^1\to X\),
\(i(s,t)=(s,0,0,t)\) is visible for the projection
\((\bm{p},\bm{q})\mapsto \bm{p}\). The Lagrangian torus
\(j\colon S^1\times S^1\to X\), \(j(s,t)=(\sin s,0,s,t)\) is
also visible\footnote{Technically, it is not visible itself
because the projection map is not a bundle, rather it is a
union of two visible cylinders. We will tolerate this and
related abuses of terminology.}, and projects to the line
segment \([-1,1]\times\{0\}\) (the preimage of each point in
\((-1,1)\times\{0\}\) is a pair of circles).

\end{Example}
\begin{Remark}\label{rmk:visible_lag_to_aff_str}
Apart from giving a useful way to visualise and construct
Lagrangian submanifolds, this theorem also gives us a way to
figure out the integral affine structure on the base of a
Lagrangian torus fibration if we don't already know it. If we
can find a Lagrangian submanifold whose image under our
Hamiltonian system is a submanifold \(K\subset\RR^n\) then we
know that the image of \(K\) under action coordinates is
supposed to be affine linear. We will use this observation in
the proof of Lemma \ref{lma:zaffexample} later.

\end{Remark}
\section{Hitting a vertex}

Suppose now that we have a Hamiltonian torus action (and toric
critical points) with moment map \(\mu\colon X\to\RR^n\) and
address the question of what visible Lagrangian surfaces look
like when the affine linear subspace \(\mu(L)\) intersects the
boundary strata of the moment polytope. For simplicity, we will
focus on the case \(\dim X=4\), \(\dim\mu(L)=1\).

\begin{Example}\label{exm:hittingvertex_1}
Consider the Lagrangian plane \(L:=\{(z,\bar{z})\ :\
z\in\CC\}\subset\CC^2\). The projection \(\mu(L)\) is the diagonal
ray \(\{(t,t)\ :\ t\in[0,\infty)\}\subset\RR^2\), so \(L\) is a
visible Lagrangian surface. See Figure \ref{fig:hittingvertex}.

\end{Example}
\begin{Example}\label{exm:antidiagonal}
Consider the Lagrangian {\em antidiagonal
sphere}\index{antidiagonal sphere}
\[\bar{\Delta}:=\{((x,y,z),(-x,-y,-z))\in S^2\times S^2\ :
(x,y,z)\in S^2\}.\] Here, we have equipped \(S^2\times S^2\)
with the equal-area symplectic form from Example
\ref{exm:products_of_spheres}. The moment map is
\(\mu((x_1,y_1,z_1),(x_2,y_2,z_2))=(z_1,z_2)\) so the
projection of the antidiagonal sphere along \(\mu\) is the
antidiagonal line \(\{(z,-z)\in[-1,1]^2\,:\,z\in[-1,1]\}\)
(see Figure \ref{fig:hittingvertex}). This is therefore a
visible Lagrangian whose projection hits two vertices, where
it is locally modelled on Example \ref{exm:hittingvertex_1}.

\end{Example}
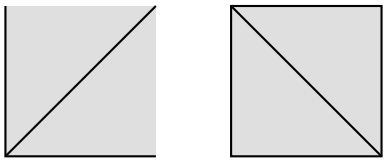
\begin{figure}
\begin{center}
\begin{tikzpicture}
\filldraw[fill=lightgray,opacity=0.5,draw=none] (0,0) -- (2,0) -- (2,2) -- (0,2) -- cycle;
\draw[thick,black] (0,2) -- (0,0) -- (2,0);
\draw[thick,black] (0,0) -- (2,2);
\begin{scope}[shift={(3,0)}]
\filldraw[fill=lightgray,opacity=0.5,draw=none] (0,0) -- (2,0) -- (2,2) -- (0,2) -- cycle;
\draw[thick,black] (0,0) -- (2,0) -- (2,2) -- (0,2) -- cycle;
\draw[thick,black] (2,0) -- (0,2);
\end{scope}
\end{tikzpicture}
\end{center}
\caption{Left: A visible Lagrangian disc (Example \ref{exm:hittingvertex_1}). Right: The antidiagonal sphere in \(S^2\times S^2\) (Example \ref{exm:antidiagonal}) is a visible Lagrangian living over the antidiagonal in the square.}
\label{fig:hittingvertex}
\end{figure}

\begin{Example}[Exercise \ref{exr:hittingvertex}]\label{exm:hittingvertex}
Fix \(m,n\in\ZZ_{> 0}\) with \(\gcd(m,n)=1\) Consider the ray
\(\{(mt,nt)\ :\ t\in[0,\infty)\}\) in the nonnegative quadrant
(Figure \ref{fig:schoenwolfson}). Above this ray is a visible
Lagrangian which I will call a\index{Schoen-Wolfson cone} {\em
Schoen-Wolfson cone}\footnote{Schoen and Wolfson
{\cite[Theorem 7.1]{SchoenWolfson}} showed that these are the
only Lagrangian cones in \(\CC^2\) which are Hamiltonian
stationary (i.e.\ critical points of the volume functional
restricted to Hamiltonian deformations).} given
parametrically by: \[(s,t) \mapsto \frac{1}{\sqrt{m+n}}
\left(t\sqrt{m}e^{is\sqrt{n/m}},
it\sqrt{n}e^{-is\sqrt{m/n}}\right),\quad
s\in[0,2\pi\sqrt{mn}],\,t\in[0,\infty)\] This cone is singular
at the origin unless \(m = n = 1\).

\begin{figure}
\begin{center}
\begin{tikzpicture}
\filldraw[fill=lightgray,opacity=0.5,draw=none] (0,2) -- (0,0) -- (2,0) -- (2,2);
\draw[thick,black] (0,2) -- (0,0) -- (2,0);
\draw[thick,black] (0,0) -- (1,2);
\node at (1,2) [above] {\((m,n)\)};
\end{tikzpicture}
\end{center}
\caption{Moment image of a Schoen-Wolfson cone.}
\label{fig:schoenwolfson}
\end{figure}
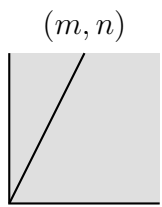

\end{Example}
\begin{Remark}
Modulo the freedom discussed in Remark \ref{rmk:qsdependence}
and Example \ref{exm:qsdependence}, this exhausts all possible
local models for visible Lagrangians living over a line which
hits the corner of a Delzant moment polygon.

\end{Remark}
\begin{Example}
If our moment polygon is a rectangle with sidelengths \(m\)
and \(n\) (positive integers) then we get a symplectic form on
\(S^2\times S^2\) which gives the factors symplectic area
\(m\) and \(n\) respectively, so the symplectic form lives in
the class \((m,n)\in H^2(S^2\times S^2)\). The diagonal line
joining opposite corners of the moment rectangle is the
projection of a visible Lagrangian
sphere\index{Lagrangian!sphere} \(L\) with two Schoen-Wolfson
singular points with parameters \((m, n)\). The homology class
of this Lagrangian is \((n,-m)\) (you can see this by
intersecting with spheres in the classes \([S^2\times\{p\}]\)
and \([\{p\}\times S^2]\)), which has symplectic area
\(0\). Indeed, the homology class \((n,-m)\) can only contain
a Lagrangian representative when \([\omega]\) is a multiple of
\((m,n)\in H^2(S^2\times S^2)\).

\end{Example}
\section{Hitting an edge}

\begin{figure}[htb]
\begin{center}
\begin{tikzpicture}
\filldraw[draw=none,fill=lightgray,opacity=0.5] (0,0) -- (3,0.5) -- ++ (0,1) -- ++ (-3,-0.5) -- cycle;
\filldraw[draw=none,fill=lightgray,opacity=0.5] (0,0) -- (3,0.5) -- ++ (-150:1) -- ++ (-3,-0.5) -- cycle;
\filldraw[draw=none,fill=lightgray,opacity=0.5] (0,0) -- (3,0.5) -- ++ (-30:1) -- ++ (-3,-0.5) -- cycle;
\draw (0,0) -- (3,0.5);
\draw[->-] (0,0) -- (90:1);
\draw[->>-] (0,0) -- (-30:1);
\draw[->>>-] (0,0) -- (-150:1);
\begin{scope}[shift={(3,0.5)}]
\draw[->>-] (0,0) -- (90:1);
\draw[->>>-] (0,0) -- (-30:1);
\draw[dotted,->-] (0,0) -- (-150:1);
\end{scope}
\end{tikzpicture}
\caption{A pinwheel core with \(n=3\) flanges: lines with arrows should be identified in pairs.}
\label{fig:pinwheel_core}
\end{center}
\end{figure}
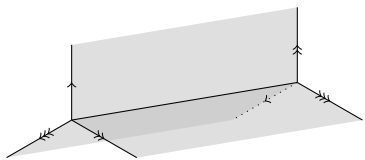

\begin{Example}\label{exm:hittingedge}
We now consider visible Lagrangians whose projection hits an
edge. For a local model, we take \(X=\RR\times S^1\times\CC\),
with coordinates \((p,q,z=x+iy)\) (\(q\in\RR/2\pi\ZZ)\) and
symplectic form \(dp\wedge dq+dx\wedge dy\). The image of the
moment map \(\mu\colon X\to\RR^2\),
\(\mu(p,q,z)=\left(p,\frac{1}{2}|z|^2\right)\) is the closed
upper half-plane \(\{(x_1,x_2)\in\RR^2\ :\ x_2\geq
0\}\). Consider the ray \(R_{m,n}=\{(ms,ns)\ :\ s\geq
0\}\). The following map is a Lagrangian immersion of the
cylinder \[i(s,t)=\left(ms,-nt,\sqrt{2ns}e^{imt}\right),\qquad
(s,t)\in[0,\infty)\times S^1\] whose projection along \(\mu\)
is the ray \(R_{m,n}\). This immersion is an embedding away
from \(s=0\), but it is \(n\)-to-\(1\) along the circle
\(s=0\) (the points \(\left(0,t+\frac{2\pi k}{n}\right)\),
\(k=0,\ldots,n-1\), all project to \((0,t\mod 2\pi,0)\)).

\begin{center}
\begin{tikzpicture}
\filldraw[fill=lightgray,opacity=0.5,draw=none] (-2,2) -- (-2,0) -- (2,0) -- (2,2) -- cycle;
\draw[thick,black] (-2,0) -- (2,0);
\draw[thick,black] (0,0) -- (1,2);
\node at (1,2) [above] {\(R_{m,n}\)};
\end{tikzpicture}
\end{center}

The image of the immersion is a Lagrangian which looks like a
collection of \(n\) flanges meeting along a circle, twisting
as they move around the circle so that the link of the circle
is an \((m,n)\)-torus knot (see Figure
\ref{fig:pinwheel_core}). For example, when \(m=1\), \(n=2\),
this is a M\"{o}bius strip. For \(n\geq 3\) it is not a
submanifold. We call the image of the immersion a {\em
Lagrangian \((n,m)\)-pinwheel
core}\index{Lagrangian!pinwheel core}.\index{pinwheel|see
{Lagrangian, pinwheel}}\index{pinwheel core|see {Lagrangian,
pinwheel, core}}

Any integral affine transformation preserving the upper
half-plane and fixing the origin acts on the set of rays
\(R_{m,n}\). These transformations are precisely the affine
shears \(\lmatrix 1 & 0 \\ k & 1 \rmatrix \), which allow us
to change \(m\) by any multiple of \(n\), so we can always
assume \(m\in\{0,\ldots,n-1\}\).

Again, modulo the freedom discussed in Remark
\ref{rmk:qsdependence} and Example \ref{exm:qsdependence},
these local models exhaust the visible Lagrangians
intersecting an edge of a moment polygon.

\end{Example}

\begin{Example}
Consider the Lagrangian
\(\rp{2}\)\index{Lagrangian!RP2@$\mathbb{RP}^2$} which is the
closure of the visible disc \(\{[z:\bar{z}:1]\ :\
z\in\CC\}\subset\cp{2}\). This projects to the diagonal
bisector in the moment
triangle\index{moment polytope!CP2@$\mathbb{CP}^2$} (see
Figure \ref{fig:rp2_pinwheel_analysis}(a)). If we use the
integral affine transformation \(\lmatrix -1 & -1 \\ 0 &
-1\rmatrix\) to make the slanted edge of the triangle
horizontal then the projection of the visible Lagrangian ends
up pointing in the \(\pm(1,2)\)-direction (Figure
\ref{fig:rp2_pinwheel_analysis}(b)), so comparison with
Example \ref{exm:hittingedge} shows that the disc is capped
off with a M\"{o}bius strip to give an
\(\rp{2}\).

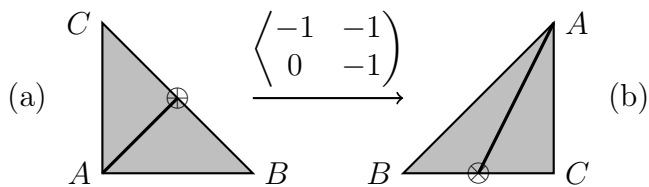
\begin{figure}[htb]
\begin{center}
\begin{tikzpicture}
\node at (-1,1) {(a)};
\node at (7,1) {(b)};
\filldraw[fill=lightgray,draw=black,thick] (0,2) node [left] {\(C\)} -- (0,0) node [left] {\(A\)} -- (2,0) node [right] {\(B\)} -- cycle;
\draw[very thick] (0,0) -- (1,1) node {\(\oplus\)};
\draw[->,thick] (2,1) -- (4,1) node [midway,above] {\(\lmatrix -1 & -1 \\ 0 & -1\rmatrix\)};
\begin{scope}[shift={(6,2)}]
\filldraw[fill=lightgray,draw=black,thick] (-2,-2) node [left] {\(B\)} -- (0,-2) node [right] {\(C\)} -- (0,0) node [right] {\(A\)} -- cycle;
\draw[very thick] (0,0) -- (-1,-2) node {\(\otimes\)};
\end{scope}
\end{tikzpicture}
\caption{(a) Visible Lagrangian \(\rp{2}\) in \(\cp{2}\). (b) The same picture after an integral affine transformation shows the line pointing in the \((1,2)\)-direction, so there is a \((2,1)\)-pinwheel core (M\"{o}bius strip) near the point marked \(\otimes\).}
\label{fig:rp2_pinwheel_analysis}
\end{center}
\end{figure}

\end{Example}
\begin{Example}[Exercise \ref{exr:kleinbot}]\label{exm:kleinbot}
The square below has vertices at \((-2,-2)\), \((-2,2)\),
\((2,-2)\), \((2,2)\). There is a smooth, closed visible
Lagrangian surface \(L\) in the corresponding toric variety,
living over the line segment connecting \((-1,-2)\) to
\((1,2)\). To which topological surface is \(L\) homeomorphic?

\begin{center}
\begin{tikzpicture}
\filldraw[fill=lightgray,opacity=0.5] (0,0) -- (2,0) -- (2,2) -- (0,2) -- (0,0);
\draw[thick] (1/2,0) -- (3/2,2);
\node at (1/2,0) [below] {\((-1,-2)\)};
\node at (3/2,2) [above] {\((1,2)\)};
\end{tikzpicture}
\end{center}

\end{Example}
Here is a higher dimensional example.

\begin{Example}[Exercise \ref{exr:nullhomologous3sphere}]\label{exm:nullhomologous3sphere}
Consider the symplectic manifold \(\cp{1}\times\CC^2\) with
the symplectic form
\(pr_1^*\omega_{\cp{1}}+pr_2^*\omega_{\CC^2}\) (here, \(pr_k\)
denotes the projection to the \(k\)th factor,
\(\omega_{\cp{1}}\) is the Fubini-Study form on \(\cp{1}\)
normalised so that
\(\frac{1}{2\pi}\int_{\cp{1}}\omega_{\cp{1}}=1\) and
\(\omega_{\CC^2}\) is the standard symplectic form). Sketch
the moment
image\index{moment polytope!CP1xC2@$\mathbb{CP}^1\times\mathbb{C}^2$}
for the \(T^3\)-action coming from the standard torus actions
on each factor. Check that the 3-sphere
\(\{([-\bar{z}_2:\bar{z}_1], z_1,
z_2)\ :\ |z_1|^2+|z_2|^2=2\}\subset\cp{1}\times\CC^2\) is
Lagrangian and sketch its projection under the moment map.

\end{Example}
\section{Solutions to inline exercises}

\begin{Exercise}[Example \ref{exm:hittingvertex}]\label{exr:hittingvertex}
Fix \(m,n\in\ZZ_{> 0}\) with \(\gcd(m,n)=1\). Verify that the
Schoen-Wolfson cone \[(s,t) \mapsto \frac{1}{\sqrt{m+n}}
\left(t\sqrt{m}e^{is\sqrt{n/m}},
it\sqrt{n}e^{-is\sqrt{m/n}}\right),\quad
s\in[0,2\pi\sqrt{mn}],\,t\in[0,\infty)\] is Lagrangian where
that makes sense (i.e.\ away from the cone point) and that its
projection under the moment map is the ray
\(\{(mt,nt)\ :\ t\in[0,\infty)\}\).
\end{Exercise}
\begin{Solution}
The parametrisation here is chosen to agree with the one from
the Schoen-Wolfson paper, but we can make our life easier by
using \(\theta=s/\sqrt{mn}\in[0,2\pi]\) and \(r=t^2/2(m+n)\)
to get the parametrisation \[(r,\theta)\mapsto (\sqrt{2mr}
e^{in\theta},\sqrt{2nr} e^{-im\theta}).\] Applying the moment
map \((|z_1|^2/2,|z_2|^2/2)\) gives us \(\{(mr,nr)\,:\,r\geq
0\}\), so the Lagrangian projects to the correct ray. The
fibre of the Lagrangian over \((mr,nr)\) is \(\{(\sqrt{2mr}
e^{in\theta},\sqrt{2nr}
e^{-im\theta})\,:\,\theta\in[0,2\pi]\}\). In the
\((\theta_1,\theta_2)\)-torus, this is a circle whose tangent
line is \((n,-m)\), which is orthogonal to the ray in the
base. Therefore this is Lagrangian by Remark
\ref{rmk:visible_orthogonal}. \qedhere

\end{Solution}
\begin{Exercise}[Example \ref{exm:kleinbot}]\label{exr:kleinbot}
The square below has vertices at \((-2,-2)\), \((-2,2)\),
\((2,-2)\), \((2,2)\). There is a smooth, closed visible
Lagrangian surface \(L\) in the corresponding toric variety,
living over the line segment connecting \((-1,-2)\) to
\((1,2)\). To which topological surface is \(L\)
homeomorphic?

\begin{center}
\begin{tikzpicture}
\filldraw[fill=lightgray,opacity=0.5] (0,0) -- (2,0) -- (2,2) -- (0,2) -- (0,0);
\draw[thick] (1/2,0) -- (3/2,2);
\node at (1/2,0) [below] {\((-1,-2)\)};
\node at (3/2,2) [above] {\((1,2)\)};
\end{tikzpicture}
\end{center}
\end{Exercise}
\begin{Solution}
There are two M\"{o}bius strips where the projection of the
visible Lagrangian meets the edge of the square. These are
joined along their common boundary, which forms a Lagrangian
Klein bottle\index{Lagrangian!Klein bottle}.\qedhere

\end{Solution}
\begin{Exercise}[Example \ref{exm:nullhomologous3sphere}]\label{exr:nullhomologous3sphere}
Consider the symplectic manifold \(\cp{1}\times\CC^2\) with
the symplectic form
\(pr_1^*\omega_{\cp{1}}+pr_2^*\omega_{\CC^2}\) (here, \(pr_k\)
denotes the projection to the \(k\)th factor,
\(\omega_{\cp{1}}\) is the Fubini-Study form on \(\cp{1}\)
normalised so that
\(\frac{1}{2\pi}\int_{\cp{1}}\omega_{\cp{1}}=1\) and
\(\omega_{\CC^2}\) is the standard symplectic form). Sketch
the moment image for the \(T^3\)-action coming from the
standard torus actions on each factor. Check that the 3-sphere
\(L:=\{([-\bar{z}_2:\bar{z}_1], z_1,
z_2)\ :\ |z_1|^2+|z_2|^2=2\}\subset\cp{1}\times\CC^2\) is
Lagrangian and sketch its projection under the moment map.
\end{Exercise}
\begin{Solution}
The moment map is \(\mu([a:b],z_1,z_2) =
\left(\frac{|b|^2}{|a|^2+|b|^2}, \frac{1}{2}|z_1|^2,
\frac{1}{2}|z_2|^2\right)\), so its image is the noncompact
polytope \(\{(x,y,z)\in\RR^3\,:\,x\in[0,1],\ y,z\geq 0\}\). To
compute the moment image of \(L\), we have
\[\mu([-\bar{z}_2:\bar{z}_1], z_1, z_2) =
\left(\frac{1}{2}|z_1|^2, \frac{1}{2}|z_1|^2,
\frac{1}{2}(2 - |z_1|^2)\right),\] where we used the fact that
\(|z_1|^2+|z_2|^2=2\). As \(|z_1|^2\) varies between \(0\) and
\(2\) (again using the constraint \(|z_1|^2+|z_2|^2=2\)) we
get the straight line segment \[t\mapsto
(t/2,t/2,(2-t)/2),\qquad t\in[0,2]\] in the \(y+z=1\) plane,
connecting \((0,0,1)\) to \((1,1,0)\).

\begin{figure}[htb]
\begin{center}
\begin{tikzpicture}
\draw (0,0) -- (10:3);
\draw (0,0) -- (0,3);
\draw (0,0) -- (-30:3);
\begin{scope}[shift={(10:3)}]
\draw[->] (0,0) -- (0,3) node [above] {\(z\)};
\draw[->] (0,0) -- (-30:3) node [below right] {\(y\)};
\end{scope}
\draw[dotted,thick] (0,2) -- ++ (10:3) coordinate (b);
\draw[dotted,thick] (-30:2) -- ++ (10:3) coordinate (a);
\draw[dotted,thick] (a) -- (b);
\draw[dotted,thick] (0,2) -- (-30:2);
\draw[very thick] (b) -- (-30:2);
\draw[->] (0,0) -- (190:1) node [left] {\(x\)};
\node at (10:1.5) [above] {\(1\)};
\node at (0,1) [left] {\(1\)};
\node at (-30:1) [below left] {\(1\)};
\end{tikzpicture}
\end{center}
\caption{The moment image (slanted line) of the Lagrangian sphere in Exercise \ref{exr:nullhomologous3sphere}. The dotted rectangle is the plane \(y+z=1\) and numbers indicate affine lengths of edge-segments.}
\label{fig:nullhomologous3sphere}
\end{figure}
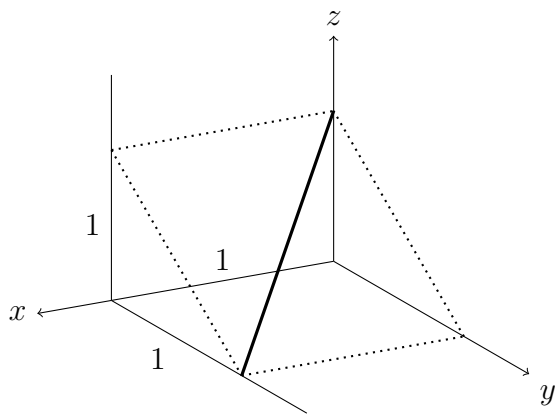

We now check that \(L\) is Lagrangian
\index{Lagrangian!submanifold!visible|)} for the symplectic
form \(pr_1^*\omega_{\cp{1}}+pr_2^*\omega_{\CC^2}\). By
definition of the Fubini-Study form, it suffices to lift the
embedding \(L\to\cp{1}\times\CC^2\) to a map \(L\to
S^3(\sqrt{2})\times\CC^2\) where \(S^3(\sqrt{2})\) is the
sphere of radius \(\sqrt{2}\) in \(\CC^2\), since
\((\cp{1},\omega_{\cp{1}})\) is obtained by the symplectic
reduction \(S^3(\sqrt{2})/S^1\). We choose the lift
\(\ell\colon(z_1,z_2)\mapsto(-\bar{z}_2,\bar{z}_1,z_1,z_2)\)
restricted to \(|z_1|^2+|z_2|^2=2\).

To begin, we work on the solid torus \(|z_1|^2\leq 1\) inside
\(|z_1|^2+|z_2|^2=2\): a similar argument holds on the
complementary solid torus \(|z_2|^2\leq 1\). Pick coordinates
\(r,\theta,\phi\) so that \(z_1=re^{i\theta}\) and
\(z_2=\sqrt{2-r^2}e^{i\phi}\). We write out \(\ell\) fully in
these coordinates (separating real and imaginary parts):
\begin{align*}\ell(r,\theta,\phi)&= \left(-\sqrt{2-r^2}\cos\phi,
\sqrt{2-r^2}\sin\phi, r\cos\theta, r\sin\theta,\right.\\
&\ \ \ \ \ \ \left. r\cos\theta,
-r\sin\theta, \sqrt{2-r^2}\cos\phi,
\sqrt{2-r^2}\sin\theta\right)\end{align*} This gives:
\begin{align*}\ell_*\partial_r &=
\left(\frac{r\cos\phi}{\sqrt{2-r^2}},
-\frac{r\sin\phi}{\sqrt{2-r^2}}, \cos\theta,
\sin\theta,\right.\\
&\ \ \ \ \ \ \left.\cos\theta, -\sin\theta,
-\frac{r\cos\phi}{\sqrt{2-r^2}},
-\frac{r\sin\phi}{\sqrt{2-r^2}}\right)\end{align*}
\begin{align*}
\ell_*\partial_\theta&=\left(0,0,-r\sin\theta,r\cos\theta,\right.\\
&\left.\ \ \ \ \ \ -r\sin\theta,-r\cos\theta,0,0\right)\end{align*}
\begin{align*}
\ell_*\partial_\phi&=\left(\sqrt{2-r^2}\sin\phi,\sqrt{2-r^2}\cos\phi,0,0,\right.\\
&\left.\ \ \ \ \ \ 0,0,-\sqrt{2-r^2}\sin\phi,\sqrt{2-r^2}\cos\phi\right)
\end{align*} from which one can check that all possible
evaluations of \(\omega\) vanish, for example
\[\omega(\ell_*\partial_r,\ell_*\partial_\phi) = r\cos^2\phi +
r\sin^2\phi - r\cos^2\phi - r\sin^2\phi=0.\qedhere\]

\end{Solution}
\chapter{Focus-focus singularities}
\label{ch:focusfocus}
\thispagestyle{cup}

So far, we have studied Lagrangian torus fibrations which have
either no critical points or else toric critical points. In this
chapter, we discuss another kind of critical point: the {\em
focus-focus critical point}\index{focus-focus!point|see
{critical point, focus-focus}}\index{critical
point!focus-focus|(}. Lagrangian torus fibrations whose critical
points are either toric or focus-focus type are called {\em
almost toric}. Allowing these critical points will drastically
expand our zoo of examples, but also make our life more
complicated. The main new features are that (a) the integral
affine structure on the base of the fibration has nontrivial
affine monodromy around the critical points, and (b) the
integral affine base no longer uniquely determines the torus
fibration.

\section{Focus-focus critical points}

\begin{Example}[Standard focus-focus system]\label{exm:ff}
Consider the following pair of
Hamiltonians\index{focus-focus!system, standard} on
\((\RR^4,dp_1\wedge dq_1+dp_2\wedge dq_2)\):
\[F_1=-p_1q_1-p_2q_2,\qquad F_2=p_2q_1-p_1q_2.\] If we
introduce complex coordinates\footnote{These complex
coordinates are not supposed to be compatible with \(\omega\),
indeed the \(p\)-plane and \(q\)-plane are both Lagrangian.}
\(p=p_1+ip_2\), \(q=q_1+iq_2\) then \(F:=F_1+iF_2=-\bar{p}q\).

\end{Example}
\begin{Lemma}[Exercise \ref{exr:ff_flow}]\label{lma:ff_flow}
The Hamiltonians \(F_1\) and \(F_2\) Poisson-commute. The
Hamiltonian \(F_1\) generates the \(\RR\)-action
\((p,q)\mapsto(e^tp,e^{-t}q)\). The Hamiltonian \(F_2\)
generates the circle action \((p,q)\mapsto(e^{it}p,e^{it}q)\).

\end{Lemma}
The orbits of the resulting \(\RR\times S^1\)-action are: the
origin (fixed point); the Lagrangian cylinders \(P:=\{(p,0)\ :\
p\neq 0\}\) and \(Q:=\{(0,q)\ :\ q\neq 0\}\); and the Lagrangian
cylinders \(\{(p,q)\ :\ \bar{p}q=c\}\) for
\(c\in\CC\setminus\{0\}\).

The diagram below represents the projection of \(\RR^4\) to
\(\RR^2\) via \[(p_1,p_2,q_1,q_2)\mapsto (|p|,|q|);\] the
projections of the \(\phi^{F_1}_t\)-flowlines are the
hyperbolae; \(\phi^{F_2}_t\)-flowlines project to points. The
Lagrangian cylinders \(P\) and \(Q\) are shown living over the
axes, the fixed point is marked with a dot at the origin.

\begin{center}
\begin{tikzpicture}
\node at (3,-0.25) {\(|p|\)};
\node at (-0.25,3) {\(|q|\)};
\filldraw[fill=lightgray,draw=none,domain=1:3.0] (0,3) -- (1,3) plot (\x,{3/(\x)}) --(3,0) -- (0,0) -- (0,3);
\draw[thick] (0,0) -- (0,3);
\node at (1,-0.25) {\(P\)};
\draw[thick] (0,0) -- (3,0);
\node at (-0.25,1) {\(Q\)};
\draw[->-,thick,domain=1:3.0] plot (\x,{3/(\x)});
\draw[->-,thick,domain=0.5:3.0] plot (\x,{1.5/(\x)});
\draw[->-,thick,domain=0.2:3.0] plot (\x,{0.6/(\x)});
\draw[->-,thick,domain=0.04:3.0] plot (\x,{0.12/(\x)});
\node at (0,0) {\(\bullet\)};
\end{tikzpicture}
\end{center}

\begin{Definition}
A {\em focus-focus chart}\index{focus-focus!chart} for an
integrable Hamiltonian system \(\bm{H}\colon X\to\RR^2\) is a
pair of embeddings \(E\colon U\to X\) and \(e\colon
V\to\RR^2\) where:
\begin{itemize}
\item \(U\subset\RR^4\) is a neighbourhood of
the origin and \(E^*\omega=\sum dp_i\wedge dq_i\),
\item \(V=\bm{F}(U)\), where \(\bm{F}\) is the Hamiltonian system in Example
\ref{exm:ff}),
\item \(\bm{H}\circ E=e\circ \bm{F}\).
\end{itemize}
We say that \(\bm{H}\colon X\to\RR^2\) has a {\em focus-focus
critical point} at \(x\in X\) if there is a focus-focus chart
\((E,e)\) with \(E(0)=x\).

\end{Definition}
\begin{Remark}\label{rmk:ffnf}
This is not the standard definition of a focus-focus critical
point: usually one specifies that \(\bm{H}\) has a critical
point at \(x\) and that the subspace of the space of quadratic
forms spanned by the Hessians\index{Hessian} of the components
\(\bm{H}\) at \(x\) agrees with the corresponding subspace for
\(\bm{F}\) at \(0\). The fact that these two definitions are
equivalent is a special case of {\em Eliasson's normal form
theorem}\index{Eliasson normal form theorem} for
non-degenerate critical points of Hamiltonian systems. For a
proof of this special case, see \cite{Chaperon}.

\end{Remark}
\begin{Lemma}[{\cite[Proposition 6.2]{NgocBohrSommerfeld}}]\label{lma:pinched_torus}
Let \(\bm{H}\colon X\to\RR^2\) be an integrable Hamiltonian
system with a focus-focus critical point \(x\) over the origin
and no other critical points. The fibre \(\bm{H}^{-1}(0)\) is
homeomorphic to a pinched torus\index{pinched torus}.
\end{Lemma}
\begin{Proof}
Recall that all our integrable systems are assumed to have
compact, connected fibres. The fibre \(\bm{H}^{-1}(0)\) is a
union of orbits \(\Orb_0=\{x\},\Orb_1,\ldots,\Orb_k\) for the
\(\RR^2\)-action generated by \(\bm{H}\). Since
\(\bm{H}^{-1}(0)\setminus\{x\}\) consists of regular
points, the orbits
\(\Orb_1,\ldots,\Orb_k\) are 2-dimensional submanifolds. By
Theorem \ref{thm:littlearnoldliouville}, there are three
possible topologies of orbit: \(\RR^2\), \(\RR\times S^1\) and
\(T^2\). The third type cannot occur because it would give a
connected component of \(\bm{H}^{-1}(0)\) not containing
\(x\), but we assume our fibres are connected. In particular,
any remaining orbits are noncompact.

Let \(\Orb_P\) be the orbit containing the Lagrangian plane
\(E(P)\) (in the focus-focus chart) and \(\Orb_Q\) be the
orbit containing \(E(Q)\). Note that it is possible that
\(\Orb_P=\Orb_Q\). Since the action of \(\RR^2\) on \(P\) (and
on \(Q\)) has stabiliser \(\ZZ\), these orbits are of the form
\(\RR\times S^1\). Moreover, these are the only orbits
containing \(\{x\}\) in their closure (such an orbit must
enter the focus-focus chart, where we can see that only
\(\Orb_P\) and \(\Orb_Q\) contain \(x\) in their
closure). There are two possibilities:
\begin{itemize}
\item \(\Orb_P\neq \Orb_Q\). In this case, the union \(\Orb_P\cup
\{x\}\cup \Orb_Q\) would be noncompact. It is impossible to
make \(\bm{H}^{-1}(0)\) compact by adding further noncompact
orbits, so because we assume \(H\) is proper, this
possibility does not occur.
\item \(\Orb_P=\Orb_Q\). In this case, the cylinder \(\Orb_P\) has
both its ends attached to the point \(x\), yielding a
(compact) pinched torus. Adding further noncompact orbits
contradicts compactness of \(\bm{H}^{-1}(0)\), so there are
no further orbits.\qedhere

\end{itemize}
\end{Proof}
The figure below shows a pinched torus fibre containing a
focus-focus critical point. The fixed point is shown with a dot,
the \(\phi^{H_2}_t\)-flowlines are the short loops going around
the fibre; the \(\phi^{H_1}_t\)-flowlines are the longer orbits
connecting the fixed point to itself.

\begin{center}
\begin{tikzpicture}
\filldraw[fill=lightgray,draw=none,opacity=0.5] (0,0) to[out=-135,in=0] (-1,-0.6) to[out=180,in=-90] (-2,0) to[out=90,in=180] (0,2) to[out=0,in=90] (2,0) to[out=-90,in=0] (1,-0.6) to[out=180,in=-45] (0,0);
\draw (0,0) to[out=-135,in=0] (-1,-0.6) to[out=180,in=-90] (-2,0) to[out=90,in=180] (0,2) to[out=0,in=90] (2,0) to[out=-90,in=0] (1,-0.6) to[out=180,in=-45] (0,0) [postaction={decorate, decoration={markings,mark=between positions 0.1 and 0.9 step 0.3 with {\arrow[line width=0.4mm]{>};}}}];
\draw (0,0) to[out=-150,in=0] (-0.75,-0.3) to[out=180,in=-90] (-1.5,0) to[out=90,in=180] (0,1.5) to[out=0,in=90] (1.5,0) to[out=-90,in=0] (0.75,-0.3) to[out=180,in=-30] (0,0) [postaction={decorate, decoration={markings,mark=between positions 0.2 and 0.8 step 0.3 with {\arrow[line width=0.4mm]{>};}}}];
\filldraw[fill=white] (0,0) to[out=-170,in=0] (-0.5,-0.1) to[out=180,in=-90] (-1,0) to[out=90,in=180] (0,1) to[out=0,in=90] (1,0) to[out=-90,in=0] (0.5,-0.1) to[out=180,in=-10] (0,0) [postaction={decorate, decoration={markings,mark=between positions 0.3 and 0.7 step 0.3 with {\arrow[line width=0.4mm]{>};}}}];
\draw (0,1) to[out=0,in=-90] (0.3,1.5) to[out=90,in=0] (0,2) [postaction={decorate, decoration={markings,mark=between positions 0.25 and 0.75 step 0.5 with {\arrow[line width=0.4mm]{>};}}}];
\draw[dotted] (0,2) to[out=180,in=90] (-0.3,1.5) to[out=-90,in=180] (0,1) [postaction={decorate, decoration={markings,mark=between positions 0.25 and 0.75 step 0.5 with {\arrow[line width=0.4mm]{>};}}}];
\draw (-1,0) to[out=90,in=0] (-1.5,0.3) to[out=180,in=90] (-2,0) [postaction={decorate, decoration={markings,mark=between positions 0.25 and 0.75 step 0.5 with {\arrow[line width=0.4mm]{>};}}}];
\draw[dotted] (-2,0) to[out=-90,in=180] (-1.5,-0.3) to[out=0,in=-90] (-1,0) [postaction={decorate, decoration={markings,mark=between positions 0.25 and 0.75 step 0.5 with {\arrow[line width=0.4mm]{>};}}}];
\draw (1,0) to[out=90,in=180] (1.5,0.3) to[out=0,in=90] (2,0) [postaction={decorate, decoration={markings,mark=between positions 0.25 and 0.75 step 0.5 with {\arrow[line width=0.4mm]{>};}}}];
\draw[dotted] (2,0) to[out=-90,in=0] (1.5,-0.3) to[out=180,in=90] (1,0) [postaction={decorate, decoration={markings,mark=between positions 0.25 and 0.75 step 0.5 with {\arrow[line width=0.4mm]{>};}}}];
\node at (0,0) {\(\bullet\)};
\end{tikzpicture}
\end{center}

\begin{Remark}
The same argument generalises to show that if
\(\bm{H}^{-1}(0)\) contains \(m>1\) focus-focus critical
points then it will form a cycle of Lagrangian spheres, each
intersecting the next transversely at a single focus-focus
point (or, if \(m=2\), two spheres intersecting transversely
at two points).

\end{Remark}
\section{Action coordinates}

Let \(\bm{H}\colon X\to\RR^2\) be an integrable Hamiltonian
system with a focus-focus critical point \(x\) over the origin
and no other critical points. Let \(E\colon U\to X\), \(e\colon
V\to\RR^2=\CC\) be a focus-focus chart centred at \(x\). Recall
that \(F=F_1+iF_2\colon U\to V\) denotes the model Hamiltonian
from Example \ref{exm:ff}. Let \(H_1 = F_1\circ E^{-1}\colon
E(U)\to \RR\) and \(H_2 = F_2\circ E^{-1}\colon E(U)\to
\RR\). By shrinking \(U\) and \(V\) if necessary, assume that
\(V=\{\bm{b}\in\RR^2\ :\ |\bm{b}|<\epsilon\}\) for some
\(\epsilon>0\); write \(B:=V\setminus\{0\}\) for the set of
regular values of \(H\). By Corollary \ref{cor:zaffine}, \(B\)
inherits an integral affine structure, coming from action
coordinates\index{action-angle coordinates!near focus-focus fibres|(} on the universal cover
\(\tilde{B}\). The next theorem identifies these action
coordinates.

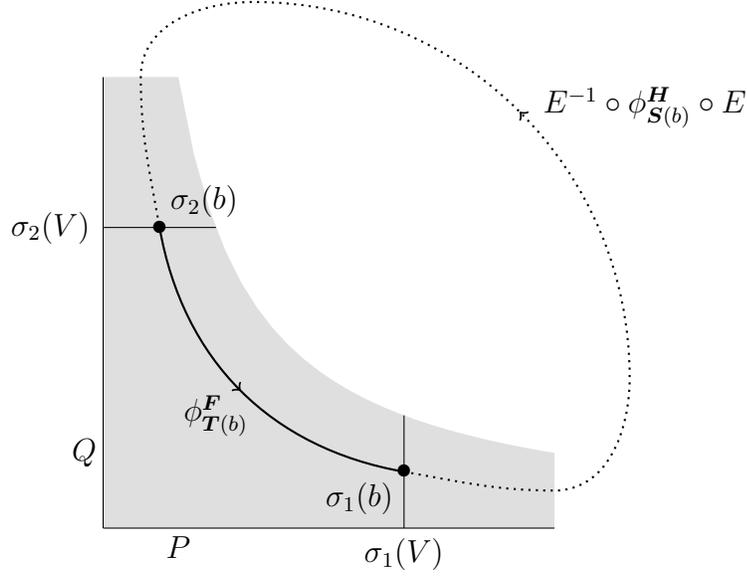
\begin{figure}[htb]
\begin{center}
\begin{tikzpicture}
\filldraw[fill=lightgray,opacity=0.5,draw=none,domain=1:6.0] (0,6) -- (1,6) plot (\x,{6/(\x)}) --(6,0) -- (0,0) -- (0,6);
\draw (0,0) -- (0,6);
\node at (1,-0.25) {\(P\)};
\node at (-0.25,1) {\(Q\)};
\draw (0,0) -- (6,0);
\draw (0,4) -- (1.5,4);
\node at (-0.7,4) {\(\sigma_2(V)\)};
\draw (4,0) -- (4,1.5);
\node at (4,-0.35) {\(\sigma_1(V)\)};
\draw[->-,dotted,thick] (4,0.75) to[out=-10,in=180] (6,0.5) to[out=0,in=-90] (7,2) to[out=90,in=0] (2,7) to[out=180,in=90] (0.5,6) to[out=-90,in=100] (0.75,4);
\node at (7.2,5.6) {\(E^{-1}\circ \phi^{\bm{H}}_{\bm{S}(b)}\circ E\)};
\node at (4,0.75) {\(\bullet\)};
\node at (0.75,4) {\(\bullet\)};
\node at (4,0.75) [below left] {\(\sigma_1(b)\)};
\node at (0.75,4) [above right] {\(\sigma_2(b)\)};
\draw[->-,thick] (0.75,4) to[out=-80,in=170] (4,0.75);
\node at (1.5,1.5) {\(\phi^{\bm{F}}_{\bm{T}(b)}\)};
\end{tikzpicture}
\end{center}
\caption{A schematic for the proof of Theorem \ref{thm:ngoc}, projected down to the \((|p|,|q|)\)-plane. The shaded region is the domain of the focus-focus chart, and we see the Lagrangian sections \(\sigma_1\) and \(\sigma_2\) intersecting the Lagrangian planes \(P\) and \(Q\) respectively. The flow \(\phi^{\bm{F}}_{\bm{T}}\) sends \(\sigma_2(b)\) to \(\sigma_1(b)\) for \(b\neq 0\). The flow \(E^{-1}\circ\phi^{\bm{H}}_{\bm{S}(b)}\circ E\) sends \(\sigma_1(b)\) to \(\sigma_2(b)\) (exiting the focus-focus chart).}\label{fig:ngoc_diagram}
\end{figure}

\begin{Theorem}[San V\~{u} Ng\d{o}c]\label{thm:ngoc}
The\index{Vu Ngoc@V\~{u} Ng\d{o}c!theorem} action map
\(\tilde{B}\to\RR^2\) has the form
\[\left(\frac{1}{2\pi}\left(S(b)+b_2\theta-b_1(\log
r-1)\right),b_2\right),\] where \(b=b_1+ib_2=re^{i\theta}\) is
the local coordinate on \(B\) and \(S(b)\) is a smooth
function.
\end{Theorem}
\begin{Proof}
The map \(\sigma_1\colon V\to\RR^4=\CC^2\),
\(\sigma_1(b)=(-1,b)\) is a Lagrangian section of \(F\) which
intersects \(Q\) at \(\sigma_1(0)\). Similarly
\(\sigma_2(b)=(-\bar{b},1)\) is a Lagrangian section which
intersects \(P\). See Figure \ref{fig:ngoc_diagram}.

For \(b\neq 0\), we can use the Hamiltonians \(F_1\) and
\(F_2\) {\em inside} our focus-focus chart to flow the point
\(\sigma_2(b)=(-\bar{b},1)\) until it hits
\(\sigma_1(b)=(-1,b)\) (see Figure \ref{fig:ngoc_diagram}). In
other words, we can find functions \(T_1(b)\) and \(T_2(b)\)
on \(V\setminus\{0\}\) with:
\begin{equation}\label{eq:sigma_2_to_sigma_1}
\phi^{F_2}_{T_2(b)}\phi^{F_2}_{T_1(b)}(-\bar{b},1) =
(-e^{T_1(b)+iT_2(b)}\bar{b}, e^{-T_1(b)+iT_2(b)}) =
(-1,b),\end{equation} namely \[T_1(b) = -\ln|b|,\quad T_2(b)
=\OP{arg}(b).\]

{\bf Claim.} After possibly shrinking \(V\), there exist
smooth functions \(S_1(b)\) and \(S_2(b)\) defined on \(V\)
such that
\begin{equation}\label{eq:sigma_1_to_sigma_2}
\phi^{H_2}_{S_2(b)}\phi^{H_1}_{S_1(b)}(E(\sigma_1(b))) =
E(\sigma_2(b))\end{equation} for all \(b\in V\) and such that
\begin{equation}\label{eq:S_1_S_2_closed} S_1 = \frac{\partial
S}{\partial b_1},\quad S_2 = \frac{\partial S}{\partial
b_2}\end{equation} for some smooth function \(S\) on \(V\).

Let us see how the claim implies the theorem. Since
\(H_k=F_k\circ E^{-1}\), we can combine Equations
\eqref{eq:sigma_2_to_sigma_1} and
\eqref{eq:sigma_1_to_sigma_2} to get
\[\phi^{\bm{H}}_{{\bm{S}}+{\bm{T}}}(E(\sigma_2(b)) =
E(\sigma_2(b))\] for all \(b\in V\), so that
\({\bm{S}}+{\bm{T}}:=(S_1(b)+T_1(b),S_2(b)+T_2(b))\) is in the
period lattice\index{period lattice!for focus-focus
singularity} (see Figure \ref{fig:ngoc_diagram}). The period
lattice is then spanned by these vectors and by \((0,2\pi)\)
(since \(H_2\) already has period \(2\pi\)). To find action
coordinates\index{action-angle coordinates!near focus-focus
fibres|)} \((G_1,G_2)\), it suffices to solve
\[\begin{pmatrix} \frac{\partial G_1}{\partial b_1} &
\frac{\partial G_1}{\partial b_2} \\ \frac{\partial
G_2}{\partial b_1} & \frac{\partial G_2}{\partial
b_2}\end{pmatrix} =
\begin{pmatrix}\frac{1}{2\pi}(S_1(b)-\ln|b|) &
\frac{1}{2\pi}(S_2(b)+arg(b)) \\ 0 & 1\end{pmatrix}.\] We can
take \(G_1(b) = \frac{1}{2\pi}(S+b_2\theta-b_1(\log
r-1))\) and \(G_2(b)=b_2\) where \(\theta=\OP{arg}(b)\)
and \(r=|b|\). This proves the theorem.

We now prove the claim. In the proof of Lemma
\ref{lma:pinched_torus}, we saw that the branch \(E(P)\) is
part of the same \(\RR^2\)-orbit as the branch
\(E(Q)\). Therefore, if we flow \(E(\sigma_1(0)\) for using
\(H_1\) for some duration \(s_1\), we will reach a point in
\(E(Q)\) at the same radius as \(\sigma_2(0)\). Further
flowing using \(H_2\) for some time \(s_2\), which preserves
the radius in the \(Q\)-plane, we can ensure that
\[\phi^{H_2}_{s_2}\phi^{H_1}_{s_1}(\sigma_1(0))=\sigma_2(0).\]
After possibly shrinking \(V\), we get local Liouville
coordinates near \(\sigma_2(0)\) using the Lagrangian section
\(\sigma_1\): \[\Psi(b,\bm{u}):=
\phi^{H_2}_{u_2}\phi^{H_1}_{u_1}(\sigma_1(b)).\] The domain of
\(\Psi\) is \(V\times I\) where \(I\) is a neighbourhood of
\((s_1,s_2)\) in \(\RR^2\). The preimage \(L :=
\Psi^{-1}(\sigma_2(V))\) is the Lagrangian submanifold of
\(V\times I\) given by \[L = \{(b,\bm{u})\in V\times I \,:\,
\phi^{H_2}_{u_2}\phi^{H_1}_{u_1}(1,-b) = (-\bar{b},1).\] We
pick the unique component of \(L\) containing
\((0,(s_1,s_2))\). This can be written as the graph of a
function \(b\mapsto (S_1(b),S_2(b))\). All that remains is to
solve the following exercise.

Exercise \ref{exr:lls}: The graph \(\{(b,(S_1(b),S_2(b)))\,:\,
b\in V\}\) is Lagrangian if and only if \(\partial
S_1/\partial b_2=\partial S_2/\partial b_1\), which holds if
and only if \(S_1 = \partial S/\partial b_1\) and \(S_2 =
\partial S/\partial b_2\) for some function \(S\).\qedhere

\end{Proof}
\begin{Remark}
In fact, any such \(S\) arises as we will show in the next
section. Moreover, V\~{u} Ng\d{o}c \cite{Ngoc}
showed\footnote{There is a subtlety here: the germ of \(S\)
can depend on the choice of focus-focus chart. This is a
finite ambiguity, and is discussed in {\cite[Section
4.3]{SepeNgoc}}: the actual V\~{u} Ng\d{o}c invariant is an
equivalence class of germs under an action of the Klein
4-group.} that the germ of \(S\) near the origin is unchanged
by any fibred symplectomorphism of the system, and that this
germ determines the (germ of the) system up to fibred
symplectomorphism in a neighbourhood of the nodal fibre. We
will write \((S)^\infty\) for the V\~{u} Ng\d{o}c
invariant\index{Vu Ngoc@V\~{u} Ng\d{o}c!invariant} of a focus-focus
critical point.

\end{Remark}
\begin{Remark}\label{rmk:base_node}
The action map has a well-defined limit point as \(r\to
0\). We call this limit point the {\em
base-node}\index{base-node} of the focus-focus critical point.

\end{Remark}
\section{Monodromy}

We\index{affine monodromy!for focus-focus singularity|(} briefly
recall the notion of affine monodromy introduced in Definition
\ref{dfn:affine_monodromy}. Let \(f\colon X\to B\) be a regular
Lagrangian fibration, let \(\tilde{B}\to B\) be the universal
cover of the base of a Lagrangian fibration and
\(\II\colon\tilde{B}\to\RR^n\) be the developing map for the
integral affine structure. Given an element \(g\in \pi_1(B)\),
we get a deck transformation \(\tilde{b}\mapsto \tilde{b}g\) of
the universal cover, and \(\II(\tilde{b}g)=\II(\tilde{b})M(g)\)
for some matrix \(M(g)\in SL(n,\ZZ)\).

\begin{Example}\label{exm:affine_monodromy_calc}
Let \(\bm{H}\) be an integrable Hamiltonian system with a
single focus-focus fibre, and let \(f\) be the restriction of
\(\bm{H}\) to the complement of the focus-focus fibre. In this
case, \(B=\RR^2\setminus\{0\}\) with polar coordinates
\(r,\theta\). The universal cover \(\tilde{B}\) is obtained by
treating \(\theta\) as a real-valued (instead of periodic
angular) coordinate. As a corollary of \ref{thm:ngoc}, we get
the developing map\index{developing map!for focus-focus
singularities} for the integral affine structure:
\[\II(r,\theta)=\left(\frac{1}{2\pi}\left(S(b) + b_2\theta -
b_1(\log r-1)\right),b_2\right),\] where
\((b_1,b_2)=(r\cos\theta,r\sin\theta)\). We have
\(\pi_1(B)=\ZZ\) and \(n\in\ZZ\) acts on \(\tilde{B}\) by
\((r,\theta)\mapsto (r,\theta+2\pi n)\).

\end{Example}
\begin{Lemma}[Exercise \ref{exr:affine_monodromy_calc}]\label{lma:affine_monodromy_calc}
The affine monodromy for \(n\in\pi_1(B)\) is
\(M(n)=\lmatrix 1 & 0 \\ n & 1\rmatrix \).

\end{Lemma}
This means that if you ``go around the loop'' in \(B\), the
action map changes by this shear matrix. Figures
\ref{fig:ngoc12} and \ref{fig:ngoc3} illustrate this by plotting
the image under \(\II\) of some different choices of fundamental
domain for the covering map \(\tilde{B}\to B\) (for the choice
\(S\equiv 0\)). We include the images under the action map of
contours of constant \(r\) (encircling the origin) and constant
\(\theta\) (pointing roughly radially outward).\index{developing
map!image}

\begin{minipage}{\linewidth}
\makebox[\linewidth]{
\includegraphics[scale=0.3,trim={12cm 2cm 9cm 2cm},clip]{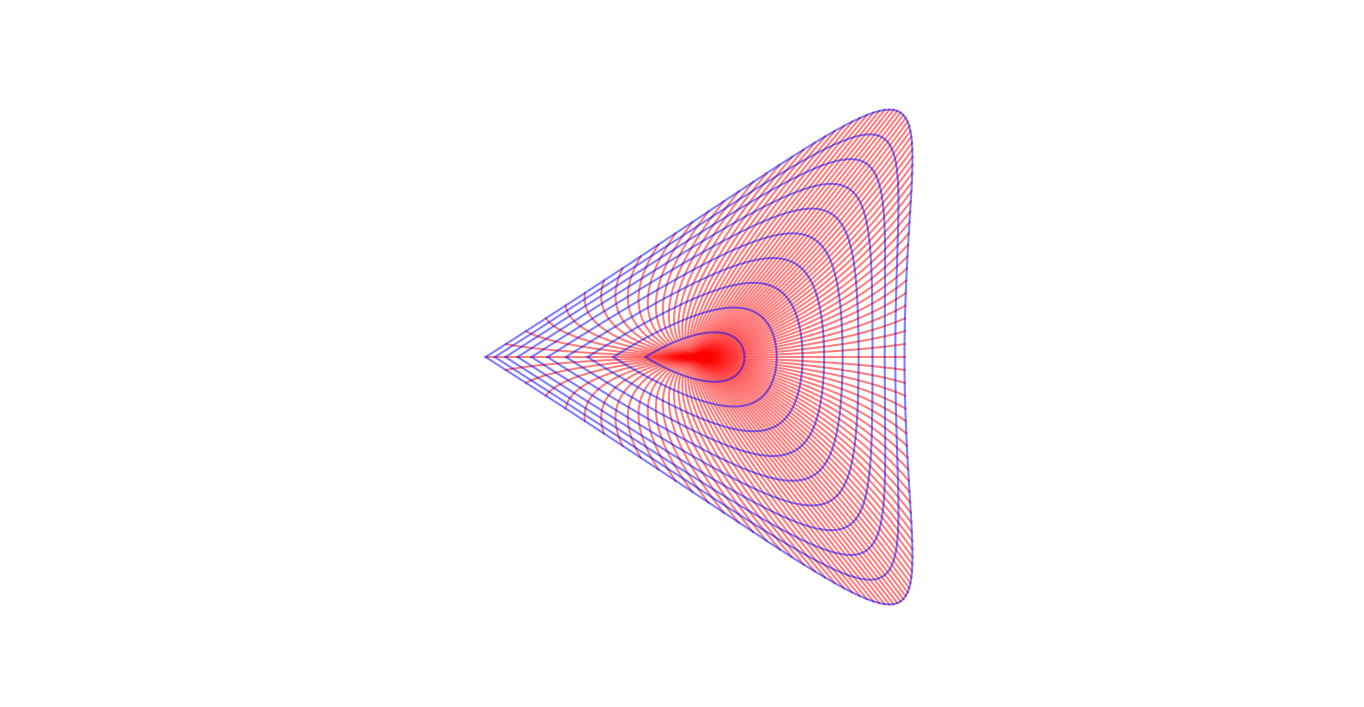}\includegraphics[scale=0.3,trim={8cm 2cm 6cm 2cm},clip]{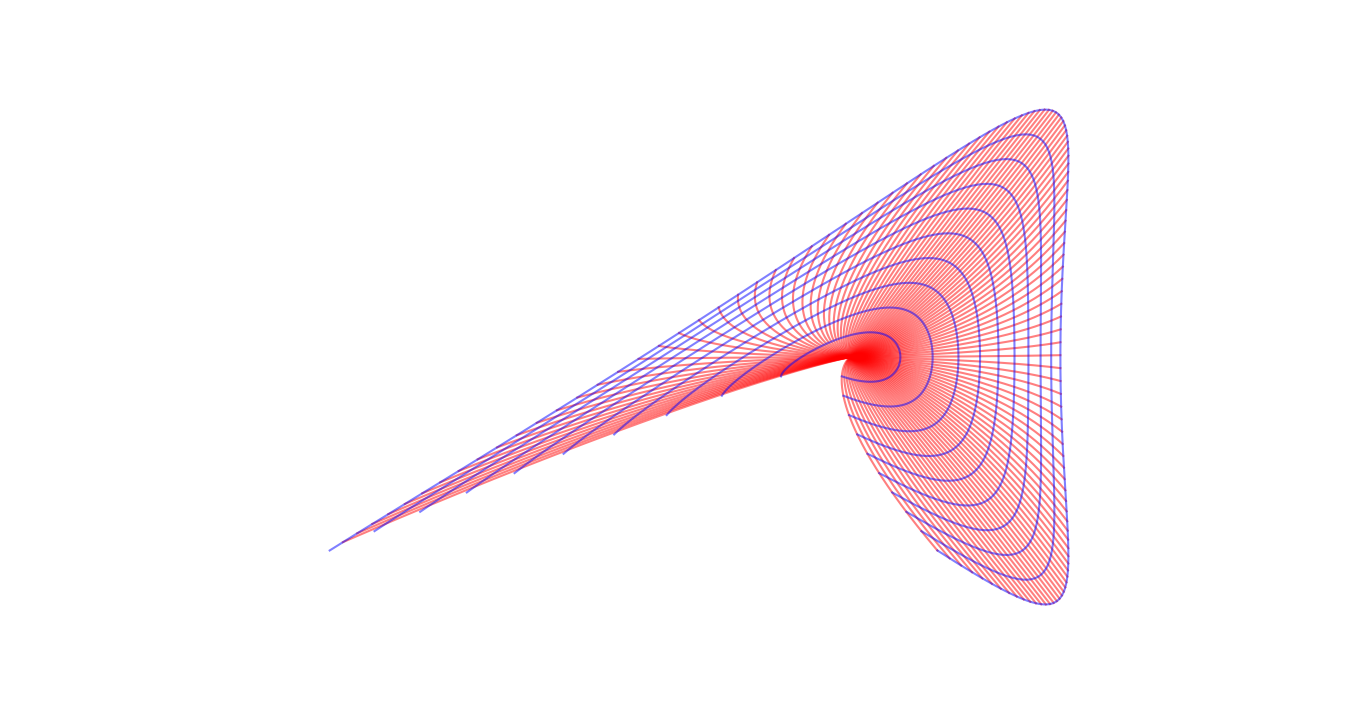}}
\captionof{figure}{Left: The image of the fundamental domain
\(\{\theta\in[-\pi,\pi)\}\). This ``closes up'' in the sense that \(\theta=-\pi\) and \(\theta=\pi\) map to the same line. This is because \(\theta=\pi\) is an eigenline of the monodromy matrix. Right: The image of the
fundamental domain \(\{\theta\in[-5\pi/7,9\pi/7)\}\). Although this
plot does not ``close up'', the
image of the radius \(\theta=-5\pi/7\) and the image of the radius
\(\theta=9\pi/7\) are related by the monodromy matrix.}\label{fig:ngoc12}
\end{minipage}

\begin{minipage}{\linewidth}
\makebox[\linewidth]{
\includegraphics[scale=0.4,trim={4cm 4cm 4cm 4cm},clip]{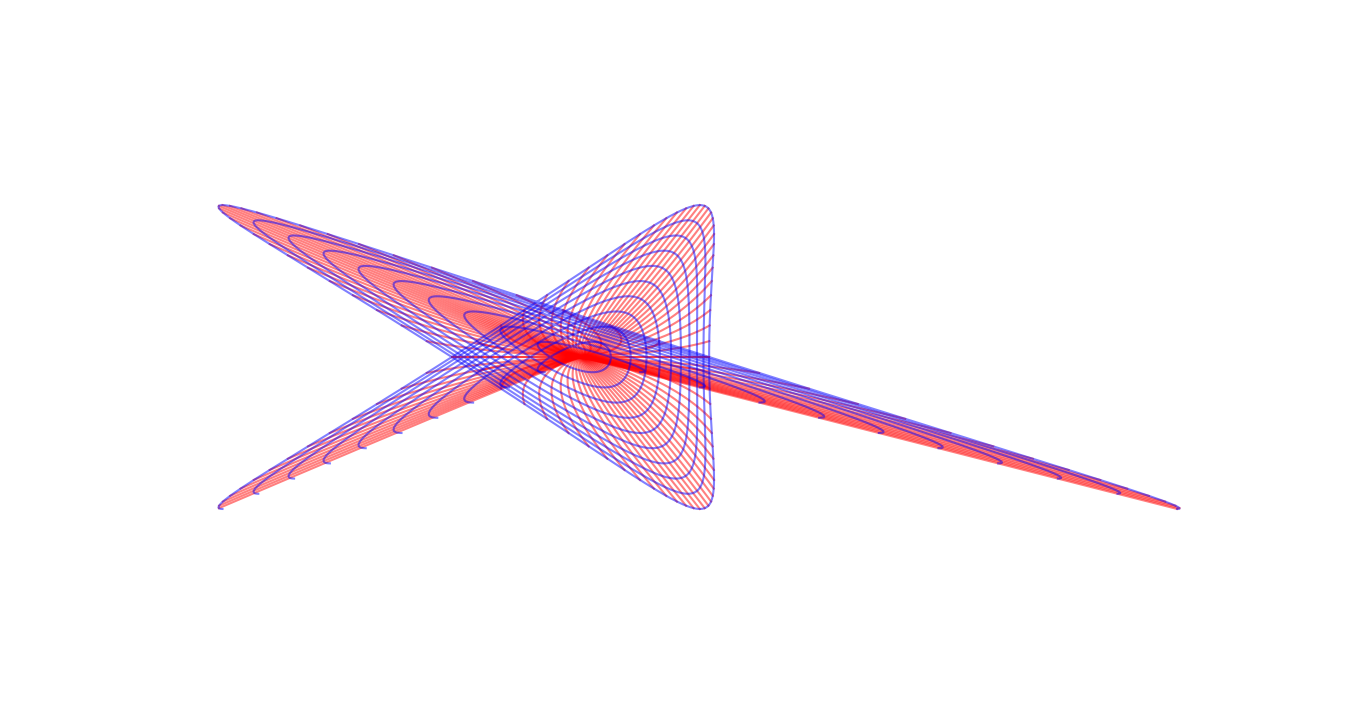}}
\captionof{figure}{In the third figure, we see the image of two fundamental domains
\(\{\theta\in[-5\pi/2,3\pi/2)\}\), related to one another by the
action of the monodromy matrix\protect\footnotemark.}\label{fig:ngoc3}
\end{minipage}

\footnotetext{Anyone who has compulsively traced out the spiral of a
raffia mat cannot fail to be moved by this image.}

\begin{Remark}\label{rmk:clockwise_vs_anticlockwise}
We think of \(B\) as being obtained from \(\tilde{B}\) as
follows. Fix a number \(\theta_0\in\RR\), take the fundamental
domain
\((0,\infty)\times[\theta_0,\theta_0+2\pi]\subset\tilde{B}\)
and identify points \((r,\theta_0)\) with
\((r,\theta_0+2\pi)\). In other words, we make a ``branch
cut''\index{branch cut|(} at \(\theta=\theta_0\). We pull back
the integral affine structure on \(\RR^2\) along the
developing map \(\II\) to get an integral affine structure on
\(\tilde{B}\), and this descends to \(B\) as in Corollary
\ref{cor:zaffine}. When we make the identification
\((r,\theta_0)\sim (r,\theta_0+2\pi)\), we need to specify how
to identify the integral affine structures. We use
\(\II(r,\theta_0+2\pi)=\II(r,\theta_0)M(1)\), that is
\(\II(r,\theta_0+2\pi)M(-1)=\II(r,\theta_0)\). In other words,
{\em when we cross the branch cut anticlockwise (direction of
increasing \(\theta\)), we apply the transformation \(M(-1)\)
to tangent vectors}.

\end{Remark}
Observe from Figures \ref{fig:ngoc12} and \ref{fig:ngoc3} that
if we use a branch cut \(\theta=0\) or \(\theta=\pi\) (parallel
to the eigenline\index{eigenline!branch
cut|(}\index{eigenray|see {eigenline}}\index{eigendirection|see
{eigenline}} of the affine monodromy) to cut out our fundamental
domain in \(\tilde{B}\) then the image of this fundamental
domain under \(\II\) ``closes up'', that is it is surjective onto
a punctured neighbourhood of the origin. If we use a different
branch cut then the image of the fundamental domain under \(\II\)
will miss out a segment of this punctured neighbourhood. For
this reason, we usually work with branch cuts\index{branch
cut|)} parallel to the eigenline\index{eigenline!branch cut|)}
of the affine monodromy.

Note that the action map from Theorem \ref{thm:ngoc} is only
unique up to post-composition by an integral affine
transformation. That is, by Lemma \ref{lma:ham_postcompose}, if
we post-compose \(\II\) by an integral affine
transformation\footnote{Here \(\bm{b}\) is a row vector and
\(A\) is acting on the right.}
\(\alpha(\bm{b})=\bm{b}A+\bm{C}\) (for some \(A\in GL(2,\ZZ)\)
and \(\bm{C}\in\RR^n\)) then we do not change the period
lattice, so we get an alternative set of action coordinates.

\begin{Lemma}[Exercise \ref{exr:monodromy_formula}]\label{lma:monodromy_formula}
If we use the action coordinates \(\II A\) for some \(A\in
SL(2,\ZZ)\) then the clockwise affine monodromy is given by
\(A^{-1}M(1)A\) and the line of
eigenvectors\index{eigenline!monodromy when crossing} points
in the \((1,0)A\)-direction. More precisely, if \((1,0)A =
(p,q)\) for some pair of coprime integers \(p,q\) and
\(\det(A) = 1\) then \(A^{-1}MA = \lmatrix 1 - pq & -q^2
\\ p^2 & 1 + pq\rmatrix \).

\end{Lemma}
\begin{Remark}
Remember that this matrix is acting {\em on the right}; if you
want to think of your action coordinates as column vectors,
you need to take the transpose matrix. Remember also that this
is the {\em clockwise} monodromy: you apply its inverse to
tangent vectors when you cross the branch cut
anticlockwise.\index{affine monodromy!for focus-focus
singularity|)}

\end{Remark}
\section{Visible Lagrangians}

The following visible Lagrangian
disc\index{Lagrangian!submanifold!visible|(} will play an
important role in our future analysis of focus-focus systems.

\begin{Lemma}\label{lma:visibleff}
Let \(\bm{H}\colon X\to\RR^2\) be an integrable Hamiltonian
system with a focus-focus critical point at \(x\in\ X\), let
\(B\) be the set of regular values and \(\tilde{B}\) its
universal cover, and let \(\II\colon\tilde{B}\to\RR^2\) be the
developing map for the integral affine structure on \(B\)
coming from action coordinates. Let \(b\in\RR^2\) be the
base-node associated to the focus-focus critical point at
\(x\). Suppose that \(\ell\) is a straight ray in \(\RR^2\)
emanating from \(b\) pointing in an
eigendirection\index{eigenline!visible Lagrangian over} for
the affine monodromy around the critical value. Then there is
a visible Lagrangian disc living over \(\ell\).
\end{Lemma}
\begin{Proof}
In the focus-focus chart we can simply use the Lagrangian disc
\(q=p\), which satisfies \(F(p,p)=-\bar{p}p\), so this lives
over the negative \(b_1\)-axis (\(b_2=0\)). By Theorem
\ref{thm:ngoc}, the image of this under \(\II\) is still the
negative \(b_1\)-axis, which is an eigenray of the affine
monodromy. \qedhere

\end{Proof}
\begin{Definition}
By analogy with a similar (but slightly different\footnote{In
Picard-Lefschetz theory, we have a {\em holomorphic} fibration
instead of a Lagrangian fibration, but the thimble is still a
Lagrangian disc.}) situation in Picard-Lefschetz theory, this
visible Lagrangian disc is called the {\em vanishing
thimble}\index{vanishing!thimble} for the focus-focus critical
point, and its intersection with any fibre over the ray
\(\ell\) is a loop in the fibre called the {\em vanishing
cycle}.\index{Lagrangian!submanifold!visible|)}\index{vanishing!cycle}

\end{Definition}
\section{Model neighbourhoods}

We now present a construction due to
V\~{u} Ng\d{o}c\index{Vu Ngoc@V\~{u} Ng\d{o}c!model|(} which, given a function
\(S\colon\RR^2\to\RR\), produces a Hamiltonian system
\(\bm{H}_S\colon X_S\to\RR^2\) with a focus-focus critical point
whose V\~{u} Ng\d{o}c invariant\index{Vu Ngoc@V\~{u} Ng\d{o}c!invariant|(} is
\((S)^\infty\). We will write \(S_i=\frac{\partial S}{\partial
b_i}\), \(i=1,2\).

Take the subset \(X:=\{(p,q)\in\RR^4\ :\ |\bar{p}q|<\epsilon\}\)
equipped with the Hamiltonian system \(F\) from Example
\ref{exm:ff}. We will construct two Liouville coordinate systems
on different regions of this space.

Recall the Lagrangian sections \(\sigma_1(b)=(-1,b)\) and
\(\sigma_2(b)=(-\bar{b},1)\). We construct a third Lagrangian
section \(\sigma_3(b) = (-e^{S_1(b)+iS_2(b)},
e^{-S_1(b)+iS_2(b)}b)=\phi^{\bm{H}}_{\bm{S}}(\sigma_1(b))\). We
can use these Lagrangian sections to construct Liouville
coordinates
\[\Psi_2(b,\bm{t})=\phi^{\bm{F}}_{\bm{t}}(\sigma_2(b))\] and
\[\Psi_3(b,\bm{t})=\phi^{\bm{F}}_{\bm{t}}(\sigma_3(b))\] with
\(t_1\in [0,\delta)\) and \(t_2\in[0,2\pi)\).

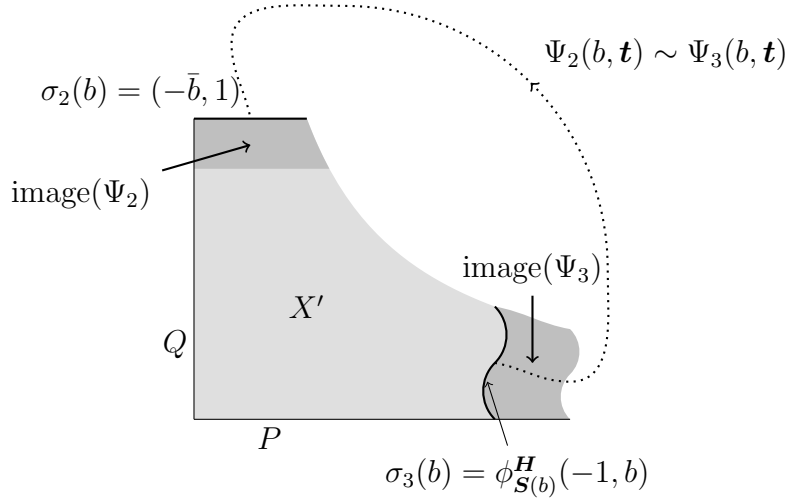
\begin{figure}[htb]
\begin{center}
\begin{tikzpicture}
\filldraw[fill=lightgray,opacity=0.5,draw=none,domain=1.5:4.0] (0,4) -- (1.5,4) plot (\x,{6/(\x)}) to[out=-45,in=45] (4,0.75) to[out=-135,in=135] (4,0) -- (0,0) -- (0,4);
\filldraw[fill=lightgray,draw=none,domain=1.5:1.8] (0,4) -- (1.5,4) plot (\x,{6/(\x)}) -- (0,6/1.8) -- (0,4);
\filldraw[fill=lightgray,draw=none,domain=1.5:1.8] (4,1.5) to[out=-45,in=45] (4,0.75) to[out=-135,in=135] (4,0) -- (5,0) to[out=135,in=-135] (5,0.6) to[out=45,in=-45] (5,1.2) to[out=170,in=-10] (4,1.5);
\node (psi_2) at (-1.5,3) {\(\OP{image}(\Psi_2)\)};
\node (psi_3) at (4.5,2) {\(\OP{image}(\Psi_3)\)};
\draw[->,thick] (psi_2) -- (0.75,3.65);
\draw[->,thick] (psi_3) -- (4.5,0.75);
\draw (0,0) -- (0,4);
\node at (1,-0.25) {\(P\)};
\node at (-0.25,1) {\(Q\)};
\draw (0,0) -- (5,0);
\draw[thick] (0,4) -- (1.5,4);
\node at (-0.7,4) [above] {\(\sigma_2(b) = (-\bar{b},1)\)};
\draw[thick] (4,1.5) to[out=-45,in=45] (4,0.75) to[out=-135,in=135] (4,0);
\node (sigma_3) at (4.3,-0.35) [below] {\(\sigma_3(b) = \phi^{\bm{H}}_{\bm{S}(b)}(-1,b)\)};
\draw[->] (sigma_3) -- (3.9,0.5);
\draw[->-,dotted,thick] (4,0.75) to[out=-10,in=180] (5,0.5) to[out=0,in=-90] (5.5,2) to[out=90,in=0] (2,5.5) to[out=180,in=90] (0.5,5) to[out=-90,in=100] (0.75,4);
\node at (4.5,4.5) [above right] {\(\Psi_2(b,\bm{t})\sim\Psi_3(b,\bm{t})\)};
\node at (1.5,1.5) {\(X'\)};
\end{tikzpicture}
\end{center}
\caption{Construction of a V\~{u} Ng\d{o}c model with invariant \(S\). The subset \(X'\) is the entire shaded region; the images of the Liouville coordinates \(\Psi_2(b,\bm{t})=\phi^{\bm{F}}_{\bm{t}}(\sigma_2(b))\) and \(\Psi_3(b,\bm{t})=\phi^{\bm F}_{\bm t}(\sigma_3(b))\) are shaded darker. The quotient \(X_S\) identifies these two darkly shaded regions.}\label{fig:ngoc_model}
\end{figure}

Let \(X'=\{(p,q)\in\RR^4\ :\ |\bar{p}q|<\epsilon,\ |q|\leq
1,\ |p|\leq e^{S_1(-\bar{p}q)}\}\) and let \(X_S\) be the
quotient \(X_S:=X'/\sim\), where \(\sim\) identifies
\(\Psi_2(b,t)\sim\Psi_3(b,t)\) (see Figure
\ref{fig:ngoc_model}). Since the domains of \(\Psi_2\) and
\(\Psi_3\) are identical and since \(\Psi_2\) and \(\Psi_3\) are
symplectomorphisms, the symplectic form on \(X\) descends to
this quotient. By construction, the map \(\bm{H}\colon
X\to\RR^2\), \(\bm{H}(p,q)=-\bar{p}q\) descends to the quotient
and produces the Hamiltonian system \(\bm{H}_S\) we want. Also
by construction, the V\~{u} Ng\d{o}c
invariant\index{Vu Ngoc@V\~{u} Ng\d{o}c!invariant|)} is
\((S)^{\infty}\).

\section{Symington's theorem on V\~{u} Ng\d{o}c models}

We now present an argument of Symington {\cite[Lemma
3.6]{Symington2}} which tells us that, although the V\~{u} Ng\d{o}c
models \(\bm{H}_{S_0}\colon
X_{S_0}\to\RR^2\), \(\bm{H}_{S_1}\colon X_{S_1}\to\RR^2\) with
\((S_0)^\infty\neq (S_1)^\infty\) are not symplectomorphic via a
fibred symplectomorphism, there is nonetheless a
symplectomorphism \(X_{S_0}\to X_{S_1}\) which is fibred outside
a compact set. Because the fibred symplectomorphism type of the
system depends only on the germ of \(S\) near the origin, we may
assume that \(S_0\) and \(S_1\) coincide outside a small
neighbourhood of the origin

\begin{Theorem}[Symington]\label{thm:symington}
Let\index{Symington's theorem!on Vu Ngoc models@on V\~{u} Ng\d{o}c
models} \(S_0\colon\RR^2\to\RR\) and \(S_1\colon\RR^2\to\RR\)
be smooth functions which coincide on the complement of a
small disc \(D\) centred at the origin and let
\(\bm{H}_{S_0}\colon X_{S_0}\to\RR^2\) and
\(\bm{H}_{S_1}\colon X_{S_1}\to\RR^2\) be the corresponding
V\~{u} Ng\d{o}c models\index{Vu Ngoc@V\~{u} Ng\d{o}c!model|)}. Then there
is a symplectomorphism \(\varphi\colon X_{S_0}\to X_{S_1}\)
which restricts to a fibred symplectomorphism
\(\bm{H}_{S_0}^{-1}(\RR^2\setminus D)\to
\bm{H}_{S_1}^{-1}(\RR^2\setminus D)\).
\end{Theorem}
\begin{Proof}
Pick a family \(S_t\) interpolating between \(S_0\) and
\(S_1\) such that \(S_t|_{\RR^2\setminus
D}=S_0|_{\RR^2\setminus D}\). Consider the family of
symplectic manifolds \(X_t:=X_{S_t}\); since the construction
depends only on \(S_t\), which is independent of \(t\) on the
complement of \(D\), the subsets \(U_t :=
\bm{H}_{S_t}^{-1}(\RR^2\setminus D)\) are
fibred-symplectomorphic (via the identity map).

We extend this identification to a isotopy of diffeomorphisms
\(\varphi_t\colon X_0\to X_t\) such that
\(\varphi_t|_{U_0}=\OP{id}\colon U_0\to U_t\). For example, we
could pick a connection on the family \(X_t\) which is trivial
on \(\bigcup_{t\in[0,1]}U_t=U_0\times[0,1]\) and take
\(\varphi_t\) to be the parallel transport of fibres.

Consider the family of symplectic forms \(\omega_t =
\varphi_t^*\omega_{S_t}\) on \(X_0\). These satisfy
\(\frac{d\omega_t}{dt} = 0\) on \(U_0\). The 2-form
\(\frac{d\omega_t}{dt}\) therefore determines a class in
\(H^2_{dR}(X_0,U_0)\). We will show that this class vanishes
for all \(t\); this will allow us to pick a family of 1-forms
\(\beta_t\) such that \(d\beta_t = d\omega_t/dt\) and
\(\beta_t|_{U_0}=0\). By Moser's trick\index{Moser argument}
(see Appendix \ref{ch:moser}), we then get diffeomorphisms
\(\phi_t\colon X_{S_0}\to X_{S_0}\), equal to the identity
outside \(U_0\), such that
\(\phi_t^*\varphi_t^*\omega_{S_t}=\omega_{S_0}\). The
symplectomorphism we want is
\(\varphi:=\varphi_1\circ\phi_1\colon X_{S_0}\to X_{S_1}\).

It remains to show that \(\frac{d[\omega_t]}{dt}=0\in
H^2_{dR}(X_0,U_0)\). Let \(V=X_0\setminus U_0\). We have
\(H^2_{dR}(X_0,U_0)=H^2_{dR}(V,\partial V)\) by excision, and
\(H^2_{dR}(V,\partial V) = H_2(V)\) by Poincar\'{e}-Lefschetz
duality. Since \(V\) deformation-retracts onto the nodal
fibre, we have \(H^2_{dR}(V)=\RR\). The cohomology class of a
closed 2-form on \(X_0\) which vanishes on \(U_0\) can
therefore be detected by its integral over a disc in \(V\)
with boundary on \(\partial V\) which intersects the nodal
fibre once transversely, e.g.\ a section of
\(\bm{H}_{S_t}\). The construction of \(X_{S_t}\) furnishes it
with an \(\omega_{S_t}\)-Lagrangian section (i.e.\
\((-\bar{b},1)\)); let us write \(\sigma_t\) for this section
viewed via \(\varphi_t\) as a submanifold of \(X_0\). Since
\(\varphi_t\) is the identity outside \(U_0\), we have
\(\sigma_t\cap U_0=\sigma_0\cap U_0\). Since \(\sigma_t\) is
\(\omega_t\)-Lagrangian, we have \(0 =
\int_{\sigma_t}\omega_t\). Fix \(T\) and let
\(\Sigma_T(\bm{b},t):=\sigma_t(\bm{b})\) be the isotopy of
sections restricted to \(t\in[0,T]\). Since \(d\omega_T=0\),
Stokes's theorem\footnote{There should be further boundary
terms corresponding to the boundary of the section, but since
the sections are all fixed over \(U_0\) these contributions
vanish.} tells us that \[0=\int\Sigma^*d\omega_T =
\int_{\sigma_T}\omega_T - \int_{\sigma_0}\omega_T.\] Since
\(\sigma_T\) is \(\omega_T\)-Lagrangian, we get
\(\int_{\sigma_0}\omega_T=0\). Therefore
\(\int_{\sigma_0}\frac{d\omega_t}{dt} =
\frac{d}{dt}\int_{\sigma_0}\omega_t = 0\), so
\(\frac{d[\omega_t]}{dt} = 0\in H^2(X_0, U_0)\). \qedhere

\end{Proof}
\section{Solutions to inline exercises}

\begin{Exercise}[Lemma \ref{lma:ff_flow}]\label{exr:ff_flow}
Verify that the Hamiltonians \(F_1=-p_1q_1-p_2q_2\) and
\(F_2=p_2q_1-p_1q_2\) Poisson-commute, that \(F_1\) generates
the \(\RR\)-action \((p,q)\mapsto(e^tp,e^{-t}q)\) and that
\(F_2\) generates the circle action
\((p,q)\mapsto(e^{it}p,e^{it}q)\).
\end{Exercise}
\begin{Solution}
We have \(-dF_1 = p_1\,dq_1+p_2\,dq_2+q_1\,dp_1+q_2\,dp_2\) which
equals \(\iota_{(p_1,-q_1,p_2,-q_2)}(dp_1\wedge
dq_1+dp_2\wedge dq_2)\). Thus \(V_{F_1}=(p_1,-q_1,p_2,-q_2)\)
and the flow satisfies \(\dot{p}=p\), \(\dot{q}=-q\), which
means \(p(t)=e^tp(0)\) and \(q(t)=e^{t}q(0)\). Similarly, we
find \(V_{F_2} = (-p_2,-q_2,p_1,q_1)\), whose flow satisfies
\(\dot{p} = ip\) and \(\dot{q} = iq\) (recall that
\(p=p_1+ip_2\) and \(q=q_1+iq_2\)) and the flow is therefore
\(p(t)=e^{it}p(0)\) and \(q(t)=e^{it}q(0)\). To see that
\(\{F_1,F_2\}=0\), we compute
\[\{F_1,F_2\}=\omega(V_{F_1},V_{F_2})=-p_1q_2-q_1p_2+p_2q_1-q_2p_1
= 0.\qedhere\]

\end{Solution}
\begin{Exercise}[From the proof of Theorem \ref{thm:ngoc}]\label{exr:lls}
Let \(V\) be a disc in \(\RR^2\). The graph
\(\{(b,(S_1(b),S_2(b)))\,:\, b\in V\}\) is Lagrangian if and
only if \(\partial S_1/\partial b_2=\partial S_2/\partial
b_1\), which holds if and only if \(S_1 = \partial S/\partial
b_1\) and \(S_2 = \partial S/\partial b_2\) for some function
\(S\).
\end{Exercise}
\begin{Solution}
The tangent space to the graph is spanned by the vectors
\begin{gather*}
(1,0,\partial S_1/\partial b_1,\partial S_2,\partial b_2), \\
(0,1,\partial S_1/\partial b_2,\partial S_2/\partial b_2),
\end{gather*}
on which the symplectic form evaluates to \[\partial
S_2/\partial b_1 - \partial S_1/\partial b_2.\] The graph is
Lagrangian if and only if this quantity vanishes. This is
equivalent to the condition that the 1-form \[S_1db_1 +
S_2db_2\] is closed. Since the disc \(V\) has zero de Rham
cohomology in degree 1, this 1-form is closed if and only if
it is exact, that is, if and only if there exists a function
\(S\) with \(\partial S/\partial b_k=S_k\) for
\(k=1,2\).\qedhere

\end{Solution}
\begin{Exercise}[Lemma \ref{lma:affine_monodromy_calc}]\label{exr:affine_monodromy_calc}
The affine monodromy for \(n\in\pi_1(B)\) in Example
\ref{exm:affine_monodromy_calc} is \(M(n)=\lmatrix 1 & 0 \\ n
& 1\rmatrix\).
\end{Exercise}
\begin{Solution}\belowdisplayskip=-12pt Since
\[\II(r,\theta)=\left(\frac{1}{2\pi}\left(S(b) + b_2\theta -
b_1(\log r-1)\right),\ b_2\right),\] we have
\begin{align*}
\II(r,\theta+2\pi n)&=\left(\frac{1}{2\pi}\left(S(b) +
b_2(\theta+2\pi n) - b_1(\log r-1)\right),\ b_2\right),\\
&=\left(\frac{1}{2\pi}\left(S(b) +
b_2\theta - b_1(\log r-1)\right)+nb_2,\ b_2\right)\\
&=\II(r,\theta)\lmatrix 1 & 0 \\ n & 1\rmatrix .
\end{align*}\qedhere
\end{Solution}

\begin{Exercise}[Lemma \ref{lma:monodromy_formula}]\label{exr:monodromy_formula}
Let \(\II\) be the developing map from Theorem
\ref{thm:ngoc}. If we use instead the developing map \(\II A\)
for some \(A\in SL(2,\ZZ)\) then the clockwise affine
monodromy is given by \(A^{-1}M(1)A\) and the line of
eigenvectors points in the \((1,0)A\)-direction. More
precisely, if \((1,0)A = (p,q)\) for some pair of coprime
integers \(p,q\) and \(\det(A)=\pm 1\) then \(A^{-1}M(1)A =
\lmatrix 1\mp pq & \mp q^2 \\ \pm p^2 & 1\pm pq\rmatrix \).
\end{Exercise}
\begin{Solution}
More generally, if \(g\in\pi_1(B)\) then we have
\(\II(\tilde{b}g)A = \II(\tilde{b})M(g)A =
\II(\tilde{b})A(A^{-1}M(g)A)\), so the affine monodromy
associated to the developing map \(\II A\) is
\(A^{-1}M(g)A\). Since \(v=(1,0)\) is an eigenvector of
\(M(g)\), so that \(vM(g)=v\), then \(vAA^{-1}M(g)A=vM(g)A =
vA\), so \(vA\) is an eigenvector of \(A^{-1}M(g)A\).

If \((1,0)A=(p,q)\) then \(A=\lmatrix p & q \\ k &
\ell\rmatrix \) for some \(k,\ell\in\ZZ\) with \(p\ell-kq=\pm
1\). Therefore\index{critical
point!focus-focus|)}
\[A^{-1}\lmatrix 1 & 0 \\ 1 & 1\rmatrix A=\pm\lmatrix \ell &
-q \\ -k & p\rmatrix \lmatrix p & q \\ p+k & q+\ell\rmatrix
=\lmatrix 1\mp pq & \mp q^2 \\ \pm p^2 & 1\pm pq\rmatrix
.\qedhere\]

\end{Solution}
\chapter{Examples of focus-focus systems}
\label{ch:focusfocus_examples}
\thispagestyle{cup}

We are now ready to introduce the notion of an {\em almost toric
manifold}: a symplectic 4-manifold with a Lagrangian torus
fibration whose critical points can be both toric and
focus-focus type. Before developing the general theory in
Chapter \ref{ch:almost_toric_manifolds}, we explore some
examples.

\section{The Auroux system}

Like many people, I first learned of the following example from
the wonderful expository article \cite{Auroux} on mirror
symmetry for Fano varieties by Denis Auroux, where it serves to
illustrate the wall-crossing phenomenon for discs.

\begin{figure}[htb]
\begin{center}
\begin{tikzpicture}
\filldraw[fill=lightgray,opacity=0.5,draw=none] (0,0) -- (0,2) -- (2,2) -- (2,0) -- cycle;
\draw[thick,black] (0,0) -- (0,2);
\draw[dotted,thick] (2,1) -- (1,1);
\node at (1,1) {\(\times\)};
\end{tikzpicture}
\end{center}
\caption{The fundamental action domain from Lemma \ref{lma:zaffexample}.}\label{fig:fadzaffexample}
\end{figure}
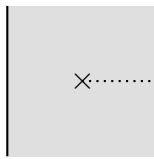

\begin{Example}[Auroux system]\label{exm:auroux}
Fix\index{Auroux system|(} a real number \(c>0\). Consider the
Hamiltonians \(\bm{H}=(H_1,H_2)\colon\CC^2\to\RR^2\) defined
by \(H_1(z_1,z_2)=|z_1z_2-c|^2\) and
\(H_2(z_1,z_2)=\frac{1}{2}\left(|z_1|^2-|z_2|^2\right)\). The
flow of \(H_2\) is
\(\phi^{H_2}_t(z_1,z_2)=(e^{it}z_1,e^{-it}z_2)\). This shows
that \(\{H_1,H_2\}=0\), because \(H_1\) is constant along the
flow of \(H_2\) (see Lemma \ref{lma:heisenberg}). The flow of
\(H_1\) is harder to compute. We can nonetheless understand
the orbits of this system geometrically.

Consider the holomorphic map \(\pi\colon\CC^2\to\CC\),
\(\pi(z_1,z_2)=z_1z_2\). This is a conic
fibration\index{conic!fibration|(}: the fibres \(\pi^{-1}(p)\)
are smooth conics except \(\pi^{-1}(0)\) which is a singular
conic (union of the \(z_1\)- and \(z_2\)-axes).

\begin{center}
\begin{tikzpicture}
\filldraw[fill=lightgray,opacity=0.5,draw=none] (0,-0.5) -- (6.5,-0.5) -- (7.5,2.5) -- (1,2.5) -- cycle;
\draw (4,1) circle [x radius=1,y radius=0.5];
\draw (4,1) circle [x radius=1.5,y radius=0.75];
\draw (4,1) circle [x radius=2,y radius=1];
\draw (4,1) circle [x radius=2.5,y radius=1.25];
\node at (6.2,0) {\(C_r\)};
\node at (6.5,-0.5) [below right] {\(\CC\)};
\filldraw[fill=white,draw=none] (1.5,1) -- (2,1.5) -- (2.5,1) -- cycle;
\filldraw[fill=white,draw=none] (1.5,2) -- (2,1.5) -- (2.5,2) -- cycle;
\filldraw[fill=white] (2,2) circle [x radius=0.5,y radius=0.25];
\filldraw[fill=white,dashed] (1.5,1) arc [x radius=0.5,y radius=0.25,start angle=180,end angle=0];
\filldraw[fill=white] (1.5,1) arc [x radius=0.5,y radius=0.25,start angle=180,end angle=360];
\draw (1.5,1) -- (2.5,2);
\draw (1.5,2) -- (2.5,1);
\node at (2,1) {\(0\)};
\filldraw[fill=white,draw=none] (3.5,1) to[out=45,in=-45] (3.5,2) -- (4.5,2) to[out=-135,in=135] (4.5,1) -- (3.5,1);
\filldraw[fill=white] (4,2) circle [x radius=0.5,y radius=0.25];
\filldraw[fill=white,dashed] (3.5,1) arc [x radius=0.5,y radius=0.25,start angle=180,end angle=0];
\filldraw[fill=white] (3.5,1) arc [x radius=0.5,y radius=0.25,start angle=180,end angle=360];
\draw (3.5,1) to[out=45,in=-45] (3.5,2);
\draw (4.5,1) to[out=135,in=-135] (4.5,2);
\node at (4,1) {\(c\)};
\filldraw[fill=white,draw=none] (5.5,1) to[out=45,in=-45] (5.5,2) -- (6.5,2) to[out=-135,in=135] (6.5,1) -- (5.5,1);
\filldraw[fill=white] (6,2) circle [x radius=0.5,y radius=0.25];
\filldraw[fill=white,dashed] (5.5,1) arc [x radius=0.5,y radius=0.25,start angle=180,end angle=0];
\filldraw[fill=white] (5.5,1) arc [x radius=0.5,y radius=0.25,start angle=180,end angle=360];
\draw (5.5,1) to[out=45,in=-45] (5.5,2);
\draw (6.5,1) to[out=135,in=-135] (6.5,2);
\node at (6,1) {\(\bullet\)};
\end{tikzpicture}
\end{center}

The Hamiltonian \(H_1\) measures the squared distance in
\(\CC\) from \(z_1z_2\) to some fixed point \(c\). The level
set \(H_1^{-1}(r^2)\) is therefore the union of all conics
living over a circle \(C_r\) of radius \(r\) centred at \(c\)
(the concentric circles in the base of the figure). The
restriction of \(H_2\) to each conic can be visualised as a
``height function'' whose level sets are circles as shown
below. The level set \(\bm{H}^{-1}(b_1,b_2)\) is therefore the
union of all circles of height \(b_2\) in conics living over
the circle \(C_{\sqrt{b_1}}\). These level sets are clearly
tori, except for the level set
\(\bm{H}^{-1}\left(c^2,0\right)\), which is a pinched
torus\index{pinched torus}.

\begin{center}
\begin{tikzpicture}
\filldraw[fill=lightgray,opacity=0.5,draw=none] (0,-0.5) -- (6.5,-0.5) -- (7.5,2.5) -- (1,2.5) -- cycle;
\draw (4,1) circle [x radius=1,y radius=0.5];
\draw (4,1) circle [x radius=1.5,y radius=0.75];
\draw (4,1) circle [x radius=2,y radius=1];
\draw (4,1) circle [x radius=2.5,y radius=1.25];
\node at (6.4,-0.1) {\(C_{\sqrt{b_1}}\)};
\draw[->] (6,0) to[out=180,in=-90] (5.6,0.3);
\node at (6.7,3.3) {\(\bm{H}^{-1}(b_1,b_2)\)};
\draw[->] (5.6,3.2) to[out=180,in=90] (5,2.8);
\draw[->,thick] (7,1.5) -- (7,2) node [right] {\(b_2\)};
\draw[dotted,thick] (6.29,1.5) -- (7,1.5);
\filldraw[fill=white,draw=none] (1.5,1) -- (2,1.5) -- (2.5,1) -- cycle;
\filldraw[fill=white,draw=none] (1.5,2) -- (2,1.5) -- (2.5,2) -- cycle;
\filldraw[fill=white,draw=black,thick] (2,2) circle [x radius=0.5,y radius=0.25];
\filldraw[fill=white,dashed] (1.5,1) arc [x radius=0.5,y radius=0.25,start angle=180,end angle=0];
\filldraw[fill=white] (1.5,1) arc [x radius=0.5,y radius=0.25,start angle=180,end angle=360];
\draw (1.5,1) -- (2.5,2);
\draw (1.5,2) -- (2.5,1);
\draw[dashed] (1.75,1.25) arc [x radius=0.25,y radius=0.125,start angle=180,end angle=0];
\draw (1.75,1.25) arc [x radius=0.25,y radius=0.125,start angle=180,end angle=360];
\draw[dashed] (1.75,1.75) arc [x radius=0.25,y radius=0.125,start angle=180,end angle=0];
\draw (1.75,1.75) arc [x radius=0.25,y radius=0.125,start angle=180,end angle=360];
\node at (2,1.5) {\(\bullet\)};
\filldraw[fill=white,draw=none] (3.5,1) to[out=45,in=-45] (3.5,2) -- (4.5,2) to[out=-135,in=135] (4.5,1) -- (3.5,1);
\filldraw[fill=white] (4,2) circle [x radius=0.5,y radius=0.25];
\filldraw[fill=white,dashed] (3.5,1) arc [x radius=0.5,y radius=0.25,start angle=180,end angle=0];
\filldraw[fill=white] (3.5,1) arc [x radius=0.5,y radius=0.25,start angle=180,end angle=360];
\draw (3.5,1) to[out=45,in=-45] (3.5,2);
\draw (4.5,1) to[out=135,in=-135] (4.5,2);
\draw[dashed] (3.66,1.25) arc [x radius=0.34,y radius=0.17,start angle=180,end angle=0];
\draw (3.66,1.25) arc [x radius=0.34,y radius=0.17,start angle=180,end angle=360];
\draw[dashed] (3.66,1.75) arc [x radius=0.34,y radius=0.17,start angle=180,end angle=0];
\draw (3.66,1.75) arc [x radius=0.34,y radius=0.17,start angle=180,end angle=360];
\draw[dashed] (3.71,1.5) arc [x radius=0.29,y radius=0.145,start angle=180,end angle=0];
\draw (3.71,1.5) arc [x radius=0.29,y radius=0.145,start angle=180,end angle=360];
\filldraw[fill=white,draw=none] (5.5,1) to[out=45,in=-45] (5.5,2) -- (6.5,2) to[out=-135,in=135] (6.5,1) -- (5.5,1);
\filldraw[fill=white,draw=violet,thick] (6,2) circle [x radius=0.5,y radius=0.25];
\filldraw[fill=white,dashed] (5.5,1) arc [x radius=0.5,y radius=0.25,start angle=180,end angle=0];
\filldraw[fill=white] (5.5,1) arc [x radius=0.5,y radius=0.25,start angle=180,end angle=360];
\draw (5.5,1) to[out=45,in=-45] (5.5,2);
\draw (6.5,1) to[out=135,in=-135] (6.5,2);
\draw[dashed] (5.66,1.25) arc [x radius=0.34,y radius=0.17,start angle=180,end angle=0];
\draw (5.66,1.25) arc [x radius=0.34,y radius=0.17,start angle=180,end angle=360];
\draw[dashed] (5.66,1.75) arc [x radius=0.34,y radius=0.17,start angle=180,end angle=0];
\draw (5.66,1.75) arc [x radius=0.34,y radius=0.17,start angle=180,end angle=360];
\draw[dashed] (5.71,1.5) arc [x radius=0.29,y radius=0.145,start angle=180,end angle=0];
\draw (5.71,1.5) arc [x radius=0.29,y radius=0.145,start angle=180,end angle=360];
\draw[thick] (4,2) circle [x radius=2.5,y radius=1.25];
\draw[thick] (4,2) circle [x radius=1.5,y radius=0.75];
\fill[even odd rule,opacity=0.25] (4,2) circle[x radius=2.5,y radius=1.25] circle[x radius=1.5,y radius=0.75];
\end{tikzpicture}
\end{center}

\begin{center}
\begin{tikzpicture}
\filldraw[fill=lightgray,opacity=0.5,draw=none] (0,-0.5) -- (6.5,-0.5) -- (7.5,2.5) -- (1,2.5) -- cycle;
\draw (4,1) circle [x radius=1,y radius=0.5];
\draw (4,1) circle [x radius=1.5,y radius=0.75];
\draw (4,1) circle [x radius=2,y radius=1];
\draw (4,1) circle [x radius=2.5,y radius=1.25];
\node at (6.7,-0.2) {\(C_{c}\)};
\draw[->] (6,0) to[out=180,in=-90] (5.6,0.3);
\node at (6.9,2.9) {\(\bm{H}^{-1}(c^2,0)\)};
\draw[->] (5.55,2.8) to[out=180,in=90] (5,2.3);
\filldraw[fill=white,draw=none] (1.5,1) -- (2,1.5) -- (2.5,1) -- cycle;
\filldraw[fill=white,draw=none] (1.5,2) -- (2,1.5) -- (2.5,2) -- cycle;
\filldraw[fill=white] (2,2) circle [x radius=0.5,y radius=0.25];
\filldraw[fill=white,dashed] (1.5,1) arc [x radius=0.5,y radius=0.25,start angle=180,end angle=0];
\filldraw[fill=white] (1.5,1) arc [x radius=0.5,y radius=0.25,start angle=180,end angle=360];
\draw (1.5,1) -- (2.5,2);
\draw (1.5,2) -- (2.5,1);
\draw[dashed] (1.75,1.25) arc [x radius=0.25,y radius=0.125,start angle=180,end angle=0];
\draw (1.75,1.25) arc [x radius=0.25,y radius=0.125,start angle=180,end angle=360];
\draw[dashed] (1.75,1.75) arc [x radius=0.25,y radius=0.125,start angle=180,end angle=0];
\draw (1.75,1.75) arc [x radius=0.25,y radius=0.125,start angle=180,end angle=360];
\node at (2,1.5) {\(\bullet\)};
\filldraw[fill=white,draw=none] (3.5,1) to[out=45,in=-45] (3.5,2) -- (4.5,2) to[out=-135,in=135] (4.5,1) -- (3.5,1);
\filldraw[fill=white] (4,2) circle [x radius=0.5,y radius=0.25];
\filldraw[fill=white,dashed] (3.5,1) arc [x radius=0.5,y radius=0.25,start angle=180,end angle=0];
\filldraw[fill=white] (3.5,1) arc [x radius=0.5,y radius=0.25,start angle=180,end angle=360];
\draw (3.5,1) to[out=45,in=-45] (3.5,2);
\draw (4.5,1) to[out=135,in=-135] (4.5,2);
\draw[dashed] (3.66,1.25) arc [x radius=0.34,y radius=0.17,start angle=180,end angle=0];
\draw (3.66,1.25) arc [x radius=0.34,y radius=0.17,start angle=180,end angle=360];
\draw[dashed] (3.66,1.75) arc [x radius=0.34,y radius=0.17,start angle=180,end angle=0];
\draw (3.66,1.75) arc [x radius=0.34,y radius=0.17,start angle=180,end angle=360];
\draw[dashed] (3.71,1.5) arc [x radius=0.29,y radius=0.145,start angle=180,end angle=0];
\draw (3.71,1.5) arc [x radius=0.29,y radius=0.145,start angle=180,end angle=360];
\filldraw[fill=white,draw=none] (5.5,1) to[out=45,in=-45] (5.5,2) -- (6.5,2) to[out=-135,in=135] (6.5,1) -- (5.5,1);
\filldraw[fill=white] (6,2) circle [x radius=0.5,y radius=0.25];
\filldraw[fill=white,dashed] (5.5,1) arc [x radius=0.5,y radius=0.25,start angle=180,end angle=0];
\filldraw[fill=white] (5.5,1) arc [x radius=0.5,y radius=0.25,start angle=180,end angle=360];
\draw (5.5,1) to[out=45,in=-45] (5.5,2);
\draw (6.5,1) to[out=135,in=-135] (6.5,2);
\draw[dashed] (5.66,1.25) arc [x radius=0.34,y radius=0.17,start angle=180,end angle=0];
\draw (5.66,1.25) arc [x radius=0.34,y radius=0.17,start angle=180,end angle=360];
\draw[dashed] (5.66,1.75) arc [x radius=0.34,y radius=0.17,start angle=180,end angle=0];
\draw (5.66,1.75) arc [x radius=0.34,y radius=0.17,start angle=180,end angle=360];
\draw[dashed,thick] (5.71,1.5) arc [x radius=0.29,y radius=0.145,start angle=180,end angle=0];
\draw[thick] (5.71,1.5) arc [x radius=0.29,y radius=0.145,start angle=180,end angle=360];
\draw[thick] (3.855,1.5) circle [x radius=1.855,y radius=0.9275];
\draw[thick] (4.145,1.5) circle [x radius=2.145,y radius=1.0725];
\fill[even odd rule,opacity=0.25] (3.855,1.5) circle [x radius=1.855,y radius=0.9275] (4.145,1.5) circle [x radius=2.145,y radius=1.0725];
\end{tikzpicture}
\end{center}

This system has a focus-focus critical point at \((0,0)\). It
also has toric critical points along the conic
\(z_1z_2=c\). Exercise \ref{exr:auroux_hessian}: Given a
Hamiltonian system \(\bm{H}\) and a critical point \(\bm{x}\)
of \(\bm{H}\), let \(Q(\bm{H},\bm{x})\) denote the subspace of
the space of quadratic forms spanned by the
Hessians\index{Hessian} of the components of \(\bm{H}\). Check
that, after a suitable symplectic change of coordinates,
\(Q(\bm{H},0)=Q(\bm{F},0)\), where \(\bm{H}\) is the Auroux
system and \(\bm{F}\) is the standard focus-focus system from
Example \ref{exm:ff}. (This is enough to guarantee the
existence of a focus-focus chart; see Remark \ref{rmk:ffnf}).

\end{Example}
\begin{Lemma}\label{lma:zaffexample}
Let \(B\) be the set of regular fibres of the Auroux system
and \(\tilde{B}\) its universal cover. There is a fundamental
domain\index{action domain, fundamental!for Auroux system} for
the deck group action on \(\tilde{B}\) whose image under
action coordinates has the form \[\{(b_1,b_2)\ :\ 0\leq
b_1\leq \phi(b_2)\}\setminus\{(b_1,0)\ :\ b_1\geq
m\}\subset\RR^2\] for some function
\(\phi\colon\RR\to(0,\infty)\) and some number \(m>0\) (see
Figure \ref{fig:fadzaffexample}). The affine
monodromy\index{affine monodromy!of Auroux system|(}, on
crossing the branch cut \(\{(b_1,0)\ :\ b_1\geq m\}\)
clockwise, is \(\lmatrix 1 & 0 \\ 1 & 1\rmatrix\).

\end{Lemma}

\begin{Remark}
Finding \(\phi\) and \(m\) precisely along with the actual map
from the fundamental domain to this subset of \(\RR^2\) is a
nontrivial task.

\end{Remark}
\begin{Proof}[Proof of Lemma \ref{lma:zaffexample}]\label{prf:zaffexample}
The image \(\bm{H}(\CC^2)\) is the closed right half-plane:
\(H_1\) is always non-negative and \(H_2\) can take on any
value. The vertical boundary of the half-plane is the image of
the toric boundary\index{toric boundary!Auroux system} (the
conic \(z_1z_2=c\)\index{conic!fibration|)}). The point
\(\left(c^2,0\right)\) is the image of the focus-focus
critical point \((0,0)\) and
\(B=\bm{H}(\CC^2)\setminus\{(c^2,0)\}\).

The Hamiltonian \(H_2\) gives a \(2\pi\)-periodic flow, so the
change of coordinates of \(\RR^2\) which gives action
coordinates\index{action-angle coordinates!for Auroux
system|(} has the form \((b_1,b_2)\mapsto (G_1(b_1,b_2),b_2)\)
for some (multiply-valued) function \(G_1\). In particular,
the monodromy of the integral affine structure around the
focus-focus critical point simply shifts amongst the branches
of \(G_1\), so has the form \(\lmatrix 1 & 0 \\ 1 & 1
\rmatrix\). We may make a branch cut\index{branch cut!Auroux
system} along the line \(R=\left\{(b_1,0)\ :\
b_1>c^2\right\}\) to get a simply-connected open set
\(U=B\setminus R\) and pick a fundamental domain \(\tilde{U}\)
lying over \(U\) in the universal cover \(p\colon\tilde{B}\to
B\).

We first compute the image \(\{(G_1(0,b_2),b_2)\ :\
b_2\in\RR\}\) of the line \(0\times\RR\) under the action
coordinates. Since this is part of the toric boundary,
Proposition \ref{prp:straight_lines} implies this image is a
straight line \(S\) with rational slope. As observed in Lemma
\ref{lma:visibleff}, there is a visible Lagrangian disc
emanating from the focus-focus critical point and living over
an eigenline\index{eigenline!visible Lagrangian over} of the
affine monodromy.\index{affine monodromy!of Auroux system|)}
Actually, we can write the disc explicitly for the Auroux
system: it is the Lagrangian disc
\(\{(z,\bar{z})\ :\ |z|^2\leq c\}\) with boundary on
\(z_1z_2=c\). This visible disc lives over the horizontal line
segment \(\{(b_1,0)\ :\ b_1\leq c^2\}\) under the map
\(\bm{H}\) and hence\footnote{See Remark
\ref{rmk:visible_lag_to_aff_str}.} over a horizontal line
segment \(\{(G_1(b_1,0),0)\ :\ b_1\leq c^2\}\) in the image of
action coordinates. This line segment connects \(S\) to the
base-node \((G_1(c^2,0),0)\). Since this visible Lagrangian is
a disc, not a pinwheel core, comparison with the local models
from Example \ref{exm:hittingedge} shows that the line \(S\)
must have slope \(1/n\) for some integer \(n\). In particular,
post-composing action coordinates with an integral affine
shear \(\lmatrix 1 & 0 \\ -n & 1\rmatrix\), we get that \(S\)
is vertical (we always have the freedom to post-compose our
action coordinates\index{action-angle coordinates!for Auroux
system|)} with an integral affine transformation, thanks to
Lemma \ref{lma:ham_postcompose}). Now it is clear that the
fundamental action domain has the required form, where
\(\phi(b_2)=\sup_{b_1\in[0,\infty)}G_1(b_1,b_2)\) and
\(m=G_1(c^2,0)\). \qedhere

\end{Proof}
\section{Different branch cuts}
\label{sct:aurouxnodaltrade}

We\index{branch cut!changing|(} can always pick a different
simply-connected domain \(U\subset B\) to get well-defined
action coordinates \(\II\), as we illustrated in Figure
\ref{fig:ngoc12} in the previous chapter. The image of \(U\)
will not in general ``close-up'': unless we take a branch cut
along the eigendirection\index{eigenline!branch cut} of the
affine monodromy\index{affine monodromy!relating branch cuts},
the boundary of \(\II(U)\) will be two branch cuts related by
the affine monodromy.

To illustrate this, we plot some of the associated pictures
below for the Auroux system as the branch cut under goes a full
rotation. It is important to emphasise that all of these are
integral affine bases for the {\em same} Hamiltonian system on
the {\em same} manifold; they differ only in the choice of a
fundamental domain for the covering space \(\tilde{B}\to B\).

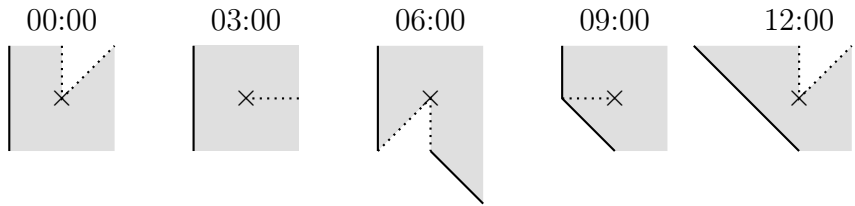
\begin{figure}[htb]
\begin{center}
\begin{tikzpicture}[baseline=1,scale=0.7]
\node at (1,2.5) {00:00};
\filldraw[fill=lightgray,opacity=0.5,draw=none] (0,0) -- (0,2) -- (1,2) -- (1,1) -- (2,2) -- (2,0) -- cycle;
\draw[thick,black] (0,0) -- (0,2);
\filldraw[fill=white,draw=black,dotted,thick] (1,2) -- (1,1) -- (2,2);
\node at (1,1) {\(\times\)};
\begin{scope}[shift={(3.5,0)}]
\node at (1,2.5) {03:00};
\filldraw[fill=lightgray,opacity=0.5,draw=none] (0,0) -- (0,2) -- (2,2) -- (2,0) -- cycle;
\draw[thick,black] (0,0) -- (0,2);
\draw[dotted,thick] (2,1) -- (1,1);
\node at (1,1) {\(\times\)};
\end{scope}
\begin{scope}[shift={(7,0)}]
\node at (1,2.5) {06:00};
\filldraw[fill=lightgray,opacity=0.5,draw=none] (0,0) -- (0,2) -- (2,2) -- (2,-1) -- (1,0) -- (1,1) -- cycle;
\draw[thick,black] (0,0) -- (0,2);
\draw[thick,black] (1,0) -- (2,-1);
\filldraw[fill=white,draw=black,dotted,thick] (0,0) -- (1,1) -- (1,0);
\node at (1,1) {\(\times\)};
\end{scope}
\begin{scope}[shift={(10.5,0)}]
\node at (1,2.5) {09:00};
\filldraw[fill=lightgray,opacity=0.5,draw=none] (0,1) -- (0,2) -- (2,2) -- (2,0) -- (1,0) -- cycle;
\draw[thick,black] (0,2) -- (0,1) -- (1,0);
\filldraw[fill=white,draw=black,dotted,thick] (0,1) -- (1,1);
\node at (1,1) {\(\times\)};
\end{scope}
\begin{scope}[shift={(14,0)}]
\node at (1,2.5) {12:00};
\filldraw[fill=lightgray,opacity=0.5,draw=none] (1,0) -- (-1,2) -- (1,2) -- (1,1) -- (2,2) -- (2,0) -- cycle;
\draw[thick,black] (1,0) -- (-1,2);
\filldraw[fill=white,draw=black,dotted,thick] (1,2) -- (1,1) -- (2,2);
\node at (1,1) {\(\times\)};
\end{scope}
\end{tikzpicture}
\end{center}
\caption{The Auroux system seen with different branch cuts; as we move from left to right in the figure, we see the branch cut rotate (from the 12-o'clock position) by 360 degrees. The final picture is related to the first by the affine monodromy.}\label{fig:branchcuts}
\end{figure}

\begin{Remark}\label{rmk:brokenlines}
In some of these pictures, the toric
boundary\index{toric boundary!appears broken} appears
``broken''. This is an artefact of the fact that it intersects
the branch cut: the two segments of the toric boundary are
related by the affine monodromy and therefore form one
straight line in the integral affine structure. If you want to
check this, the {\em anticlockwise} affine monodromy is
\(\lmatrix 1 & 0 \\ -1 & 1\rmatrix\), so, for example in the
9:00 diagram, the tangent vector \((0,-1)\) to the line above
the branch cut gets sent to \((0,-1)\lmatrix 1 & 0\\-1 &
1\rmatrix=(1,-1)\) below the branch cut, which is tangent to
the continuation of the boundary.\index{Auroux system|)}

\end{Remark}
\begin{Remark}
We can apply an integral affine transformation to any of these
diagrams. Applying the matrix \(\lmatrix 1 & 1 \\ 0 &
1\rmatrix\) to the 09:00 diagram in Figure
\ref{fig:branchcuts} yields Figure \ref{fig:nodaltrade0} which
will be important in the next chapter and which has
anticlockwise affine monodromy \(\lmatrix 2 & 1 \\ -1 &
0\rmatrix\) by Lemma \ref{lma:monodromy_formula}. The
importance of this example is that away from the branch cut,
the integral affine manifold looks like the standard Delzant
corner. We will see that this means we can always ``implant''
this local Hamiltonian system whenever we have a polygon with
a standard Delzant corner, an operation known as a {\em nodal
trade}.\index{nodal trade}\index{branch cut!changing|)}

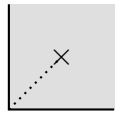
\begin{figure}[htb]
\begin{center}
\begin{tikzpicture}[baseline=1,scale=0.7]
\filldraw[fill=lightgray,opacity=0.5,draw=none] (0,0) -- (0,2) -- (2,2) -- (2,0) -- cycle;
\draw[thick,black] (2,0) -- (0,0) -- (0,2);
\draw[dotted,thick] (0,0) -- (1,1);
\node at (1,1) {\(\times\)};
\end{tikzpicture}
\end{center}
\caption{Another fundamental action domain for the Auroux system.}\label{fig:nodaltrade0}
\end{figure}

\end{Remark}
\section{Smoothing \(A_n\) singularities}

\begin{Example}\label{exm:anmilnorfibre}
Let \(P(z)\) be a polynomial of degree \(n+1\) with \(n+1\)
distinct roots and \(P(0)\neq 0\). Let
\(M_P=\{(z_1,z_2,z_3)\in\CC^3\ :\ z_1z_2+P(z_3)=0\}\). If you
allow \(P\) to vary, you get a family of such varieties; as
\(P\) approaches the degenerate polynomial \(P(z)=z^{n+1}\),
the variety \(M_P\) develops a singularity called an
\(A_n\)-singularity\index{singularity!An@$A_n$}. In other
words, for generic \(P\) (with distinct roots) the variety
\(M_P\) is the {\em Milnor fibre}\index{Milnor fibre!of An
singularity@of $A_n$ singularity} (or {\em
smoothing}\index{smoothing}\footnote{Algebraic geometers
usually say ``smoothing'' to mean the total space of a family
which smooths a singularity; some other people say
``smoothing'' to mean the smooth fibre of such a family.}) of
the \(A_n\) singularity (see Milnor's book
\cite{MilnorSingularities} for more about Milnor
fibres). Milnor fibres of singularities provide a rich class
of symplectic manifolds which have been intensively studied.

Let \(\pi\colon M_P\to\CC\) be the conic
fibration\index{conic!fibration} \(\pi(z_1,z_2,z_3)=z_3\). By
analogy with the Auroux system, we define
\[\bm{H}(z_1,z_2,z_3) = \left(|z_3|^2,
\frac{1}{2}\left(|z_1|^2-|z_2|^2\right)\right).\] Again, these
Hamiltonians commute with one another, but only \(H_2\)
generates a circle action.

The subvariety \(z_3=0\), \(z_1z_2+P(0)=0\) is a conic along
which the Hamiltonian system has toric critical points; this
projects to the line \(\{(0,b_2)\,:\,b_2\in\RR\}\) under
\(\bm{H}\).

The level sets \(\bm{H}^{-1}(b_1,b_2)\) for \(b_2\neq 0\) are
Lagrangian tori, and the level sets \(\bm{H}^{-1}(b_1,0)\) are
Lagrangian tori unless the circle \(|z_3|^2=b_1\) contains a
root of \(P\). If this circle contains \(k\) roots of \(P\)
then the fibre \(\bm{H}^{-1}(b_1,0)\) is a Lagrangian torus
with \(k\) pinches.
\begin{itemize}
\item For example, if \(P(z)=z^{n+1}-1\) then the fibre
\(\bm{H}^{-1}(1,0)\) is the only critical fibre; it has
\(n+1\) focus-focus critical points (see Figure
\ref{fig:anmilnorfibre_1}).
\item If \(0<a_1<a_2<\ldots<a_{n+1}\) are real numbers and \(P(z)
= (z-a_1)(z-a_2)\cdots(z-a_{n+1})\) then there are \(n+1\)
focus-focus fibres which project via \(\bm{H}\) to the
points \(\{(a_1^2,0),(a_2^2,0),\ldots,(a_{n+1}^2,0)\}\) (see
Figure \ref{fig:anmilnorfibre_2}).

\end{itemize}
\begin{figure}[htb]
\begin{center}
\begin{tikzpicture}[scale=0.8]
\SmoothConic{0}{0}
\SingularConic{-3}{0}
\SingularConic{0}{2}
\SingularConic{3}{0}
\SingularConic{0}{-2}
\filldraw[opacity=0.25] (-3,0) to[out=-90,in=180] (0,-2) to[out=-90,in=-90] (-4,-1.6) to[out=90,in=180] (-3,0);
\draw[thick] (-3,0) to[out=-90,in=180] (0,-2);
\draw[thick] (0,-2) to[out=-90,in=-90] (-4,-1.6);
\draw[thick] (-4,-1.6) to[out=90,in=180] (-3,0);
\draw (-4,-1.6) to[out=-45,in=-90] (-2.65,-1);
\draw[dotted,thick] (-2.65,-1) to[out=110,in=105] (-4,-1.6);
\draw (-2,-2.98) to[out=0,in=-70] (-1,-1.9);
\draw[dotted,thick] (-1,-1.9) to[out=180,in=105] (-2,-2.98);
\begin{scope}[xscale=-1,yscale=1]
\filldraw[opacity=0.25] (-3,0) to[out=-90,in=180] (0,-2) to[out=-90,in=-90] (-4,-1.6) to[out=90,in=180] (-3,0);
\draw[thick] (-3,0) to[out=-90,in=180] (0,-2);
\draw[thick] (0,-2) to[out=-90,in=-90] (-4,-1.6);
\draw[thick] (-4,-1.6) to[out=90,in=180] (-3,0);
\draw (-4,-1.6) to[out=-45,in=-90] (-2.65,-1);
\draw[dotted,thick] (-2.65,-1) to[out=110,in=105] (-4,-1.6);
\draw (-2,-2.98) to[out=0,in=-70] (-1,-1.9);
\draw[dotted,thick] (-1,-1.9) to[out=180,in=105] (-2,-2.98);
\end{scope}
\begin{scope}[xscale=1,yscale=-1]
\filldraw[opacity=0.25] (-3,0) to[out=-90,in=180] (0,-2) to[out=-90,in=-90] (-4,-1.6) to[out=90,in=180] (-3,0);
\draw[thick] (-3,0) to[out=-90,in=180] (0,-2);
\draw[thick] (0,-2) to[out=-90,in=-90] (-4,-1.6);
\draw[thick] (-4,-1.6) to[out=90,in=180] (-3,0);
\draw (-4,-1.6) to[out=-45,in=-90] (-2.65,-1);
\draw[dotted,thick] (-2.65,-1) to[out=110,in=105] (-4,-1.6);
\draw (-2,-2.98) to[out=0,in=-70] (-1,-1.9);
\draw[dotted,thick] (-1,-1.9) to[out=180,in=105] (-2,-2.98);
\end{scope}
\begin{scope}[xscale=-1,yscale=-1]
\filldraw[opacity=0.25] (-3,0) to[out=-90,in=180] (0,-2) to[out=-90,in=-90] (-4,-1.6) to[out=90,in=180] (-3,0);
\draw[thick] (-3,0) to[out=-90,in=180] (0,-2);
\draw[thick] (0,-2) to[out=-90,in=-90] (-4,-1.6);
\draw[thick] (-4,-1.6) to[out=90,in=180] (-3,0);
\draw (-4,-1.6) to[out=-45,in=-90] (-2.65,-1);
\draw[dotted,thick] (-2.65,-1) to[out=110,in=105] (-4,-1.6);
\draw (-2,-2.98) to[out=0,in=-70] (-1,-1.9);
\draw[dotted,thick] (-1,-1.9) to[out=180,in=105] (-2,-2.98);
\end{scope}
\end{tikzpicture}
\end{center}
\caption{The Hamiltonian system from Example \ref{exm:anmilnorfibre} with \(P=z^4+1\) (\(n=3\)). There is a single fibre with 4 focus-focus critical points and a smooth conic which consists of toric critical points.}
\label{fig:anmilnorfibre_1}
\end{figure}
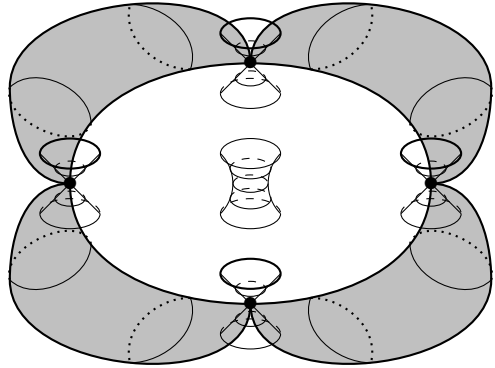

\begin{figure}[htb]
\begin{center}
\begin{tikzpicture}[scale=0.7]
\SmoothConic{-2}{0}
\SingularConic{0}{0}
\SingularConic{2}{0}
\SingularConic{4}{0}
\filldraw[fill=lightgray,opacity=0.5] (0,0) to[out=-90,in=-90] (-4,0) to[out=90,in=90] (0,0) to[out=90,in=90] (-5,0) to[out=-90,in=-90] (0,0);
\draw (0,0) to[out=-90,in=-90] (-4,0) to[out=90,in=90] (0,0) to[out=90,in=90] (-5,0) to[out=-90,in=-90] (0,0);
\filldraw[fill=lightgray,opacity=0.5] (2,0) to[out=-90,in=-90] (-6,0) to[out=90,in=90] (2,0) to[out=90,in=90] (-7,0) to[out=-90,in=-90] (2,0);
\draw (2,0) to[out=-90,in=-90] (-6,0) to[out=90,in=90] (2,0) to[out=90,in=90] (-7,0) to[out=-90,in=-90] (2,0);
\filldraw[fill=lightgray,opacity=0.5] (4,0) to[out=-90,in=-90] (-8,0) to[out=90,in=90] (4,0) to[out=90,in=90] (-9,0) to[out=-90,in=-90] (4,0);
\draw (4,0) to[out=-90,in=-90] (-8,0) to[out=90,in=90] (4,0) to[out=90,in=90] (-9,0) to[out=-90,in=-90] (4,0);
\filldraw[fill=lightgray,opacity=0.5] (1,0) circle [x radius = 1,y radius = 0.3];
\draw (1,0) circle [x radius = 1,y radius = 0.3];
\filldraw[fill=lightgray,opacity=0.5] (3,0) circle [x radius = 1,y radius = 0.3];
\draw (3,0) circle [x radius = 1,y radius = 0.3];
\draw[thick] (1,0.3) arc [x radius = 0.3, y radius = 0.3, start angle = 90, end angle = -90];
\draw[dotted,thick] (1,0.3) arc [x radius = 0.3, y radius = 0.3, start angle = 90, end angle = 270];
\draw[thick] (3,0.3) arc [x radius = 0.3, y radius = 0.3, start angle = 90, end angle = -90];
\draw[dotted,thick] (3,0.3) arc [x radius = 0.3, y radius = 0.3, start angle = 90, end angle = 270];
\filldraw[fill=lightgray,opacity=0.5] (-10,0.3) to[out=0,in=180] (-2,0.15) to[out=0,in=150] (0,0) to[out=-150,in=0] (-2,-0.15) to[out=180,in=0] (-10,-0.3);
\draw (-10,0.3) to[out=0,in=180] (-2,0.15) to[out=0,in=150] (0,0) to[out=-150,in=0] (-2,-0.15) to[out=180,in=0] (-10,-0.3);
\draw[thick] (-2,0) circle [x radius=0.3,y radius=0.15];
\draw[thick,dotted] (-4.5,0) circle [x radius=0.5,y radius=0.15];
\draw[thick] (-4.5,0.15) arc [x radius=0.5,y radius=0.15,start angle=90,end angle=-90];
\draw[thick,dotted] (-6.5,0) circle [x radius=0.5,y radius=0.2];
\draw[thick] (-6.5,0.2) arc [x radius=0.5,y radius=0.2,start angle=90,end angle=-90];
\draw[thick,dotted] (-8.5,0) circle [x radius=0.5,y radius=0.3];
\draw[thick] (-8.5,0.3) arc [x radius=0.5,y radius=0.3,start angle=90,end angle=-90];
\filldraw[fill=lightgray,opacity=0.5] (8,0.3) -- (6,0.3) to[out=180,in=30] (4,0) to[out=-30,in=180] (6,-0.3) -- (8,-0.3);
\draw (8,0.3) -- (6,0.3) to[out=180,in=30] (4,0) to[out=-30,in=180] (6,-0.3) -- (8,-0.3);
\draw[thick] (6,0.3) to[out=-45,in=45] (6,-0.3);
\draw[dotted,thick] (6,0.3) to[out=-135,in=135] (6,-0.3);
\end{tikzpicture}
\caption{The Hamiltonian system from Example \ref{exm:anmilnorfibre} with \(P=(z-a_1)(z-a_2)(z-a_3)\). There is a smooth conic consisting of toric critical points, and three focus-focus fibres which encircle it. We also show, horizontally across the figure, the visible Lagrangian submanifolds described in Remark \ref{rmk:vis_lag_milnor_fibre}.}
\label{fig:anmilnorfibre_2}
\end{center}
\end{figure}
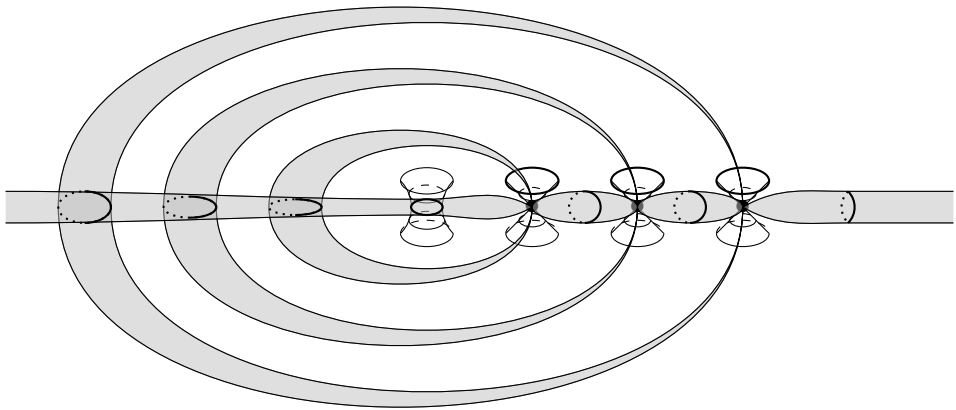

\clearpage

In any case, the image of the developing map can be analysed
as in the Auroux system. We get the diagram shown in Figure
\ref{fig:anmilnor_atbd_1}.\index{almost toric base
diagram!Milnor fibre of An singularity@Milnor fibre of $A_n$
singularity|(}

\begin{figure}[htb]
\begin{center}
\begin{tikzpicture}
\filldraw[fill=lightgray,opacity=0.5,draw=none] (0,-2) -- (0,2) -- (4,2) -- (4,-2) -- cycle;
\draw (0,-2) -- (0,2);
\draw[dotted,thick] (4,0) -- (3,0) node {\(\times\)} -- (2,0) node {\(\times\)} -- (1,0) node {\(\times\)};
\end{tikzpicture}
\caption{Almost toric diagram for \(A_2\) Milnor fibre.}
\label{fig:anmilnor_atbd_1}
\end{center}
\end{figure}
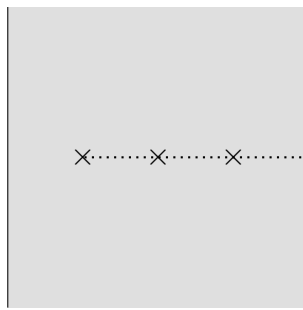

The number of base-nodes here is the number of values of
\(b_1\) for which \(P(z)\) has a zero with \(|z|^2=b_1\). We
see that the affine monodromy\index{affine monodromy!for An
Milnor fibre@for $A_n$ Milnor fibre} as we cross the branch
cut at position \((t,0)\) is the product of all the individual
affine monodromies for base-nodes with \(b_1<t\), all of which
are \(\lmatrix 1 & 0 \\ 1 & 1\rmatrix \) (just as in the
Auroux system). In particular, the ``total monodromy'' as we
cross the branch cut far to the right is \(\lmatrix 1 & 0
\\ n+1 & 1\rmatrix \). If we change the branch cut by \(180\)
degrees clockwise then we get the diagram in Figure
\ref{fig:anmilnor_atbd_2} (drawn in the case \(n=2\)):

\begin{figure}[htb]
\begin{center}
\begin{tikzpicture}
\filldraw[fill=lightgray,opacity=0.5,draw=none] (0,0) -- (0,2) -- (4,2) -- (4,-1) -- (3,-1) -- cycle;
\draw[thick,->] (0,2) -- (0,0) -- (3,-1) node [below] {\((n+1,-1)\)};
\draw[dotted,thick] (0,0) -- (1,0) node {\(\times\)} -- (2,0) node {\(\times\)} -- (3,0) node {\(\times\)};
\end{tikzpicture}
\caption{Another view of Figure \ref{fig:anmilnor_atbd_1}. Here \(n=2\) and there are \(n+1\) singularities.}
\label{fig:anmilnor_atbd_2}
\end{center}
\end{figure}
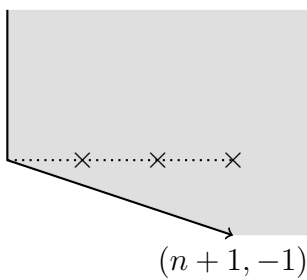

If we apply the integral affine transformation \(\lmatrix 1 &
1\\ 0& 1\rmatrix\) to this diagram, we get Figure
\ref{fig:anmilnor_atbd_3}.

\begin{figure}[htb]
\begin{center}
\begin{tikzpicture}
\filldraw[fill=lightgray,opacity=0.5,draw=none] (0,0) -- (0,10/3) -- (5,10/3) -- cycle;
\draw[thick,->] (0,10/3) -- (0,0) -- (5,10/3) node [right] {\((n+1,n)\)};
\draw[dotted,thick] (0,0) -- (1,1) node {\(\times\)} -- (2,2) node {\(\times\)} -- (3,3) node {\(\times\)};
\end{tikzpicture}
\caption{A third view of Figure \ref{fig:anmilnor_atbd_1}.}
\label{fig:anmilnor_atbd_3}
\end{center}
\end{figure}
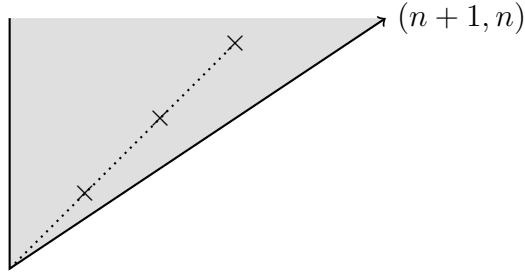

Compare this with the moment polygon\index{moment
polytope!cyclic quotient singularity} \(\pi(n+1,n)\) from
Example \ref{exm:cycquot} for the cyclic quotient singularity
\(\frac{1}{n+1}(1,n)\). Indeed, this is precisely the
\(A_n\)-singularity\index{singularity!An@$A_n$} mentioned
above. As the polynomial \(P(z)\) degenerates to \(z^{n+1}\),
the base-nodes move in the diagram along the dotted line
towards the non-Delzant corner\index{polytope!non-Delzant}.

\end{Example}
\begin{Remark}\label{rmk:vis_lag_milnor_fibre}
There are some visible
Lagrangians\index{Lagrangian!submanifold!visible} in the
Milnor fibre \(M_P\) when \(P(z)=
(z-a_1)(z-a_2)\cdots(z-a_{n+1})\). Namely, consider the
antisymplectic involution\index{antisymplectic involution}
\((z_1,z_2,z_3)\mapsto (\bar{z}_2,\bar{z}_1,\bar{z}_3)\). The
fixed locus consists of points
\(\{(z_1,\bar{z}_1,z_3)\,:\,|z_1|^2=-P(z_3),\,z_3\in
\RR\}\). The fixed locus of an antisymplectic involution is
always a Lagrangian submanifold; in this case, it consists of:
\begin{itemize}
\item Lagrangian spheres\index{Lagrangian!sphere} \(z_3\in[a_k,a_{k+1}]\), \(k\)
even,
\item the Lagrangian plane \(z_3\in(-\infty,a_1]\) if \(n\) is even.
\end{itemize}
The involution \((z_1,z_2,z_3)\mapsto
(-\bar{z}_2,\bar{z}_1,\bar{z}_3)\) gives more Lagrangians
\(\{(z_1,-\bar{z}_1,z_3)\,:\,|z_1|^2=P(z_3),\,z_3\in\RR\}\),
which consists of:
\begin{itemize}
\item Lagrangian spheres \(z_3\in[a_k,a_{k+1}]\), \(k\)
odd,
\item the Lagrangian plane \(z_3\in[a_{n+1},\infty)\),
\item the Lagrangian plane \(z_3\in(-\infty,a_1]\) if \(n\) is
odd.
\end{itemize}
These are all visible Lagrangians mapping to the line
\(H_2=0\): the Lagrangian spheres project to the compact
segments connecting focus-focus fibres; the plane
\(z_3\in[a_{n+1},\infty)\) projects to the segment connecting
the right-most focus-focus fibre to infinity. The plane
\(z_3\in(-\infty,a_1]\) has a more singular projection: the
disc \(z_3\in[0,a_1]\) projects to the segment connecting the
left-most focus-focus fibre to the toric
boundary\index{toric boundary!of An Milnor fibre@of $A_n$
Milnor fibre}; the annulus \(z_3\in[-\infty,0]\) projects to
the whole ray \(H_2=0\) emanating from the toric boundary. In
other words, this visible Lagrangian ``folds over itself'' at
the toric boundary. See Figure \ref{fig:vis_lag_milnor_fibre}
for the images of these visible Lagrangians under the action
map, and Figure \ref{fig:anmilnorfibre_2} to see how they look
in the total space of the conic fibration. In what follows,
the visible disc \(z_3\in[0,1]\) will be more important; we
will denote it by \(\Delta\).\index{almost toric base
diagram!Milnor fibre of An singularity@Milnor fibre of
$A_n$ singularity|)}

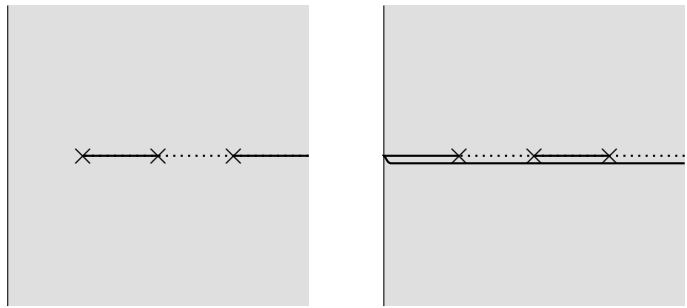
\begin{figure}[htb]
\begin{center}
\begin{tikzpicture}
\filldraw[fill=lightgray,opacity=0.5,draw=none] (0,-2) -- (0,2) -- (4,2) -- (4,-2) -- cycle;
\draw (0,-2) -- (0,2);
\draw[dotted,thick] (4,0) -- (3,0) node {\(\times\)} -- (2,0) node {\(\times\)} -- (1,0) node {\(\times\)};
\draw[thick] (4,0) -- (3,0);
\draw[thick] (2,0) -- (1,0);
\begin{scope}[shift={(5,0)}]
\filldraw[fill=lightgray,opacity=0.5,draw=none] (0,-2) -- (0,2) -- (4,2) -- (4,-2) -- cycle;
\draw (0,-2) -- (0,2);
\draw[dotted,thick] (4,0) -- (3,0) node {\(\times\)} -- (2,0) node {\(\times\)} -- (1,0) node {\(\times\)};
\draw[thick] (3,0) -- (2,0);
\draw[thick] (1,0) -- (0,0) to[out=-45,in=180] (0.1,-0.1) -- (4,-0.1);
\end{scope}
\end{tikzpicture}
\end{center}
\caption{The visible Lagrangians described in Remark \ref{rmk:vis_lag_milnor_fibre} in the case \(n=2\). On the left, we see the components of the fixed locus of \((z_1,z_2,z_3)\mapsto (\bar{z}_2,\bar{z}_1,\bar{z}_3)\): a Lagrangian sphere and a Lagrangian plane. On the right we see the components of the fixed locus of \((z_1,z_2,z_3)\mapsto (-\bar{z}_2,\bar{z}_1,\bar{z}_3)\): a Lagrangian sphere and a Lagrangian plane whose projection ``folds over itself''. In reality, this projection is contained in a single horizontal line; we have separated it for clarity.}
\label{fig:vis_lag_milnor_fibre}
\end{figure}

\end{Remark}
\section{Smoothing cyclic quotient T-singularities}

\begin{Example}\label{exm:bpq}
Let \(d\geq 1\) be an integer, \(p,q\) be coprime positive
integers with \(1\leq q<p\), and \(0<a_1<\ldots<a_d\) be real
numbers. Let \(P\) be the polynomial
\(P(z)=(z^p-a_1)(z^p-a_2)\cdots(z^p-a_d)\). Consider the
action of the group \(\mu_{p}\) of \(p\)th roots of unity on
the variety \(M_P\) from Example \ref{exm:anmilnorfibre} given
by \(\mu\cdot(z_1,z_2,z_3)=(\mu z_1,\mu^{-1}z_2,\mu^qz_3)\),
\(\mu\in\mu_p\). This action is free and
\(\pi(\mu\cdot(z_1,z_2,z_3))=\mu\pi(z_1,z_2,z_3)\).

The Hamiltonian system \(\bm{H}\) on \(M_P\) from Example
\ref{exm:anmilnorfibre} has a line of toric critical points
along \(H_1=0\) and \(d\) isolated critical fibres with
\(H_1=1,2,\ldots,d\), each of which has \(p\) focus-focus
critical points. The \(\mu_p\)-action preserves the critical
fibres: the \(p\) focus-focus critical points in each fibre
form a \(\mu_p\)-orbit. The Hamiltonian system \(\bm{H}\)
descends to give a system \(\bm{G}\colon B_{d,p,q}\to\RR^2\)
on the quotient space\index{Bdpq@$B_{d,p,q}$|(}
\(B_{d,p,q}:=M_P/\mu_p\) with \(d\) focus-focus critical
points and \(\bm{H}(M_P)=\bm{G}(B_{p,q})\). However, the
action coordinates are different: quotienting by the
\(\mu_p\)-action changes the period lattice (compare with
Example \ref{exm:cycquot}). In fact, a fundamental action
domain for \(\bm{G}\) is the polygon \(\pi(dp^2,dpq-1)\), and
there are \(d\) base-nodes. The branch cut\index{branch
cut!Bdpq@$B_{d,p,q}$} is along a line pointing in the
\((p,q)\)-direction, which is an
eigenvector\index{eigenline!for Bdpq@for $B_{d,p,q}$} of the
affine monodromy\index{affine monodromy!for Bdpq@for
$B_{d,p,q}$}; the base-nodes all lie on this branch cut. See
Figure \ref{fig:b_dpq_diagram}\index{almost toric base
diagram!Bdpq@$B_{d,p,q}$|(}.

\begin{figure}[htb]
\begin{center}
\begin{tikzpicture}
\filldraw[fill=lightgray,opacity=0.5,draw=none] (0,3) -- (0,0) -- (8,3) -- cycle;
\draw[->] (0,3) -- (0,0) -- (8,3) node [right] {\((dp^2,dpq-1)\)};
\draw[dotted,thick] (0,0) -- (2,1) node {\(\times\)} -- (3,1.5) node {\(\times\)} -- (4,2) node {\(\times\)};
\end{tikzpicture}
\end{center}
\caption{The fundamental action domain for \(B_{d,p,q}\), shown in the case \(d=2,p=2,q=1\).}
\label{fig:b_dpq_diagram}
\end{figure}
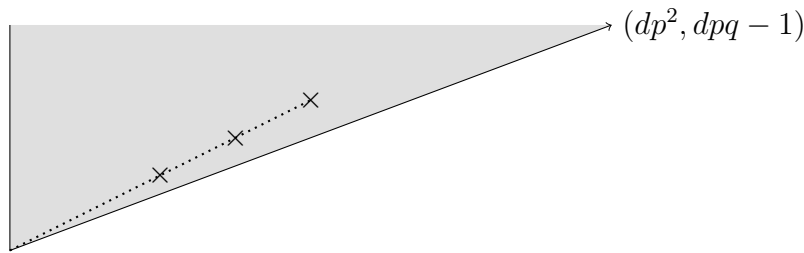

\end{Example}
\begin{Remark}[Advertising]\label{rmk:advertising}
The manifold \(B_{d,p,q}\) is the Milnor fibre\index{Milnor
fibre!of cyclic quotient T-singularity} of the cyclic quotient
singularity\index{singularity!cyclic quotient|(}
\(\frac{1}{dp^2}(1,dpq-1)\). Cyclic quotient singularities of
this form are called {\em cyclic quotient T-singularities},
and are the most general cyclic quotient surface singularities
admitting a
\(\QQ\)-Gorenstein\index{smoothing!Q-Gorenstein@$\mathbb{Q}$-Gorenstein}
smoothing\footnote{This means that the total space of the
smoothing is \(\QQ\)-Gorenstein; this condition picks out a
distinguished deformation class of smoothings {\cite[Theorem
3.9]{KSB}}. Other cyclic quotient singularites can be
smoothed, but the total space of the smoothing is not
\(\QQ\)-Gorenstein.} see {\cite[Proposition
5.9]{LooijengaWahl}} or {\cite[Proposition
3.10]{KSB}}. Perhaps the case that has attracted the most
attention is the case \(B_{1,p,q}\), often abbreviated to
\(B_{p,q}\), because in that case the Milnor fibre has
\(H_*(B_{1,p,q};\QQ)=H_*(B^4;\QQ)\), i.e.\ it is a rational
homology ball. We will show this below. This makes the
manifold \(B_{1,p,q}\) a useful building block for
constructing exotic 4-dimensional manifolds\index{small exotic
4-manifolds} with small homology groups, for example using the
rational blow-down\index{rational blow-down} construction
\cite{FintushelStern}.

The symplectic geometry of the manifolds \(B_{1,p,q}\) has
also been studied. Lekili and Maydanskiy
\cite{LekiliMaydanskiy} showed that \(B_{1,p,q}\) contains no
compact exact Lagrangian submanifolds despite having nonzero
symplectic cohomology\footnote{The ``standard way'' to rule
out exact Lagrangian submanifolds is to show that the
symplectic cohomology vanishes.}. Karabas \cite{Karabas}
showed that the Kontsevich cosheaf conjecture holds for
\(B_{1,p,1}\), i.e.\ that the wrapped Fukaya category of
\(B_{1,p,1}\) can be calculated using microlocal sheaf theory
on the Lagrangian
skeleton\index{skeleton!Bdpq@$B_{d,p,q}$|(}\index{Lagrangian
skeleton|see {skeleton}} discussed in Lemma
\ref{lma:lag_skel_b_dpq} below. Evans and Smith
\cite{EvansSmith1,EvansSmith2}, building on ideas of
Khodorovskiy \cite{Khodorovskiy}, used obstructions to
symplectic embeddings\index{symplectic embeddings of rational
homology balls} of \(B_{1,p,q}\) to obtain
restrictions\footnote{Compare with the algebro-geometric
approaches to these problems in the work of Hacking and
Prokhorov \cite{HackingProkhorov} and the work of Rana and
Urz\'{u}a \cite{RanaUrzua}.} on which cyclic quotient
singularities can occur under stable degenerations of complex
surfaces.\index{singularity!cyclic quotient|)}

\end{Remark}
\begin{Remark}\label{rmk:vis_lag_b_dpq}
As in Remark \ref{rmk:vis_lag_milnor_fibre}, the manifold
\(M_P\) with \(P(z) = (z^p-a_1) (z^p-a_2) \cdots (z^p-a_d)\)
contains \(p(d-1)\) visible Lagrangian spheres \[S_{1,1},\
S_{1,2},\ \ldots,\ S_{1,p},\ S_{2,1},\ \ldots,\ S_{d-1,p}.\]
These can be obtained by taking fixed loci as before and
applying the \(\mu_p\) action to the result. When we quotient
by \(\mu_p\), the spheres descend to give \(d-1\) visible
Lagrangian spheres \(\overline{S}_1, \ldots,
\overline{S}_{d-1}\) in \(B_{d,p,q}\). We also obtain \(p\)
visible Lagrangian discs \(\Delta_1, \ldots, \Delta_p\) with
common boundary along the toric
boundary\index{toric boundary!of Bdpq@of $B_{d,p,q}$}. These
discs descend to a visible Lagrangian CW-complex
\(\overline{\Delta}\subset B_{d,p,q}\) which we can think of
as a quotient of the unit disc by the equivalence relation
which identifies points \(z\sim e^{2\pi iq/p}z\) in its
boundary. Where \(\overline{\Delta}\) meets the toric
boundary, it does so along a visible Lagrangian
\((p,q)\)-pinwheel core\index{Lagrangian!pinwheel core}. We
call such a Lagrangian CW-complex a {\em Lagrangian
\((p,q)\)-pinwheel}.\index{Lagrangian!pinwheel} See Figure
\ref{fig:b_dpq_diagram_visibles} for the projections of these
visible Lagrangians.

\begin{figure}[htb]
\begin{center}
\begin{tikzpicture}
\filldraw[fill=lightgray,opacity=0.5,draw=none] (0,3) -- (0,0) -- (8,3) -- cycle;
\draw[->] (0,3) -- (0,0) -- (8,3) node [right] {\((dp^2,dpq-1)\)};
\draw[dotted,thick] (0,0) -- (2,1) node {\(\times\)} -- (3,1.5) node {\(\times\)} -- (4,2) node {\(\times\)};
\draw[thick] (0,0) -- (2,1) -- (3,1.5) -- (4,2);
\node at (1,0.5) [above] {\(\overline{\Delta}\)};
\node at (2.5,1.25) [above] {\(\overline{S}_1\)};
\node at (3.5,1.75) [above] {\(\overline{S}_2\)};
\end{tikzpicture}
\end{center}
\caption{The visible spheres and pinwheel in \(B_{d,p,q}\), shown in the case \(d=2,p=2,q=1\) (in this case, \(\overline{\Delta}\) is an \(\rp{2}\)).}
\label{fig:b_dpq_diagram_visibles}
\end{figure}
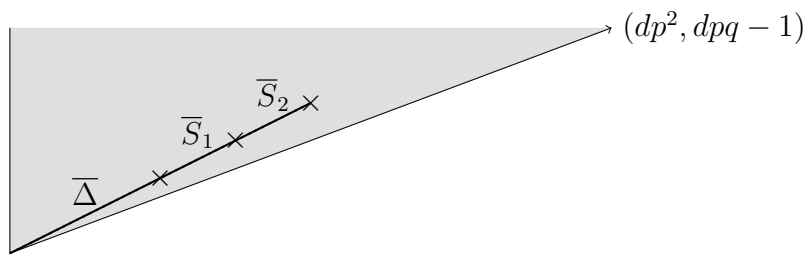

\end{Remark}
\begin{Lemma}\label{lma:lag_skel_b_dpq}
This union of visible Lagrangians
\(\overline{\Delta} \cup
\bigcup_{i=1}^{d-1} \overline{S}_i\) is homotopy equivalent to
\(B_{d,p,q}\). In particular, \(\pi_1(B_{d,p,q})=\ZZ/p\),
\(H_1(B_{d,p,q};\ZZ) = \ZZ/p\) and \(H_2(B_{d,p,q};\ZZ) =
\ZZ^{d-1}\).
\end{Lemma}
\begin{Proof}
The manifold \(B_{d,p,q}\)\index{Bdpq@$B_{d,p,q}$|)}
deformation retracts onto the preimage of the line segment
\(\ell\) shown in the fundamental action domain in Figure
\ref{fig:b_dpq_fad}.

\begin{figure}[htb]
\begin{center}
\begin{tikzpicture}
\filldraw[fill=lightgray,opacity=0.5,draw=none] (-1,-1) -- (-1,2) -- (4,2) -- (4,-1) -- cycle;
\draw (-1,-1) -- (-1,2);
\draw[dotted,thick] (4,2) -- (3,1.5) node {\(\times\)} -- (2,1) node {\(\times\)} -- (1,0.5) node {\(\times\)};
\draw[thick] (-1,-0.5) -- (3,1.5);
\draw [decorate,decoration={brace,amplitude=10pt},xshift=0pt,yshift=0pt] (-1,-0.5) -- (3,1.5) node [midway,above left=0.2cm] {\(\ell\)};
\node at (0,0) {\(\bullet\)};
\end{tikzpicture}
\caption{A fundamental action domain for \(B_{d,p,q}\). The manifold \(B_{d,p,q}\) deformation-retracts onto a Lagrangian CW-complex which projects to the line \(\ell\).}
\label{fig:b_dpq_fad}
\end{center}
\end{figure}
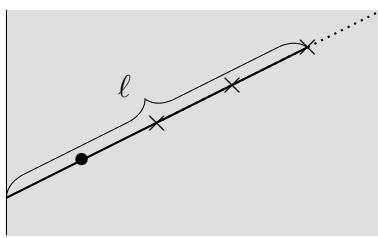

Let \(\ell_-\) and \(\ell_+\) be the segments of \(\ell\) to
the left and right (respectively) of the marked point
\(\bullet\) in the diagram. The preimage of \(\ell_-\) is a
solid torus \(T_-\); our convention will be that the loop
\((1,0)\) in \(\partial T_-\) bounds a disc in \(T_-\). The
preimage of \(\ell_+\) can be understood as follows. If there
were no base-nodes, it would be \(T^2\times[0,1]\). Each
base-node means that we pinch the torus above it along a loop
in the homology class \((-q,p)\). Up to homotopy equivalence,
this is the same as attaching a disc to \(T^2\times[0,1]\)
along a loop in this homology class. We can further homotope
these attaching maps so that they attach to loops in the
boundary of \(T_-\). Therefore the preimage of \(\ell\) is
homotopy equivalent to a solid torus with \(d\) discs attached
along its boundary along \(d\) parallel copies of the loop
\((-q,p)\). Now by a homotopy equivalence we can collapse
\(T_-\) to its core circle. The result is a CW-complex built
up from the core circle by adding \(d\) 2-cells using the
attaching map which winds the boundaries \(p\) times around
the core circle. Let \(\overline{\Delta}\) be the result of
attaching the first of these 2-cells to the core circle. Since
all these attaching maps are homotopic, we can homotope the
remaining \(d-1\) to a point in \(\overline{\Delta}\) and we
see that the resulting CW-complex is homotopy equivalent to
\(\overline{\Delta}\) wedged with \(d-1\) spheres. The
homology and fundamental group can be calculated using this
CW-decomposition.\index{almost toric base
diagram!Bdpq@$B_{d,p,q}$|)}\index{skeleton!Bdpq@$B_{d,p,q}$|)}\qedhere

\end{Proof}
\section{Solutions to inline exercises}

\begin{Exercise}\label{exr:auroux_hessian}
Given a Hamiltonian system \(\bm{H}\) and a critical point
\(\bm{x}\) of \(\bm{H}\), let \(Q(\bm{H},\bm{x})\) denote the
subspace of the space of quadratic forms spanned by the
Hessians\index{Hessian|(} of the components of
\(\bm{H}\). Check that, after a suitable symplectic change of
coordinates, \(Q(\bm{H},0)=Q(\bm{F},0)\) where \(\bm{H}\) is
the Auroux system and \(\bm{F}\) is the standard focus-focus
system from Example \ref{exm:ff}.
\end{Exercise}
\begin{Solution}
If we set \(z_k=x_k+iy_k\) then \(H_1(z_1,z_2)=|z_1z_2-c|^2 =
2c(y_1y_2-x_1x_2)+\cdots\), where the dots stand for terms of
higher order in the Taylor expansion, and
\(H_2(z_1,z_2)=\frac{1}{2}(|z_1|^2-|z_2|^2) =
\frac{1}{2}(x_1^2+y_1^2-x_2^2-y_2^2)\). Recall that the
standard focus-focus Hamiltonians are \[F_1 =
-p_1q_1-p_2q_2,\quad F_2 = p_2q_1-p_1q_2.\] If we make the
symplectic change of coordinates:
\begin{align*}
p_1&=\frac{1}{\sqrt{2}}(x_2-y_1)&
p_2&=\frac{1}{\sqrt{2}}(x_1-y_2)\\
q_1&=\frac{1}{\sqrt{2}}(x_1+y_2)&
q_2&=\frac{1}{\sqrt{2}}(x_2+y_1)
\end{align*}
then we get \[-(p_1q_1+p_2q_2)=y_1y_2-x_1x_2,\quad
p_2q_1-p_1q_2 = \frac{1}{2}(x_1^2-y_1^2-x_2^2+y_2^2)\] so
that, to second order, \(H_1(z_1,z_2) = 2cF_1(p,q)\) and
\(H_2(z_1,z_2)=F_2(p,q)\). Therefore the Hessians of \(H_1\)
and \(H_2\) span the same subspace of quadratic forms as the
Hessians\index{Hessian|)} of \(F_1\) and \(F_2\) after this
coordinate change.\qedhere

\end{Solution}
\chapter{Almost toric manifolds}
\label{ch:almost_toric_manifolds}
\thispagestyle{cup}

We have now seen some examples of Hamiltonian systems with
focus-focus critical points; in particular, we have seen what
their fundamental action domains look like. We now introduce a
definition ({\em almost toric fibrations}) which covers all of
these examples, and use it to develop some general theory for
manipulating and interpreting their fundamental action domains
({\em almost toric base diagrams}).

\section{Almost toric fibration}

\begin{Definition}
An {\em almost toric fibration}\index{fibration!almost
toric} is a Lagrangian torus fibration \(f\colon X\to B\) on
a 4-dimensional symplectic manifold such that the discriminant
locus comprises a collection of 0- and 1-dimensional strata
such that the smooth structure on \(B\) extends over these
strata, \(f\) is smooth with respect to this extended smooth
structure and has either toric or focus-focus critical points
there.

\end{Definition}
\begin{Remark}
We remark that the smooth structure mentioned in the
definition plays a somewhat auxiliary role: as in Remark
\ref{rmk:smooth_atlas}, the regular locus of the base inherits
a, possibly different, smooth structure from its canonical
integral affine structure, and this may not extend.

\end{Remark}
Let \(B^{reg}\subset B\) be the set of regular values of an
almost toric fibration \(f\), let \(\tilde{B}^{reg}\) be its
universal cover, and let \(\II\colon\tilde{B}^{reg}\to\RR^2\) be
the flux map. Let \(D\subset \tilde{B}^{reg}\) be a fundamental
domain for the action of \(\pi_1(B^{reg})\). Recall from Remark
\ref{rmk:base_node} that if \(b_1,b_2,\ldots\in B^{reg}\) is a
sequence tending to a focus-focus critical point then
\(\lim_{k\to\infty}\II(b_k)\) is a well-defined point in
\(\RR^2\) called the {\em base-node}\index{base-node} of that
critical point. There is also an affine monodromy associated to
loops in \(B^{reg}\) that go around the base-node. Lemma
\ref{lma:monodromy_formula} tells us that this monodromy is
completely determined by specifying, at each base-node, a
primitive integral eigenvector\index{eigenline} \((p,q)\) with
eigenvalue \(1\) for the monodromy.

\begin{Definition}
The {\em almost toric base diagram}\index{almost toric base
diagram} associated to these choices is the fundamental action
domain\index{action domain, fundamental}
\(\II(D)\subset\RR^2\) decorated with the positions of the
base-nodes and the eigenvector at each base node.

\end{Definition}
\begin{Remark}
We will usually (though not always) choose our fundamental
domain \(D\) by making branch cuts connecting the base nodes
to the boundary along eigenlines\index{eigenline!branch cut}.

\end{Remark}
Although the almost toric base diagram does not determine
\(f\colon X\to B\) up to {\em fibred} symplectomorphism, the
next theorem guarantees that it does determine \(X\) up to
symplectomorphism.

\begin{Theorem}[Symington {\cite[Corollary 5.4]{Symington}}]\label{thm:atbd}
Suppose\index{Symington's theorem!on almost toric base
diagrams} that \(f\colon X\to B\) and \(g\colon Y\to B\) are
almost toric fibrations whose almost toric base diagrams are
the same. If \(B\) is a punctured 2-dimensional surface then
\(X\) and \(Y\) are symplectomorphic.
\end{Theorem}
\begin{Proof}
Let \(N\subset B\) be the set of base-nodes. By Theorem
\ref{thm:symington}, there is a neighbourhood \(U\) of \(N\)
together with a symplectomorphism \(\Phi\colon f^{-1}(U)\to
g^{-1}(U)\). Although this symplectomorphism is not fibred, it
is fibred near the boundary of \(U\). Choose Lagrangian
sections of \(f\) and \(g\) over \(U\) (for example, we can
use the section \(\sigma_1(b)=(-\bar{b},1)\) in each
focus-focus chart, see the proof of Theorem
\ref{thm:ngoc}). If we can find global Lagrangian sections
over \(B\setminus N\) which match with the chosen Lagrangian
section over \(U\) along \(U\setminus N\) then Theorem
\ref{thm:uniqueness} gives us a fibred symplectomorphism
\(f^{-1}(B\setminus U)\to g^{-1}(B\setminus U)\) extending the
symplectomorphism \(\Phi\colon f^{-1}(U)\to g^{-1}(U)\). Since
\(B\setminus N\) has the homotopy type of a punctured surface
and \(N\) is a strict subset of the punctures (by assumption),
the relative cohomology group \(H^2(B\setminus N,U\setminus
N)\) vanishes. Therefore Corollary \ref{cor:extend_section_2}
tells us we can extend the Lagrangian section as
required. (This was the strategy we alluded to in Remark
\ref{rmk:rel_coh}.)\qedhere

\end{Proof}
\section{Operation I: nodal trade}

We now introduce some tools for manipulating and constructing
almost toric base diagrams which will give us a wealth of
examples. The first of these is Symington's {\em nodal
trade}\index{nodal trade|(}. It allows us to ``trade'' a Delzant
corner\index{Delzant!vertex|(} for a base-node.

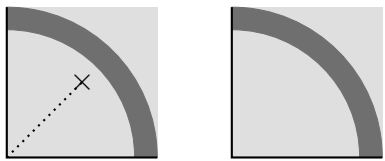
\begin{figure}[htb]
\begin{center}
\begin{tikzpicture}[baseline=1]
\filldraw[fill=lightgray,opacity=0.5,draw=none] (0,0) -- (0,2) -- (2,2) -- (2,0) -- cycle;
\draw[thick,black] (2,0) -- (0,0) -- (0,2);
\draw[dotted,thick] (0,0) -- (1,1);
\node at (1,1) {\(\times\)};
\filldraw[opacity=0.5] (1.7,0) arc [start angle=0,end angle=90,radius=1.7] to (0,2) arc [start angle=90,end angle=0,radius=2] to (1.7,0);
\end{tikzpicture}
\qquad
\begin{tikzpicture}[baseline=1]
\filldraw[fill=lightgray,opacity=0.5,draw=none] (0,0) -- (0,2) -- (2,2) -- (2,0) -- cycle;
\draw[thick,black] (2,0) -- (0,0) -- (0,2);
\filldraw[opacity=0.5] (1.7,0) arc [start angle=0,end angle=90,radius=1.7] to (0,2) arc [start angle=90,end angle=0,radius=2] to (1.7,0);
\end{tikzpicture}
\end{center}
\caption{(a) A fundamental action domain for the Auroux system on \(\CC^2\). (b) The moment image for the standard torus action on \(\CC^2\). The shaded subsets in both diagrams are integral affine equivalent to one another.}\label{fig:nodaltrade1}
\end{figure}

Recall from Figure \ref{fig:nodaltrade0} that there is an almost
toric structure on \(\CC^2\) which admits a fundamental action
domain as drawn on the left in Figure \ref{fig:nodaltrade1}(a)
below. The shaded region is integral affine equivalent to the
shaded region in Figure \ref{fig:nodaltrade1}(b), which is a
subset of the moment polygon for the standard torus action on
\(\CC^2\). This means that the preimages of these two regions
are fibred-symplectomorphic.

In particular, whenever we see a Delzant corner, we can excise
it and glue in a copy of the Auroux system, using this fibred
symplectomorphism to make identifications. Since the
identifications are fibred, this operation yields a new
Lagrangian torus fibration on the same manifold\footnote{To see
that the manifold does not change, observe that we are excising
a symplectic ball and gluing in another symplectic ball with the
same boundary (a contact 3-sphere\index{contact!-type
hypersurface}). The contactomorphism group of the 3-sphere is
connected \cite{Eliashberg}, so there is a unique way to glue up
to isotopy.}. In fact, there are many different operations, one
for each V\~{u} Ng\d{o}c model, but the results are all (non-fibred)
symplectomorphic to one another by Theorem
\ref{thm:symington}. We call an operation like this a {\em nodal
trade}.

\begin{figure}[htb]
\begin{center}
\includegraphics[scale=0.6]{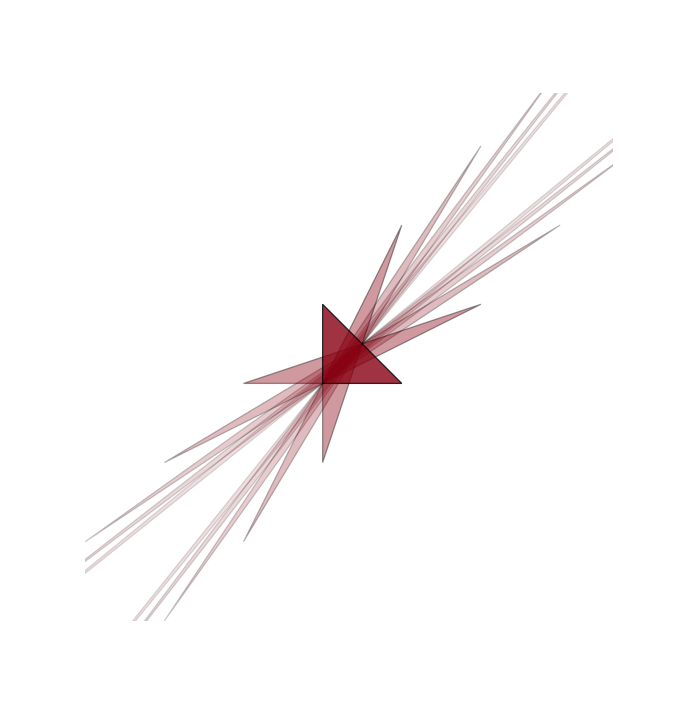}
\caption{The image of the developing map for an almost toric structure on \(\cp{2}\) obtained from the standard moment triangle by a single nodal trade in the lower left corner.}
\label{fig:dev1}
\end{center}
\end{figure}

\begin{figure}[htb]
\begin{center}
\includegraphics[scale=0.8]{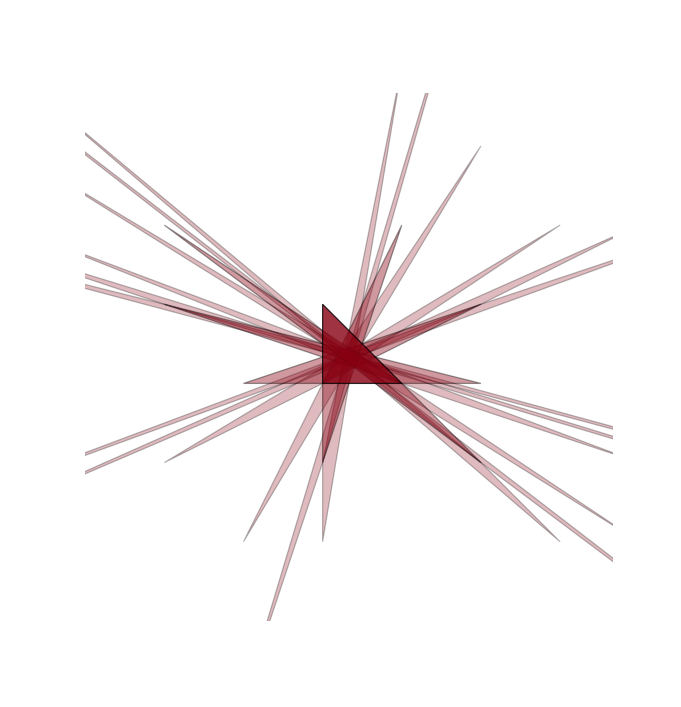}
\caption{The image of the developing map for an almost toric
structure on \(\cp{2}\) obtained from the standard moment triangle by
a nodal trade in the lower left corner and a nodal trade in the lower
right corner.}\label{fig:dev2}
\end{center}
\end{figure}

\begin{figure}[htb]
\begin{center}
\includegraphics[scale=0.8]{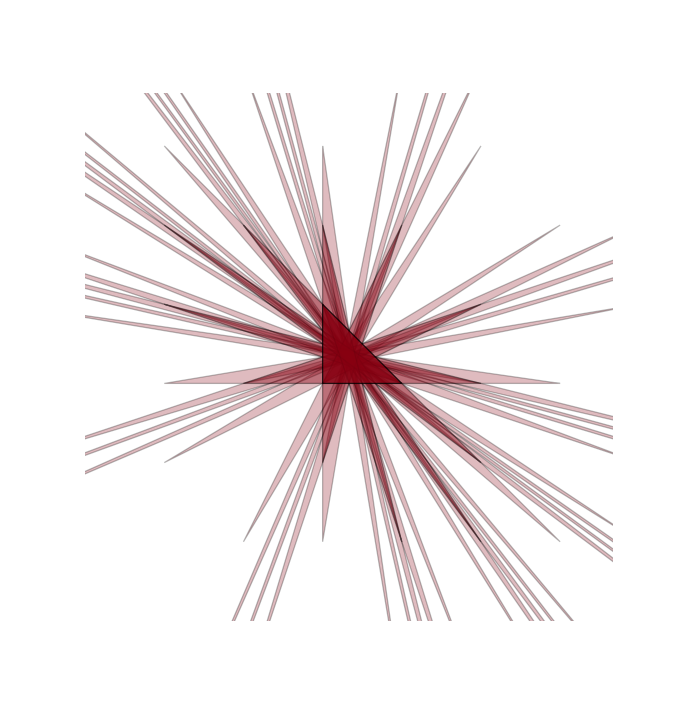}
\caption{The image of the developing map for an almost toric structure on \(\cp{2}\) obtained from the standard moment triangle by nodal trades in all three corners.}
\label{fig:dev3}
\end{center}
\end{figure}

\begin{Remark}
The toric boundary\index{toric boundary!after nodal trade}
near a Delzant corner comprises two symplectic discs meeting
transversally at the vertex. When you perform a nodal trade,
the toric boundary becomes a symplectic annulus which is a
smoothing\index{smoothing!and nodal trades|(} of this pair of
discs. For example, in the Auroux system this is the smoothing
from \(z_1z_2=0\) to \(z_1z_2=c\).

\end{Remark}
\begin{Example}\label{exm:atbd_cp2}
Here are some Lagrangian torus fibrations on
\(\cp{2}\):\index{almost toric base
diagram!CP2@$\mathbb{CP}^2$}

\begin{center}
\begin{tikzpicture}
\filldraw[fill=lightgray,opacity=0.5,draw=black,thick] (0,0) -- (0,3) -- (3,0) -- cycle;
\draw[thick] (0,0) -- (0,3) -- (3,0) -- cycle;
\draw[dotted,thick] (0,0) -- (1/2,1/2);
\node at (1/2,1/2) {\(\times\)};
\end{tikzpicture}
\qquad
\begin{tikzpicture}
\filldraw[fill=lightgray,opacity=0.5,draw=black,thick] (0,0) -- (0,3) -- (3,0) -- cycle;
\draw (0,0) -- (0,3) -- (3,0) -- cycle;
\draw[dotted,thick] (0,0) -- (1/2,1/2);
\draw[dotted,thick] (3,0) -- (2,0.5);
\node at (2,1/2) {\(\times\)};
\node at (1/2,1/2) {\(\times\)};
\end{tikzpicture}
\qquad
\begin{tikzpicture}
\filldraw[fill=lightgray,opacity=0.5,draw=black,thick] (0,0) -- (0,3) -- (3,0) -- cycle;
\draw (0,0) -- (0,3) -- (3,0) -- cycle;
\draw[dotted,thick] (0,0) -- (1/2,1/2);
\draw[dotted,thick] (0,3) -- (0.5,2);
\draw[dotted,thick] (3,0) -- (2,0.5);
\node at (1/2,2) {\(\times\)};
\node at (2,1/2) {\(\times\)};
\node at (1/2,1/2) {\(\times\)};
\end{tikzpicture}
\end{center}

The nodal trade in the lower left corner should look familiar;
we call this a {\em standard Delzant corner}. To find the
eigendirection for {\em any} Delzant corner \(p\), if \(A\) is
the unique integral affine transformation which maps the
standard Delzant corner to \(p\) then then the eigendirection
at \(p\) is \((1,1)A\). For example, the top left corner is
the image of the standard Delzant
corner\index{Delzant!vertex|)} under \(\lmatrix 0 & -1 \\ 1 &
-1\rmatrix\), so the eigendirection is \((1,-2)\), as
shown.\index{nodal trade|)}

As noted in Remark \ref{rmk:brokenlines}, although the toric
boundary\index{toric boundary!CP2 with nodal
trades@$\mathbb{CP}^2$ with nodal trades|(} looks like three
line segments, every time it crosses a branch cut you have to
apply the affine monodromy\index{toric boundary!appears
broken} to its tangent vector, so the apparent breaks in the
line when it crosses a branch cut are just an illusion: it is
really an uninterrupted straight line in the affine
structure. In the three examples above, the toric boundary
comprises:
\begin{itemize}
\item a conic\index{conic} and a line (two spheres intersecting
transversely at two points, one having twice the symplectic
area of the other),
\item a nodal cubic curve\index{cubic curve!nodal} (pinched
torus\index{pinched torus} having symplectic area three),
\item a smooth cubic curve\index{cubic curve!smooth} (torus having
symplectic area three).
\end{itemize}
This should make sense: the toric boundary for the usual toric
picture of \(\cp{2}\) comprises three lines and these
configurations above are obtained by
smoothing\index{smoothing!and nodal trades|)} one or more
intersections between these lines. Although I have used the
terminology ``line'', ``conic'', and ``cubic'' from algebraic
geometry, it is not clear for these new integrable Hamiltonian
systems whether the toric boundary is actually a subvariety
for the standard complex structure. It is, at least, a
symplectic submanifold\index{symplectic submanifold!of CP2@of
$\mathbb{CP}^2$} (immersed, where there are double points),
and it is known that low-degree symplectic surfaces in
\(\cp{2}\) are isotopic amongst symplectic surfaces to
subvarieties (see Gromov \cite{Gromov} for smooth surfaces of
degrees 1 and 2, Sikorav \cite{Sikorav} for smooth surfaces of
degree 3, Shevchishin \cite{ShevNodal} for nodal surfaces of
genus at most 4, and Siebert-Tian
\cite{SiebertTian1,SiebertTian2} for smooth surfaces in
degrees less than or equal to 17), hence the abuse of
terminology.\index{toric boundary!CP2 with nodal
trades@$\mathbb{CP}^2$ with nodal trades|)}

\end{Example}
The diagrams we have drawn rely on a specific choice of
fundamental domain \(D\) in the universal cover of the regular
locus of the almost toric fibration. If we simply plot the image
of the developing map\index{developing map!image} (flux map) on
the whole of \(\tilde{B}^{reg}\), we get some very beautiful
pictures which are symmetric under the action of
\(\pi_1(B^{reg})\) via affine monodromy. In Figures
\ref{fig:dev1}-\ref{fig:dev3} you can see what this looks like
for the three examples in Example \ref{exm:atbd_cp2}. These
pictures are closely related to the idea of {\em mutation} we
will meet in Section \ref{sct:mutation}, and the {\em Vianna
triangles} that we will see later (Theorem \ref{thm:vianna_tri}
and Appendix \ref{ch:markov_triples}).

\clearpage

\begin{Remark}\label{rmk:generalised_nodal_trade}
A closely related operation takes an almost toric base diagram
with a (non-Delzant)
corner\index{polytope!non-Delzant}\index{non-Delzant|see
{polytope, non-Delzant}} modelled on the polygon
\(\pi(dp^2,dpq-1)\) and replaces it with the diagram in Figure
\ref{fig:b_dpq_diagram}. We will call this a {\em generalised
nodal trade}\index{nodal trade!generalised|(}. One can think
of this as smoothing\index{smoothing!and generalised nodal
trades} the cyclic quotient
T-singularity\index{singularity!cyclic
quotient}\index{singularity!T-}\index{T-singularity|see
{singularity, T-}} in the original almost toric
orbifold\index{orbifold}.

\end{Remark}
\begin{Example}
Take the cubic surface\index{cubic surface} in \(\cp{3}\)
given in homogeneous coordinates \([z_1:z_2:z_3:z_4]\) by
\(z_1z_2z_3=z_4^3\). This has three
\(A_2\)-singularities\index{singularity!An@$A_n$} at
\([1:0:0:0]\), \([0:1:0:0]\) and \([0:0:1:0]\). It is toric,
with the torus action given by \((z_1,z_2,z_3,z_4)\mapsto
(e^{3i\theta_1}z_1, e^{3i\theta_2}z_2, z_3,
e^{i(\theta_1+\theta_2)}z_4)\). The Hamiltonians \(H_1 =
\frac{3|z_1|^2+|z_4|^2}{|z|^2}\) and
\(H_2=\frac{3|z_2|^2+|z_4|^2}{|z|^2}\) (with
\(|z|^2=\sum_{i=1}^4|z_i|^2\)) generate this action and their
image is the triangle \(\{(b_1,b_2)\in\RR^2\,:\,b_1\geq
0,b_2\geq 0,b_1+b_2\leq 3\}\). However, the period lattice is
not standard, for example the element \(\theta_1 = 2\pi/3\)
\(\theta_2 = 4\pi/3\) acts trivially. If we use the
Hamiltonians \(H_2\) and \((H_1+2H_2)/3\), whose period
lattice is standard, then the moment
polygon\index{moment polytope!singular cubic surface} becomes
the triangle with vertices at \((0,0)\), \((0,1)\) and
\((3,2)\). This has three corners each modelled on the polygon
\(\pi(3,2)\), corresponding to the cyclic quotient singularity
\(\frac{1}{3}(1,2)\) (also known as \(A_2\)). Performing three
generalised nodal trades\index{nodal trade!generalised|)}
gives the almost toric base diagram\index{almost toric base
diagram!smooth cubic surface} for the smooth cubic surface
shown in Figure \ref{fig:cubic_surface}.

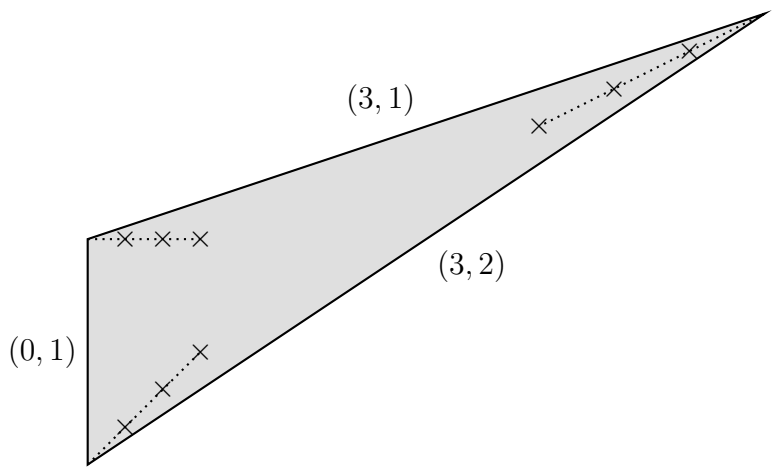
\begin{figure}[htb]
\begin{center}
\begin{tikzpicture}[baseline=1]
\filldraw[fill=lightgray,opacity=0.5,draw=black] (0,0) -- (0,3) -- (9,6) -- cycle;
\draw[thick] (0,0) -- (0,3) -- (9,6) -- cycle;
\draw[dotted,thick] (0,0) -- (0.5,0.5) node {\(\times\)} -- (1,1) node {\(\times\)} -- (1.5,1.5) node {\(\times\)};
\draw[dotted,thick] (0,3) -- (0.5,3) node {\(\times\)} -- (1,3) node {\(\times\)} -- (1.5,3) node {\(\times\)};
\draw[dotted,thick] (9,6) -- (8,5.5) node {\(\times\)} -- (7,5) node {\(\times\)} -- (6,4.5) node {\(\times\)};
\node at (4.5,3) [below right] {\((3,2)\)};
\node at (4.5,4.5) [above left] {\((3,1)\)};
\node at (0,1.5) [left] {\((0,1)\)};
\end{tikzpicture}
\end{center}
\caption{An almost toric base diagram for a cubic surface; the edges are labelled with primitive integer vectors pointing along them.}
\label{fig:cubic_surface}
\end{figure}

\end{Example}
\section{Operation II: nodal slide}

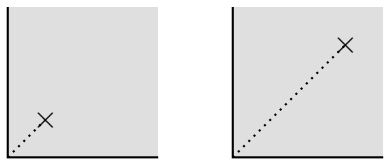
\begin{figure}[htb]
\begin{center}
\begin{tikzpicture}[baseline=1]
\filldraw[fill=lightgray,opacity=0.5,draw=none] (0,0) -- (0,2) -- (2,2) -- (2,0) -- cycle;
\draw[thick,black] (2,0) -- (0,0) -- (0,2);
\draw[dotted,thick] (0,0) -- (0.5,0.5);
\node at (0.5,0.5) {\(\times\)};
\end{tikzpicture}
\qquad
\begin{tikzpicture}[baseline=1]
\filldraw[fill=lightgray,opacity=0.5,draw=none] (0,0) -- (0,2) -- (2,2) -- (2,0) -- cycle;
\draw[thick,black] (2,0) -- (0,0) -- (0,2);
\draw[dotted,thick] (0,0) -- (1.5,1.5);
\node at (1.5,1.5) {\(\times\)};
\end{tikzpicture}
\end{center}
\caption{A nodal slide.}\label{fig:nodalslide}
\end{figure}

Note that there is a free parameter \(c>0\) in the Auroux
system. As this parameter varies, we obtain a family of
Lagrangian torus fibrations in which the focus-focus critical
point moves in the direction of the
eigenvector\index{eigenline!sliding along|see {nodal, slide}}
for its affine monodromy (see Figure \ref{fig:nodalslide}). Such
a family of fibrations is called a {\em nodal
slide}\index{nodal slide|(}.

\begin{Theorem}[Symington's theorem on nodal slides {\cite[Theorem 6.5]{Symington}}]\label{thm:slide_symplecto}
Suppose\index{Symington's theorem!on nodal slides} that \(X\)
and \(X'\) are almost toric manifolds whose almost toric base
diagrams \(B\) and \(B'\) are related by a nodal
slide. Suppose that the base diagrams have the homotopy type
of a punctured 2-dimensional surface. Then \(X\) and \(X'\)
are symplectomorphic.
\end{Theorem}
\begin{Proof}
This argument is a Moser argument\index{Moser argument} very
similar to the one used to prove Theorem \ref{thm:symington}.

Let \(b\in B\) be the base node which slides and let \(b'\)
denote the corresponding base node in \(B'\). For simplicity,
we consider only the case where \(B\) and \(B'\) are obtained
by taking branch cuts\index{branch cut!nodal slide} along an
eigenray for \(b\). There is a family of almost toric base
diagrams \(B_t\) interpolating between \(B\) and \(B'\). These
diagrams all coincide outside a contractible neighbourhood
\(K\) of the sliding line.

\begin{center}
\begin{tikzpicture}
\filldraw[fill=lightgray,opacity=0.5] (0,0) -- (1,1) -- (5,1) -- (5,-1) -- (1,-4) -- (0,-4) -- cycle;
\draw[thick] (0,0) -- (1,1) -- (5,1) -- (5,-1) -- (1,-4) -- (0,-4) -- cycle;
\filldraw[fill=gray,opacity=0.5] (2,0) circle [x radius = 1.5, y radius = 0.5];
\node at (2,0) {\(K\)};
\draw[dotted,thick] (0,0) -- (1,0) node {\(\times\)};
\draw[dotted,thick] (5,-1) -- (3,-2) node {\(\times\)};
\draw[dotted,thick] (0,-4) -- (1,-3) node {\(\times\)};
\draw[dotted,thick] (5,1) -- (4,0) node {\(\times\)};
\begin{scope}[shift={(6,0)}]
\filldraw[fill=lightgray,opacity=0.5] (0,0) -- (1,1) -- (5,1) -- (5,-1) -- (1,-4) -- (0,-4) -- cycle;
\draw[thick] (0,0) -- (1,1) -- (5,1) -- (5,-1) -- (1,-4) -- (0,-4) -- cycle;
\filldraw[fill=gray,opacity=0.5] (2,0) circle [x radius = 1.5, y radius = 0.5];
\draw[dotted,thick] (0,0) -- (3,0) node {\(\times\)};
\draw[dotted,thick] (5,-1) -- (3,-2) node {\(\times\)};
\draw[dotted,thick] (0,-4) -- (1,-3) node {\(\times\)};
\draw[dotted,thick] (5,1) -- (4,0) node {\(\times\)};
\end{scope}
\end{tikzpicture}
\end{center}

This gives a family of almost toric manifolds \(X_t\) together
with almost toric fibrations \(f_t\colon X_t\to B_t\). Since
the diagrams \(B_0\) and \(B\) (respectively \(B_1\) and
\(B'\)) are identical, \(X_0\) and \(X\) (respectively \(X_1\)
and \(X'\)) are symplectomorphic by Theorem
\ref{thm:atbd}. (Here is the first place where we use the
assumption on the topology of the base).

It suffices to show that \(X_0\) and \(X_1\) are
symplectomorphic. We use the fact that the almost toric bases
coincide outside \(K\) to identify the subsets
\(U_t=f_t^{-1}(B_t\setminus K)\), extend this to a family of
diffeomorphisms \(\varphi_t\colon X_0\to X_t\), and show that
the symplectic forms \(\omega_t=\varphi_t^*\omega_{X_t}\) on
\(X_0\) satisfy \(d[\omega_t]/dt=0\). The result will then
follow from Moser's argument (Appendix \ref{ch:moser}). The
subset on which the symplectic forms differ is
\(f_t^{-1}(K)\), which is just a neighbourhood of the nodal
fibre. As we saw in the proof of Theorem \ref{thm:symington},
it suffices to show that \(\int_{\sigma_t}d\omega_t/dt=0\) for
some family of submanifolds \(\sigma_t\) such that
\(\varphi_t(\sigma_t)\) intersects the nodal fibre of \(f_t\)
once transversely and \(\sigma_t=\sigma_0\) on
\(\varphi_t^{-1}(U_t)\). As in that proof, we can use a chosen
family of Lagrangian sections (here again we use the topology
of the base to guarantee the existence of these
sections). \qedhere

\end{Proof}
As a result, nodal sliding does not change the symplectic
manifold, but it certainly changes the Lagrangian
fibration. Here is an example.

\begin{Example}\label{exm:vianna_1}
Start with the almost toric fibration\index{almost toric base
diagram!CP2@$\mathbb{CP}^2$|(} on \(\cp{2}\) from Example
\ref{exm:atbd_cp2} with three base nodes. Pick the
bottom-right node (labelled \(B\) in Figure
\ref{fig:distinguishing_ltfs}, and slide it towards the
opposite edge, beyond the barycentre of the triangle. We get
two almost toric fibrations on \(\cp{2}\) which we can
distinguish as follows. Consider the Lagrangian vanishing
thimble\index{vanishing!thimble} emanating from the base-node
labelled \(A\): this is a visible Lagrangian disc with centre
at the focus-focus critical point, and its projection is a ray
in the base diagram pointing in the eigendirection of this
node, which is \((1,1)\).

Before the nodal slide, this ray hits the slanted
edge. Exercise \ref{exr:which_visible_lagrangians_vianna}: The
visible Lagrangian\index{Lagrangian!RP2@$\mathbb{RP}^2$} meets
this edge along a \((2,1)\)-pinwheel
core\index{Lagrangian!pinwheel core}, i.e.\ a M\"{o}bius
strip\footnote{Indeed, this visible Lagrangian is just
\(\rp{2}\subset\cp{2}\)}. After the nodal
slide\index{nodal slide|)}, this ray crosses the sliding
branch cut; when it emerges its direction has changed by the
affine monodromy \(\lmatrix -1 & 1 \\ -4 & 3\rmatrix \) so
that the ray now points in the \((-5,4)\)-direction. Exercise
\ref{exr:which_visible_lagrangians_vianna}: This now
intersects the slanted edge of the base diagram in a
\((5,4)\)-pinwheel core. This difference in topology
distinguishes the torus fibrations.

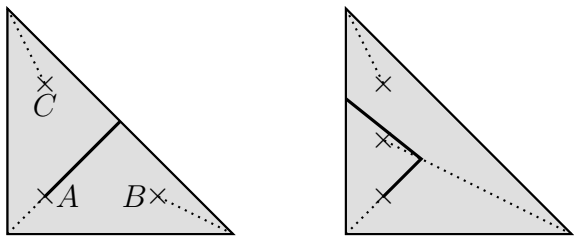
\begin{figure}[htb]
\begin{center}
\begin{tikzpicture}
\filldraw[fill=lightgray,opacity=0.5,draw=black,thick] (0,0) -- (0,3) -- (3,0) -- cycle;
\draw[thick] (0,0) -- (0,3) -- (3,0) -- cycle;
\draw[very thick] (1/2,1/2) -- (1.5,1.5);
\draw[dotted,thick] (0,0) -- (1/2,1/2);
\draw[dotted,thick] (0,3) -- (0.5,2);
\draw[dotted,thick] (3,0) -- (2,0.5);
\node at (1/2,2) {\(\times\)};
\node at (2,1/2) {\(\times\)};
\node at (1/2,1/2) {\(\times\)};
\node at (1/2,1/2) [right] {\(A\)};
\node at (2,1/2) [left] {\(B\)};
\node at (1/2,2) [below] {\(C\)};
\begin{scope}[shift={(4.5,0)}]
\filldraw[fill=lightgray,opacity=0.5,draw=black,thick] (0,0) -- (0,3) -- (3,0) -- cycle;
\draw[thick] (0,0) -- (0,3) -- (3,0) -- cycle;
\draw[very thick] (1/2,1/2) -- (1,1) -- (0,9/5);
\draw[dotted,thick] (0,0) -- (1/2,1/2);
\draw[dotted,thick] (0,3) -- (0.5,2);
\draw[dotted,thick] (3,0) -- (1/2,5/4);
\node at (1/2,2) {\(\times\)};
\node at (1/2,5/4) {\(\times\)};
\node at (1/2,1/2) {\(\times\)};
\end{scope}
\end{tikzpicture}
\end{center}
\caption{Left: Before the nodal slide, the vanishing thimble of base-node \(A\) is part of a visible Lagrangian \(\rp{2}\). Right: After the nodal slide, it is part of a visible Lagrangian \((5,4)\)-pinwheel.}
\label{fig:distinguishing_ltfs}
\end{figure}

\end{Example}
\begin{Remark}
Visible Lagrangians obtained by capping a \((p,q)\)-pinwheel
core\index{Lagrangian!pinwheel core} with a disc are quite
common in this context, and we call them {\em Lagrangian
\((p,q)\)-pinwheels}\index{Lagrangian!pinwheel}.

\end{Remark}
The curious reader may wonder what happens if we try and slide a
node in a direction which is not its eigenline\index{nodal
slide!transverse to eigenline}. If we do so, we obtain a family
of almost toric base diagrams, and so a family of symplectic
manifolds. However, the cohomology class of the symplectic form
can vary\index{symplectic form!cohomology
class|(}\index{symplectic form!deformation of}, so they are not
all symplectomorphic. Here is an example.

\begin{Example}
Take the toric diagram for \(\mathcal{O}(-1)\) giving the
compact edge affine length \(2\) and make two nodal trades.

\begin{center}
\begin{tikzpicture}
\filldraw[fill=lightgray,opacity=0.5] (0,2) -- (0,0) -- (4,0) -- (6,2);
\draw[thick] (0,2) -- (0,0) -- (4,0) -- (6,2);
\draw[dotted,thick] (0,0) -- (0.5,0.5) node {\(\times\)};
\draw[dotted,thick] (4,0) -- (4,0.5) node {\(\times\)};

\end{tikzpicture}
\end{center}
Now we attempt to slide the right-hand node to the left. As it
moves, more and more of what used to be the compact edge
passes through the branch cut and ends up as part of the
slanted non-compact edge on the right.

\begin{center}
\begin{tikzpicture}
\filldraw[fill=lightgray,opacity=0.3] (0,2) -- (0,0) -- (4,0) -- (6,2);
\filldraw[fill=lightgray,opacity=0.5] (0,2) -- (0,0) -- (2,0) -- (4,2);
\draw[thick] (0,2) -- (0,0) -- (2,0) -- (4,2);
\draw[dotted,thick] (0,0) -- (0.5,0.5) node {\(\times\)};
\draw[dotted,thick] (2,0) -- (2,0.5) node {\(\times\)};
\draw[dotted,thick,opacity=0.3] (4,0) -- (4,0.5) node {\(\times\)};

\end{tikzpicture}
\end{center}
If we ``undo'' the nodal trades, we see that the compact edge
has shrunk. In other words, the symplectic
area\index{symplectic area} of the zero-section has decreased,
and the cohomology class\index{symplectic form!cohomology
class|)} of \(\omega\) has changed.

\end{Example}
Indeed, one can use this to ``flip'' the sign of the symplectic
area of curves. This is closely related to the theory of flips
in algebraic geometry; see \cite{EvansUrzua} for a full
discussion with examples.

\section{Operation III: Mutation}
\label{sct:mutation}

In Example \ref{exm:vianna_1}, we needed to keep track of a
visible Lagrangian which crossed a branch cut. This can get very
tricky if there are several branch cuts. Sometimes, it is more
convenient to change the choice of fundamental domain in
\(\tilde{B}^{reg}\). We will do this by rotating the branch
cut\index{branch cut!mutation|(} as we did in Section
\ref{sct:aurouxnodaltrade} and Figure \ref{fig:branchcuts}.

In fact, in most examples, we will start with a branch cut which
emanates from a base-node in the direction \(v\) of an
eigenvector for the affine monodromy\index{affine monodromy!and
mutation|(}, and rotate by 180 degrees anticlockwise (or
clockwise) to get a branch cut in the \(-v\)-direction.

The fundamental action domain (almost toric base diagram)
transforms in the following way. The eigenline\index{eigenline}
bisects the diagram; call the two pieces \(D_1\) and \(D_2\) (we
adopt the convention that \(D_1\) lies clockwise of the branch
cut and \(D_2\) lies anticlockwise). Let \(M\) be the affine
monodromy around the base-node as we cross the branch cut in the
anticlockwise direction. When the branch cut rotates
anticlockwise 180 degrees, we replace \(D_2\) with
\((D_2)M^{-1}\) to get the new almost toric base diagram
\(D_1\cup (D_2)M^{-1}\). When the branch cut rotates clockwise
180 degrees, we get \((D_1)M\cup D_2\) instead.

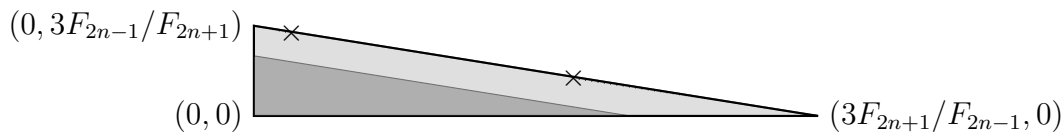
\begin{figure}[htb]
\begin{center}
\begin{tikzpicture}
\filldraw[fill=lightgray,opacity=0.5,draw=black,thick] (0,0) -- (0,6/5) -- (15/2,0) -- cycle;
\draw[thick] (0,0) -- (0,6/5) -- (15/2,0) -- cycle;
\filldraw[fill=gray,opacity=0.5] (0,0) -- (0,12/15) -- (30/6,0) -- cycle;
\draw[dotted] (0,6/5) -- (1/2,11/10) node {\(\times\)};
\draw[dotted] (15/2,0) -- (15/2-13/4,1/2) node {\(\times\)};
\node at (15/2,0) [right] {\((3F_{2n+1}/F_{2n-1},0)\)};
\node at (0,6/5) [left] {\((0,3F_{2n-1}/F_{2n+1})\)};
\node at (0,0) [left] {\((0,0)\)};
\end{tikzpicture}
\end{center}
\caption{The general almost toric base diagram in the Fibonacci staircase. The shaded region is a solid ellipsoid which can be rescaled to fill an arbitrarily large fraction of the volume of \(\cp{2}\) provided we slide the base-nodes closer to the corners.}
\label{fig:fibonacci_staircase}
\end{figure}

\begin{Remark}
Note that changing branch cut has no effect on the symplectic
manifold nor on the Lagrangian torus fibration. All that
changes is the picture: the picture is the image of a
fundamental action domain under the developing map, and the
change of branch cut amounts to a different choice of
fundamental action domain.

\end{Remark}
\begin{Example}\label{exm:atbd_cp_2_mutations}
Take\index{mutation!CP2 example@$\mathbb{CP}^2$ example|(} the
almost toric base diagram shown below and let \(x\) be the
base-node marked \(B\). The anticlockwise affine
monodromy\index{affine monodromy!and mutation|)} is
\(M=\lmatrix -1 & 1 \\ -4 & 3\rmatrix \) and the branch cut
points in the \((2,-1)\)-direction. We have indicated the
coordinates of the corners of the triangle. The affine lengths
of the three edges are \(3\); this choice corresponds to the
pullback of the Fubini-Study form along the anticanonical
embedding of \(\cp{2}\) in \(\cp{9}\).

\begin{center}
\begin{tikzpicture}
\filldraw[fill=lightgray,opacity=0.5,draw=black,thick] (0,0) -- (0,3) -- (3,0) -- cycle;
\draw[thick] (0,0) -- (0,3) -- (3,0) -- cycle;
\draw[dotted,thick] (3,0) -- (1/2,5/4);
\node at (1/2,5/4) {\(\times\)};
\node at (1/2,5/4) [left] {\(B\)};
\node at (3,0) [right] {\((3,0)\)};
\node at (0,3) [left] {\((0,3)\)};
\node at (0,0) [left] {\((0,0)\)};
\end{tikzpicture}
\end{center}

If we rotate the associated branch cut 180 degrees
anticlockwise, then the result is:

\begin{center}
\begin{tikzpicture}
\filldraw[fill=lightgray,opacity=0.5,draw=black,thick] (0,0) -- (0,3/2) -- (6,0) -- cycle;
\draw[thick] (0,0) -- (0,3/2) -- (6,0) -- cycle;
\draw[dotted,thick] (0,3/2) -- (1/2,5/4);
\node at (1/2,5/4) {\(\times\)};
\node at (1/2,5/4) [below] {\(B\)};
\node at (6,0) [right] {\((6,0)\)};
\node at (0,3/2) [left] {\((0,3/2)\)};
\node at (0,0) [left] {\((0,0)\)};
\end{tikzpicture}
\end{center}

We superimpose the two pictures for easier comparison.

\begin{center}
\begin{tikzpicture}
\filldraw[fill=lightgray,opacity=0.5,draw=black,very thick] (0,0) -- (0,3) -- (3,0) -- cycle;
\filldraw[fill=lightgray,opacity=0.5,draw=black,very thick] (0,0) -- (0,3/2) -- (6,0) -- cycle;
\draw[dotted,thick] (3,0) -- (1/2,5/4);
\node at (1/2,5/4) {\(\times\)};
\draw[dotted,thick] (0,3/2) -- (1/2,5/4);
\node at (1/2,1/2) {\(D_1\)};
\node at (1/2,2) {\(D_2\)};
\node at (4,1) {\((D_2)M^{-1}\)};
\end{tikzpicture}
\end{center}

After performing this change of branch cut, the visible
Lagrangian \((5,4)\)-pinwheel\index{Lagrangian!pinwheel|(} from
Example \ref{exm:vianna_1} is easier to see, as its projection
does not cross any branch cuts (see Figure
\ref{fig:vianna_after_mutation}).

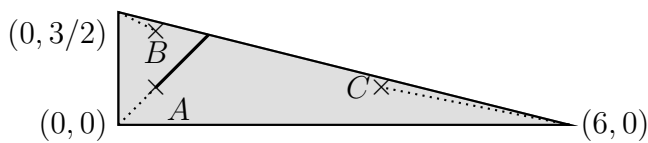
\begin{figure}[htb]
\begin{center}
\begin{tikzpicture}
\filldraw[fill=lightgray,opacity=0.5,draw=black,thick] (0,0) -- (0,3/2) -- (6,0) -- cycle;
\draw[thick] (0,0) -- (0,3/2) -- (6,0) -- cycle;
\draw[dotted,thick] (0,0) -- (1/2,1/2) node {\(\times\)};
\draw[dotted,thick] (0,3/2) -- (1/2,5/4) node {\(\times\)};
\draw[dotted,thick] (6,0) -- (6-5/2,1/2) node {\(\times\)};
\draw[very thick] (1/2,1/2) -- (6/5,6/5);
\node at (1/2,1/2) [below right] {\(A\)};
\node at (1/2,5/4) [below] {\(B\)};
\node at (6-5/2,1/2) [left] {\(C\)};
\node at (6,0) [right] {\((6,0)\)};
\node at (0,6/5) [left] {\((0,3/2)\)};
\node at (0,0) [left] {\((0,0)\)};
\end{tikzpicture}
\end{center}
\caption{The visible \((5,4)\)-pinwheel is easier to see after applying a mutation to Figure \ref{fig:distinguishing_ltfs}. Note that the direction of the branch cut at \(C\) is obtained from the original cut by applying the clockwise \(B\)-monodromy. These cuts are hard to see, as they are almost parallel to the edges. The vectors indicate the coordinates of the vertices.}
\label{fig:vianna_after_mutation}
\end{figure}

\end{Example}
\begin{Definition}
A (clockwise or anticlockwise) {\em
mutation}\index{mutation!definition} at a base-node is a
combination of a nodal slide followed by a change of the same
branch cut by 180 degrees (clockwise or anticlockwise).

\end{Definition}
\begin{Remark}
I have used the word ``mutation'' in earlier papers to mean
just a change of branch cut, e.g.\ \cite{EvansUrzua}. I am
changing my preference here to avoid overusing the phrase
``combination of nodal slide and mutation''.

\end{Remark}
\begin{Example}
Let us continue by mutating anticlockwise at the base-node
labelled \(C\) in Figure \ref{fig:vianna_after_mutation}. We
now find a visible Lagrangian \((13,2)\)-pinwheel living over
the eigenray emanating from the base-node \(C\).

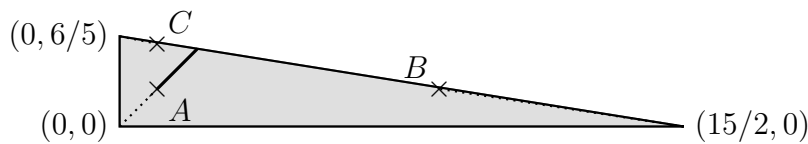
\begin{figure}[htb]
\begin{center}
\begin{tikzpicture}
\filldraw[fill=lightgray,opacity=0.5,draw=black,thick] (0,0) -- (0,6/5) -- (15/2,0) -- cycle;
\draw[thick] (0,0) -- (0,6/5) -- (15/2,0) -- cycle;
\draw[dotted,thick] (0,0) -- (1/2,1/2) node {\(\times\)};
\draw[dotted,thick] (0,6/5) -- (1/2,11/10) node {\(\times\)};
\draw[dotted,thick] (15/2,0) -- (15/2-13/4,1/2) node {\(\times\)};
\draw[very thick] (1/2,1/2) -- (6*15/87,6*15/87);
\node at (1/2,1/2) [below right] {\(A\)};
\node at (1/2,11/10) [above right] {\(C\)};
\node at (15/2-13/4,1/2) [above left] {\(B\)};
\node at (15/2,0) [right] {\((15/2,0)\)};
\node at (0,6/5) [left] {\((0,6/5)\)};
\node at (0,0) [left] {\((0,0)\)};
\end{tikzpicture}
\end{center}
\caption{After a further mutation, we see a visible \((13,2)\)-pinwheel. By now, the \(B\) and \(C\) branch cuts are indistinguishably close to the edges.}
\label{fig:vianna_after_further_mutation}
\end{figure}

\end{Example}
\begin{Remark}
One useful trick for figuring out the directions of the branch
cuts is the following. In the original \(\cp{2}\) triangle,
the eigenlines\index{eigenline!meeting at common point|(} all
intersect at the barycentre of the triangle \((1,1)\). This
remains true after mutation. So the directions of the branch
cuts can be found by taking the vectors from \((1,1)\) to the
corners. For example, the \(B\) branch cut in Figure
\ref{fig:vianna_after_further_mutation} points in the
\((13,-2)\)-direction. This trick works when all the
eigenlines intersect at a common point. This is a
cohomological condition; see Section
\ref{sct:coh_class}.\index{eigenline!meeting at common
point|)}

\end{Remark}
One could continue in this fashion, mutating at \(B\), then
\(C\), then \(B\), etc. This would give a sequence of Lagrangian
\((p_k,q_k)\)-pinwheels in \(\cp{2}\) where \(p_k\) runs over
the odd-indexed Fibonacci\index{Fibonacci|see {ellipsoid,
Fibonacci staircase}} numbers\footnote{Our convention is that
\(F_1 = 1\), \(F_3 = 2\), \(F_5=5\), etc.}. The general almost
toric base diagram\index{almost toric base diagram!Fibonacci
staircase} in this sequence is shown in Figure
\ref{fig:fibonacci_staircase}; to declutter the diagram we have
``undone'' the nodal trade at the vertex which remains
unmutated. The shaded region in this figure is a symplectically
embedded copy of the solid ellipsoid\index{ellipsoid!solid}
\[\left\{(z_1, z_2)\in \CC^2 \,:\, \frac{1}{2}\left(\frac{F_{2n-
1}}{3F_{2n+ 1}} |z_1|^2+ \frac{F_{2n+1}}{3F_{2n-
1}}|z_2|^2\right)\leq \lambda\right\}\] for \(\lambda < 1\). We
can get an embedding with \(\lambda\) arbitrarily close to \(1\)
by nodally sliding the two base-nodes very close to their
vertices and pushing the slanted edge of the shaded region
towards the toric boundary. In other words, we can fill up an
arbitrarily large fraction of the volume of \(\cp{2}\) by a
symplectically embedded ellipsoid\index{ellipsoid!Fibonacci
staircase}\index{ellipsoid!embedding
problems}\index{symplectic embedding of ellipsoid} with this
Fibonacci-ratio of radii. This is related to the famous {\em
Fibonacci staircase} pattern for full ellipsoid embeddings in
\(\cp{2}\) observed by McDuff and Schlenk
\cite{McDuffSchlenk}. For papers which discuss this and other
examples from an almost toric point of view, see
\cite{BertozziEtAl,CasalsVianna,CGHMP,MagillMcDuff}.

In fact, one could also allow sequences of mutations involving
base-node \(A\), to obtain an infinite trivalent tree of almost
toric base diagrams all representing almost toric fibrations on
\(\cp{2}\). This tree is closely related to another famous
infinite trivalent tree arising from the theory of Diophantine
approximation: the {\em Markov tree}\index{Markov tree}. For the
\((p,q)\)-pinwheels which appear, \(p\) is a Markov
number\index{Markov number}\footnote{The odd-indexed Fibonacci
numbers are a subset of the Markov numbers.}.

\begin{Definition}
A {\em Markov triple}\index{Markov triple} is a positive
integer solution \(p_1,p_2,p_3\) to the {\em Markov equation}
\[p_1^2+p_2^2+p_3^2 = 3p_1p_2p_3.\] A Markov number is a
number which appears in a Markov triple.

\end{Definition}
\begin{Theorem}[Vianna \cite{Vianna2}]\label{thm:vianna_tri}
For every Markov triple \(p_1,p_2,p_3\), there is an almost
toric diagram\index{almost toric base
diagram!CP2@$\mathbb{CP}^2$|)} \(D(p_1,p_2,p_3)\) (a {\em
Vianna triangle}) with the following
properties.\index{Vianna!triangle}
\begin{itemize}
\item The diagram \(D(1,1,1)\) is obtained from the standard toric
diagram of \(\cp{2}\) by performing three nodal trades.
\item The diagram \(D(p_1,p_2,p_3)\) is a triangle with three
base-nodes \(n_1,n_2,n_3\), obtained by iterated mutation on
\(D(1,1,1)\) (in particular, the associated almost toric
manifold is \(\cp{2}\)).
\item For \(k=1,2,3\), there is an integer \(q_k\) and a
Lagrangian pinwheel of type \((p_k,q_k)\) living over the
branch cut which connects \(n_k\) to a corner \(P_k\).
\item The affine length\index{affine length!Vianna triangles} of
the edge opposite the corner \(P_k\) is
\(3p_k/(p_{k+1}p_{k+2})\) where indices are taken modulo
\(3\).
\end{itemize}
\end{Theorem}
\begin{Proof}
We prove this in Appendix \ref{ch:markov_triples}, where we
also remind the reader of the basic properties of Markov
triples.\qedhere

\end{Proof}
\begin{Remark}
The diagrams in Figures \ref{fig:vianna_after_mutation} and
\ref{fig:vianna_after_further_mutation} are \(D(1,1,2)\) and
\(D(1,2,5)\).

\end{Remark}
\begin{Remark}
Superimposing all these almost toric base diagrams yields
Figure \ref{fig:dev3}, that is the image of the developing
map of the integral affine
structure on the universal
cover of the complement of the base-nodes.

\end{Remark}
\begin{Remark}
For each diagram \(D(p_1,p_2,p_3)\), let \(T(p_1,p_2,p_3)\)
denote\footnote{Vianna \cite{Vianna2}
uses the notation \(T(p_1^2,p_2^2,p_3^2)\).} the Lagrangian
torus fibre over the
barycentre\index{Vianna!tori}. Before the work of Vianna, the
tori \(T(1,1,1)\) (the {\em Clifford torus})\index{Clifford
torus} and \(T(1,1,2)\) (the {\em Chekanov
torus})\index{Chekanov torus} had been constructed and
Chekanov \cite{Chekanov} had shown that they were not
Hamiltonian isotopic. Vianna's truly remarkable contribution,
besides constructing \(T(p_1,p_2,p_3)\), was to show that if
\(p_1,p_2,p_3\) and \(p'_1,p'_2,p'_3\) are distinct as
unordered Markov triples then \(T(p_1,p_2,p_3)\) and
\(T(p'_1,p'_2,p'_3)\) are not Hamiltonian isotopic\footnote{He
had an earlier paper \cite{Vianna1} which distinguished
\(T(1,2,5)\) from \(T(1,1,1)\) and \(T(1,2,5)\), which was
already a major breakthrough. At the time it appeared, I was
convinced there should be only two Hamiltonian isotopy classes
of monotone\index{monotone!Lagrangian torus} Lagrangian tori
in \(\cp{2}\).}.\index{mutation!CP2 example@$\mathbb{CP}^2$
example|)}

\end{Remark}
\begin{Remark}
Vianna \cite{Vianna3} has also studied mutations of other
triangular almost toric base diagrams, and used this to
construct exotic Lagrangian tori in other symplectic
4-manifolds.\index{branch cut!mutation|)}

\end{Remark}
\begin{Remark}\label{rmk:open_mutation}
As far as I know, it is a completely open question to
characterise which quadrilaterals arise as mutations of a
square or a rectangle.\index{mutation!quadrilateral}

\end{Remark}
\section{More examples}

\begin{Example}\label{exm:cp_2_rp_2}
In Example \ref{exm:atbd_cp_2_mutations}, we constructed the
almost toric base diagram\index{almost toric base
diagram!CP2@$\mathbb{CP}^2$} for \(\cp{2}\) shown in Figure
\ref{fig:atbd_cp_2_rp_2} by performing one nodal trade and a
mutation (we have nodally-slid the base node to make it
clearer). The visible Lagrangian
pinwheel\index{Lagrangian!pinwheel|)} over the branch cut is
\(\rp{2}\) and the symplectic sphere living over the opposite
edge is isotopic to a conic curve (its symplectic area is
twice that of a complex line in \(\cp{2}\), so it inhabits the
homology class of a conic, and symplectic curves in this class
are known to be isotopic \cite{Gromov}). If we excise the
conic, what is left is a subset of the almost toric base
diagram for \(B_{1,2,1}\) from Example \ref{exm:bpq}. This is
actually the cotangent\index{cotangent bundle!of RP2@of
$\mathbb{RP}^2$} bundle of \(T^*\rp{2}\). Thus we see that the
complement of a conic\index{conic!complement of} in \(\cp{2}\)
is symplectomorphic to a neighbourhood of \(\rp{2}\) in its
cotangent bundle.

\begin{figure}[htb]
\begin{center}
\begin{tikzpicture}
\filldraw[fill=lightgray,opacity=0.5,draw=black,thick] (0,0) -- (0,3/2) -- (6,0) -- cycle;
\draw[thick] (0,0) -- (0,3/2) -- (6,0) -- cycle;
\draw[dotted,thick] (0,3/2) -- (1,1);
\node at (3,0) [below] {conic};
\node at (3,3/2) {\(\rp{2}\)};
\draw[thick,->] (2.5,3/2) -- (0.5,1.25);
\node at (1,1) {\(\times\)};
\node at (6,0) [right] {\((6,0)\)};
\node at (0,3/2) [left] {\((0,3/2)\)};
\node at (0,0) [left] {\((0,0)\)};
\end{tikzpicture}
\caption{}
\label{fig:atbd_cp_2_rp_2}
\end{center}
\end{figure}

\end{Example}
\begin{Remark}[Exercise \ref{exr:anti_diag_diag}]\label{rmk:anti_diag_diag}
Using\index{almost toric base diagram!S2xS2@$S^2\times S^2$}
the same idea, we find that if \(X = S^2\times S^2\) with its
equal-area symplectic form and \(C\subset X\) is a symplectic
sphere isotopic to the diagonal then \(X\setminus C\) is a
neighbourhood of a Lagrangian sphere\index{Lagrangian!sphere}
in its cotangent bundle\index{cotangent bundle!of sphere}. Of
course these results can be proved directly and explicitly,
but almost toric diagrams give us a way to see them and
generalise them to less obvious examples.

\end{Remark}
\begin{Example}
The complement of a cubic curve\index{cubic curve!complement}
in \(\cp{2}\) can be given the almost toric
diagram\index{almost toric base diagram!complement of
cubic curve} shown in Figure \ref{fig:cp_2_3_trades}(a). This
is a Weinstein domain which retracts onto a visible Lagrangian
cell complex, coloured in the figure: we take the fibre over
the barycentre, together with the three Lagrangian vanishing
thimbles\index{vanishing!thimble} coming from Lemma
\ref{lma:visibleff} living over the three lines connecting the
base-nodes to the barycentre. The attaching maps for these
discs are loops in the barycentric torus living in the
homology classes \((1,2)\), \((1,2)\) and \((1,-1)\). In
Figure \ref{fig:cp_2_3_trades}(b), we draw these three loops
in a square-picture of the barycentric torus. This Lagrangian
cell complex is called the {\em Lagrangian
skeleton}\index{skeleton!complement of plane cubic} of the
complement of the cubic. The paper
\cite{WeinsteinHandlebodies} explores in more detail how to
read Weinstein handlebody decompositions off from almost toric
base diagrams, and the paper \cite{ExactCombinatorics}
explains how a torus with Lagrangian discs attached is all you
need to recover the full complexity of Vianna's
tori\index{Vianna!tori}.

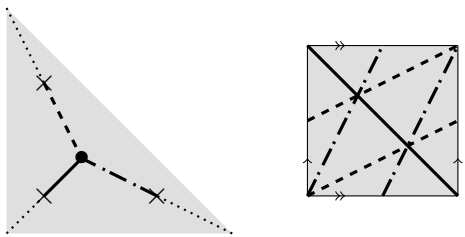
\begin{figure}[htb]
\begin{center}
\begin{tikzpicture}
\filldraw[fill=lightgray,opacity=0.5,draw=none] (0,0) -- (0,3) -- (3,0) -- cycle;
\draw[dotted,thick] (0,0) -- (1/2,1/2);
\draw[dotted,thick] (0,3) -- (0.5,2);
\draw[dotted,thick] (3,0) -- (2,0.5);
\node at (1/2,2) {\(\times\)};
\node at (2,1/2) {\(\times\)};
\node at (1/2,1/2) {\(\times\)};
\draw[black,very thick] (1/2,1/2) -- (1,1);
\draw[dashed,very thick] (1/2,2) -- (1,1);
\draw[dash pattern={on 7pt off 2pt on 1pt off 3pt},very thick] (2,1/2) -- (1,1);
\node at (1,1) {\(\bullet\)};
\begin{scope}[shift={(5,1.5)}]
\filldraw[lightgray,opacity=0.5] (-1,-1) -- (1,-1) -- (1,1) -- (-1,1) -- cycle;
\draw[dashed,very thick] (-1,-1) --+ (2,1);
\draw[dashed,very thick] (-1,0) --+ (2,1);
\draw[dash pattern={on 7pt off 2pt on 1pt off 3pt},very thick] (-1,-1) --+ (1,2);
\draw[dash pattern={on 7pt off 2pt on 1pt off 3pt},very thick] (0,-1) --+ (1,2);
\draw[black,very thick] (-1,1) -- (1,-1);
\draw[->-] (-1,-1) -- (-1,0);
\draw (-1,0) -- (-1,1);
\draw[->-] (1,-1) -- (1,0);
\draw (1,0) -- (1,1);
\draw[->>-] (-1,-1) -- (0,-1);
\draw (0,-1) -- (1,-1);
\draw[->>-] (-1,1) -- (0,1);
\draw (0,1) -- (1,1);
\end{scope}
\end{tikzpicture}
\caption{(a) An almost toric fibration on the complement of a cubic in \(\cp{2}\). (b) The boundaries of the three visible discs on the barycentric torus (thought of as a square with its sides identified in opposite pairs).}
\label{fig:cp_2_3_trades}
\end{center}
\end{figure}

\end{Example}
\section{Cohomology class of the symplectic form}
\label{sct:coh_class}

We have seen that the symplectic area of a curve in the toric
boundary of a toric manifold is given by \(2\pi\) times the
affine length of the edge to which it projects. Since the
cohomology class of the symplectic form is determined by its
integrals over curves, this is enough to determine the
cohomology class\index{symplectic form!cohomology class|(} of
the symplectic form on a toric manifold. We now come up with a
prescription for finding the cohomology class of the symplectic
form on an almost toric manifold. We focus on a restricted class
of base diagrams where the almost toric base is \(\RR^2\) with
some base-nodes; one can easily extend our analysis, e.g.\ to
the case where there are toric critical points.

\begin{Remark}
The case where the base is a closed 2-manifold like the sphere
is trickier because it is impossible to ``read off'' the
symplectic areas of sections. For example, one can change the
symplectic area of a section without changing the almost toric
base diagram by pulling back a non-exact 2-form from the base
of the fibration and adding it to the symplectic form.

\end{Remark}
Here is the class of almost toric base diagrams on which we will
focus. Suppose \(f\colon X\to B\) is an almost toric fibration
where the almost toric base diagram is \(\RR^2\) with a
collection of base-nodes \(n_1,\ldots,n_k\) and branch
cuts. Since we are only interested in the integral affine
structure, we can choose which point we consider to be the
origin \(0\). We can ensure that \(0\not\in\{n_1,\ldots,n_k\}\)
and that the straight line segments \(s_i\) connecting \(0\) to
\(n_i\) intersect only at \(0\). We write the positions of the
base-nodes as \((x_i,y_i)\) relative to this origin. Suppose
that the branch cut at \(n_i\) points in the
\((a_i,b_i)\)-direction, and suppose that this is a primitive
integer eigenvector of the affine monodromy at \(n_i\). We will
suppose for simplicity that none of the line segments \(s_i\)
cross the branch cuts.

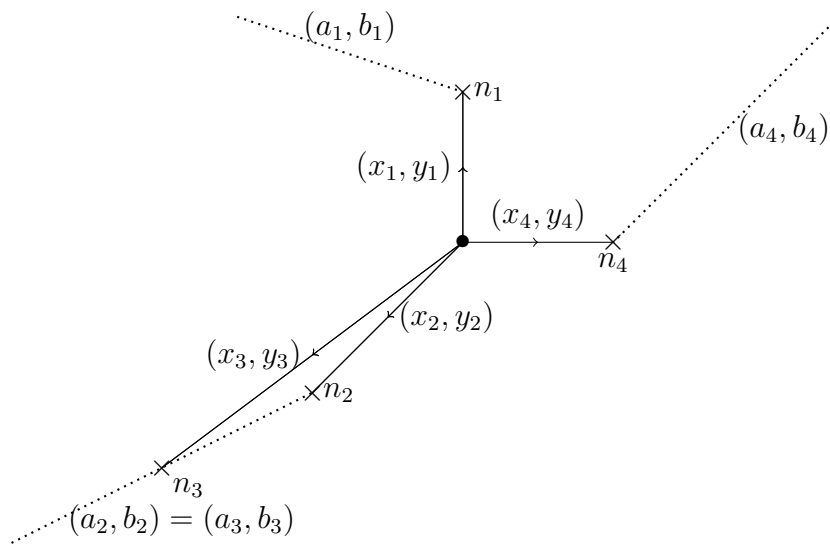
\begin{figure}[htb]
\begin{center}
\begin{tikzpicture}
\draw[thick,dotted] (-3,3) -- (0,2) node {\(\times\)} node [midway,above] {\((a_1,b_1)\)};
\draw[thick,dotted] (-6,-4) -- node [pos=0.3,right] {\((a_2,b_2)=(a_3,b_3)\)} (-4,-3) node {\(\times\)} -- (-2,-2) node {\(\times\)};
\draw[thick,dotted] (5,3) -- (2,0) node {\(\times\)} node [midway,right] {\((a_4,b_4)\)};
\node at (0,0) {\(\bullet\)};
\node at (0,2) [right] {\(n_1\)};
\node at (-2,-2) [right] {\(n_2\)};
\node at (-4,-3) [below right] {\(n_3\)};
\node at (2,0) [below] {\(n_4\)};
\draw[->-] (0,0) -- (0,2);
\draw[->-] (0,0) -- (-2,-2);
\draw[->-] (0,0) -- (-4,-3);
\draw[->-] (0,0) -- (0,2) node [midway,left] {\((x_1,y_1)\)};
\draw[->-] (0,0) -- (-2,-2) node [midway,right] {\((x_2,y_2)\)};
\draw[->-] (0,0) -- (-4,-3) node [midway,left] {\((x_3,y_3)\)};
\draw[->-] (0,0) -- (2,0) node [midway,above] {\((x_4,y_4)\)};
\end{tikzpicture}
\end{center}
\caption{A typical almost toric base diagram of the kind we consider.}
\label{fig:typical_atbd}
\end{figure}

\begin{Theorem}
The second cohomology (with \(\RR\)-coefficients) of \(X\) is
the cokernel of the map \(\partial\colon\RR^2\to\RR^{k+1}\)
defined by \[\partial(v,w)=(0,b_1v-a_1w,\ldots,b_kv-a_kw).\]
If we write \([z_0,z_1,\ldots,z_k]\) for the equivalence class
of \((z_0,z_1,\ldots,z_k)\) in this cokernel, then the
symplectic form lives in the cohomology class \([0,
2\pi(b_1x_1 - a_1y_1), \ldots, 2\pi(b_kx_k - a_ky_k)]\).
\end{Theorem}
\begin{Proof}
The space \(X\) deformation retracts onto the following CW
complex \(W\). Let \(F\) be the torus fibre over the origin
and let \(s_i\) be the line segment connecting \(0\) to
\(n_i\). Inside \(f^{-1}(s_i)\) there is a disc \(\Delta_i\)
with boundary on \(F\) and with centre at the focus-focus
critical point over \(n_i\). The boundary of this disc is a
loop in \(F\) which inhabits the class \((-b_i,a_i)\in
H_1(F;\ZZ)\). Let \(W = F\cup \bigcup_{i=1}^k\Delta_i\).

The fact that \(X\) retracts onto \(W\) can be proved using
Morse theory (though we only sketch it here). The Morse-Bott
function \(|f|^2\) has a global minimum along \(F\) and index
2 critical points at the focus-focus points. The discs
\(\Delta_i\) are the downward manifolds emanating from these
critical points, which gives a handle-decomposition of \(X\)
with \(W\) as the union of cores of handles.

We pick a CW structure on \(F\) with one 0-cell, two 1-cells,
and one 2-cell \(\Delta_0\). To get a CW structure on \(W\) we
need to homotope the attaching maps for the 2-cells
\(\Delta_i\) so that they attach to the 1-skeleton of \(F\)
(which we can do using cellular approximation). Now the
cellular chain complex computing \(H^2(X;\RR)\) is
\[C_2(X;\RR) = \RR^{k+1} \stackrel{\partial}{\to} C_1(X;\RR) =
\RR^2 \stackrel{\partial}{\to} C_0(X;\RR) = \RR.\] The map
\(\partial\colon C_2\to C_1\) is given by looking at the
boundaries of the 2-cells: \(\partial\Delta_0 = 0\) and
\(\partial \Delta_i = (-b_i,a_i)\). We can think of this as a
\(2\)-by-\((k+1)\) matrix \[\begin{pmatrix} 0 & -b_1 & \cdots
& -b_k \\ 0 & a_1 & \cdots & a_k\end{pmatrix}.\] When we take
the dual complex (to compute cohomology), the differential
\(\partial\colon C^1\to C^2\) is given by the transpose of
this matrix. The cohomology is the cokernel of this map (by
definition).

We can compute the integrals of the symplectic form over the
2-cells \(\Delta_i\), which then tells us its cohomology
class. In action-angle coordinates away from the focus-focus
point, \(\Delta_i\) is the cylinder \([0,1)\times
S^1\ni(s,t)\mapsto(x_is,y_is,-b_it,a_it)\), which has
symplectic area \[\int_0^1\int_{S^1}(a_iy_i-b_ix_i)ds\wedge dt
= 2\pi(a_iy_i - b_ix_i)\] as required.\qedhere

\end{Proof}
\begin{Corollary}
The symplectic form on \(X\) is
exact\index{symplectic form!exact} if there is a common point
of intersection of all the eigenlines\index{eigenline!meeting
at common point|(} from the base-nodes.
\end{Corollary}
\begin{Proof}
If we were to use this common point of intersection as our
origin, it would give all the discs \(\Delta_i\) area zero
because each vector \((x_i,y_i)\) would be proportional to
\((a_i,b_i)\).\index{symplectic form!cohomology
class|)}\qedhere

\end{Proof}
\begin{Remark}
For readers who are interested in disc classes with boundary
on Lagrangian tori\index{Lagrangian!torus}, this also tells us
that if there is a common intersection point of the eigenlines
then the Lagrangian torus fibre over this point is
exact\index{Lagrangian!submanifold!exact}. More generally, if
there is a toric boundary divisor, this
Lagrangian torus will be {\em
monotone}\index{monotone!Lagrangian torus} provided the
ambient manifold is monotone\index{monotone!symplectic
manifold}. In particular, Vianna's tori\index{Vianna!tori} are
all monotone.\index{eigenline!meeting at common point|)}

\end{Remark}
\section{Solutions to inline exercises}

\begin{Exercise}\label{exr:which_visible_lagrangians_vianna}
Show that the visible
Lagrangians\index{Lagrangian!submanifold!visible} in the
pictures below have respectively \((2,1)\)- and
\((5,4)\)-pinwheel cores where they meet the edge. In the
first diagram, the direction of the line in the base is
\((1,1)\). In the second diagram, the line points in the
\((-5,4)\)-direction after crossing the branch cut.

\begin{center}
\begin{tikzpicture}
\filldraw[fill=lightgray,opacity=0.5,draw=black,thick] (0,0) -- (0,3) -- (3,0) -- cycle;
\draw[thick] (0,0) -- (0,3) -- (3,0) -- cycle;
\draw[very thick] (1/2,1/2) -- (1.5,1.5);
\draw[dotted,thick] (0,0) -- (1/2,1/2);
\draw[dotted,thick] (0,3) -- (0.5,2);
\draw[dotted,thick] (3,0) -- (2,0.5);
\node at (1/2,2) {\(\times\)};
\node at (2,1/2) {\(\times\)};
\node at (1/2,1/2) {\(\times\)};
\begin{scope}[shift={(4.5,0)}]
\filldraw[fill=lightgray,opacity=0.5,draw=black,thick] (0,0) -- (0,3) -- (3,0) -- cycle;
\draw[thick] (0,0) -- (0,3) -- (3,0) -- cycle;
\draw[very thick] (1/2,1/2) -- (1,1) -- (0,9/5);
\draw[dotted,thick] (0,0) -- (1/2,1/2);
\draw[dotted,thick] (0,3) -- (0.5,2);
\draw[dotted,thick] (3,0) -- (1/2,5/4);
\node at (1/2,2) {\(\times\)};
\node at (1/2,5/4) {\(\times\)};
\node at (1/2,1/2) {\(\times\)};
\end{scope}
\end{tikzpicture}
\end{center}
\end{Exercise}
\begin{Solution}
In each case, let \(\ell\) be the primitive integral vector
pointing along the line of the visible Lagrangian where it
intersects the edge, and let \(e\) be the primitive integral
vector pointing along the edge. Make an integral affine change
of coordinates \(M\) so that \(eM=(1,0)\). By comparison with
Example \ref{exm:hittingedge}, we see that if \(\ell
M=\pm(q,p)\) then we have a visible Lagrangian \((p,q)\)-core
(here \(q\) is determined only up to adding a multiple of
\(p\)).

In the first case, we have \(\ell=(-1,-1)\) and \(e=(-1,1)\),
so we use \(M=\lmatrix -1 &-1\\ 0&-1\rmatrix \) which
gives \(\ell M=(1,2)\).

In the second case, we have \(\ell=(5,-4)\) and \(e=(0,-1)\),
so we use \(M=\lmatrix 0 & 1 \\ -1 & 0 \rmatrix \) which gives
\(\ell M=(4,5)\).\qedhere

\end{Solution}
\begin{Exercise}[Remark \ref{rmk:anti_diag_diag}]\label{exr:anti_diag_diag}
Using the same idea as in Example \ref{exm:cp_2_rp_2}, show
that if \(X = S^2\times S^2\) with its equal-area symplectic
form and \(C\subset X\) is a symplectic sphere isotopic to the
diagonal then \(X\setminus C\) is a neighbourhood of a
Lagrangian sphere in its cotangent bundle.
\end{Exercise}
\begin{Solution}
Figure \ref{fig:atbd_s_2_anti_diag} shows what you get when
you start with the standard square moment map picture of
\(S^2\times S^2\), perform nodal trades at two opposite
corners, and then mutate one of them.

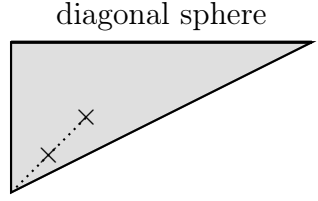
\begin{figure}[htb]
\begin{center}
\begin{tikzpicture}
\filldraw[fill=lightgray,opacity=0.5,draw=black,thick] (0,0) -- (4,2) -- (0,2) -- cycle;
\draw[thick] (0,0) -- (4,2) -- (0,2) -- cycle;
\draw[dotted,thick] (0,0) -- (0.5,0.5) node {\(\times\)} -- (1,1) node {\(\times\)};
\draw[very thick] (0,2) -- (4,2) node [midway,above] {diagonal sphere};
\end{tikzpicture}
\caption{An almost toric picture of \(S^2\times S^2\). The complement of the diagonal sphere is a neighbourhood of the visible Lagrangian sphere living over the dotted line connecting the two base-nodes.}
\label{fig:atbd_s_2_anti_diag}
\end{center}
\end{figure}

After the nodal trades, the toric boundary consists of two
symplectic spheres each isotopic to the diagonal. The
complement of the one marked as ``diagonal sphere'' in the
picture is an open subset of \(B_{2,1,1}\), which is the
Milnor fibre of a \(\QQ\)-Gorenstein smoothing of a
\(\frac{1}{2}(1,1)\) (or \(A_1\)) singularity. But
\(B_{2,1,1}\) is symplectomorphic to the cotangent
bundle\index{Milnor fibre!of A1 singularity@of $A_1$
singularity|(}\index{cotangent bundle!of sphere|(} \(T^*S^2\)
(Exercise \ref{exr:quadric_cotangent}) and the visible
Lagrangian sphere living over the dotted line connecting the
base-nodes (made up of two vanishing thimbles) is the
zero-section.\qedhere

\end{Solution}
\begin{Exercise}[From the solution to Exercise \ref{exr:anti_diag_diag}]\label{exr:quadric_cotangent}
Show that the Milnor fibre of the \(A_1\) singularity is
symplectomorphic to the cotangent bundle of \(S^2\).
\end{Exercise}
\begin{Solution}
The Milnor fibre is an affine quadric\index{quadric
hypersurface!affine} \(z_1z_2 + z_3^2 = 1\). We will write an
explicit symplectomorphism between the affine quadric and
\(T^*S^2\). First we make a change of coordinates
\(z_1=\xi_1+i\xi_2\), \(z_2=\xi_1-i\xi_2\), \(z_3=\xi_3\) to
get the quadric in the form \[\xi_1^2 + \xi_2^2 + \xi_3^2 =
1.\] If we think of \(S^2\) as sitting inside \(\RR^3\) as the
unit sphere \(q_1^2+q_2^2+q_3^2=1\) then its cotangent bundle
sits inside \(\RR^3\times\RR^3\) and consists of pairs
\((\bm{p},\bm{q})\in\RR^3\times S^2\) such that \(\sum_{k=1}^3
p_k q_k = 0\). Here, we give \(\RR^3\times\RR^3\) the
symplectic structure \(\sum_{k=1}^3 dp_k\wedge dq_k\). Now the
map sending \(\bm{\xi}=\bm{x}+i\bm{y}\) to
\((\bm{p},\bm{q})=(-|\bm{x}|\bm{y},\bm{x}/|\bm{x}|)\) is a
symplectomorphism. To see this, write \(r=|\bm{x}|\) and
pullback the 2-form \(\sum dp_k\wedge dq_k\):
\begin{align*} dp_k\wedge dq_k &=
-\left(y_k\,dr+r\,dy_k\right)\wedge
\left(\frac{dx_k}{r}-\frac{x_k\,dr}{r^2}\right)\\ & =
dx_k\wedge dy_k +
(x_k\,dy_k+y_k\,dx_k)\wedge\frac{dr}{r},\end{align*} and use
the fact that \(1=\sum (x_k+iy_k)^2 = \sum(x_k^2-y_k^2)+2i\sum
x_ky_k\), so \(\sum d(x_ky_k)=\sum(x_k\,dy_k+y_k\,dx_k)\)
vanishes on the quadric. Overall, we get \[\sum dp_k \wedge
dq_k = -\sum dx_k\wedge dy_k.\] Note that this shows more
generally that \(T^*S^n\) is symplectomorphic to a smooth
affine quadric in \(n+1\) complex variables.\index{Milnor
fibre!of A1 singularity@of $A_1$
singularity|)}\index{cotangent bundle!of sphere|)}\qedhere

\end{Solution}
\chapter{Surgery}
\label{ch:surgery}
\thispagestyle{cup}

In this chapter, we describe almost toric pictures of some of
the most important surgery operations in 4-dimensional
topology. We first revisit blow-up\index{blow-up!non-toric|(}
but allow ourselves to blow-up at an {\em edge point} on the
toric boundary\index{toric boundary!blowing up a point on|(}
rather than a vertex. Then we discuss rational
blow-up/blow-down. Finally, we use these ideas to explore the
symplectic fillings of lens spaces.

\section{Non-toric blow-up}

Let \(X\) be an almost toric manifold and let \(x\in X\) be a
toric fixed point, i.e.\ a point lying over a Delzant vertex
\(b\) in the almost toric base diagram. We have seen (Example
\ref{exm:symp_cut_blow_up}) that performing the symplectic cut
corresponding to truncating the moment
polygon\index{moment polytope!truncation} at \(b\) has the
result of symplectically blowing-up \(X\) at \(x\). This is
often called {\em toric blow-up}\index{toric blow-up} because it
can be understood purely in terms of toric geometry. But what if
we want to blow-up a point \(x'\in X\) which does not live over
a Delzant vertex? In this section, we will explain how to
blow-up a point living in the toric boundary over an edge of the
base diagram. We begin by describing the local picture; we use a
strategy similar to what we used to analyse the Auroux system
(Example \ref{exm:auroux}).

\begin{Example}
Consider the manifold\index{almost toric base
diagram!non-toric blow-up|(}
\(\mathcal{O}(-1)\)\index{O-1@$\mathcal{O}(-1)$|(} from
Example \ref{exm:taut_bun}. Recall that this is the variety
\[\{(z_1, z_2, [z_3: z_4]) \in \CC^2 \times \cp{1}\,:\,z_1z_4
= z_2z_3\}\] and that this is the blow-up of \(\CC^2\) at the
origin. Pick a real number \(c>0\) and let \(X=\mathcal{O}(-1)
\setminus \{z_1=c\}\). We will write down an almost toric
fibration on \(X\) whose toric boundary is the cylinder
\(\{(\xi,0,[1:0])\,:\,\xi\neq c\}\) and which has one
focus-focus fibre. Namely, we take \[\bm{H}(z_1, z_2, [z_3 :
z_4]) = \left(|z_1-c|^2, \frac{1}{2}|z_2|^2-
\frac{|z_3|^2}{|z_3|^2 + |z_4|^2}\right).\] The function
\(H_2\) satisfies \(H_2\geq -1\) with equality if and only if
\(z_2=z_4=0\). The image of \(\bm{H}\) is
\(\{(b_1,b_2)\in\RR^2\,:\, b_1 > 0, b_2\geq -1\}\). The toric
boundary is the cylinder \(\{(z_1,0,[1:0])\,z_1\neq c\}\)
living over the bottom edge \(b_2=-1\), and the focus-focus
fibre is \(\bm{H}^{-1}(c^2, 0)\) (see Figure
\ref{fig:nontoric_blow_up}).

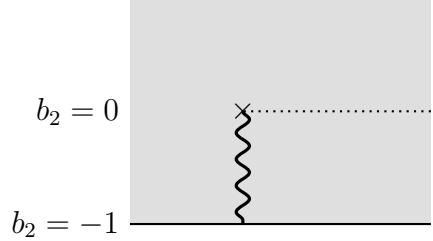
\begin{figure}[htb]
\begin{center}
\begin{tikzpicture}
\filldraw[fill=lightgray,opacity=0.5,draw=none] (0,0) -- (4,0) -- (4,3) -- (0,3) -- cycle;
\draw[thick,black] (0,0) -- (4,0);
\draw[dotted,thick] (4,1.5) -- (1.5,1.5) node {\(\times\)};
\node at (0,1.5) [left] {\(b_2=0\)};
\node at (0,0) [left] {\(b_2=-1\)};
\draw[very thick,snake it] (1.5,1.5) -- (1.5,0);
\end{tikzpicture}
\caption{The local model for nontoric blow-up at a point over an edge in the toric boundary. The wiggly line indicates roughly where the exceptional sphere projects to (finding the precise image would require a nontrivial calculation of action coordinates).}
\label{fig:nontoric_blow_up}
\end{center}
\end{figure}

As with the Auroux system, the Hamiltonian \(H_2\) generates a
circle action, so the conversion to action coordinates has the
form \((b_1,b_2)\mapsto (G_1(b_1,b_2),b_2)\). This means that
(a) the toric boundary in action coordinates is still given by
the horizontal line \(b_2=-1\), (b) the affine
monodromy\index{affine monodromy!non-toric blow-up} around the
base-node is \(\lmatrix 1 & 0 \\ 1 & 1\rmatrix \), so (c) both
the toric boundary and the eigenline of the affine monodromy
are horizontal, and choosing a horizontal branch
cut\index{branch cut!for non-toric blow-up}, we get an almost
toric base diagram of the form shown in Figure
\ref{fig:nontoric_blow_up}. The projection of the exceptional
sphere \(\{(0,0)\}\times\cp{1}\) under \(\bm{H}\) is a
vertical line connecting the base-node to the toric boundary;
I have not checked whether this remains vertical in action
coordinates, so have drawn it as a wiggly line in Figure
\ref{fig:nontoric_blow_up}. Note that regardless of this
projection, the symplectic area of the exceptional sphere is
the affine length between the edge and the
base-node.\index{O-1@$\mathcal{O}(-1)$|)}

Finally, we rotate the branch cut by \(90^{\circ}\) clockwise,
to get the diagram shown in Figure
\ref{fig:nontoric_blow_up_2}. This picture can now be
implanted near to an edge point in any almost toric base
diagram\index{almost toric base diagram!rational elliptic
surface}, see example Figure \ref{fig:rational_elliptic} for
an example of a nine-point blow-up of \(\cp{2}\) (a {\em
rational elliptic surface}).\index{elliptic surface!rational}

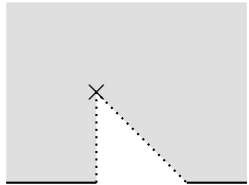
\begin{figure}[htb]
\begin{center}
\begin{tikzpicture}[scale=0.8]
\filldraw[fill=lightgray,opacity=0.5,draw=none] (0,0) -- (1.5,0) -- (1.5,1.5) -- (3,0) -- (4,0) -- (4,3) -- (0,3) -- cycle;
\draw[thick,black] (0,0) -- (1.5,0);
\draw[thick,black] (3,0) -- (4,0);
\draw[dotted,thick] (1.5,0) -- (1.5,1.5) node {\(\times\)} -- (3,0);
\end{tikzpicture}
\caption{A different picture of the same local model as in Figure \ref{fig:nontoric_blow_up}, obtained by rotating the branch cut \(90^{\circ}\) clockwise.}
\label{fig:nontoric_blow_up_2}
\end{center}
\end{figure}

\begin{figure}[htb]
\begin{center}
\begin{tikzpicture}
\filldraw[fill=lightgray,opacity=0.5,draw=none] (0,9) -- (0,6) -- (1,6) -- (0,5) -- (0,4) -- (1,4) -- (0,3) -- (0,2) -- (1,2) -- (0,1) -- (0,0) -- (2,0) -- (2,1) -- (3,0) -- (4,0) -- (4,1) -- (5,0) -- (6,0) -- (6,1) -- (7,0) -- (9,0) -- (7,2) -- (6,2) -- (6,3) -- (5,4) -- (4,4) -- (4,5) -- (3,6) -- (2,6) -- (2,7) -- cycle;
\draw[thick,black] (0,9) -- (0,6);
\draw[thick,black] (0,5) -- (0,4);
\draw[thick,black] (0,3) -- (0,2);
\draw[thick,black] (0,1) -- (0,0);
\draw[thick,black] (9,0) -- (7,0);
\draw[thick,black] (6,0) -- (5,0);
\draw[thick,black] (4,0) -- (3,0);
\draw[thick,black] (2,0) -- (0,0);
\draw[thick,black] (9,0) -- (7,2);
\draw[thick,black] (6,3) -- (5,4);
\draw[thick,black] (4,5) -- (3,6);
\draw[thick,black] (2,7) -- (0,9);
\draw[dotted,thick] (2,0) -- (2,1) -- (3,0);
\node at (2,1) {\(\times\)};
\draw[dotted,thick] (4,0) -- (4,1) -- (5,0);
\node at (4,1) {\(\times\)};
\draw[dotted,thick] (6,0) -- (6,1) -- (7,0);
\node at (6,1) {\(\times\)};
\draw[dotted,thick] (0,1) -- (1,2) -- (0,2);
\node at (1,2) {\(\times\)};
\draw[dotted,thick] (0,3) -- (1,4) -- (0,4);
\node at (1,4) {\(\times\)};
\draw[dotted,thick] (0,5) -- (1,6) -- (0,6);
\node at (1,6) {\(\times\)};
\draw[dotted,thick] (2,7) -- (2,6) -- (3,6);
\node at (2,6) {\(\times\)};
\draw[dotted,thick] (4,5) -- (4,4) -- (5,4);
\node at (4,4) {\(\times\)};
\draw[dotted,thick] (6,3) -- (6,2) -- (7,2);
\node at (6,2) {\(\times\)};
\end{tikzpicture}
\caption{A nontoric blow-up of \(\cp{2}\) in nine balls.}
\label{fig:rational_elliptic}
\end{center}
\end{figure}
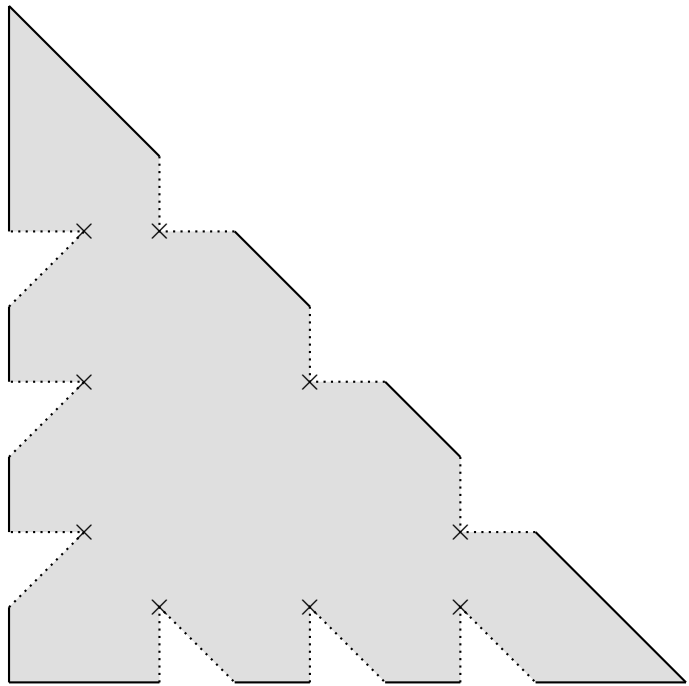

\end{Example}
\begin{Example}\label{exm:atbd_torus_nbhd}
Consider the product symplectic manifold \(X =
T^2\times\CC\). We use \(2\pi\)-periodic coordinates
\((\theta_1,\theta_2)\) on \(T^2\) and Cartesian coordinates
\(x+iy\) on \(\CC\), and equip \(X\) with the symplectic form
\(d\theta_1\wedge d\theta_2 + dx\wedge dy\). The function
\(f\colon T^2\times\CC\to \RR/2\pi\ZZ\times\RR\) defined by
\(f(\theta_1,\theta_2,x+iy) = (\theta_1,(x^2+y^2)/2)\) is a
Lagrangian torus fibration, which induces the product integral
affine structure on \((\RR/2\pi\ZZ)\times\RR\). The torus
\(T:=T^2\times\{0\}\) has self-intersection
zero\index{self-intersection!of torus}. We draw the integral
affine base in Figure \ref{fig:atbd_torus_blow_up}(a); the
dotted lines with arrows indicate that the two sides of the
picture should be identified.

\begin{figure}[htb]
\begin{center}
\begin{tikzpicture}
\node at (3.5,0) {(a)};
\filldraw[lightgray,opacity=0.5] (0,2) -- (0,0) -- (3,0) -- (3,2) -- cycle;
\draw[thick] (0,0) -- (3,0);
\draw[dotted,->-,thick] (0,0) -- (0,2);
\draw[dotted,->-,thick] (3,0) -- (3,2);
\node at (1.5,0) [below] {\(T\)};
\begin{scope}[shift = {(5,0)}]
\node at (3.5,0) {(b)};
\filldraw[lightgray,opacity=0.5] (0,2) -- (0,0) -- (1,0) -- (1,1) -- (2,0) -- (3,0) -- (3,2) -- cycle;
\draw[thick] (0,0) -- (1,0);
\draw[thick,dotted] (1,0) -- (1,1) node {\(\times\)} -- (2,0);
\draw[thick] (2,0) -- (3,0);
\draw[dotted,->-,thick] (0,0) -- (0,2);
\draw[dotted,->-,thick] (3,0) -- (3,2);
\begin{scope}[shift={(5,0)}]
\filldraw[lightgray,opacity=0.5] (0,2) -- (0,0) -- (2,0) -- (3,1) -- (3,2) -- cycle;
\draw[thick] (0,0) -- (2,0);
\draw[thick,dotted] (1,1) node {\(\times\)} -- (3,1);
\draw[dotted,->-,thick] (0,0) -- (0,2);
\draw[dotted,->-,thick] (2,0) -- (3,1) -- (3,2);
\node at (3,0) {(c)};
\node at (1,0) [below] {\(\tilde{T}\)};
\end{scope}
\end{scope}
\end{tikzpicture}
\caption{(a) The base of a Lagrangian torus fibration on \(X = T^2\times\CC\). (b) Non-toric blow-up of \(X\). (c) Same diagram after a change of branch cut.}
\label{fig:atbd_torus_blow_up}
\end{center}
\end{figure}
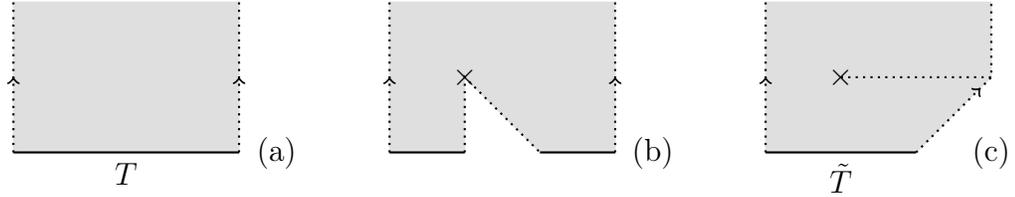

Now perform non-toric blow-up at a point on \(T\); we get the
almost toric manifold whose base diagram is shown in Figure
\ref{fig:atbd_torus_blow_up}(b). Write \(E\) for the
exceptional sphere and \(\tilde{T}\) for the proper transform
of \(T\) (which is the torus living over the toric boundary in
Figure \ref{fig:atbd_torus_blow_up}(b) or (c)). Since
\(\tilde{T}+E\) is homologous to a pushoff of \(T\), we have
\(0 = (\tilde{T}+E)^2 = \tilde{T}^2 +2\tilde{T}\cdot E +
E^2\). Since \(\tilde{T}\cdot E = 1\) and \(E^2 = -1\), we
have \(\tilde{T}^2 = -1\). If we perform a change of branch
cut on Figure \ref{fig:atbd_torus_blow_up}(b), we get the
diagram in Figure \ref{fig:atbd_torus_blow_up}(c). By focusing
on a neighbourhood of \(\tilde{T}\), we get a local (almost
toric) model\index{almost toric base diagram!neighbourhood
of symplectic torus} for a symplectic manifold in a
neighbourhood of a symplectic
torus\index{symplectic torus!neighbourhood of} with
self-intersection \(-1\) (Figure
\ref{fig:atbd_torus_nbhd}(a)). Note that the matrix \(\lmatrix
1 & 0 \\ 1 & 1\rmatrix\) is now used to identify the dotted
edges with arrows. A similar argument gives a picture for tori
with self-intersection \(-n\) (Figure
\ref{fig:atbd_torus_nbhd}(b)).\index{blow-up!non-toric|)}\index{non-toric
blow-up|see {blow-up, non-toric}}\index{almost toric base
diagram!non-toric blow-up|)}\index{toric boundary!blowing
up a point on|)}

\begin{figure}[htb]
\begin{center}
\begin{tikzpicture}
\node at (4,0) {(a)};
\filldraw[lightgray,opacity=0.5] (0,2) -- (0,0) -- (2,0) -- (4,2) -- cycle;
\draw[thick] (0,0) -- (2,0);
\draw[dotted,->-,thick] (0,0) -- (0,2);
\draw[dotted,->-,thick] (2,0) -- (4,2);
\node at (3,1) [right] {\((1,1)\)};
\begin{scope}[shift={(6,0)}]
\node at (6,0) {(b)};
\filldraw[lightgray,opacity=0.5] (0,2) -- (0,0) -- (2,0) -- (6,2) -- cycle;
\draw[thick] (0,0) -- (2,0);
\draw[dotted,->-,thick] (0,0) -- (0,2);
\draw[dotted,->-,thick] (2,0) -- (6,2);
\node at (4,1) [right] {\((n,1)\)};
\end{scope}
\end{tikzpicture}
\caption{(a) Torus with self-intersection \(-1\). (b) Torus with self-intersection \(-n\).}
\label{fig:atbd_torus_nbhd}
\end{center}
\end{figure}
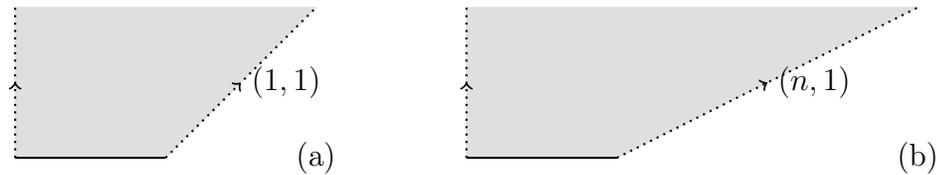

\end{Example}
\begin{Example}
Consider the almost toric diagram of \(\cp{2}\) obtained
from\index{almost toric base diagram!5-point blow up of
CP2@5-point blow-up of $\mathbb{CP}^2$} the standard diagram
by a single nodal trade and mutation:

\begin{center}
\begin{tikzpicture}[scale=2]
\filldraw[draw=none,fill=lightgray,opacity=0.5] (0,0) -- (6,1.5) -- (0,1.5) -- cycle;
\draw[thick] (0,0) -- (6,1.5) -- (0,1.5) -- cycle;
\draw[dotted] (0,0) -- (0.5,0.25) node {\(\times\)};
\end{tikzpicture}
\end{center}

This has a visible Lagrangian \(\rp{2}\) over the branch cut
and the preimage of the top edge is a conic. We can perform
five non-toric blow-ups at points on the conic:

\begin{center}
\begin{tikzpicture}[scale=2]
\filldraw[draw=none,fill=lightgray,opacity=0.5] (0,0) -- (0,1.5) -- (0.16,1.5) -- (0.16,0.5) node {\(\times\)} -- (1.16,1.5) -- (1.32,1.5) -- (0.32,0.5) node {\(\times\)} -- (2.32,1.5) -- (2.48,1.5) -- (0.48,0.5) node {\(\times\)} -- (3.48,1.5) -- (3.64,1.5) -- (0.64,0.5) node {\(\times\)} -- (4.64,1.5) -- (4.8,1.5) -- (0.8,0.5) node {\(\times\)} -- (5.8,1.5) -- (6,1.5) -- cycle;
\draw[thick] (0,0) -- (0,1.5) -- (0.16,1.5);
\draw[dotted] (0.16,1.5) -- (0.16,0.5) -- (1.16,1.5);
\draw[thick] (1.16,1.5) -- (1.32,1.5);
\draw[dotted] (1.32,1.5) -- (0.32,0.5) -- (2.32,1.5);
\draw[thick] (2.32,1.5) -- (2.48,1.5);
\draw[dotted] (2.48,1.5) -- (0.48,0.5) -- (3.48,1.5);
\draw[thick] (3.48,1.5) -- (3.64,1.5);
\draw[dotted] (3.64,1.5) -- (0.64,0.5) -- (4.64,1.5);
\draw[thick] (4.64,1.5) -- (4.8,1.5);
\draw[dotted] (4.8,1.5) -- (0.8,0.5) -- (5.8,1.5);
\draw[thick] (5.8,1.5) -- (6,1.5) -- (0,0);
\draw[dotted] (0,0) -- (0.5,0.25) node {\(\times\)};
\end{tikzpicture}
\end{center}

The ``bites'' we have taken from the edge look a bit different
to the usual picture in Figure \ref{fig:nontoric_blow_up_2},
but are related to it by integral affine transformations. To
be a little more precise, we take the affine length of the top
edge to be \(6\) before the blow-up, and each bite takes out
affine length \(1\) (so the corresponding exceptional sphere
has area \(2\pi\)). The result is a monotone\footnote{i.e.\
the cohomology class of \(\omega\) is a positive multiple of
the first Chern class.} symplectic form on a Del Pezzo surface
obtained by blowing-up \(\cp{2}\) at five points. Changing
branch cuts to make them parallel to their eigenlines, we
obtain:

\begin{center}
\begin{tikzpicture}[scale=2]
\filldraw[draw=none,fill=lightgray,opacity=0.5] (0,0) -- (0,1.5) -- (1,1.5) -- (2,0.5) -- cycle;
\draw[thick] (0,0) -- (0,1.5) -- (1,1.5) -- (2,0.5) -- cycle;
\draw[dotted] (0,0) -- (0.5,0.25) node {\(\times\)};
\draw[dotted] (0.16,0.5) -- (2,0.5);
\node at (0.16,0.5) {\(\times\)};
\node at (0.32,0.5) {\(\times\)};
\node at (0.48,0.5) {\(\times\)};
\node at (0.64,0.5) {\(\times\)};
\node at (0.8,0.5) {\(\times\)};
\end{tikzpicture}
\end{center}

\end{Example}
\section{Rational blow-down/rational blow-up}

Recall\index{rational blow-down|(} that the lens
space\index{lens space!symplectic fillings|(}\index{lens space}
\(L(n,a)\) is the quotient of the standard
contact\index{hypersurface, contact-type|see {contact, -type
hypersurface}}\index{contact!structure on lens space} 3-sphere
in \(\CC^2\) by the action of the cyclic group \(\bm{\mu}_n\) of
\(n\)th roots of unity where \(\mu\) acts by \((z_1,z_2)\mapsto
(\mu z_1,\mu^a z_2)\). This lens space appears as the boundary
of the symplectic orbifold\index{orbifold} \(\CC^2/\bm{\mu}_n\),
which is toric with moment polygon \(\pi(n,a)\) (Example
\ref{exm:cycquot}). We can make symplectic cuts to get the {\em
minimal resolution}\index{singularity!resolution!minimal|(} of
this orbifold singularity as in Example
\ref{exm:minimal_resolution}: we get a smooth symplectic
``filling'' of \(L(n,a)\) which retracts onto a chain of
symplectic spheres\index{symplectic chain of spheres} whose
self-intersection numbers \(b_1,\ldots,b_k\) satisfy
\[\frac{n}{a} = b_1 - \frac{1}{b_2- \frac{1}{\cdots-
\frac{1}{b_k}}}.\] Suppose that \(n=p^2\) and \(a=pq-1\) for
some coprime integers \(p,q\). We have seen another symplectic
filling of the lens space \(L(p^2,pq-1)\), namely
\(B_{1,p,q}\)\index{Bdpq@$B_{d,p,q}$} from Example
\ref{exm:bpq}.\index{lens space!symplectic fillings|)}

\begin{Definition}
Suppose that \(U\subset X\) is a symplectically embedded copy
of (a neighbourhood of the chain of spheres in) the minimal
resolution of \(\frac{1}{p^2}(1,pq-1)\) inside a symplectic
manifold \(X\). The {\em rational
blow-down}\footnote{Originally, rational blow-down was
reserved for the case \(q=1\) and the more general procedure
was called {\em generalised rational blow-down}. One could
generalise still further using \(B_{d,p,q}\) for \(d\geq 2\),
but the beauty of using \(B_{1,p,q}\) is the drastic reduction
in second Betti number that can be achieved.} of \(X\) along
\(U\) is the symplectic manifold obtained by replacing \(U\)
with (an open set in) \(B_{1,p,q}\). Rational blow-up is the
inverse operation.

\end{Definition}
\begin{Example}
The toric diagram in Figure \ref{fig:qbd_1} shows a Hirzebruch
surface containing a symplectic sphere (over the short edge)
with self-intersection \(-4\). We can rationally blow-down
along this sphere and we obtain an almost toric
manifold\index{Lagrangian!RP2@$\mathbb{RP}^2$} containing a
symplectically embedded \(B_{1,2,1}\). In fact, this is
symplectomorphic to \(\cp{2}\) (you can get back to the
standard picture of \(\cp{2}\) by mutating).

\begin{figure}[htb]
\begin{center}
\begin{tikzpicture}
\filldraw[fill=lightgray,opacity=0.5] (0,0) -- (1,0) -- (5,1) -- (0,1) -- cycle;
\draw[thick] (0,0) -- (1,0) -- (5,1) -- (0,1) -- cycle;
\node at (3,1/2) [below right] {\((4,1)\)};
\begin{scope}[shift={(6,0)}]
\filldraw[fill=lightgray,opacity=0.5] (0,-1/4) -- (5,1) -- (0,1) -- cycle;
\draw[thick] (0,-1/4) -- (5,1) -- (0,1) -- cycle;
\draw[dotted,thick] (0,-1/4) -- (1,1/2) node {\(\times\)};
\end{scope}
\end{tikzpicture}
\caption{Rationally blow-down the Hirzebruch surface \(\mathbb{F}_4\) along a neighbourhood of the \(-4\)-sphere to get \(\cp{2}\) (the branch cut points in the \(-(2,1)\)-direction).}
\label{fig:qbd_1}
\end{center}
\end{figure}
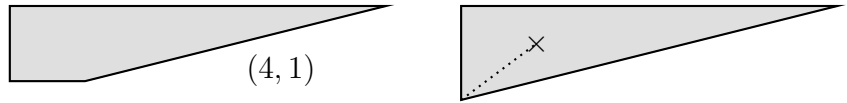

\end{Example}
\begin{Remark}
Rational blow-up/down was introduced by Fintushel and Stern
\cite{FintushelStern}, who showed that certain {\em elliptic
surfaces}\index{elliptic surface} (minimal complex surfaces
with Kodaira dimension 1) could be related by
blow-ups and rational blow-downs, and used this
to calculate their Donaldson
invariants. Symington\index{Symington's theorem!on rational
blow-down} \cite{Symington2} used almost toric methods to show
that these surgeries can be performed symplectically. Rational
blow-down has since been used extensively to construct exotic
4-manifolds\index{small exotic 4-manifolds} with small second
Betti number. This technique was pioneered by Park
\cite{ParkQBD}, who constructed a 4-manifold homeomorphic but
not diffeomorphic to \(\cp{2}\sharp
7\overline{\mathbb{CP}}^2\) by perfoming rational blow-down on
a certain rational surface; the literature on ``small exotic
4-manifolds'' has grown significantly since then.

\end{Remark}
Rather than focusing on exotica, we will content ourselves with
constructing a symplectic filling of a lens space.\index{lens
space!symplectic fillings|(}

\begin{figure}[htb]
\begin{center}
\begin{tikzpicture}
\filldraw[fill=lightgray,opacity=0.5,draw=none] (0,0) node {\(\bullet\)} --++ (1,0) node {\(\bullet\)} --++ (3,1) node {\(\bullet\)} --++ (14/5,5/5) node {\(\bullet\)} --++ (25/9,9/9) node {\(\bullet\)} --++ (36/13,13/13) -- (0,1+5/5+9/9+13/13) -- cycle;
\draw[thick] (0,1+5/5+9/9+13/13) -- (0,0) --++ (1,0) node [midway,below] {\((1,0)\)} --++ (3,1) node [midway,below right] {\((3,1)\)} --++ (14/5,5/5) node [midway,below right] {\((14,5)\)} --++ (25/9,9/9) node [midway,below right] {\((25,9)\)} --++ (36/13,13/13) node [midway,below right] {\((36,13)\)};
\path (0,0) --++ (1,0) node [midway,above] {\(-3\)} --++ (3,1) node [midway,above] {\(-5\)} --++ (14/5,5/5) node [midway,above] {\(-2\)} --++ (25/9,9/9) node [midway,above] {\(-2\)};
\end{tikzpicture}
\caption{The minimal resolution of \(\frac{1}{36}(1,13)\). Delzant vertices are marked with dots and self-intersections of the spheres in the toric boundary are indicated.}
\label{fig:qbd_2}
\end{center}
\end{figure}
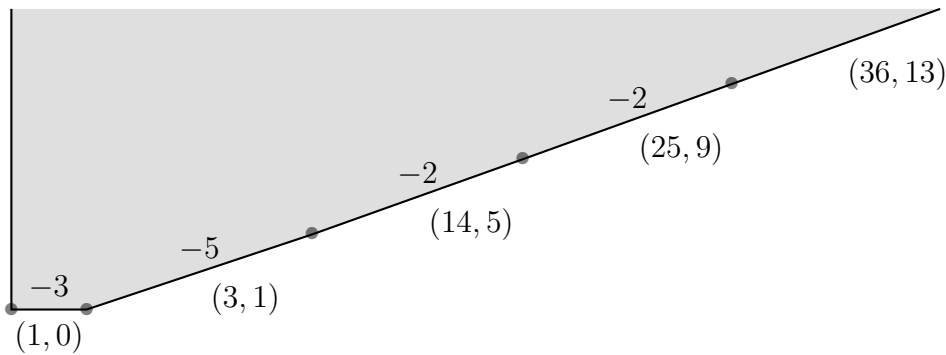

\begin{Example}\label{exm:fillings_36_13}
Consider the singularity \(\frac{1}{36}(1,13)\) and its
minimal resolution\index{singularity!resolution!minimal|)},
which contains a chain of spheres of self-intersections
\(-3,-5,-2,-2\); see Figure \ref{fig:qbd_3}. This is a
symplectic filling of \(L(36,13)\) whose second Betti number
is \(4\). One can rationally blow-down the sub-chain
\(-5,-2\), because \[5-\frac{1}{2} = \frac{9}{2},\] which
corresponds to \(p^2/(pq-1)\) for \(p=3\), \(q=1\). This gives
a different symplectic filling of \(L(36,13)\) whose second
Betti number is 1; see Figure \ref{fig:qbd_2}.

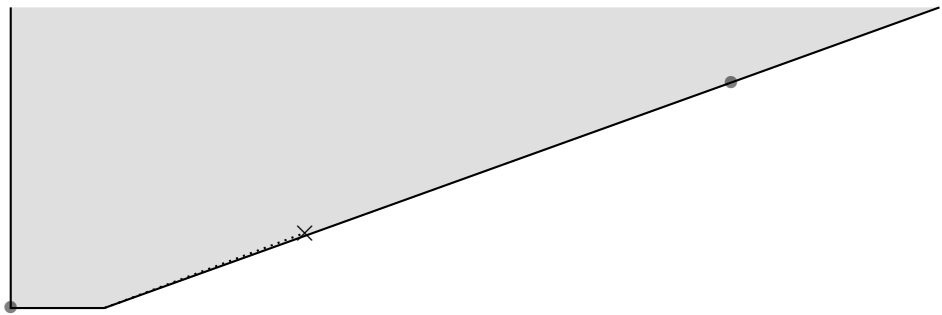
\begin{figure}[htb]
\begin{center}
\begin{tikzpicture}
\filldraw[fill=lightgray,opacity=0.5,draw=none] (0,0) node {\(\bullet\)} -- (56/45,0) -- (1+3+14/5+25/9,3) node {\(\bullet\)} --++ (36/13,13/13) -- (0,1+5/5+9/9+13/13) -- cycle;
\draw[thick] (0,1+5/5+9/9+13/13) -- (0,0) -- (56/45,0) -- (1+3+14/5+25/9,3) --++ (36/13,13/13);
\draw[thick,dotted] (56/45,0) --++ (8/3,3/3) node {\(\times\)};
\end{tikzpicture}
\caption{The result of rationally-blowing down Figure \ref{fig:qbd_2} along the \(-5,-2\) sub-chain of spheres; this contains a copy of \(B_{1,3,1}\). The branch cut points in the \(-(8,3)\)-direction (very close to the edge). The second Betti number of this filling is \(1\) (there is only one compact edge).}
\label{fig:qbd_3}
\end{center}
\end{figure}

Here is a different symplectic filling of \(L(36,13)\) which
also has second Betti number equal to \(1\). First blow-up the
minimal resolution at the intersection point between the
\(-3\) and \(-5\)-spheres. This yields a chain of spheres with
self-intersections
\(-4,-1,-6,-2,-2\). We can rationally blow-down along both the
\(-4\)-sphere and the \(-6,-2,-2\)-subchain, replacing them
with \(B_{1,2,1}\) and \(B_{1,4,1}\) respectively (Figure
\ref{fig:qbd_4}).\index{rational blow-down|)}

\begin{figure}[htb]
\begin{center}
\begin{tikzpicture}
\node at (7723/585,0) {(a)};
\filldraw[fill=lightgray,opacity=0.5,draw=none] (0,0) node {\(\bullet\)} --++ (3-9/4,0) node {\(\bullet\)} --++ (1,0.25) node {\(\bullet\)} --++ (3*0.75,1*0.75) node {\(\bullet\)} --++ (14/5,5/5) node {\(\bullet\)} --++ (25/9,9/9) node {\(\bullet\)} --++ (36/13,13/13) -- (0,1+5/5+9/9+13/13) -- cycle;
\draw[thick] (0,1+5/5+9/9+13/13) -- (0,0) --++ (3-9/4,0) node [midway,below] {\((1,0)\)} --++ (1,0.25) node [midway,below right] {\((4,1)\)} --++ (3*0.75,1*0.75) node [midway,below right] {\((3,1)\)} --++ (14/5,5/5) node [midway,below right] {\((14,5)\)} --++ (25/9,9/9) node [midway,below right] {\((25,9)\)} --++ (36/13,13/13) node [midway,below right] {\((36,13)\)};
\path (0,0) --++ (3-9/4,0) node [midway,above] {\(-4\)} --++ (1,0.25) node [midway,above] {\(-1\)} --++ (3*0.75,1*0.75) node [midway,above] {\(-6\)} --++ (14/5,5/5) node [midway,above] {\(-2\)} --++ (25/9,9/9) node [midway,above] {\(-2\)};
\begin{scope}[shift={(0,-5)}]
\node at (7723/585,0) {(b)};
\filldraw[fill=lightgray,opacity=0.5,draw=none] (0,-3/16) --++ (4*1757/2880,1*1757/2880) -- (7223/585,4) -- (0,4) -- cycle;
\draw[thick] (0,4) -- (0,-3/16) --++ (4*1757/2880,1*1757/2880) node [midway,below right] {\((4,1)\)} -- (7223/585,4) node [midway,below right] {\((36,13)\)};
\draw[thick,dotted] (0,0-3/16) --++ (2,1) node {\(\times\)};
\draw[thick,dotted] (4*1757/2880,1*1757/2880-3/16) --++ (8/3,3/3) node {\(\times\)};
\end{scope}
\end{tikzpicture}
\caption{(a) Blow-up the minimal resolution of \(\frac{1}{36}(1,13)\), (b) then rationally blow-down along the \(-4\) and \(-6,-2,-2\) subchains. The branch cuts point in the \((2,1)\) and \((8,3)\)-directions.}
\label{fig:qbd_4}
\end{center}
\end{figure}
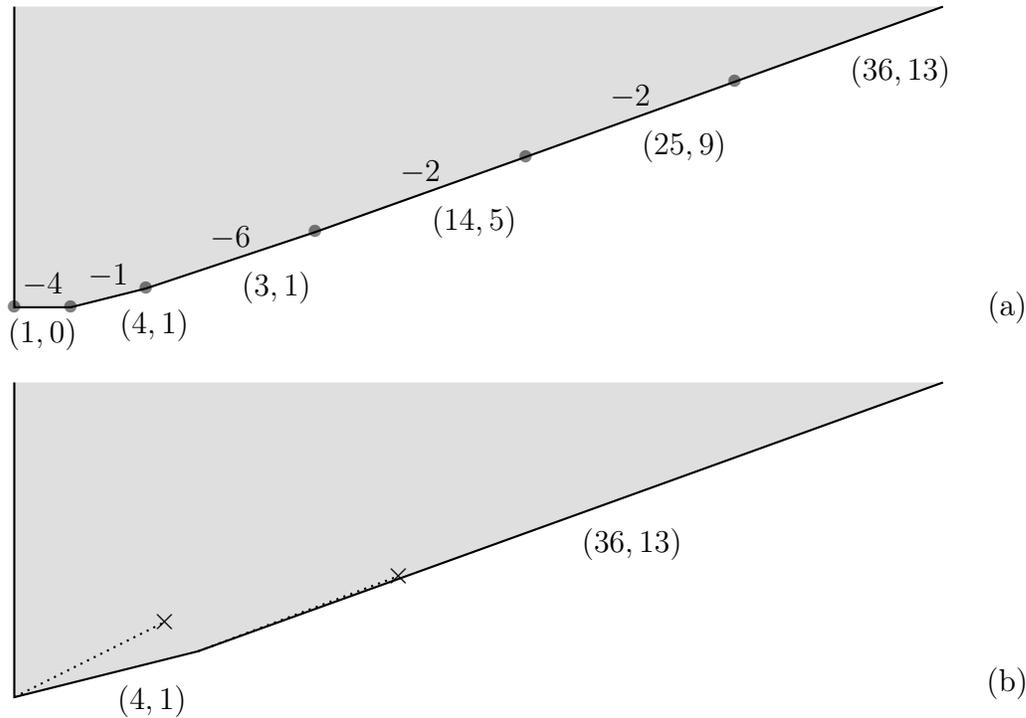

\end{Example}
\begin{Remark}[Exercise \ref{exr:non_diffeo_fillings}]\label{rmk:non_diffeo_fillings}
These two fillings are non-diffeomorphic: the first is
simply-connected, while the second has fundamental group
isomorphic to \(\ZZ/2\).

\end{Remark}
\section{Symplectic fillings of lens spaces}

In Example \ref{exm:fillings_36_13}, we constructed two
non-diffeomorphic fillings of the lens space \(L(36,13)\) with
second Betti number \(1\). In fact, symplectic fillings of lens
spaces\index{contact!structure on lens space}\footnote{There are
many contact structures on lens spaces. We always equip
\(L(n,a)\) with its ``standard'' contact structure descended
from the tight contact structure on \(S^3\).} are completely
classified: Lisca \cite{Lisca} proved the
classification\index{Lisca classification} up to diffeomorphism,
and this was later strengthened to give a classification up to
deformation/symplectomorphism by Bhupal and Ono
\cite{BhupalOno}.

There is an almost toric recipe for constructing all of Lisca's
fillings, which we now explain. We first need to introduce some
ingredients.

\begin{figure}[htb]
\begin{center}
\begin{tikzpicture}
\node at (4,0) {(a)};
\filldraw[lightgray,opacity=0.5] (0,4) -- (0,0) -- (2,0) -- (4,2) -- (4,4) -- cycle;
\draw[thick] (0,4) -- (0,0) node {\(\bullet\)} -- (2,0) node {\(\bullet\)} -- (4,2) node {\(\bullet\)} -- (4,4);
\begin{scope}[shift={(6,0)}]
\node at (4,0) {(b)};
\filldraw[lightgray,opacity=0.5] (0,4) -- (0,0) -- (1,0) -- (3,1) -- (4,2) -- (4,4) -- cycle;
\draw[thick] (0,4) -- (0,0) node {\(\bullet\)} -- (1,0) node {\(\bullet\)} -- (3,1) node {\(\bullet\)} -- (4,2) node {\(\bullet\)} -- (4,4);
\end{scope}
\end{tikzpicture}
\caption{(a) Toric diagram corresponding to the ZCF \([1,1]\). (b) Toric diagram corresponding to a blow-up of \([1,1]\) (in this case, \([2,1,2]\)).}
\label{fig:zcf_toric}
\end{center}
\end{figure}
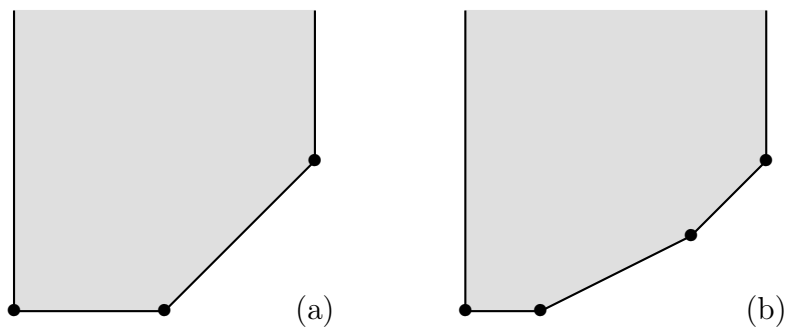

\begin{Definition}
A continued fraction\index{continued fraction!zero (ZCF)|(}
\[[c_1, \ldots, c_m]= c_1- \frac{1}{c_2- \frac{1}{\cdots-
\frac{1}{c_m}}}\] is called a {\em zero continued fraction}
(ZCF) if it evaluates to zero and its evaluation does not
involve dividing by zero at any stage. For example, \([1,1]\)
is a ZCF because \(1-\frac{1}{1}=0\), while
\([1,1,1,1,1]=[1,1,1,0]=[1,1,\infty]=[1,1]=0\) is not.

\end{Definition}
\begin{Example}[Exercise \ref{exr:zcf_blow_up}]\label{exm:zcf_blow_up}
If \([c_1,\ldots,c_m]\) is a ZCF then so are
\[[1,c_1+1,\ldots,c_m],\ \ [c_1,\ldots,c_m+1,1]\mbox{ and }
[c_1,\ldots,c_i+1,1,c_{i+1}+1,\ldots,c_m]\] for any
\(i\in\{1,\ldots,m-1\}\). We call this {\em blowing-up} a
ZCF\index{blow-up!zero continued fraction|(} because it
captures the combinatorics behind the following geometric
procedure. Suppose we have a chain of spheres with
self-intersections \(-c_1,\ldots,-c_m\):
\begin{itemize}
\item If we blow-up a (non-intersection) point on the first sphere
then we get a chain of spheres with self-intersections
\[-1,-c_1-1,-c_2,\ldots,-c_m.\]
\item If we blow-up a (non-intersection) point on the final sphere
then we get a chain of spheres with self-intersections
\[-c_1,-c_2,\ldots,-c_m-1,-1.\]
\item If we blow-up the intersection between the \(i\)th and
\((i+1)\)st sphere then we get a chain of spheres with
self-intersections
\[-c_1,\ldots,-c_i-1,-1,-c_{i+1}-1,\ldots,-c_m.\]
\end{itemize}
We call the obvious inverse procedure {\em blowing-down}.

\end{Example}
\begin{Lemma}\label{lma:zcf_iterated_blow_up}
Any ZCF is obtained from \([1,1]\) by iterated blow-up.
\end{Lemma}
\begin{Proof}
Lemma \ref{lma:zcf_no_1} below implies that any ZCF of length
at least \(2\) contains an entry equal to \(1\), so can be
blown-down to get a shorter continued fraction. Therefore it
is sufficient to prove that the only ZCF of length \(2\) is
\([1,1]\). If \([c_1,c_2]=0\) then \((c_1c_2-1)/c_2=0\), so
\(c_1c_2\) are positive integers whose product is \(1\), hence
\(c_1=c_2=1\).\index{blow-up!zero continued fraction|)}
\qedhere

\end{Proof}
\begin{Lemma}[Exercise \ref{exr:zcf_no_1}]\label{lma:zcf_no_1}
If \([c_1,\ldots,c_m]\) is a continued fraction with \(c_i\geq
2\) for all \(i\) then it is not a ZCF. In fact,
\([c_1,\ldots,c_m] > 1\).

\end{Lemma}
\begin{Corollary}\label{cor:zcf_toric}
If \([c_1,\ldots,c_m]\) is a ZCF then there is a toric
manifold\index{toric manifold!corresponding to zero continued
fraction} containing a chain of spheres with
self-intersections \(-c_1,\ldots,-c_m\) such that the moment
polygon\index{moment polytope!corresponding to zero continued
fraction} is an iterated truncation of Figure
\ref{fig:zcf_toric}(a).
\end{Corollary}
\begin{Proof}
Since \([c_1,\ldots,c_m]\) is obtained by blowing-up
\([1,1]\), we simply follow the geometric procedure outlined
in Example \ref{exm:zcf_blow_up} starting with Figure
\ref{fig:zcf_toric}(a), which contains a chain of spheres with
self-intersections \(-1,-1\). We can take all the blow-ups to
be toric (i.e.\ at vertices of the moment
polygon).\index{continued fraction!zero (ZCF)|)}\qedhere

\end{Proof}

We need one final ingredient.

\begin{Example}
The polygon in Figure \ref{fig:lisca_orbifold}(a) defines a
toric variety \(X\) with two cyclic quotient singularities
\(A\) and \(B\). The singularity \(A\) is
\(\frac{1}{n}(1,a)\); the singularity \(B\) is isomorphic to
\(\frac{1}{n}(1,n-a)\): we can use the matrix\footnote{This
matrix has determinant \(-1\), which is responsible for the
reversed ordering of the exceptional spheres later.}
\(\lmatrix -1 & -1 \\ 0 & 1\rmatrix\) to identify a
neighbourhood of this vertex with a neighbourhood of the
vertex in \(\pi(n,n-a)\). Let \(\tilde{X}\) be the toric
variety obtained by taking the minimal resolution of \(X\) at
\(B\); this has the moment polygon shown in Figure
\ref{fig:lisca_orbifold}(b). Note that if
\(\frac{n}{n-a}=[b_1,\ldots,b_m]\) then the self-intersections
of the spheres in the minimal resolution appear as
\(-b_m,-b_{m-1},\ldots,-b_1\) as we traverse the boundary
anticlockwise.

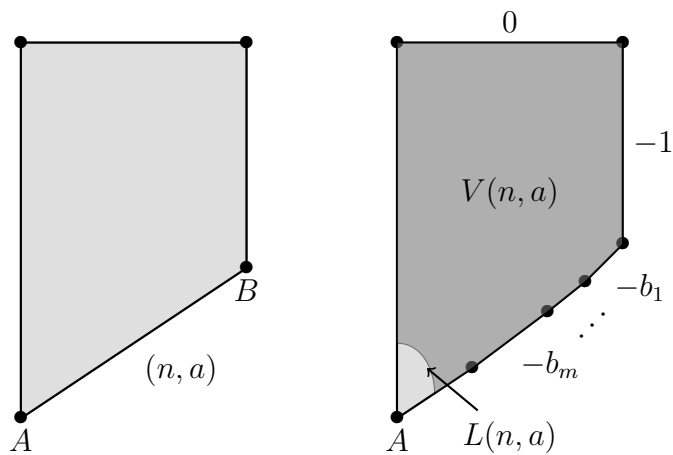
\begin{figure}[htb]
\begin{center}
\begin{tikzpicture}
\filldraw[fill=lightgray,opacity=0.5,draw=black,thick] (0,5) node {\(\bullet\)} -- (0,0) node {\(\bullet\)} -- (3,2) node {\(\bullet\)} -- (3,5) node {\(\bullet\)} -- cycle;
\draw[thick] (0,5) node {\(\bullet\)} -- (0,0) node {\(\bullet\)} -- (3,2) node {\(\bullet\)} -- (3,5) node {\(\bullet\)} -- cycle;
\node at (0,0) [below] {\(A\)};
\node at (3,2) [below] {\(B\)};
\node at (1.5,1) [below right] {\((n,a)\)};
\begin{scope}[shift={(5,0)}]
\filldraw[fill=lightgray,opacity=0.5,draw=black,thick] (0,5) node {\(\bullet\)} -- (0,0) node {\(\bullet\)} -- (1,2/3) node {\(\bullet\)} --++ (1,3/4) node {\(\bullet\)} --++ (0.5,4/5*0.5) node {\(\bullet\)} --++ (0.5,0.5) node {\(\bullet\)} -- (3,5) node {\(\bullet\)} -- cycle;
\draw[thick] (0,5) node {\(\bullet\)} -- (0,0) node {\(\bullet\)} -- (1,2/3) node {\(\bullet\)} --++ (1,3/4) node {\(\bullet\)} --++ (0.5,4/5*0.5) node {\(\bullet\)} --++ (0.5,0.5) node {\(\bullet\)} -- (3,5) node {\(\bullet\)} -- cycle;
\filldraw[fill=gray,opacity=0.5] (0,5) -- (0,1) to[out=0,in=90] (0.5,1/3) -- (1,2/3) --++ (1,3/4) --++ (0.5,4/5*0.5) --++ (0.5,0.5) -- (3,5) -- cycle;
\path (1,2/3) --++ (1,3/4) node [midway,below right] {\(-b_m\)} --++ (0.5,4/5*0.5) node (p) [midway,below right] {} --++ (0.5,0.5) node [midway,below right] {\(-b_1\)} -- (3,5) node [midway,right] {\(-1\)} -- (0,5) node [midway,above] {\(0\)};
\node at (0,0) [below] {\(A\)};
\node at (1.5,3) {\(V(n,a)\)};
\node (L) at (1.5,-0.25) {\(L(n,a)\)};
\draw[->,thick] (L) -- (0.4,2/3);
\node at (p) [below right] {\(\cdot\)};
\node at (p.north east) [below right] {\(\cdot\)};
\node at (p.south west) [below right] {\(\cdot\)};
\end{scope}
\end{tikzpicture}
\caption{(a) A toric variety \(X\) with two orbifold singularities: \(A\) of type \(\frac{1}{n}(1,a)\) and \(B\) of type \(\frac{1}{n}(1,n-a)\). (b) The variety \(\tilde{X}\) obtained by taking the minimal resolution of \(X\) at \(B\). The self-intersections of curves in the toric boundary are indicated. The shaded region is the submanifold \(V(n,a)\) whose concave contact boundary is \(L(n,a)\).}
\label{fig:lisca_orbifold}
\end{center}
\end{figure}

We define the symplectic manifold \(V(n,a)\) to be the
complement of a neighbourhood of \(A\) in \(\tilde{X}\),
shaded in Figure \ref{fig:lisca_orbifold}(b). This submanifold
has {\em concave} contact boundary\footnote{See Definition
\ref{dfn:liouville_vector_field}.} \(L(n,a)\).

\end{Example}
We will now give a recipe for constructing symplectic fillings
of lens spaces.

\begin{Recipe}\label{recipe:lisca}
Let\index{Lisca classification} \([b_1,\ldots,b_m]\) be the
continued fraction\index{continued fraction}\index{continued
fraction!zero (ZCF)} expansion of \(n/(n-a)\). Suppose
that\footnote{The reverse ordering is not a typo!}
\([c_m,\ldots,c_1]\) is a ZCF with \(c_i\leq b_i\) for
\(i=1,\ldots,m\). Let \(Y\) be the corresponding toric
variety\index{toric manifold!corresponding to zero continued
fraction} from Corollary \ref{cor:zcf_toric}. If we perform
\(b_i-c_i\) non-toric blow-ups on the edge with
self-intersection \(-c_i\) then we obtain an almost toric
manifold\index{almost toric base diagram!filling of lens
space|(} containing a chain of spheres with self-intersections
\(-b_m,\ldots,-b_1,-1,0\). A neighbourhood of this chain is
symplectomorphic to \(V(n,a)\), so the complement of a
neighbourhood of this chain gives a symplectic filling of
\(L(n,a)\).

\end{Recipe}
We illustrate this recipe with some simple examples.

\begin{Example}\label{exm:mcduff_filling}
The fillings of \(L(4,1)\) were classified earlier by McDuff
{\cite[Theorem 1.7]{McDuffRatRuled}}. Up to deformation, there
are two: \(B_{1,2,1}\) and \(\mathcal{O}(-4)\). We will
construct these using Recipe \ref{recipe:lisca}. We have
\(\frac{4}{4-1}=[2,2,2]\). There are two possible ZCFs for use
in Recipe \ref{recipe:lisca}: \([2,1,2]\) and \([1,2,1]\).

{\bf Case \([2,1,2]\):} We start with the toric manifold shown
in Figure \ref{fig:mcduff_filling_1}(a) and perform a
non-toric blow-up on the middle edge, yielding Figure
\ref{fig:mcduff_filling_2}(a). The red region is \(V(4,1)\)
and its complement is a symplectic filling of \(L(4,1)\). If
we perform a change of branch cut (Figure
\ref{fig:mcduff_filling_3}(a)), we see that this is precisely
the almost toric diagram of \(B_{1,2,1}\) from Example
\ref{exm:bpq}.

{\bf Case \([1,2,1]\):} We start with the toric manifold shown
in Figure \ref{fig:mcduff_filling_1}(b) and perform a
non-toric blow-up on the two outer edges, yielding Figure
\ref{fig:mcduff_filling_2}(b). The shaded region is \(V(4,1)\)
and its complement is a symplectic filling of \(L(4,1)\). If
we perform changes of branch cut (Figure
\ref{fig:mcduff_filling_3}(b)), we see that this is precisely
the almost toric diagram of \(\mathcal{O}(-4)\) from Example
\ref{exm:o_n}.

\begin{figure}[htb]
\begin{center}
\begin{tikzpicture}[scale=0.8]
\node at (4,0) {(a)};
\filldraw[fill=lightgray,opacity=0.5,draw=black,thick] (0,3) node {\(\bullet\)} -- (0,0) node {\(\bullet\)} -- (1,0) node {\(\bullet\)} -- (3,1) node {\(\bullet\)} -- (4,2) node {\(\bullet\)} -- (4,3) node {\(\bullet\)} -- cycle;
\draw[thick] (0,3) node {\(\bullet\)} -- (0,0) node {\(\bullet\)} -- (1,0) node {\(\bullet\)} -- (3,1) node {\(\bullet\)} -- (4,2) node {\(\bullet\)} -- (4,3) node {\(\bullet\)} -- cycle;
\path (0,0) -- (1,0) node [midway,above] {\(-2\)} -- (3,1) node [midway,above left] {\(-1\)} -- (4,2) node [midway,above left] {\(-2\)};
\path (0,0) -- (1,0) node [midway,below] {\((1,0)\)} -- (3,1) node [midway,below right] {\((2,1)\)} -- (4,2) node [midway,below right] {\((1,1)\)};
\begin{scope}[shift={(6,0)}]
\filldraw[fill=lightgray,opacity=0.5,draw=black,thick] (0,4) node {\(\bullet\)} -- (0,0) node {\(\bullet\)} -- (1,0) node {\(\bullet\)} -- (2,1) node {\(\bullet\)} -- (3,3) node {\(\bullet\)} -- (3,4) node {\(\bullet\)} -- cycle;
\draw[thick] (0,4) node {\(\bullet\)} -- (0,0) node {\(\bullet\)} -- (1,0) node {\(\bullet\)} -- (2,1) node {\(\bullet\)} -- (3,3) node {\(\bullet\)} -- (3,4) node {\(\bullet\)} -- cycle;
\path (0,0) -- (1,0) node [midway,above] {\(-1\)} -- (2,1) node [midway,above left] {\(-2\)} -- (3,3) node [midway,above left] {\(-1\)};
\path (0,0) -- (1,0) node [midway,below] {\((1,0)\)} -- (2,1) node [midway,below right] {\((1,2)\)} -- (3,3) node [midway,below right] {\((1,1)\)};
\node at (3,0) {(b)};
\end{scope}
\end{tikzpicture}
\caption{(a) The toric variety associated with the ZCF \([2,1,2]\). (b) The toric variety associated with the ZCF \([1,2,1]\).}
\label{fig:mcduff_filling_1}
\end{center}
\end{figure}
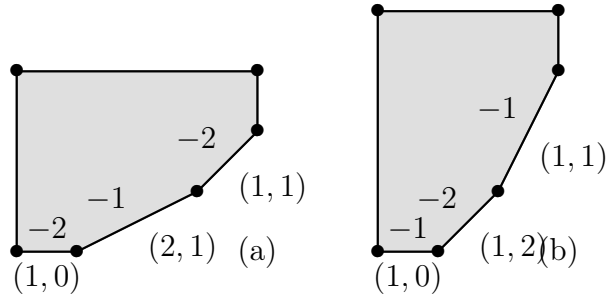

\begin{figure}[htb]
\begin{center}
\begin{tikzpicture}[scale=0.8]
\node at (4,0) {(a)};
\filldraw[fill=lightgray,opacity=0.5,draw=none] (0,3) -- (0,0) -- (1,0) -- (1.25,0.125) -- (1.75,0.625) -- (2.25,0.625) -- (3,1) -- (4,2) -- (4,3) -- cycle;
\filldraw[gray,opacity=0.5] (0,0.3) to[out=0,in=-175] (1.5,0.375) -- (1.25,0.125) -- (1,0) -- (0,0) -- cycle;
\filldraw[gray,opacity=0.5] (2.25,0.625) -- (3,1) -- (4,2) -- (4,3) -- (0,3) -- (0,2.7) to[out=0,in=45] (2,0.625) -- cycle;
\draw[draw=black,thick] (2.25,0.625) -- (3,1) node {\(\bullet\)} -- (4,2) node {\(\bullet\)} -- (4,3) node {\(\bullet\)} -- (0,3) node {\(\bullet\)} -- (0,0) node {\(\bullet\)} -- (1,0) node {\(\bullet\)} -- (1.25,0.125);
\draw[thick,dotted] (1.25,0.125) -- (1.75,0.625) node {\(\times\)} -- (2.25,0.625);
\node at (0.5,0) [below] {\(-2\)};
\node at (2.75,0.8) [below] {\(-2\)};
\node at (3.5,1.5) [below] {\(-2\)};
\node at (4,2.5) [right] {\(-1\)};
\node at (2,3) [above] {\(0\)};
\begin{scope}[shift={(6,0)}]
\filldraw[fill=lightgray,opacity=0.5,draw=none] (0,4) -- (0,0) -- (0.25,0) -- (0.25,0.5) -- (0.75,0) -- (1,0) -- (2,1) -- (2.25,1.5) -- (2.25,2) -- (2.75,2.5) -- (3,3) -- (3,4) -- cycle;
\filldraw[gray,opacity=0.5] (0,0.3) -- (0.25,0.3) -- (0.25,0) -- (0,0) -- cycle;
\filldraw[gray,opacity=0.5] (0.45,0.3) to[out=0,in=-116] (2.25,1.5+0.3) -- (2.25,1.5) -- (2,1) -- (1,0) -- (0.75,0) -- cycle;
\filldraw[gray,opacity=0.5] (2.45,2.2) to[out=64,in=0] (0,3.7) -- (0,4) -- (3,4) -- (3,3) -- (2.75,2.5) -- cycle;
\draw[thick,black] (2.75,2.5) -- (3,3) node {\(\bullet\)} -- (3,4) node {\(\bullet\)} -- (0,4) node {\(\bullet\)} -- (0,0) node {\(\bullet\)} -- (0.25,0);
\draw[thick,black] (0.75,0) -- (1,0) node {\(\bullet\)} -- (2,1) node {\(\bullet\)} -- (2.25,1.5);
\draw[thick,dotted] (0.25,0) -- (0.25,0.5) node {\(\times\)} -- (0.75,0);
\draw[thick,dotted] (2.25,1.5) -- (2.25,2) node {\(\times\)} -- (2.75,2.5);
\node at (3,0) {(b)};
\node at (0.5,0) [below] {\(-2\)};
\node at (1.5,0.5) [below] {\(-2\)};
\node at (2.25,1.25) [below right] {\(-2\)};
\node at (3,3.5) [right] {\(-1\)};
\node at (1.5,4) [above] {\(0\)};
\end{scope}
\end{tikzpicture}
\caption{(a) Non-toric blow-up of Figure \ref{fig:mcduff_filling_1}(a). (b) Non-toric blow-up of Figure \ref{fig:mcduff_filling_1}(b). The submanifold \(V(4,1)\) is shaded in both diagrams, and the numbers along the edges indicate the self-intersections of spheres in \(V(4,1)\). The remaining grey regions are the fillings of \(L(4,1)\). In all cases, the eigenlines for the affine monodromy are parallel to the blown-up edges.}
\label{fig:mcduff_filling_2}
\end{center}
\end{figure}
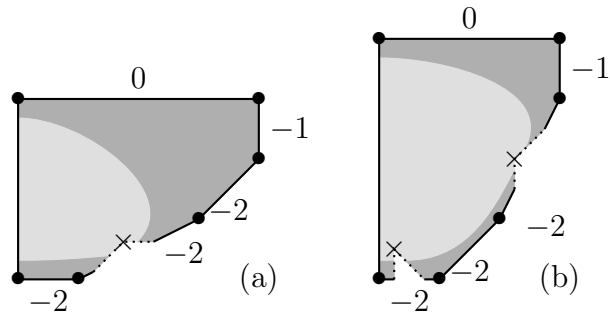

\clearpage

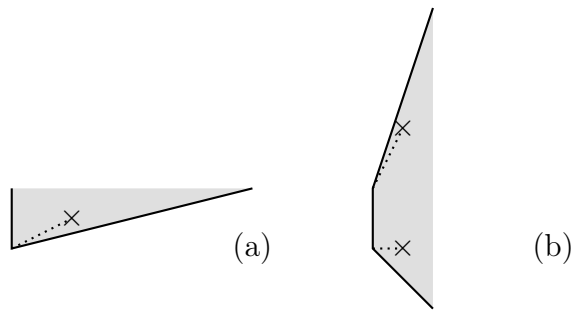
\begin{figure}[htb]
\begin{center}
\begin{tikzpicture}[scale=0.8]
\node at (4,0) {(a)};
\filldraw[fill=lightgray,opacity=0.5,draw=none] (0,1) -- (0,0) -- (4,1) -- cycle;
\draw[thick] (0,1) -- (0,0) -- (4,1);
\draw[thick,dotted] (0,0) -- (1,0.5) node {\(\times\)};
\begin{scope}[shift={(6,0)}]
\filldraw[fill=lightgray,opacity=0.5,draw=none] (1,-1) -- (0,0) -- (0,1) --++ (1,3) -- cycle;
\draw[thick] (1,-1) -- (0,0) -- (0,1) --++ (1,3);
\draw[thick,dotted] (0,0) -- (0.5,0) node {\(\times\)};
\draw[thick,dotted] (0,1) --++ (0.5,1) node {\(\times\)};
\node at (3,0) {(b)};
\end{scope}
\end{tikzpicture}
\caption{(a) The filling from Figure \ref{fig:mcduff_filling_2}(a) after a change of branch cut; this is \(B_{1,2,1}\) from Example \ref{exm:bpq}. (b) The filling from Figure \ref{fig:mcduff_filling_2}(b) after a change of branch cut; this is the standard toric picture of \(\mathcal{O}(-4)\) after two nodal trades and a shear. In both cases we have also made deformations of the symplectic form to shift the base-nodes around.}
\label{fig:mcduff_filling_3}
\end{center}
\end{figure}

\end{Example}
\begin{Theorem}[Lisca \cite{Lisca}, Bhupal-Ono \cite{BhupalOno}]\label{thm:lisca}
Any\index{Lisca classification} symplectic filling of a lens
space is symplectomorphic to a deformation of a filling
constructed using Recipe \ref{recipe:lisca}.

\end{Theorem}
We will not prove this theorem, as it uses the theory of
holomorphic curves in a nontrivial way. Roughly, the idea is to
cap off a symplectic filling \(U\) of \(L(n,a)\) using
\(V(n,a)\); the result contains an embedded symplectic sphere of
self-intersection \(1\) because \(V(n,a)\) does (a smoothing of
the \(0\) and \(-1\)-spheres in the chain). By a result of
McDuff {\cite[Theorem 1.1(i) + Theorem 1.4]{McDuffRatRuled}},
this implies that \(U\cup V(n,a)\) is rational (an iterated
symplectic blow-up of \(\cp{2}\)), so one reduces to studying
different ways that \(V(n,a)\) can embed in a rational
4-manifold.

\begin{Remark}
If we remove all of the base-nodes from one
of these almost toric fillings by performing inverse
generalised nodal trades\index{nodal trade!generalised}, we
obtain a (possibly) singular toric variety with cyclic
quotient T-singularities\index{singularity!T-}. This is a partial
resolution\index{singularity!resolution!partial} of
\(\frac{1}{n}(1,a)\). Note that Koll\'{a}r and Shepherd-Barron
used techniques from Mori theory in \cite{KSB} to show that
any smoothing\index{smoothing} of \(\frac{1}{n}(1,a)\) can be
obtained as a \(\QQ\)-Gorenstein
smoothing\index{smoothing!Q-Gorenstein@$\mathbb{Q}$-Gorenstein}
of a {\em P-resolution}\index{singularity!resolution!P-} of
\(\frac{1}{n}(1,a)\), that is a partial resolution with at
worst T-singularities and such that all exceptional curves
pair nonnegatively with the canonical class. Stevens
\cite{Stevens} and Christophersen \cite{Christophersen} showed
that P-resolutions of \(\frac{1}{n}(1,a)\) are in bijection
with ZCFs \([c_1,\ldots,c_m]\) with \(c_i\leq b_i\) where
\(n/(n-a) = [b_m,\ldots,b_1]\). The aforementioned
prescription gives us a way to go directly between the Lisca
description of the filling and the Koll\'{a}r-Shepherd-Barron
P-resolution. It also follows from {\cite[Lemma 3.14]{KSB}}
that all symplectic fillings of \(L(n,a)\) are obtained by a
combination of blow-downs and rational blow-downs from a
certain {\em maximal resolution}\index{maximal resolution|see
{singularity, resolution,
maximal}}\index{singularity!resolution!maximal} which is a
P-resolution dominating the minimal resolution.

\end{Remark}
\section{Solutions to inline exercises}

\begin{Exercise}[Remark \ref{rmk:non_diffeo_fillings}]\label{exr:non_diffeo_fillings}
The two fillings from Example \ref{exm:fillings_36_13} are
non-diffeomorphic: the first is simply-connected, while the
second has fundamental group isomorphic to \(\ZZ/2\).
\end{Exercise}
\begin{Solution}
We can deformation retract these almost toric manifolds onto
the union of the almost toric boundary and the vanishing
thimbles living over the branch cuts.

In the first case, the almost toric boundary is a sphere (over
the compact edge) with two planes (over the non-compact edges)
attached. We can ignore the planes as they are
contractible. Thus the almost toric manifold deformation
retracts onto a CW complex obtained by attaching a 2-cell to a
sphere; the fundamental group is trivial by Van Kampen's
theorem.

In the second case, the almost toric boundary is a cylinder
\(C\) (with \(\pi_1(C)=\ZZ\)), and there are two vanishing
thimbles. The boundaries of the vanishing thimbles attach to
the cylinder in homotopy classes which correspond to the even
numbers \(2\) and \(4\) in \(\pi_1(C)\). You can see this
because the vanishing thimbles are part of visible Lagrangian
\((2,1)\)- and \((4,1)\)-pinwheels, so the thimble ``caps
off'' a loop which wraps twice (respectively four times)
around the cylinder. By Van Kampen's theorem again, this means
that the fundamental group of this CW complex is the quotient
of \(\ZZ\) by the subgroup generated by \(2\) and \(4\), that
is \(\ZZ/2\).\index{lens space!symplectic
fillings|)}\index{almost toric base diagram!filling of
lens space|)}\qedhere

\end{Solution}
\begin{Exercise}[Example \ref{exm:zcf_blow_up}]\label{exr:zcf_blow_up}
If \([c_1,\ldots,c_m]\) is a zero continued
fraction\index{continued fraction!zero (ZCF)|(} then so are
\([1,c_1+1,\ldots,c_m]\), \([c_1,\ldots,c_m+1,1]\) and
\([c_1,\ldots,c_i+1,1,c_{i+1}+1,\ldots,c_m]\) for any
\(i\in\{1,\ldots,m-1\}\).
\end{Exercise}
\begin{Solution}
We deal with these three cases in order. Let
\(x=[c_2,\ldots,c_m]\). We have \[[1,c_1+1,c_2,\ldots,c_m] =
1-\frac{1}{c_1+1-1/x} =
\left(c_1-1/x\right)/\left(c_1+1-1/x\right).\] Since
\(c_1-1/x=0\), this equals zero, proving the first case.

Blowing-up at the end of a continued fraction does not change
its value: \[c_m = c_m+1 - \frac{1}{1},\] proving the second
case.

Finally, we blow-up in the middle of the chain. Let
\(x=[c_{i+2},\ldots,c_m]\). We need to show that
\[c_i+1-\frac{1}{1-\frac{1}{c_{i+1}+1-1/x}} = c_i -
\frac{1}{c_{i+1}-1/x}.\] We have
\begin{align*}
c_i+1-\frac{1}{1-\frac{1}{c_{i+1}+1-\frac{1}{x}}} &= c_i + 1-\frac{c_{i+1}+1-1/x}{c_{i+1}-1/x}\\
&= c_i + \frac{c_{i+1}-1/x-(c_{i+1}+1-1/x)}{c_{i+1}-1/x}\\
&= c_i - \frac{1}{c_{i+1}-1/x}.
\end{align*}
This shows that blowing up the zero continued fraction does
not change its value.\qedhere

\end{Solution}
\begin{Exercise}[Lemma \ref{lma:zcf_no_1}]\label{exr:zcf_no_1}
If \([c_1,\ldots,c_m]\) is a continued fraction with \(c_i\geq
2\) for all \(i\) then it is not a ZCF. In fact,
\([c_1,\ldots,c_m] > 1\).
\end{Exercise}
\begin{Proof}
We prove this by induction on the length of the continued
fraction. It is clearly true if \(m=1\). Assume it is true for
all continued fractions with length \(m\) and all entries
\(\geq 2\); let \([c_1,\ldots,c_{m+1}]\) be a continued
fraction of length \(m+1\).\index{continued fraction!zero
(ZCF)|)} Then \([c_2,\ldots,c_m] > 1\) and \(c_1\geq 2\) by
assumption, so \[[c_1,\ldots,c_{m+1}] =
c_1-1/[c_2,\ldots,c_{m+1}] > 2-1/1 =1.\qedhere\]

\end{Proof}
\chapter{Elliptic and cusp singularities}
\label{ch:cusps}
\thispagestyle{cup}

I first learned about the following pictures in conversation
with Paul Hacking, and then finally understood them by reading
the paper \cite{Engel} by Philip Engel.

\section{Another picture of \(\cp{2}\)}

Recall that to obtain an almost toric diagram for an almost
toric fibration \(f\colon X\to B\), we picked a simply-connected
fundamental domain for the deck group action in the universal
cover \(\widetilde{B}^{reg}\). Since this fundamental domain is
simply-connected, we can find single-valued action coordinates
on the whole domain, and we took the almost toric diagram
(fundamental action domain) to be the image of the fundamental
domain under the action coordinates. In this section, we will
allow ourselves something more exotic: we will take a branch
cut\index{branch cut!more exotic|(} in \(B\) whose complement is
not simply-connected, but rather has fundamental group
\(\ZZ\). The action coordinates will be multi-valued, but
related by a \(\ZZ\)-action which will be given by iterated
application of an integral affine matrix. As a result, our
pictures will have a high degree of redundancy (we are
superimposing infinitely many fundamental action domains) but
this can be fixed by quotienting them by the \(\ZZ\)-action. The
result will be a ``conical'' almost toric diagram\index{almost
toric base diagram!conical} rather than a planar diagram.

\begin{Example}\label{exm:cp2_elliptic}
Consider the standard almost toric picture of \(\cp{2}\) where
we have made three nodal trades at the corners; there are
three branch cuts extending from the nodes to the corners
(Figure \ref{fig:tripod1}). The almost toric boundary is a
cubic curve with self-intersection 9. We will redraw this
picture as follows. Let \(A,B,C\) be the three nodes and let
\(O\) be the barycentre of the triangle; let \(T\) be the
tripod of lines \(OA,OB,OC\) (shown dashed in Figure
\ref{fig:tripod1}). Let \(\alpha,\beta,\gamma\) be the three
triangular regions labelled in Figure \ref{fig:tripod1}.

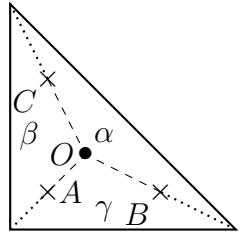
\begin{figure}[htb]
\begin{center}
\begin{tikzpicture}
\draw[thick] (0,0) -- (3,0) -- (0,3) -- cycle;
\draw[dotted,thick] (0,0) -- (1/2,1/2);
\draw[dotted,thick] (3,0) -- (2,1/2);
\draw[dotted,thick] (0,3) -- (1/2,2);
\node (A) at (1/2,1/2) {\(\times\)};
\node (B) at (2,1/2) {\(\times\)};
\node (C) at (1/2,2) {\(\times\)};
\node at (A) [right] {\(A\)};
\node at (B) [below left] {\(B\)};
\node at (C) [below left] {\(C\)};
\node at (1,1) {\(\bullet\)};
\node at (1,1) [left] {\(O\)};
\draw[dashed] (1/2,1/2) -- (1,1);
\draw[dashed] (1/2,2) -- (1,1);
\draw[dashed] (2,1/2) -- (1,1);
\node at (1.25,0.25) {\(\gamma\)};
\node at (1.25,1.25) {\(\alpha\)};
\node at (0.25,1.25) {\(\beta\)};
\end{tikzpicture}
\caption{Changing branch cuts, reprise.}
\label{fig:tripod1}
\end{center}
\end{figure}

In Figure \ref{fig:tripod2}, we draw the image of the
developing map on the complement of \(T\) in such a way that
the almost toric boundary unwraps as a horizontal line. The
image of the developing map is a strip (closed at the bottom,
open at the top). We have shaded some fundamental action
domains alternately light and dark. The fundamental group of
the complement of \(T\) is \(\ZZ\), which acts on the strip by
powers of \(\lmatrix 1 & 0 \\ 9 & 1\rmatrix \) (treating \(O\)
as the origin). In each translate of the fundamental action
domain, you can ``see'' the hole left by excising \(T\). For
example, take the central light-coloured fundamental action
domain. This is obtained by taking the region labelled
\(\gamma\) in Figure \ref{fig:tripod1} and appending the
images of regions \(\alpha\) and \(\beta\) under the
monodromies\index{affine monodromy!around boundary of
CP2@around boundary of $\mathbb{CP}^2$|(} around the nodes
\(B\) and \(A\) respectively so that the almost toric boundary
becomes straight.

\begin{figure}[htb]
\begin{center}
\begin{tikzpicture}
\clip (-4,-1) rectangle (10,3);
\draw[thick] (-4,0) -- (10,0);
\node (O) at (1,1) {};
\node at (O) [above] {\(O\)};
\foreach \x in {-3,...,3}
{\filldraw[gray,opacity=0.5] (1,1) -- (-12+18*\x,0) -- (-3+18*\x,0) -- cycle;}
\foreach \y in {-3,...,3}
{\filldraw[lightgray,opacity=0.5] (1,1) -- (-3+18*\y,0) -- (6+18*\y,0) -- cycle;}
\foreach \x in {-3,...,3}
{\node at (1/2-9/2+9*\x,1/2) {\(\times\)};
\node at (2-9/2+9*\x,1/2) {\(\times\)};
\node at (-1-9/2+9*\x,1/2) {\(\times\)};
\node at (1/2-9/2+9*\x,1/2) [below] {\(A\)};
\node at (2-9/2+9*\x,1/2) [below] {\(B\)};
\node at (-1-9/2+9*\x,1/2) [below] {\(C\)};
\draw[dashed,thick] (1/2-9/2+9*\x,1/2) -- (O.center);
\draw[dashed,thick] (2-9/2+9*\x,1/2) -- (O.center);
\draw[dashed,thick] (-1-9/2+9*\x,1/2) -- (O.center);};
\foreach \y in {-3,...,3}
{\node at (1/2+9*\y,1/2) {\(\times\)};
\node at (2+9*\y,1/2) {\(\times\)};
\node at (-1+9*\y,1/2) {\(\times\)};
\node at (1/2+9*\y,1/2) [below] {\(A\)};
\node at (2+9*\y,1/2) [below] {\(B\)};
\node at (-1+9*\y,1/2) [below] {\(C\)};
\draw[dashed,thick] (1/2+9*\y,1/2) -- (O.center);
\draw[dashed,thick] (2+9*\y,1/2) -- (O.center);
\draw[dashed,thick] (-1+9*\y,1/2) -- (O.center);}
\node at (1.25,1/2) [below] {\(\gamma\)};
\node at (2.75,1/2) [below] {\(\alpha\)};
\node at (-0.25,1/2) [below] {\(\beta\)};
\end{tikzpicture}
\caption{The image of the complement of \(T\) under the developing map for the integral affine structure; alternately shaded regions are fundamental action domains, tiling the strip. The dashed lines show where \(T\) has been excised.}
\label{fig:tripod2}
\end{center}
\end{figure}
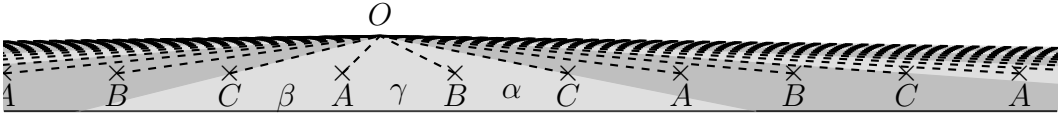

\end{Example}
\begin{Remark}[Exercise \ref{exr:boundary_monodromy}]\label{rmk:boundary_monodromy}
The significance of the matrix \(\lmatrix 1 & 0 \\ 9 &
1\rmatrix \) is that it is the total monodromy (anticlockwise)
around the boundary loop in Figure \ref{fig:tripod1},
considered as starting in region \(\gamma\). Note that \(9\)
is also the self-intersection of the almost toric boundary
curve. This is not a coincidence: compare with Example
\ref{exm:atbd_torus_nbhd}.

\end{Remark}
We will now be more explicit about monodromies. Let \[M_A =
\lmatrix 2 & 1 \\ -1 & 0 \rmatrix ,\quad M_B = \lmatrix -1 & 1
\\ -4 & 3\rmatrix ,\quad M_C = \lmatrix -1 & 4 \\ -1 & 3
\rmatrix.\] be the anticlockwise monodromies around the nodes
\(A,B,C\) in Figure \ref{fig:tripod1} (we can calculate these
using Lemma \ref{lma:monodromy_formula}). Let
\(T^o:=T\setminus\{0\}\) and call the three components ``branch
cuts'' (for now we ignore the interesting point \(0\)). A path
which crosses a branch cut lifts to a path in the strip but, as
usual, the affine structure is twisted by a monodromy matrix
when you cross the branch cut. In Figure
\ref{fig:tripod_monodromies}, we label the branch cuts by their
monodromies; these are obtained by conjugating the matrices
\(M_A,M_B,M_C\). For example, the nodes in the region \(\gamma\)
are identical to \(A\) and \(B\) in Figure \ref{fig:tripod1}, so
they have monodromies \(M_A\) and \(M_B\). The node just to the
right of \(B\) in Figure \ref{fig:tripod_monodromies} is related
to node \(C\) in Figure \ref{fig:tripod1} by crossing the branch
cut \(B\) anticlockwise, so its monodromy is
\(M_BM_CM_B^{-1}\). The next node to the right is related to
node \(A\) by crossing the branch cuts \(B\) and then \(C\)
anticlockwise, so its monodromy is \(M_B M_C M_A M_C^{-1}
M_B^{-1}\).\index{affine monodromy!around boundary of CP2@around
boundary of $\mathbb{CP}^2$|)}

\begin{figure}[htb]
\begin{center}
\begin{tikzpicture}
\clip (-4,-1) rectangle (10,3);
\node (O) at (1,1) {};
\foreach \x in {-3,...,3}
{\node at (1/2-9/2+9*\x,1/2) {\(\bullet\)};
\node at (2-9/2+9*\x,1/2) {\(\bullet\)};
\node at (-1-9/2+9*\x,1/2) {\(\bullet\)};
\draw[dashed,thick] (1/2-9/2+9*\x,1/2) -- (O.center);
\draw[dashed,thick] (2-9/2+9*\x,1/2) -- (O.center);
\draw[dashed,thick] (-1-9/2+9*\x,1/2) -- (O.center);};
\foreach \y in {-3,...,3}
{\node at (1/2+9*\y,1/2) {\(\bullet\)};
\node at (2+9*\y,1/2) {\(\bullet\)};
\node at (-1+9*\y,1/2) {\(\bullet\)};
\draw[dashed,thick] (1/2+9*\y,1/2) -- (O.center);
\draw[dashed,thick] (2+9*\y,1/2) -- (O.center);
\draw[dashed,thick] (-1+9*\y,1/2) -- (O.center);}
\node at (1/2,0) {{\tiny\(M_A\)}};
\node at (2,0) {{\tiny\(M_B\)}};
\node at (-1,0) {{\tiny\(M_A^{-1}M_CM_A\)}};
\node at (1/2+9/2,0) {{\tiny\(M_BM_CM_AM_C^{-1}M_B^{-1}\)}};
\node at (2+9/2,-0.5) {{\tiny\(M_BM_CM_AM_BM_A^{-1}M_C^{-1}M_B^{-1}\)}};
\draw (2+9/2,-0.25) -- (2+9/2,0.25);
\node at (-1+9/2,-0.5) {{\tiny\(M_BM_CM_B^{-1}\)}};
\draw (-1+9/2,-0.25) -- (-1+9/2,0.25);
\node at (-1+9,0) {\(\cdots\)};
\node at (1-9/2,0) {\(\cdots\)};
\node at (2-9/2,-0.5) {{\tiny\(M_A^{-1}M_C^{-1}M_BM_CM_A\)}};
\draw (2-9/2,-0.25) -- (2-9/2,0.25);
\node at (2-9/2,1/2) {{\(\circlearrowleft\)}};
\node at (1/2,1/2) {{\(\circlearrowleft\)}};
\node at (2,1/2) {{\(\circlearrowleft\)}};
\node at (-1,1/2) {{\(\circlearrowleft\)}};
\node at (1/2+9/2,1/2) {{\(\circlearrowleft\)}};
\node at (2+9/2,1/2) {{\(\circlearrowleft\)}};
\node at (-1+9/2,1/2) {{\(\circlearrowleft\)}};
\end{tikzpicture}
\caption{Monodromies for Figure \ref{fig:tripod2} (all acting from the right).}
\label{fig:tripod_monodromies}
\end{center}
\end{figure}
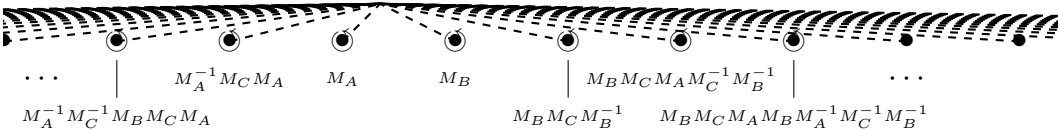

We now try to understand what is happening near the point
\(O\). If we look at Figure \ref{fig:tripod1}, the point \(O\)
looks less intimidating: it is just a point in \(B^{reg}\) away
from the branch cuts. The integral affine
structure\index{integral affine!structure!extend over cone
point} around \(O\) is easy to understand: it is just an open
ball in \(\RR^2\). The only reason it looks so interesting in
Figure \ref{fig:tripod2} is that the three legs of \(T\) pass
through \(O\). It becomes important to ask: if we have a diagram
like Figure \ref{fig:tripod2}, how can we tell if the point
\(O\) is a perfectly innocuous point in disguise? To answer this
question, we first introduce the language of integral affine
cones.\index{branch cut!more exotic|)}

\section{Integral affine cones}

\begin{figure}[htb]
\begin{center}
\begin{tikzpicture}
\filldraw[fill=lightgray,opacity=0.5,draw=none] (0,3) -- (0,0) -- (6,3) -- cycle;
\draw[dotted,thick,->-] (0,0) -- (0,3);
\draw[dotted,thick,->-] (0,0) -- (6,3);
\node at (0,1.5) [left] {\(\ell\)};
\node at (3,1.5) [below right] {\(\ell M\)};
\node at (1,1.5) {\(W_{M,\ell}\)};
\begin{scope}[shift={(10,3)}]
\draw (-3,0) -- (0,-3) -- (3,0);
\filldraw[fill=white,draw=black,dotted] (0,0) circle [x radius = 2.81,y radius = 1];
\draw (2.64,-0.35) arc [start angle = -20, end angle = -160, x radius = 2.81, y radius = 1];
\draw[dashed] (0,-3) -- (2.81*0.5,-3+2.13397459622);
\draw[->] (1,1.2-3) -- (0.2,1.2-3) node [left] {\(M^{-1}\)};
\end{scope}
\end{tikzpicture}
\caption{The cone \(B_{M,\ell}\) obtained by ``wrapping up'' \(W_{M,\ell}\) using \(M\) to glue across the (dashed) branch cut.}
\label{fig:cone}
\end{center}
\end{figure}
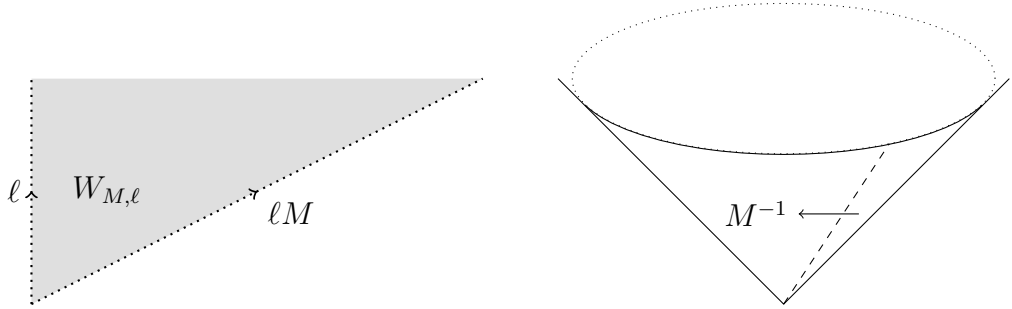

\begin{Definition}
Let \(M\in SL(2,\ZZ)\) and let \(\langle M\rangle\) be the
subgroup of \(SL(2,\ZZ)\) generated by \(M\). Let
\(\ell\subset\RR^2\) be a ray emanating from the
origin. Consider the wedges \(W_{M,\ell}\subset\RR^2\)
(respectively \(W'_{M,\ell}\subset\RR^2\)) which are swept
out\footnote{If \(\ell M=\ell\) then both wedges are the whole
of \(\RR^2\).} as \(\ell M\neq \ell\) rotates anticlockwise
(respectively clockwise) back to \(\ell\). We can equip the
quotient space \(B_{M,\ell} = W_{M,\ell}\setminus\{0\}/\langle
M\rangle\) (respectively \(B'_{M,\ell}\)) with the structure
of an integral affine manifold, by identifying \(x\in \ell M\)
with \(xM^{-1}\in \ell\) at the level of points and \(v\in
T_x\RR^2\) with \(vM^{-1}\in T_{xM^{-1}}\RR^2\) at the level
of tangent vectors. We call such an integral affine manifold a
{\em punctured cone}; we obtain a singular integral affine
manifold\index{integral affine!cone} called a {\em cone} by
adding in the cone point \(0\). The matrix \(M\) is the {\em
affine monodromy around the cone point}.\index{affine
monodromy!around cone point}

\end{Definition}
\begin{Remark}
We can visualise \(B_{M,\ell}\) as a cone obtained by wrapping
\(W_{M,\ell}\) up so that \(\ell\) and \(\ell M\) are
identified by the map \(x\mapsto xM\); see Figure
\ref{fig:cone}.

\end{Remark}

\begin{Remark}
Note that \(W'_{M,\ell}=W_{M^{-1},\ell M}\) and
\(B'_{M,\ell}=B_{M^{-1},\ell M}\).

\end{Remark}
\begin{Exercise}\label{exr:paper_model}
Make a paper model of the integral affine manifold
\(B_{M,\ell}\) for \[M = \lmatrix 0 & 1 \\ -1 & 0 \rmatrix \]
and \(\ell\) the positive \(x\)-axis. Is it possible to make
paper models for \(B_{M,\ell}\) when \(M\) is not conjugate to
this matrix?

\end{Exercise}
\begin{Lemma}\label{lma:integral_affine_smooth_cone_point}
The integral affine structure\index{integral
affine!structure!extend over cone} on the punctured cone
extends over the cone point when \(M=I\).
\end{Lemma}
\begin{Proof}
In this case, \(W_{M,\ell}=\RR^2\setminus\ell\) and we obtain
the cone by identifying both sides of the cut using the
identity, so the cone is just \(\RR^2\).\qedhere

\end{Proof}
\section{Back to the example}

\begin{figure}[htb]
\begin{center}
\begin{tikzpicture}
\draw (0,0) -- (3,0) -- (0,3) -- cycle;
\draw[dotted,thick] (0,0) -- (1/2,1/2);
\draw[dotted,thick] (3,0) -- (2,1/2);
\draw[dotted,thick] (0,3) -- (1/2,2);
\node (A) at (1/2,1/2) {\(\bullet\)};
\node (B) at (2,1/2) {\(\bullet\)};
\node (C) at (1/2,2) {\(\bullet\)};
\node[blue] at (1,0.65) {{\tiny\(\bullet\)}};
\draw[dashed] (1/2,1/2) -- (1,1);
\draw[dashed] (1/2,2) -- (1,1);
\draw[dashed] (2,1/2) -- (1,1);
\draw[blue,thick,->-] (0.65,0.65) -- (1.7,0.65) -- (0.65,1.7)-- cycle;
\draw[->,thick] (1,4) -- (1,2.5);
\begin{scope}[shift={(0,4.5)}]
\clip (-4,-1) rectangle (10,3);
\draw (-4,0) -- (10,0);
\node (O) at (1,1) {};
\node at (O) [above] {\(O\)};
\foreach \x in {-3,...,3}
{\node at (1/2-9/2+9*\x,1/2) {\(\times\)};
\node at (2-9/2+9*\x,1/2) {\(\times\)};
\node at (-1-9/2+9*\x,1/2) {\(\times\)};
\draw[dashed,thick] (1/2-9/2+9*\x,1/2) -- (O.center);
\draw[dashed,thick] (2-9/2+9*\x,1/2) -- (O.center);
\draw[dashed,thick] (-1-9/2+9*\x,1/2) -- (O.center);};
\foreach \y in {-3,...,3}
{\node at (1/2+9*\y,1/2) {\(\times\)};
\node at (2+9*\y,1/2) {\(\times\)};
\node at (-1+9*\y,1/2) {\(\times\)};
\draw[dashed] (1/2+9*\y,1/2) -- (O.center);
\draw[dashed] (2+9*\y,1/2) -- (O.center);
\draw[dashed] (-1+9*\y,1/2) -- (O.center);}
\draw[blue,thick,->-] (-4,0.65) -- (10,0.65);
\node[blue] at (1,0.65) {{\tiny\(\bullet\)}};
\node[blue] at (1+9*0.35,0.65) {{\tiny\(\bullet\)}};
\end{scope}
\end{tikzpicture}
\caption{Below: A triangular path with trivial monodromy. Above: Its lift to the strip.}
\label{fig:tripod_blue_line}
\end{center}
\end{figure}
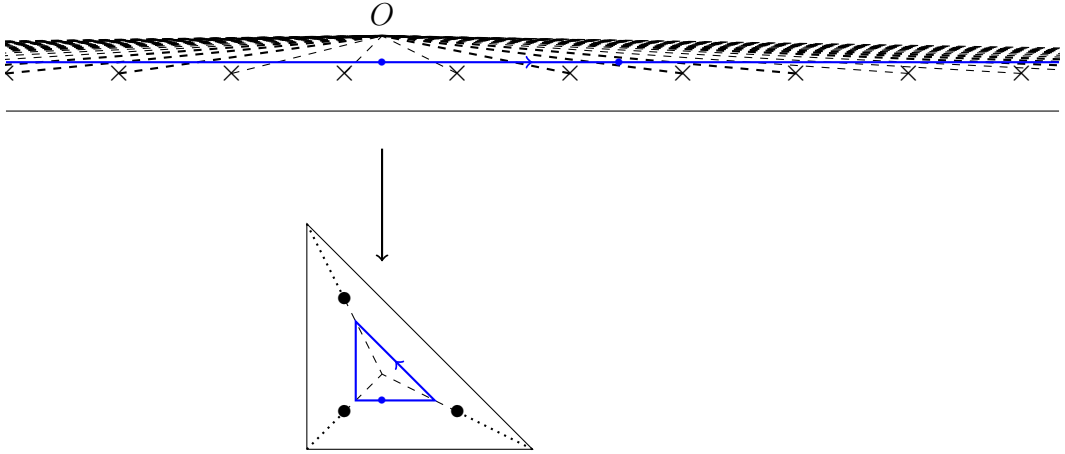

Now take the horizontal line shown crossing the branch cuts in
the top part of Figure \ref{fig:tripod_blue_line}. The portion
of the almost toric base lying above this line is an integral
affine cone\index{integral affine!cone}: it is made up of
infinitely many triangular segments which are glued together
using the affine monodromies. To check that the integral affine structure extends
over the cone point \(O\), it suffices to check that the total
monodromy as we traverse the boundary of this cone is the
identity.

This horizontal line projects to a triangle in the lower part of
Figure \ref{fig:tripod_blue_line} (with corners where it hits
the branch cuts\index{branch cut!more exotic|(}). Since this
triangle does not cross any of the dotted branch cuts, the
monodromy around it is the identity, but we can perform the same
calculation ``upstairs'' by multiplying the clockwise
monodromies of the three branch cuts crossed by the horizontal
line, and then multiplying by \(\lmatrix 1 & 0 \\ 9 & 1 \rmatrix
\) to send the lifted end-point back to the lifted
start-point:\index{branch cut!more exotic|)}
\begin{gather*}
M_B^{-1}\cdot M_BM_C^{-1}M_B^{-1} \cdot
M_BM_CM_A^{-1}M_C^{-1}M_B^{-1} \cdot\lmatrix 1 & 9 \\ 0 & 1
\rmatrix
\\= \cdot
M_A^{-1}M_C^{-1}M_B^{-1}\cdot \lmatrix 1 & 9 \\ 0 & 1 \rmatrix
\\=\lmatrix 1 & 0 \\ 0 &
1\rmatrix .\end{gather*}

\section{Developing map for cones}

We now discuss the developing map\index{developing map!for
cones} for the integral affine structure, and the dependence of
\(B_{M,\ell}\) on \(\ell\). Let \(\tilde{B}_{M,\ell}\) be the
universal cover of \(B_{M,\ell}\). This can be constructed as
follows. Take infinitely many copies \(W_i\), \(i\in\ZZ\), of
\(W_{M,\ell}\) and write \(\partial_\ell W_i\) and
\(\partial_{\ell M} W_i\) for the two boundary rays. Let
\(\tilde{B}_{M,\ell}\) be the quotient of
\(\coprod_{i\in\ZZ}W_i\) by identifying \(x\in\partial_\ell
W_i\) with \(xM\in\partial_{\ell M}W_{i-1}\) for all
\(i\in\ZZ\). The image of \(W_i\) under the developing map for
the integral affine structure on \(\tilde{B}_{M,\ell}\) is then
\(M^iW_{M,\ell}\).

\begin{Lemma}\label{lma:ell_dependence}
If \(\ell'\subset M^iW_{M,\ell}\) for some \(i\) then
\(B_{M,\ell}\) and \(B_{M,\ell'}\) are isomorphic as integral
affine manifolds.
\end{Lemma}
\begin{Proof}
If \(\ell'=\ell M^i\) then the map \(M^i\colon W_{M,\ell}\to
W_{M,\ell'}\) descends to give an isomorphism \(B_{M,\ell}\to
B_{M,\ell'}\). Using this, we may assume that \(\ell'\subset
W_{M,\ell}\).

\begin{figure}[htb]
\begin{center}
\begin{tikzpicture}
\draw (0,0) -- (0,3) node [above] {\(\ell\)};
\draw (0,0) -- (3,3) node [above] {\(\ell'\)};
\draw (0,0) -- (6,3) node [above] {\(\ell M\)};
\draw (0,0) -- (9,3) node [above] {\(\ell' M\)};
\node at (1,2) {\(S\)};
\node at (3,2) {\(T\)};
\node at (5,2) {\(SM\)};
\node at (3,3.7) [above] {\(W_{M,\ell}\)};
\node at (6,4.2) [above] {\(W_{M,\ell'}\)};
\draw[decorate,decoration={brace,amplitude=4pt}] (0,3.7) -- (6,3.7);
\draw[decorate,decoration={brace,amplitude=4pt}] (3,4.2) -- (9,4.2);
\end{tikzpicture}
\caption[Constructing an isomorphism \(B_{M,\ell}\to B_{M,\ell'}\).]{Constructing an isomorphism \(B_{M,\ell}\to B_{M,\ell'}\). In this example, \(M=\lmatrix 1 & 0 \\ 2 & 1 \rmatrix \).}
\label{fig:lma_ell_proof}
\end{center}
\end{figure}

Let \(S\) be the sector in \(W_{M,\ell}\) swept out as
\(\ell'\) moves anticlockwise to \(\ell\) and \(T\) the other
sector. We see that \(W_{M,\ell'}=T\cup SM\). The piecewise
linear map \[W_{M,\ell}\to
W_{M,\ell'},\qquad x\mapsto \begin{cases}xM\mbox{ if }x\in
S\\ x\mbox{ if }x\in T,\end{cases}\] descends to give the
desired isomorphism \(B_{M,\ell}\to B_{M,\ell'}\). \qedhere

\end{Proof}
\begin{Lemma}\label{lma:symmetry_ell_dependence}
If \(K\in GL(2,\ZZ)\) then \(B_{K^{-1}MK,\ell K}\) is
isomorphic to \(B_{M,\ell}\) via the map \(K\). In particular,
if \(K=-I\) then we see that \(B_{M,\ell}=B_{M,-\ell}\).
\end{Lemma}
\begin{Proof}
If we consider \(K\) as a change of coordinates then, in the
new coordinates, \(M\) is represented by \(K^{-1}MK\) and
\(\ell\) is sent to \(\ell K\). In particular, \(W_{M,\ell}\)
is sent to \(W_{K^{-1}MK,\ell K}\), and the recipe for gluing
together \(\tilde{B}_{M,\ell}\) transforms into the recipe for
gluing together \(\tilde{B}_{K^{-1}MK,\ell K}\). \qedhere

\end{Proof}
\section{Examples}

We will focus on examples where the ray \(\ell\) is not an
eigenray of \(M\). We will also focus on examples where the
wedge \(W_{M,\ell}\) or \(W'_{M,\ell}\) subtends an angle less
then \(\pi\) radians.

\begin{Example}\label{exm:parabolic}
Consider \(M=\lmatrix 1 & 0 \\ n & 1 \rmatrix \) for
some integer \(n>0\). If we take \(\ell\) to be a ray pointing
vertically up, then \(W_{M,\ell}\) is shown in Figure
\ref{fig:parabolic} for \(n=3\). The image of
\(\tilde{B}_{M,\ell}\) is the strict upper half-plane, which
is tiled by the domains \(W_{M,\ell}M^i\).

\begin{figure}[htb]
\begin{center}
\begin{tikzpicture}
\clip (-4,-1) rectangle (10,3);
\foreach \x in {0,...,20}
{\draw[opacity=(20-\x)/20] (0,0) -- (0+9*\x,3);}
\foreach \x in {-20,...,-1}
{\draw[opacity=(20+\x)/20] (0,0) -- (0+9*\x,3);}
\filldraw[fill=lightgray,opacity=0.5,draw=none] (0,0) -- (9,3) -- (0,3) -- cycle;
\node at (2,2) {\(W_{M,\ell}\)};
\node at (-2,2) {\(W_{M,\ell}M^{-1}\)};
\node at (8,2) {\(W_{M,\ell}M\)};
\end{tikzpicture}
\caption[The image of \(\tilde{B}_{M,\ell}\) under the developing map.]{The image of \(\tilde{B}_{M,\ell}\) under the developing map for \(M=\lmatrix 1 & 0 \\ 2 & 1\rmatrix \), tiled by domains \(W_{M,\ell}M^i\).}
\label{fig:parabolic}
\end{center}
\end{figure}
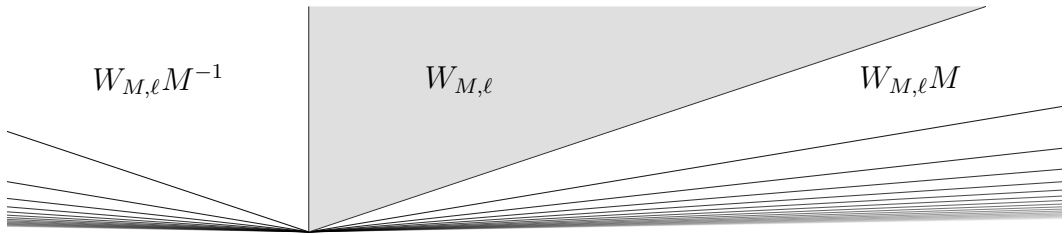

Any choice of \(\ell\) in the strict upper half-plane will
give the same \(B_{M,\ell}\) by Lemma
\ref{lma:ell_dependence}. If we take \(\ell\) in the lower
half-plane then we get the same integral affine manifold
\(B_{M,\ell}\) by Lemma \ref{lma:symmetry_ell_dependence}.

\end{Example}
\begin{Remark}
The only other option is to take \(\ell\) to be horizontal (an
eigenray). If \(n=1\), the resulting \(B_{M,\ell}\) is the
integral affine base for a single focus-focus fibre.

\end{Remark}
\begin{Example}\label{exm:cat_map}
Consider\footnote{For 2-by-2 matrix connaisseurs, this is the
Arnold cat map.} \(M=\lmatrix 2 & 1 \\ 1 &
1\rmatrix \). This has eigenvalues
\(\frac{3\pm\sqrt{5}}{2}\) and eigenvectors \(v_{\pm}=(x,y)\)
with \(y=-\frac{1}{2}(1\mp\sqrt{5})x\). The eigendirections
divide the plane into four open quadrants. Let us call the
quadrants \(Q,Q',-Q,-Q'\) where \(Q\) and \(Q'\) have the
property that if \(\ell\subset Q\) (respectively \(\ell\subset
Q'\)) then \(W_{M,\ell}\subset Q\) (respectively
\(W'_{M,\ell}\subset Q'\)).

If we pick \(\ell\subset Q\) or \(\ell\subset -Q\)
(respectively \(\ell\subset Q'\) or \(\ell\subset -Q'\)) then
the images \(W_{M,\ell}M^i\) (respectively
\(W'_{M,\ell}M^i\)), \(i\in\ZZ\) tile the chosen quadrant. By
Lemma \ref{lma:ell_dependence}, this means that, other than
choosing an eigenray, the choice of \(\ell\) amounts to a
choice of quadrant. Moreover, rays from opposite quadrants
yield the same integral affine manifold by Lemma
\ref{lma:symmetry_ell_dependence}. Thus, there are two {\em
essentially different} choices of ray: \(\ell_1\subset Q\) or
\(\ell_2\subset Q'\); see Figure \ref{fig:hyperbolic}. The
same phenomenon occurs whenever \(M\) has two distinct
positive real eigenvalues (which is equivalent to
\(\OP{Tr}(M)\geq 3\)).

\begin{figure}[htb]
\begin{center}
\begin{tikzpicture}
\filldraw[fill=lightgray,opacity=0.5,draw=none] (0,3) -- (0,0) -- (3,3) -- cycle;
\filldraw[fill=lightgray,opacity=0.5,draw=none] (3,0) -- (0,0) -- (3,1.5) -- cycle;
\node at (0,3) [above] {\(\ell_1\)};
\node at (3,3) [above] {\(\ell_1 M\)};
\node at (3,0) [right] {\(\ell_2\)};
\node at (3,1.5) [right] {\(\ell_2 M\)};
\node at (1,2) {\(W_{M,\ell_1}\)};
\node at (2.3,0.5) {\(W'_{M,\ell_2}\)};
\node at (-1,3) [above] {\(Q\)};
\node at (3,-2) [right] {\(Q'\)};
\node at (-1.3,0.2) {\(-Q'\)};
\node at (-0.5,-1) {\(-Q\)};
\clip (-3,-3) rectangle (3,3);
\draw[very thick] (-3,3*1.61803398875) -- (3,-3*1.61803398875);
\draw[very thick] (-3,-3*0.61803398875) -- (3,3*0.61803398875);
\draw[opacity=0.5] (0,0) -- (-24,39);
\draw (0,0) -- (-9,15);
\draw (0,0) -- (-3,6);
\draw (0,0) -- (0,3);
\draw (0,0) -- (3,3);
\draw (0,0) -- (9,6);
\draw[opacity=0.5] (0,0) -- (24,15);
\draw[opacity=0.5] (0,0) -- (15,-24);
\draw (0,0) -- (6,-9);
\draw (0,0) -- (3,-3);
\draw (0,0) -- (3,0);
\draw (0,0) -- (6,3);
\draw (0,0) -- (15,9);
\draw[opacity=0.5] (0,0) -- (39,24);
\end{tikzpicture}
\caption{Some choices of rays \(\ell_1\) and \(\ell_2\) in Example \ref{exm:cat_map} which give different manifolds \(B_{M,\ell_1}\) and \(B'_{M,\ell_2}\). The thick lines are the eigenlines.}
\label{fig:hyperbolic}
\end{center}
\end{figure}
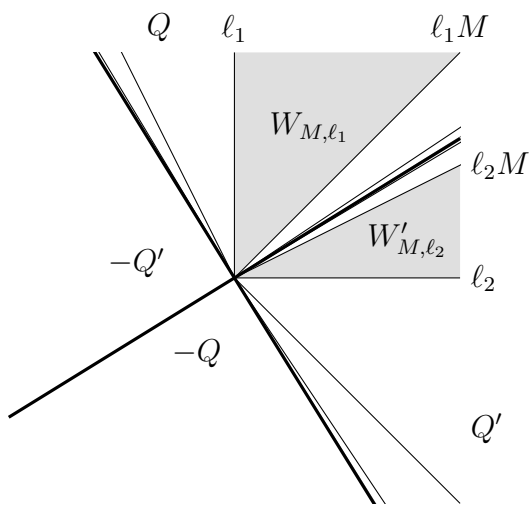

\end{Example}
\section{Symplectic manifolds}

\begin{Definition}
We can associate to \(B_{M,\ell}\) a symplectic manifold
\(X_{M,\ell}\) together with a Lagrangian torus fibration
\(X_{M,\ell}\to B_{M,\ell}\) compatible with the integral
affine structure. To construct \(X_{M,\ell}\), we take
\(W_{M,\ell}\times T^2\) with the symplectic form \(\sum
dp_i\wedge dq_i\) and identify \((\bm{p},\bm{q})\)
with\footnote{Compare with the map \(\Phi\) in the proof of
Lemma \ref{lma:globcoords}.} \((\bm{p}M^{-1},M\bm{q})\) for
\(\bm{p}\in \ell M\).

\end{Definition}
\begin{Lemma}[Exercise \ref{exr:find_path}]\label{lma:find_path}
There exists a smooth path \(\gamma\colon[0,1]\to W_{M,\ell}\)
which is transverse to all rays emanating from the origin such
that \(\gamma(0)\in\ell\), \(\gamma(1)=\gamma(0)M\), and
\(\frac{d^n\gamma}{dt^n}(1)=\frac{d^n\gamma}{dt^n}(0)M\) for
all \(n\geq 1\).

\end{Lemma}
\begin{Lemma}\label{lma:identify_b_m_ell}
Let \(Y_{M,\ell}\) be the quotient of \([0,1]\times T^2\)
which identifies \((0,\bm{q})\) with \((1,M^T\bm{q})\). This is a torus
bundle over the circle. The symplectic manifold \(X_{M,\ell}\)
is diffeomorphic to \((0,\infty)\times Y_{M,\ell}\).
\end{Lemma}
\begin{Proof}
Let \(\gamma\colon[0,1]\to W_{M,\ell}\) be a path given by
Lemma \ref{lma:find_path}. This descends to a closed smooth
loop in \(B_{M,\ell}\). Moreover, because \(\gamma\) is
transverse to all rays through the origin, we can foliate
\(B_{M,\ell}\) by loops \(s\gamma\) for
\(s\in(0,\infty)\). Define a map \((0,\infty)\times
[0,1]\times T^2\to W_{M,\ell}\times T^2\) by \((s,t,\bm{q})\mapsto
(s\gamma(t),\bm{q})\). This descends to give a diffeomorphism
\((0,\infty)\times Y_{M,\ell}\to X_{M,\ell}\).\qedhere

\end{Proof}
In fact, we can be more precise.

\begin{Lemma}\label{lma:symplectically_identify_b_m_ell}
\(X_{M,\ell}\) is symplectomorphic to \((0,\infty)\times
Y_{M,\ell}\) with the symplectic form \(\omega=d(s\alpha)\),
where \(s\in(0,\infty)\) and \(\alpha\) is a contact form on
\(Y_{M,\ell}\). This is called the {\em
symplectisation}\index{symplectisation} of the contact
manifold \(Y_{M,\ell}\).
\end{Lemma}
\begin{Proof}
Consider the map \((0,\infty)\times Y_{M,\ell}\to
X_{M,\ell}\), \((s,t,\bm{q})\mapsto (s\gamma(t),\bm{q})\) from
Lemma \ref{lma:identify_b_m_ell}. The symplectic
form\index{symplectic form!on symplectisation} on
\(X_{M,\ell}\) is \(\sum dp_i\wedge dq_i\) with
\(\bm{p}=s\gamma(t)\), so the pullback of this form to
\((0,\infty)\times Y_{M,\ell}\) is \(\sum_i
d(s\gamma_i(t))\wedge dq_i=d\left(\sum_i
s\gamma_i(t)\,dq_i\right)\). The 1-form \(\alpha\) on
\(Y_{M,\ell}\) is \(\sum_i \gamma_i(t)\,dq_i\), so we see the
pullback of \(\omega\) along our diffeomorphism is
\(d(s\alpha)\) on \((0,\infty)\times Y_{M,\ell}\). \qedhere

\end{Proof}
\section{Elliptic and cusp singularities}

The integral affine manifold \(B_{M,\ell}\) has a natural
partial compactification, \(\overline{B}_{M,\ell}:=
W_{M,\ell}/\langle M\rangle\), which adds in the cone point
\(\overline{b}\). If \(M\) is the identity matrix, we observed
in Lemma \ref{lma:integral_affine_smooth_cone_point} that the
integral affine structure on \(B_{M,\ell}\) extends over the
cone point, and so \(X_{M,\ell}\) naturally sits inside a larger
symplectic manifold \(\overline{X}_{M,\ell}\) which fibres over
\(\overline{B}_{M,\ell}\) with a regular torus fibre over
\(\overline{b}\).

If \(M\) is not the identity, we simply take
\(\overline{X}_{M,\ell}\) to be the partial compactification of
\(X_{M,\ell}\) which adds in a single point \(\overline{x}\)
over \(\overline{b}\). We will think of this as a ``singular
symplectic manifold'', i.e.\ a singular space with a symplectic
form defined away from the singularity. In fact, this singular
space can be equipped with the structure of a complex algebraic
variety, with a singularity at \(\overline{x}\) called an {\em
elliptic}\index{elliptic singularity|see {singularity,
elliptic}}\index{singularity!elliptic} or {\em cusp
singularity}\index{singularity!cusp} depending on the matrix
\(M\); see Hirzebruch's seminal paper {\cite[Section
2.2--3]{HirzebruchHilbert}}. The contact manifold \(Y_{M,\ell}\)
from Lemma \ref{lma:symplectically_identify_b_m_ell} is known as
the {\em link}\index{singularity!link} of the singularity.

The minimal resolution\index{singularity!resolution!minimal} of
this algebraic variety replaces \(\overline{x}\) with either a
smooth elliptic curve (in the elliptic case), a nodal elliptic
curve, or a cycle of rational curves\index{cycle of rational
curves}. We will discuss the symplectic version of this, where
we resolve the singularity by symplectic cuts.

\begin{Example}[Parabolic matrix]\label{exm:parabolic_ctd}
Let \(M=\lmatrix 1 & 0 \\ n & 1 \rmatrix \) and take \(\ell\)
to point in the \((0,1)\)-direction, as in Example
\ref{exm:parabolic}. If we perform a symplectic cut on
\(\overline{X}_{M,\ell}\) at some positive height then we
obtain the following integral affine base:

\begin{center}
\begin{tikzpicture}
\filldraw[fill=lightgray,opacity=0.5,draw=none] (0,4) -- (0,2) -- (2,2) -- (4,4) -- cycle;
\draw[->-,thick,dotted] (0,2) -- (0,4);
\draw[->-,thick,dotted] (2,2) -- (4,4);
\node at (0,3) [left] {\(\ell=(0,1)\)};
\node at (3,3) [right=0.3cm] {\((n,1)=\ell M\)};
\draw[very thick] (0,2) -- (2,2);
\end{tikzpicture}
\end{center}

where we are identifying the edges labelled with arrows using
the matrix \(M\). The corresponding symplectic manifold is
smooth, and lying over the compact edge we have a symplectic
torus with self-intersection \(-n\) (compare with Example
\ref{exm:atbd_torus_nbhd}). This corresponds to the case of an
elliptic singularity\index{singularity!elliptic}, where the
minimal resolution\index{singularity!resolution!elliptic}
introduces a smooth elliptic
curve.\index{toric boundary!resolution of elliptic
singularity}

\end{Example}
\begin{Example}[Hyperbolic matrix]\label{exm:cat_map_ctd}
Let\index{singularity!cusp|(} \(M=\lmatrix 2 & 1 \\ 1 & 1
\rmatrix \) and take \(\ell\) to point in the
\((0,1)\)-direction as in Example \ref{exm:cat_map}. If we
make a symplectic cut at some positive height then we obtain
the following integral affine base:

\begin{center}
\begin{tikzpicture}
\filldraw[fill=lightgray,opacity=0.5,draw=none] (0,4) -- (0,2) -- (2,2) -- (4,4) -- cycle;
\draw[->-,thick,dotted] (0,2) -- (0,4);
\draw[->-,thick,dotted] (2,2) -- (4,4);
\node at (0,3) [left] {\(\ell=(0,1)\)};
\node at (3,3) [right=0.3cm] {\((1,1)=\ell M\)};
\node at (0,2) {\(\bullet\)};
\node at (2,2) {\(\bullet\)};
\draw[very thick] (0,2) -- (2,2);
\end{tikzpicture}
\end{center}

where the edges with arrows are identified using \(M\). The
point\footnote{Because of the identifications, there is only
one point!} we have marked with a dot is actually a Delzant
vertex. To see this, imagine how the left- and right-pointing
edges look from the point of view of the left-most dot in the
diagram. The right-pointing edge points in the
\((1,0)\)-direction. The left-pointing edge points in the
\((-1,0)M^{-1}= (-1,1)\)-direction. We can make this clearer
by shifting the branch cut \(\ell\) to \(\ell'\) parallel to
\((-1,3)\):

\begin{center}
\begin{tikzpicture}
\filldraw[fill=lightgray,opacity=0.5,draw=none] (-5/3,5) -- (-1,3) -- (0,2) -- (1,2) -- (2,4) -- cycle;
\draw[->-,thick,dotted] (-1,3) -- (-5/3,5);
\draw[->-,thick,dotted] (1,2) -- (2,4);
\node at (-4/3,4) [left] {\(\ell'=(-1,2)\)};
\node at (1.5,3) [right=0.3cm] {\((1,2)=\ell' M\)};
\node at (0,2) {\(\bullet\)};
\draw[very thick] (-1,3) -- (0,2) -- (1,2);
\end{tikzpicture}
\end{center}

The toric boundary\index{toric boundary!resolution of cusp
singularity} is therefore a nodal elliptic curve. This
singularity is therefore a {\em cusp singularity}.

\end{Example}
If this first symplectic cut had not fully resolved our
singularity (had the corner not been Delzant), we could have
continued in the manner of Example \ref{exm:minimal_resolution},
making more cuts, until all vertices were Delzant. The result
would be a cycle of rational curves; this is the typical
behaviour when \(M\) is hyperbolic\footnote{i.e.\ has two
distinct real eigenvalues.}.

\begin{Lemma}
If these rational curves have
self-intersection
\(s_1,s_2,\ldots,s_k\) then the infinite periodic continued
fraction\index{continued fraction!infinite periodic|(}
\[s_1-\frac{1}{s_2- \frac{1}{\cdots - \frac{1}{s_k-
\frac{1}{\cdots}}}}\] converges to the slope of the
dominant\footnote{i.e.\ the
eigenline\index{eigenline!irrational slope|(} corresponding to
the largest positive eigenvalue.} eigenline of \(M\).
\end{Lemma}
\begin{Proof}
Suppose we make a sequence of
cuts to the fundamental
domain \(W_{M,\ell}\) to get a new Delzant
polygonal domain
\(\tilde{W}_{M,\ell}\). Using the action of \(\langle
M\rangle\), we get infinitely many translates of this domain.

If necessary, change the branch cut\index{branch cut!integral
affine cone} \(\ell\) to ensure that none of the vertices of
\(\tilde{W}_{M,\ell}\) are on the branch cut. Let \(v_k\) be
the primitive integer vectors pointing rightwards along the
left-most edge of \(\tilde{W}_{M,\ell}M^k\) By construction,
we have \(v_kM=v_{k+1}\). The argument from Example
\ref{exm:minimal_resolution} shows that \(v_kS=v_{k+1}\) with
\(S=\lmatrix 0 &-1 \\ 1 & s_k\rmatrix \cdots \lmatrix 0 &-1
\\ 1 & s_1\rmatrix \). Since \(v_k=v_0M^k\), we know that the
slope of \(v_k\) converges to the slope of the dominant
eigenline of \(M\) as \(k\to \infty\). But as in Example
\ref{exm:minimal_resolution}, the recursion \(v_kS=v_{k+1}\)
tells us that if the slope of \(v_k\) is \(\alpha_k\) then the
slope of \(v_{k+1}\) is \[s_1-\frac{1}{s_2-\frac{1}{\cdots
s_k-\frac{1}{\alpha_k}}}.\] Therefore the infinite periodic
continued fraction\index{continued fraction!infinite
periodic|)} with coefficients
\[s_1,\ldots,s_k,s_1,\ldots,s_k,\ldots\] converges to the
slope of the dominant eigenline\index{eigenline!irrational
slope|)} of \(M\).\qedhere

\end{Proof}
\begin{Remark}
Recall that for every hyperbolic matrix \(M\) there were two
possible integral affine manifolds \(B_{M,\ell}\) up to
isomorphism, depending on the choice of \(\ell\). The
corresponding cusp singularities are said to be {\em dual} to
one another. Dual pairs of cusps were the subject of a
long-standing conjecture of Looijenga\index{Looijenga
conjecture}, which was resolved by Gross, Hacking and Keel
\cite{GHK} using ideas from mirror symmetry. A different, more
direct, proof of this conjecture, which uses these almost
toric pictures in an essential way was given by Engel
\cite{Engel}.\index{singularity!cusp|)}

\end{Remark}
\section{K3 surfaces from fibre sum}

\begin{figure}[htb]
\begin{center}
\begin{tikzpicture}[scale=1.5]
\node at (-3,1) {(a)};
\filldraw[lightgray,opacity=0.5] (0,1) -- (6,0) -- (-3,0) -- cycle;
\node (O) at (0,1) {\(\bullet\)};
\node at (O) [above] {\(O\)};
\draw[thick,dotted,->-] (O.center) --++ (-3,-1);
\draw[thick,dotted,->-] (O.center) --++ (6,-1);
\draw[very thick,black] (-3,0) -- (6,0);
\draw[thick,dotted] (O.center) --++ (1/2,-1/2) node {\(\times\)};
\draw[thick,dotted] (O.center) --++ (2,-1/2) node {\(\times\)};
\draw[thick,dotted] (O.center) --++ (-1,-1/2) node {\(\times\)};
\begin{scope}[shift={(0,-2)}]
\node at (-3,1) {(b)};
\filldraw[lightgray,opacity=0.5] (0,1) -- (-3,0) --++ (1.5,0) --++ (0,0.3) --++ (0.3,-0.3) --++ (0.3,0) --++ (0,0.3) --++ (0.3,-0.3) --++ (0.3,0) --++ (0,0.3) --++ (0.3,-0.3) --++ (0.3,0) --++ (0,0.3) --++ (0.3,-0.3) --++ (0.3,0) --++ (0,0.3) --++ (0.3,-0.3) --++ (0.3,0) --++ (0,0.3) --++ (0.3,-0.3) --++ (0.3,0) --++ (0,0.3) --++ (0.3,-0.3) --++ (0.3,0) --++ (0,0.3) --++ (0.3,-0.3) --++ (0.3,0) --++ (0,0.3) --++ (0.3,-0.3) --++ (0.3,0) -- (6,0) -- cycle;
\draw[thick,dotted] (-3,0) --++ (1.5,0) --++ (0,0.3) node {\(\times\)} --++ (0.3,-0.3) --++ (0.3,0) --++ (0,0.3) node {\(\times\)} --++ (0.3,-0.3) --++ (0.3,0) --++ (0,0.3) node {\(\times\)} --++ (0.3,-0.3) --++ (0.3,0) --++ (0,0.3) node {\(\times\)} --++ (0.3,-0.3) --++ (0.3,0) --++ (0,0.3) node {\(\times\)} --++ (0.3,-0.3) --++ (0.3,0) --++ (0,0.3) node {\(\times\)} --++ (0.3,-0.3) --++ (0.3,0) --++ (0,0.3) node {\(\times\)} --++ (0.3,-0.3) --++ (0.3,0) --++ (0,0.3) node {\(\times\)} --++ (0.3,-0.3) --++ (0.3,0) --++ (0,0.3) node {\(\times\)} --++ (0.3,-0.3) --++ (0.3,0) -- (6,0);
\node (O) at (0,1) {\(\bullet\)};
\node at (O) [above] {\(O\)};
\draw[thick,dotted,->-] (O.center) --++ (-3,-1);
\draw[thick,dotted,->-] (O.center) --++ (6,-1);
\draw[very thick,black] (-3,0) -- (-1.5,0);
\draw[very thick,black] (-1.2,0) --++ (0.3,0);
\draw[very thick,black] (-0.6,0) --++ (0.3,0);
\draw[very thick,black] (0.0,0) --++ (0.3,0);
\draw[very thick,black] (0.6,0) --++ (0.3,0);
\draw[very thick,black] (1.2,0) --++ (0.3,0);
\draw[very thick,black] (1.8,0) --++ (0.3,0);
\draw[very thick,black] (2.4,0) --++ (0.3,0);
\draw[very thick,black] (3,0) --++ (0.3,0);
\draw[very thick,black] (3.6,0) -- (6,0);
\draw[thick,dotted] (O.center) --++ (1/2,-1/2) node {\(\times\)};
\draw[thick,dotted] (O.center) --++ (2,-1/2) node {\(\times\)};
\draw[thick,dotted] (O.center) --++ (-1,-1/2) node {\(\times\)};
\end{scope}
\begin{scope}[shift={(0,-4)}]
\node at (-3,1) {(c)};
\filldraw[lightgray,opacity=0.5] (0,1) -- (-3*0.7,0.3) -- (-3*0.7,0) -- (6*0.7,0) -- (6*0.7,0.3) -- cycle;
\node at (-1.5,0.3) {\(\times\)};
\node at (-0.9,0.3) {\(\times\)};
\node at (-0.3,0.3) {\(\times\)};
\node at (0.3,0.3) {\(\times\)};
\node at (0.9,0.3) {\(\times\)};
\node at (1.5,0.3) {\(\times\)};
\node at (2.1,0.3) {\(\times\)};
\node at (2.7,0.3) {\(\times\)};
\node at (3.3,0.3) {\(\times\)};
\node (O) at (0,1) {\(\bullet\)};
\node at (O) [above] {\(O\)};
\draw[thick,dotted,->-] (O.center) -- (-3*0.7,0.3);
\draw[thick,dotted,->>-] (-3*0.7,0.3) -- (-3*0.7,0);
\draw[thick,dotted,->-] (O.center) -- (6*0.7,0.3);
\draw[thick,dotted,->>-] (6*0.7,0.3) -- (6*0.7,0);
\draw[very thick,black] (-3*0.7,0) -- (6*0.7,0);
\draw[thick,dotted] (O.center) --++ (1/2,-1/2) node {\(\times\)};
\draw[thick,dotted] (O.center) --++ (2,-1/2) node {\(\times\)};
\draw[thick,dotted] (O.center) --++ (-1,-1/2) node {\(\times\)};
\draw[thick,dotted] (0.3,0.3) -- (6*0.7,0.3);
\draw[thick,dotted] (-0.3,0.3) -- (-3*0.7,0.3);
\end{scope}
\end{tikzpicture}
\caption[(a) The almost toric base diagram of \(\cp{2}\) from Example \ref{exm:cp2_elliptic}. (b) Perform nine non-toric blow-ups along the toric boundary. (c) Make the branch cuts horizontal.]{(a) The almost toric base diagram of \(\cp{2}\) from Example \ref{exm:cp2_elliptic}.
The dotted edges are identified using the matrix \(M = \lmatrix 1 & 0 \\ 9 & 1\rmatrix\).
(b) Perform nine non-toric blow-ups along the toric boundary. (c) Make the branch cuts
horizontal. The opposite edges are identified using \(M\) above the horizontal cut
and using the identity below the cut.}
\label{fig:rational_elliptic_unwrapped}
\end{center}
\end{figure}
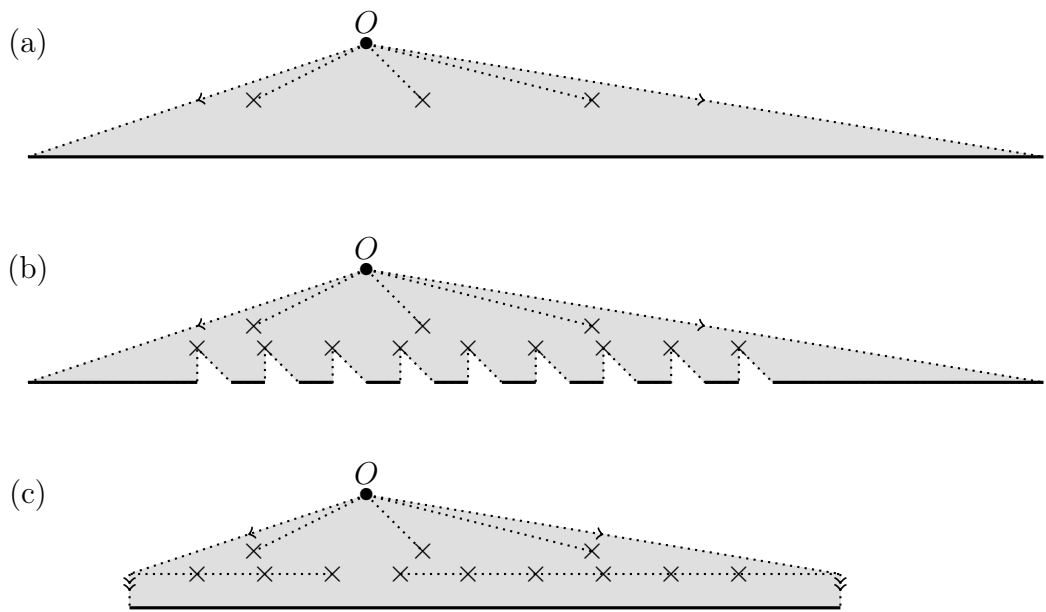

Take\index{almost toric base diagram!K3 surface|(} the almost
toric picture of \(\cp{2}\) from Example \ref{exm:cp2_elliptic};
we redraw a single fundamental action domain in Figure
\ref{fig:rational_elliptic_unwrapped}(a), with boundary
identifications indicated. The toric boundary is a symplectic
torus with self-intersection \(9\), as we can see by comparing
with Example \ref{exm:parabolic_ctd}. If we make \(9\) non-toric
blow-ups along the boundary as in Figure
\ref{fig:rational_elliptic} then our picture changes: see Figure
\ref{fig:rational_elliptic_unwrapped}(b). We can make a change
of branch cuts\index{branch cut!for K3 surface} to make all of
these branch cuts horizontal (Figure
\ref{fig:rational_elliptic_unwrapped}(c)); the opposite edges of
the fundamental action domain are identified using the identity
matrix or \(\lmatrix 1 & 0 \\ 9 & 1\rmatrix\) depending on
whether they are below (respectively above) the horizontal
branch cut. The associated manifold \(\cp{2}\sharp
9\overline{\CC\mathbb{P}}^2\) is called a {\em rational elliptic
surface}\index{elliptic surface!rational} and is often written
\(E(1)\)\index{En@$E(n)$|see {elliptic surface}} by
low-dimensional topologists; its toric boundary is a torus with
self-intersection \(0\).

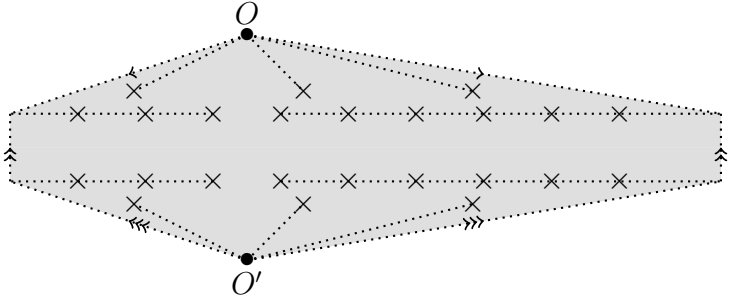
\begin{figure}[htb]
\begin{center}
\begin{tikzpicture}[scale=1.5]
\filldraw[lightgray,opacity=0.5,draw=none] (0,1) -- (-3*0.7,0.3) -- (-3*0.7,0) -- (6*0.7,0) -- (6*0.7,0.3) -- cycle;
\node at (-1.5,0.3) {\(\times\)};
\node at (-0.9,0.3) {\(\times\)};
\node at (-0.3,0.3) {\(\times\)};
\node at (0.3,0.3) {\(\times\)};
\node at (0.9,0.3) {\(\times\)};
\node at (1.5,0.3) {\(\times\)};
\node at (2.1,0.3) {\(\times\)};
\node at (2.7,0.3) {\(\times\)};
\node at (3.3,0.3) {\(\times\)};
\node (O) at (0,1) {\(\bullet\)};
\node at (O) [above] {\(O\)};
\draw[thick,dotted,->-] (O.center) -- (-3*0.7,0.3);
\draw[thick,dotted] (-3*0.7,0.3) -- (-3*0.7,0);
\draw[thick,dotted,->-] (O.center) -- (6*0.7,0.3);
\draw[thick,dotted] (6*0.7,0.3) -- (6*0.7,0);
\draw[thick,dotted] (O.center) --++ (1/2,-1/2) node {\(\times\)};
\draw[thick,dotted] (O.center) --++ (2,-1/2) node {\(\times\)};
\draw[thick,dotted] (O.center) --++ (-1,-1/2) node {\(\times\)};
\draw[thick,dotted] (0.3,0.3) -- (6*0.7,0.3);
\draw[thick,dotted] (-0.3,0.3) -- (-3*0.7,0.3);
\begin{scope}[xscale=1,yscale=-1]
\filldraw[lightgray,opacity=0.5,draw=none] (0,1) -- (-3*0.7,0.3) -- (-3*0.7,0) -- (6*0.7,0) -- (6*0.7,0.3) -- cycle;
\node at (-1.5,0.3) {\(\times\)};
\node at (-0.9,0.3) {\(\times\)};
\node at (-0.3,0.3) {\(\times\)};
\node at (0.3,0.3) {\(\times\)};
\node at (0.9,0.3) {\(\times\)};
\node at (1.5,0.3) {\(\times\)};
\node at (2.1,0.3) {\(\times\)};
\node at (2.7,0.3) {\(\times\)};
\node at (3.3,0.3) {\(\times\)};
\node (O) at (0,1) {\(\bullet\)};
\draw[thick,dotted,->>>-] (O.center) -- (-3*0.7,0.3);
\draw[thick,dotted,->>] (-3*0.7,0.3) -- (-3*0.7,0);
\draw[thick,dotted,->>>-] (O.center) -- (6*0.7,0.3);
\draw[thick,dotted,->>] (6*0.7,0.3) -- (6*0.7,0);
\draw[thick,dotted] (O.center) --++ (1/2,-1/2) node {\(\times\)};
\draw[thick,dotted] (O.center) --++ (2,-1/2) node {\(\times\)};
\draw[thick,dotted] (O.center) --++ (-1,-1/2) node {\(\times\)};
\draw[thick,dotted] (0.3,0.3) -- (6*0.7,0.3);
\draw[thick,dotted] (-0.3,0.3) -- (-3*0.7,0.3);
\end{scope}
\node at (0,-1) [below] {\(O'\)};
\end{tikzpicture}
\caption{An almost toric fibration on a K3 surface. Edges are identified in pairs as indicated by the arrows.}
\label{fig:k3}
\end{center}
\end{figure}

The almost toric picture in Figure \ref{fig:k3} is obtained by
reflecting Figure \ref{fig:rational_elliptic_unwrapped}(c)
horizontally and ignoring the toric boundary. The associated
almost toric manifold is obtained from a pair of rational
elliptic surfaces by performing a {\em Gompf sum} \cite{Gompf}
on the square-zero tori: in other words, we have excised a
neighbourhood of the toric boundary in each copy of \(E(1)\) and
glued the results together along their common boundary
\(T^3\). This construction yields the elliptic surface \(E(2)\),
otherwise known as a K3 surface\index{K3 surface}. The integral
affine base is a sphere with 24 base-nodes: we have drawn a
cylinder in Figure \ref{fig:k3}, and by Lemma
\ref{lma:integral_affine_smooth_cone_point}, the integral affine
structure extends over the sphere we get by adding in the points
\(O\) and \(O'\) to the cylinder.

\index{almost toric base diagram!K3 surface|)}

\section{Solutions to inline exercises}

\begin{Exercise}[Remark \ref{rmk:boundary_monodromy}]\label{exr:boundary_monodromy}
Show that the matrix \(\lmatrix 1 & 0 \\ 9 & 1\rmatrix \) is
the total monodromy (anticlockwise) around the boundary loop
in Figure \ref{fig:tripod1} starting at a point in the region
\(\gamma\).
\end{Exercise}
\begin{Solution}
By Lemma \ref{lma:monodromy_formula}, if the branch cut at a
base-node points in the \((p,q)\)-direction then the clockwise
monodromy around that node is \[\lmatrix 1-pq & -q^2\\ p^2 &
1+pq\rmatrix.\] The branch cuts are in the directions
\((-1,-1),\quad (2,-1),\quad (-1,2)\) at \(A\), \(B\), and
\(C\) respectively. This gives anticlockwise monodromy
matrices \[M_A = \lmatrix 2 & 1 \\ -1 & 0 \rmatrix ,\quad M_B
= \lmatrix -1 & 1 \\ -4 & 3\rmatrix ,\quad M_C = \lmatrix -1 &
4 \\ -1 & 3 \rmatrix.\] If we start with a vector \(v\) in the
region \(\gamma\) then it crosses the branch cuts emanating
from \(B\), \(C\) and \(A\) in that order, so the monodromy is
\[M_B M_C M_A = \lmatrix 1 & 0 \\ 9 & 1\rmatrix .\] Note that
if you start with a vector in the region \(\beta\) then you
end up with the matrix \(M_A M_B M_C\), which is different
(though conjugate by \(M_A\)).\qedhere

\end{Solution}
\begin{Exercise}
Make a paper model of the integral affine manifold
\(B_{M,\ell}\) for \[M = \lmatrix 0 & 1 \\ -1 & 0 \rmatrix \]
and \(\ell\) the positive \(x\)-axis. Is it possible to make
paper models for \(B_{M,\ell}\) when \(M\) is not conjugate to
this matrix?
\end{Exercise}
\begin{Solution}
You can only make a paper model of the cone if \(M\) is an
orthogonal (distance-preserving) map as well as being
\(\ZZ\)-linear. The group \(O(2)\cap SL(2,\ZZ)\) consists of
the four matrices \(\begin{pmatrix} \pm 1 & 0 \\ 0 & \pm
1\end{pmatrix},\quad \begin{pmatrix} 0 & \mp 1 \\ \pm 1 &
0\end{pmatrix}\). These four are all possible to
make. \qedhere

\end{Solution}
\begin{Exercise}[Lemma \ref{lma:find_path}]\label{exr:find_path}
There exists a smooth path \(\gamma\colon[0,1]\to W_{M,\ell}\)
which is transverse to all rays emanating from the origin such
that \(\gamma(0)\in\ell\), \(\gamma(1)=\gamma(0)M\), and
\(\frac{d^n\gamma}{dt^n}(1)=\frac{d^n\gamma}{dt^n}(0)M\) for
all \(n\geq 1\).
\end{Exercise}
\begin{Solution}
Let \(\Theta\) be the angle between \(\ell\) and \(\ell M\)
subtended by \(W_{M,\ell}\). Pick a segment of
\(\gamma\colon(-\epsilon,\epsilon)\to\RR^2\) in a
neighbourhood of \(\ell\) so that \(\gamma(0)\in\ell\) and so
that \(\gamma\) is transverse to all rays emanating from the
origin. In fact, we can assume by reparametrising \(\gamma\)
that the ray through \(\gamma(t)\) makes an angle \(\Theta t\)
with \(\ell\) for all \(t\in(-\epsilon,\epsilon)\). Apply
\(M\) to get a path-segment \(\delta:=\gamma M\) passing
through \(\ell M\). Again, by reparametrising, we can assume
that \(\delta(1)\in \ell M\) and \(\delta(1+t)\) makes an
angle \(\Theta t\) with \(\ell M\) for all
\(t\in(-\epsilon,\epsilon)\). We write these two segments in
polar coordinates \((r,\theta)\) as \((r(t),\theta_0+\Theta
t)\) where \(\theta_0\) is the argument of the ray
\(\ell\). We can now extend \(r\) in a completely arbitrary
(smooth) way to the interval \([0,1]\) and the graph will give
a path with the desired properties. \qedhere

\end{Solution}
\appendix
\part{Appendices}

\chapter{Symplectic linear algebra}
\label{ch:sla}
\thispagestyle{cup}

\section{Symplectic vector spaces}

\begin{Definition}\label{dfn:symplectic_vector_space}
Let \(V\) be a finite-dimensional vector space over \(\RR\)
and \(\omega\) be a bilinear map \(V\times V\to\RR\). We say
that \(\omega\) is a {\em linear symplectic
2-form}\index{symplectic form!linear} if \(\omega(v,v)=0\) for
all \(v\in V\) and for all nonzero \(v\in V\) there exists
\(w\in V\) such that \(\omega(v,w)\neq 0\) (nondegeneracy).
We say that the pair \((V,\omega)\) is a symplectic vector
space\index{symplectic vector space|(}.

\end{Definition}
\begin{Lemma}\label{lma:nondegeneracy_map}
The map\footnote{Recall that \(\iota_v\omega\) denotes the
1-form \(\omega(v,-)\).} \(v\mapsto \iota_v\omega\) gives an
isomorphism \(V\to V^*\).
\end{Lemma}
\begin{Proof}
If \(\iota_v\omega=0\) then \(\omega(v,w)=0\) for all \(w\), so
\(v=0\) by nondegeneracy. Therefore this map is an injective map
between vector spaces of the same dimension, hence it is an
isomorphism. \qedhere

\end{Proof}
\begin{Definition}\label{dfn:symplectic_orthogonal_complement}
If \(W\subset V\) then we define the {\em symplectic
orthogonal complement}\index{symplectic orthogonal complement}
\(W^{\omega}\subset V\) to be the subspace \[W^{\omega} :=
\{v\in V\,:\,\omega(v,w)=0\mbox{ for all }w\in W\}.\]

\end{Definition}
\begin{Lemma}\label{lma:complementary}
\(\dim(W^{\omega})=\dim(V)-\dim(W)\).
\end{Lemma}
\begin{Proof}
Under the isomorphism \(v\mapsto \iota_v\omega\), the subspace
\(W^{\omega}\) is identified with the annihilator
\(W^\circ=\{f\in W^*\,:\,f(w)=0\mbox{ for all }w\in W\}\subset
W^*\). The annihilator has dimension \(\dim(V)-\dim(W)\)
{\cite[\S 16, Theorem 1]{Halmos}}.\qedhere

\end{Proof}
\begin{Lemma}\label{lma:symplectic_basis_exists}
Any nonzero finite-dimensional symplectic vector space \(V\)
admits a basis \(e_1,\ldots,e_n,f_1,\ldots,f_n\) such that
\(\omega(e_i,f_j)=\delta_{ij}\) and
\(\omega(e_i,e_j)=\omega(f_i,f_j)=0\) for all \(i,j\). Such a
basis is called a {\em symplectic
basis}\index{symplectic basis} for \(V\). As a corollary,
\(\dim(V)\) is even.
\end{Lemma}
\begin{Proof}
Suppose that we have constructed a (possibly empty) partial
symplectic basis \(e_1,\ldots,e_k,f_1,\ldots,f_k\) of size
\(2k\), i.e.\ a linearly independent set of \(2k\) vectors
satisfying the conditions in the statement of the lemma. Write
\(W\subset V\) for the span of this partial basis. Note that
the restriction of \(\omega\) to \(W\) is nondegenerate, so
\(W^{\omega}\cap W=0\) and \(W^{\omega}\) is a complement to
\(W\) by Lemma \ref{lma:complementary}.

If \(V\neq W\), pick \(e_{k+1}\in W^{\omega}\). By
construction, \(\omega(e_{k+1},e_i)=\omega(e_{k+1},f_i)=0\)
for \(i\leq k\). By nondegeneracy, we can find \(f'_{k+1}\in
W\) with \(\omega(e_{k+1},f'_{k+1})=1\). With respect to the
splitting \(V=W\oplus W^{\omega}\) we have
\(f'_{k+1}=g+f_{k+1}\) for uniquely determined vectors \(g\in
W\) and \(f_{k+1}\in W^{\omega}\). Since \(e_{k+1}\in
W^{\omega}\), we have
\(\omega(e_{k+1},g+f_{k+1})=\omega(e_{k+1},f_{k+1})\). Now
\(e_1,\ldots,e_{k+1},f_1,\ldots,f_{k+1}\) is a partial
symplectic basis of size \(2(k+1)\). At some point this
construction terminates because \(V\) is
finite-dimensional.\qedhere

\end{Proof}
\begin{Definition}\label{dfn:isotropic_symplectic}
Let \((V,\omega)\) be a symplectic vector space and \(W\subset V\) a
subspace. We say
\begin{itemize}
\item \(W\) is isotropic\index{subspace!isotropic} if \(W\subset
W^{\omega}\),
\item \(W\) is coisotropic\index{subspace!coisotropic} if
\(W^{\omega}\subset W\),
\item \(W\) is Lagrangian\index{subspace!Lagrangian|(} if it is
both isotropic and coisotropic,
\item \(W\) is symplectic\index{subspace!symplectic|(} if \(W\cap
W^{\omega}=0\).

\end{itemize}
\end{Definition}
\begin{Lemma}\label{lma:dimension_bound}
If \(W\) is isotropic then \(2\dim(W)\leq \dim(V)\). If \(W\) is
coisotropic then \(\dim(V)\leq 2\dim(W)\). In particular, we see
that Lagrangian subspaces satisfy \(2\dim(W)=\dim(V)\).
\end{Lemma}
\begin{Proof}
If \(W\) is isotropic then we have \[\dim(W)\leq
\dim(W^{\omega})=\dim(V)-\dim(W),\] so \(2\dim(W)\leq \dim(V)\). If
\(W\) is coisotropic then we have
\[\dim(V)-\dim(W)=\dim(W^{\omega})\leq \dim(W),\] so \(\dim(V)\leq
2\dim(W)\). \qedhere

\end{Proof}
\section{Complex structures}

\begin{Definition}
Let \(V\) be a vector space. A linear map \(J\colon V\to V\) is
called a {\em complex structure}\index{complex structure|(} if \(J^2 = -I\).

\end{Definition}
\begin{Definition}\label{dfn:cpt_j}
If \((V,\omega)\) is symplectic vector space then we say that:
\begin{itemize}
\item \(J\) tames\index{complex structure!tame} \(\omega\) if
\(\omega(v,Jv)>0\) for any \(v\neq 0\),
\item \(J\) is \(\omega\)-compatible\index{complex
structure!compatible} if it tames \(\omega\) and
\(\omega(Jv,Jw)=\omega(v,w)\) for all \(v,w\in V\).

\end{itemize}
\end{Definition}
\begin{Lemma}
Let \((V,\omega)\) be a symplectic vector space. If \(J\) is a
complex structure on \(V\) taming \(\omega\) and \(W\subset
V\) is a \(J\)-complex subspace\index{subspace!complex}
(i.e.\ \(JW=W\)) then \(W\) is a symplectic
subspace\index{subspace!symplectic|)}.
\end{Lemma}
\begin{Proof}
The subspace \(W\cap W^\omega\) consists of vectors \(v\in W\)
such that \(\omega(v,w)=0\) for all \(w\in W\). However, if
\(v\in W\) then \(Jv\in W\), so \(\omega(v,Jv)=0\) and
tameness implies \(v=0\). Thus \(W\cap W^{\omega}=0\) and
\(W\) is symplectic. \qedhere

\end{Proof}
\begin{Lemma}\label{lma:J_metric}
If \((V,\omega)\) is a symplectic vector space and \(J\) is an
\(\omega\)-compatible complex structure on \(V\) then
\(g_J(v,w)=\omega(v,Jw)\) defines a positive-definite
symmetric bilinear form on \(V\).
\end{Lemma}
\begin{Proof}
Bilinearity follows from bilinearity of \(\omega\). Symmetry
follows from \[g(w,v) = \omega(w,Jv) = \omega(Jw,J^2v) =
\omega(Jw,-v) = \omega(v,Jw) = g(v,w).\] Positive-definiteness
follows from the fact that \(g(v,v)=\omega(v,Jv)>0\) if
\(v\neq 0\). \qedhere

\end{Proof}
\begin{Lemma}\label{lma:j_orth}
Let \(J\) be an \(\omega\)-compatible complex structure on a
symplectic vector space \((V,\omega)\) and let \(W\) be a
subspace. We have \(W^{\omega}=(JW)^{\perp}\), where \(\perp\)
denotes the orthogonal complement with respect to \(g_J\).
\end{Lemma}
\begin{Proof}\belowdisplayskip=-12pt
Since \(g_J(v,Jw)=-\omega(v,w)\) we have
\begin{align*}
W^\omega&=\{v\in V\,:\,\omega(v,w)=0\mbox{ for all }w\in W\} \\
&= \{v\in V\,:\,g_J(v,Jw)=0\mbox{ for all }w\in W\}\\
&= (JW)^\perp.
\end{align*}\qedhere
\end{Proof}

\begin{Lemma}\label{lma:lag_orth}
Let \(J\) be an \(\omega\)-compatible complex structure on a
symplectic vector space \((V,\omega)\) and let \(L\subset V\) be a
subspace. The following are equivalent.
\begin{enumerate}
\item [(a)] \(L\) is Lagrangian\index{subspace!Lagrangian|)};
\item [(b)] \(L \perp JL\).
\item [(c)] \(JL\) is Lagrangian;

\end{enumerate}
\end{Lemma}
\begin{Proof}
We have \(L^\omega = (JL)^\perp\), so \(L=L^\omega\) if and
only if \(L\perp JL\). Thus (a) is equivalent to (b). Since
\(J^2=-I\), \(J^2L=L\), which means that (b) is symmetric in
\(L\) and \(JL\). Thus (b) is equivalent to (c).\index{complex
structure|)}\index{symplectic vector space|)}\qedhere

\end{Proof}
\chapter{Lie derivatives}
\label{ch:lie_ds}
\thispagestyle{cup}

The background we assume on differential geometry can be found
in many places, for example Lee's compendious book on smooth
manifolds \cite{Lee}, Warner's terser book on manifolds and Lie
theory \cite{Warner}, or Arnold's wonderful introduction to
differential forms {\cite[Chapter 7]{Arnold}}. However, the
theory of Lie derivatives can be difficult to swallow on a first
encounter. The philosophy behind this book has been to give a
complementary perspective rather than rehashing what can be
found written better elsewhere. In this appendix we give a quick
and high-level conceptual overview of Lie derivatives from the
point of view of Lie groups and Lie algebras, in the hopes that
the reader will find this viewpoint helpful alongside a more
traditional treatment. The final goal is to give a proof of the
``magic formulas'' relating the Lie derivative, interior product
and exterior derivative.

\section{Recap on Lie groups}

We start by giving a lightning review of Lie groups (see
\cite{Warner} for a more thorough introduction from the ground
up). A {\em Lie group}\index{Lie group} \(G\) is a
finite-dimensional manifold which is also a group in such a way
that the multiplication and inversion maps are smooth. You will
lose nothing by imagining that it is a group of matrices with
real entries. The Lie algebra\index{Lie algebra}
\(\mathfrak{g}\) of \(G\) is the tangent space of \(G\) at the
identity; in other words, if \(\phi_t\) is a smooth path in
\(G\) with \(\phi_0=\OP{id}\) then
\(\left.\frac{d\phi_t}{dt}\right|_{t=0}\) is an element of
\(\mathfrak{g}\). Going back the other way, each \(V\in
\mathfrak{g}\) arises as the tangent vector
\(\left.\frac{d\phi_t}{dt}\right|_{t=0}\) of a unique {\em
1-parameter subgroup} \(\phi_t\) of \(G\), i.e. a path
satisfying \(\phi_0=\OP{id}\) and \(\phi_{s+t}=\phi_s\phi_t\)
for all \(s,t\). This 1-parameter subgroup is usually written as
\(\exp(tV)\).

The Lie algebra \(\mathfrak{g}\) is a much simpler object than
\(G\): it is a vector space instead of a manifold, so it has no
interesting topology. It retains some knowledge of the group
structure of \(G\): it is equipped with an antisymmetric
bilinear operation \([,]\) called {\em Lie bracket}\index{Lie
bracket!on Lie algebra} which is defined as follows. If
\(V=\left.\frac{d\psi_s}{ds}\right|_{s=0}\) and
\(W=\left.\frac{d\phi_t}{dt}\right|_{t=0}\) then \[[V,W] :=
\left.\frac{d}{ds}\right|_{s=0} \left.\frac{d}{dt}\right|_{t=0}
\psi_s\phi_t\psi_s^{-1}.\] If \(\psi_s\) and \(\phi_t\) commute
for all \(s,t\) then clearly \([V,W]=0\). Somewhat miraculously,
a kind of converse holds: if \([V,W]=0\) then \(\exp(sV)\) and
\(\exp(tW)\) commute. Indeed, the full multiplication structure
of \(G\) in a neighbourhood of the identity can be determined
from the Lie bracket.

The easiest example of a Lie group is the group \(GL(n,\RR)\) of
invertible \(n\)-by-\(n\) real matrices; indeed its subgroups
can be understood in the same way. Its Lie algebra is the space
\(\mathfrak{gl}(n,\RR)\) of all \(n\)-by-\(n\) real matrices,
and the 1-parameter subgroup associated to a matrix \(V\) is
\[\exp(tV)=\sum_{n=0}^\infty\frac{1}{n!}V^n.\] The Lie bracket
is simply the commutator \([V,W]=VW-WV\).

One good way to understand a given Lie group \(G\) is to map it
(smoothly and homomorphically) to a subgroup of
\(GL(n,\RR)\). Such a smooth homomorphism \(R\colon G\to
GL(n,\RR)\) is called a {\em representation}\index{Lie
group!representation of} of \(G\). Its differential at the
identity matrix is a linear map
\[\rho:=d_1R\colon\mathfrak{g}\to\mathfrak{gl}(n,\RR)\] which is
a {\em representation} of the Lie algebra\index{Lie
algebra!representation of}, in the sense that
\[[\rho(V),\rho(W)]=\rho([V,W])\mbox{ for all
}V,W\in\mathfrak{g}.\] For example, \(GL(n,\RR)\) acts on
\(\mathfrak{gl}(n,\RR)\) by conjugation, which gives a
representation \[\OP{Ad}\colon GL(n,\RR)\to
GL(\mathfrak{gl}(n,\RR)),\qquad\OP{Ad}(g)(V)=gVg^{-1}\] and its
differential is the representation
\[\OP{ad}\colon\mathfrak{gl}(n,\RR)\to
\mathfrak{gl}(\mathfrak{gl}(n,\RR)),\qquad
\OP{ad}(V)(W)=[V,W].\] These are both called the {\em adjoint}
representation.\index{representation!adjoint}

\section{Diffeomorphism groups}

Let \(M\) be a smooth manifold and \(G=\OP{Diff}(M)\) be the
group of diffeomorphisms \(M\to M\). This is not a Lie group,
because it is not finite-dimensional, but many of the same ideas
apply.

If we take a point \(p\in M\) and a 1-parameter family of
diffeomorphisms \(\phi_t\) then we get a path
\(\phi_t(p)\). This gives a tangent vector
\(\left.\frac{d\phi_t(p)}{dt}\right|_{t=0}\) at \(p\) and hence
a vector field on \(M\) whose value at \(p\) is \[V(p) =
\left.\frac{d\phi_t(p)}{dt}\right|_{t=0}.\] Conversely, given a
vector field \(V\) on \(M\) we get a 1-parameter subgroup of
\(\OP{Diff}(M)\) given by the flow\footnote{It is of course
possible that the flow is only locally defined, or defined for
small \(t\). For example, if \(M=\RR\) and \(V(x)=-x^2\) then
the flow is \(\phi_t(x) = x/(xt+1)\), and we see that
\(\lim_{t\to -1/x}\phi_t(x)=\infty\). We avoid this kind of
behaviour if \(M\) is compact or if \(\phi_t\) preserves the
(compact) level sets of some proper function, e.g. if \(\phi_t\)
is the Hamiltonian flow of a proper Hamiltonian.} \(\phi_t\) of
\(V\), which is the family of diffeomorphisms satisfying
\[\frac{d}{dt}\phi_t(p) = V(\phi_t(p))\mbox{ for all
}t\in\RR,\,p\in M.\] For this reason, we will think of the space
\(\Vect(M)\) of vector fields on \(M\) as the Lie algebra of
\(\OP{Diff}(M)\).\index{diffeomorphism group!Lie algebra} We
will figure out the Lie bracket in a moment.

The easiest way to study \(\OP{Diff}(M)\) and its Lie algebra is
via its representations. It comes with a natural and plentiful
supply. We will write \(\Omega^k(M)\) for the space of smooth
differential \(k\)-forms (so \(\Omega^0(M)\) means smooth
functions).\index{diffeomorphism group!representations}

\begin{Example}
\(\OP{Diff}(M)\) acts on functions by pullback: \[f\mapsto
\phi^*f :=f\circ\phi.\] This is a {\em right action} in the
sense that \((\phi_1\circ\phi_2)^* = \phi_2^*\circ\phi_1^*\)
so it gives an\index{antirepresentation} {\em
antirepresentation}\footnote{If the notation
\(GL(\Omega^0(M))\) is giving you a headache, it just means
``invertible linear maps
\(\Omega^0(M)\to\Omega^0(M)\)''. Similarly
\(\mathfrak{gl}(\Omega^0(M))\) will mean ``linear maps
\(\Omega^0(M)\to\Omega^0(M))\)''.} \[R\colon \OP{Diff}(M)\to
GL(\Omega^0(M)),\qquad\phi\mapsto\phi^*.\] This therefore
gives an antirepresentation \(\rho\colon\Vect(M)\to
\mathfrak{gl}(\Omega^0(M))\) of the Lie algebra \(\Vect(M)\)
by defining \[\rho(V)(f) =
\left.\frac{d}{dt}\right|_{t=0}\phi_t^*f\] where \(\phi_t\) is
the flow of \(V\). This is better known as the {\em
directional derivative} of \(f\) in the \(V\)-direction, and
written \(V(f)\). In fact, since pullback of functions is
multiplicative: \[\phi^*(fg)=(\phi^*f)(\phi^*g),\] the
representation \(R\) lands in the subgroup
\(\OP{Aut}(\Omega^0(M))\subset GL(\Omega^0(M))\) of
automorphisms of the ring of functions, and so \(\rho\) lands
in the subalgebra
\(\mathfrak{der}(\Omega^0(M))\subset\mathfrak{gl}(\Omega^0(M))\)
of derivations of the ring of functions. In fact, this
representation is injective: we can really think of vector
fields as derivations on functions without losing information,
and many expositions {\em define} vectors as derivations.

We can now identify the Lie bracket on \(\OP{Vect}(M)\):
because \(\rho\) is an {\em anti}repr\-esen\-ta\-tion, it
should be {\em minus} the commutator bracket on
derivations.\index{Lie bracket!on vector fields} To avoid
sign-clashes with the rest of the literature, we actually
write \([,]\) for the commutator bracket \[[V,W](f) =
V(W(f)) - W(V(f)).\] and call it ``Lie bracket of vector
fields'', even though it is out by a minus sign.

\end{Example}
\begin{Example}
Pullback of differential forms gives another natural
antirepresentation \[\OP{Diff}(M)\to
GL(\Omega^*(M)),\qquad\phi\mapsto\phi^*.\] This preserves
wedge product of differential forms: \[\phi^*(\eta_1\wedge
\eta_2) = \phi^*\eta_1 \wedge \phi^*\eta_2.\] We write
\(V\mapsto \Lie_V\) for the corresponding Lie algebra
antirepresentation \(\Vect(M)\to \mathfrak{gl}(\Omega^*(M))\):
\[\Lie_V\eta := \left.\frac{d}{dt}\right|_{t=0}\phi_t^*\eta.\]
By the Leibniz rule, \(\Lie_V\) acts as a {\em derivation} of
the algebra \(\Omega^*(M)\): \[\Lie_V(\eta_1\wedge\eta_2) =
(\Lie_V\eta_1)\wedge \eta_2+ \eta_1\wedge\Lie_V\eta_2,\qquad
\Lie_Vd\eta = d\Lie_V\eta.\] Moreover, since this is an {\em
anti}representation, and since we have grudgingly accepted a
historical minus sign in our bracket, we get
\begin{equation}\label{eq:lie_alg_rep}\Lie_{[V,W]}\eta =
\Lie_V\Lie_W\eta - \Lie_W\Lie_V\eta.\end{equation}

\end{Example}
\section{Cartan's magic formulas}

The algebra of differential forms admit further natural
operations which have played a key role in this book. Our goal
here is to prove the ``magic formulas'' which govern the
interplay between these operations and \(\Lie_V\).

The operations in question are the exterior
derivative\footnote{The most conceptually well-motivated
exposition of \(d\) that I know is due to Arnold: {\cite[Chapter
7]{Arnold}}.}: \[d\colon\Omega^*(M)\to \Omega^{* +1}(M),\] and
the interior product: \[\iota_V\colon \Omega^*(M)\to\Omega^{*
-1}(M),\qquad (\iota_V\eta)(V_1,\ldots,V_{k-1}) =
\eta(V,V_1,\ldots,V_{k-1}),\] where \(V\) is a choice of vector
field. Both of these operations are {\em antiderivations},
i.e.\, \[d(\eta_1\wedge \eta_2) = (d\eta_1)\wedge \eta_2 +
(-1)^{|\eta_1|}d\eta_2.\] First, we prove some lemmas.

\begin{Lemma}\label{lma:lie_d_commute}
\[d\Lie_V\eta=\Lie_V d\eta.\]
\end{Lemma}
\begin{Proof}
Let \(\phi_t\) be the flow along \(V\). We have
\(\phi_t^*d\eta=d\phi_t^*\eta\). Differentiating with respect
to \(t\) at \(t=0\) gives the required identity.\qedhere

\end{Proof}
\begin{Lemma}\label{lma:lie_func}
If \(f\) is a function (0-form) then \[\Lie_Vf=\iota_Vdf.\]
\end{Lemma}
\begin{Proof}
Let \(\phi_t\) be the flow of \(V\). In local coordinates, for
small \(t\), we have\footnote{As usual, \(\oo(t)\) denotes a
quantity such that \(\lim_{t\to 0}\oo(t)/t=0\).} \(\phi_t(p)
= p+tV(p)+\oo(t)\) and therefore \[(\phi_t^*f)(p) =
f(\phi_t(p))=f(p)+t\,df(V(p))+\oo(t).\] This means that
\[\left.\frac{d}{dt}\right|_{t=0}(\phi_t^*f)(p) = df(V(p)),\]
or \[\Lie_Vf = \iota_Vdf.\qedhere\]

\end{Proof}
\begin{Theorem}[Cartan's magic formulas]\label{thm:mfs}
We have\index{Cartan's formula}
\begin{align}
\label{eq:cartan_mf}\iota_Vd+d\iota_V&= \Lie_V\\
\label{eq:other_mf}\Lie_V\iota_W-\iota_W\Lie_V &= \iota_{[V,W]}.
\end{align}
\end{Theorem}
\begin{Proof}
The operators \(d\) and \(\iota_V\) are antiderivations and
\(\Lie_V\) is a derivation. This implies that the operator
\(\iota_Vd+d\iota_V\) is a derivation and that
\(\Lie_V\iota_W-\iota_W\Lie_V\) is an
antiderivation. Derivations and antiderivations on
\(\Omega^k(M)\) are determined by their effect on functions
and on exact 1-forms. This is easy to see in local
coordinates: if \[\eta = \sum \eta_{i_1\cdots
i_k}dx_{i_1}\wedge\cdots\wedge dx_{i_k}\] and \(D\) is a
derivation, for example, then \begin{align*}D\eta &= \sum
(D\eta_{i_1\cdots i_k})dx_{i_1}\wedge\cdots\wedge
dx_{i_k}+\sum \eta_{i_1\cdots
i_k}(D(dx_{i_1}))\wedge\cdots\wedge dx_{i_k}+\\&
\qquad+\cdots+\sum \eta_{i_1\cdots
i_k}dx_{i_1}\wedge\cdots\wedge (D(dx_{i_k})),\end{align*} so
\(D\) is determined completely by its action on functions
(like the coefficients \(\eta_{i_1\cdots i_k}\)) and exact
1-forms (like the local coordinate 1-forms \(dx_i\)).

Therefore it suffices to check Equation \eqref{eq:cartan_mf}
and \eqref{eq:other_mf} for functions and for exact 1-forms.

{\bf Equation \eqref{eq:cartan_mf}.} For functions \(f\),
Equation \eqref{eq:cartan_mf} is simply the identity
\(\Lie_Vf=df(V)\) which we proved in Lemma
\ref{lma:lie_func}. Now suppose we have an exact 1-form
\(df\). We have \((d\iota_V+\iota_Vd)df=d\iota_Vdf\) because
\(d^2=0\) and \[\Lie_Vdf=d\Lie_Vf=d\iota_Vdf,\] using Lemmas
\ref{lma:lie_d_commute} and \ref{lma:lie_func}, so both sides
of Equation \eqref{eq:cartan_mf} agree when applied to
\(df\). This proves Equation \eqref{eq:cartan_mf}.

{\bf Equation \eqref{eq:other_mf}.} For functions \(f\),
Equation \eqref{eq:other_mf} reduces to \(0=0\). For exact
1-forms \(df\), using \(\iota_Wdf=\Lie_Wf\), the left-hand
side of the Equation \eqref{eq:other_mf} becomes
\[\Lie_V\Lie_Wf-\iota_W\Lie_Vdf.\] Since
\(\iota_W\Lie_Vdf=\iota_Wd\Lie_V=\Lie_W\Lie_Vf\), this becomes
\(\left(\Lie_V\iota_W-\iota_W\Lie_V\right)df =
\left(\Lie_V\Lie_W-\Lie_W\Lie_V\right)f\), which becomes
\(\Lie_{[V,W]}f\) using Equation
\eqref{eq:lie_alg_rep}. Therefore
\[\left(\Lie_V\iota_W-\iota_W\Lie_V\right)df =
\Lie_{[V,W]}f=\iota_{[V,W]}df.\] This proves Equation
\eqref{eq:other_mf} for exact 1-forms and hence in
general.\qedhere

\end{Proof}
\newpage

\chapter{Complex projective spaces}
\label{ch:complex_projective_spaces}
\thispagestyle{cup}

\section{\(\cp{n}\)}

\begin{Definition}
The complex projective space \(\cp{n}\)
\index{projective space!complex|(}
\index{CPn@$\mathbb{CP}^n$|see {projective, space, complex}}
is the space of complex lines in \(\CC^{n+1}\) passing through
the origin.

\end{Definition}
Recall that a complex line passing through the origin is a
subspace\index{subspace!complex} of the form
\(\CC\cdot \bm{z}:=\{(\lambda z_1,\ldots,\lambda
z_{n+1})\,:\,\lambda\in\CC\}\subset\CC^{n+1}\) for some complex
vector \(\bm{z}=(z_1,\ldots,z_{n+1})\neq 0\).

\begin{Lemma}[Exercise \ref{exr:complex_lines}]\label{lma:complex_lines}
If \(\bm{z},\bm{z}'\neq 0\) are complex vectors then
\(\CC\cdot \bm{z}=\CC\cdot \bm{z}'\) if and only if
\(\bm{z}=\mu \bm{z}'\) for some complex number \(\mu\neq 0\).

\end{Lemma}
\begin{Lemma}\label{lma:cpn_quotient}
The complex projective space \(\cp{n}\) is the quotient of
\(\CC^{n+1}\setminus\{0\}\) by the equivalence relation
\(\bm{z}\sim \bm{z}'\) if and only if \(\bm{z}'=\mu \bm{z}\)
for some complex number \(\mu\neq 0\).
\end{Lemma}
\begin{Proof}
Each \(\bm{z}\in\CC^{n+1}\setminus\{0\}\) gives us a line
\(\CC\cdot \bm{z}\) and every complex line has this form, so
the map \[Q\colon\CC^{n+1}\setminus\{0\}\to\cp{n},\qquad
Q(\bm{z})=\CC\cdot \bm{z}\] is a surjection. By Lemma
\ref{lma:complex_lines}, the fibres of \(\CC\,\cdot\) are the
stated equivalence classes. \qedhere

\end{Proof}
We equip \(\cp{n}\) with the quotient topology induced by
\(Q\). We often write \([\bm{z}]\) or \([z_1:\cdots:z_{n+1}]\)
for \(\CC\cdot \bm{z}\), and call the \(z_i\) {\em homogeneous
coordinates}\index{homogeneous coordinates|(} on
\(\cp{n}\). Homogeneous coordinates are not like Cartesian
coordinates: Cartesian coordinates have the property that if
\(p\) and \(q\) have different coordinates then \(p\neq q\), but
the homogeneous coordinates \([1:1]\) and \([2:2]\) specify the
same point in \(\cp{1}\). We remedy this redundancy, at the cost
of missing some points, by passing to {\em affine
charts}\index{affine chart|(}.

\begin{Definition}
Let \(A_k = \{(z_1,\ldots,z_{n+1})\in\CC^{n+1}\, :
\,z_k=1\}\subset\CC^{n+1}\).

\end{Definition}
\begin{Lemma}[Exercise \ref{exr:aff_chart_cpn}]\label{lma:aff_chart_cpn}
The restriction \(Q|_{A_k}\colon A_k\to\cp{n}\) is an
embedding. We call its image an {\em affine chart} in
\(\cp{n}\).

\end{Lemma}
\begin{Lemma}\label{lma:cpn_atlas}
The topological space \(\cp{n}\) is a complex manifold (in
fact an algebraic variety).
\end{Lemma}
\begin{Proof}
If \(\bm{z}\in \CC^{n+1}\) has \(z_k\neq 0\) then
\(\bm{z}/z_k\in A_k\) and \(Q(\bm{z})=Q(\bm{z}/z_k)\), so
\(Q(\bm{z})\in Q(A_k)\). The space \(\cp{n}\) is therefore
covered by the \(n+1\) affine charts
\(Q(A_1),\ldots,Q(A_{n+1})\). The transition map
\[\varphi_{k\ell}:=Q|_{A_k}^{-1}\circ
Q|_{A_\ell}\colon\{\bm{z}\in A_\ell\, :\, z_k\neq
0\}\to\{\bm{z}\in A_k\, :\, z_\ell\neq 0\}\] is given by
\[\varphi_{k\ell}(z_1, \ldots, z_k, \ldots, z_{\ell}=1,
\ldots, z_{n+1}) = \left(\frac{z_1}{z_k}, \ldots, 1, \ldots,
\frac{1}{z_k}, \ldots, \frac{z_{n+1}}{z_k}\right),\] which is
an algebraic isomorphism, so \(\cp{n}\) is an algebraic
variety\footnote{Just as one can define a manifold as a
collection of charts glued together by transition maps, an
algebraic variety can be defined as a collection of affine
varieties glued together by algebraic transition maps. Affine
varieties are just subsets of affine space cut out by
polynomial equations. In this case, the affine
charts\index{affine chart|)} are copies of \(\CC^n\) rather
than something more exotic.} locally modelled on \(\CC^n\)
and hence a complex manifold. \qedhere

\end{Proof}
\begin{Example}
The complex projective 1-space \(\cp{1}\) is diffeomorphic to
the 2-sphere. To see this, let \((r,\theta,h)\) be cylindrical
polar coordinates on \(S^2\) and define \(\sigma_1\) and
\(\sigma_2\) by
\[\sigma_1(r,\theta,h)=\frac{r}{1-h}e^{i\theta},\qquad
\sigma_2(r,\theta,h)=\frac{r}{1+h}e^{i\theta}\] for \(h<1\)
and \(h>-1\) respectively. These are the
stereographic\index{stereographic projection} projections from
the North and South poles respectively. Since \(h^2+r^2=1\),
we have \(\sigma_1\overline{\sigma}_2 = 1\), so if we use
\(\sigma_1\) and \(\overline{\sigma}_2\) as coordinate charts
on \(S^2\) then the transition map is \(\sigma_1\mapsto
\bar{\sigma}_2=1/\sigma_1\). This is the same as the
transition map \(\varphi_{1,2}\) for \(\cp{1}\) from the proof
of Lemma \ref{lma:cpn_atlas}, so these manifolds are
diffeomorphic.

\end{Example}
\begin{figure}[htb]
\begin{center}
\begin{tikzpicture}
\draw (0,0) circle [radius = 2cm];
\draw (-2,0) arc [x radius = 2cm, y radius = 0.5cm, start angle = 180, end angle = 360];
\draw[dashed] (2,0) arc [x radius = 2cm, y radius = 0.5cm, start angle = 0, end angle = 180];
\draw (-4,-1.5) -- (8,-1.5) -- (6,1.5) -- (-2,1.5) -- cycle;
\draw[dotted,thick] (0,-2) -- (45:2cm);
\node (s_2) at (0.828427124746190,0) {\(\bullet\)};
\node (s_3) at (s_2) [above left=10pt] {\(\sigma_2(r,h,\theta)\)};
\draw[->] (s_3) -- (s_2);
\node (s_4) at (45:2cm) {\(\bullet\)};
\node at (s_4) [above right] {\((r,h,\theta)\)};
\node (s_1) at (4.828427124746190,0) {\(\bullet\)};
\node at (s_1) [right] {\(\sigma_1(r,h,\theta)\)};
\draw[dotted,thick] (0,2) -- (s_1.center);
\draw (0,0) -- (s_1.center);
\draw[->] (0,0) -- (-120:1.5cm);
\draw (-120:0.45cm) arc [radius = 0.45, start angle = -120, end angle = 0];
\node at (-45:0.25cm) {\(\theta\)};
\draw[dashed] (45:2cm) -- (1.41421356237310,0);
\node at (1.41421356237310,0.5) [right] {\(h\)};
\node at (1.41421356237310,0) [below] {\(r\)};
\end{tikzpicture}
\end{center}
\caption{Stereographic projections}
\label{fig:stereographic_projections}
\end{figure}

\begin{Lemma}
The complex projective space \(\cp{n}\) is the quotient of the
unit sphere \(S^{2n+1}\subset\CC^{n+1}\) by the
\(U(1)\)-action where \(u\in U(1)\) acts on \(\CC^{n+1}\) by
\(\bm{z}\mapsto u\bm{z}\). Here, \(U(1)\) is the multiplicative group of
unit complex numbers.
\end{Lemma}
\begin{Proof}
The restriction of \(Q\) to \(S^{2n+1} \subset \CC^{n+1}
\setminus \{0\}\) is still surjective because every complex
line contains a circle of vectors of unit length. The fibres
of \(Q|_{S^{2n+1}}\) are precisely these unit circles in each
complex line, which are also the orbits of the \(U(1)\)-action
in the statement of the
lemma\index{projective space!complex|)}.\qedhere

\end{Proof}
\section{Projective varieties}

A polynomial \(F(z_1,\ldots,z_{n+1})\) is called {\em
homogeneous of degree \(d\)}\index{homogeneous polynomial|(} if
\[F(\lambda z_1,\ldots, \lambda z_{n+1})=\lambda^d
F(z_1,\ldots,z_{n+1})\] for all \(\lambda\in\CC\). For example,
\(z_1z_2+z_3z_4\) is homogeneous of degree \(2\), whereas
\(z_1+z_2^2\) is not homogeneous of any degree. It is not hard
to show that a homogeneous polynomial of degree \(d\) is
precisely a linear combination of monomials each of which has
degree precisely \(d\).

The advantage of working with homogeneous polynomials is that if
\(F(\bm{z})=0\) and \([\bm{z}]=[\bm{z}']\) then
\(F(\bm{z}')=0\). In other words, it makes sense to write
\(F([\bm{z}])=0\): this condition doesn't depend on the choice
of the homogeneous coordinate \(\bm{z}\).

\begin{Definition}
We define the {\em projective subvariety cut out by \(F\)},
\index{subvariety} to be the subset
\(\VV(F):=\{[\bm{z}]\in\cp{n}\, :\,
F([\bm{z}])=0\}\). Similarly, we can define subvarieties cut
out by a (possibly empty) set of homogeneous polynomials
\(\{F_1,\ldots,F_s\}\). Note that \(\VV(\emptyset)=\cp{n}\).

\end{Definition}
\begin{Example}
Let \(F(z_1,z_2)=z_1z_2\). Then \(\VV(F)\subset\cp{1}\) consists
of the points \([1:0]\) and \([0:1]\).

\end{Example}
\begin{Remark}
In affine algebraic geometry, an affine variety is the subset
cut out of an affine space by a collection of (not necessarily
homogeneous) polynomials. The intersection of the projective
variety \(\VV(F)\) with the affine chart \(Q(A_k)\) is defined
by \(F(z_1,\ldots,z_k=1,\ldots,z_{n+1})=0\), which is a (not
necessarily homogeneous) polynomial, so \(\VV(F)\cap Q(A_k)\)
is an affine variety in the chart \(Q(A_k)\).

\end{Remark}
\begin{Remark}
Unless \(F\) is constant, the expression \(F([\bm{z}])\) does
not make sense as a complex-valued function on \(\cp{n}\). For
example, if \(F(\bm{z})\neq 0\) and \(d\geq 1\) then
\(F(2\bm{z})=2^d F(\bm{z})\neq F(\bm{z})\), but
\([2\bm{z}]=[\bm{z}]\).

\end{Remark}
\section{Zariski-closure}

We will use the notion of Zariski-closure\index{Zariski-closure|(} in Chapter
\ref{ch:toric_varieties}.

\begin{Definition}
The {\em Zariski-closure} of a subset \(S\subset\cp{n}\) is
the smallest subvariety\index{subvariety|(} containing \(S\). If
\(V\subset\cp{n}\) is a subvariety then a subset \(S\subset
V\) is called {\em Zariski-dense}\index{Zariski-dense|(} in
\(V\) if its Zariski-closure is \(V\).

\end{Definition}
\begin{Example}
Consider the set of points \(S=\{[n:1]\, :\,
n\in\ZZ\}\subset\cp{1}\). The Zariski-closure of this set is
\(\cp{1}=\VV(0)\). This is because if there is a homogeneous
polynomial \(F=\sum_{m=0}^d a_mz_1^mz_2^{d-m}\) with
\(F([n:1])=0\) for all \(n\in\ZZ\) then \(\sum_{m=0}^da_mz^m\)
has infinitely many zeros (every integer) and hence vanishes
identically.

\end{Example}
\begin{Example}
Consider the set of points \(S=\{[z:1:0]\in\cp{2}\, :\,
z\in\CC\}\). This is contained in the subvariety \(\VV(z_3)\)
defined by the vanishing of the \(z_3\)-coordinate. Moreover,
it is Zariski-dense inside that subvariety. To see this,
suppose that \(F=\sum_{m_1+m_2+m_3=d} a_{m_1,m_2,m_3}
z_1^{m_1} z_2^{m_2} z_3^{m_3}\) is a homogeneous polynomial
with \(F(z,1,0)=0\) for all \(z\in\CC\). Then
\(\sum_{m_1,m_2}a_{m_1,m_2,0}z^{m_1}=0\) has infinitely many
solutions (any \(z\in\CC\)), so \(a_{m_1,m_2,0}=0\) for all
\(m_1,m_2\). Thus \(F\) is divisible by \(z_3\). Therefore
\(\VV(F)\) contains \(\VV(z_3)\). Thus \(\VV(z_3)\) is the smallest
subvariety\index{subvariety|)} containing
\(S\).\index{Zariski-closure|)}\index{Zariski-dense|)}\index{homogeneous
coordinates|)}\index{homogeneous polynomial|)}

\end{Example}
\section{Solutions to inline exercises}

\begin{Exercise}[Lemma \ref{lma:complex_lines}]\label{exr:complex_lines}
If \(\bm{z},\bm{z}'\neq 0\) are complex vectors then
\(\CC\cdot\bm{z}=\CC\cdot\bm{z}'\) if and only if \(\bm{z}=\mu
\bm{z}'\) for some complex number \(\mu\neq 0\).
\end{Exercise}
\begin{Solution}
If \(\bm{z}=\mu \bm{z}'\) then \(\CC\cdot\bm{z}=\{\lambda
\bm{z}\,\:\,\lambda\in\CC\}=\{\lambda\mu
\bm{z}'\,:\,\lambda\in\CC\}\) and as \(\lambda\) varies over
\(\CC\), so \(\lambda\mu\) varies over \(\CC\), so
\(\CC\cdot\bm{z}=\CC\cdot\bm{z}'\). Conversely, if
\(\CC\cdot\bm{z}=\CC\cdot\bm{z}'\) then \(\bm{z}'\in
\CC\cdot\bm{z}\) so \(\bm{z}'=\mu \bm{z}\) for some
\(\mu\in\CC\). Since \(\bm{z}'\neq 0\), \(\mu\neq
0\). \qedhere

\end{Solution}
\begin{Exercise}[Lemma \ref{lma:aff_chart_cpn}]\label{exr:aff_chart_cpn}
The restriction \(Q|_{A_k}\colon A_k\to\cp{n}\) is an
embedding.
\end{Exercise}
\begin{Proof}
If \(\bm{z},\bm{z}'\in A_k\) then \(z_k=z'_k=1\). If
\(Q(\bm{z})=Q(\bm{z}')\) then \(\bm{z}'=\mu \bm{z}\) and so
\(1 = z'_k = \mu z_k=\mu\). Thus \(\mu=1\) and
\(\bm{z}'=\bm{z}\). \qedhere

\end{Proof}
\chapter{Cotangent bundles}
\label{ch:cot_bun}
\thispagestyle{cup}

The simplest symplectic manifolds beyond
\((\RR^{2n},\sum_{i=1}^n dp_i\wedge dq^i)\) are the {\em
cotangent bundles}\index{cotangent bundle|(} of manifolds. We
start by defining the symplectic structure on cotangent bundles,
and use it to give a formulation of Noether's famous theorem
relating symmetries and conserved quantities. Next, we introduce
a class of Hamiltonians on cotangent bundles which generate
(co)geodesic flows. This is all intended to illustrate the power
and scope of the Hamiltonian formalism, but also serves as a
source of examples in Chapter \ref{ch:symp_cut} (Examples
\ref{exm:periodic_geodesics} and
\ref{exm:periodic_geodesics_ctd}). Finally, we give a geometric
interpretation of the Hamilton-Jacobi equation.

Note: In this appendix, we use up-indices for components of
vectors, down-indices for components of
covectors\index{covector}. This is because of the appearance of
metric tensors and Christoffel symbols in the section on
cogeodesic flow.

\section{Cotangent bundles}

Let \(Q\) be a manifold and \(q\in Q\) be a point. Recall that a
\index{covector}covector\footnote{I like to use \(\eta\) (eta)
because \(\eta\) {\em eats} a vector and outputs a number.}
\(\eta\) at \(q\) is a linear map \(T_qQ\to \RR\). We write
\(T^*_qQ\) for the space of covectors at \(q\) and \(\pi\colon
T^*Q\to Q\) for the {\em cotangent bundle} of covectors over
\(Q\). Recall that a 1-form is a section of the cotangent
bundle, that is a covector at each point.

\begin{Definition}\label{dfn:canonical_1_form}
The cotangent bundle \(T^*Q\) carries a {\em canonical
1-form}\index{canonical 1-form} \(\lambda\), defined as
follows. Let \(q\in Q\) and \(\eta\in T^*_qQ\). If \(v\in
T_{\eta}T^*Q\) is a tangent vector to the cotangent bundle at
the point \(\eta\) then \(\lambda(v):=\eta(\pi_*(v))\).

\end{Definition}
\begin{Remark}
If we pick local coordinates \((q^1,\ldots,q^n)\) on \(Q\) and use
the coordinates \(p=\sum p_i\, dq^i\) for \(p\in T_q^*Q\) then
\(\lambda=\sum p_i\,dq^i\).

\end{Remark}
\begin{Definition}\label{dfn:canonical_symplectic_form}
We call the 2-form \(\omega = d\lambda\) the {\em canonical
symplectic structure}\index{symplectic form!canonical}
\index{canonical symplectic form|see {symplectic form,
canonical}}\index{symplectic form!on cotangent bundle|see
{symplectic form, canonical}} on \(T^*Q\); in local
coordinates this is \(\sum_{i=1}^n dp_i\wedge dq^i\), which
makes it clear why it is symplectic.

\end{Definition}
The coordinates \(p_i,q^i\) are called {\em canonical
coordinates}\index{canonical coordinates} on \(T^*Q\). They are
canonical in the sense that, once the local coordinates
\(q^1,\ldots,q^n\) are chosen on \(Q\), we get a basis
\(dq^1,\ldots,dq^n\) for the fibres \(T_q^*Q\), and so we get
fibre coordinates \(p_i\) for free.

\begin{Remark}[Exercise \ref{exr:cot_as_ham_sys}]\label{rmk:cot_as_ham_sys}
Pick local coordinates \(q^i\) on a patch in \(Q\) and
consider the Hamiltonian system
\((q^1\circ\pi,\ldots,q^n\circ\pi)\) on the \(\pi\)-preimage
of this patch. Show that the canonical coordinates \(p_i\) are
{\em minus} the Liouville coordinates\index{Liouville
coordinates} associated with the global Lagrangian section
given by the zero-section. Does the zero-section inherit an
integral affine structure?

\end{Remark}
The fact that the \(p_i\) are canonical in the above sense has
the consequence that changing coordinates on \(Q\) induces a
{\em symplectic} change of canonical coordinates on
\(T^*Q\). Namely:

\begin{Lemma}\label{lma:diff_action_on_cotangent}
If \(\psi\colon Q\to Q\) is a diffeomorphism then
\[\psi_*=(\psi^{-1})^*\colon T^*Q\to
T^*Q,\quad\psi_*(\eta)=\eta\circ (d\psi^{-1})\in
T_{\psi(q)}^*Q\mbox{ for }\eta\in T_q^*Q\] is a symplectomorphism of
the cotangent bundle.
\end{Lemma}
\begin{Proof}
This is immediate from the fact that the canonical 1-form (and
hence \(\omega\)) are defined without reference to
coordinates, but you can see it explicitly as follows. Suppose
\(\psi\colon Q\to Q\) is a diffeomorphism (change of
coordinates); we will write the new coordinates as
\(\psi^1(q),\ldots,\psi^n(q)\). The basis \(d\psi^i\) is given
by \(d\psi^i=\sum\frac{\partial\psi^i}{\partial
q^j}\,dq^j\). Let us write \(\Psi:=d\psi\) for the matrix
\(\frac{\partial\psi^i}{\partial q^j}\). The new fibre
coordinates \(p'_i\) are chosen so that \(\sum
p'_i\,d\psi^i=\sum p_i\,dq^i\), so \(p'_i=p_j(\Psi^{-1})^j_i\)
(here, we think of \(p\) as a row vector so the matrix
\(\Psi^{-1}\) multiplies it on the right). The change of
coordinates \((p_i,q^j)\mapsto (p'_i,\psi^j)\) is symplectic
because \[\sum_i dp'_i\wedge d\psi^i = \sum_{i,j,k}
dp_j(\Psi^{-1})^j_i \wedge \Psi^i_kdq^k=\sum_j dp_j\wedge
dq^j.\qedhere\]

\end{Proof}
\begin{Remark}[Exercise \ref{exr:parallels_differences}]\label{rmk:parallels_differences}
The observant reader will detect parallels, but also notice subtle
differences, between Lemma \ref{lma:globcoords} and Lemma
\ref{lma:diff_action_on_cotangent}. Contemplate these parallels and
differences, and then turn to Exercise
\ref{exr:parallels_differences} and its solution to read more.

\end{Remark}
\section{Cogeodesic flow}

Suppose that \(g\) is a metric on \(Q\). Write \(|\eta|\) for
the length (with respect to \(g\)) of a covector \(\eta\in
T^*_qQ\).

\begin{Definition}\label{dfn:cogeodesic_flow}
Consider the function \(H\colon T^*Q\to \RR\) defined by
\(H(\eta)=\frac{1}{2}|\eta|^2\). The Hamiltonian flow
generated by \(H\) is called the {\em cogeodesic
flow}.\index{flow!cogeodesic}

\end{Definition}
\begin{Lemma}
If \((p(t),q(t))\) is a flowline of the cogeodesic flow then
\(q(t)\) is a geodesic on \(Q\) for the metric \(g\) and \(p(t)\) is
\(g\)-dual to \(\dot{q}(t)\), that is \(p(t)(w)=g(\dot{q}(t),w)\)
for all \(w\in T_{q(t)}Q\).
\end{Lemma}
\begin{Proof}
In local coordinates \(q^i\), let \(g_{ij}\) be the metric
(symmetric in \(i\) and \(j\)) and \(g^{ij}\) be its inverse
(i.e.\ \(\sum_j g^{ij}g_{jk}=\delta^i_k\)), so that if \(\eta=\sum
p_i\,dq^i\) then \(H=\frac{1}{2}|\eta|^2=\frac{1}{2}\sum_{i,j}
g^{ij}p_ip_j\). Thus \[\dot{q}^k = \frac{\partial H}{\partial
p_k}=\sum_j g^{kj}p_j\,\qquad \dot{p}_k = \frac{\partial H}{\partial
q^k}=\frac{1}{2}\sum_{i,j}\frac{\partial g^{ij}}{\partial
q^k}p_ip_j.\] The first equation tells us that \(p_j=\sum_k
g_{jk}\dot{q}^k\) as desired. It remains to show that \(q(t)\) is a
geodesic. Substituting \(p_k=\sum_j g_{jk}\dot{q}^j\) into the
second equation
gives \begin{equation}\label{eq:geod_1}\dot{p_k}=\sum_{k,\ell}\frac{\partial
g_{jk}}{\partial q^\ell}\dot{q}^j\dot{q}^\ell +
\sum_{j}g_{jk}\ddot{q}^j=-\frac{1}{2}\sum_{i,j,\ell,m}
\frac{\partial g^{ij}}{\partial q^k}g_{i\ell}g_{jm}\dot{q}^\ell
\dot{q}^m.\end{equation} By differentiating \(\sum_j
g^{ij}g_{jm}=\delta^i_m\) we get \(\sum_{i,j}\frac{\partial
g^{ij}}{\partial q^k}g_{i\ell}g_{jm}=-\frac{\partial g_{\ell
m}}{\partial q^k}\). Rearranging Equation \eqref{eq:geod_1} gives
\[\ddot{q}^j = \sum_{k,\ell,m}g^{jk}\left(\frac{1}{2}\frac{\partial
g_{\ell m}}{\partial q^k}-\frac{\partial g_{mk}}{\partial
q^\ell}\right)\dot{q}^\ell \dot{q}^m.\] Because we are summing over
\(\ell,m\) and \(\dot{q}^\ell \dot{q}^m\) is symmetric in
\(\ell,m\), we can rewrite this as \[\ddot{q}^j =
\frac{1}{2}\sum_{k,\ell,m}g^{jk}\left(\frac{\partial g_{\ell
m}}{\partial q^k} - \frac{\partial g_{mk}}{\partial q^\ell} -
\frac{\partial g_{\ell k}}{\partial q^m}\right)\dot{q}^\ell
\dot{q}^m=-\sum_{\ell,m}\Gamma^j_{\ell m}\dot{q}^\ell \dot{q}^m,\]
where \(\Gamma^j_{\ell m}\) are the Christoffel symbols. This is
precisely the geodesic equation for the path \(q(t)\). \qedhere

\end{Proof}
\begin{Remark}\label{rmk:arc_length}
Note that \(q(t)\) moves with speed \(|\dot{q}(t)|=|p(t)|\). Since
\(H=\frac{1}{2}|p|^2\) is conserved along the Hamiltonian flow, this
speed is constant. In other words, the geodesic is parametrised
proportionally to arc-length, where the constant of proportionality
depends on the level set of \(H\).

\end{Remark}
\section{Noether's theorem}
\label{sct:noether}

The Hamiltonian formalism assigns a Hamiltonian flow
\(\phi^H_t\) to a function \(H\) on a symplectic manifold; the
function \(H\) is conserved along the flow in the sense that
\(H(\phi^H_t(x))=H(x)\) for all \(t\). This is responsible for
Noether's famous correspondence\index{Noether's theorem|(} between
symmetries and conserved quantities. We reiterate here that
\(H\) does not need to correspond to ``energy'': it can be any
function.

\begin{Example}[Translation]\label{exm:translation}
Consider the symplectic manifold \(\RR^2\) with coordinates
\((p,q)\) and symplectic form \(dp\wedge dq\). This is the
phase space\index{phase space} of a particle on a line, where
\(q\) denotes the position of the particle and \(p\) its
momentum\index{momentum}. The Hamiltonian vector field
associated with the function \(p\) is \[\dot{p}=0,\quad
\dot{q}=1,\] which generates the flow
\(\phi^p_t(p,q)=(p,q+t)\). This is a
translation\index{translation} in the \(q\)-direction.

\end{Example}
More generally, if \(\xi\) is a vector field on \(Q\) and
\(\psi_t\colon Q\to Q\) is its flow, Noether's theorem gives an
explicit Hamiltonian on \(T^*Q\) generating the 1-parameter
family of symplectomorphisms \((\psi_t)_*\colon T^*Q\to T^*Q\)
defined in Lemma \ref{lma:diff_action_on_cotangent}:

\begin{Theorem}[Noether's theorem]\label{thm:noeth}
Let \(\xi\) be a vector field on
\(Q\) and let \(\psi_t\) be its flow. Define \(H_\xi\colon
T^*Q\to\RR\) by \(H_\xi(\eta)=\eta(\xi(q))\) for \(\eta\in
T^*_qQ\). If \(\phi^{H_\xi}_t\) is the Hamiltonian flow of
\(H_\xi\) then \(\phi^{H_\xi}_t=(\psi_t)_*\), where
\((\psi_t)_*\) is the action of \(\psi_t\) on \(T^*Q\) defined
in Lemma \ref{lma:diff_action_on_cotangent}.

\end{Theorem}
\begin{Remark}
Let \(\mathfrak{diff}(Q)\) denote the space of vector fields
on \(Q\) (the Lie algebra of \(\OP{Diff}(Q)\)) and define the
map \(\mu\colon T^*Q\to\left(\mathfrak{diff}(Q)\right)^*\) by
\(\mu(\eta)(\xi)=H_\xi(\eta)\). In the language of Chapter
\ref{ch:glob}, Theorem \ref{thm:noeth} can be phrased by
saying that \(\mu\) is a moment map\index{moment map!for
diffeomorphism group action on cotangent bundle} for the
Hamiltonian \(\OP{Diff}(Q)\) action on
\(T^*Q\)\index{diffeomorphism group!action on cotangent
bundle}.

\end{Remark}
\begin{Proof}
Write \(\eta = \sum p_i\, dq^i\). We have
\(H_\xi(p,q)=\eta(\xi(q))=\sum p_i\xi^i\). Both \(\phi^{H_\xi}_t\)
and \((\psi_t)_*\) are generated by vector fields. It suffices to
check that these vector fields coincide.

By definition, \(\phi^{H_\xi}_t\) is generated by the
Hamiltonian vector field \(V_{H_\xi}\). We claim that
\[V_{H_\xi}:=-\sum_{i,j} p_i\frac{\partial\xi^i}{\partial
q^j}\frac{\partial}{\partial
p_j}+\sum_j\xi^j\frac{\partial}{\partial q^j}.\] To see this,
observe that
\begin{align*}
\iota_{V_{H_\xi}}\sum dp_k\wedge dq^k&=\sum dp_k(V_{H_\xi})\,dq^k-\sum
dq^k(V_{H_\xi})\,dp_k\\ &=-\sum p_i\frac{\partial\xi^i}{\partial q^k}\,dq^k-\sum
\xi^k\,dp_k = -dH_\xi.
\end{align*}
Now let us calculate the infinitesimal action of \(\xi\) on
\(T^*Q\). Suppose that \(\psi_t\) is the flow of \(\xi\); in
coordinates, \((\psi_t)_*(p_i,q^i)=\left(\sum_j
p_j\frac{\partial(\psi_t^{-1})^j}{\partial
q^i},\psi_t^i(q)\right)\). The infinitesimal action of \(\xi\)
on \(T^*Q\) is given by
\begin{align*}
\left.\frac{d}{dt}\right|_{t=0}(\psi_t)_*(p_i,q^i)&=\left.\frac{d}{dt}\right|_{t=0}\left(\sum_j p_j\frac{\partial(\psi_t^{-1})^j}{\partial q^i},\psi_t^i(q)\right)\\
&=\left(-\sum p_j\frac{\partial\xi^j}{\partial q^i},\xi^i\right)
\end{align*}
This is just another way of writing \(V_{H_\xi}\), so the theorem is
proved. In getting to the final line, we used the fact that
\(\left.\frac{d}{dt}\right|_{t=0}\psi_t^{-1} = -\xi\), i.e.\
the inverse of \(\psi_t\) is obtained by flowing backwards
along \(\xi\) for time \(t\). \qedhere

\end{Proof}
\begin{Example}[Angular momentum, Exercise \ref{exr:angular_momentum}]\label{exm:angular_momentum}
Suppose \(Q=\RR^3\) with coordinates \(q^1,q^2,q^3\), and
consider the 1-parameter family of diffeomorphisms
\[\psi_t(q^1,q^2,q^3)=(q^1\cos t- q^2\sin t,q^1\sin t+q^2\cos
t,q^3)\] given by rotating around the \(q^3\)-axis. Find the
Hamiltonian on \(T^*Q\) which generates
\((\psi_t)_*\).\index{angular momentum}

\end{Example}
\begin{Remark}
For those who have encountered Noether's theorem in the
context of classical field theory, the field theory version is
proved in the same way, where we take \(Q\) to be the space of
fields and \(T^*Q\) to be the phase space.\index{Noether's theorem|)}

\end{Remark}
\section{Lagrangian submanifolds and the Hamilton-Jacobi equation}

There are some easy examples of Lagrangian
submanifolds\index{Lagrangian!submanifold|(} in cotangent
bundles.

\begin{Example}
The {\em zero-section}\index{zero-section} is the submanifold
which intersects every cotangent fibre at the zero
covector. This is Lagrangian: in local canonical coordinates
\((p,q)\), with \(\omega=\sum_i dp_i\wedge dq^i\), it is given
by \(p=0\). Dually, the cotangent fibres\index{cotangent
fibre} \(q=\mbox{const}\) are also Lagrangian submanifolds.

\end{Example}
\begin{Example}
Suppose \(H\colon T^*Q\to\RR\) is the Hamiltonian from
Definition \ref{dfn:cogeodesic_flow} generating the cogeodesic
flow on \(T^*Q\) for some metric. What are the geodesics
connecting \(x\in Q\) to \(y\in Q\) in time \(t\)? They are
in bijection with the intersection points between the
Lagrangian submanifolds \(\phi^H_t(T^*_xQ)\) and
\(T^*_yQ\). In this way, Lagrangian submanifolds can be used
to impose initial/terminal conditions on geodesics or other
Hamiltonian systems. The utility of this stems from the fact
there is a variational interpretation for Hamilton's equations
with Lagrangian boundary conditions (the Hamiltonian
trajectories are critical points for the {\em action
functional}). This is the point of departure for applications
of Floer theory\index{Floer theory} to symplectic geometry.

\end{Example}
Recall that a {\em section} of the cotangent bundle is a map
\(\eta\colon Q\to T^*Q\) such that \(\eta(q)\in T_q^*Q\). This
is the same thing as a 1-form on \(Q\). We define the {\em
graph} of a 1-form \(\eta\) to be the image of the corresponding
section \(\eta(Q)=\{(q,\eta(q)\in T^*Q\,:\,q\in Q\}\subset
T^*Q\). This is a submanifold diffeomorphic to \(Q\).

\begin{Lemma}
The graph\index{graph of a closed 1-form} of a 1-form \(\eta\)
is Lagrangian if and only if \(d\eta=0\), i.e.\ \(\eta\) is
closed.
\end{Lemma}
\begin{Proof}
Let \(q^i\) be local coordinates on \(Q\). The tangent space to
\(\eta(Q)\) at a point living over this coordinate patch is spanned
by the vectors
\[\eta_*(\partial_{q^i})=\frac{\partial\eta_j}{\partial
q^i}\partial_{p_j}+\partial_{q^i}.\] We have \[\left(\sum_m
dp_m\wedge
dq^m\right)(\eta_*\partial_{q^k},\eta_*\partial_{q^\ell})=
\frac{\partial\eta_\ell}{\partial q^k} -
\frac{\partial\eta_k}{\partial q^\ell}.\] This is the \(dq^k\wedge
dq^\ell\)-component of \(d\eta\). (Compare this with the proof of
Theorem \ref{thm:actionangle}.) \qedhere

\end{Proof}
Note that if \(L\subset T^*Q\) is transverse to the cotangent
fibres near some point \(x\in L\) then, locally near \(x\),
\(L\) is the graph of some section (by the inverse function
theorem). Moreover, locally, any closed 1-form \(\eta\) admits
an antiderivative, that is a function \(S\) such that \(\eta =
dS\). We often call such a function \(S\) a (local) {\em
generating function}\index{generating function} for
\(\eta(Q)\). So to describe a Lagrangian submanifold of
\(T^*Q\), away from points where it is tangent to cotangent
fibres, it is sufficient to give a collection of local
generating functions.

If a Lagrangian submanifold is allowed to evolve under a
Hamiltonian flow, then (up to a time-dependent constant shift)
its local generating functions evolve according to a
differential equation called the {\em Hamilton-Jacobi
equation}\index{Hamilton-Jacobi equation}. We state this in its
simplest form for Lagrangians that admit a global generating
function, i.e.\ Lagrangians which are the graph of an exact
1-form\index{graph of an exact 1-form}.

\begin{Theorem}[Hamilton-Jacobi equation]\label{thm:HJ}
Let \(L=\OP{graph}(dS)\subset T^*Q\) be a Lagrangian
submanifold which is the graph of an exact 1-form \(dS\). Let
\(H_t\colon T^*Q\to\RR\) be a time-dependent Hamiltonian. If
\(S_t\) is a solution to the {\em Hamilton-Jacobi
equation}\index{Hamilton-Jacobi equation} \begin{equation}\label{eq:HJ}\frac{\partial
S_t}{\partial t}=-H_t\left(\frac{\partial S_t}{\partial
\bm{q}},\bm{q}\right),\qquad S_0=S,\end{equation} then
\[\OP{graph}(dS_t)=\phi^{H_t}_t(\OP{graph}(dS)).\] Conversely,
if \(\phi^{H_t}_t(\OP{graph}(dS))=\OP{graph}(dF_t)\) then
\(F_t=S_t+c(t)\) where \(S_t\) solves Equation \eqref{eq:HJ}
and \(c(t)\) is a time-dependent constant.
\end{Theorem}
\begin{Proof}
Let \(dS\colon Q\to T^*Q\) be the section corresponding to the
1-form \(dS\) and let \(i_t\colon Q\to T^*Q\) be the
Lagrangian inclusion\index{Lagrangian!submanifold|)} of
\(\phi^{H_t}_t(dS(Q))\) defined by
\(i_t(q)=\phi^{H_t}_t(dS(q))\). Pick local canonical
coordinates \((\bm{p},\bm{q})\) and write \(i_{t +
\epsilon}(\bm{q}) = i_t(\bm{q}) + \epsilon v_t(q)\) for some
vector field \(v_t\) along \(i_t(Q)\). To first order in
\(\epsilon\), \(v_t = V_{H_t}\), that is:
\[v_t=\left(-\frac{\partial H}{\partial \bm{q}} +
\mathfrak{o}(\epsilon), \frac{\partial H}{\partial
\bm{p}}+\mathfrak{o}(\epsilon)\right)\] where we write
\(\mathfrak{o}(\epsilon)\) for any terms such that
\(\lim_{\epsilon\to 0}|\mathfrak{o}(\epsilon)| = 0\).

\begin{figure}[htb]
\begin{center}
\begin{tikzpicture}
\draw (0,0) to[out=80,in=180] (3,3) to[out=0,in=180] (6,-1) to[out=0,in=-135] (8,1);
\draw (-0.5,0) to[out=80,in=180] (2.5,1.5) to[out=0,in=180] (6,-0.5) to[out=0,in=-135] (8,2);
\draw[thick,->] (2.5,1.5) -- (4,2.5);
\draw[dotted] (2.5,1.5) -- (2.5,-1);
\draw[dotted] (4,2.5) -- (4,-1);
\draw[dotted] (5,1.5) -- (0,1.5);
\draw[dotted] (5,2.5) -- (0,2.5);
\node at (2.5,-1) [below] {\(\bm{q}\)};
\node at (4.5,-1) [below] {\(\bm{q}+\epsilon\frac{\partial H}{\partial\bm{p}}+\epsilon\mathfrak{o}(\epsilon)\)};
\node at (0,1.5) [left] {\(\frac{\partial S_t}{\partial\bm{q}}\)};
\node at (0,2.5) [left] {\(\frac{\partial S_{t+\epsilon}}{\partial\bm{q}}(\bm{q}+\epsilon\frac{\partial H_t}{\partial\bm{p}}+\epsilon\mathfrak{o}(\epsilon))\)};
\node at (5,2) [right] {\(-\epsilon\frac{\partial H_t}{\partial\bm{q}}+\epsilon\mathfrak{o}(\epsilon)\)};
\draw[decorate,decoration={brace,amplitude=4pt}] (4.8,2.5) -- (4.8,1.5);
\node at (8,2) [right] {\(\OP{graph}(dS_t)\)};
\node at (8,1) [right] {\(\OP{graph}(dS_{t+\epsilon})\)};
\node at (3,2.2) {\(\epsilon v_t\)};
\end{tikzpicture}
\end{center}
\caption{The graphs of \(dS_t\) and \(dS_{t+\epsilon}\) differ by \(\epsilon v_t =\epsilon (V_{H_t}+\mathfrak{o}(\epsilon))\).}
\label{fig:hj}
\end{figure}
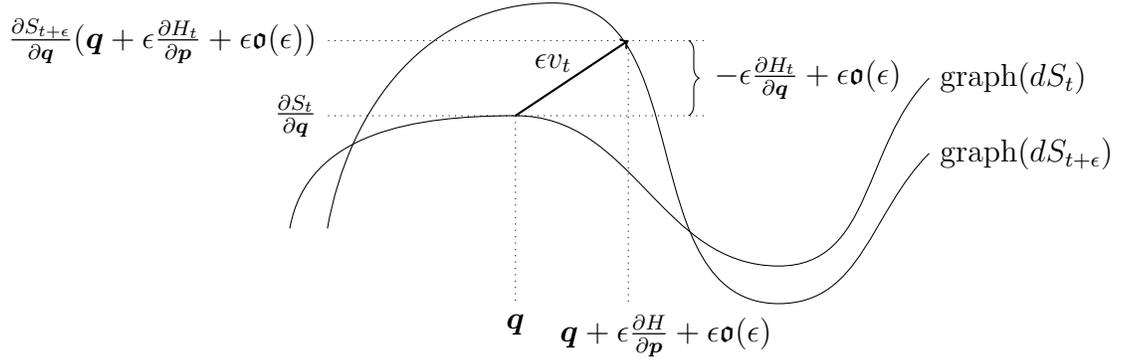

We have \[\frac{\partial S_{t+\epsilon}}{\partial
q^i}\left(\bm{q}+\epsilon\frac{\partial H}{\partial
\bm{p}}+\mathfrak{o}(\epsilon)\right) = \frac{\partial S_t}{\partial
q^i}(\bm{q}) - \epsilon\frac{\partial H_t}{\partial
q^i}\left(\frac{\partial
S_t}{\partial\bm{q}}(\bm{q}),\bm{q}\right)+\mathfrak{o}(\epsilon)\]
and, by Taylor expanding, we also have \[\frac{\partial
S_{t+\epsilon}}{\partial q^i}\left(\bm{q}+\epsilon\frac{\partial
H}{\partial \bm{p}}+\mathfrak{o}(\epsilon)\right) = \frac{\partial
S_t}{\partial q^i}(\bm{q}) +
\epsilon\left(\frac{\partial^2S_{t}}{\partial t\partial
q^i}(\bm{q})+\sum_j\frac{\partial^2 S_t}{\partial q^i\partial
q^j}(\bm{q})\frac{\partial H_t}{\partial
p_j}\right)+\mathfrak{o}(\epsilon).\] Comparing terms of order
\(\epsilon\), we get \[\frac{\partial^2 S_t}{\partial q^i\partial
t} + \frac{\partial H_t}{\partial q^i} + \sum_j\frac{\partial
S_t}{\partial q^j}\frac{\partial H_t}{\partial p_j} = 0,\] where the
derivatives of \(H_t\) are evaluated at \((\partial S_t/\partial
\bm{q},\bm{q})\). In particular, this means that \(\frac{\partial
H_t}{\partial q^i} + \sum_j\frac{\partial S_t}{\partial
q^j}\frac{\partial H_t}{\partial p_j}=\frac{\partial}{\partial
q^i}((dS_t)^*H_t)\), so the equation is telling us that
\[\frac{\partial S_t}{\partial t} + (dS_t)^*H_t\] is constant on
\(Q\), say equal to \(C(t)\). This means that \(S_t\) satisfies the
equation \[\frac{\partial S_t}{\partial t} =
-H_t\left(\frac{\partial S_t}{\partial \bm{q}},\bm{q}\right) +
C(t),\] which means that \(S_t-c(t)\) satisfies the Hamilton-Jacobi
equation provided \(\dot{c}(t)=C(t)\). \qedhere

\end{Proof}
\begin{Remark}
The proof is purely local, and therefore also works when the
generating function is local, but it is trickier to state in that
case because the domain of the local generating function changes.

\end{Remark}
\section{Solutions to inline exercises}

\begin{Exercise}[Remark \ref{rmk:cot_as_ham_sys}]\label{exr:cot_as_ham_sys}
Pick local coordinates \(q^i\) on a patch in \(Q\) and
consider the Hamiltonian system
\((q^1\circ\pi,\ldots,q^n\circ\pi)\) on the \(\pi\)-preimage
of this patch. Show that the canonical coordinates \(p_i\) are
{\em minus} the Liouville coordinates associated with the
global Lagrangian section given by the zero-section. Does the
zero-section inherit an integral affine structure?
\end{Exercise}
\begin{Solution}
In local coordinates \((\bm{p},\bm{q})\),
\(\bm{p}=(p_1,\ldots,p_n)\) and \(\bm{q}=(q^1,\ldots,q^n)\),
with \(\omega=\sum dp_i\wedge dq^i\), the Hamiltonian flow of
\(q^i\) is translation in the \(-p_i\) direction, and the
zero-section is given by \(\sigma(\bm{q})=(\bm{q},0)\). We
have \(\phi^{\bm{q}}_{\bm{t}}(\bm{q},0)=(\bm{q},-\bm{t})\),
which shows that the \(-p_i\) are Liouville coordinates. Since
the fibres of \(\pi\) are \(\RR^n\), there are no periodic
orbits, so the period lattice is \(0\) in each
fibre. Therefore there is no natural integral affine
structure: that construction would need the period lattice to
have full rank. \qedhere

\end{Solution}
\begin{Exercise}[Remark \ref{rmk:parallels_differences}]\label{exr:parallels_differences}
Explain the parallels and differences between Lemmas
\ref{lma:globcoords} and \ref{lma:diff_action_on_cotangent}.
\end{Exercise}
\begin{Solution}
By Exercise \ref{exr:cot_as_ham_sys}, we can think of
\(\pi\colon T^*Q\to Q\) as a Hamiltonian system. In both
lemmas, we assume the existence of a diffeomorphism between
the images of our Hamiltonian systems: for Lemma
\ref{lma:globcoords} we have \(\phi\colon \bm{F}(X)\to
\bm{G}(X)\) and for Lemma \ref{lma:diff_action_on_cotangent}
we have \(\psi\colon Q\to Q\). In both cases, we obtain a
symplectomorphism between the total spaces: respectively
\(\Phi\colon X\to Y\) and \((\psi)_*\colon T^*Q\to
T^*Q\). Moreover, these symplectomorphisms are given by the
same formula in Liouville coordinates.

The difference is that \(\phi\) is required to be an integral
affine transformation, whereas \(\psi\) can be any
diffeomorphism. This is because the Hamiltonian systems
\(\bm{F}\) and \(\bm{G}\) have period lattices of rank \(n\),
and the derivative of \(\phi\) is required to preserve these
period lattices, which tells us that \(\phi\) is integral
affine. The period lattice for \(\pi\) is trivial, so there is
no constraint on \(d\psi\). \qedhere

\end{Solution}
\begin{Exercise}[Angular momentum, Example \ref{exm:angular_momentum}]\label{exr:angular_momentum}
Suppose \(Q=\RR^3\) with coordinates \(q^1,q^2,q^3\), and consider
the 1-parameter family of diffeomorphisms
\[\psi_t(q^1,q^2,q^3)=(q^1\cos t- q^2\sin t,q^1\sin t+q^2\cos
t,q^3)\] given by rotating around the \(q^3\)-axis. Find the
Hamiltonian on \(T^*Q\) which generates
\((\psi_t)_*\).\index{angular momentum}
\end{Exercise}
\begin{Solution}
The flow \(\psi_t\) is generated by the vector field
\(\xi=(-q_2,q_1,0)\), so by Theorem \ref{thm:noeth}, the
induced symplectomorphism on \(T^*Q\) is generated by the
Hamiltonian \[H_\xi(p,q)=p_2q_1-p_1q_2.\] This is the usual
formula for the component of angular momentum around the
\(q^3\)-axis.\index{cotangent bundle|)}\qedhere

\end{Solution}
\chapter{Moser's argument}
\label{ch:moser}
\thispagestyle{cup}

At various points in the book, we have appealed to {\em the
Moser argument}\index{Moser argument} \index{Moser trick|see
{Moser argument}}. This is a famous and extremely useful trick,
first introduced by Moser \cite{MoserVolume}. We include a proof
here for completeness. When we say a ``family of \(k\)-forms'',
we mean a \(k\)-form whose coefficients (with respect to any
local coordinate system) depend continuously-differentiably on a
parameter \(t\).

\begin{Theorem}[Moser's argument]\label{thm:moser}
Suppose that \(X\) is a manifold and \(\omega_t\) is a family
of symplectic forms. If \(d\omega_t/dt = d\sigma_t\) for some
family of compactly-supported 1-forms \(\sigma_t\) then there
is a family of diffeomorphisms \(\phi_t\) with
\(\phi_0=\OP{id}\) and \(\phi_t^*\omega_t=\omega_0\).
\end{Theorem}
\begin{Proof}
Let \(V_t\) be the vector field \(\omega_t\)-dual to
\(-\sigma_t\), that is \(\iota_{V_t}\omega_t =
-\sigma_t\). This is a compactly-supported vector field, so we
can define the flow along \(V_t\). The flow is a 1-parameter
family of diffeomorphisms \(\phi_t\) satisfying
\(\phi_0=\OP{id}\) and \(\frac{d\phi_t(x)}{dt} =
V_t(\phi_t(x))\). We will differentiate \(\phi_t^*\omega_t\)
with respect to \(t\) and show that the result is zero. This
will imply that \(\phi_t^*\omega_t\) is independent of \(t\),
and hence equal to \(\omega_0\).
\begin{align*}
\frac{d}{dt}\phi_t^*\omega_t &= \phi_t^*\left(\Lie_{V_t}\omega_t\right) +
\phi_t^*\frac{d\omega_t}{dt}\\
&= \phi_t^*(d\iota_{V_t}\omega_t) - \phi_t^*d\sigma_t\\
&= \phi_t^*(d\sigma_t - d\sigma_t) = 0
\end{align*}
where we used Cartan's formula\index{Cartan's formula}
\(\Lie_{V_t}\omega_t =
d\iota_{V_t}\omega_t+\iota_{V_t}d\omega_t\) and the fact that
\(d\omega_t = 0\).\qedhere

\end{Proof}
\chapter{Toric varieties revisited}
\label{ch:toric_varieties}
\thispagestyle{cup}

In this appendix, we will construct the toric
variety\index{toric variety|(} associated to a convex rational
polytope\index{polytope!convex rational} using only algebraic
geometry (no symplectic cuts\index{symplectic cut}). Since most
expositions of toric geometry (for example, Danilov
\cite{Danilov} or Fulton \cite{Fulton}) start from the dual
(fan) picture, and we are aiming to give alternative viewpoints
wherever possible, we will confine ourselves to work only with
the moment polytope\index{moment polytope}. Throughout this
appendix we will make use of homogeneous coordinates; see
Appendix \ref{ch:complex_projective_spaces} for a rapid
overview.

\section{Construction}

Let \(\Delta\subset\RR^n\) be a compact Delzant
polytope\index{polytope!Delzant}. We will focus on the special
case where the vertices of \(\Delta\) have integer coordinates
and explain how to construct the manifold \(X_\Delta\) whose
existence is guaranteed by Delzant's existence
theorem\index{Delzant!existence theorem}, Theorem
\ref{thm:convexity}(2).

\begin{Theorem}\label{thm:toric_varieties}
Suppose that \(\bm{p}_1,\ldots,\bm{p}_N\in\ZZ^n\) are the
integer lattice points contained in a Delzant polytope
\(\Delta\), and write \(\bm{p}_i=(p_{i1},\ldots,p_{in})\). Let
\(\bm{z}^{\bm{p}_i}\) be the monomial
\(z_1^{p_{i1}}z_2^{p_{i2}}\cdots z_n^{p_{in}}\). Consider the
map \[F_\Delta\colon(\CC^*)^n\to\cp{N-1},\qquad
F_\Delta(\bm{z})=[\bm{z}^{\bm{p}_1} : \cdots :
\bm{z}^{\bm{p}_N}].\] Let \(X_\Delta\) be the
Zariski-closure\index{Zariski-closure} of the image of
\(F_\Delta\). Then \(X_\Delta\) is a smooth projective
variety. Let \(P\colon\RR^N\to\RR^n\) be the linear projection
given by right-multiplication with the matrix
\[\scaleleftright[1.75ex]{<}{\begin{matrix} p_{11} & p_{21}&
\cdots & p_{n1}\\ p_{12} & \ddots& & p_{n2}\\ \vdots & &
\ddots & \vdots \\ p_{1N} & \cdots & \cdots &
p_{Nn}\end{matrix}}{)}.\] If \(\mu\colon\cp{N}\to\RR^N\) is
the moment map\index{moment map} for the standard
\(T^N\)-action\index{Hamiltonian!torus action} then
\(\mu|_{X_{\Delta}}\cdot P\colon X_\Delta\to\RR^n\) is the
moment map for a \(T^n\)-action on \(X_\Delta\) whose moment
image is \(\Delta\).

\end{Theorem}
\begin{Definition}
The projective variety \(X_\Delta\) is called the {\em
projective toric variety}\index{projective variety!toric}
associated to the polytope \(\Delta\).

\end{Definition}
Before proving this theorem, we will work out some examples.

\section{Examples}

\begin{Example}\label{exm:segre}
Suppose \(\Delta\) is the square with vertices \((0,0)\),
\((1,0)\), \((0,1)\), and \((1,1)\). Since this is a square,
Delzant's uniqueness theorem tells us that \(X_\Delta\) will
be \(S^2\times S^2\). We will confirm that this is the output
of Theorem \ref{thm:toric_varieties}.

\begin{center}
\begin{tikzpicture}
\filldraw[fill=lightgray,opacity=0.5,draw=none] (0,0) -- (0,2) -- (2,2) -- (2,0) -- cycle;
\draw[black,thick] (0,0) -- (0,2) -- (2,2) -- (2,0) -- cycle;
\node at (0,0) {\(\bullet\)};
\node at (2,0) {\(\bullet\)};
\node at (0,2) {\(\bullet\)};
\node at (2,2) {\(\bullet\)};

\end{tikzpicture}
\end{center}
Theorem \ref{thm:toric_varieties} tells us to consider the map
\[F_\Delta\colon(\CC^*)^2\to\cp{3},\qquad
F_\Delta(z_1,z_2)=[1:z_1:z_2:z_1z_2].\] If
\([x_1:x_2:x_3:x_4]\) are our homogeneous coordinates on
\(\cp{3}\) then we see that the image of \(F_\Delta\) is
contained (as a Zariski-dense subset) in the
subvariety\index{subvariety} \(V=\{x_1x_4=x_2x_3\}\). This
subvariety is a smooth quadric\index{quadric hypersurface!Segre} surface and it
is the Zariski-closure of the image of \(F_\Delta\). Note that
\(F_\Delta\) is the restriction of the Segre embedding
\[\cp{1}\times\cp{1}\to\cp{3},\qquad ([a:b],[c:d])\mapsto
[ac:bc:ad:bd]\] to the affine chart \(a=c=1\), and \(V\) is
the image of the Segre embedding\index{Segre embedding}. Since
\(\cp{1}\cong S^2\), this confirms that \(X_\Delta=S^2\times
S^2\).

The matrix \(P\) is
\[\scaleleftright[1.75ex]{<}{\begin{matrix} 0 & 0\\ 1 & 0 \\ 0
& 1 \\ 1 & 1\end{matrix}}{)}\] and the moment map \(\mu\) is
\[\left(\frac{1}{2}\frac{|x_1|^2}{|x|^2},\ \
\frac{1}{2}\frac{|x_2|^2}{|x|^2},\ \
\frac{1}{2}\frac{|x_3|^2}{|x|^2},\ \
\frac{1}{2}\frac{|x_4|^2}{|x|^2}\right),\] so \[\mu(x)\cdot
P=\left(\frac{1}{2}\frac{|x_2|^2+|x_4|^2}{|x|^2},\ \
\frac{1}{2}\frac{|x_3|^2+|x_4|^2}{|x|^2}\right).\]
Precomposing with the Segre embedding to get a function on
\(X_\Delta = \cp{1}\times \cp{1}\), we get
\[\left(\frac{1}{2}\frac{|b|^2}{|a|^2+|b|^2},\ \
\frac{1}{2}\frac{|d|^2}{|c|^2+|d|^2}\right).\] The first
(respectively second) component is the Hamiltonian generating
the standard circle action on
the first (respectively second) factor \(\cp{1}\) (see Example
\ref{exm:cpn}).

The moment image is the convex hull of the moment images of
the fixed points. The fixed points are \(([1:0],[1:0])\),
\(([1:0],[0:1])\), \(([0:1],[1:0])\), and \(([0:1],[0:1])\),
whose images are \((0,0)\), \((1,0)\), \((0,1)\) and \((1,1)\)
respectively. Therefore the moment image of \(X_\Delta\) is
\(\Delta\).

\end{Example}
\begin{Example}\label{exm:cp2_triangle}
If \(\Delta\) is the triangle with vertices \((0,0)\),
\((1,0)\) and \((0,1)\) then
\(F_\Delta\colon(\CC^*)^2\to\cp{2}\) is the map
\(F_\Delta(z_1,z_2)=[1:z_1:z_2]\). The image of \(F_\Delta\)
is dense in \(\cp{2}\), so \(X_\Delta=\cp{2}\). By Lemmas
\ref{lma:sphere_area} and \ref{lma:selfint}, the preimage of
an edge is a symplectic sphere with area \(2\pi\) and
self-intersection \(1\); this is a line in \(\cp{2}\).

\end{Example}
\begin{figure}
\begin{center}
\begin{tikzpicture}
\filldraw[fill=lightgray,opacity=0.5,draw=none] (0,0) -- (0,2) -- (2,0) -- cycle;
\draw[black,thick] (0,0) -- (0,2) -- (2,0) -- cycle;
\node at (0,0) {\(\bullet\)};
\node at (2,0) {\(\bullet\)};
\node at (0,2) {\(\bullet\)};
\begin{scope}[shift={(4,0)}];
\filldraw[fill=lightgray,opacity=0.5,draw=none] (0,0) -- (0,4) -- (4,0) -- cycle;
\draw[black,thick] (0,0) -- (0,4) -- (4,0) -- cycle;
\node at (0,0) {\(\bullet\)};
\node at (2,0) {\(\bullet\)};
\node at (4,0) {\(\bullet\)};
\node at (0,2) {\(\bullet\)};
\node at (2,2) {\(\bullet\)};
\node at (0,4) {\(\bullet\)};
\end{scope}
\end{tikzpicture}
\end{center}
\caption{The polygons \(\Delta\) and \(2\Delta\) for Examples \ref{exm:cp2_triangle} and \ref{exm:cp2_triangle_2}.}
\end{figure}
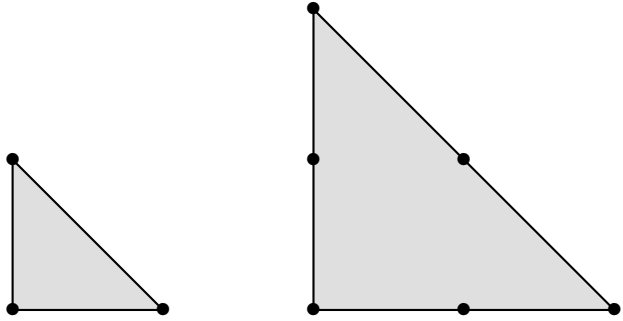

If we rescale both the Fubini-Study \index{Fubini-Study form|see
{symplectic form, Fubini-Study}} form on \(\cp{2}\) and the
moment map\index{moment polytope!rescaling} for the torus action
by a factor of \(2\) then we get a moment map whose image is the
triangle \(2\Delta\) with vertices \((0,0)\), \((2,0)\) and
\((0,2)\). Delzant's uniqueness theorem tells us that
\(X_{2\Delta}\cong (\cp{2},2\omega_{FS})\). But the isomorphism
is not obvious from the construction:

\begin{Example}\label{exm:cp2_triangle_2}
If \(2\Delta\) is the triangle with vertices \((0,0)\),
\((2,0)\) and \((0,2)\) then there are now 6 integer points in
\(2\Delta\), and we get
\[F_\Delta\colon(\CC^*)^2\to\cp{5},\qquad
F_\Delta(z_1,z_2)=[1:z_1:z_1^2:z_2:z_1z_2:z_2^2].\] The map
\(F_\Delta\) factors through the {\em quadratic Veronese
embedding}\index{Veronese embedding}
\[\mathcal{V}\colon\cp{2}\to\cp{5},\qquad \mathcal{V}([a:b:c])
= [a^2:ab:b^2:ac:bc:c^2]\] by taking \(a=1\), \(b=z_1\),
\(c=z_2\), and the image of \(F_\Delta\) is dense inside
\(\mathcal{V}(\cp{2})\). Thus \(X_\Delta=\cp{2}\). Again, the
preimage of an edge is a line in \(\cp{2}\), but it has
symplectic area \(4\pi\) because the pullback of the
Fubini-Study form along the quadratic Veronese embedding is
symplectomorphic to {\em twice} the Fubini-Study form on
\(\cp{2}\) (a hyperplane of \(\cp{5}\) intersects
\(\mathcal{V}(\cp{2})\) in a conic, not a line).

\end{Example}
Rescaling the polytope by a factor of \(k\) always corresponds
to reimbedding via a Veronese map of degree
\(k\).

Finally, let us try to apply the construction from Theorem
\ref{thm:toric_varieties} when \(\Delta\) is not
Delzant\index{polytope!non-Delzant|(}. The corresponding toric
variety will have singularities living over the non-Delzant
points of \(\Delta\).

\begin{Example}
Let \(\Delta\) be the non-Delzant polygon with vertices
\((0,0)\), \((0,1)\) and \((2,1)\) from Figure
\ref{fig:non_delzant}. This additionally contains the integer
point \((1,1)\). We therefore get
\[F_\Delta\colon(\CC^*)^2\to\cp{3},\qquad F_\Delta(z_1,z_2) =
[1:z_2:z_1z_2:z_1^2z_2].\] In homogeneous coordinates
\([x_1:x_2:x_3:x_4]\) this satisfies the equation
\(x_2x_4=x_3^2\). This is a singular
quadric\index{quadric hypersurface!nodal} surface with an ordinary double
point at \([1:0:0:0]\). Under the moment
map\index{moment polytope!nodal quadric surface} \(\mu\cdot
P\), this point projects to the origin, which is precisely the
point where \(\Delta\) fails to be Delzant.

By Lemma \ref{lma:selfint}, the preimage of the horizontal
edge is a symplectic sphere with square \(2\). In homogeneous
coordinates, this is the conic\index{conic} \(x_2x_4=x_3^2\)
in the plane \(x_1=0\).

\end{Example}
\begin{Remark}
The ordinary double point is the cyclic quotient singularity
\(\frac{1}{2}(1,1)\). The germ of our non-Delzant polygon near
the origin agrees with the germ of the non-Delzant polygon
from Example \ref{exm:cycquot} with \(n=2\), \(a=1\). This is
a general fact: one can read off the singularities of
\(X_\Delta\) from the non-Delzant points in
\(\Delta\).\index{polytope!non-Delzant|)}

\end{Remark}
\section{Proof of Theorem \ref{thm:toric_varieties}}

Consider the \(T^n\)-actions
\begin{align*}
e^{i\bm{t}}\bm{z} &= (e^{it_1}z_1,\ldots,e^{it_n}z_n)\\
e^{i\bm{t}}\star[Z_1:\ldots:Z_N] &= \left[e^{i\left(p_{11}t_1+ \cdots+p_{1n}t_n\right)}Z_1 : \cdots :
e^{i\left(p_{N1}t_1+\cdots+p_{Nn}t_n\right)}Z_N\right]
\end{align*}
on \((\CC^*)^n\) and \(\cp{N}\) respectively. The action denoted
by \(\star\) is generated by the Hamiltonian \(\mu\cdot P\). The
map \(F_\Delta\) intertwines the actions in the sense that
\(F_\Delta(e^{i\bm{t}}\bm{z})=e^{i\bm{t}}\star
F_{\Delta}(\bm{z})\); this means that \(e^{i\bm{t}}\star\)
preserves the image of \(F_\Delta\), and hence its
Zariski-closure \(X_\Delta\).

It remains to show that \(X_\Delta\) is smooth and that the
moment image agrees with \(\Delta\). We will start by writing
down equations for \(X_\Delta\). Let \([Z_1:\cdots:Z_N]\) be
homogeneous coordinates on \(\cp{N-1}\). Note that each
coordinate \(Z_i\) corresponds to an integer lattice point
\(\bm{p}_i\in\Delta\).

\begin{Lemma}\label{lma:equations}
Let \(a_1,\ldots,a_N\) be integers. If the relation \(\sum_i
a_i\bm{p}_i = 0\) holds then the equation \[\prod_{a_i\geq
0}Z_i^{a_i}=\prod_{a_i<0}Z_i^{-a_i}\] holds on \(X_\Delta\).
\end{Lemma}
\begin{Proof}
This holds on the image of \(F_\Delta\) because it translates
to \(\prod_i z^{\sum a_i\bm{p}_i}=z^0=1\). It therefore holds
on the Zariski-closure of the image of \(F_\Delta\), which is
\(X_\Delta\) by definition. \qedhere

\end{Proof}
Let \(V_\Delta\) be the subvariety cut out by the equations
coming from Lemma \ref{lma:equations}. The lemma shows that
\(X_\Delta\subset V_\Delta\). We will show that \(V_\Delta\) is
a smooth variety containing the image of \(F_\Delta\) as a
Zariski-open set, which will show that \(V_\Delta=X_\Delta\) (in
particular, it will show that \(X_\Delta\) is smooth).

Let \(\vtx = \{v \in \{1,\ldots,N\}\,:\,\bm{p}_v \mbox{
is a vertex of }\Delta\}\). Note that every integer lattice
point in \(\Delta\) can be written as a linear combination
\(\sum_{j\in\vtx} a_j\bm{p}_j\) with \(a_j\in\QQ_{\geq
0}\) for all \(j\in\vtx\).

\begin{Corollary}
For each \(i\in\{1,\ldots,N\}\), write
\(\bm{p}_i=\sum_{j\in\vtx(\Delta)}a_j\bm{p}_j\) with
\(a_j\in\QQ_{\geq 0}\). Let \(\vtx_i=\{j\in\vtx\,:\,a_j\neq
0\}\subset\vtx\). The open set \(V_{\Delta}\cap\{Z_i\neq
0\}\) is contained in the intersection
\[V_\Delta\cap\bigcap_{j\in\vtx_i}\{Z_{j}\neq 0\}.\] In
particular, \(V_\Delta\) is covered by the open sets
\(V_\Delta\cap\{Z_j\neq 0\}\), \(j\in\vtx\).
\end{Corollary}
\begin{Proof}
Let \(b\in\ZZ_{>0}\) be such that \(c_j:=ba_j\in\ZZ_{\geq
0}\). The equation
\(Z_i^b=\prod_{j\in\vtx_i}Z_j^{c_j}\) holds on
\(V_\Delta\) by Lemma \ref{lma:equations}. If \(Z_i\neq 0\)
then this means \(Z_j\neq 0\) for all
\(j\in\vtx_i\).\qedhere

\end{Proof}
\begin{Lemma}[Exercise \ref{exr:automorphisms}]\label{lma:automorphisms}
Let \(A\) be an integer matrix with rows
\(\bm{A}_i\). Consider the morphism
\(\tilde{A}\colon(\CC^*)^n\to(\CC^*)^n\) defined by
\[\tilde{A}(\bm{z}) = (\bm{z}^{\bm{A}_1}, \ldots,
\bm{z}^{\bm{A}_n}).\] If we write
\(\bm{w} :=\tilde{A}(\bm{z})\) then \(\bm{w}^{\bm{q}} =
\bm{z}^{\bm{q}A}\) for any integer row vector \(\bm{q}\). The
morphism \(\tilde{A}\) is invertible if and only if \(A\in
GL(n,\ZZ)\).

\end{Lemma}
\begin{Lemma}\label{lma:transform}
Suppose \(T\colon\RR^n\to\RR^n\) is a map of the form
\(T(\bm{x})=\bm{x}A+\bm{c}\) for some \(A\in GL(n,\ZZ)\) and
\(\bm{c}\in\ZZ^n\). Let
\(\tilde{A}^{-1}\colon(\CC^*)^n\to(\CC^*)^n\) be the morphism
given by the matrix \(A^{-1}\) as in Lemma
\ref{lma:automorphisms}. Then
\(F_{T(\Delta)}\circ\tilde{A}^{-1} = F_{\Delta}\).
\end{Lemma}
\begin{Proof}\belowdisplayskip=-12pt First note that the
constant term \(c\) has no effect on the image of \(F_\Delta\):
it introduces an overall scale factor \(\bm{z}^{\bm{c}}\) into
every homogeneous coordinate. We therefore assume without loss
of generality that \(\bm{c}=0\). Let \(\bm{q}_i=\bm{p}_iA\) be
the vertices of \(T(\Delta)\). We have: \[F_{T(\Delta)}(\bm{w})
= [\bm{w}^{\bm{q}_1}:\cdots:\bm{w}^{\bm{q}_N}],\] so
\begin{align*}
F_{T(\Delta)}\circ\tilde{A}^{-1}(\bm{z})&=[\bm{z}^{\bm{q}_1A^{-1}}:\cdots:\bm{z}^{\bm{q}_NA^{-1}}]\\
&=[\bm{z}^{\bm{p}_1}:\cdots:\bm{z}^{\bm{p}_N}]=F_\Delta(\bm{z}).
\end{align*}\qedhere \end{Proof}

\begin{Remark}
Note that although the variety \(X_\Delta\) is unchanged by \(T\),
the moment map is changed by \(T\) because the projection \(P\)
from Theorem \ref{thm:toric_varieties} changes in such a way that
the moment image is \(T(\Delta)\).

\end{Remark}
\begin{Lemma}\label{lma:equivariant_local_charts}
If \(i\in\vtx\) then \(V_{\Delta}\cap\{z_i\neq 0\}\) is
\(T^n\)-equivariantly biholomorphic to \(\CC^n\) with its
standard torus action.
\end{Lemma}
\begin{Proof}
Because our polytope is Delzant, we can apply a transformation
as in Lemma \ref{lma:transform} so that \(\bm{p}_i\) is at the
origin. By making a further transformation, we can assume that
if \(\bm{p}_{j_1},\ldots,\bm{p}_{j_n}\) are the closest
lattice points to \(\bm{p}_i\) along the \(n\) edges meeting
at \(\bm{p}_i\) then these sit at the points
\((1,0,\ldots,0),\ldots,(0,\ldots,0,1)\). Now any lattice
point \(\bm{p}_k\in\Delta\) can be written as a nonnegative
integer linear combination of these basis vectors, so
\(Z_k=\prod_{s=1}^{n}Z_{j_s}^{a_{j_s}}\) with
\(a_{j_s}\in\ZZ_{\geq 0}\). This means that on \(Z_i\neq 0\)
we can take \(Z_i=1\) and use \(Z_{j_1},\ldots,Z_{j_n}\) as
global coordinates on \(V_\Delta\cap\{Z_i\neq 0\}\). Since
\(\bm{p}_{j_s}\) is the \(s\)th basis vector, the torus action
rotates \(Z_{j_s}\) by \(e^{it_s}\), which shows that the
biholomorphism we have chosen is equivariant with the standard
torus action. \qedhere

\end{Proof}
As a consequence, we see that \(V_\Delta\) is smooth because we
have covered \(V_\Delta\) by smooth coordinate charts. We also
see that \(F_{\Delta}((\CC^*)^n)\subset V_\Delta\) is
Zariski-dense in \(V_\Delta\) because it intersects each chart
\(V_\Delta\cap\{Z_i\neq 0\}\cong\CC^n\) in the Zariski-dense
subset \((\CC^*)^n\). We deduce that \(X_\Delta=V_\Delta\) and
that \(X_\Delta\) is smooth.

Finally, we need to check that the moment image of \(X_\Delta\)
is \(\Delta\). Recall from Theorem \ref{thm:convexity}(1) that
the moment image is the convex hull of the moment images of the
fixed points, so it suffices to show that the fixed points are
sent by the moment map to the vertices of \(\Delta\).

For each \(i\in\vtx\) (i.e.\ \(\bm{p}_i\) is a vertex of
\(\Delta\)), let \(e_i\in\cp{N-1}\) be the point whose
homogeneous coordinates are \(Z_i=1\) and \(Z_j=0\) if \(j\neq
i\). In the \(T^n\)-equivariant local chart \(Z_i\neq 0\) from
Lemma \ref{lma:equivariant_local_charts}, \(e_i\) is sent to the
origin, which is a \(T^n\)-fixed point and the only
\(T^n\)-fixed point in that chart. This shows that the
\(T^n\)-fixed points in \(X_\Delta\) are precisely the points
\(e_i\). We have \(\mu(e_i)\cdot P=\bm{p}_i\), so we deduce that
the \(T^n\)-fixed points map under the moment map to the
vertices of \(\Delta\), as required.\index{toric variety|)}

\section{Solutions to inline exercises}

\begin{Exercise}[Lemma \ref{lma:automorphisms}]\label{exr:automorphisms}
Let \(A\) be an integer matrix with rows
\(\bm{A}_i\). Consider the morphism
\(\tilde{A}\colon(\CC^*)^n\to(\CC^*)^n\) defined by
\[\tilde{A}(\bm{z}) = (\bm{z}^{\bm{A}_1}, \ldots,
\bm{z}^{\bm{A}_n}).\] If we write
\(\bm{w} :=\tilde{A}(\bm{z})\) then show that
\(\bm{w}^{\bm{q}} = \bm{z}^{\bm{q}A}\) for any integer row
vector \(\bm{q}\). Prove that the morphism \(\tilde{A}\) is
invertible if and only if \(A\in GL(n,\ZZ)\).
\end{Exercise}
\begin{Solution}
We have \(w_i = z_1^{A_{i1}}\cdots z_n^{A_{in}}\),
so
\begin{align*}\bm{w}^{\bm{q}} &= (z_1^{A_{11}q_1}\cdots
z_n^{A_{1n}q_1})\cdots (z_1^{A_{n1}q_n}\cdots
z_n^{A_{nn}q_n})\\
&= z_1^{A_{11}q_1+\cdots+A_{n1}q_n}\cdots
z_n^{A_{1n}q_1+\cdots+A_{nn}q_n}\\
&= \bm{z}^{\bm{q}A}.\end{align*}
In particular, this shows that
\(\widetilde{AB}=\tilde{A}\tilde{B}\).

In particular, if \(A\in GL(n,\ZZ)\) then
\(\widetilde{A^{-1}}\) gives an inverse for
\(\tilde{A}\).

More generally, let \(S,T\) be invertible integer matrices
such that \(D:=S A T\) is in Smith normal form. Since
\(\tilde{D}=\tilde{S}\tilde{A}\tilde{T}\), we see that
\(\tilde{D}\) is invertible if and only if \(\tilde{A}\) is
invertible. But \(\tilde{D}(z_1, \ldots, z_n) = (z_1^{D_{11}},
\ldots, z_n^{D_{nn}})\) which is invertible only if
\(D_{ii}=\pm 1\) for all \(i\), which holds only if \(D\in
GL(n,\ZZ)\), which holds only if \(A\in GL(n,\ZZ)\). \qedhere

\end{Solution}
\chapter{Visible contact hypersurfaces and Reeb flows}
\label{ch:visible_contact}
\thispagestyle{cup}

We have focused a lot on visible Lagrangian submanifolds, but
one can also ``see'' other sorts of submanifolds using
Lagrangian torus fibrations. In this section, we discuss visible
submanifolds of codimension 1 and {\em contact geometry}.

\section{Hypersurfaces}

Suppose we have a Lagrangian fibration \(f\colon X\to B\). One
easy way to produce real codimension 1 hypersurfaces in \(X\) is
to take the preimage of a codimension 1 submanifold of \(B\).

\begin{Example}\label{exm:ellipsoids}
Let \(X=\CC^2\), \(B=\RR^2\) and \[f\colon X\to \RR^2, \quad
f(z_1,z_2) = \left(\frac{1}{2}|z_1|^2,
\frac{1}{2}|z_2|^2\right)\] be the moment map for the standard
torus action. Let \(a,b,c>0\) be positive constants and
consider the line \(\ell_{a,b,c}\subset B\) defined by the
equation \(a x+b y=c\) (see Figure \ref{fig:ellipsoid}). The
preimage in \(X\) is
the\index{ellipsoid!contact-type|(}\index{contact!-type
hypersurface} ellipsoid\footnote{Note that {\em ellipsoid} is
also used to mean the compact region bounded by this
hypersurface (sometimes called a {\em solid ellipsoid}).}
\[Y_{a,b,c} := \{a|z_1|^2+b|z_2|^2=2c\}.\]

\begin{figure}[htb]
\begin{center}
\begin{tikzpicture}
\filldraw[lightgray,opacity=0.5,draw=none] (0,2) -- (0,0) -- (2,0) -- (2,2) -- cycle;
\draw (0,2) -- (0,0) -- (2,0);
\draw (0,0.8) -- (1.6,0);
\end{tikzpicture}
\caption{The preimage of the line segment is an ellipsoid.}
\label{fig:ellipsoid}
\end{center}
\end{figure}
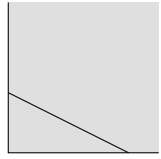

\end{Example}
\begin{Definition}\label{dfn:liouville_vector_field}
A vector field \(Z\) on a symplectic manifold \((X,\omega)\)
is a {\em Liouville}\index{vector field!Liouville|(} or {\em
symplectically dilating} \index{vector field!symplectically
dilating|see {vector field, Liouville}} vector field if
\(\Lie_Z\omega=\omega\). A hypersurface \(Y\subset X\) is said to
be of {\em contact-type}\index{contact!-type hypersurface|(} if
there is a Liouville vector field defined in a neighbourhood
of \(Y\) which is everywhere transverse to \(Y\). If
\(Y=\partial X\) then we say \(Y\) is {\em
convex}\index{boundary!convex} or {\em
concave}\index{boundary!concave} if \(Z\) points respectively
out of or into \(X\).

\end{Definition}
\begin{Example}\label{exm:contact_ellipsoids}
Continuing Example \ref{exm:contact_ellipsoids}, let
\((p_1,p_2,q_1,q_2)\) be action-angle coordinates on \(\CC^2\)
for the standard torus action (so \(p_i=\frac{1}{2}|z_i|^2\)
and \(q_i\) is the argument of \(z_i\)) and let \(Z\) be the
vector field given in action-angle coordinates by
\(p_1\frac{\partial}{\partial p_1}+
p_2\frac{\partial}{\partial p_2}\). This is a Liouville vector
field: \[\Lie_Z\omega= d\iota_Z\omega = d(p_1 dq_1 + p_2
dq_2)=\sum dp_i\wedge dq_i.\] Moreover, \(f_*Z\) is the radial
vector field in (the positive quadrant of) \(\RR^2\) (see
Figure \ref{fig:ellipsoid_liouville_vector_field}). Since
\(a,b,c>0\), this is transverse to \(\ell_{a,b,c}\) and hence
\(Z\) is transverse to \(Y_{a,b,c}\). Thus our
ellipsoids\index{ellipsoid!contact-type|)} are contact-type
hypersurfaces\index{contact!-type hypersurface|)}.

\begin{figure}[htb]
\begin{center}
\begin{tikzpicture}
\filldraw[lightgray,opacity=0.5,draw=none] (0,2) -- (0,0) -- (2,0) -- (2,2) -- cycle;
\draw (0,2) -- (0,0) -- (2,0);
\draw (0,0.8) -- (1.6,0);
\foreach \x in {1,...,8}
\foreach \y in {1,...,5}
{\draw[->] (90-10*\x:0.5+\y/10) -- (90-10*\x:0.5+2*\y/10);}
\end{tikzpicture}
\caption{There is a Liouville vector field transverse to the ellipsoid which projects to the radial vector field in action-coordinates.}
\label{fig:ellipsoid_liouville_vector_field}
\end{center}
\end{figure}
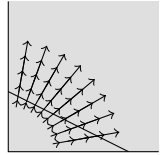

\end{Example}
More generally, this proves:

\begin{Lemma}\label{lma:contact_hypersurface}
Let \((X,\omega)\) be a symplectic manifold and suppose that
\((p,q)\) are local action-angle coordinates on a chart
\(U\subset X\). The vector field \(\sum_i
p_i\frac{\partial}{\partial p_i}\) is a Liouville vector field
on \(U\). If \(H\colon U\to\RR\) is a function which depends
only on \(p\) and \(c\) is a regular value then \(H^{-1}(c)\)
is a contact-type hypersurface\index{contact!-type
hypersurface} if and only if \(\sum_ip_i\frac{\partial
H}{\partial p_i}\) is nowhere vanishing on \(H^{-1}(c)\).
\end{Lemma}
\begin{Proof}
The calculation from Example \ref{exm:contact_ellipsoids}
shows that \(Z\) is Liouville. To understand when
\(H^{-1}(c)\) is contact-type, we see that \(Z\) is transverse
to \(H^{-1}(c)\) if and only if \(dH(Z)\neq 0\) everywhere
along \(H^{-1}(c)\). We compute: \[dH(Z)=\sum_i\frac{\partial
H}{\partial p_i}\,dp_i\left(\sum_jp_j \frac{\partial}{\partial
p_j}\right)=\sum p_i\frac{\partial H}{\partial p_i},\] which
proves the result.\qedhere

\end{Proof}
\section{Contact forms and Reeb flows}

\begin{Definition}
A 1-form \(\alpha\) on a \((2n-1)\)-dimensional manifold \(Y\)
is called a {\em contact form}\index{contact!form|(} if
\(\alpha\wedge (d\alpha)^{n-1}\) is nowhere zero.

\end{Definition}
\begin{Lemma}
Let \(i\colon Y\to X\) be the inclusion map for a contact-type
hypersurface\index{contact!-type hypersurface} with transverse
Liouville field\index{vector field!Liouville|)} \(Z\). The
1-form \(\alpha:=i^*\iota_Z\omega\) is a contact form on
\(Y\).
\end{Lemma}
\begin{Proof}
The \(2n\)-form \(\omega^n\) is a nowhere-vanishing volume
form on \(X\), so \(i^*\iota_Z(\omega^n)\) is a
nowhere-vanishing volume form on \(Y\), since \(Z\) is
transverse to \(Y\). We have \[\iota_Z(\omega^n) =
(\iota_Z\omega)\wedge(\omega^{n-1}) + \cdots +
(\omega^{n-1})\wedge (\iota_Z\omega) =
n(\iota_Z\omega)\wedge(\omega^{n-1}).\] Therefore
\[\frac{i^*\iota_Z\omega^n}{n} =
i^*\left(\iota_Z\omega\wedge\omega^{n-1}\right).\] Set
\(\alpha = i^*\iota_Z\omega\). We have
\(d\alpha=i^*d\iota_Z\omega=i^*\Lie_Z\omega=i^*\omega\), so
\(\alpha\wedge d\alpha^{n-1} =
i^*\left(\iota_Z\omega\wedge\omega^{n-1}\right)\). As this is
a nowhere-vanishing volume form on \(Y\), this proves the
result. \qedhere

\end{Proof}
\begin{Definition}
If \(\alpha\) is a contact form\index{contact!form|)} on
\(Y\), we define the {\em Reeb vector field}\index{vector
field!Reeb|(} \(R_\alpha\) to be the unique vector field
satisfying \[\iota_{R_\alpha}d\alpha=0,\quad
\iota_{R_\alpha}\alpha=1.\]

\end{Definition}
\begin{Remark}
The equation \(\iota_{R_\alpha}d\alpha=0\) says that the Reeb
field points along the line field\footnote{This is called the
{\em characteristic line field}. Recall that \(^\omega\)
denotes the symplectic orthogonal complement; see Definition
\ref{dfn:symplectic_orthogonal_complement}.}
\((TY)^\omega\subset TY\). The equation
\(\iota_{R_\alpha}\alpha=1\) is simply a normalisation
condition which picks out a specific vector \((TY)^\omega\).

\end{Remark}
\begin{Example}[Exercise \ref{exr:reeb_flow}]\label{exm:reeb_flow}
Let \(X\) be a symplectic manifold and \(U\subset X\) be the
domain of an action-angle chart with action-angle coordinates
\((p,q)\). Let \(H\colon U\to\RR\) be a function depending
only on \(p\) and satisfying the condition
\(\sum_ip_i\frac{\partial H}{\partial p_i}\neq 0\) along a
regular level set \(H^{-1}(c)\). Let \(\alpha=i^*(\sum_k p_k\,
dq_k)\) be the contact form\index{contact!form} guaranteed by
Lemma \ref{lma:contact_hypersurface}. The Reeb vector field is
given by \[R_\alpha = \left(\sum_i\frac{\partial H}{\partial
p_i}\frac{\partial}{\partial q_i}\right) \bigg/
\left(\sum_jp_j\frac{\partial H}{\partial p_j}\right).\]

\end{Example}
\begin{Example}[Exercise \ref{exr:reeb_ellipsoids}]\label{exm:reeb_ellipsoids}
Continuing Example \ref{exm:contact_ellipsoids}, we take \(H =
a p_1 + b p_2\), which gives the contact-type
ellipsoid\index{ellipsoid!contact-type}
\index{contact!-type hypersurface}
\(Y_{a,b,c}\subset\CC^2\). By Example
\ref{exm:reeb_flow}, the Reeb field is
\[c^{-1}\left(a\frac{\partial}{\partial
q_1}+b\frac{\partial}{\partial q_2}\right).\] The dynamics
of the flow along the Reeb field now depend on the constants
\(a\) and \(b\).

If the ratio \(b/a\) is irrational then the orbits of
\(R_\alpha\) are lines of irrational slope in the
\((q_1,q_2)\)-torus. The exception is when \(p_1=0\) or
\(p_2=0\): here the action-angle coordinates are degenerate in
the sense that the fibre is a circle parametrised by \(q_2\)
respectively \(q_1\). These two circles are closed orbits of
the Reeb field (see Figure \ref{fig:ellipsoid_reeb_orbits}).

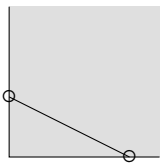
\begin{figure}[htb]
\begin{center}
\begin{tikzpicture}
\filldraw[lightgray,opacity=0.5,draw=none] (0,2) -- (0,0) -- (2,0) -- (2,2) -- cycle;
\draw (0,2) -- (0,0) -- (2,0);
\draw (0,0.8) -- (1.6,0);
\node at (0,0.8) {\(\circ\)};
\node at (1.6,0) {\(\circ\)};
\end{tikzpicture}
\caption{If the slope is irrational then the only closed Reeb orbits are the circles living over the points marked \(\circ\).}
\label{fig:ellipsoid_reeb_orbits}
\end{center}
\end{figure}

If the ratio \(b/a\) is rational then all the Reeb orbits are
closed. More precisely, if \(a=\rho m\) and \(b = \rho n\) for
coprime integers \(m,n\) then the orbits away from \(p_1=0\)
and \(p_2=0\) have period \(2\pi c/\rho\). There are still two
exceptional orbits at \(p_1=0\) and \(p_2=0\), with periods
\(2\pi c/n\rho\) and \(2\pi c/m\rho\) respectively. In this
case, we can take the symplectic quotient of \(\CC^2\) with
respect to this Hamiltonian and obtain a symplectic sphere
with two orbifold\index{orbifold} points (compare with Example
\ref{exm:wps}, where we obtained weighted projective
spaces\index{projective space!weighted} in this way). The
quotient map \(Y_{a,b,c}\to Y_{a,b,c}/S^1\) is an example of a
{\em Seifert fibration}\index{Seifert fibration}; for more
about the topology of Seifert fibred manifolds, see
\cite{JankinsNeumann,OrlikSeifert} or Seifert's appendix to
\cite{SeifertThrelfall}.

\end{Example}
\section{Solutions to inline exercises}

\begin{Exercise}[Example \ref{exm:reeb_flow}]\label{exr:reeb_flow}
Show that the Reeb vector field in Example \ref{exm:reeb_flow}
is given by \[R_\alpha = \left(\sum_i\frac{\partial
H}{\partial p_i}\frac{\partial}{\partial q_i}\right) \bigg/
\left(\sum_jp_j\frac{\partial H}{\partial p_j}\right).\]
\end{Exercise}
\begin{Solution}
We recall that our contact manifold is the level set
\(Y:=H^{-1}(c)\) of a function \(H(\bm{p},\bm{q})\) which is
independent of \(\bm{q}\). Since \(R_\alpha\) has no
\(\partial_{p_j}\) components, we have \(dH(R_\alpha)=0\), so
\(R_\alpha\) is tangent to \(Y\). Next, we have
\(d\alpha=i^*\sum dp_k\wedge dq_k\), where \(i\) is the
inclusion of the level set \(Y\). Therefore
\[\iota_{R_\alpha}d\alpha = \left(\sum_k\frac{\partial
H}{\partial p_k}dp_k\right)\bigg/\left(\sum_j p_j\frac{\partial
H}{\partial p_j}\right) \propto dH,\] and \(dH\) vanishes on
\(Y\). Finally, observe
that \begin{align*}\iota_{R_\alpha}\alpha &= \left(\sum_k
p_k\,dq_k\right)\left(\left(\sum_i\frac{\partial H}{\partial
p_i}\frac{\partial}{\partial q_i}\right) \bigg/
\left(\sum_jp_j\frac{\partial H}{\partial p_j}\right)\right)
\\ &= \left(\sum_ip_i \frac{\partial H}{\partial p_i}\right)
\bigg/ \left(\sum_jp_j\frac{\partial H}{\partial
p_j}\right)=1\end{align*}

\end{Solution}
\begin{Exercise}[Example \ref{exm:reeb_ellipsoids}]\label{exr:reeb_ellipsoids}
Show that in Example \ref{exm:reeb_ellipsoids}, the Reeb field
is \[c^{-1}\left(a\frac{\partial}{\partial
q_1}+b\frac{\partial}{\partial q_2}\right).\]
\end{Exercise}
\begin{Solution}
We have \(\partial H/\partial p_1=a\) and \(\partial
H/\partial p_2=b\), so the Reeb field\index{vector
field!Reeb|)} from Example \ref{exm:reeb_flow} becomes
\[\left(a\frac{\partial}{\partial
q_1}+b\frac{\partial}{\partial
q_2}\right)\bigg/\left(a p_1+b p_2\right).\] Since
\(H(p,q)=a p_1+b p_2=c\) along
\(Y_{a,b,c}\), this is
\(c^{-1}\left(a\frac{\partial}{\partial
q_1}+b\frac{\partial}{\partial q_2}\right)\) as
required. \qedhere

\end{Solution}
\chapter{Tropical Lagrangian submanifolds}
\label{ch:tropical_lag}
\thispagestyle{cup}

While visible Lagrangians\index{Lagrangian!submanifold!visible}
are associated with straight line segments or affine subspaces
in the base of a Lagrangian torus fibration, {\em tropical
Lagrangians}\index{Lagrangian!submanifold!tropical|(} are
associated with certain {\em piecewise} linear subsets (like
trivalent graphs). More precisely, a tropical
Lagrangian\index{tropical Lagrangian|see {Lagrangian,
submanifold, tropical}} \index{visible Lagrangian|see
{Lagrangian, submanifold, visible}} is an immersed Lagrangian
whose image under the fibration is a small thickening of a {\em
tropical subvariety} in the base of the fibration. We will focus
on the case of tropical curves.

Mikhalkin \cite{Mikhalkin} and Matessi \cite{Matessi1, Matessi2}
have given constructions of tropical Lagrangians associated with
tropical curves and tropical hypersurfaces respectively. We will
focus on the 4-dimensional case where these constructions
coincide, and we will use Mikhalkin's conventions.

\section{A Lagrangian pair-of-pants}

We consider \(\CC^*\times\CC^*\) equipped with complex
coordinates \((z_1,z_2)\). We identify this with \(\RR^2\times
T^2\) with coordinates \((p_1,p_2)\in\RR^2\) and \((q_1,q_2)\in
T^2=(\RR/2\pi\ZZ)^2\) using the identification \[z_k =
\exp(p_k+iq_k).\]

\begin{Theorem}[Mikhalkin]\label{thm:mikhalkin}
Let \(R_1\), \(R_2\), \(R_3\) be three rays with rational
slope in the \(\bm{p}\)-plane emanating from the origin, and
let \(v_1,v_2,v_3\) be the primitive integer vectors pointing
along these rays. Suppose that the {\em balancing
condition} \begin{equation}\label{eq:balancing} v_1+ v_2+ v_3=
0\end{equation} holds and that any two of these vectors form a
\(\ZZ\)-basis for the integer lattice\footnote{These three
vectors form what Conway \cite{Conway} calls a {\em
superbase}.}. Let \(L_1,L_2,L_3\) be the visible
Lagrangian
half-cylinders living over \(R_1,R_2,R_3\) and fix
\(\epsilon>0\). Let
\(U:=\{\bm{p}\,:\,|\bm{p}|>\epsilon\}\). There is an embedded
Lagrangian submanifold\index{Lagrangian!pair-of-pants|(}
\(L\subset \RR^2\times T^2\), diffeomorphic to the
pair-of-pants, such that \(U\cap L = U\cap(L_1\cup L_2\cup
L_3)\).
\end{Theorem}
\begin{Proof}
We will focus on the case \(v_1=(-1,0)\), \(v_2=(0,-1)\),
\(v_3=(1,1)\). This implies the general case: if \(v_1,v_2\)
is obtained from this basis by an element of \(GL(2,\ZZ)\)
then Lemma \ref{lma:globcoords} gives us a fibred
symplectomorphism living over this integral affine
transformation of the \(\bm{p}\)-plane, and we can apply this
symplectomorphism to the Lagrangian obtained for
\(v_1=(-1,0)\), \(v_2=(0,-1)\).

The starting point of the construction is the following exercise.

\begin{Lemma}[Exercise \ref{exr:hyperkahler_twist}]\label{lma:hyperkahler_twist}
Consider the {\em hyperK\"{a}hler twist}\index{hyperKaehler
twist@hyperK\"{a}hler twist} \[\RR^2\times T^2\to\RR^2\times
T^2,\quad (p_1,p_2,q_1,q_2)\mapsto (p_1,p_2,-q_2,q_1)\] If
\(C\subset\RR^2\times T^2\) is a complex curve with respect
to the complex coordinates \(z_k=e^{p_k+iq_k}\) then the
image of \(C\) under the hyperK\"{a}hler twist is Lagrangian
for the symplectic form \(\sum dp_i\wedge dq_i\).

\end{Lemma}
We will apply this lemma to the complex curve \(C = \{z_2 = 1+
z_1\} \subset\CC^*\times \CC^*\). This is diffeomorphic to a
pair-of-pants (3-punctured sphere); we can parameterise it as
\(z\mapsto (z,1+z)\) with \(z\in\CC\setminus\{0,-1\}\). Let
\(L\subset \RR^2\times T^2\) be the hyperK\"{a}hler twist of
\(C\). In coordinates, this is given parametrically by
\begin{align*}
p_1(z)&=\ln|z|&p_2(z)&=\ln|1+z|\\
q_1(z)&=-\arg(1+z)&q_2(z)&=\arg(z).
\end{align*}
We now write \(z=re^{i\theta}\) and see what happens as \(r\to
0\). We have \(\lim_{r\to 0}p_2(z) = \ln(1)=0\) and \(\lim_{r\to
0}q_1(z)=-\arg(1)=0\mod 2\pi\). This means that near the
puncture \(0\), our Lagrangian \(L\) is asymptotic to the
cylinder \[L_1 :=\{(p_1,0,0,q_2)\,:\,p_1 < 0,q_2\in[0,2\pi]\},\]
which is a visible Lagrangian cylinder associated to the
negative \(p_1\)-axis (note that \(\lim_{r\to
0}p_1(z)=-\infty\)). We write \(R_1\) for the negative
\(p_1\)-axis. A similar analysis near the punctures
\(z\to -1\) and \(z\to \infty\) shows that \(L\) has asymptotes
along the visible Lagrangian cylinders
\begin{align*}
L_2 &:= \{(0,p_2,q_1,0)\,:\,p_2<0,q_1\in[0,2\pi]\},\\
L_3 &:=\{(p,p,q,-q)\,:\,p>0,q\in[0,2\pi]\}
\end{align*}
associated to the rays \(R_2\) and \(R_3\) shown in Figure
\ref{fig:tropical_rays}.

\begin{figure}[htb]
\begin{center}
\begin{tikzpicture}
\filldraw[gray] (180:2) to[out=0,in=-135] (45:2) to[out=-135,in=90] (-90:2) to[out=90,in=0] (180:2);
\begin{scope}[shift={(7,0)}]
\draw[gray] (45:1.5) -- (45:2) node [right] {\(R_3\)};
\draw[gray] (-1.5,0) -- (-2,0) node [left] {\(R_1\)};
\draw[gray] (0,-1.5) -- (0,-2) node [below] {\(R_2\)};
\draw (0,0) circle [radius = 1];
\node at (-45:0.3) [right] {\small \(K_1\)};
\draw (0,0) circle [radius = 1.5];
\node at (20:1.23) {\small \(K_2\)};
\filldraw[gray,opacity=0.5] (180:1.5) -- (175.5:1) -- (197:1) -- cycle;
\filldraw[gray,opacity=0.5] (-90:1.5) -- (-85.5:1) -- (-107:1) -- cycle;
\filldraw[gray,opacity=0.5] (49.5:1) -- (45:1.5) -- (40.5:1) -- cycle;
\begin{scope}
\clip (0,0) circle [radius = 1];
\filldraw[gray] (180:2) to[out=0,in=-135] (45:2) to[out=-135,in=90] (-90:2) to[out=90,in=0] (180:2);
\end{scope}
\end{scope}
\end{tikzpicture}
\caption{(a) The projection of \(L\) to the \(\bm{p}\)-plane. (b) \(L\) is asymptotic to the visible Lagrangians cylinders \(L_1\), \(L_2\), \(L_3\) which live over the three rays shown. We modify it outside the region \(K_1\) so that it coincides with these visible Lagrangians outside \(K_2\).}
\label{fig:tropical_rays}
\end{center}
\end{figure}
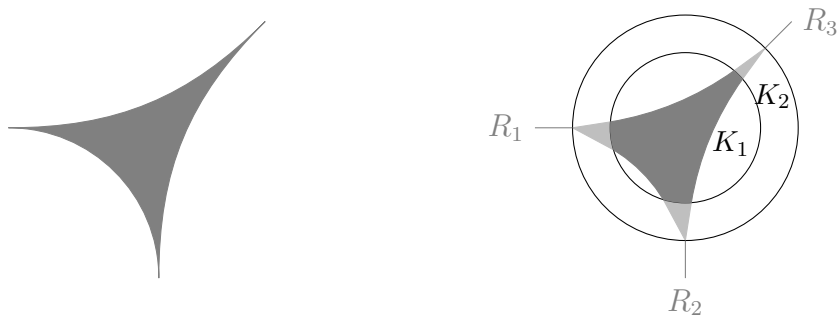

Given \(0<k_1<k_2\), define the compact regions
\(K_i=\{p_1^2+p_2^2\leq k_i\}\), \(i=1,2\), and let \(U_i =
(\CC^*\times\CC^*)\setminus K_i\) be their complements. We
will modify \(L\) to obtain a Lagrangian pair-of-pants \(L'\)
with \[K_1 \cap L' = K_1\cap L\mbox{ and }U_2 \cap L' =
U_2\cap (L_1\cup L_2\cup L_3).\]
We will explain this modification for the puncture asymptotic to
\(L_1\); the other cases are similar. We can identify a
neighbourhood of \(L_1\) with a neighbourhood of the
zero-section in \(T^*L_1\). More precisely, we think of
\(p_1=\ln r\) and \(q_2=\theta\) as coordinates on the cylinder
and \(-q_1\) and \(p_2\) as dual momenta\footnote{The symplectic
form is \(-dq_1\wedge dp_1+dp_2\wedge dq_2\). Because \(q_1\) is
circle-valued, the identification of \(q_1\) with a coordinate
on the fibre of \(T^*L_1\) only makes sense if \(q_1\approx
0\).}.

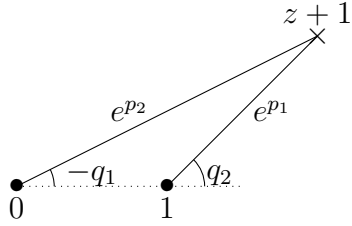
\begin{figure}[htb]
\begin{center}
\begin{tikzpicture}
\node (A) at (0,0) {\(\bullet\)};
\node at (A) [below] {\(0\)};
\node (B) at (2,0) {\(\bullet\)};
\node at (B) [below] {\(1\)};
\node (C) at (4,2) {\(\times\)};
\node at (C) [above] {\(z+1\)};
\draw (A.center) -- (C.center);
\draw (B.center) -- (C.center);
\draw[dotted] (0,0) -- (3,0);
\draw (2.5,0) arc [radius = 0.5,start angle = 0, end angle=45];
\node at (2.7,0.15) {\(q_2\)};
\node at (3,1) [right] {\(e^{p_1}\)};
\draw (0.5,0) arc [radius = 0.5,start angle = 0, end angle=27];
\node at (1,0.17) {\(-q_1\)};
\node at (1.5,1) {\(e^{p_2}\)};
\end{tikzpicture}
\caption{Geometric picture behind formulae for \(q_1\) and \(p_2\).}
\label{fig:linkage}
\end{center}
\end{figure}

From Figure \ref{fig:linkage}, you can extract equations for
our section: \[q_1 = -\arctan\left(\frac{e^{p_1}\sin
q_2}{1+e^{p_1}\cos q_2}\right)\qquad p_2 =
\frac{1}{2}\ln(1+2e^{p_1}\cos q_2 + e^{2p_1}),\] defined on
the subset \(p_1 < -\ln 2\). In other words, it is the graph
of the 1-form \[\beta := \frac{1}{2}\ln(1+2e^{p_1}\cos q_2 +
e^{2p_1})dp_1+ \arctan\left(\frac{e^{p_1}\sin
q_2}{1+e^{p_1}\cos q_2}\right)dq_2.\] This 1-form is closed
(this is equivalent to the Lagrangian condition, but you can
also check it directly by differentiating), but it is also
exact: the obstruction to exactness\footnote{on the
cylindrical end \(p_1 < -\ln 2\) whose de Rham cohomology has
rank 1.} is the integral \(\int\beta\) around the loop
\(p_1=0\), \(q_2\in [-\pi,\pi]\), i.e.\footnote{This integral
vanishes because the integrand is an odd function.}
\[\int_{-\pi}^{\pi}\arctan\left(\frac{\sin q_2}{1+\cos
q_2}\right)dq_2 = 0.\] Exactness means there is a function
\(\varphi(p_1,q_2)\) such that \(\beta = d\varphi\). Pick
\(\epsilon>0\) and let \(\rho(p_1)\) be a cut-off function
equal to \(0\) for \(p_1\leq -\ln 2-\epsilon\) and equal to
\(1\) for \(-\ln 2 + \epsilon \leq p_1\). If we take
\(k_1=(\ln 2-\epsilon)^2\) and \(k_2=(\ln 2+\epsilon)^2\) and
define \(K_1,K_2\) as above, then the Lagrangian cylinder
given by the graph of \(d(\rho\varphi)\) coincides with the
cylinder \(L_1\) outside \(K_2\) and coincides with \(L\) in
the compact region \(K_1\).

We perform a similar modification near each of the three
punctures. The result is a Lagrangian
pair-of-pants\index{Lagrangian!pair-of-pants|)} which
coincides with the three Lagrangian cylinders \(L_1\),
\(L_2\), and \(L_3\) outside the compact set \(K_2\). This
does not quite prove Theorem \ref{thm:mikhalkin}, because we
cannot take \(K_2\) arbitrarily small in this
construction. However, notice that the radial vector field in
the \(\bm{p}\)-plane is a Liouville vector field, by Lemma
\ref{lma:contact_hypersurface}. Our visible Lagrangian
cylinders \(L_1,L_2,L_3\) are preserved by the flow of this
Liouville field, and if we flow \(L\) backwards along this
Liouville field, we ensure that it agrees with \(L_1,L_2,L_3\)
on a larger and larger region (Figure
\ref{fig:tropical_liouville_flow}). This completes the proof
of Theorem \ref{thm:mikhalkin}.\qedhere

\end{Proof}
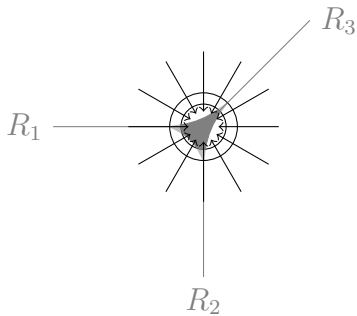
\begin{figure}[htb]
\begin{center}
\begin{tikzpicture}
\draw[gray] (45:0.3) -- (45:2) node [right] {\(R_3\)};
\draw[gray] (-0.3,0) -- (-2,0) node [left] {\(R_1\)};
\draw[gray] (0,-0.3) -- (0,-2) node [below] {\(R_2\)};
\draw (0,0) circle [radius = 0.3];
\draw (0,0) circle [radius = 0.45];
\filldraw[gray,opacity=0.5] (180:0.45) -- (175.5:0.3) -- (197:0.3) -- cycle;
\filldraw[gray,opacity=0.5] (-90:0.45) -- (-85.5:0.3) -- (-107:0.3) -- cycle;
\filldraw[gray,opacity=0.5] (49.5:0.3) -- (45:0.45) -- (40.5:0.3) -- cycle;
\begin{scope}[scale=0.3]
\clip (0,0) circle [radius = 1];
\filldraw[gray] (180:2) to[out=0,in=-135] (45:2) to[out=-135,in=90] (-90:2) to[out=90,in=0] (180:2);
\end{scope}
\foreach \x in {1,...,12}
{\draw[->] ({30*\x}:1) -- ({30*\x}:0.2);}
\end{tikzpicture}
\caption{We can find a Lagrangian whose projection to the \(\bm{p}\)-plane is arbitrarily close to \(R_1\cup R_2\cup R_3\) by flowing backwards along a Liouville field.}
\label{fig:tropical_liouville_flow}
\end{center}
\end{figure}

\section{Immersed Lagrangians}

Now suppose that in the statement of Theorem
\ref{thm:mikhalkin}, the primitive vectors \(v_1,v_2,v_3\) still
satisfy the balancing condition \(v_1+v_2+v_3\), but that each
pair fails to form a \(\ZZ\)-basis for \(\ZZ^2\). If we write
the matrix \(M\) whose rows are \(v_1\) and \(v_2\) then
applying this matrix (on the right) gives an integer matrix
sending \((1,0)\) and \((0,1)\) to \(v_1\) and \(v_2\)
respectively. This map is not induced by a symplectic map
\(\RR^2\times T^2\to\RR^2\times T^2\), but it is induced by a
holomorphic covering map \(h\colon \CC^*\times\CC^*\to
\CC^*\times\CC^*\) of degree \(\det(M)\), namely
\(h(\exp(\bm{p}+i\bm{q})) = \exp(\bm{p}M+i\bm{q}M)\). The
hyperK\"{a}hler twist of \(h(C)\) is another Lagrangian
submanifold with three punctures asymptotic to Lagrangian
cylinders living over the rays pointing in the \(v_1\)-,
\(v_2\)-, and \(v_3\)-directions. As in Theorem
\ref{thm:mikhalkin}, we can modify this Lagrangian so that it
actually coincides with these Lagrangian cylinders outside a
compact set.

The main difference is that the resulting Lagrangian
pair-of-pants\index{Lagrangian!immersed} is not embedded: it is
only immersed.

\begin{Example}\label{exm:mikhalkin_delta}
Consider the case \(v_1 = (2,-1)\), \(v_2 = (-1,2)\). We have
\[h(z_1,z_2) = (z_1^2/z_2, z_2^2/z_1).\] A
self-intersection\index{self-intersection!of immersed
Lagrangian|(} of \(h(C)\) corresponds to a pair of points
\(\xi,\xi'\in\CC\setminus\{0,-1\}\) with
\(h(\xi,\xi+1)=h(\xi',\xi'+1)\). In this case, we can show
there is precisely one self-intersection. Suppose that \(\xi\)
and \(\xi'\) are distinct solutions to \(h(z,z+1) = (u,v)\)
for some \(u,v\in\CC^*\). Then the quadratic equations
\[z^2=u(z+1),\quad (z+1)^2 = vz\] have \(\xi\) and \(\xi'\) as
roots. But then \(\xi+\xi' = -u = v-2\) and \(\xi\xi' = -u =
1\). In particular, \(u=-1\) and \(v = 3\) and \(\xi,\xi'\)
are roots of \(z^2+z+1\), so \(\xi\) and \(\xi'\) are
\(\frac{-1\pm i\sqrt{3}}{2}\). Thus there is precisely one
self-intersection at \((-1,3)\).

\end{Example}
Write \(|v\wedge w|\) for the absolute value of the determinant
of the 2-by-2 matrix whose rows are \(v\) and \(w\). In Example
\ref{exm:mikhalkin_delta}, we have \(|v_1\wedge v_2|= 3\).

\begin{Theorem}\label{thm:mikhalkin_delta}
Let \(\Delta\) be the absolute value of the determinant of the
matrix whose rows are \(v_1\) and \(v_2\). The number of
self-intersections of the Lagrangian
pair-of-pants\index{Lagrangian!pair-of-pants} is \(\delta =
\frac{\Delta-1}{2}\).

\end{Theorem}
\begin{Remark}[Exercise \ref{exr:indep_order}]\label{rmk:indep_order}
If \(v_1,v_2,v_3\) are primitive integer vectors with
\(v_1+v_2+v_3=0\) then \(|v_k\wedge v_{\ell}|\) is an odd
number and is independent of \(k,\ell\).

\end{Remark}
We will not prove Theorem \ref{thm:mikhalkin_delta}, and refer
the interested reader to {\cite[Corollary 4.3]{Mikhalkin}}. We
leave the following related lemma as an exercise:

\begin{Lemma}[Exercise \ref{exr:tropical_double_points}]\label{lma:tropical_double_points}
Suppose we have several straight lines of rational slope in
\(\RR^2\) incident on a point \(b\in B\). Let \(v_1,\ldots,k\)
be primitive integer vectors pointing along these lines. Show
that the visible Lagrangian cylinders above these lines have a
total of \(\delta(b)\) transverse
intersections\index{self-intersection!of immersed
Lagrangian|)}, where \[\delta(b)=\sum_{i< j}|v_i\wedge v_j|.\]

\end{Lemma}
\section{Lagrangians from tropical curves}

We have already seen how to construct Lagrangian submanifolds
living over straight lines of rational slope. Provided we are
willing to allow pinwheel core and Schoen-Wolfson singularities,
these straight lines are allowed to terminate on the toric
boundary. Thanks to Lemma \ref{lma:visibleff}, we also have
visible Lagrangian discs terminating on base-nodes of almost
toric fibrations, providing they live over eigenlines. Now,
courtesy of Theorem \ref{thm:mikhalkin} we have immersed
Lagrangian pairs-of-pants living over trivalent vertices
satisfying the balancing condition \eqref{eq:balancing}. Since
this pair-of-pants coincides with the three visible Lagrangian
cylinders over the edges of the graph (except in a small
neighbourhood of the vertex), we can combine all of these
constructions to get a Lagrangian submanifold (possibly singular
and immersed) living over a trivalent graph whose edges have
rational slope and whose vertices satisfy the balancing
condition \eqref{eq:balancing}. These graphs are called {\em
tropical curves}\index{tropical curve} and the associated
Lagrangians are called {\em tropical Lagrangians}. In summary, a
tropical Lagrangian is made up of:
\begin{itemize}
\item an immersed pair-of-pants over every trivalent vertex \(b\) of
the graph (with \(\delta(b)\) self-intersections),
\item a visible Lagrangian cylinder over every
edge,
\item a \((p,q)\)-pinwheel core\index{Lagrangian!pinwheel core} over
every point where an edge terminates on an edge of the almost
toric base diagram,
\item a disc or Schoen-Wolfson cone\index{Schoen-Wolfson cone} over
every point where an edge terminates at a vertex of the almost
toric base diagram,
\item a disc over every base-node at which an edge terminates,
providing the edge points in the
eigendirection\index{eigenline!visible Lagrangian over} for
the affine monodromy of the base-node.
\end{itemize}
Rather than writing this carefully and formally, it is easier to
explain in some examples.

\begin{figure}[htb]
\begin{center}
\begin{tikzpicture}
\node at (-1,-1) [left] {(a)};
\filldraw[fill=lightgray,opacity=0.5,draw=black] (-1,-1) -- (-1,2) -- (2,-1) -- cycle;
\draw[thick] (-1,-1) -- (-1,2) -- (2,-1) -- cycle;
\draw[very thick] (-1,-1) -- (0,0);
\draw[very thick] (-1,2) -- (0,0);
\draw[very thick] (2,-1) -- (0,0);
\node at (-1,-1) {\(\bullet\)};
\node at (-1,0) {\(\bullet\)};
\node at (-1,1) {\(\bullet\)};
\node at (-1,2) {\(\bullet\)};
\node at (0,-1) {\(\bullet\)};
\node at (0,0) {\(\bullet\)};
\node at (0,1) {\(\bullet\)};
\node at (1,-1) {\(\bullet\)};
\node at (1,0) {\(\bullet\)};
\node at (2,-1) {\(\bullet\)};
\begin{scope}[shift={(4.5,0)}]
\node at (-1,-1) [left] {(b)};
\filldraw[fill=lightgray,opacity=0.5,draw=black] (-1,-1) -- (-1,4) -- (4,-1) -- cycle;
\draw[thick] (-1,-1) -- (-1,4) -- (4,-1) -- cycle;
\draw[very thick] (-1,-1) -- (0,0);
\draw[very thick] (-1,4) -- (0,2);
\draw[very thick] (4,-1) -- (2,0);
\draw[very thick] (0,0) -- (2,0) -- (0,2) -- cycle;
\foreach \x in {-1,...,4}
\foreach \y in {-1,...,\x}
{\node at (3-\x , \y ) {\(\bullet\)};}
\end{scope}
\end{tikzpicture}
\caption{(a) A tropical Lagrangian immersed sphere in \(\cp{2}\). (b) An embedded Lagrangian torus.}
\label{fig:trop_1}
\end{center}
\end{figure}
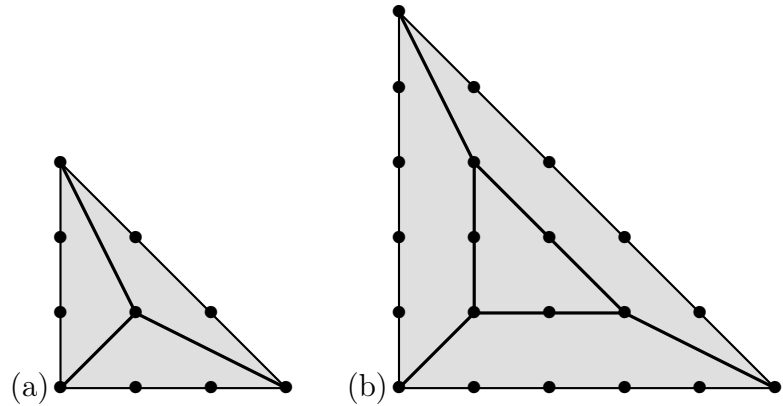

\begin{Example}\label{exm:trop_1a}
Consider the tropical curve shown in Figure
\ref{fig:trop_1}(a). The three corners are all \(\ZZ\)-affine
equivalent; for example, if we think of the trivalent point as
the origin then the matrix \(\lmatrix -1 & 1 \\ -1 & 0
\rmatrix \) sends the top vertex to the bottom left and
preserves the tropical curve. The bottom left corner is the
local model from Example \ref{exm:hittingvertex_1}, so the
three corners give us visible Lagrangian discs with which to
cap off the Lagrangian pair-of-pants coming from the 3-valent
point. The outgoing vectors at the 3-valent point are:
\[v_1=(-1,-1),\quad v_2 = (2,-1),\quad v_3 = (-1,2)\] so
\(|v_1\wedge v_2|=\left|\det\begin{pmatrix} -1 & -1 \\ 2 &
-1\end{pmatrix} \right|=3\), so the \(\delta\)-invariant is
\((3-1)/2=1\). This tropical Lagrangian is therefore an
immersed sphere\index{Lagrangian!immersed} in \(\cp{2}\) with
one transverse double point.

\end{Example}
\begin{Example}\label{exm:trop_1b}
Consider the tropical curve shown in Figure
\ref{fig:trop_1}(b). The corners are all the same as in
Example \ref{exm:trop_1a}. The trivalent vertices all have
\(\delta=1\), so the result is an embedded Lagrangian. By
inspection, it is topologically a
torus\index{Lagrangian!torus}.

\end{Example}
\begin{Example}[Exercise \ref{exr:wtf_is_it}]\label{exm:wtf_is_it}
Nemirovski and Shevchishin proved, independently and in very
different ways, that there is no embedded Lagrangian Klein
bottle\index{Lagrangian!Klein bottle} in \(\cp{2}\). Why is
Figure \ref{fig:wtf_is_it} not a counterexample to their
theorem?

\begin{figure}[htb]
\begin{center}
\begin{tikzpicture}
\filldraw[fill=lightgray,opacity=0.5,draw=black] (-2,-2) -- (-2,4) -- (4,-2) -- cycle;
\draw[thick] (-2,-2) -- (-2,4) -- (4,-2) -- cycle;
\draw[very thick] (-2,-2) -- (0,0);
\draw[very thick] (-2,3) -- (0,0);
\draw[very thick] (3,-2) -- (0,0);
\node at (-2,-2) {\(\bullet\)};
\node at (-2,-1) {\(\bullet\)};
\node at (-2,0) {\(\bullet\)};
\node at (-2,1) {\(\bullet\)};
\node at (-2,2) {\(\bullet\)};
\node at (-2,3) {\(\bullet\)};
\node at (-2,4) {\(\bullet\)};
\node at (-1,-2) {\(\bullet\)};
\node at (-1,-1) {\(\bullet\)};
\node at (-1,0) {\(\bullet\)};
\node at (-1,1) {\(\bullet\)};
\node at (-1,2) {\(\bullet\)};
\node at (-1,3) {\(\bullet\)};
\node at (0,-2) {\(\bullet\)};
\node at (0,-1) {\(\bullet\)};
\node at (0,0) {\(\bullet\)};
\node at (0,1) {\(\bullet\)};
\node at (0,2) {\(\bullet\)};
\node at (1,-2) {\(\bullet\)};
\node at (1,-1) {\(\bullet\)};
\node at (1,0) {\(\bullet\)};
\node at (1,1) {\(\bullet\)};
\node at (2,-2) {\(\bullet\)};
\node at (2,-1) {\(\bullet\)};
\node at (2,0) {\(\bullet\)};
\node at (3,-2) {\(\bullet\)};
\node at (3,-1) {\(\bullet\)};
\node at (4,-2) {\(\bullet\)};
\end{tikzpicture}
\end{center}
\caption{Another tropical Lagrangian in \(\cp{2}\).}
\label{fig:wtf_is_it}
\end{figure}

\end{Example}
\begin{Example}\label{exm:shevsmi}
The almost toric base diagrams\index{almost toric base
diagram!symplectic triangle inequality|(} in Figure
\ref{fig:shevsmi} represent blow-ups of the standard
symplectic ball in three smaller balls, where the symplectic
areas\index{symplectic area} of the exceptional spheres after
blowing up are \(a\), \(b\) and \(c\) (you can get these
diagrams by performing one toric blow-up at the origin and two
non-toric blow-ups and changing branch cuts). In Figure
\ref{fig:shevsmi}(a), the tropical Lagrangian associated with
the tropical curve is a Lagrangian \(\rp{2}\) representing the
\(\ZZ/2\)-homology class \(E_1+E_2+E_3\), meeting the toric
boundary along a \((2,1)\)-pinwheel core (M\"{o}bius
strip). In Figure \ref{fig:shevsmi}(b), the tropical
Lagrangian is a disc with boundary on the left-hand toric
boundary. Shevchishin and Smirnov \cite{ShevSmi} showed that
the \(\ZZ/2\)-homology class \(E_1+E_2+E_3\) can be
represented by a
Lagrangian\index{Lagrangian!RP2@$\mathbb{RP}^2$} \(\rp{2}\) if
and only if the {\em symplectic triangle
inequalities}\index{symplectic triangle inequalities} hold:
\[a<b+c,\quad b<c+a,\quad c<a+b.\] Whenever these inequalities
hold (e.g.\ Figure \ref{fig:shevsmi}(a)), the tropical curve
modelled on this tripod furnishes us with a Lagrangian
\(\rp{2}\), and whenever they fail, it furnishes us with a
disc. For more tropical discussions along these lines, see
\cite{EvansKB}.

\begin{figure}[htb]
\begin{center}
\begin{tikzpicture}
\filldraw[fill=gray!25,draw=none] (0,4/3) -- (0,4) -- (4,0) -- (4/3,0) -- cycle;
\draw[thick] (0,4) -- (0,4/3) -- (4/3,0) -- (4,0);
\draw[dashed] (1,3) -- (1,2) node {\(\times\)};
\draw[dashed] (3,1) -- (2,1) node {\(\times\)};
\draw[very thick] (2,1) -- (1,1) -- (1,2);
\draw[very thick] (1,1) -- (2/3,2/3);
\draw[dotted,thick] (4/3,0) -- (0,0) -- (0,4/3);
\node at (2/3,0) [below] {\(c\)};
\node at (0,2/3) [left] {\(c\)};
\draw[dotted,thick] (0,2) -- (1,2) node [midway,above] {\(a\)};
\draw[dotted,thick] (2,0) -- (2,1) node [midway,right] {\(b\)};
\begin{scope}[shift={(5,0)}]
\filldraw[fill=gray!25,draw=none] (0,2/3) -- (0,4) -- (4,0) -- (2/3,0) -- cycle;
\draw[thick] (0,4) -- (0,2/3) -- (2/3,0) -- (4,0);
\draw[dashed] (2/3,10/3) -- (2/3,2) node {\(\times\)};
\draw[dashed] (7/3,5/3) -- (2,5/3) node {\(\times\)};
\draw[very thick] (2,5/3) -- (2/3,5/3) -- (2/3,2);
\draw[very thick] (0,3/3) -- (2/3,5/3);
\draw[dotted,thick] (2/3,0) -- (0,0) -- (0,2/3);
\node at (1/3,0) [below] {\(c\)};
\node at (0,1/3) [left] {\(c\)};
\draw[dotted,thick] (0,2) -- (2/3,2) node [midway,above] {\(a\)};
\draw[dotted,thick] (2,0) -- (2,5/3) node [midway,right] {\(b\)};
\end{scope}
\end{tikzpicture}
\caption{Two almost toric blow-ups of a symplectic ball in three smaller balls (Example \ref{exm:shevsmi}). (a) The symplectic triangle inequalities hold, and we find a tropical Lagrangian \(\rp{2}\). (b) The symplectic triangle inequalities fail and we find a tropical Lagrangian disc instead. (Figure taken from {\cite[Figure 1]{EvansKB}}).}
\label{fig:shevsmi}
\end{center}
\end{figure}
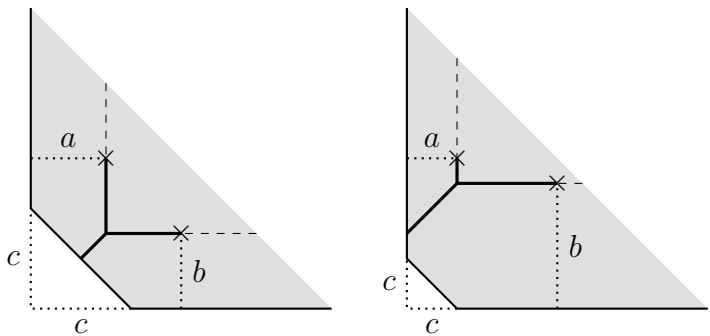

\end{Example}
\index{almost toric base diagram!symplectic triangle
inequality|)}

\section{Solutions to inline exercises}

\begin{Exercise}[Lemma \ref{lma:hyperkahler_twist}]\label{exr:hyperkahler_twist}
Consider the {\em hyperK\"{a}hler twist}\index{hyperKaehler
twist@hyperK\"{a}hler twist} \[\RR^2\times T^2\to\RR^2\times
T^2,\quad (p_1,p_2,q_1,q_2)\mapsto (p_1,p_2,-q_2,q_1)\] If
\(C\subset\RR^2\times T^2\) is a complex curve with respect to
the complex coordinates \(z_k=e^{p_k+iq_k}\) then show that
the image of \(C\) under this hyperK\"{a}hler twist is
Lagrangian for the symplectic form \(\sum dp_i\wedge dq_i\).
\end{Exercise}
\begin{Solution}
Let \(z=x+iy\) be a local complex coordinate on \(C\) and
suppose \((p_1(z),p_2(z),q_1(z),q_2(z))\) is a local
parametrisation of \(C\), holomorphic with respect to the
complex coordinates \(z_k=e^{p_k+iq_k}\) Since (any branch of)
the logarithm is holomorphic, this is equivalent to the
Cauchy-Riemann equations \[\frac{\partial p_k}{\partial
x}=\frac{\partial q_k}{\partial y},\quad \frac{\partial
p_k}{\partial y}=-\frac{\partial q_k}{\partial x},\mbox{ for
}k=1,2.\] After applying the twist, we get a submanifold which
satisfies \[\frac{\partial p_1}{\partial x}=\frac{\partial
q_2}{\partial y},\quad \frac{\partial p_1}{\partial
y}=-\frac{\partial q_2}{\partial x},\quad\frac{\partial
p_2}{\partial x}=-\frac{\partial q_1}{\partial y},\quad
\frac{\partial p_2}{\partial y}=\frac{\partial q_1}{\partial
x}.\] This means that if \(\partial_x\) pushes forward along
the twisted embedding to \((a,b,c,d)\) then \(\partial_y\)
pushes forward to \((-d,c,-b,a)\), and \[\omega((a,b,c,d),
(-d,c,-b,a))=-ab+cd+ab-cd = 0,\] so this twisted embedding is
Lagrangian.\qedhere

\end{Solution}
\begin{Exercise}[Example \ref{exm:wtf_is_it}]\label{exr:wtf_is_it}
What is the tropical Lagrangian in Figure \ref{fig:wtf_is_it}?
\end{Exercise}
\begin{Solution}
Where the tropical curve meets the corner and edges, it gives
us a Lagrangian disc and two Lagrangian M\"{o}bius
strips. These are used to cap off the pair-of-pants coming
from the trivalent vertex. Therefore this tropical Lagrangian
is an immersed Klein bottle\index{Lagrangian!Klein bottle}
\index{Lagrangian!immersed}. The number of transverse double
points is given by the \(\delta\)-invariant of the trivalent
vertex. The outgoing vectors at this point are: \[v_1 =
(-2,3),\quad v_2 = (3,-2),\quad v_3 = (-1,-1),\] so
\(|v_1\wedge v_2|=\left|\det\begin{pmatrix} -2 & 3 \\ 3 &
-2\end{pmatrix}\right|=5\) and \(\delta =
(5-1)/2=2\). Therefore this Klein bottle has two double
points. \qedhere

\end{Solution}
\begin{Exercise}[Remark \ref{rmk:indep_order}]\label{exr:indep_order}
If \(v_1,v_2,v_3\) are primitive integer vectors with
\(v_1+v_2+v_3=0\) then \(|v_k\wedge v_{\ell}|\) is an odd
number and is independent of \(k,\ell\).
\end{Exercise}
\begin{Solution}
If we switch the rows \(v_1\) and \(v_2\) then we just change
the sign of the determinant. If we replace \(v_2\) by \(v_3\)
then the balancing condition means \(v_3=-v_1-v_2\), so
\[|v_1\wedge v_3| = |v_1\wedge (v_1+v_2)| = |v_1\wedge v_2|.\]
Similarly the result is unchanged if we replace \(v_1\) by
\(v_3\).

To see that it is odd, write \(v_1=(a,b)\) and
\(v_2=(c,d)\). Suppose for contradiction that \(|v_1\wedge
v_2|=|ad-bc|\) is even. We claim that \(a,b,c,d\) are then all
odd. This will give a contradiction because it implies
\(v_3=(-a-c,-b-d)\) is divisible by \(2\) and so not
primitive.

To prove that \(|ad-bc|\) even implies \(a,b,c,d\) are all
odd, suppose that \(a\) is even (a similar argument works for
\(b,c,d\)). Then \(ad-bc=0\mod 2\) implies that either \(b\)
or \(c\) is even. But if \(a\) is even then \(b\) is odd by
primitivity of \(v_1\), so \(c\) is even. This means \(d\) is
odd. But then \(v_3=(-a-c,-b-d)\) has both components even,
which contradicts primitivity.\qedhere

\end{Solution}
\begin{Exercise}[Lemma \ref{lma:tropical_double_points}]\label{exr:tropical_double_points}
Suppose we have several straight lines of rational slope in
\(\RR^2\) incident on a point \(b\in B\). Let \(v_1,\ldots,k\)
be primitive integer vectors pointing along these lines. Show
that the visible Lagrangian cylinders above these lines have a
total of \(\delta(b)\) transverse intersections, where
\[\delta(b)=\sum_{i< j}|v_i\wedge v_j|.\]
\end{Exercise}
\begin{Solution}
If \(v_j=(m_j,n_j)\) for \(j=1,\ldots,j\) then the visible
Lagrangian cylinders intersect the fibre over \(b\) in the
circles (in \(\RR^2/\ZZ^2\)) with slopes \((-n_i,m_i)\). The
number of intersections between two of these circles is
\(|m_2n_1-m_1n_2|\). Summing this over pairs of circles gives
\(\delta(b)\).

\begin{figure}[htb]
\begin{center}
\begin{tikzpicture}
\filldraw[lightgray,opacity=0.5] (0,0) -- (3,0) -- (3,3) -- (0,3) -- cycle;
\draw[->-] (0,0) -- (1.5,0);
\draw (1.5,0) -- (3,0);
\draw[->-] (0,3) -- (1.5,3);
\draw (1.5,3) -- (3,3);
\draw[->>-] (0,0) -- (0,1.5);
\draw (0,1.5) -- (0,3);
\draw[->>-] (3,0) -- (3,1.5);
\draw (3,1.5) -- (3,3);
\draw[very thick] (0,0) -- (3,3);
\draw[very thick,dash dot] (1.5,3) -- (3,1.5);
\draw[very thick,dash dot] (0,1.5) -- (1.5,0);
\end{tikzpicture}
\caption{Intersections between visible Lagrangian cylinders in the fibre over a point where their projections meet.}
\label{fig:exr_circles_intersect}
\end{center}
\end{figure}
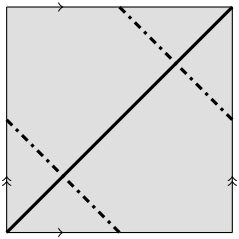

For example, suppose the lines have slopes \(1\) and
\(-1\). Then the circles in the fibre have slopes \(-1\) and
\(1\) respectively. If we draw these as lines in the square
with opposite sides identified (see Figure
\ref{fig:exr_circles_intersect}), we see that these intersect
twice.\index{Lagrangian!submanifold!tropical|)} \qedhere

\end{Solution}
\chapter{Markov triples}
\label{ch:markov_triples}
\thispagestyle{cup}

\section{The Markov equation}

The\index{Markov triple|(} Diophantine equation
\[a^2+b^2+c^2=3abc\] is called {\em the Markov equation}. It
occurs in the theory of Diophantine approximation and continued
fractions \cite{CasselsDiophantine}, the theory of quadratic
forms \cite{CasselsDiophantine}, the study of exceptional
collections \cite{Rudakov}, the hyperbolic geometry of a
punctured torus \cite{SeriesMarkoff}, in governing the
\(\QQ\)-Gorenstein degenerations\index{Q-Gorenstein
degeneration@$\mathbb{Q}$-Gorenstein degeneration} of \(\cp{2}\)
\cite{Manetti,HackingProkhorov}, and elsewhere; a wonderful
exposition can be found in Aigner's book \cite{Aigner}. A triple
of positive integers solving Markov's equation is called a {\em
Markov triple}.

\begin{Lemma}
If \(a,b,c\) is a Markov triple then so is \(a,b,3ab-c\).
\end{Lemma}
\begin{Proof}
Fix \(a\) and \(b\) and consider the quadratic function
\(f(x):=x^2-3abx+a^2+b^2\). This quadratic has \(c\) as a
root: \(f(c)=0\). A quadratic equation \(x^2+\beta
x+\gamma=0\) has two roots (counted with multiplicity) which
sum to \(-\beta\). In our case, \(\beta=-3ab\), so this means
that \(3ab-c\) is another root. To see that both roots are
positive, note that \(f(0)=a^2+b^2>0\) and
\(f'(0)=-3ab<0\). This means that \(f(x)>0\) for all \(x<0\)
(see Figure \ref{fig:graph_markov}) so there are no negative
roots. \qedhere

\end{Proof}
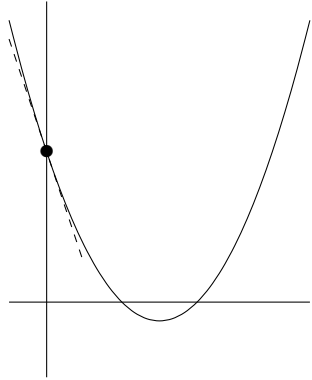
\begin{figure}[htb]
\begin{center}
\begin{tikzpicture}
\draw (-1/2,0) -- (3.5,0);
\draw (0,-1) -- (0,4);
\draw[black] plot[smooth,domain=-0.5:3.5] (\x, {\x*\x-3*\x+2});
\node at (0,2) {\(\bullet\)};
\draw[black,dashed] plot[smooth,domain=-0.5:0.5] (\x,{-3*\x+2});
\end{tikzpicture}
\end{center}
\caption{The graph of \(f\) when \(a=b=1\). Positive \(y\)-intercept and negative gradient at \(x=0\) ensures both roots are positive.}
\label{fig:graph_markov}
\end{figure}

\begin{Definition}
The operation of replacing the Markov triple \((a,b,c)\) with
\((a,b,3ab-c)\) is called a {\em mutation}\index{mutation!of
Markov triple} on \(c\). The graph whose vertices are Markov
triples and whose edges connect triples related by a mutation
is called the {\em Markov graph}. In fact, this graph is a
connected tree and we will usually refer to it as the {\em
Markov tree}\index{Markov tree|(}.

\end{Definition}
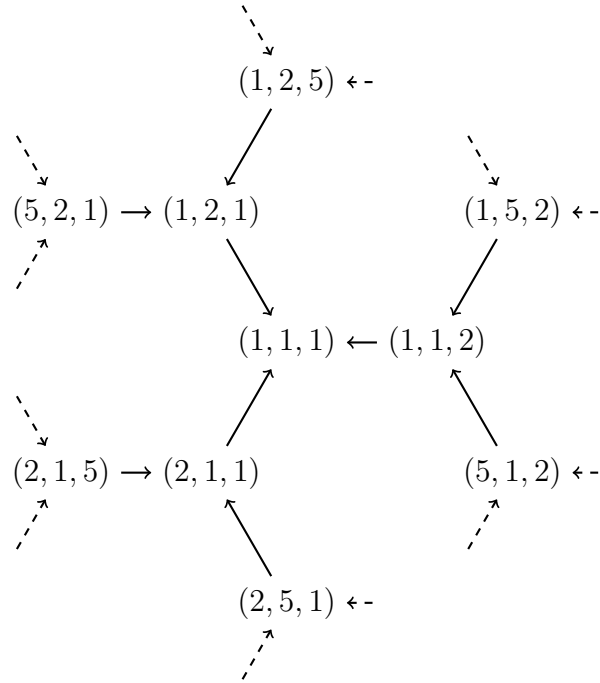
\begin{figure}[htb]
\begin{center}
\begin{tikzpicture}
\node (111) at (0,0) {\((1,1,1)\)};
\node (112) at (0:2) {\((1,1,2)\)};
\begin{scope}[shift={(0:2)}]
\node (152) at (60:2) {\((1,5,2)\)};
\draw[thick,<-,dashed] (152) -- ++ (0:1.2);
\draw[thick,<-,dashed] (152) -- ++ (120:1.2);
\node (512) at (-60:2) {\((5,1,2)\)};
\draw[thick,<-,dashed] (512) -- ++ (0:1.2);
\draw[thick,<-,dashed] (512) -- ++ (240:1.2);
\end{scope}
\node (121) at (120:2) {\((1,2,1)\)};
\begin{scope}[shift={(120:2)}]
\node (125) at (60:2) {\((1,2,5)\)};
\draw[thick,<-,dashed] (125) -- ++ (0:1.2);
\draw[thick,<-,dashed] (125) -- ++ (120:1.2);
\node (521) at (180:2) {\((5,2,1)\)};
\draw[thick,<-,dashed] (521) -- ++ (120:1.2);
\draw[thick,<-,dashed] (521) -- ++ (240:1.2);
\end{scope}
\node (211) at (240:2) {\((2,1,1)\)};
\begin{scope}[shift={(240:2)}]
\node (215) at (180:2) {\((2,1,5)\)};
\draw[thick,<-,dashed] (215) -- ++ (240:1.2);
\draw[thick,<-,dashed] (215) -- ++ (120:1.2);
\node (251) at (-60:2) {\((2,5,1)\)};
\draw[thick,<-,dashed] (251) -- ++ (0:1.2);
\draw[thick,<-,dashed] (251) -- ++ (240:1.2);
\end{scope}
\draw[thick,<-] (111) -- (112);
\draw[thick,<-] (111) -- (121);
\draw[thick,<-] (111) -- (211);
\draw[thick,<-] (112) -- (152);
\draw[thick,<-] (112) -- (512);
\draw[thick,<-] (121) -- (125);
\draw[thick,<-] (121) -- (521);
\draw[thick,<-] (211) -- (215);
\draw[thick,<-] (211) -- (251);
\end{tikzpicture}
\caption{A small part of the Markov tree (dashed lines indicate it continues). The arrows are explained in Lemma \ref{lma:markov_orientation}.}
\label{fig:markov_tree}
\end{center}
\end{figure}

\begin{Theorem}\label{thm:markov_connected}
The Markov graph is connected.
\end{Theorem}
\begin{Proof}
Suppose we are given a Markov triple \((a,b,c)\) with
\(c>a,b\). Perform a mutation
on the largest element \(c\). This decreases the value of the
largest element (Lemma \ref{lma:markov_decrease} below). If
there is still a unique largest element, repeat this
procedure. Since the values in the triple are decreasing, but
always positive, this process terminates, and we find a Markov
triple with a repeated largest element. The only triple with a
repeated largest element is \((1,1,1)\) (Lemma
\ref{lma:markov_unique_lowest_triple} below). This shows that
the graph is connected. \qedhere

\end{Proof}
\begin{Lemma}\label{lma:markov_decrease}
Let \(a,b,c\) be a Markov triple with \(a\leq b\leq c\). Then
\(b\) lies in the closed interval between \(c\) and
\(3ab-c\). If \(a<b\) then \(b\) lies in the interior of this
interval. In particular, if \(b<c\) then \(3ab-c\leq b\).
\end{Lemma}
\begin{Proof}
If \(f(x)=x^2-3abx+a^2+b^2\) then
\(f(b)=b^2-3ab^2+a^2+b^2=(2-3a)b^2+a^2\leq a^2-b^2\leq
0\). This means that \(b\) lies in the region between the two
roots of \(f\), and strictly between if \(a^2-b^2<0\). These
roots are \(c\) and \(3ab-c\). \qedhere

\end{Proof}
\begin{Lemma}\label{lma:markov_unique_lowest_triple}
The only Markov triple with no unique largest element is
\((1,1,1)\).
\end{Lemma}
\begin{Proof}
If \((a,b,c)\) is a Markov triple with \(a\leq b=c\) then
substituting \(b=c\) in the Markov equation gives
\(a^2+2b^2=3ab^2\), or \(a^2=(3a-2)b^2\). Since \(a\geq 1\),
\(3a-2\geq 1\) and hence \(a^2\geq b^2\). Since \(a\leq b\),
we get \(a=b=c\). The Markov equation now reduces to
\(a^2=(3a-2)a^2\), so \(3a-2=1\) and \(a=b=c=1\).\qedhere

\end{Proof}
\begin{Remark}
We briefly recall the definition of a tree. A {\em path} in a
graph is a sequence of oriented edges such that the endpoint
of edge \(i\) is the start-point of edge \(i+1\) for all
\(i\). A path is {\em non-simple} if the sequence of edges
contains a subsequence of the form \(e\bar{e}\) where \(e\) is
an edge and \(\bar{e}\) is the same edge with the opposite
orientation. Otherwise, a path is called {\em simple}. A graph
is a {\em tree} if any pair of vertices can be connected by a
unique simple path.

\end{Remark}
\begin{Lemma}\label{lma:markov_orientation}
The following prescription defines a global choice of
orientations on the edges of the Markov graph. If an edge
connects two triples \((a,b,c)\) and \((a,b,3ab-c)\) then
orient it so that it points towards the triple whose maximal
element is smaller. For example \((1,1,2)\to (1,1,1)\).
\end{Lemma}
\begin{Proof}
To show this is well-defined, we need to show that the maximal
elements in \((a,b,c)\) and \((a,b,3ab-c)\) are different
(otherwise there is no way to decide the orientation of the
edge). Without loss of generality, suppose \(a\leq b\) and
\(c<3ab-c\). Lemma \ref{lma:markov_decrease} tells us that
\(c\leq b\leq 3ab-c\). Therefore \(\max(a,b,c)=b\) and
\(\max(a,b,3ab-c)=3ab-c\). The only problem is if
\(b=3ab-c\). This is only possible if \(b=a\) by Lemma
\ref{lma:markov_decrease}, but in this case the triple
\((a,b,c)\) has a repeated maximum \(a=b\), so
\((a,b,c)=(1,1,1)\) and all the neighbouring triples are
permutations of \((1,1,2)\), so it is clear how to orient the
graph around this vertex. \qedhere

\end{Proof}
\begin{Remark}
The vertex \((1,1,1)\) has three incoming edges. The proof of
Lemma \ref{lma:markov_orientation} shows that every other
vertex has two incoming edges (corresponding to mutations on
the smaller elements of the triple) and one outgoing edge
(corresponding to mutation on the largest element).

\end{Remark}
\begin{Theorem}
The Markov graph is a tree.
\end{Theorem}
\begin{Proof}
We have seen that the Markov graph is connected. It remains to
show that two triples can be connected by a {\em unique}
simple path. Orient the Markov graph according to Lemma
\ref{lma:markov_orientation}. We call a path {\em downwards}
if it follows the orientation. The proof of Theorem
\ref{thm:markov_connected} can be interpreted as saying that
there is a unique downwards path from every vertex to
\((1,1,1)\). Define the {\em height} of a triple to be the
length of this unique downwards path. Similarly, we can define
the height of an edge to be the height of its start-point.

Fix two Markov triples \(t_1\) and \(t_2\). Follow the paths
down from \(t_1\) and from \(t_2\). Eventually these paths
meet up at a vertex (this could be at \((1,1,1)\) or somewhere
higher up). Let \(m\) be the first vertex where these paths
meet; we get a simple path \(P\) by going down from \(t_1\) to
\(m\) and then up from \(m\) to \(t_2\).

Suppose there is another simple path, \(Q\). If \(P\neq Q\)
then there is an edge in \(Q\) which does not appear in
\(P\). Amongst the edges of \(Q\) which do not appear in
\(P\), let \(e_i\) be one which maximises height. There are
three possible cases, each leading to a contradiction:
\begin{itemize}
\item Suppose \(e_i=e_1\) is the first edge in \(Q\). Since it is
not in \(P\), it cannot be the downward arrow from \(t_1\),
so it must go up. Since it is a highest edge, the next edge
\(e_2\) must return along \(\bar{e}_1\) as there is only one
way down. This contradicts the fact that \(Q\) is simple.
\item A similar argument applies when \(e_i\) is the final edge in
\(Q\).
\item If \(e_i\) is neither the first nor last edge in \(Q\) then
there are edges \(e_{i-1}\) and \(e_{i+1}\) adjacent in the
path. At least one of these must start at the highest point
of \(e_i\). Since it cannot go higher, it must be
\(\bar{e}_i\), which contradicts the assumption that \(Q\)
is simple.
\end{itemize}
Therefore there is a unique simple path between two vertices
and the Markov graph is a tree.\index{Markov
triple|)}\index{Markov tree|)}\qedhere

\end{Proof}
\section{Triangles}

A {\em Vianna triangle}\index{Vianna!triangle|(} is an almost
toric diagram whose edges are \(v_1,v_2,v_3\) with affine
lengths\index{affine length!Vianna triangles|(}
\(\ell_1,\ell_2,\ell_3\) and whose vertices \(P_1,P_2,P_3\) are
modelled on the \(T\)-singularities
\(\frac{1}{d_kp_k^2}(1,d_kp_kq_k-1)\) for \(k=1,2,3\). The
relative orientations and positions of the edges and vertices
are shown below. We call the numbers \(d_k,p_k,q_k,\ell_k\),
\(k=1,2,3\), the {\em Vianna data}\index{Vianna!data} of the
triangle.

\begin{center}
\begin{tikzpicture}
\node (a) at (0,0) {\(\bullet\)};
\node (b) at (3,1) {\(\bullet\)};
\node (c) at (1,2) {\(\bullet\)};
\draw[->-,thick] (a.center) -- (b.center);
\draw[->-,thick] (b.center) -- (c.center);
\draw[->-,thick] (c.center) -- (a.center);
\node at (1.5,0.5) [below right] {\(v_1\)};
\node at (2,1.5) [above right] {\(v_2\)};
\node at (0.5,1) [above left] {\(v_3\)};
\node at (c) [above] {\(P_1\)};
\node at (a) [below left] {\(P_2\)};
\node at (b) [right] {\(P_3\)};
\end{tikzpicture}
\end{center}

Our indices take values in the cyclic group \(\ZZ/3\).

\begin{Lemma}\label{lma:useful_eq}
If \(\hat{v}_k\) denotes the primitive integer vector along
\(v_k\) and \(\ell_k\) denotes the affine length of \(v_k\)
then we have the following relation\footnote{Recall that
\((a,b)\wedge (c,d) := ad-bc\).}: \[\hat{v}_k\wedge
\hat{v}_{k+1} = d_{k+2}p_{k+2}^2.\]
\end{Lemma}
\begin{Proof}
Since the vertex \(P_{k+2}\) is modelled on the moment polygon
of the
\(\frac{1}{d_{k+2}p_{k+2}^2}(1,d_{k+2}p_{k+2}q_{k+2}-1)\)
singularity, there is an integral affine transformation making
\(\hat{v}_k=(0,-1)\) and \(\hat{v}_{k+1} = (d_{k+2}p_{k+2}^2,
d_{k+2}p_{k+2}q_{k+2} - 1)\). Thus \(\hat{v}_k\wedge
\hat{v}_{k+1} = (0,-1)\wedge
(d_{k+2}p_{k+2}^2,d_{k+2}p_{k+2}q_{k+2}-1) =
d_{k+2}p_{k+2}^2\).\qedhere

\end{Proof}
\begin{Corollary}\label{cor:common_value_vianna}
We have \[\ell_1\ell_2d_3p_3^2 = \ell_2\ell_3d_1p_1^2 =
\ell_3\ell_1d_2p_2^2.\]
\end{Corollary}
\begin{Proof}
Since \(v_1+v_2+v_3=0\), we get \[0=(v_1+v_2+v_3)\wedge v_k =
v_{k-1}\wedge v_k+0+v_{k+1}\wedge v_k=v_{k-1}\wedge v_k-
v_k\wedge v_{k+1}.\] The quantity \(v_k\wedge v_{k+1}\) is
therefore independent of \(k\). Since \(v_k=\ell_k\hat{v}_k\)
and \(v_{k+1}=\ell_{k+1}\hat{v}_{k+1}\), the previous lemma
implies this is \(\ell_k\ell_{k+1}d_{k+2}p_{k+2}^2\). \qedhere

\end{Proof}
\begin{Definition}
We write \(K\) for the common value
\(\ell_k\ell_{k+1}d_{k+2}p_{k+2}^2\) and \(L\) for the total
affine length \(\ell_1 + \ell_2 + \ell_3\).

\end{Definition}
\begin{Corollary}\label{cor:formula_for_aff_length}
We have \(\ell_k = \frac{p_k}{p_{k+1}p_{k+2}}
\sqrt{\frac{Kd_k}{d_{k+1}d_{k+2}}}\).
\end{Corollary}
\begin{Proof}
For concreteness, take \(k=3\). By Corollary
\ref{cor:common_value_vianna}, We have
\[\ell_1=\frac{K}{\ell_3d_2p_2^2},\quad
\ell_2=\frac{K}{\ell_3d_1p_1^2},\quad \ell_1\ell_2d_3p_3^2 =
K,\] so \[\frac{K^2d_3p_3^2}{\ell_3^2d_1d_2p_1^2p_2^2}=K,\]
which gives \(\ell_3^2 = (Kd_3/d_1d_2)(p_3^2/p_1^2p_2^2)\) as
required.\qedhere

\end{Proof}
\begin{Corollary}\label{cor:markov_type_eq}
We have \[d_1p_1^2+d_2p_2^2+d_3p_3^2 =
\frac{L\sqrt{d_1d_2d_2}}{\sqrt{K}}p_1p_2p_3.\]
\end{Corollary}
\begin{Proof}
This follows from Corollary \ref{cor:formula_for_aff_length}
and the fact that \(\ell_1+\ell_2+\ell_3=L\).\qedhere

\end{Proof}
We will now prove a sequence of lemmas to show that the
constants \(K\) and \(L\) are unchanged by
mutation\index{mutation!of Vianna triangle|(}.

\begin{Lemma}
The eigenline\index{eigenline!for Vianna triangle|(} at vertex
\(P_{k+2}\) points in the direction
\(\frac{\hat{v}_{k+1}-\hat{v}_k}{d_{k+2}p_{k+2}}\).
\end{Lemma}
\begin{Proof}
Again, making an integral affine transformation we can assume
\(\hat{v}_k=(0,-1)\) and
\(\hat{v}_{k+1}=(d_{k+2}p_{k+2}^2,d_{k+2}p_{k+2}q_{k+2}-1)\). In
these coordinates, the eigenline points in the
\((p_{k+2},q_{k+2})\)-direction. which is
\((\hat{v}_{k+1}-\hat{v}_k)/d_{k+2}p_{k+2}\).\qedhere

\end{Proof}
\begin{Lemma}
If we perform a mutation on the vertex \(P_3\) then we obtain
a new Vianna triangle with data (omitting the \(q_k\)'s):
\begin{align*}
d'_1&=d_1,& d'_2&=d_2,&d'_3&=d_3\\
p'_1&=p_1,& p'_2&=p_2,&p'_3 &= (dp_1^2+d_2p_2^2)/d_3p_3\\
\ell'_1 &=\frac{\ell_3d_2p_2^2}{d_1p_1^2+d_2p_2^2},& \ell'_2 &=
\frac{\ell_3d_1p_1^2}{d_1p_1^2+d_2p_2^2},&\ell'_3 &=
\ell_1+\ell_2.\end{align*}
\end{Lemma}
\begin{Proof}
When we mutate, a new vertex is introduced at the point
\(P'_3\) where the edge \(v_3\) intersects the eigenline
emanating out of \(P_3\) in the
\(\frac{\hat{v}_2-\hat{v}_1}{d_3p_3}\) direction. We find the
new \(p'_3\) by taking
\(\frac{\hat{v}_2-\hat{v}_1}{d_3p_3}\wedge \hat{v}_3\), which
gives \[p'_3=\frac{1}{d_3p_3}(\hat{v}_2\wedge
\hat{v}_3-\hat{v}_1\wedge\hat{v}_3) =
\frac{d_1p_1^2+d_2p_2^2}{d_3p_3}.\] This intersection point
\(P'_3\) lives on the line joining \(P_1\) and \(P_2\). If we
choose coordinates where \(P_3\) is the origin, \(P_1\) and
\(P_2\) are the vectors \(v_2\) and \(-v_1\) respectively. Let
\(s\in[0,1]\) be the number such that
\(P_3=-sv_1+(1-s)v_2\). This lives on the
eigenline\index{eigenline!for Vianna triangle|)}, so can be
written as \(t(\hat{v}_2-\hat{v}_1)\) for some \(t\). Since
\(v_1\) and \(v_2\) are linearly independent, we deduce that
\(-sv_1=-t\hat{v}_1\) and \((1-s)v_2=t\hat{v}_2\). Eliminating
\(t\) gives \(s/(1-s) = \ell_2/\ell_1\). We have
\(\ell_2/\ell_1=d_2p_2^2/d_1p_1^2\) by Corollary
\ref{cor:formula_for_aff_length}, so we get
\(s=\frac{d_2p_2^2}{d_1p_1^2+d_2p_2^2}\). Thus
\(P'_3=\frac{d_1p_1^2v_2-d_2p_2^2v_1}{d_1p_1^2+d_2p_2^2}\). The
vector from \(P_1=v_2\) to \(P'_3\) is therefore
\(-\frac{d_2p_2^2}{d_1p_1^2+d_2p_2^2}(v_1+v_2) =
\frac{d_2p_2^2}{d_1p_1^2+d_2p_2^2}v_3\). Thus, the edge
\(v_3\) is subdivided so that its affine length \(\ell_3\) is
split in the ratio \(d_1p_1^2:d_2p_2^2\). \qedhere

\end{Proof}
\begin{Corollary}\label{cor:mut_inv}
The common value \(\ell_1\ell_2d_3p_3^2= \ell_2\ell_3d_1p_1^2
= \ell_3\ell_1d_2p_2^2\) is unchanged by mutation at \(P_3\)
(or, by symmetry, any other vertex). Similarly, the sum
\(\ell_1+\ell_2+\ell_3\) is unchanged by mutation.
\end{Corollary}
\begin{Proof}
We compute the value \(\ell'_1\ell'_2d'_3(p'_3)^2\) after
mutation: \[\ell'_1\ell'_2d'_3(p'_3)^2 =
\frac{\ell_3d_2p_2^2\cdot\ell_3d_1p_1^2}{(d_1p_1^2+d_2p_2^2)^2}\
d_3\left(\frac{d_1p_1^2 + d_2p_2^2}{d_2p_3}\right)^2 =
\frac{\ell_3^2d_1p_1^2d_2p_2^2}{d_3p_3^2}.\] By Corollary
\ref{cor:common_value_vianna},
\(\ell_3d_1p_1^2=\ell_1d_3p_3^2\) and
\(\ell_3d_2p_2^2=\ell_2d_3p_3^2\), so this reduces to
\(\ell_1\ell_2d_3p_3^2\), which is the same as before
mutation.

The sum of the affine lengths after mutation is
\[(\ell_1+\ell_2) + \frac{\ell_3d_2p_2^2}{d_1p_1^2+d_2p_2^2} +
\frac{\ell_3d_1p_1^2}{d_1p_1^2+d_2p_2^2} =
\ell_1+\ell_2+\ell_2.\qedhere\]

\end{Proof}
\begin{Example}[Proof of Theorem \ref{thm:vianna_tri}]\label{exm:vianna_mutations}
The triangle \(D(1,1,1)\) in Figure \ref{fig:diag_111} is a
Vianna triangle with data \(d_1=d_2=d_3=1\), \(p_1=p_2=p_3=1\)
and \(\ell_1=\ell_2=\ell_3=3\). This has \(K=9\) and \(L=9\),
so\index{affine length!Vianna triangles|)}
Corollaries \ref{cor:markov_type_eq} and \ref{cor:mut_inv}
tell us that if \(D\) is obtained from \(D(1,1,1)\) by
iterated mutation and has Vianna data \(d_k,p_k,\ell_k\) then
\[p_1^2+p_2^2+p_3^2 = 3p_1p_2p_3\] and \(\ell_k =
3p_k/(p_{k+1}p_{k+2})\). Moreover, if we write
\(D(p_1,p_2,p_3)\) for the triangle associated with the Markov
triple\index{Markov triple} \(p_1,p_2,p_3\), mutation at
vertex \(3\) gives the Markov triple
\(p_1,p_2,p'_3=3p_1p_2-p_3\). One can check this from the
formulas, or simply observe that this mutation leaves \(p_1\)
and \(p_2\) unchanged, and there are only two Markov triples
containing both \(p_1\) and \(p_2\).\index{mutation!of Vianna
triangle|)}\index{Vianna!triangle|)}

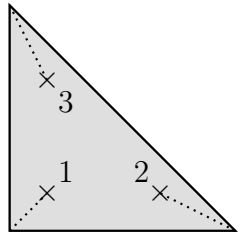
\begin{figure}[htb]
\begin{center}
\begin{tikzpicture}
\filldraw[fill=lightgray,opacity=0.5,draw=black,thick] (0,0) -- (0,3) -- (3,0) -- cycle;
\draw[thick] (0,0) -- (0,3) -- (3,0) -- cycle;
\draw[dotted,thick] (0,0) -- (1/2,1/2);
\draw[dotted,thick] (0,3) -- (0.5,2);
\draw[dotted,thick] (3,0) -- (2,0.5);
\node (a) at (1/2,2) {\(\times\)};
\node (b) at (2,1/2) {\(\times\)};
\node (c) at (1/2,1/2) {\(\times\)};
\node at (a) [below right] {\(3\)};
\node at (b) [above left] {\(2\)};
\node at (c) [above right] {\(1\)};
\end{tikzpicture}
\end{center}
\caption{The almost toric diagram \(D(1,1,1)\). The edges all have affine length \(3\).}
\label{fig:diag_111}
\end{figure}

\end{Example}
\chapter{Open problems}
\label{ch:open_problems}
\thispagestyle{cup}

In the preface, I mentioned that this book is intended to give
you the tools to explore further in symplectic geometry. In this
final appendix, you will find a handful of open\footnote{at the
time of writing.} problems which you might like to
explore. These are some of my pet problems: not particularly
venerable or widely known, but they have intrigued and
frustrated me over the years. The first few are hopefully purely
``combinatorial'' or geometric in that they should not require
any Floer theory\index{Floer theory} to solve. Some will
probably require Floer theory, and I have not hesitated to use
terminology from beyond this book in discussing these problems,
because you will need to learn it to solve them.

\section{Mutation of quadrilaterals}

\begin{Problem}\label{prb:mutation}
Which\index{mutation!quadrilateral} quadrilaterals arise as
mutations of a square? Or a rectangle? For a rectangle, how
does the ratio of the side-lengths affect the answer?

\end{Problem}
It is easy to generate such quadrilaterals by iterated mutation,
but I know of no succinct description of this class of
quadrilaterals analogous to the Markov-triple\index{Markov
triple} description of mutants of the \(\cp{2}\) triangle.

\section{Fillings of lens spaces}

Equip the lens space\index{lens space|(}\index{lens
space!symplectic fillings|(} \(L(n,a)\) with its standard
contact structure. It is known that \(L(n,a)\) has a filling
with second Betti number zero if and only if \(n=p^2\) and
\(a=pq-1\) for some coprime integers \(p,q\); moreover, such a
filling is unique up to symplectomorphism. It is also known
({\cite[Theorem 4.3]{HTU}}, {\cite[Theorem 2.8]{UrzuaVilches}})
that \(L(n,a)\) admits at most two fillings (up to
deformation/symplectomorphism) with \(b_2=1\). For example, in
Remark \ref{rmk:non_diffeo_fillings} we saw that \(L(36,13)\)
admits two non-diffeomorphic fillings.

\begin{Problem}
Characterise which (standard contact) lens spaces admit
precisely two symplectic fillings with \(b_2=1\).

\end{Problem}
This is a purely combinatorial problem, thanks to Lisca's
classification. Again, it is easy to generate examples (even
infinite families of examples) and there are many tantalising
patterns, but I do not know a complete answer. These examples
arise very naturally in the context of algebraic geometry; the
papers \cite{UrzuaNeighbours} and \cite{UrzuaVilches} have many
examples of surfaces of general type containing a lens space
hypersurface which are related (topologically) by the surgery
which interchanges the two fillings. Urz\'{u}a and Vilches
\cite{UrzuaVilches} call these {\em wormholes}\index{wormhole}
in the moduli space of surfaces.\index{lens space|)}\index{lens
space!symplectic fillings|)}

\section{Minimal genus problem}

Let \((X,\omega_\lambda)\) be the symplectic manifold
\(S^2\times S^2\) with a product symplectic form giving the
factors areas \(1\) and \(\lambda\). Let \(\beta\in
H_2(S^2\times S^2;\ZZ/2)\) be the homology class of
\(S^2\times\{\OP{pt}\}\).

\begin{Problem}\label{prb:genus}
As a function of \(\lambda\), what is the minimal genus of a
(nonorientable) embedded Lagrangian
submanifold\index{Lagrangian!submanifold} in \(X\) inhabiting
the \(\ZZ/2\)-homology class \(\beta\)? Is this genus
uniformly bounded in \(\lambda\)?

\end{Problem}
Here, the genus of the nonorientable surface \(\sharp_m\rp{2}\)
is defined to be \(m\). One can find \cite{EvansKB} explicit
tropical Lagrangians\index{Lagrangian!submanifold!tropical} in
this homology class whose genus goes to infinity with \(m\), and
it seems difficult (to me) to do significantly better than this
with tropical Lagrangians.

Proving a lower bound on the minimal genus will be difficult and
most likely require techniques of Floer theory\index{Floer
theory}, but you can give an upper bound just by constructing
examples.

Obviously this problem admits many generalisations (just pick a
different manifold \(X\)) but this seems to be the simplest
version which is open.

\section{Nodal slides in higher dimensions}

In this book, we focused on symplectic 4-manifolds fibring over
2-dimensional integral affine bases. Almost toric fibrations
make sense in higher dimensions too, and the discriminant
locus\index{discriminant locus|(} has codimension 2. For
example, if the base is 3-dimensional then the discriminant
locus will be some kind of knotted graph.

A nodal slide\index{nodal slide} is a kind of isotopy of the
discriminant locus. If the base is 2-dimensional, you can
usually slide base-nodes out of the way of another sliding
node\footnote{unless their eigenlines are collinear, in which
case they can harmlessly slide over one another anyway because
the vanishing cycles can be made disjoint.} which gives great
flexibility in constructing and modifying Lagrangian torus
fibrations. But in dimension 3, this is no longer true: in a
generic 1-parameter family of Lagrangian torus fibrations, you
expect a discrete set of fibrations where the discriminant locus
fails to be embedded. It is not clear how to continue the nodal
slide beyond this point.

\begin{Problem}
Is there a theory of nodal slides in dimension 3? If the
discriminant locus hits itself at some point along the slide,
can we modify the fibration to allow the discriminant
locus\index{discriminant locus|)} to pass through itself?

\end{Problem}
See the work of Groman and Varolgunes \cite{GromanVarolgunes}
for a different and enlightening perspective on nodal slides.

\section{Lagrangian rational homology spheres}

Lagrangian submanifolds\index{Lagrangian!rational homology
sphere|(} \(L\) with \(H^1(L;\RR)=0\) are ``rigid'' in the sense
that any family of Lagrangian submanifolds \(L_t\) with
\(L_0=L\) come from a Hamiltonian isotopy in the sense that
\(L_t = \phi^{H_t}_t(L)\) for some time-dependent Hamiltonian
\(H_t\). If \(L\) is 3-dimensional then \(H^1(L;\RR)=0\) if and
only if \(L\) is a {\em rational homology sphere}, that is
\(H^*(L;\QQ) \cong H^*(S^3;\QQ)\).

In \(\cp{3}\), we have two well-known Lagrangian rational
homology spheres: \(\rp{3}\)
\index{Lagrangian!RP3@$\mathbb{RP}^3$} and the {\em Chiang
Lagrangian}\index{Lagrangian!Chiang} (see
\cite{Chiang,EvansLekili,JackPlatonic}).

\begin{Problem}\label{prb:chiang}
Are there any other Lagrangian rational homology 3-spheres in
\(\cp{3}\)?

\end{Problem}
Konstantinov \cite{Momchil} used Floer theory\index{Floer
theory|(} to give strong restrictions on which homology spheres
could occur. Namely, in {\cite[Corollary 3.2.12]{Momchil}}, he
showed that the only possibility is that \(L\) is a quotient of
\(S^3\) by a finite subgroup of \(SO(4)\) which is either cyclic
of order \(4k\) for some \(k\geq 1\) or else isomorphic to a
product \(D_{2^k(2n+1)}\times C_m\) where \(k\geq 2\), \(n\geq
1\), \(\gcd(2^k(2n+1),m)=1\) and \(D_\ell\) is the dihedral
group of order \(\ell\). The case \(D_{12}\times C_1\) is the
Chiang Lagrangian. \index{Lagrangian!rational homology sphere|)}

Ruling the other cases out would probably require Floer
theoretic machinery, but maybe one can construct examples using
tropical techniques\index{Lagrangian!submanifold!tropical}:
Mikhalkin's theory works for tropical curves in any dimension.
In particular, Mikhalkin has a construction of a Chiang-like
Lagrangian {\cite[Example 6.20]{Mikhalkin}} and he raises
Problem \ref{prb:chiang} as {\cite[Question 6.22]{Mikhalkin}}.

\section{Disjoint pinwheels}

Evans and Smith \cite{EvansSmith1} showed that one can find at
most three pairwise-disjointly embedded Lagrangian
pinwheels\index{Lagrangian!pinwheel|(} in \(\cp{2}\). The
argument is somewhat convoluted and mysterious; in particular,
we {\em first} showed that any triple of disjoint pinwheels
satisfy Markov's equation\index{Markov triple}, and then deduced
from this that there can be at most three. Since we do not know
what the analogue of the Markov equation should be for other Del
Pezzo surfaces like \(S^2\times S^2\), it is not clear how to
generalise the result.

Our result was inspired by a theorem of Hacking and Prokhorov
\cite{HackingProkhorov} which bounded the number of cyclic
quotient singularities\index{singularity!cyclic quotient} on a
singular degeneration of \(\cp{2}\). Prokhorov has a more
general bound for the number of cyclic quotient singularities in
a \(\QQ\)-Gorenstein degeneration\index{Q-Gorenstein
degeneration@$\mathbb{Q}$-Gorenstein degeneration} of a Del
Pezzo surface in terms of the Picard rank \(\rho\) of the
singular fibre (there should be at most \(\rho+2\)). This in
turn is bounded above by the rank of the quantum cohomology of
the smooth fibre, which motivates the following question:

\begin{Problem}\label{prb:pinwheels}
Let \(X\) be a monotone\index{monotone!symplectic manifold}
symplectic 4-manifold with semisimple quantum
cohomology\footnote{See \cite{McDuffSalamonJQH} for an
introduction to quantum cohomology, and the references below
for what it might have to do with disjoint Lagrangian
embeddings. Technically we should specify which ground field
we are working over for this assumption to make sense.} and
let \(r\) be the rank of its quantum cohomology. Is it true
that there are at most \(r\) pairwise-disjointly embedded
Lagrangian pinwheels in \(X\)?\index{Lagrangian!pinwheel|)}

\end{Problem}
Upper bounds like this hold for pairwise-disjointly embedded
Lagrangian submanifolds whose Floer cohomology is non-vanishing
over some field when the quantum cohomology is semisimple over
that field, see {\cite[Theorem 1.25]{EntovPolterovich}},
{\cite[Theorem 1.3]{SeidelDilation}}. However, Lagrangian
pinwheels do not have well-defined Floer cohomology, and even if
they did it would be likely that one must work over fields of
different characteristics to get well-defined Floer
cohomology\index{Floer theory|)} (even to get a fundamental
class). Perhaps one can make some headway using recent
developments in relative symplectic cohomology due to McLean
\cite{McLeanBirational}, Varolgunes \cite{Varolgunes}, and
Venkatesh \cite{Venkatesh}?

\section{Pinwheel content and deformations of \(\omega\)}

\begin{figure}[htb]
\begin{center}
\begin{tikzpicture}
\draw (0,2) -- (2,2) -- (2,1) -- (0,-3) -- cycle;
\draw[dashed] (2,2) -- (1.7,1.7) node {\(\times\)};
\draw[dashed] (0,-3) -- (0.6,1.8-3) node {\(\times\)};
\draw (0.6,1.8-3) -- (10/6,2);
\node[align=left] (lb) at (4,0) {visible\\ \((3,1)\)-pinwheel};
\draw[->] (lb.west) -- (0+1,0);
\node at (2,1.5) [right] {\(\lambda-1\)};
\node at (0,0) [left] {\(1+\lambda\)};
\node at (1,2) [above] {\(1\)};
\begin{scope}[shift={(6,0)}]
\draw (0,2) -- (2,2) -- (2,0) -- (0,-4) -- cycle;
\draw[dashed] (2,2) -- (1.7,1.7) node {\(\times\)};
\draw[dashed] (0,-4) -- (0.6,1.8-4) node {\(\times\)};
\draw (0.6,1.8-4) -- (2,2);
\node[align=left] (lb) at (4,-1) {visible \\ disc};
\draw[->] (lb.west) -- (0+1,-1);
\node at (2,1) [right] {\(\lambda=2\)};
\node at (0,-2) [left] {\(3\)};
\node at (1,2) [above] {\(1\)};
\end{scope}
\end{tikzpicture}
\caption{(a) A \((3,1)\)-pinwheel in \((S^2\times S^2,\omega_\lambda)\) for \(\lambda < 2\). (b) What becomes of this when \(\lambda\geq 2\).}
\label{fig:b_31_prob}
\end{center}
\end{figure}
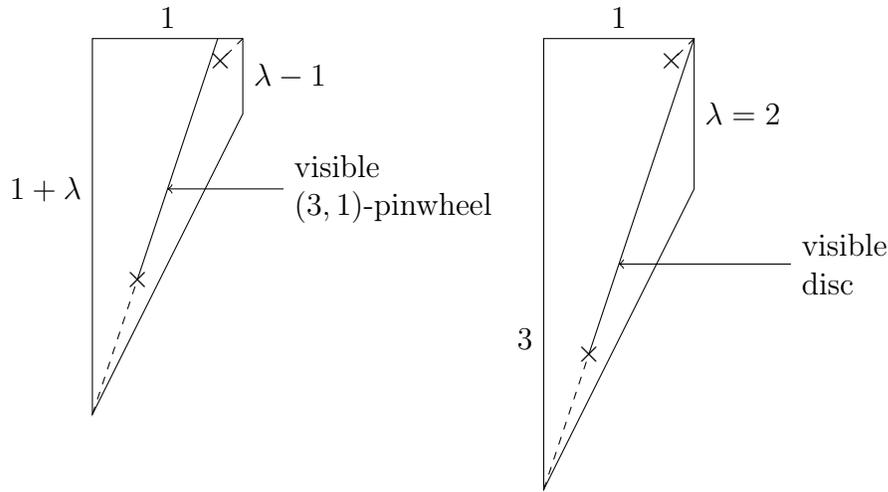

Given a symplectic 4-manifold \((X,\omega)\), define its {\em
pinwheel content} to be the set of pairs \((p,q)\) such that
some neighbourhood of the Lagrangian pinwheel in \(B_{1,p,q}\)
embeds symplectically in \(X\). Combining Problems
\ref{prb:mutation}, \ref{prb:genus} and \ref{prb:pinwheels}, we
can ask how the pinwheel content depends on the cohomology class
of \(\omega\). The simplest example is \(X=S^2\times S^2\) with
the symplectic form \(\omega_\lambda\) giving the factors areas
\(1\) and \(\lambda\). For \(\lambda<2\), it is possible to
construct a symplectically embedded \(B_{1,3,1}\subset X\): this
becomes visible after mutating the standard moment rectangle
(see Figure \ref{fig:b_31_prob}(a)). For \(\lambda\geq 2\), this
construction fails (Figure \ref{fig:b_31_prob}(b)).

\begin{Problem}
Is there a symplectically embedded \(B_{1,3,1}\subset
S^2\times S^2\) when \(\lambda\geq 2\)?

\end{Problem}
A much harder open question along these lines (about symplectic
embeddings of \(B_{1,p,q}\)s in surfaces of general type) was
raised in {\cite[Remark 1.4]{EvansUrzua}}.

\section{Big balls}

The final problem I want to state is motivated by a result of
Ein, K\"{u}chle and Lazarsfeld \cite{EKL} which gives a lower
bound on the Seshadri constants of projective
varieties. Rather than explaining this theorem, here is how
Lazarsfeld {\cite[Remark 5.2.7]{Lazarsfeld}} reformulates
their result in terms of symplectic topology:

\begin{Theorem}\label{thm:lazarsfeld}
Let \(X\subset\cp{N}\) be a smooth projective variety of
complex dimension \(n\) and let \(\omega\) be the symplectic
form on \(X\) pulled back from the Fubini-Study form on
\(\cp{N}\). Then \(X\) contains a symplectically embedded ball
of radius \(1/\sqrt{n\pi}\).

\end{Theorem}
All symplectic manifolds contain symplectic balls of some
(possibly very small) radius, but this result is saying you can
always find a ``big ball'' (of size at least \(1/\sqrt{n\pi}\))
if \(X\) is the symplectic manifold underlying a complex
\(n\)-dimensional projective variety.

The proof makes heavy use of techniques in algebraic geometry
which are not available in symplectic topology. However, one
might hope to find a purely symplectic construction of this big
ball. More generally:

\begin{Problem}[Big Ball conjecture]\label{prb:big_ball_conjecture}
Prove that there exists a universal constant \(R_n\) depending
only on \(n\) such that any integral\footnote{``Integral''
here means \([\omega]/2\pi\) is an integral cohomology
class. You might like to think about why this condition is
necessary.} symplectic manifold of dimension \(2n\) admits a
symplectically embedded ball of radius \(R_n\).

\end{Problem}
Almost toric fibrations provide a means of constructing balls,
for example taking the preimage of a suitable neighbourhood of a
Delzant vertex. Schlenk's book \cite{Schlenk} provides a wealth
of material on how to construct balls and packings by balls.

\bibliography{ltf}

\begin{thebibliography}{100}

\bibitem{WeinsteinHandlebodies}
B.~Acu, O.~Capovilla-Searle, A.~Gadbled, A.~Marinkovic, E.~Murphy,
  L.~Starkston, and A.~Wu.
\newblock An introduction to {W}einstein handlebodies for complements of
  smoothed toric divisors.
\newblock In {\em Research directions in symplectic and contact geometry and
  topology}, volume~27 of {\em Assoc. Women Math. Ser.}, pages 217--243.
  Springer, Cham, 2021.

\bibitem{Aigner}
M.~Aigner.
\newblock {\em Markov's theorem and 100 years of the uniqueness conjecture}.
\newblock Springer, Cham, 2013.
\newblock A mathematical journey from irrational numbers to perfect matchings.

\bibitem{Arnold}
V.~I. Arnol'd.
\newblock {\em Mathematical methods of classical mechanics}, volume~60 of {\em
  Graduate Texts in Mathematics}.
\newblock Springer-Verlag, New York, 1989.
\newblock Translated from the 1974 Russian original by K. Vogtmann and A.
  Weinstein, second edition.

\bibitem{Atiyah}
M.~F. Atiyah.
\newblock Convexity and commuting {H}amiltonians.
\newblock {\em Bull. London Math. Soc.}, 14(1):1--15, 1982.

\bibitem{Audin}
M.~Audin.
\newblock {\em Torus actions on symplectic manifolds}, volume~93 of {\em
  Progress in Mathematics}.
\newblock Birkh\"{a}user Verlag, Basel, revised edition, 2004.

\bibitem{AudinSkeletons}
M.~Audin.
\newblock Lagrangian skeletons, periodic geodesic flows and symplectic
  cuttings.
\newblock {\em Manuscripta Math.}, 124(4):533--550, 2007.

\bibitem{Auroux}
D.~Auroux.
\newblock Mirror symmetry and {$T$}-duality in the complement of an
  anticanonical divisor.
\newblock {\em J. G\"{o}kova Geom. Topol. GGT}, 1:51--91, 2007.

\bibitem{BertozziEtAl}
M.~Bertozzi, T.~Holm, E.~Maw, D.~McDuff, G.~T. Mwakyoma, A.~R. Pires, and
  M.~Weiler.
\newblock Infinite staircases for {H}irzebruch surfaces.
\newblock In {\em Research directions in symplectic and contact geometry and
  topology}, volume~27 of {\em Assoc. Women Math. Ser.}, pages 47--157.
  Springer, Cham, 2021.

\bibitem{BhupalOno}
M.~Bhupal and K.~Ono.
\newblock Symplectic fillings of links of quotient surface singularities.
\newblock {\em Nagoya Math. J.}, 207:1--45, 2012.

\bibitem{BottTu}
R.~Bott and L.~W. Tu.
\newblock {\em Differential forms in algebraic topology}, volume~82 of {\em
  Graduate Texts in Mathematics}.
\newblock Springer-Verlag, New York-Berlin, 1982.

\bibitem{CasalsVianna}
R.~Casals and R.~Vianna.
\newblock Full ellipsoid embeddings and toric mutations.
\newblock {\em Selecta Math. (N.S.)}, 28(3):Paper No. 61, 62, 2022.

\bibitem{CasselsDiophantine}
J.~W.~S. Cassels.
\newblock {\em An introduction to {D}iophantine approximation}.
\newblock Cambridge Tracts in Mathematics and Mathematical Physics, No. 45.
  Cambridge University Press, New York, 1957.

\bibitem{CastanoBernardMatessi2}
R.~Casta\~{n}o Bernard and D.~Matessi.
\newblock Lagrangian 3-torus fibrations.
\newblock {\em J. Differential Geom.}, 81(3):483--573, 2009.

\bibitem{CastanoBernardMatessi1}
R.~Casta\~{n}o Bernard and D.~Matessi.
\newblock Semi-global invariants of piecewise smooth {L}agrangian fibrations.
\newblock {\em Q. J. Math.}, 61(3):291--318, 2010.

\bibitem{Chaperon}
M.~Chaperon.
\newblock Normalisation of the smooth focus-focus: a simple proof.
\newblock {\em Acta Math. Vietnam.}, 38(1):3--9, 2013.

\bibitem{Chekanov}
Yu.~V. Chekanov.
\newblock Lagrangian tori in a symplectic vector space and global
  symplectomorphisms.
\newblock {\em Math. Z.}, 223(4):547--559, 1996.

\bibitem{CheungVianna}
M.-W.~M. Cheung and R.~Vianna.
\newblock Algebraic and symplectic viewpoint on compactifications of
  two-dimensional cluster varieties of finite type.
\newblock In D.~R. Wood, J.~de~Gier, C.~E. Praeger, and T.~Tao, editors, {\em
  2019-20 MATRIX Annals}, volume~4 of {\em MATRIX Book Series}. Springer, 2021.

\bibitem{Chiang}
R.~Chiang.
\newblock New {L}agrangian submanifolds of {$\Bbb{CP}^n$}.
\newblock {\em Int. Math. Res. Not.}, (45):2437--2441, 2004.

\bibitem{Christophersen}
J.~A. Christophersen.
\newblock On the components and discriminant of the versal base space of cyclic
  quotient singularities.
\newblock In {\em Singularity theory and its applications, {P}art {I}
  ({C}oventry, 1988/1989)}, volume 1462 of {\em Lecture Notes in Math.}, pages
  81--92. Springer, Berlin, 1991.

\bibitem{Cicero}
M.~T. Cicero.
\newblock {\em Tusculanae {D}isputationes V}.
\newblock 45 BC.

\bibitem{Conway}
J.~H. Conway.
\newblock {\em The sensual (quadratic) form}, volume~26 of {\em Carus
  Mathematical Monographs}.
\newblock Mathematical Association of America, Washington, DC, 1997.
\newblock With the assistance of Francis Y. C. Fung.

\bibitem{CGHMP}
D.~Cristofaro-Gardiner, T.~Holm, A.~Mandini, and A.~R. Pires.
\newblock On infinite staircases in toric symplectic four-manifolds.
\newblock {\em ar{X}iv:2004.07829}, 2020.

\bibitem{Danilov}
V.~I. Danilov.
\newblock The geometry of toric varieties.
\newblock {\em Uspekhi Mat. Nauk}, 33(2(200)):85--134, 247, 1978.

\bibitem{Delzant}
T.~Delzant.
\newblock Hamiltoniens p\'{e}riodiques et images convexes de l'application
  moment.
\newblock {\em Bull. Soc. Math. France}, 116(3):315--339, 1988.

\bibitem{DufourMolino}
J.-P. Dufour and P.~Molino.
\newblock Compactification d'actions de {${\bf R}^n$} et variables action-angle
  avec singularit\'{e}s.
\newblock In {\em Travaux du {S}\'{e}minaire {S}ud-{R}hodanien de
  {G}\'{e}om\'{e}trie, {I}}, volume~88 of {\em Publ. D\'{e}p. Math. Nouvelle
  S\'{e}r. B}, pages 161--183. Univ. Claude-Bernard, Lyon, 1988.

\bibitem{Duistermaat}
J.~J. Duistermaat.
\newblock On global action-angle coordinates.
\newblock {\em Comm. Pure Appl. Math.}, 33(6):687--706, 1980.

\bibitem{EbertMO}
J.~Ebert.
\newblock Relative de rham cohomologies.
\newblock MathOverflow.
\newblock URL:https://mathoverflow.net/q/111063 (version: 2012-10-30).

\bibitem{EKL}
L.~Ein, O.~K\"{u}chle, and R.~Lazarsfeld.
\newblock Local positivity of ample line bundles.
\newblock {\em J. Differential Geom.}, 42(2):193--219, 1995.

\bibitem{Eliashberg}
Y.~Eliashberg.
\newblock Contact {$3$}-manifolds twenty years since {J}. {M}artinet's work.
\newblock {\em Ann. Inst. Fourier (Grenoble)}, 42(1-2):165--192, 1992.

\bibitem{Eliasson}
L.~H. Eliasson.
\newblock Normal forms for {H}amiltonian systems with {P}oisson commuting
  integrals---elliptic case.
\newblock {\em Comment. Math. Helv.}, 65(1):4--35, 1990.

\bibitem{Engel}
P.~Engel.
\newblock Looijenga's conjecture via integral-affine geometry.
\newblock {\em J. Differential Geom.}, 109(3):467--495, 2018.

\bibitem{EntovPolterovich}
M.~Entov and L.~Polterovich.
\newblock Rigid subsets of symplectic manifolds.
\newblock {\em Compos. Math.}, 145(3):773--826, 2009.

\bibitem{EvansKB}
J.~D. Evans.
\newblock A {L}agrangian {K}lein bottle you can't squeeze.
\newblock {\em J. Fixed Point Theory Appl.}, 24(2):Paper No. 47, 10, 2022.

\bibitem{EvansLekili}
J.~D. Evans and Y.~Lekili.
\newblock Floer cohomology of the {C}hiang {L}agrangian.
\newblock {\em Selecta Math. (N.S.)}, 21(4):1361--1404, 2015.

\bibitem{EvansMauri}
J.~D. Evans and M.~Mauri.
\newblock Constructing local models for {L}agrangian torus fibrations.
\newblock {\em Ann. H. Lebesgue}, 4:537--570, 2021.

\bibitem{EvansSmith1}
J.~D. Evans and I.~Smith.
\newblock Markov numbers and {L}agrangian cell complexes in the complex
  projective plane.
\newblock {\em Geom. Topol.}, 22(2):1143--1180, 2018.

\bibitem{EvansSmith2}
J.~D. Evans and I.~Smith.
\newblock Bounds on {W}ahl singularities from symplectic topology.
\newblock {\em Algebr. Geom.}, 7(1):59--85, 2020.

\bibitem{EvansUrzua}
J.~D. Evans and G.~Urz\'{u}a.
\newblock Antiflips, mutations, and unbounded symplectic embeddings of rational
  homology balls.
\newblock {\em Ann. Inst. Fourier (Grenoble)}, 71(5):1807--1843, 2021.

\bibitem{FintushelStern}
R.~Fintushel and R.~J. Stern.
\newblock Rational blowdowns of smooth {$4$}-manifolds.
\newblock {\em J. Differential Geom.}, 46(2):181--235, 1997.

\bibitem{Fulton}
W.~Fulton.
\newblock {\em Introduction to toric varieties}, volume 131 of {\em Annals of
  Mathematics Studies}.
\newblock Princeton University Press, Princeton, NJ, 1993.
\newblock The William H. Roever Lectures in Geometry.

\bibitem{Gompf}
R.~E. Gompf.
\newblock A new construction of symplectic manifolds.
\newblock {\em Ann. of Math. (2)}, 142(3):527--595, 1995.

\bibitem{Greenhill}
A.~G. Greenhill.
\newblock {\em The applications of elliptic functions}.
\newblock Dover Publications, Inc., New York, 1959.

\bibitem{GromanVarolgunes}
Y.~Groman and U.~Varolgunes.
\newblock Locality of relative symplectic cohomology for complete embeddings.
\newblock {\em arXiv:2110.08891}, 2021.

\bibitem{Gromov}
M.~Gromov.
\newblock Pseudo holomorphic curves in symplectic manifolds.
\newblock {\em Invent. Math.}, 82(2):307--347, 1985.

\bibitem{Gross3}
M.~Gross.
\newblock Examples of special {L}agrangian fibrations.
\newblock In {\em Symplectic geometry and mirror symmetry ({S}eoul, 2000)},
  pages 81--109. World Sci. Publ., River Edge, NJ, 2001.

\bibitem{Gross1}
M.~Gross.
\newblock Special {L}agrangian fibrations. {I}. {T}opology.
\newblock In {\em Winter {S}chool on {M}irror {S}ymmetry, {V}ector {B}undles
  and {L}agrangian {S}ubmanifolds ({C}ambridge, {MA}, 1999)}, volume~23 of {\em
  AMS/IP Stud. Adv. Math.}, pages 65--93. Amer. Math. Soc., Providence, RI,
  2001.

\bibitem{Gross2}
M.~Gross.
\newblock Special {L}agrangian fibrations. {II}. {G}eometry. {A} survey of
  techniques in the study of special {L}agrangian fibrations.
\newblock In {\em Winter {S}chool on {M}irror {S}ymmetry, {V}ector {B}undles
  and {L}agrangian {S}ubmanifolds ({C}ambridge, {MA}, 1999)}, volume~23 of {\em
  AMS/IP Stud. Adv. Math.}, pages 95--150. Amer. Math. Soc., Providence, RI,
  2001.

\bibitem{GHK}
M.~Gross, P.~Hacking, and S.~Keel.
\newblock Mirror symmetry for log {C}alabi-{Y}au surfaces {I}.
\newblock {\em Publ. Math. Inst. Hautes \'{E}tudes Sci.}, 122:65--168, 2015.

\bibitem{GS}
V.~Guillemin and S.~Sternberg.
\newblock Convexity properties of the moment mapping.
\newblock {\em Invent. Math.}, 67(3):491--513, 1982.

\bibitem{GuilleminSternbergBirational}
V.~Guillemin and S.~Sternberg.
\newblock Birational equivalence in the symplectic category.
\newblock {\em Invent. Math.}, 97(3):485--522, 1989.

\bibitem{HackingProkhorov}
P.~Hacking and Y.~Prokhorov.
\newblock Smoothable del {P}ezzo surfaces with quotient singularities.
\newblock {\em Compos. Math.}, 146(1):169--192, 2010.

\bibitem{HTU}
P.~Hacking, J.~Tevelev, and G.~Urz\'{u}a.
\newblock Flipping surfaces.
\newblock {\em J. Algebraic Geom.}, 26(2):279--345, 2017.

\bibitem{Halmos}
P.~R. Halmos.
\newblock {\em Finite {D}imensional {V}ector {S}paces}.
\newblock Annals of Mathematics Studies, No. 7. Princeton University Press,
  Princeton, N.J., 1942.

\bibitem{HirzebruchHilbert}
F.~E.~P. Hirzebruch.
\newblock Hilbert modular surfaces.
\newblock {\em Enseign. Math. (2)}, 19:183--281, 1973.

\bibitem{JankinsNeumann}
M.~Jankins and W.~D. Neumann.
\newblock {\em Lectures on {S}eifert manifolds}, volume~2 of {\em Brandeis
  Lecture Notes}.
\newblock Brandeis University, Waltham, MA, 1983.

\bibitem{JoyceSpLagSYZ}
D.~Joyce.
\newblock Singularities of special {L}agrangian fibrations and the {SYZ}
  conjecture.
\newblock {\em Comm. Anal. Geom.}, 11(5):859--907, 2003.

\bibitem{Karabas}
D.~Karabas.
\newblock Microlocal sheaves on pinwheels.
\newblock {\em ar{X}iv:1810.09021}, page 109, 2018.

\bibitem{Khodorovskiy}
T.~Khodorovskiy.
\newblock Bounds on embeddings of rational homology balls in symplectic
  4-manifolds.
\newblock {\em ar{X}iv:1307.4321}, 2013.

\bibitem{KSB}
J.~Koll\'{a}r and N.~I. Shepherd-Barron.
\newblock Threefolds and deformations of surface singularities.
\newblock {\em Invent. Math.}, 91(2):299--338, 1988.

\bibitem{Momchil}
M.~Konstantinov.
\newblock {\em Symplectic Topology of Projective Space: {L}agrangians, Local
  Systems and Twistors}.
\newblock PhD thesis, University College London, 2019.

\bibitem{KS}
M.~Kontsevich and Y.~Soibelman.
\newblock Homological mirror symmetry and torus fibrations.
\newblock In {\em Symplectic geometry and mirror symmetry ({S}eoul, 2000)},
  pages 203--263. World Sci. Publ., River Edge, NJ, 2001.

\bibitem{Lazarsfeld}
R.~Lazarsfeld.
\newblock {\em Positivity in algebraic geometry. {I}}, volume~48 of {\em
  Ergebnisse der Mathematik und ihrer Grenzgebiete. 3. Folge. A Series of
  Modern Surveys in Mathematics [Results in Mathematics and Related Areas. 3rd
  Series. A Series of Modern Surveys in Mathematics]}.
\newblock Springer-Verlag, Berlin, 2004.
\newblock Classical setting: line bundles and linear series.

\bibitem{Lee}
J.~M. Lee.
\newblock {\em Introduction to smooth manifolds}, volume 218 of {\em Graduate
  Texts in Mathematics}.
\newblock Springer, New York, second edition, 2013.

\bibitem{LekiliMaydanskiy}
Y.~Lekili and M.~Maydanskiy.
\newblock The symplectic topology of some rational homology balls.
\newblock {\em Comment. Math. Helv.}, 89(3):571--596, 2014.

\bibitem{Lerman}
E.~Lerman.
\newblock Symplectic cuts.
\newblock {\em Math. Res. Lett.}, 2(3):247--258, 1995.

\bibitem{SymingtonLeung}
N.~C. Leung and M.~Symington.
\newblock Almost toric symplectic four-manifolds.
\newblock {\em J. Symplectic Geom.}, 8(2):143--187, 2010.

\bibitem{TJLiRuan}
T.-J. Li and Y.~Ruan.
\newblock Symplectic birational geometry.
\newblock In {\em New perspectives and challenges in symplectic field theory},
  volume~49 of {\em CRM Proc. Lecture Notes}, pages 307--326. Amer. Math. Soc.,
  Providence, RI, 2009.

\bibitem{Lisca}
P.~Lisca.
\newblock On symplectic fillings of lens spaces.
\newblock {\em Trans. Amer. Math. Soc.}, 360(2):765--799, 2008.

\bibitem{LooijengaWahl}
E.~Looijenga and J.~Wahl.
\newblock Quadratic functions and smoothing surface singularities.
\newblock {\em Topology}, 25(3):261--291, 1986.

\bibitem{MagillMcDuff}
N.~Magill and D.~McDuff.
\newblock Staircase symmetries in {H}irzebruch surfaces.
\newblock {\em ar{X}iv:2106.09143}, 2021.

\bibitem{Manetti}
M.~Manetti.
\newblock Normal degenerations of the complex projective plane.
\newblock {\em J. Reine Angew. Math.}, 419:89--118, 1991.

\bibitem{Matessi1}
D.~Matessi.
\newblock Lagrangian pairs of pants.
\newblock {\em Int. Math. Res. Not. IMRN}, (15):11306--11356, 2021.

\bibitem{Matessi2}
D.~Matessi.
\newblock Lagrangian submanifolds from tropical hypersurfaces.
\newblock {\em Internat. J. Math}, 32(7):Paper No. 2150046, 63, 2021.

\bibitem{McDuffRatRuled}
D.~McDuff.
\newblock The structure of rational and ruled symplectic {$4$}-manifolds.
\newblock {\em J. Amer. Math. Soc.}, 3(3):679--712, 1990.

\bibitem{McDuffPolterovich}
D.~McDuff and L.~Polterovich.
\newblock Symplectic packings and algebraic geometry.
\newblock {\em Invent. Math.}, 115(3):405--434, 1994.
\newblock With an appendix by Yael Karshon.

\bibitem{McDuffSalamonJQH}
D.~McDuff and D.~Salamon.
\newblock {\em {$J$}-holomorphic curves and quantum cohomology}, volume~6 of
  {\em University Lecture Series}.
\newblock American Mathematical Society, Providence, RI, 1994.

\bibitem{McDuffSalamon}
D.~McDuff and D.~Salamon.
\newblock {\em {Introduction to symplectic topology. 2nd ed}}.
\newblock New York, NY: Oxford University Press, 2nd ed. edition, 1998.

\bibitem{McDuffSchlenk}
D.~McDuff and F.~Schlenk.
\newblock The embedding capacity of 4-dimensional symplectic ellipsoids.
\newblock {\em Ann. of Math. (2)}, 175(3):1191--1282, 2012.

\bibitem{McLeanBirational}
M.~McLean.
\newblock Birational {C}alabi-{Y}au manifolds have the same small quantum
  products.
\newblock {\em Ann. of Math. (2)}, 191(2):439--579, 2020.

\bibitem{Mikhalkin}
G.~Mikhalkin.
\newblock {Examples of tropical-to-Lagrangian correspondence}.
\newblock {\em {Eur. J. Math.}}, 5(3):1033--1066, 2019.

\bibitem{MilnorSingularities}
J.~Milnor.
\newblock {\em Singular points of complex hypersurfaces}.
\newblock Annals of Mathematics Studies, No. 61. Princeton University Press,
  Princeton, N.J.; University of Tokyo Press, Tokyo, 1968.

\bibitem{MoserVolume}
J.~Moser.
\newblock On the volume elements on a manifold.
\newblock {\em Trans. Amer. Math. Soc.}, 120:286--294, 1965.

\bibitem{OrlikSeifert}
P.~Orlik.
\newblock {\em Seifert manifolds}.
\newblock Lecture Notes in Mathematics, Vol. 291. Springer-Verlag, Berlin-New
  York, 1972.

\bibitem{ParkQBD}
J.~Park.
\newblock Simply connected symplectic 4-manifolds with {$b^+_2=1$} and
  {$c^2_1=2$}.
\newblock {\em Invent. Math.}, 159(3):657--667, 2005.

\bibitem{TonkonogPascaleff}
J.~Pascaleff and D.~Tonkonog.
\newblock The wall-crossing formula and {L}agrangian mutations.
\newblock {\em Adv. Math.}, 361:106850, 67, 2020.

\bibitem{RanaUrzua}
J.~Rana and G.~Urz\'{u}a.
\newblock Optimal bounds for {T}-singularities in stable surfaces.
\newblock {\em Adv. Math.}, 345:814--844, 2019.

\bibitem{Ruan1}
W.-D. Ruan.
\newblock Lagrangian torus fibration of quintic hypersurfaces. {I}. {F}ermat
  quintic case.
\newblock In {\em Winter {S}chool on {M}irror {S}ymmetry, {V}ector {B}undles
  and {L}agrangian {S}ubmanifolds ({C}ambridge, {MA}, 1999)}, volume~23 of {\em
  AMS/IP Stud. Adv. Math.}, pages 297--332. Amer. Math. Soc., Providence, RI,
  2001.

\bibitem{Ruan2}
W.-D. Ruan.
\newblock Lagrangian torus fibration of quintic {C}alabi-{Y}au hypersurfaces.
  {II}. {T}echnical results on gradient flow construction.
\newblock {\em J. Symplectic Geom.}, 1(3):435--521, 2002.

\bibitem{Ruan3}
W.-D. Ruan.
\newblock Lagrangian torus fibration of quintic {C}alabi-{Y}au hypersurfaces.
  {III}. {S}ymplectic topological {SYZ} mirror construction for general
  quintics.
\newblock {\em J. Differential Geom.}, 63(2):171--229, 2003.

\bibitem{Rudakov}
A.~N. Rudakov.
\newblock Markov numbers and exceptional bundles on {${\bf P}^2$}.
\newblock {\em Izv. Akad. Nauk SSSR Ser. Mat.}, 52(1):100--112, 240, 1988.

\bibitem{Schlenk}
F.~Schlenk.
\newblock {\em Embedding problems in symplectic geometry}, volume~40 of {\em De
  Gruyter Expositions in Mathematics}.
\newblock Walter de Gruyter GmbH \& Co. KG, Berlin, 2005.

\bibitem{SchoenWolfson}
R.~Schoen and J.~Wolfson.
\newblock Minimizing area among {L}agrangian surfaces: the mapping problem.
\newblock {\em J. Differential Geom.}, 58(1):1--86, 2001.

\bibitem{SeidelDilation}
P.~Seidel.
\newblock Disjoinable {L}agrangian spheres and dilations.
\newblock {\em Invent. Math.}, 197(2):299--359, 2014.

\bibitem{SeifertThrelfall}
H.~Seifert and W.~Threlfall.
\newblock {\em Seifert and {T}hrelfall: a textbook of topology}, volume~89 of
  {\em Pure and Applied Mathematics}.
\newblock Academic Press, Inc. [Harcourt Brace Jovanovich, Publishers], New
  York-London, 1980.
\newblock Translated from the German edition of 1934 by Michael A. Goldman,
  With a preface by Joan S. Birman, With ``Topology of $3$-dimensional fibered
  spaces'' by Seifert, Translated from the German by Wolfgang Heil.

\bibitem{SepeNgoc}
D.~Sepe and S.~V\~{u}~Ng\d{o}c.
\newblock Integrable systems, symmetries, and quantization.
\newblock {\em Lett. Math. Phys.}, 108(3):499--571, 2018.

\bibitem{SeriesMarkoff}
C.~Series.
\newblock The geometry of {M}arkoff numbers.
\newblock {\em Math. Intelligencer}, 7(3):20--29, 1985.

\bibitem{ExactCombinatorics}
V.~Shende, D.~Treumann, and H.~Williams.
\newblock On the combinatorics of exact {L}agrangian surfaces.
\newblock {\em ar{X}iv:1603.07449}, 2016.

\bibitem{ShevSmi}
V.~Shevchishin and G.~Smirnov.
\newblock Symplectic triangle inequality.
\newblock {\em Proc. Amer. Math. Soc.}, 148(4):1389--1397, 2020.

\bibitem{ShevNodal}
V.~V. Shevchishin.
\newblock On the local {S}everi problem.
\newblock {\em Int. Math. Res. Not.}, (5):211--237, 2004.

\bibitem{SiebertTian1}
B.~Siebert and G.~Tian.
\newblock On the holomorphicity of genus two {L}efschetz fibrations.
\newblock {\em Ann. of Math. (2)}, 161(2):959--1020, 2005.

\bibitem{SiebertTian2}
B.~Siebert and G.~Tian.
\newblock Lectures on pseudo-holomorphic curves and the symplectic isotopy
  problem.
\newblock In {\em Symplectic 4-manifolds and algebraic surfaces}, volume 1938
  of {\em Lecture Notes in Math.}, pages 269--341. Springer, Berlin, 2008.

\bibitem{Sikorav}
J.-C. Sikorav.
\newblock The gluing construction for normally generic {$J$}-holomorphic
  curves.
\newblock In {\em Symplectic and contact topology: interactions and
  perspectives ({T}oronto, {ON}/{M}ontreal, {QC}, 2001)}, volume~35 of {\em
  Fields Inst. Commun.}, pages 175--199. Amer. Math. Soc., Providence, RI,
  2003.

\bibitem{JackPlatonic}
J.~Smith.
\newblock Floer cohomology of {P}latonic {L}agrangians.
\newblock {\em J. Symplectic Geom.}, 17(2):477--601, 2019.

\bibitem{Stevens}
J.~Stevens.
\newblock The versal deformation of cyclic quotient singularities.
\newblock In {\em Deformations of surface singularities}, volume~23 of {\em
  Bolyai Soc. Math. Stud.}, pages 163--201. J\'{a}nos Bolyai Math. Soc.,
  Budapest, 2013.

\bibitem{Symington2}
M.~Symington.
\newblock Generalized symplectic rational blowdowns.
\newblock {\em Algebr. Geom. Topol.}, 1:503--518, 2001.

\bibitem{Symington}
M.~Symington.
\newblock Four dimensions from two in symplectic topology.
\newblock In {\em Topology and geometry of manifolds ({A}thens, {GA}, 2001)},
  volume~71 of {\em Proc. Sympos. Pure Math.}, pages 153--208. Amer. Math.
  Soc., Providence, RI, 2003.

\bibitem{Thurston}
W.~P. Thurston.
\newblock Some simple examples of symplectic manifolds.
\newblock {\em Proc. Amer. Math. Soc.}, 55(2):467--468, 1976.

\bibitem{UrzuaNeighbours}
G.~Urz\'{u}a.
\newblock Identifying neighbors of stable surfaces.
\newblock {\em Ann. Sc. Norm. Super. Pisa Cl. Sci. (5)}, 16(4):1093--1122,
  2016.

\bibitem{UrzuaVilches}
G.~Urz\'{u}a and N.~Vilches.
\newblock On wormholes in the moduli space of surfaces.
\newblock {\em Algebr. Geom.}, 9(1):39--68, 2022.

\bibitem{UsherMO}
M.~Usher.
\newblock Symplectic blow-up.
\newblock MathOverflow.
\newblock https://mathoverflow.net/q/133709 (version: 2013-06-14).

\bibitem{Ngoc}
S.~V\~{u}~Ng\d{o}c.
\newblock On semi-global invariants for focus-focus singularities.
\newblock {\em Topology}, 42(2):365--380, 2003.

\bibitem{NgocBohrSommerfeld}
San V\~{u}~Ng\d{o}c.
\newblock Bohr-{S}ommerfeld conditions for integrable systems with critical
  manifolds of focus-focus type.
\newblock {\em Comm. Pure Appl. Math.}, 53(2):143--217, 2000.

\bibitem{Varolgunes}
U.~Varolgunes.
\newblock Mayer-{V}ietoris property for relative symplectic cohomology.
\newblock {\em Geom. Topol.}, 25(2):547--642, 2021.

\bibitem{Venkatesh}
S.~Venkatesh.
\newblock Rabinowitz {F}loer homology and mirror symmetry.
\newblock {\em J. Topol.}, 11(1):144--179, 2018.

\bibitem{Vianna1}
R.~Vianna.
\newblock {On exotic Lagrangian tori in \(\mathbb{CP}^{2}\)}.
\newblock {\em {Geom. Topol.}}, 18(4):2419--2476, 2014.

\bibitem{Vianna2}
R.~Vianna.
\newblock {Infinitely many exotic monotone Lagrangian tori in
  \(\mathbb{C}\mathbb{P}^{2}\)}.
\newblock {\em {J. Topol.}}, 9(2):535--551, 2016.

\bibitem{Vianna3}
R.~Vianna.
\newblock {Infinitely many monotone Lagrangian tori in del Pezzo surfaces}.
\newblock {\em {Sel. Math., New Ser.}}, 23(3):1955--1996, 2017.

\bibitem{Warner}
F.~W. Warner.
\newblock {\em Foundations of differentiable manifolds and {L}ie groups},
  volume~94 of {\em Graduate Texts in Mathematics}.
\newblock Springer-Verlag, New York-Berlin, 1983.
\newblock Corrected reprint of the 1971 edition.

\bibitem{Whittaker}
E.~T. Whittaker.
\newblock {\em A treatise on the analytical dynamics of particles and rigid
  bodies: {W}ith an introduction to the problem of three bodies}.
\newblock Cambridge University Press, New York, 1959.
\newblock 4th ed.

\bibitem{Zung1}
N.~T. Zung.
\newblock Symplectic topology of integrable {H}amiltonian systems. {I}.
  {A}rnold-{L}iouville with singularities.
\newblock {\em Compositio Math.}, 101(2):179--215, 1996.

\bibitem{Zung2}
N.~T. Zung.
\newblock Symplectic topology of integrable {H}amiltonian systems. {II}.
  {T}opological classification.
\newblock {\em Compositio Math.}, 138(2):125--156, 2003.

\end{thebibliography}
\bibliographystyle{plain}

\cleardoublepage
\addcontentsline{toc}{chapter}{Index}
\printindex
\end{document}